\documentclass[12pt]{amsart}
\usepackage[a4paper]{anysize}
\usepackage{amsthm}
\usepackage{amssymb}
\usepackage{times}
\usepackage{mathrsfs}

\input xy
\xyoption{v2}
\xyoption{all}

    \newtheorem{theorem}{Theorem}[subsection]
    \newtheorem{proposition}[theorem]{Proposition}
    \newtheorem{lemma}[theorem]{Lemma}
    \newtheorem{corollary}[theorem]{Corollary}
    
    \theoremstyle{definition}
    \newtheorem{example}[theorem]{Example}
    \newtheorem{definition}[theorem]{Definition}
    \newtheorem{remark}[theorem]{Remark}
    \theoremstyle{remark}
    \newtheorem{claim}[theorem]{Claim}

    \def\set{\setcounter{equation}
             {\value{theorem}}\addtocounter{theorem}{1}}
    \numberwithin{equation}{subsection}
    \def\sset{\setcounter{subsubsection}
             {\value{theorem}}\addtocounter{theorem}{1}}

\def\A{{\mathbb A}}
\def\cA{{{\mathscr A}}}

\def\B{{\mathbb B}}
\def\cB{{{\mathscr B}}}

\def\cC{{{\mathscr C}}}  
\def\sC{{\text{\sf C}}}
\def\cD{{{\mathscr D}}}

\def\fD{{\mathfrak D}}
\def\sD{{\text{\sf D}}}

\def\cE{{{\mathscr E}}}
\def\F{{\mathbb F}}
\def\cF{{{\mathscr F}}}
\def\G{{\mathbb G}}
\def\cG{{{\mathscr G}}}

\def\cI{{{\mathscr I}}}
\def\cJ{{{\mathscr J}}}

\def\cL{{{\mathscr L}}}
\def\L{{\mathbb L}}
\def\cM{{{\mathscr M}}}

\def\N{{\mathbb N}}
\def\cO{{{\mathscr O}}}

\def\cP{{{\mathscr P}}}
\def\fp{{\mathfrak p}}
\def\P{{\mathbb P}}
\def\bP{{\mathbf P}}
\def\cQ{{{\mathscr Q}}}
\def\fq{{\mathfrak q}}
\def\Q{{\mathbb Q}}
\def\R{{\mathbb R}}
\def\cS{{{\mathscr S}}}
\def\cT{{{\mathscr T}}}

\def\Z{{\mathbb Z}}

\def\cW{{{\mathscr W}}}

\def\cU{{{\mathscr U}}}

\newcommand{\indref}[2]{\ref{#1}}

\newcommand{\limdir}[1]
{\underset{#1}{\lim}}

\newcommand{\colim}[1]
{\underset{#1}{\mathrm{colim}\,}}

\newcommand{\liminv}[1]
{\underset{#1}{\lim}}

\newcommand{\derotimes}
{\overset{\mathbf{L}}{\otimes}}

\renewcommand{\labelenumi}{(\roman{enumi})}

\newenvironment{pfclaim}
{\noindent {\em Proof of the claim:\/}}
{{\smallskip}}

\newenvironment{acknowledgement}{\vskip .5cm
\setlength{\baselineskip}{0mm}
\noindent \footnotesize\rmfamily
{\em Acknowledgements\ \ \ }}

\def\Tr{\mathrm{Tr}}
\def\tr{\mathrm{tr}}
\def\gr{\mathrm{gr}}

\def\Ann{\mathrm{Ann}}

\def\flip{\mathrm{flip}}

\def\Cone{\mathrm{Cone}}
\def\Gal{\mathrm{Gal}}
\def\Spec{\mathrm{Spec}}
\def\Alhom{\mathrm{alHom}}
\def\Hom{\mathrm{Hom}}
\def\Hot{\mathrm{Hot}}
\def\Der{\mathrm{Der}}
\def\AlExt{\mathrm{alExt}}
\def\Ext{\mathrm{Ext}}
\def\hExt{{\mathbb{E}\mbox{\sf xt}}}
\def\Exal{\mathrm{Exal}}
\def\Exmon{\mathrm{Exmon}}

\def\bExun{\mathbf{Exun}}
\def\bExal{\mathbf{Exal}}
\def\bExmon{\mathbf{Exmon}}
\def\Tor{\mathrm{Tor}}
\def\End{\mathrm{End}}
\def\Idemp{\mathrm{Idemp}}

\def\Spf{\mathrm{Spf}}
\def\Coker{\mathrm{Coker}}
\def\Ker{\mathrm{Ker}}
\def\Img{\mathrm{Im}}
\def\Aut{\mathrm{Aut}}
\def\Alg{{{\text{-}}\mathbf{Alg}}}
\def\Sym{\mathrm{Sym}}
\def\AlgMorph{{\text{-}\mathbf{Alg.Morph}}}
\def\AlgMod{{\text{-}\mathbf{Alg.Mod}}}
\def\Aliso{{\text{-}\mathbf{al.Iso}}}
\def\Desc{\mathbf{Desc}}
\def\Et{\text{-}\mathbf{\acute{E}t}}
\def\fibEt{\mathbf{\acute{E}t}}
\def\uEt{\text{-}\mathbf{u.\acute{E}t}}
\def\fibuEt{\mathbf{u.\acute{E}t}}
\def\wEt{\text{-}\mathbf{w.\acute{E}t}}
\def\fibwEt{\mathbf{w.\acute{E}t}}
\def\Mod{\text{-}\mathbf{Mod}}
\def\Mon{\text{-}\mathbf{Mon}}
\def\UniMod{\text{-}\mathbf{Uni.Mod}}
\def\Set{\mathbf{Set}}
\def\chara{\mathrm{char}}
\def\nil{\mathrm{nil}}
\def\pr{\mathrm{pr}}
\def\rk{\mathrm{rk}}
\def\one{\mathbf{1}}
\def\bar#1{\overline{#1}}
\def\hat{\widehat} \def\tilde{\widetilde} 
\def\fm{\mathfrak m}
\def\fp{\mathfrak p}
\def\eps{\varepsilon}

\def\ev{\mathrm{ev}}
\def\frk{\mathrm{f.rk}}
\def\Split{\mathrm{Split}}
\def\Iperv{I_{!*}}
\def\thetaperv{$\theta_{M|N}$}
\def\Mperv{M_!}
\def\Bperv{B_{!!}}

\def\bmu{\boldsymbol\mu}
\def\prelog{\mathbf{pre}\text{-}\mathbf{log}}

\def\fG{\mathfrak{G}}
\def\fH{\mathfrak{H}}
\def\sep{\mathrm{s}}
\def\tame{\mathrm{t}}

\def\fSet{\mathbf{f.Set}}
\renewcommand\emptyset{\varnothing}
\def\rad{\mathrm{rad}}

\font\prelim=pagko

\makeindex

\begin{document}
\title{Almost ring theory}

\author{Ofer Gabber}
\author{Lorenzo Ramero}

\maketitle

\centerline{\prelim sixth (and final) release}

\vskip 1cm


\vskip 2cm

    Ofer Gabber

    I.H.E.S.

    Le Bois-Marie

    35, route de Chartres

    F-91440 Bures-sur-Yvette

    {\em e-mail address:} {\ttfamily gabber@ihes.fr}

\vskip 1cm

    Lorenzo Ramero 

    Universit{\'e} de Bordeaux I 

    Laboratoire de Math{\'e}matiques 

    351, cours de la Liberation 

    F-33405 Talence Cedex

{\em e-mail address:} {\ttfamily ramero@math.u-bordeaux.fr}

{\em web page:} {\ttfamily http://www.math.u-bordeaux.fr/$\sim$ramero}
\newpage

\tableofcontents

\newpage

\hfill{\small\rm Le bruit des vagues {\'e}tait encore plus paresseux, plus}

\hfill{\small\rm {\'e}tale qu'a midi. C'{\'e}tait le m{\^e}me soleil, la m{\^e}me}

\hfill{\small\rm lumi{\`e}re sur le m{\^e}me sable qui se prolongeait ici.}

\smallskip
\hfill{A. Camus - \it{L'{\'e}tranger}}
\bigskip

\section{Introduction}

From a pragmatic standpoint, one can describe the theory
of almost rings as a useful tool for performing calculations
of Galois cohomology groups. Indeed, this is the main 
application of Faltings' ``almost purity theorem", which
is the technical heart of \cite{Fa2}.

Though almost ring theory is developed here as an independent 
branch of mathematics, stretching somewhere in between commutative
algebra and category theory, the original applications to
Galois cohomology  still provide the main motivation and
influence largely the evolution of the subject.

It is therefore fitting to introduce the present
work by reviewing briefly the main ideas behind these calculations. 
Let us consider first a complete discretely valued field $K$ of 
zero characteristic, with perfect residue field of characteristic 
$p>0$, and uniformizer $\pi$; we denote by $K^+$ the ring of integers 
of $K$. The valuation 
$v$ of $K$ extends uniquely to any algebraic extension, and we want 
to normalize the value group in such a way that $v(p)=1$ in every 
such extension.

Let $E$ be a finite Galois extension of $K$, with Galois group
$G$. Typically, one is given a discrete $E^+[G]$-module $M$
(such that the $G$-action on $M$ is {\em semilinear}, that is, 
compatible with the $G$-action on $E^+$), and is interested
in studying the (modified) Tate  cohomology $\hat H^i:=\hat H^i(G,M)$ 
(for $i\in\Z$). (Recall that $\hat H^i(M)$ agrees with
Galois cohomology $R^i\Gamma^GM$ for $i>0$, with Galois
homology for $i<-1$, and for $i=0$ it equals $M^G/\Tr_{E/K}(M)$,
the $G$-invariants divided by the image of the trace map).

In such a situation, the scalar multiplication map 
$E^+\otimes_\Z M\to M$ induces natural cup product pairings
$\hat H^i(G,E^+)\otimes_\Z\hat H^j\to\hat H^{i+j}$. Especially, the
action of $(E^+)^G=K^+$ on $\hat H^i$ factors through 
$K^+/\Tr_{E/K}(E^+)$; in other words, the image of $E^+$ under 
the trace map annihilates the modified Tate cohomology.

If now the extension $E$ is {\em tamely ramified\/} over $K$,
then $\Tr_{E/K}(E^+)=K^+$, so the $\hat H^i$ vanish for all $i\in\Z$.
Even sharper results can be achieved when the extension is 
{\em unramified}. Indeed, in
such case $E^+$ is a $G$-torsor for the {\'e}tale topology of $K^+$,
hence, some basic descent theory tells us that the natural map
$$
E^+\otimes_{K^+}R\Gamma^GM\to M[0]
$$
is an isomorphism in the derived category of the category of 
$E^+[G]$-modules (where we have denoted by $M[0]$ the complex 
consisting of $M$ placed in degree zero).

In Tate's paper \cite{Ta} there occurs a variant of the above 
situation : instead of the finite extension $E$ one considers 
the algebraic closure $K^\mathrm{a}$ of $K$, so that $G$ is 
the absolute Galois group of $K$, and the discrete $G$-module 
$M$ is replaced by the {\em topological\/} module $C(\chi)$, 
where $C$ is the $p$-adic completion of $K^\mathrm{a}$, whose
natural $G$-action we ``twist" by a continuous character 
$\chi:G\to K^\times$. 
Then the relevant $H^\bullet$ is the {\em continuous\/} Galois 
cohomology $H_{\mathrm{cont}}^\bullet(G,C(\chi))$, which is defined
in general as the homology of a complex of continuous cochains.
Under the present assumptions, $H^i$ can be computed by the formula:
$$
H_{\mathrm{cont}}^i(G,C(\chi)):=
(\lim_{\substack{\longleftarrow \\
                 n}}~
H^i(G,K^{\mathrm{a}+}(\chi)\otimes_\Z\Z/p^n\Z))\otimes_\Z\Q.
$$
Let now $K_\infty$ be a totally ramified Galois extension
with Galois group $H$ isomorphic to $\Z_p$.
Tate realized that, for cohomological purposes,
the extension $K_\infty$ plays the role of a maximal
totally ramified Galois extension of $K$. More precisely,
let $L$ be any finite extension of $K$, and set $L_n:=L\cdot K_n$,
where $K_n$ is the subfield of $K_\infty$ fixed by 
$H^{p^n}\simeq p^n\cdot\Z_p$. The extension $K_n\subset L_n$
is unramified if and only if the different ideal 
$\cD_n:=\cD_{L_n^+/K_n^+}$ equals $L_n^+$. In case this fails, 
the valuation $v(\delta_n)$ of a generator $\delta_n$ of $\cD_n$ 
will be a strictly positive rational number, giving a quantitative
measure for the ramification of the extension. With this notation,
\cite[\S 3.2, Prop.9]{Ta} reads 
\set\begin{equation}\label{eq_Tate-reads}
\lim_{n\to\infty}v(\delta_n)=0
\end{equation}
(indeed, $v(\delta_n)$ approaches zero about as fast as $p^{-n}$).
In this sense, one can say that the extension 
$K_\infty\subset L_\infty:=L\cdot K_\infty$ is {\em almost unramified}.
One immediate consequence is that the maximal ideal $\fm$
of $K^+_\infty$ is contained in $\Tr_{L_\infty/K_\infty}(L_\infty^+)$.
If, additionally, $L$ is a Galois extension of $K$, we can consider the
subgroup 
$$
G_\infty:=\Gal(L_\infty/K_\infty)\subset\Gal(L/K)
$$ 
and the foregoing implies that $\fm$ annihilates
$H^i_\mathrm{cont}(G_\infty,M)$, for every $i>0$, and every
topological $L_\infty^+[G_\infty]$-module $M$. More precisely,
the homology of the cone of the natural morphism
\set\begin{equation}\label{eq_repeat}
L_\infty^+\otimes_{K^+_\infty}R\Gamma^{G_\infty}M\to M[0]
\end{equation}
is annihilated by $\fm$ in all degrees, {\em i.e.} it is 
{\em almost zero}. Equivalently, one says that the maps on 
homology induced by \eqref{eq_repeat} are {\em almost isomorphisms\/} 
in all degrees.

Tate goes on to apply these cohomological vanishings to the
study of $p$-divisible groups; in turns, this study enables
him to establish a comparison between the {\'e}tale and the 
Hodge cohomology of an abelian scheme over $K^+$, which has 
become the prototype for all subsequent investigation of $p$-adic 
Hodge theory.

A first generalization of \eqref{eq_Tate-reads} can be
found in the work \cite{Fr-Ma} by Fresnel and Matignon;
one interesting aspect of this work is that it does away
with any consideration of local class field theory (which
was used to get the main estimates in \cite{Ta}); instead,
Fresnel and Matignon write a general extension $L$  as a
tower of monogenic subextensions, whose structure is sufficiently
well understood to allow a direct and very explicit analysis. 
The main tool in \cite{Fr-Ma} is a notion of different
ideal $\cD_{E^+/K^+}$ for a possibly infinite algebraic field extension
$K\subset E$; then the extension $K_\infty$ considered in \cite{Ta}
is replaced by any extension $E$ of $K$ such that $\cD_{E^+/K^+}=(0)$,
and \eqref{eq_Tate-reads} is generalized by the claim that
$\cD_{F^+/E^+}=F^+$, for every finite extension $E\subset F$.

In some sense, the arguments of \cite{Fr-Ma} anticipate
those used by Faltings in the first few paragraphs of his 
fundamental article \cite{Fa1}. There we find, first of all,
a further extension of \eqref{eq_Tate-reads}: the
residue field of $K$ is now not necessarily perfect, instead
one assumes only that it admits a finite $p$-basis; then the 
relevant $K_\infty$ is an extension whose residue field is perfect,
and whose value group is $p$-divisible. This generalization
paves the way to the almost purity theorem, of which it
represents the one-dimensional case. In order to state and
prove the higher dimensional case, Faltings invents the 
method of ``almost {\'e}tale extensions", and indeed sketches
in a few pages a whole program of ``almost commutative algebra",
with the aim of transposing to the almost context as much as
possible of the classical theory. So, for instance, if $A$
is a given $K^+_\infty$-algebra, and $M$ is an $A$-module,
one says that $M$ is {\em almost flat\/} if, for every $A$-module
$N$, the natural map of complexes 
$$
M\derotimes_AN\to M\otimes_AN[0]
$$
induces almost isomorphisms on homology in all degrees.
Similarly, $M$ is {\em almost projective\/} if the same
holds for the map of complexes $\Hom_A(M,N)[0]\to R\Hom_A(M,N)$.
Then, according to \cite{Fa1}, a map $A\to B$ of $K^+_\infty$-algebras 
is called {\em almost {\'e}tale\/} if $B$ is almost projective as 
an $A$-module and as a $B\otimes_AB$-module (moreover, $B$ is required 
to be {\em almost finitely generated} : the discussion of finiteness
conditions in almost ring theory is a rather subtle business,
and we dedicate the better part of section \ref{sec_unif.spaces} to 
its clarification).

With this new language, the almost
purity theorem should be better described as an almost version
of Abhyankar's lemma, valid for morphisms $A\to B$ of $K^+$-algebras
that are {\'e}tale in characteristic zero and possibly wildly 
ramified on the locus of positive characteristic.  The actual
statement goes as follows. Suppose that
$A$ admits global {\'e}tale coordinates, that is, there exists
an {\'e}tale map $K^+[T_1^{\pm 1},...,T_d^{\pm 1}]\to A$; whereas in
the tamely ramified case a finite ramified base change
$K^+\to K^+[\pi^{1/n}]$ (with $(p,n)=1$) suffices to kill all 
ramification, the infinite extension 
$A\to A_\infty:=A[T_1^{\pm 1/p^\infty},...,T_d^{\pm 1/p^\infty}]
\otimes_{K^+}K^+_\infty$ is required in the wildly
ramified case, to kill {\em almost all\/} ramification, 
which means that the normalization  $B_\infty$ of $A_\infty\otimes_AB$ 
is almost {\'e}tale over $A_\infty$.

Faltings has proposed two distinct strategies for the proof
of his theorem : the first one, presented in \cite{Fa1}, consists
in adapting Grothendieck's proof of Zariski-Nagata's 
purity\footnote{At the time of writing, there are still some
obscure points in this proof};
a more recent one (\cite{Fa2}) uses the action of Frobenius on some 
local cohomology modules, and is actually valid under more general 
assumptions (one does not require the existence of {\'e}tale 
coordinates, but only a weaker semi-stable reduction hypothesis
on the special fibre).

As a corollary, one deduces cohomological vanishings generalizing
the foregoing : indeed, suppose that the extension of fraction fields
$\mathrm{Frac}(A)\subset\mathrm{Frac}(B)$ is Galois with group $G$;
then, granting almost purity, $B_\infty$ is an ``almost $G$-torsor"
over $A_\infty$, therefore, for any $B_\infty[G]$-module 
$M$, the natural map of complexes 
$B_\infty\otimes_{A_\infty}R\Gamma^GM\to M[0]$ induces 
almost isomorphisms on homology. Finally, these results can be used 
(together with a lot of hard work) to deduce comparison theorems 
between $p$-adic {\'e}tale cohomology and deRham (or other kinds of) 
cohomology, for arbitrary smooth projective varieties over $K$. 
This method can even be extended to treat cohomology with not 
necessarily constant coefficients (see \cite{Fa2}), thereby providing 
the most comprehensive approach to $p$-adic Hodge theory found so far.

The purpose of our text is to fully work out the foundations of
``almost commutative algebra" outlined by Faltings; in the process
we generalize and simplify considerably the theory, and also 
extend it in directions that were not explored in \cite{Fa1}, 
\cite{Fa2}. 

It turns out that most of almost ring theory can be built up
satisfactorily  from a very slim and general set of 
assumptions: our basic setup, introduced in section
\ref{sec_ring.prel}, consists of a ring $V$ and an ideal 
$\fm\subset V$ such that $\fm=\fm^2$; starting from  
\eqref{subsec_so-much} we also assume that $\fm\otimes_V\fm$ 
is a flat $V$-module : simple considerations show this to be 
a natural hypothesis, often verified in practice. 

The $V$-modules killed by $\fm$ are the objects of a (full) Serre
subcategory $\Sigma$ of the category $V\Mod$ of all $V$-modules,
and the quotient $V^a\Mod:=V\Mod/\Sigma$ is an abelian
category which we call the category of {\em almost $V$-modules}.
It is easy to check that the usual tensor product of $V$-modules
descends to a bifunctor $\otimes$ on almost $V$-modules, so that
$V^a\Mod$ is a monoidal abelian category in a natural way.
Then an {\em almost ring\/} is just an almost $V$-module $A$
endowed with a ``multiplication" morphism $A\otimes A\to A$ satisfying
certain natural axioms. Together with the obvious morphisms, these
gadgets form a category $V^a\Alg$. Given any almost $V$-algebra $A$,
one can then define the notion of $A$-module and $A$-algebra, just
like for usual rings. The purpose of the game is to reconstruct in 
this new framework as much as possible (and useful) of classical
linear and commutative algebra.
Essentially, this is the same as the
ideology informing Deligne's paper \cite{De}, which sets out to
develop algebraic geometry in the context of abstract tannakian
categories. We could also claim an even earlier ancestry, in that
some of the leading motifs resonating throughout our text, can
be traced as far back as Gabriel's memoir \cite{Ga} ``Des cat{\'e}gories
ab{\'e}liennes".

In evoking Deligne's and Gabriel's works, we have unveiled
another source of motivation whose influence has steadily
grown throughout the long gestation of our paper. Namely,
we have come to view almost ring theory as a contribution to
that expanding body of research of still uncertain range
and shifting boundaries, that we could call 
``abstract algebraic geometry''. We would like to encompass 
under this label several heterogeneous developments: notably,
it should include various versions of non-commutative geometry
that have been proposed in the last twenty years, but also the
relative schemes of \cite{Hak}, as well as Deligne's ideas for
algebraic geometry over symmetric monoidal categories.

The common thread loosely unifying these works is the
realization that "geometric spaces" do not necessarily
consist of set-theoretical points, and - perhaps more
importantly - functions on such "spaces" do not necessarily
form (sheaves of) commutative rings. Much effort has been devoted
to extending the reach of geometric intuition to non-commutative
algebras; alternatively, one can retain commutativity, but
allow ``structure sheaves'' which take values in tensor categories
other than the category of rings. As a case in point, to any
given almost ring $A$ one can attach its {\em spectrum}
$\Spec\,A$, which is just $A$ viewed as an object of the
opposite of the category $V^a\Alg$. $\Spec\,A$ has even a
natural {\em flat topology\/}, which allows to define
more general {\em almost schemes\/} by gluing ({\em i.e.}
taking colimits of) diagrams of affine spectra; all this
is explained in section \ref{sec_qu-proj-schs}, where we
also introduce {\em quasi-projective} almost schemes and
investigate some basic properties of the {\em smooth locus}
of a quasi-projective almost scheme. By way of illustration,
these generalities are applied in section \ref{sec_lift-tors}
in order to solve a deformation problem for torsors over affine
almost group schemes; let us stress that the problem in question
is stated purely in terms of affine objects ({\em i.e.} almost
rings and ``almost Hopf algebras''), but the solution requires
the introduction of certain auxiliary almost schemes that are not
affine.

Having soared into the thin air of abstract geometry, we
come back to earth in the last two chapters, which deal
with applications to valuation theory and to $p$-adic
analytic geometry : especially the reader will find there our
own contributions to almost purity. This general presentation
has thus come full circle, and we defer to the introductory
remarks at the beginning of the respective chapters for a
more detailed description of our results.

\begin{acknowledgement} The second author is very much indebted
to Gerd Faltings for many patient explanations on the 
method of almost  {\'e}tale extensions.
Next he would like to acknowledge several interesting 
discussions with Ioannis Emmanouil. He is also much obliged 
to Pierre Deligne, for a useful list of critical remarks.
Finally, he owes a special thank to Roberto Ferretti, who has 
read the first tentative versions of this work, has corrected 
many slips and has made several valuable suggestions.

This project began in 1997 while the second author was
supported the IHES, and has been brought to completion
while the second author was a guest of the Laboratoire de
Th\'eorie des Nombres of the University of Paris VI.
\end{acknowledgement}

\newpage

\section{Homological theory}\label{ch_homol-th}

As explained in the introduction, in order to define a category
of almost modules one requires a pair $(V,\fm)$ consisting of
a ring $V$ and an ideal $\fm\subset V$ such that $\fm=\fm^2$.
In section \ref{sec_ring.prel} we collect a few useful ring-theoretic
preliminaries concerning such pairs. In section \ref{sec_alm.cat}
we introduce the category $V^a\Mod$ of {\em almost modules} : it is
a quotient $V\Mod/\Sigma$ of the category of $V$-modules, where
$\Sigma$ is the thick subcategory of the $V$-modules killed by $\fm$.
$V^a\Mod$ is an abelian tensor category and its commutative
unitary monoids, called {\em almost algebras}, are the chief
objects of study in this work. 
The first useful observation is that the localization functor 
$V\Mod\to V^a\Mod$ admits both left and right adjoints. Taken 
together, these functors exhibit the kind of exactness properties 
that one associates to open embeddings of topoi, perhaps a hint 
of some deeper geometrical structure, still to be unearthed.

After these generalities, we treat in section \ref{sec_unif.spaces}
the question of finiteness conditions for almost modules.
Let $A$ denote an almost algebra, fixed for the rest of
this introduction.
It is certainly possible to define as usual a notion of finitely
generated $A$-module, however this turns out to be too restrictive 
a class for applications.
The main idea here is to define a uniform structure on the 
set of equivalence classes of $A$-modules ; then we will say that 
an $A$-module is {\em almost finitely generated\/} if its 
isomorphism class lies in the topological closure of the subspace 
of finitely generated $A$-modules. Similarly we define
{\em almost finitely presented\/} $A$-modules. The uniform
structure also comes handy when we want to construct operators
on almost modules : if one can show that the operator in question
is uniformly continuous on a class $\cC$ of almost modules, then
its definition extends right away by continuity to the topological
closure $\bar\cC$ of $\cC$. This is exemplified by the construction
of the (almost) Fitting ideals for $A$-modules, at the end of section 
\ref{sec_unif.spaces}.

In section \ref{sec_homol} we introduce the basic toolkit of homological
algebra, beginning with the notion of flat almost module, which poses
no problem, since we do have a tensor product in our category.
The notion of projectivity is more subtle : it turns out that
the category of $A$-modules usually does {\em not} have enough 
projectives. The useful notion is {\em almost projectivity}: simply
one uses the standard definition, except that the role of the 
$\Hom$ functor is played by the internal $\Alhom$ functor.
The scarcity of projectives should not be regarded as surprising
or pathological: it is quite analogous to the lack of
enough projective objects in the category of quasi-coherent 
$\cO_X$-modules on a non-affine scheme $X$. 

Section \ref{sec_hot} introduces the cotangent complex of a morphism
of almost rings, and establishes its usual properties, such as
transitivity and Tor-independent base change theorems.
These foundations will be put to use in chapter \ref{ch_alm.ring.th},
to study infinitesimal deformations of almost algebras.

\subsection{Some ring-theoretic preliminaries}
\label{sec_ring.prel}

Unless otherwise stated, every ring is commutative with unit.
This section collects some results of general nature that will
be used throughout this work.

\sset\subsubsection{}\label{subsec_basic.setup}
\index{$(V,\fm)$, $\tilde\fm$ : Basic setup|indref{subsec_basic.setup}}
Our {\em basic setup\/} consists of a fixed base ring $V$ containing 
an ideal $\fm$ such that $\fm^2=\fm$. Starting from 
\eqref{subsec_so-much}, we will also assume that 
$\tilde\fm:=\fm\otimes_V\fm$ is a flat $V$-module. 

\begin{example}\label{ex_rings}
(i) The main example is given by a non-discrete valuation
ring $(V,|\cdot|)$ of rank one; in this case $\fm$ will be 
the maximal ideal. 

(ii) Take $\fm:=V$. This is the ``classical limit''.  
\index{Classical limit|indref{ex_rings}}
In this case almost ring theory reduces to usual ring theory. 
Thus, all the discussion that follows specialises to, and 
sometimes gives alternative proofs for, statements about rings 
and their modules.
\end{example}

\sset\subsubsection{}\label{subsec_almost_zero}
\index{Almost zero module|indref{subsec_almost_zero}}
\index{Almost isomorphism|indref{subsec_almost_zero}}
Let $M$ be a given $V$-module. We say that $M$ is 
{\em almost zero\/} if $\fm M=0$. A map $\phi$ of 
$V$-modules is an {\em almost isomorphism\/} if both 
$\Ker\,\phi$ and $\Coker\,\phi$ are almost zero $V$-modules.

\begin{remark}\label{rem_almost.zero}
(i) It is easy to check that a $V$-module $M$ is almost zero if and
only if $\fm\otimes_VM=0$. Similarly, a map $M\to N$ of $V$-modules 
is an almost isomorphism if and only if the induced map 
$\tilde\fm\otimes_VM\to\tilde\fm\otimes_VN$ is an isomorphism. 
Notice also that, if $\fm$ is flat, then $\fm\simeq\tilde\fm$.

(ii) Let $V\to W$ be a ring homomorphism. For a $V$-module $M$
set $M_W:=W\otimes_VM$. We have an exact sequence 
\set\begin{equation}\label{eq_economy}
0\to K\to\fm_W\to\fm W\to 0
\end{equation} 
where $K:=\Tor^V_1(V/\fm,W)$ is an almost zero $W$-module. By (i) 
it follows that $\fm\otimes_VK\simeq(\fm W)\otimes_WK\simeq 0$. 
Then, applying $\fm_W\otimes_W-$ and $-\otimes_W(\fm W)$
to \eqref{eq_economy} we derive 
$$\fm_W\otimes_W\fm_W\simeq\fm_W\otimes_W(\fm W)\simeq
(\fm W)\otimes_W(\fm W)$$ 
{\em i.e.\/} $\tilde\fm_W\simeq(\fm W)^\sim$. In particular, 
if $\tilde\fm$ is a flat $V$-module, then $\tilde\fm_W$ is a flat 
$W$-module. This means that our basic assumptions on the pair 
$(V,\fm)$ are stable under arbitrary base extension. Notice that 
the flatness of $\fm$ does not imply the flatness of $\fm W$. 
This partly explains why we insist that $\tilde\fm$, rather than 
$\fm$, be flat. 
\end{remark}

\sset\subsubsection{}\label{subsec_conditions.A.B}
\index{Conditions {\bf A} and {\bf B}|indref{subsec_conditions.A.B}}
Before moving on, we want to analyze in some detail how
our basic assumptions relate to certain other natural 
conditions that can be postulated on the pair $(V,\fm)$.
Indeed, let us consider the following two hypotheses :
\smallskip

\noindent ({\bf A})\ \  $\fm=\fm^2$ and $\fm$ is a filtered 
union of principal ideals.
\smallskip

\noindent ({\bf B})\ \  $\fm=\fm^2$ and, for all integers $k>1$, 
the $k$-th powers of elements of $\fm$ generate $\fm$.
\smallskip

Clearly ({\bf A}) implies ({\bf B}). Less obvious is
the following result.
\begin{proposition}\label{prop_less.obv}
{\em(i)} {\em({\bf A})} implies that $\tilde\fm$ is flat.
\begin{enumerate}
\addtocounter{enumi}{1}
\item
If $\tilde\fm$ is flat then {\em({\bf B})} holds. 
\end{enumerate}
\end{proposition}
\begin{proof} Suppose that ({\bf A}) holds, so that 
$\fm=\colim{\alpha\in I}Vx_\alpha$, where $I$ is a directed set
parametrizing elements $x_\alpha\in\fm$ (and 
$\alpha\leq\beta\Leftrightarrow Vx_\alpha\subset Vx_\beta$).
For any $\alpha\in I$ we have natural isomorphisms
\set\begin{equation}\label{eq_pi.alpha}
Vx_\alpha\simeq V/\Ann_V(x_\alpha)\simeq
(Vx_\alpha)\otimes_V(Vx_\alpha).
\end{equation}
For $\alpha\leq\beta$, let 
$j_{\alpha\beta}:Vx_\alpha\hookrightarrow Vx_\beta$ be the imbedding;
we have a commutative diagram
$$\xymatrix{
V \ar[rr]^{\mu_{z^2}} \ar[d]_{\pi_\alpha} & & V \ar[d]^{\pi_\beta} \\
(Vx_\alpha)\otimes_V(Vx_\alpha) 
\ar[rr]^{j_{\alpha\beta}\otimes j_{\alpha\beta}} & &
(Vx_\beta)\otimes_V(Vx_\beta)
}$$
where $z\in V$ is such that $x_\alpha=z\cdot x_\beta$, $\mu_{z^2}$
is multiplication by $z^2$ and $\pi_\alpha$ is the projection
induced by \eqref{eq_pi.alpha} (and similarly for $\pi_\beta$).
Since $\fm=\fm^2$, for all $\alpha\in I$ we can find $\beta$
such that $x_\alpha$ is a multiple of $x_\beta^2$. Say 
$x_\alpha=y\cdot x^2_\beta$; then we can take $z:=y\cdot x_\beta$,
so $z^2$ is a multiple of $x_\alpha$ and in the above diagram
$\Ker\,\pi_\alpha\subset\Ker\,\mu_{z^2}$. Hence one can define
a map $\lambda_{\alpha\beta}:(Vx_\alpha)\otimes_V(Vx_\alpha)\to V$
such that $\pi_\beta\circ\lambda_{\alpha\beta}=
j_{\alpha\beta}\otimes j_{\alpha\beta}$ and 
$\lambda_{\alpha\beta}\circ\pi_\alpha=\mu_{z^2}$. It now
follows that for every $V$-module $N$, the induced morphism
$\Tor_1^V(N,(Vx_\alpha)\otimes_V(Vx_\alpha))\to
\Tor_1^V(N,(Vx_\beta)\otimes_V(Vx_\beta))$
is the zero map. Taking the colimit we derive that $\tilde\fm$
is flat. This shows (i).
In order to show (ii) we consider, for any prime number $p$, 
the following condition 
\medskip

\noindent($*_p$)\qquad $\fm/p\cdot\fm$ is generated (as a $V$-module)
by the $p$-th powers of its elements.
\medskip

Clearly ({\bf B}) implies ($*_p$) for all $p$. In fact we have :
\begin{claim}
({\bf B}) holds if and only if ($*_p$) holds for every prime 
$p$. 
\end{claim}
\begin{pfclaim} Suppose that ($*_p$) holds for every prime $p$.
The polarization identity 
$$k!\cdot x_1\cdot x_2\cdot...\cdot x_k=
\sum_{I\subset\{1,2,...,k\}}(-1)^{k-|I|}\cdot
\left(\sum_{i\in I}x_i\right)^k$$
shows that if $N:=\sum_{x\in\fm}Vx^k$ then 
$k!\cdot\fm\subset N$.
To prove that $N=\fm$ it then suffices to show that for every
prime $p$ dividing $k!$ we have $\fm=p\cdot\fm+N$. Let 
$\phi:V/pV\to V/pV$ be the Frobenius 
($x\mapsto x^p$); 
we can denote by $(V/p V)^\phi$ the ring $V/pV$ 
seen as a $V/pV$-algebra via the homomorphism $\phi$. 
Also set $\phi^*M:=M\otimes_{V/pV}(V/pV)^\phi$ 
for a $V/pV$-module $M$. Then the map 
$\phi^*(\fm/p\cdot\fm)\to(\fm/p\cdot\fm)$ (defined by raising 
to $p$-th power) is surjective by ($*_p$). Hence so is
$(\phi^r)^*(\fm/p\cdot\fm)\to(\fm/p\cdot\fm)$ for every $r>0$,
which says that $\fm=p\cdot\fm+N$ when $k=p^r$, hence for every
$k$.\end{pfclaim}

Next recall (see \cite[Exp. XVII 5.5.2]{SGA4-3}) that, if $M$
is a $V$-module, the module of symmetric tensors 
$\mathrm{TS}^k(M)$
is defined as $(\otimes_V^kM)^{S_k}$, the invariants under the
natural action of the symmetric group $S_k$ on $\otimes_V^kM$.
We have a natural map $\Gamma^k(M)\to\mathrm{TS}^k(M)$ that 
is an isomorphism when $M$ is flat (see 
{\em loc. cit.\/} 5.5.2.5;
here $\Gamma^k$ denotes the $k$-th graded piece of the divided 
power algebra).

\begin{claim} The group $S_k$ acts trivially on 
$\otimes_V^k\fm$ and the map $\tilde\fm\otimes_V\fm\to\tilde\fm$ 
($x\otimes y\otimes z\mapsto x\otimes yz$) is an isomorphism.
\end{claim}
\begin{pfclaim} The first statement is reduced to the case
of transpositions and to $k=2$. There we can compute :
$x\otimes yz=xy\otimes z=y\otimes xz=yz\otimes x$.
For the second statement note that the imbedding 
$\fm\hookrightarrow V$ is an almost isomorphism, 
and apply remark \ref{rem_almost.zero}(i).
\end{pfclaim}

Suppose now that $\tilde\fm$ is flat and pick a prime $p$. 
Then $S_p$ acts trivially on $\otimes^p_V\tilde\fm$. Hence 
\set\begin{equation}\label{eq_iso.Gamma}
\Gamma^p(\tilde\fm)\simeq\otimes^p_V\tilde\fm\simeq\tilde\fm.
\end{equation}
But $\Gamma^p(\tilde\fm)$ is spanned as a $V$-module by
the products $\gamma_{i_1}(x_1)\cdot...\cdot\gamma_{i_k}(x_k)$
(where $x_i\in\tilde\fm$ and $\sum_ji_j=p$). Under the isomorphism
\eqref{eq_iso.Gamma} these elements map to 
$\binom{p}{i_1,...,i_k}\cdot x_1^{i_1}\cdot...\cdot x_k^{i_k}$;
but such an element vanishes in $\tilde\fm/p\cdot\tilde\fm$ 
unless $i_k=p$ for some $k$. Therefore $\tilde\fm/p\cdot\tilde\fm$
is generated by $p$-th powers, so the same is true for $\fm/p\cdot\fm$,
and by the above, ({\bf B}) holds, which shows (ii).
\end{proof}

\begin{theorem}\label{th_countably.pres} 
Let $(\eps_i~|~i\in I)$ be a family of generators of 
$\fm$ and, for every subset $S\subset I$, denote by 
$\fm_S\subset\fm$ the subideal generated by 
$(\eps_i~|~i\in S)$. Then we have:
\begin{enumerate}
\item
Every countable subset $S\subset I$ is contained in
another countable $S'\subset I$ such that:
\begin{enumerate}
\item
$\fm_{S'}^2=\fm_{S'}$.
\item
If\/ $\tilde\fm$ is a flat $V$-module, the same holds
for $\fm_{S'}\otimes_V\fm_{S'}$.
\end{enumerate}
\item
Suppose that $\fm$ is countably generated as a $V$-module. 
Then :
\begin{enumerate}
\item
$\tilde\fm$ is countably presented as a $V$-module.
\item
If\/ $\tilde\fm$ is a flat $V$-module, then it is of 
homological dimension $\leq 1$.
\end{enumerate}
\end{enumerate}
\end{theorem}
\begin{proof}
For every $i\in I$, we can write $\eps_i=\sum_jx_{ij}\eps_j$, 
for certain $x_{ij}\in\fm$.  For any $i,j\in I$ such that 
$x_{ij}\neq 0$, let us write 
$x_{ij}=\sum_kx_{ijk}\eps_k$ for some $x_{ijk}\in V$.
We say that a subset $S\subset I$ is {\em saturated\/} 
if the following holds: whenever $i\in S$ and $x_{ijk}\neq 0$,
we have $j,k\in S$. 
\begin{claim}\label{cl_saturation} 
Every countable subset $S\subset I$ is contained
in a countable saturated subset $S_\infty\subset I$.
\end{claim}
\begin{pfclaim} We let $S_0:=S$ and define recursively
$S_n$ for every $n>0$ as follows. Suppose $S_{n-1}$ has
already been given; then we set
$$
S_n:=S_{n-1}\cup\{i\in I~|~\text{there exists $a\in S_{n-1}$,
$b\in I$ such that either $x_{aib}\neq 0$ or $x_{abi}\neq 0$}\}.
$$
Notice that, for every $i\in I$ we have $x_{ij}=0$ for all
but finitely many $j\in I$, hence $x_{ijk}=0$ for all
but finitely many $j,k\in I$. It follows easily that 
$S_\infty:=\bigcup_{n\in\N}S_n$ will do.
\end{pfclaim}

(i.a) is now straightforward, when one remarks that
$\fm_S=\fm_S^2$ for every saturated subset $S$.

Next, let $S\subset I$ be any saturated subset.
Clearly the family $(\eps_i\otimes\eps_j~|~i,j\in S)$ 
generates $\tilde\fm_S:=\fm_S\otimes_V\fm_S$ 
and $\eps_i\cdot\eps_j\cdot(\eps_k\otimes\eps_l)=
\eps_k\cdot\eps_l\cdot(\eps_i\otimes\eps_j)$ for all 
$i,j,k,l\in S$. 
Let $F(S)$ be the $V$-module defined by
generators $(e_{ij})_{i,j\in S}$, subject to the 
relations:
$$\eps_i\cdot\eps_j\cdot e_{kl}=
\eps_k\cdot\eps_l\cdot e_{ij}\qquad 
e_{ik}=\sum_jx_{ij}e_{jk}\qquad\text{for all $i,j,k,l\in S$.}$$
We get an epimorphism $\pi_S:F(S)\to\tilde\fm_S$ by 
$e_{ij}\mapsto\eps_i\otimes\eps_j$. The relations 
imply that, if $x:=\sum_{k,l}y_{kl}e_{kl}\in\Ker\,\pi_S$, 
then $\eps_i\cdot\eps_j\cdot x=0$, so $\fm_S\cdot\Ker\,\pi_S=0$.
Whence $\fm_S\otimes_V\Ker\,\pi_S=0$ and 
$\one_{\fm_S}\otimes_V\pi_S$ is an isomorphism. We consider 
the diagram 
$$\xymatrix{\fm_S\otimes_VF(S) \ar[r]^-{\sim} \ar[d]_\phi & 
\fm_S\otimes_V\tilde\fm_S \ar[d]^\psi \\ F(S) \ar[r]^-{\pi_S} & 
\tilde\fm_S}$$
where $\phi$ and $\psi$ are induced by scalar multiplication.
We already know that $\psi$ is an isomorphism, and since 
$F(S)=\fm_S\cdot F(S)$, we see that $\phi$ is an epimorphism, 
so $\pi_S$ is an isomorphism. If now $I$ is countable, this 
shows that (ii.a) holds. Now (ii.b) follows from (ii.a) and 
the following lemma \ref{lem_flat.aleph-pres}. In order to 
show (i.b) we will use the following well known criterion.
\begin{claim}\label{cl_Lazard.crit} (\cite[Ch.I, Th.1.2]{La}).
Let $R$ be a ring, $M$ an $R$-module. Then
$M$ is $R$-flat if and only if, for every finitely presented
$R$-module $N$, every morphism $N\to M$ factors as
a composition $N\to L\to M$ where $L$ is a free $R$-module
of finite rank. 
\end{claim}

In order to apply claim \ref{cl_Lazard.crit} we show:
\begin{claim} Let $S\subset I$ be a countable saturated
subset. Then there exists a countable saturated subset 
$\sigma(S)\subset I$ containing $S$, with the following 
property. For every finitely presented $V$-module $N$ and 
every morphism $f:N\to\tilde\fm_S$, we can find a
commutative diagram:
$$
\xymatrix{ N \ar[r]^-f \ar[d] & \tilde\fm_S \ar[d]^-j \\
L \ar[r] & \tilde\fm_{\sigma(S)}
}$$
where $L$ is  free $V$-module of finite rank and
$j$ is the natural map.
\end{claim}
\begin{pfclaim} In view of (ii.a), we can write 
$\tilde\fm_S\simeq\colim{\alpha\in A}M_\alpha$,
where $A$ is some filtered countable set and every
$M_\alpha$ is a finitely presented $V$-module.
Given $f$ as in the claim, we can find $\alpha\in A$
such that $f$ factors thorugh the natural map
$\iota_\alpha:M_\alpha\to\tilde\fm_S$, so we are 
reduced to prove the claim for $N=M_\alpha$ and
$f=\iota_\alpha$. However, by assumption $\tilde\fm$
is flat, so by claim \ref{cl_Lazard.crit} the composition 
$M_\alpha\stackrel{\iota_\alpha}{\to}\tilde\fm_S\to\tilde\fm$
factors through some map $L_\alpha\to\tilde\fm$, with
$L_\alpha$ free of finite rank over $V$.
Furthermore, thanks to claim \ref{cl_saturation} we have 
$\tilde\fm=\colim{J}{\tilde\fm_J}$, where $J$ runs over 
the family of all countable saturated subsets of $I$.
It follows that, for some countable saturated 
$S_\alpha\supset S$, the map $L_\alpha\to\tilde\fm$ 
factors through $\tilde\fm_{S_\alpha}$. Clearly 
$\sigma(S):=\bigcup_{\alpha\in A}S_\alpha$ will do.
\end{pfclaim}

Finally, for any countable subset $S\subset I$,
let us set $S':=\bigcup_{n\in\N}\sigma^n(S_\infty)$,
where $S_\infty$ is the saturation of $S$ as in claim
\ref{cl_saturation}. One verifies easily that (i.b) 
holds for this choice of $S'$.
\end{proof}

The proof of following well known lemma is due to D.Lazard
(\cite[Ch.I, Th.3.2]{La}), up to some slight imprecisions 
which were corrected in \cite[pp.49-50]{Mit}. For the
convenience of the reader we reproduce the argument.

\begin{lemma}\label{lem_flat.aleph-pres} Let $R$ be any
ring. A flat countably presented $R$-module has homological 
dimension $\leq 1$.
\end{lemma}
\begin{proof} Let $M$ be flat and countable presented over $R$,
and choose a presentation $F_1\stackrel{\phi}{\to}F_0\to M\to 0$
with $F_0$ and $F_1$ free $R$-modules of (infinite) countable rank.
Let $(e_j~|~j\in\N)$ be a basis of $F_0$ and $(f_i~|~i\in\N)$
a basis of $F_1$. Say that $\phi(f_i)=\sum_ix_{ij}e_j$ with 
$x_{ij}\in R$, for every $i\in\N$. For every $i\in\N$ we
define the finite set $S_i:=\{j\in\N~|~x_{ij}\neq 0\}$; 
also, let $S'_n:=\bigcup_{i\leq n}(\{i\}\cup S_i)$.
Define $R$-modules $G_n:=\oplus_{i\leq n}f_iR$, 
$H_n:=\oplus_{j\in S'_n}e_jR$; the restriction of $\phi$
induces maps $\phi_n:G_n\to H_n$ and we have 
$M\simeq\colim{n\in\N}\,\Coker\,\phi_n$.
By claim \ref{cl_Lazard.crit}, the natural map $\Coker\,\phi_n\to M$
factors as a composition
$\Coker\,\phi_n\stackrel{\alpha_n}{\to}L_n\stackrel{\psi_n}{\to}M$, 
where $L_n$ is a free $R$-module of finite rank. Then $\psi_n$ factors
through a map $\psi'_n:L_n\to\Coker\,\phi_{t(n)}$ for some 
$t(n)\in\N$. We further define $\psi''_n:L_n\to\Coker\,\phi_{k(n)}$
as the composition of $\psi'_n$ with the natural transition
map $\Coker\,\phi_{t(n)}\to\Coker\,\phi_{k(n)}$, where $k(n)$
is suitably chosen, so that $k(n)>\max(n,t(n))$ and the composition
$\psi''_n\circ\alpha_n:\Coker\,\phi_n\to\Coker\,\phi_{k(n)}$ is
the natural transition map.
We define by induction on $n\in\N$ a direct
system of maps $\tau_n:L_{h_n}\to L_{h_{n+1}}$, as follows. Set 
$h_0:=0$; next, suppose that $h_n\in\N$ has already been
given up to some $n\in\N$; we let $h_{n+1}:=k(h_n)$ and
$\tau_n:=\alpha_{h_{n+1}}\circ\psi''_{h_n}$. Clearly
$M\simeq\colim{n\in\N}L_{h_n}$; set $L:=\oplus_{n\in\N}L_{h_n}$
and let $\tau:L\to L$ be the map given by the rule
$\tau(x_n~|~n\in\N):=(x_n-\tau_n(x_{n-1})~|~n\in\N)$.
We derive a short exact sequence:
$0\to L\stackrel{\tau}{\to}L\to M\to 0$, whence the claim.
\end{proof}

\subsection{Categories of almost modules and algebras}
\label{sec_alm.cat}
\index{$\cC$ : Category(ies)! $\cC^{\,o}$ : opposite of a|indref{sec_alm.cat}}
\index{$\cC$ : Category(ies)! $s.\cC$ : simplicial objects of a|indref{sec_alm.cat}} 
\index{$\cC$ : Category(ies)! $\one_X$ : identity morphism of an object 
in a|indref{sec_alm.cat}} 
 \index{$\cC$ : Category(ies)! $\sD(\cC)$, $\sD^+(\cC)$, $\sD^-(\cC)$ : 
derived categories of a|indref{subsec_deriv.cat.notate}}
If $\cC$\ \ is a category, and $X,Y$ two objects of $\cC$, we will
usually denote by $\Hom_\cC(X,Y)$ the set of morphisms in $\cC$\ \
from $X$ to $Y$ and by $\one_X$ the identity morphism of $X$. 
Moreover we denote by $\cC^{\,o}$ the opposite category of $\cC$ 
and by $s.\cC$ the category of simplicial objects over $\cC$, 
that is, functors $\Delta^o\to\cC$, where $\Delta$ is the category 
whose objects are the ordered sets $[n]:=\{0,...,n\}$ for each integer 
$n\ge 0$ and where a morphism $\phi:[p]\to[q]$ is a non-decreasing
map. A morphism $f:X\to Y$ in $s.\cC$\ \ is a sequence of morphisms
$f_{[n]}:X[n]\to Y[n]$, $n\ge 0$ such that the obvious diagrams commute. 
We can imbed $\cC$\ \ in $s.\cC$\ \ by sending each object $X$ to the 
``constant'' object $s.X$ such that $s.X[n]=X$ for all $n\ge 0$
and $s.X[\phi]=\one_X$ for all morphisms $\phi$ in $\Delta$. 

\sset\subsubsection{}\label{subsec_deriv.cat.notate}
\index{$R\Mod$, $R\Alg$ (for a ring $R$)|indref{subsec_deriv.cat.notate}}
\index{$\Set$ : category of sets|indref{subsec_deriv.cat.notate}}
If $\cC$ is an abelian category, $\sD(\cC)$ will denote the 
derived category of $\cC$. As usual we have also the full subcategories 
$\sD^+(\cC),\sD^-(\cC)$ of complexes of objects of $\cC$ that are 
exact for sufficiently large negative (resp. positive) degree.
If $R$ is a ring, the category of $R$-modules (resp. $R$-algebras) will 
be denoted by $R\Mod$ (resp. $R\Alg$). Most of the times we will write 
$\Hom_R(M,N)$ instead of $\Hom_{R\Mod}(M,N)$. 

We denote by $\Set$ the category of sets. 
The symbol $\N$ denotes the set of non-negative 
integers; in particular $0\in\N$.

\sset\subsubsection{}\label{subsec_define.V^a-mod}
\index{Almost module(s)|indref{subsec_define.V^a-mod}}
\index{$V^a\Mod$, $(V,\fm)^a\Mod$|indref{subsec_define.V^a-mod}}
\index{$M^a$|indref{subsec_define.V^a-mod}}
The full subcategory $\Sigma$ of $V\Mod$ consisting of all $V$-modules
that are almost isomorphic to $0$ is clearly a Serre subcategory and
hence we can form the quotient category $V\Mod/\Sigma$.
There is a localization functor 
$$V\Mod\to V\Mod/\Sigma \qquad M\mapsto M^a$$ 
that takes a $V$-module $M$ to the same module, seen as an object 
of $V\Mod/\Sigma$. In particular, we have the object $V^a$
associated to $V$; it seems therefore natural to use the notation 
$V^a\Mod$ for the category $V\Mod/\Sigma$, and an object of 
$V^a\Mod$ will be indifferently referred to as ``a $V^a$-module''
or ``an almost $V$-module''. In case we need to stress the
dependance on the ideal $\fm$, we can write $(V,\fm)^a\Mod$.

Since the almost isomorphisms form a multiplicative system (see 
{\em e.g.\/} \cite[Exerc.10.3.2]{We}), it is possible to describe 
the morphisms in $V^a\Mod$ via a calculus of fractions, as follows. 
Let $V\Aliso$ be the category that has the same objects as $V\Mod$,
but such that $\Hom_{V\Aliso}(M,N)$ consists of all almost
isomorphisms $M\to N$. If $M$ is any object of $V\Aliso$ we write 
$(V\Aliso/M)$ for the category of objects of $V\Aliso$ over $M$ 
({\em i.e.\/} morphisms $\phi:X\to M$). If $\phi_i:X_i\to M$ $(i=1,2)$ are
two objects of $(V\Aliso/M)$ then $\Hom_{(V\Aliso/M)}(\phi_1,\phi_2)$
consists of all morphisms $\psi:X_1\to X_2$ in $V\Aliso$ such that
$\phi_1=\phi_2\circ\psi$. For any two $V$-modules $M,N$ we define a 
functor $\cF_N:(V\Aliso/M)^o\to V\Mod$ by
associating to an object $\phi:P\to M$ the $V$-module $\Hom_V(P,N)$ 
and to a morphism $\alpha:P\to Q$ the map 
$\Hom_V(Q,N)\to\Hom_V(P,N)~:~\beta\mapsto\beta\circ\alpha$.
Then we have
\set\begin{equation}\label{eq_alhom}
\Hom_{V^a\Mod}(M^a,N^a)=\colim{(V\Aliso/M)^o}\cF_N.
\end{equation}
However, formula \eqref{eq_alhom} can be simplified considerably
by remarking that for any $V$-module $M$, the natural morphism 
$\tilde\fm\otimes_VM\to M$ is an initial object of $(V\Aliso/M)$.
Indeed, let $\phi:N\to M$ be an almost isomorphism; the diagram
$$\xymatrix{
\tilde\fm\otimes_VN \ar[r]^\sim \ar[d] & 
\tilde\fm\otimes_VM \ar[d] \\
N \ar[r]^\phi & M
}$$
(cp. remark \ref{rem_almost.zero}(i)) allows one to define a 
morphism $\psi:\tilde\fm\otimes_VM\to N$ over $M$. We need to
show that $\psi$ is unique. But if   
$\psi_1,\psi_2:\tilde\fm\otimes_VM\to N$ are two maps over
$M$, then $\Img(\psi_1-\psi_2)\subset\Ker(\phi)$ is almost 
zero, hence $\Img(\psi_1-\psi_2)=0$, since 
$\tilde\fm\otimes_VM=\fm\cdot(\tilde\fm\otimes_VM)$. 
Consequently, \eqref{eq_alhom} boils down to
\set\begin{equation}\label{eq_alm.morph}
\Hom_{V^a\Mod}(M^a,N^a)=\Hom_V(\tilde\fm\otimes_VM,N).
\end{equation}
In particular $\Hom_{V^a\Mod}(M,N)$ has a natural structure of $V$-module
for any two $V^a$-modules $M,N$, {\em i.e.\/} $\Hom_{V^a\Mod}(-,-)$
is a bifunctor that takes values in the category $V\Mod$.

\sset\subsubsection{}\label{subsec_comm.constr}
\index{$\cC$ : Category(ies)!tensor|indref{subsec_comm.constr}}
\index{$\cC$ : Category(ies)!tensor!monoid, algebra in
  a|indref{subsec_comm.constr}}
\index{$\cC$ : Category(ies)!tensor!(left, right, unitary) module in a
|indref{subsec_comm.constr}}
\index{Almost algebra(s)|indref{subsec_comm.constr}}
\index{Almost algebra(s)!ideal in an|indref{subsec_comm.constr}}
\index{Almost module(s)!$\Ann_A(M)$ : annihilator of an|indref{subsec_comm.constr}}
\index{$\cC$ : Category(ies)!tensor!\thetaperv : switch operator
for two modules in a|indref{subsec_comm.constr}}
\index{$\cC$ : Category(ies)!tensor!$\mu_A$ : multiplication
of a monoid in a|indref{subsec_comm.constr}}
\index{$\cC$ : Category(ies)!tensor!$\sigma_{M/A}$ : scalar
multiplication of a module in a|indref{subsec_comm.constr}}
One checks easily (for instance using \eqref{eq_alm.morph}) that 
the usual tensor product induces a bifunctor $-\otimes_V-$ on almost
$V$-modules, which, in the jargon of \cite{DeM} makes of $V^a\Mod$ an 
{\em abelian tensor category\/}. Then an {\em almost $V$-algebra} 
is just a commutative unitary monoid in the tensor category $V^a\Mod$. 
Let us recall what this means. Quite generally, let $(\cC,\otimes,U)$ 
be any abelian tensor category, so that $\otimes:\cC\times\cC\to\cC$\ \ 
is a biadditive functor, $U$ is the identity object
of $\cC$ (see \cite[p.105]{DeM}) and for any two objects $M$ 
and $N$ in $\cC$\ \ we have a ``commutativity constraint'' ({\em i.e.\/}
a functorial isomorphism $\theta_{M|N}:M\otimes N\to N\otimes M$ 
that ``switches the two factors'') and a functorial isomorphism
$\nu_M:U\otimes M\to M$. Then a $\cC$-monoid $A$ is an
object of $\cC$\ \ endowed with a morphism $\mu_A:A\otimes A\to A$ 
(the ``multiplication'' of $A$) satisfying the associativity condition
$$\mu_A\circ(\one_A\otimes\mu_A)=\mu_A\circ(\mu_A\otimes\one_A).$$
We say that $A$ is {\em unitary\/} if additionally $A$ is endowed 
with a ``unit morphism'' $\underline 1_A:U\to A$ satisfying the 
(left and right) unit property : 
$$\mu_A\circ(\underline 1_A\otimes\one_A)=\nu_A\qquad
\mu_A\circ(\underline 1_A\otimes\one_A)\circ\theta_{A|U}=
\mu_A\circ(\one_A\otimes \underline 1_A).$$
Finally $A$ is {\em commutative\/} if $\mu_A=\mu_A\circ\theta_{A|A}$
(to be rigorous, in all of the above one should indicate the 
associativity constraints, which we have omitted : see \cite{DeM}).
A commutative unitary monoid will also be simply called
an {\em algebra\/}. 
With the morphisms defined in the obvious way, the $\cC$-monoids 
form a category; furthermore, given a $\cC$-monoid $A$, a 
{\em left $A$-module} is an object $M$ of $\cC$ 
endowed with a morphism $\sigma_{M/A}:A\otimes M\to M$ 
such that $\sigma_{M/A}\circ(\one_A\otimes\sigma_{M/A})=
\sigma_{M/A}\circ(\mu_A\otimes\one_M)$.
Similarly one defines right $A$-modules and $A$-bimodules. 
In the case of bimodules we have left and right morphisms
$\sigma_{M,l}:A\otimes M\to M$, $\sigma_{M,r}:M\otimes A\to M$
and one imposes that they ``commute'', {\em i.e.\/} that
$$\sigma_{M,r}\circ(\sigma_{M,l}\otimes\one_A)=
\sigma_{M,l}\circ(\one_A\otimes\sigma_{M,r}).$$
Clearly the (left resp. right) $A$-modules (and the 
$A$-bimodules) form an additive category with 
{\em $A$-linear morphisms\/} defined as one expects.
One defines the notion of a submodule as an equivalence class of
monomorphisms $N\to M$ such that the composition 
$A\otimes N\to A\otimes M\to M$ factors through $N$. 
Especially, a {\em two-sided ideal\/} of $A$ is an
$A$-sub-bimodule $I\to A$. For given submodules $I,J$ 
of $A$ one denotes 
$IJ:=\Img(I\otimes J\to A\otimes A\stackrel{\mu_A}{\to}A)$.
For an $A$-module $M$, the {\em annihilator\/} $\Ann_A(M)$
of $M$ is the largest ideal $j:I\to A$ of $A$ such that 
the composition
$I\otimes M\stackrel{j\otimes\one_M}{\longrightarrow} 
A\otimes M\stackrel{\sigma_{M/A}}{\longrightarrow}M$ is 
the zero morphism.

\sset\subsubsection{} 
If $f:M\to N$ is a morphism of left $A$-modules, then 
$\Ker(f)$ exists in the underlying abelian category $\cC$ 
and one checks easily that it has a unique structure of left 
$A$-module which makes it a submodule of $M$. 
{\em If moreover $\otimes$ is right exact\/} when either 
argument is fixed, then also $\Coker\,f$ has a unique 
$A$-module structure for which $N\to\Coker\,f$ is $A$-linear. 
In this case the category of left $A$-modules is abelian. 
Similarly, if $A$ is a unitary $\cC$-monoid, then one 
defines the notion of {\em unitary\/} left $A$-module by 
requiring that 
$\sigma_{M/A}\circ(\underline 1_A\otimes\one_M)=\nu_M$ 
and these form an abelian category when $\otimes$ is right exact.

\sset\subsubsection{}\label{subsec_specialising}
\index{Almost algebra(s)|indref{subsec_specialising}}
\index{$V^a\Alg$, $A\Mod$, $A\Alg$ (for an almost algebra $A$)
|indref{subsec_specialising}}
Specialising to our case we obtain the category $V^a\Alg$ of 
almost $V$-algebras and, for every almost $V$-algebra $A$, the 
category $A\Mod$ of unitary left $A$-modules. Clearly the 
localization functor restricts to 
a functor $V\Alg\to V^a\Alg$ and for any $V$-algebra $R$ 
we have a localization functor $R\Mod\to R^a\Mod$.

Next, if $A$ is an almost $V$-algebra, we can define the category
$A\Alg$ of $A$-algebras. It consists of all the morphisms 
$A\to B$ of almost $V$-algebras. 

\sset\subsubsection{}\label{subsec_def.functor_*}
Let again $(\cC,\otimes,U)$ be any abelian tensor category.
By \cite[p.119]{DeM}, the endomorphism ring $\End_\cC(U)$ of $U$
is commutative. For any object $M$ of $\cC$, denote 
$M_*=\Hom_\cC(U,M)$; then $M\mapsto M_*$ defines a functor
$\cC\to\End_\cC(U)\Mod$. Moreover, if $A$ is a $\cC$-monoid,
$A_*$ is an associative $\End_\cC(U)$-algebra, with multiplication 
given as follows. For $a,b\in A_*$ let 
$a\cdot b:=\mu_A\circ(a\otimes b)\circ\nu_U^{-1}$.
Similarly, if $M$ is an $A$-module, $M_*$ is an $A_*$-module
in a natural way, and in this way we obtain a functor 
from $A$-modules and $A$-linear morphisms to $A_*$-modules
and $A_*$-linear maps. Using \cite[Prop.1.3]{DeM}, one 
can also check that $\End_\cC(U)=U_*$ as $\End_\cC(U)$-algebras,
where $U$ is viewed as a $\cC$-monoid using $\nu_U$. 

\sset\subsubsection{}\label{subsec_almost_element}
\index{Almost module(s)!$M_*$ : almost element(s) of
an|indref{subsec_def.functor_*}{}, \indref{subsec_almost_element}}
All this applies especially to our categories of almost
modules and almost algebras : in this case we call $M\mapsto M_*$
the {\em functor of almost elements\/}. So, if $M$ is an 
almost module, an almost element of $M$ is just an honest 
element of $M_*$. Using \eqref{eq_alm.morph} one can show 
easily that for every $V$-module $M$ the natural map 
$M\to(M^a)_*$ is an almost isomorphism.

\sset\subsubsection{}\label{subsec_almost_homomorph}
\index{$\Alhom_A(M,N)$ : Almost homomorphisms|indref{subsec_almost_homomorph}}
Let $A$ be a $V^a$-algebra. For any two $A$-modules $M,N$, 
the set $\Hom_{A\Mod}(M,N)$ has a natural structure of 
$A_*$-module and we obtain an internal Hom functor by letting
$$\Alhom_A(M,N)=\Hom_{A\Mod}(M,N)^a.$$
This is the functor of {\em almost homomorphisms\/} from $M$ to $N$.

\sset\subsubsection{}\label{subsec_tens.prod.a-mods}
\index{Almost module(s)!$M\otimes_AN$ : tensor product of|indref{subsec_tens.prod.a-mods}}
For any $A$-module $M$ we have also a functor 
of tensor product $M\otimes_A-$ on $A$-modules which, in view
of the following proposition \ref{prop_adjoint} can be shown to be a 
left adjoint to the functor $\Alhom_A(M,-)$. It can be defined as
$M\otimes_AN:=(M_*\otimes_{A_*}N_*)^a$.

With this tensor product, $A\Mod$ is an abelian tensor category 
as well, and $A\Alg$ could also be described as the category of 
($A\Mod$)-algebras. Under this equivalence, a morphism 
$\phi:A\to B$ of almost $V$-algebras becomes the unit morphism 
$\underline 1_B:A\to B$ of the corresponding monoid. We will 
sometimes drop the subscript and write simply $\underline 1$.

\begin{remark}\label{rem_base.comp} Let $V\to W$ be a map of 
base rings, $W$ taken with the extended ideal $\fm W$. 
Then $W^a$ is an almost $V$-algebra so we have defined the 
category $W^a\Mod$ using base ring $V$ and the category 
$(W,\fm W)^a\Mod$ using base $W$. One shows easily 
that they are equivalent: we have an obvious functor 
$(W,\fm W)^a\Mod\to W^a\Mod$ and a quasi-inverse 
is provided by $M\mapsto M_*$. 
Similar base comparison statements hold for the categories 
of almost algebras.
\end{remark} 
\begin{proposition}\label{prop_adjoint}
Let $A$ be a $V^a$-algebra, $R$ a $V$-algebra.
\begin{enumerate}
\item
There is a natural isomorphism $A\simeq A_*^a$ of almost 
$V$-algebras.
\item
The functor $M\mapsto M_*$ from $R^a\Mod$ to $R\Mod$ 
(resp. from $R^a\Alg$ to $R\Alg$) is right adjoint to the 
localization functor $R\Mod\to R^a\Mod$ (resp. $R\Alg\to R^a\Alg$).
\item
The counit of the adjunction $M_*^a\to M$ is a natural 
isomorphism from the composition of the two functors to the 
identity functor $\one_{R^a\Mod}$ (resp. $\one_{R^a\Alg}$).
\end{enumerate}
\end{proposition}

\begin{proof} (i) has already been remarked. We show (ii). 
In light of remark \ref{rem_base.comp} (applied with 
$W=R$) we can assume that $V=R$. Let $M$ be a $V$-module
and $N$ an almost $V$-module; we have natural bijections
$$\begin{array}{r@{\:\simeq\:}l}
\Hom_{V^a\Mod}(M^a,N) & \Hom_{V^a\Mod}(M^a,(N_*)^a) 
\simeq \Hom_V(\tilde\fm\otimes_VM,N_*) \\
& \Hom_V(M,\Hom_V(\tilde\fm,N_*)) 
\simeq \Hom_V(M,\Hom_{V^a\Mod}(V,(N_*)^a)) \\
& \Hom_V(M,N_*)
\end{array}$$
which proves (ii). Now (iii) follows by inspecting 
the proof of (ii), or by \cite[Ch.III Prop.3]{Ga}.
\end{proof}

\begin{remark}\label{rem_tricky}
(i) Let $M_1,M_2$ be two $A$-modules. By proposition 
\ref{prop_adjoint}(iii) it is clear that a morphism 
$\phi:M_1\to M_2$ of $A$-modules is uniquely determined 
by the induced morphism $M_{1*}\to M_{2*}$. On this basis, 
we will very often define morphisms of $A$-modules (or 
$A$-algebras) by saying how they act on almost elements.

(ii) It is a bit tricky to deal with preimages of almost elements 
under morphisms: for instance, if $\phi:M_1\to M_2$ is an 
epimorphism (by which we mean that $\Coker\,\phi\simeq 0$) and 
$m_2\in M_{2*}$, then it is not true in general that we can 
find an almost element $m_1\in M_{1*}$ such that 
$\phi_*(m_1)=m_2$. What remains true is that for arbitrary 
$\eps\in\fm$ we can find $m_1$ such that 
$\phi_*(m_1)=\eps\cdot m_2$. 

(iii) The existence of the right adjoint $M\mapsto M_*$ follows 
also directly from \cite[Chap.III \S 3 Cor.1 or Chap.V \S2]{Ga}. 
\end{remark}

\begin{corollary}\label{cor_coco}
The categories $A\Mod$ and $A\Alg$ are both complete and cocomplete.
\end{corollary}
\begin{proof} We recall that the categories $A_*\Mod$ and 
$A_*\Alg$ are both complete and cocomplete. Now let $I$ be 
any small indexing category and $M:I\to A\Mod$ be any functor. 
Denote by $M_*:I\to A_*\Mod$ the composed functor 
$i\mapsto M(i)_*$. We claim that
$\colim{I} M=(\colim{I} M_*)^a$.
The proof is an easy application of proposition 
\ref{prop_adjoint}(iii). A similar argument also works for limits 
and for the category $A\Alg$.
\end{proof}

\sset\subsubsection{} 
For any $V^a$-algebra $A$, The abelian category $A\Mod$ satisfies 
axiom (AB5) (see {\em e.g.\/} \cite[\S A.4]{We}) and it has a 
generator, namely the object $A$ itself. It then follows by a 
general result that $A\Mod$ has enough injectives.

\begin{corollary} The functor $M\mapsto M_*$ from $R^a\Mod$ 
to $R\Mod$ sends injectives to injectives and injective envelopes
to injective envelopes.
\end{corollary}
\begin{proof} The functor $M\mapsto M_*$ is right adjoint to 
an exact functor, hence it preserves injectives. Now, let $J$ 
be an injective envelope of $M$; to show that $J_*$ is an 
injective envelope of $M_*$, it suffices to show that $J_*$ 
is an essential extension of $M_*$. However, if 
$N\subset J_*$ and $N\cap M_*=0$, then $N^a\cap M=0$, hence 
$\fm N=0$, but $J_*$ does not contain $\fm$-torsion, 
thus $N=0$.
\end{proof}

\sset\subsubsection{}
Note that the essential image of $M\mapsto M_*$ is closed under 
limits. Next recall that the forgetful functor 
$A_*\Alg\to\Set$ 
(resp. $A_*\Mod\to\Set$) has a left adjoint 
$A_*[-]:\Set\to A_*\Alg$ (resp. $A^{(-)}:\Set\to A_*\Mod$)
that assigns to a set $S$ the free $A_*$-algebra $A_*[S]$ 
(resp. the free $A_*$-module $A_*^{(S)}$) generated by $S$. 
If $S$ is any set, it is natural to write $A[S]$ (resp. $A^{(S)}$) 
for the $A$-algebra $(A_*[S])^a$ (resp. for the $A$-module 
$(A_*^{(S)})^a$. This yields a left adjoint, called the 
{\em free $A$-algebra\/} functor $\Set\to A\Alg$ (resp. the 
{\em free $A$-module\/} functor $\Set\to A\Mod$) to the 
``forgetful'' functor $A\Alg\to\Set$ (resp. $A\Mod\to\Set$) 
$B\mapsto B_*$.

\sset\subsubsection{}\label{subsec_define.Mperv}
\index{$\Mperv$|indref{subsec_define.Mperv}}
Now let $R$ be any $V$-algebra; we want to construct a left adjoint 
to the localisation functor $R\Mod\to R^a\Mod$. For a given 
$R^a$-module $M$, let
\set\begin{equation}\label{eq_left.adj}
M_!:=\tilde\fm\otimes_V(M_*).
\end{equation}
We have the natural map (unit of adjunction) $R\to R^a_*$, 
so that we can view $M_!$ as an $R$-module.

\begin{proposition}\label{prop_left.adj} Let $R$ be a $V$-algebra.
\begin{enumerate}
\item
The functor $R^a\Mod\to R\Mod$ defined by \eqref{eq_left.adj} is left 
adjoint to localisation.
\item
The unit of the adjunction $M\to M^a_!$ is a natural isomorphism
from the identity functor $\one_{R^a\Mod}$ to the composition of the 
two functors.
\end{enumerate}
\end{proposition}
\begin{proof} (i) follows easily from \eqref{eq_alm.morph} 
and (ii) follows easily from (i).
\end{proof}

\begin{corollary}\label{cor_left.exact} Suppose that $\tilde\fm$ 
is a flat $V$-module. Then we have : 
\begin{enumerate}
\item
the functor $M\mapsto M_!$ is exact;
\item
the localisation functor $R\Mod\to R^a\Mod$ sends injectives
to injectives.
\end{enumerate}
\end{corollary}
\begin{proof} By proposition \ref{prop_left.adj} it follows
that $M\mapsto M_!$ is right exact. To show that it is also left 
exact when $\tilde\fm$ is a flat $V$-module, it suffices to remark 
that $M\mapsto M_*$ is left exact. Now, by (i), the functor
$M\mapsto M^a$ is right adjoint to an exact functor, so (ii) 
is clear.
\end{proof}

\sset\subsubsection{}\label{subsec_define_Bperv}
\index{$\Bperv$|indref{subsec_define_Bperv}{}, \indref{subsec_adj-for-schemes}}
Let $B$ be any $A$-algebra. The multiplication on
$B_*$ is inherited by $B_!$, which is therefore a non-unital
ring in a natural way. We endow the $V$-module
$V\oplus B_!$ with the ring structure determined by the rule:
$(v,b)\cdot(v',b'):=(v\cdot v',v\cdot b'+v'\cdot b+b\cdot b')$
for all $v,v'\in V$ and $b,b'\in B_!$. Then $V\oplus B_!$ is a 
(unital) ring. We notice that the $V$-submodule generated 
by all the elements of the form 
$(x\cdot y,-x\otimes y\otimes\underline 1)$ (for arbitrary 
$x,y\in\fm$) forms an ideal $I$ of $V\oplus B_!$. 
Set $B_{!!}:=(V\oplus B_!)/I$.
Thus we have a sequence of $V$-modules
\set\begin{equation}\label{eq_left.adj.algebras}
0\to\tilde\fm\to V\oplus B_!\to B_{!!}\to 0
\end{equation}
which in general is only right exact.
\begin{definition}\label{def_exact_algebra}
\index{Almost algebra(s)!exact|indref{def_exact_algebra}}
We say that $B$ is an {\em exact
$V^a$-algebra\/} if the sequence \eqref{eq_left.adj.algebras}
is exact.
\end{definition}
\begin{remark}\label{rem_exact.alg} Notice that if 
$\tilde\fm\stackrel{\sim}{\to}\fm$ ({\em e.g.\/} when $\fm$ is flat), 
then all $V^a$-algebras are exact. In the general case, 
if $B$ is any $A$-algebra, then $V^a\times B$ is always exact.
Indeed, we have $(V^a\times B)_*\simeq V^a_*\times B_*$
and, by remark \ref{rem_almost.zero}(i), 
$\tilde\fm\otimes_VV^a_*\simeq\tilde\fm$.
\end{remark}
Clearly we have a natural isomorphism $B\simeq B_{!!}^a$.
\begin{proposition} The functor $B\mapsto B_{!!}$ is left
adjoint to the localisation functor $A_{!!}\Alg\to A\Alg$.
\end{proposition}
\begin{proof} Let $B$ be an $A$-algebra, $C$ an
$A_{!!}$-algebra and $\phi:B\to C^a$ a morphism of 
$A$-algebras. By proposition \ref{prop_left.adj} we obtain 
a natural $A_*$-linear morphism $B_!\to C$. Together with
the structure morphism $V\to C$ this yields a map 
$\tilde\phi:V\oplus B_!\to C$ which is easily seen to be 
a ring homomorphism. It is equally clear that the ideal 
$I$ defined above is mapped to zero by $\tilde\phi$, hence 
the latter factors through a map of $A_{!!}$-algebras $B_{!!}\to C$.
Conversely, such a map induces a morphism of $A$-algebras
$B\to C^a$ just by taking localisation. It is easy to check
that the two procedures are inverse to each other, which shows
the assertion.
\end{proof}

\begin{remark}\label{rem_adjoints}
(i) The functor of almost elements commutes with arbitrary limits, 
because all right adjoints do. It does not in general commute with 
colimits, not even with arbitrary infinite direct sums.
Dually, the functors $M\mapsto M_!$ and $B\mapsto B_{!!}$ 
commute with all colimits. In particular, the latter commutes
with tensor products.

(ii) Resume the notation of remark \ref{rem_base.comp}.
The change of setup functor $(W,\fm W)^a\Mod\to W^a\Mod$
commutes with the operations $M\mapsto M_*$ and $M\mapsto M_!$.
The corresponding functor $F:(W,\fm W)^a\Alg\to W^a\Alg$
satisfies the identity: $F(B)_{!!}\otimes_{W^a_{!!}}W=B_{!!}$
for every $(W,\fm W)^a$-algebra $B$.
\end{remark}

\subsection{Uniform spaces of almost modules}\label{sec_unif.spaces}
\index{Almost module(s)!$\cM_c(A)$ : $A$-modules generated by
at most $c$ almost elements|indref{sec_unif.spaces}}
Let $A$ be a $V^a$-algebra.
For any cardinal number $c$, we let $\cM_c(A)$ 
be the set of isomorphism classes of $A$-modules which admit
a set of generators of cardinality $\leq c$. In the following
we fix some (very) large infinite cardinality $\omega$,
and suppose that the isomorphism classes of all our $A$-modules
lie in $\cM_\omega(A)$. The choice of $\omega$ is required to 
avoid set-theoretical inconsistencies, but it is immaterial for 
our purposes, so we will henceforth just write $\cM(A)$ instead
of $\cM_\omega(A)$.

\begin{definition}\label{def_unif.on.A}
\index{$\cI_A(M)$, $E_M(\fm_0)$, $E_\cM(\fm_0)$|indref{def_unif.on.A}}
Let $A$ be a $V^a$-algebra and $M$ an $A$-module.
\begin{enumerate}
\item
We define a uniform structure on the set $\cI_A(M)$ of 
$A$-submodules of $M$, as follows. For every finitely 
generated ideal $\fm_0\subset\fm$, the subset of 
$\cI_A(M)\times\cI_A(M)$ given by 
$E_M(\fm_0):=\{(M_0,M_1)~|~\fm_0 M_0\subset M_1~
\text{and}~\fm_0 M_1\subset M_0\}$
is an entourage for the uniform structure, and the subsets 
of this kind form a fundamental system of entourages.
\item
We define a uniform structure on $\cM(A)$ as follows.
For every finitely generated ideal $\fm_0\subset\fm$
and every integer $n\geq 0$ we define the entourage 
$E_\cM(\fm_0)\subset\cM(A)\times\cM(A)$, which consists 
of all pairs of $A$-modules $(M,M')$ such that there 
exist a third module $N$ and morphisms $\phi:N\to M$, 
$\psi:N\to M'$, such that $\fm_0$ annihilates the kernel 
and cokernel of $\phi$ and $\psi$. We declare that the 
$E_\cM(\fm_0)$ form a fundamental system of entourages 
for the uniform structure of $\cM(A)$.
\end{enumerate}
\end{definition}

\begin{remark}\label{rem_alternative} 
Notice that the entourage $E_\cM(\fm_0)$
can be defined equivalently by all the pairs of $A$-modules
$(M,M')$ such that there exists a third module $L$ and
morphisms $\phi':M\to L$, $\psi:M'\to L$ such that $\fm_0$
annihilates the kernel and cokernel of $\phi$ and $\phi$.
Indeed, given a pair $(M,M')\in E_\cM(\fm_0)$, and a datum 
$(N,\phi,\psi)$ as in definition \ref{def_unif.on.A}(ii),
a datum $(L,\phi',\psi')$ satisfying the above condition
is obtained from the push out diagram
\set\begin{equation}\label{eq_pushpull}
{\diagram
N \ar[r]^\phi \ar[d]_\psi & M \ar[d]^{\phi'} \\
M' \ar[r]^{\psi'} & L.
\enddiagram}
\end{equation}
Conversely, given a datum $(L,\phi',\psi')$, one obtains
another diagram as \eqref{eq_pushpull}, by letting $N$
be the fibred product of $M$ and $M'$ over $L$. 
\end{remark}

\sset\subsubsection{}\label{subsec_Cauchy.prod}
\index{Cauchy product|indref{subsec_Cauchy.prod}}
We will also need occasionally a notion of ``Cauchy product'' : 
let $\prod^\infty_{n=0}I_n$ be a formal infinite product of ideals
$I_n\subset A$. We say that the formal product {\em satisfies the 
Cauchy condition\/} (or briefly : {\em is a Cauchy product\/})
if, for every neighborhood $\cU$ of $A$ in $\cI_A(A)$ there exists 
$n_0\ge 0$ such that $\prod_{m=n}^{n+p}I_m\in\cU$\ \ for all 
$n\ge n_0$ and all $p\ge 0$.

\begin{lemma}\label{lem_separated.complete}
Let $M$ be an $A$-module.
\begin{enumerate}
\item 
$\cI_A(M)$ with the uniform structure of definition 
{\em\ref{def_unif.on.A}} is complete and separated.
\item
The following maps are uniformly continuous :
\begin{enumerate}
\item
$\cI_A(M)\times\cI_A(M)\to\cI_A(M)~:~(M',M'')\mapsto M'\cap M''$.
\item
$\cI_A(M)\times\cI_A(M)\to\cI_A(M)~:~(M',M'')\mapsto M'+M''$.
\item
$\cI_A(A)\times\cI_A(A)\to\cI_A(A)~:~(I,J)\mapsto IJ$.
\end{enumerate}
\item
For any $A$-linear morphism $\phi:M\to N$, the following
maps are uniformly continuous:
\begin{enumerate}
\item
$\cI_A(M)\to\cI_A(N)~:~M'\mapsto\phi(M')$.
\item
$\cI_A(N)\to\cI_A(M)~:~N'\mapsto\phi^{-1}(N')$.
\end{enumerate}
\end{enumerate}
\end{lemma}
\begin{proof} (i) : The separation property is easily verified. 
We show that $\cI_A(M)$ is complete. Therefore, suppose that 
$\cF$ is some Cauchy filter of $\cI_A(M)$. Concretely, this 
means that for every finitely generated $\fm_0\subset\fm$, 
there exists $F(\fm_0)\in\cF$ such that 
$\fm_0I\subset J$ for every $I,J\in F(\fm_0)$. Let 
$L:=\bigcup_{F\in\cF}(\bigcap_{I\in F}I)$.
We claim that $L$ is the limit of our filter. Indeed,
for a given finitely generated $\fm_0\subset\fm$, 
we have $\fm_0I\subset\bigcap_{J\in F(\fm_0)}J$, 
for every $I\in F(\fm_0)$, whence 
$\fm_0I\subset L$. On the other hand,
if $I\in F\subset F(\fm_0)$, we can write:
$\fm_0L=\bigcup_{F'\subset F}
\fm_0(\bigcap_{J\in F'}J)\subset
\bigcup_{F'\subset F}(\bigcap_{J\in F'}\fm_0 J)\subset
\bigcup_{F'\subset F}I=I$ (where $F'$ runs over all the
subsets $F'\in\cF$ such that $F'\subset F$). This shows 
that $(L,I)\in E_M(\fm_0)$ whenever $I\in F(\fm_0)$, which 
implies the claim. (ii) and (iii) are easy and will be 
left to the reader.
\end{proof}

\begin{remark} In general, the uniform space $\cM(A)$
is not separated. In view of proposition \ref{prop_best}, 
a counterexample is provided by remark \ref{rem_counter}.
\end{remark}

\begin{lemma}\label{lem_ancora}
Let $\phi:M\to N$ be an $A$-linear morphism, $B$ an $A$-algebra. 
The following maps are uniformly continuous :
\begin{enumerate}
\item 
$\cM(A)\to\cI_A(A)~:~M\mapsto\Ann_A(M)$.
\item
$\cI_A(M)\times\cI_A(N)\to\cM(A)~:~
(M',N')\mapsto(\phi(M')+N')/\phi(M')$.
\item
$\cM(A)\times\cM(A)\to\cM(A)~:~
(M',M'')\mapsto\Alhom_A(M',M'')$.
\item
$\cM(A)\times\cM(A)\to\cM(A)~:~
(M',M'')\mapsto M'\otimes_AM''$.
\item
$\cM(A)\to\cM(B)~:~M\mapsto B\otimes_AM$.
\item
$\cM(A)\to\cM(A)~:~M\mapsto\Lambda_A^rM$\quad for any $r\geq 0$,
provided {\em ({\bf B})} holds.
\end{enumerate}
\end{lemma}
\begin{proof} We show (iv) and leave the others to the reader.
By symmetry, we reduce to verifying that, if 
$(M',M'')\in E_\cM(\fm_0)$ and $N$ is an arbitrary $A$-module,
then $(N\otimes_AM',N\otimes_AM'')\in E_\cM(\fm_0^2)$.
Then we can further assume that there is a morphism 
$\phi:M'\to M''$ with 
$\fm_0\cdot\Ker\,\phi=\fm_0\cdot\Coker\,\phi=0$.
We factor $\phi_!$ as an epimorphism followed by a monomorphism 
$M'_!\stackrel{\phi_1}{\to}\Img\,\phi_!\stackrel{\phi_2}{\to}M''_!$, 
and then we reduce to checking that the kernels and cokernels 
of both $\one_N\otimes_A\phi_1$ and $\one_N\otimes_A\phi_2$ are 
killed by $\fm_0$. This is clear for $\one_N\otimes_A\phi_1$, and 
it follows easily for $\one_N\otimes_A\phi_2$ as well, 
by using the Tor sequences.
\end{proof}

\begin{definition}\label{def_finit_conditions}
\index{Almost module(s)!(almost) finitely generated
|indref{def_finit_conditions}}
\index{Almost module(s)!(almost) finitely presented
|indref{def_finit_conditions}}
\index{Almost module(s)!$\cU_n(A)$ : uniformly almost finitely
generated|indref{def_finit_conditions}}
For a subset $S$ of a topological space
$T$, let $\bar S$ denote the adherence of $S$ in $T$. Let
$M$ be an $A$-module.
\begin{enumerate}
\item
$M$ is said to be {\em finitely generated\/}
if its isomorphism class lies in $\bigcup_{n\in\N}\cM_n(A)$.
\item
$M$ is said to be {\em almost finitely generated\/}
if its isomorphism class lies in $\bar{\bigcup_{n\in\N}\cM_n(A)}$.
\item
$M$ is said to be {\em uniformly almost finitely
generated\/} if there exists an integer $n\geq 0$ such that
the isomorphism class of $M$ lies in $\cU_n(A):=\bar{\cM_n(A)}$. 
Then we will say that $n$ is a {\em uniform bound\/} for $M$.
\item
$M$ is said to be {\em finitely presented\/} if it is
isomorphic to the cokernel of a morphism of free finitely
generated $A$-modules. We denote by $\cF\cP(A)\subset\cM(A)$ 
the subset of the isomorphism classes of finitely presented
$A$-modules.
\item
$M$ is {\em almost finitely presented\/} if its
isomorphism class lies in $\bar{\cF\cP(A)}$.
\end{enumerate}
\end{definition}

\begin{remark} Under condition ({\bf A}), an $A$-module
$M$ lies in $\cU_n(A)$ if and only if, for every $\eps\in\fm$ 
there exists an $A$-linear morphism $A^n\to M$ whose cokernel
is killed by $\eps$.
\end{remark}

\begin{proposition}\label{prop_equiv.cond} 
Let $M$ be an $A$-module.
\begin{enumerate}
\item
$M$ is almost finitely generated if and only if for every 
finitely generated ideal $\fm_0\subset\fm$ there exists a 
finitely generated submodule $M_0\subset M$ such that 
$\fm_0 M\subset M_0$.
\item 
The following conditions are equivalent:
\begin{enumerate}
\item
$M$ is almost finitely presented.
\item
for arbitrary $\eps,\delta\in\fm$ there exist positive integers 
$n=n(\eps)$, $m=m(\eps)$ and a three term complex 
$A^m\stackrel{\psi_\eps}{\to}A^n\stackrel{\phi_\eps}{\to}M$ 
with $\eps\cdot\Coker(\phi_\eps)=0$ and
$\delta\cdot\Ker\,\phi_\eps\subset\Img\,\psi_\eps$.
\item
For every finitely generated ideal $\fm_0\subset\fm$ there 
is a complex $A^m\stackrel{\psi}{\to}A^n\stackrel{\phi}{\to}M$ 
with $\fm_0\cdot\Coker\,\phi=0$ and 
$\fm_0\cdot\Ker\,\phi\subset\Img\,\psi$.
\end{enumerate}
\end{enumerate}
\end{proposition}
\begin{proof}(i): Let $M$ be an almost finitely generated 
$A$-module, and $\fm_0\subset\fm$ a finitely generated 
subideal. Choose a finitely generated subideal 
$\fm_1\subset\fm$ such that $\fm_0\subset\fm_1^3$; by hypothesis, 
there exist $A$-modules $M'$ and $M''$, where $M''$ is finitely 
generated, and morphisms $f:M'\to M$, $g:M'\to M''$ whose kernels 
and cokernels are annihilated by $\fm_1$. We get morphisms
$\fm_1\otimes_V M''\to\Img(g)$ and $\fm_1\otimes_V\Img(g)\to M'$,
hence a composed morphism 
$\phi:\fm_1\otimes_V\fm_1\otimes_V M''\to M'$; it is easy to check
that $\Coker(f\circ\phi)$ is annihilated by $\fm_1^3$, hence 
$M_0:=\Img(f\circ\phi)$ will do.

To show (ii) we will need the following :
\begin{claim}\label{cl_fin.pres}
Let $F_1$ be a finitely generated $A$-module
and suppose that we are given $a,b\in V$ and a (not necessarily 
commutative) diagram
$$\xymatrix{F_1 \ar[r]^p \ar@<.5ex>[d]^\phi & M \\
F_2 \ar@<.5ex>[u]^\psi \ar[ur]_q
}$$
such that $q\circ\phi=a\cdot p$, $p\circ\psi=b\cdot q$. Let
$I\subset V$ be an ideal such that $\Ker\,q$ has a finitely
generated submodule containing $I\cdot\Ker\,q$. Then
$\Ker\,p$ has a finitely generated submodule containing
$ab\cdot I\cdot\Ker\,p$.
\end{claim}
\begin{pfclaim} Let $R$ be the submodule of $\Ker\,q$
given by the assumption. We have 
$\Img(\psi\circ\phi-ab\cdot\one_{F_1})\subset\Ker\,p$
and $\psi(R)\subset\Ker\,p$. We take 
$R_1:=\Img(\psi\circ\phi-ab\cdot\one_{F_1})+\psi(R)$.
Clearly $\phi(\Ker\,p)\subset\Ker\,q$, so 
$I\cdot\phi(\Ker\,p)\subset R$, hence
$I\cdot\psi\circ\phi(\Ker\,p)\subset\psi(R)$ and finally
$ab\cdot I\cdot\Ker\,p\subset R_1$.
\end{pfclaim}

\begin{claim}\label{cl_fin.pres.was.lemma}
If $M$ satisfies condition (b) of the proposition, and 
$\phi:F\to M$ is a morphism with $F\simeq A^n$, then for 
every finitely generated ideal 
$\fm_1\subset\fm\cdot\Ann_V(\Coker\,\phi)$ there is a 
finitely generated submodule of $\Ker\,\phi$ containing 
$\fm_1\cdot\Ker\,\phi$.
\end{claim}
\begin{pfclaim}
Now, let $\delta\in\Ann_V(\Coker\,\phi)$ and 
$\eps_1,\eps_2,\eps_3,\eps_4\in\fm$. By assumption
there is a complex $A^r\stackrel{t}{\to}A^s\stackrel{q}{\to}M$
with $\eps_1\cdot\Coker\,q=0$, 
$\eps_2\cdot\Ker\,q\subset\Img\,t$. Letting $F_1:=F$, $F_2:=A^s$, 
$a:=\eps_1\cdot\eps_3$, $b:=\eps_4\cdot\delta$, one checks easily 
that $\psi$ and $\phi$ can be given such that all the assumptions
of claim \ref{cl_fin.pres} are fulfilled. So, with 
$I:=\eps_2\cdot V$ we see that 
$\eps_1\cdot\eps_2\cdot\eps_3\cdot\eps_4\cdot\delta\cdot
\Ker\,\phi$
lies in a finitely generated submodule of $\Ker\,\phi$.
But $\fm_1$ is contained in an ideal generated by finitely
many such products 
$\eps_1\cdot\eps_2\cdot\eps_3\cdot\eps_4\cdot\delta$.
\end{pfclaim}

Now, it is clear that (c) implies (a) and (b). To show
that (b) implies (c), take a finitely generated ideal 
$\fm_1\subset\fm$ such that $\fm_0\subset\fm\cdot\fm_1$, 
pick a morphism $\phi:A^n\to M$ whose cokernel is annihilated 
by $\fm_1$, and apply claim \ref{cl_fin.pres.was.lemma}.
We show that (a) implies (c). For a given finitely generated 
subideal $\fm_0\subset\fm$, pick another finitely generated
$\fm_1\subset\fm$ such that $\fm_0\subset\fm_1^3$; 
find morphisms $f:M'\to M$ and $g:M'\to M''$ whose 
kernels and cokernels are annihilated by $\fm_1$, and such 
that $M''$ is finitely presented. 
Let $\underline\eps:=(\eps_1,...,\eps_r)$ be a finite sequence 
of generators of $\fm_1$, and denote by 
$K_\bullet:=K_\bullet(\underline\eps)$
the Koszul complex of $V$-modules associated to the sequence
$\underline\eps$. Set $\fm_1':=\Coker(K_2\to K_1)$; we derive
a natural surjection $\partial:\fm_1'\to\fm_1$ and, for every $i=1,...,r$,
maps $e_i:V\to\fm_1'$ such that the compositions
$$
\fm'_1\stackrel{\partial}{\to}\fm_1\hookrightarrow
V\stackrel{e_i}{\to}\fm_1'
\qquad
V\stackrel{e_i}{\to}\fm_1'\stackrel{\partial}{\to}\fm_1
\hookrightarrow V
$$
are both scalar multiplication by $\eps_i$. Hence, for every 
$V^a$-module $N$, the kernel and cokernel of the natural morphism
$\fm_1'\otimes_VN\to N$ are annihilated by $\fm_1$.
Let now $\phi$ be as in the proof of (i); notice that
the diagram:
$$
\xymatrix{
\fm_1\otimes_V\fm_1\otimes_VM' \ar[r] \ar[d] & M' \ar[d] \\
\fm_1\otimes_V\fm_1\otimes_VM'' \ar[r] \ar[ru]^\phi & M''
}
$$
commutes. It follows that the composed morphism 
$$
\fm_1'\otimes_V\fm_1'\otimes_VM''\to
\fm_1\otimes_V\fm_1\otimes_VM''\stackrel{\phi}{\to} M'\to M
$$
has kernel and cokernel annihilated by $\fm_1^3$, so the
claim follows.
\end{proof}

\sset\subsubsection{}\label{subsec_reduce.to.count}
Suppose that $\fm=\bigcup_{\lambda\in\Lambda}\fm_\lambda$,
where $(\fm_\lambda~|~\lambda\in\Lambda)$ is a filtered family of
subideals such that $\fm_\lambda^2=\fm_\lambda$ for every 
$\lambda\in\Lambda$. Let $R$ be a $V$-algebra, $M$ an $R$-module,
and denote by $R^a_\lambda$ (resp. $R^a$) the 
$(V,\fm_\lambda)^a$-algebra (resp. $(V,\fm)$-algebra)
associated to $R$; define similarly $M^a$ and $M^a_\lambda$, 
for all $\lambda\in\Lambda$. 

\begin{lemma} With the notation of \eqref{subsec_reduce.to.count},
the following conditions are equivalent:
\begin{enumerate}
\item
$M^a$ is an almost finitely generated (resp. almost finitely
presented) $R^a$-module.
\item
$M^a_\lambda$ is an almost finitely generated (resp. almost
finitely presented) $R^a_\lambda$-module for all $\lambda\in\Lambda$.
\end{enumerate}
\end{lemma}
\begin{proof} We only give the proof for the case of almost finitely
presented modules, and let the reader spell out the analogous
(and easier) argument for almost finitely generated modules.

(i) $\Rightarrow$ (ii) : for given $\lambda\in\Lambda$
and a finitely generated subideal $\fm_0\subset\fm_\lambda$, 
pick a complex $(R^a)^m\to(R^a)^n\to M^a$ as in proposition
\ref{prop_equiv.cond}(ii.c); after applying termwise the natural
functor $R^a\Mod\to R^a_\lambda\Mod$, we obtain another
complex that satisfies again the condition of 
proposition \ref{prop_equiv.cond}(ii.c), so the claim follows.

(ii) $\Rightarrow$ (i) : for finitely generated $\fm_0\subset\fm$
find $\lambda\in\Lambda$ such that $\fm_0\subset\fm_\lambda$;
then pick a complex 
$(R^a_\lambda)^m\stackrel{\psi}{\to}(R^a_\lambda)^n
\stackrel{\phi}{\to} M^a_\lambda$ as in the foregoing.
Set $N:=M^a_{\lambda*}$, let $N^a$ be the image of $N$
in $R^a\Mod$ and notice that $(M^a,N^a)\in E_{\cM}(\fm_0)$.
Let $\eps_1,...,\eps_k$ be a set of generators for $\fm_0$;
let $\alpha:R^m\to(R^a_{\lambda*})^n$ be the composition of the
adjunction map $R^m\to(R^a_{\lambda*})^m$ and the map $\psi_*$.
Let $e_1,...,e_k$ be a basis of $V^k$ and $f_1,...f_m$ a
basis of $R^m$; we define a map $\beta:V^k\otimes_VR^m\to R^m$ 
by the rule: $e_i\otimes f_j\mapsto\eps_i\cdot f_j$.
The map $\alpha\circ\beta$ lifts to an $R$-linear map
$\gamma:V^k\otimes_VR^m\to R^n$ and $\phi_*$ induces an
$R$-linear map $\Coker\,\gamma\to N$ whose kernel and cokernel
are annihilated by $\fm_0^2\cdot\fm$. All in all, this shows that
$(\Coker\,\gamma^a,M^a)\in E_{\cM}(\fm_0^4)$, so the claim 
follows.
\end{proof}

The following proposition generalises a well-known 
characterization of finitely presented modules over usual 
rings.
\begin{proposition}\label{prop_al.small} 
Let $M$ be an $A$-module.
\begin{enumerate}
\item
$M$ is almost finitely generated if and only if, for every 
filtered system $(N_\lambda,\phi_{\lambda\mu})$ 
(indexed by a directed set $\Lambda$) the natural morphism
\set\begin{equation}\label{eq_was.nu}
M:\colim{\Lambda}\Alhom_A(M,N_\lambda)\to
\Alhom_A(M,\colim{\Lambda}N_\lambda)
\end{equation}
is a monomorphism.
\item
$M$ is almost finitely presented if and only if for every
filtered inductive sytem as above, \eqref{eq_was.nu} is an 
isomorphism.
\end{enumerate}
\end{proposition}
\begin{proof} The ``only if'' part in (i) (resp. (ii)) is first 
checked when $M$ is finitely generated (resp. finitely presented)
and then extended to the general case. We leave the details
to the reader and we proceed to verify the ``if'' part.
For (i), choose a set $I$ and an epimorphism $p:A^{(I)}\to M$.
Let $\Lambda$ be the directed set of finite subsets of $I$, ordered
by inclusion. For $S\in\Lambda$, let $M_S:=p(A^S)$. Then 
$\colim{\Lambda}(M/M_S)=0$, so the assumption gives 
$\colim{\Lambda}\Alhom_A(M,M/M_S)=0$, {\em i.e.\/} 
$\colim{\Lambda}\Hom_A(M,M/M_S)=0$ is almost zero, so, for
every $\eps\in\fm$, the image of $\eps\cdot\one_M$ in the 
above colimit is $0$, {\em i.e.\/} there exists $S\in\Lambda$
such that $\eps M\subset M_S$, which proves the contention. 
For (ii), we present $M$ as a filtered colimit 
$\colim{\Lambda}M_\lambda$, where each $M_\lambda$ is finitely
presented (this can be done {\em e.g.\/} by taking such a
presentation of the $A_*$-module $M_*$ and applying $N\mapsto N^a$).
The assumption of (ii) gives that 
$\colim{\Lambda}\Hom_A(M,M_\lambda)\to\Hom_A(M,M)$ is an almost
isomorphism, hence, for every $\eps\in\fm$ there is 
$\lambda\in\Lambda$ and $\phi_\eps:M\to M_\lambda$ such that 
$p_\lambda\circ\phi_\eps=\eps\cdot\one_M$, where 
$p_\lambda:M_\lambda\to M$ is the natural morphism to the
colimit. If such a $\phi_\eps$ exists for $\lambda$, then 
it exists for every $\mu\geq\lambda$. Hence, if $\fm_0\subset\fm$
is a finitely generated subideal, say $\fm_0=\sum_j^kV\eps_j$,
then there exist $\lambda\in\Lambda$ and $\phi_i:M\to M_\lambda$
such that $p_\lambda\circ\phi_i=\eps_i\cdot\one_M$ for $i=1,...,k$.
Hence $\Img(\phi_i\circ p_\lambda-\eps_i\cdot\one_{M_\lambda})$
is contained in $\Ker\,p_\lambda$ and contains 
$\eps_i\cdot\Ker\,p_\lambda$. Hence $\Ker\,p_\lambda$ has a 
finitely generated submodule $L$ containing
$\fm_0\cdot\Ker\,p_\lambda$. Choose a presentation 
$A^m\to A^n\stackrel{\pi}{\to} M_\lambda$. Then one can 
lift $\fm_0 L$ to a finitely generated submodule 
$L'$ of $A^n$. Then $\Ker(\pi)+L'$ is a finitely generated 
submodule of $\Ker(p_\lambda\circ\pi)$ containing
$\fm_0^2\cdot\Ker(p_\lambda\circ\pi)$. Since we also have 
$\fm_0\cdot\Coker(p_\lambda\circ\pi)=0$ and $\fm_0$ is
arbitrary, the conclusion follows from proposition 
\ref{prop_equiv.cond}.
\end{proof}

\begin{lemma}\label{lem_long.forgotten} 
Let\ $0\to M'\to M\to M''\to 0$ be an exact sequence of 
$A$-modules. Then:
\begin{enumerate}
\item
If $M'$, $M''$ are almost finitely generated (resp. presented)
then so is $M$.
\item
If $M$ is almost finitely presented, then $M''$ is almost
finitely presented if and only if $M'$ is almost finitely
generated.
\end{enumerate}
\end{lemma}
\begin{proof} These facts can be deduced from proposition
\ref{prop_al.small} and remark \ref{rem_flproj}(iii), or proved 
directly.
\end{proof}

\begin{lemma}\label{lem_limfingen}
Let $(M_n~;~\phi_n:M_n\to M_{n+1}~|~n\in\N)$ 
be a direct system of $A$-modules and suppose there exist
sequences $(\eps_n~|~n\in\N)$ and $(\delta_n~|~n\in\N)$ 
of ideals of $V$ such that
\begin{enumerate}
\item
$\limdir{n\to\infty}\eps_n^a=V^a$ (for the uniform 
structure of definition {\em\ref{def_unif.on.A}}) and 
$\prod^\infty_{j=0}\delta_j$ is a Cauchy product (see 
\eqref{subsec_Cauchy.prod};
\item
for all $n\in\N$ there exist integers $N(n)$ and morphisms 
of $A$-modules $\psi_n:A^{N(n)}\to M_n$ such that 
$\eps_n\cdot\Coker\,\psi_n=0$;
\item
$\delta_n\cdot\Coker\,\phi_n=0$ for all $n\in\N$.
\end{enumerate}
Then $\colim{n\in\N}M_n$ is an almost finitely generated 
$A$-module.
\end{lemma}
\begin{proof} Let $M:=\colim{n\in\N}M_n$. For any $n\in\N$ let
$a_n=\bigcap_{m\ge 0}
(\prod_{j=n}^{n+m}\delta_j)$.
Then $\limdir{n\to\infty}a_n=V$.
For $m>n$ set $\phi_{n,m}=
\phi_m\circ...\circ\phi_{n+1}\circ\phi_n:M_n\to M_{m+1}$ and let 
$\phi_{n,\infty}:M_n\to M$ be the natural morphism. An easy 
induction shows that 
$\prod^m_{j=n}\delta_j\cdot\Coker\,\phi_{n,m}=0$
for all $m>n\in\N$. Since 
$\Coker\,\phi_{n,\infty}=\colim{m\in\N}\Coker\,\phi_{n,m}$ we
obtain $a_n\cdot\Coker\,\phi_{n,\infty}=0$ for all $n\in\N$. 
Therefore 
$\eps_n\cdot a_n\cdot\Coker(\phi_{n,\infty}\circ\psi_n)=0$
for all $n\in\N$. Since $\limdir{n\to\infty}\eps_n\cdot a_n=V$, 
the claim follows.
\end{proof}

In the remaining of this section {\em we assume that condition\/}
({\bf B}) {\em of \eqref{subsec_conditions.A.B} holds}. We wish to 
define the Fitting ideals of an arbitrary uniformly almost finitely 
generated $A$-module $M$. This will be achieved in two steps:
first we will see how to define the Fitting ideals
of a finitely generated module, then we will deal with
the general case. We refer to \cite[Ch.XIX]{La} for the 
definition of the Fitting ideals $F_i(M)$ of a finitely 
generated module over an arbitrary ring $R$.

\begin{lemma}\label{lem_Fit.fin.gen} 
Let $R$ be a $V$-algebra and $M$, $N$ two finitely
generated $R$-modules with an isomorphism of
$R^a$-modules $M^a\simeq N^a$.
Then $F_i(M)^a=F_i(N)^a$ for every $i\geq 0$.
\end{lemma}
\begin{proof} By the usual arguments, for every 
$\eps\in\fm$ we have morphisms $\alpha:M\to N$, 
$\beta:N\to M$ with kernels and cokernels 
killed by $\eps^2$. Then we have: $F_i(N)\supset 
F_i(\Img\,\alpha)\cdot F_0(\Coker\,\alpha)$.
If $N$ is generated by $k$ elements, then the same
holds for $\Coker\,\alpha$, whence 
$\Ann_R(\Coker\,\alpha)^k\subset F_0(\Coker\,\alpha)$,
therefore $\eps^{2k}R\subset F_0(\Coker\,\alpha)$, and 
consequently 
$F_i(N)\supset\eps^{2k}F_i(\Img\,\alpha)$.
Since $\Img(\alpha)$ is a quotient of $M$, it is
clear that $F_i(M)\subset F_i(\Img\,\alpha)$, so
finally $\eps^{2k}F_i(M)\subset F_i(N)$.
Arguing symmetrically with $\beta$ one has 
$\eps^{2k}F_i(N)\subset F_i(M)$.
Since we assume {(\bf B)}, the claim follows.
\end{proof}

\sset\subsubsection{}
Let $M$ be a finitely generated $A$-module. In light
of lemma \ref{lem_Fit.fin.gen}, the Fitting ideals
$F_i(M)$ are well defined as ideals in $A$. 

\begin{lemma}\label{lem_uniform.estimate} 
Let $\fm_0\subset\fm$ be a finitely generated subideal 
and $n\in\N$. Pick $\eps_1,...,\eps_k\in\fm$ 
such that $\fm_0\subset(\eps_1^{3n},...,\eps_k^{3n})$ and 
set $\fm_1:=(\eps_1,...,\eps_k)$. Then
$(F_i(M),F_i(M'))\in E_A(\fm_0)$ for every 
$(M,M')\in E_\cM(\fm_1)$ such that $M$ and $M'$ are 
generated by at most $n$ of their almost elements.
\end{lemma}
\begin{proof} Let $M,M'$ be as in the lemma. By
hypothesis, there exist an $A$-module $N$ and morphisms 
$\phi:N\to M$ and $\psi:N\to M'$ such that $\fm_1$ 
annihilates the kernel and cokernel of $\phi$ and $\psi$. 
By symmetry, it suffices to show that 
$\eps_i^{3n}F_i(M)\subset F_i(M')$
for every $i=1,...,k$. Now, for every $i\leq k$, the 
morphism $M\to M$ : $x\mapsto\eps_i\cdot x$ factors
through a morphism $\alpha:M\to\phi(N)$, and similarly, 
scalar multiplication by $\eps_i$ on $N$ factors through
a morphism $\beta:\phi(N)\to N$. Then 
$\eta:=\psi\circ\beta\circ\alpha:M\to M'$ has kernel
and cokernel annihilated by $\eps_i^3$. Pick finitely
generated $A_*$-modules $L\subset M_*$, $L'\subset M'_*$ 
such that $L^a=M$ and $L^{\prime a}=M'$. Replacing
$L'$ by $L'+\eta_*(L)$ we can assume that 
$\eta_*(L)\subset L'$. Then $F_i(M)=F_i(L)^a$, 
$F_i(L')^a=F_i(M')$ and 
$F_i(L')\supset F_i(L/(L\cap\Ker~\eta_*))\cdot 
F_0(L'/\eta_*L)$.
Since $L'/\eta_*(L)$ is generated by at most $n$
elements and is annihilated by $\eps_i^3\cdot\fm$, we have 
$\eps^{3n}_i\cdot\fm\subset F_0(L'/(\eta_*L))$. Furthermore
$F_i(L/(L\cap\Ker~\eta_*))\supset F_i(L)$, so the
claim follows.
\end{proof}

\begin{proposition}\label{prop_uniform.cont.F} 
For every $i,n\geq 0$, the map $F_i:\cM_n(A)\to\cI_A(A)$ 
is uniformly continuous and therefore it extends uniquely 
to a uniformly continuous map $F_i:\cU_n(A)\to\cI_A(A)$.
\end{proposition}
\begin{proof} The uniform continuity follows readily 
from lemma \ref{lem_uniform.estimate}. Since $\cI_A(A)$ 
is complete, it follows that $F_i$ extends to the whole 
of $\cU_n(A)$. Finally, the extension is unique because 
$\cI_A(A)$ is separated.
\end{proof}
\index{Almost module(s)!$\cU_n(A)$ : uniformly almost finitely
generated!$F_i(M)$ : Fitting ideals of a|indref{def_Fitting-ideal}}
\begin{definition}\label{def_Fitting-ideal} Let $M$ be 
a uniformly almost finitely generated $A$-module. We call 
$F_i(M)$ the {\em $i$-th Fitting ideal\/} of $M$.
\end{definition}

\begin{proposition}\label{prop_Fitting.short} {\em (i)} Let 
$0\to M'\stackrel{\phi}{\to}M\stackrel{\psi}{\to}M''\to 0$
be a short exact sequence of uniformly almost finitely
generated $A$-modules. Then:
$\sum_{j+k=i}F_j(M')\cdot F_k(M'')\subset F_i(M)$
for every $i\geq 0$.
\begin{enumerate}
\addtocounter{enumi}{1}
\item
For every uniformly almost finitely generated $A$-module
$M$, any $A$-algebra $B$ and any $i\geq 0$ we have 
$F_i(B\otimes_AM)=F_i(M)\cdot B$.
\end{enumerate}
\end{proposition}
\begin{proof} (i): Let $n$ be uniform bound for $M$ and $M'$;
by remark \ref{rem_alternative} we can find, for every subideal 
$\fm_0\subset\fm$, $A$-modules $M_0$, $M'_0$, $L$, $L'$ and 
morphisms $M\stackrel{\alpha}{\to} L\stackrel{\beta}{\leftarrow} M_0$, 
$M'\stackrel{\alpha'}{\to} L'\stackrel{\beta'}{\leftarrow} M'_0$
whose kernels and cokernels are annihilated by $\fm_0$, and
such that $M_0$ and $M'_0$ are generated by $n$ almost elements.
Let $N$ be defined by the push-out diagram 
$$
\xymatrix{
M' \ar[r]^\phi \ar[d]_{\alpha'} & M \ar[r]^\alpha & L \ar[d]^\gamma \\
L' \ar[rr]^{\gamma'} & & N.}
$$
Furthermore set $M'_1:=\Img(\gamma'\circ\beta':M'_0\to N)$,
$M_1:=\Img((\gamma\circ\beta)\oplus(\gamma'\circ\beta'):
M_0\oplus M'_0\to N)$ and let $M''_1$ be the cokernel of
the induced monomorphism $M'_1\to M_1$. We deduce a commutative
diagram with short exact rows:
\set\begin{equation}\label{eq_triplerow}{
\diagram
0 \ar[r] & M' \ar[r]^\phi \ar[d] & M \ar[r]^\psi \ar[d] &
M'' \ar[d] \ar[r] & 0 \\
0 \ar[r] & \Img\,\gamma' \ar[r] & N \ar[r] & \Coker\,\gamma' 
\ar[r] & 0 \\
0 \ar[r] & M'_1 \ar[u] \ar[r] & M_1 \ar[u] \ar[r] & M''_1 
\ar[u] \ar[r] & 0.
\enddiagram
}\end{equation}
One checks easily that the kernels and cokernels of all the
vertical arrows in \eqref{eq_triplerow} are annihilated by 
$\fm_0^2$, {\em i.e.} $(M,M_1), (M',M'_1), 
(M'',M''_1)\in E_\cM(\fm_0^2)$. Let $x_1,...,x_n\in M_{0*}$ 
(resp. $x'_1,...,x'_n\in M'_{0*}$) be a set of generators for 
$M_0$ (resp. for $M'_0$). 
For every $i=1,...,n$, let $z_i:=\gamma\circ\beta(x_i)$ 
and $z'_i:=\gamma'\circ\beta'(x'_i)$. Let $Q'\subset N_*$
(resp $Q\subset N_*$) be the $A_*$-module generated by
the $z'_i$ (resp. and by the $z_i$). It is clear that
the bottom row of \eqref{eq_triplerow} is naturally
isomorphic to the short exact sequence 
$(0\to Q'\to Q\to Q/Q'\to 0)^a$. It is well-known that
$F_i(Q')\cdot F_j(Q/Q')\subset F_{i+j}(Q)$ for every $i,j\in\N$;
by lemma \ref{lem_separated.complete}(ii.c) and proposition 
\ref{prop_uniform.cont.F} all the operations under considerations
are uniformly continuous, so we deduce 
$F_i(M')\cdot F_j(M/M')\subset F_{i+j}(M)$, which is (i).

(ii): since the identity is known for usual finitely generated
modules over rings, the claim follows easily from proposition 
\ref{prop_uniform.cont.F} and lemma \ref{lem_ancora}(v).
\end{proof}

\subsection{Almost homological algebra}\label{sec_homol} 
In this section we fix an almost $V$-algebra $A$ and we consider 
various constructions in the category of $A$-modules. 

\sset\subsubsection{}
By corollary \ref{cor_coco} any inverse system $(M_n~|~n\in\N)$ of 
$A$-modules has an (inverse) limit $\liminv{n\in\N}M$. As usual,
we denote by $\lim^1$ the right derived functor of the inverse limit
functor. Notice that \cite[Cor.~3.5.4]{We} holds in the almost case 
since axiom (AB4*) holds in $A\Mod$ (on the other hand,
\cite[Lemma 3.5.3]{We} does not hold under (AB4*), (the
proof given there uses elements : for a counterexample in an
exotic abelian category, see \cite{Nee}).

\begin{lemma}\label{lem_invlim}
Let $(M_n~;~\phi_n:M_n\to M_{n+1}~|~n\in\N)$ 
(resp. $(N_n~;~\psi_n:N_{n+1}\to N_n~|~n\in\N)$)
be a direct (resp. inverse) system of $A$-modules and morphisms
and $(\eps_n~|~n\in\N)$ a sequence of ideals of $V^a$ converging 
to $V^a$ (for the uniform structure of definition 
{\em\ref{def_unif.on.A}}).
\begin{enumerate}
\item
If $\eps_n\cdot M_n=0$ for all $n\in\N$ then 
$\colim{n\in\N}M_n\simeq 0$.
\item
If $\eps_n\cdot N_n=0$ for all $n\in\N$ then 
$\liminv{n\in\N}N_n\simeq 0\simeq\liminv{n\in\N}^1N_n$.
\item
If $\eps_n\cdot\Coker\,\psi_n=0$ for all $n\in\N$ and 
$\prod_{j=0}^\infty\eps_j$ is a Cauchy product, then 
$\liminv{n\in\N}^1N_n\simeq 0$.
\end{enumerate}
\end{lemma}

\begin{proof} (i) and (ii) : we remark only that 
$\liminv{n\in\N}^1N_n\simeq\liminv{n\in\N}^1N_{n+p}$ for all
$p\in\N$ and leave the details to the reader.
We prove (iii). From \cite[Cor. 3.5.4]{We} it follows easily
that $(\liminv{n\in\N}^1N_{n*})^a\simeq\liminv{n\in\N}^1N_n$.
It then suffices to show that $\liminv{n\in\N}^1N_{n*}$
is almost zero. We have $\eps^2_n\cdot\Coker\,\psi_{n*}=0$
and the product $\prod_{j=0}^\infty\eps_j^2$ is again a Cauchy 
product. Next let $N'_n:=\bigcap_{p\ge 0}\Img(N_{n+p*}\to N_{n*})$. 
If $J_n:=
\bigcap_{p\ge 0}(\eps_n\cdot\eps_{n+1}\cdot...\cdot\eps_{n+p})^2$
then $J_n N_{n*}\subset N'_n$ and 
$\limdir{n\to\infty}J_n^a=V^a$. In view of (ii),
$\liminv{n\in\N}^1N_{n*}/N'_n$ is almost zero, hence we reduce 
to showing that $\liminv{n\in\N}^1N'_n$ is almost zero.
But
$$
J_{n+p+q}\cdot N'_n\subset\Img(N'_{n+p+q}\to N'_n)
\subset\Img(N'_{n+p}\to N'_n)
$$
for all $n,p,q\in\N$. On the other hand, since the ideals $J_n^a$
converge to $V^a$, we get 
$\bigcup_{q=0}^\infty\fm\cdot J_{n+p+q}=\fm$, hence
$\fm N'_n\subset\Img(N'_{n+p}\to N'_n)$ and finally
$\fm N'_n=\fm^2 N'_n\subset\Img(\fm N'_{n+p}\to\fm N'_n)$
which means that $\{\fm N'_n\}$ is a surjective inverse system,
so its $\liminv{}^1$ vanishes and the result follows.
\end{proof}

\begin{example} Let $(V,\fm)$ be as in example \ref{ex_rings}.
Then every finitely generated ideal in $V$ is principal, so in 
the situation of the lemma we can write $\eps_j=(x_j)$ for some 
$x_j\in V$. Then the hypothesis in (iii) can be stated by saying that
there exists $c\in\N$ such that $x_j\neq 0$ for all $j\geq c$
and the sequence $n\mapsto\prod_{j=c}^n|x_j|$ is Cauchy in 
$\Gamma$.
\end{example}

\begin{definition}\label{def_flat_modules}
\index{Almost module(s)!(faithfully) flat|indref{def_flat_modules}}
\index{Almost module(s)!almost projective|indref{def_flat_modules}}
\index{Almost module(s)!$A\Mod_\mathrm{afpr}$ : almost finitely generated
projective|indref{def_flat_modules}}
Let $M$ be an $A$-module. 
\begin{enumerate}
\item 
We say that $M$ is {\em flat\/} (resp. {\em faithfully
flat\/}) if the functor $N\mapsto M\otimes_A N$, from the 
category of $A$-modules to itself is exact (resp.
exact and faithful). 
\item
We say that $M$ is {\em almost projective\/} if 
the functor $N\mapsto\Alhom_A(M,N)$ is exact.
\end{enumerate}
For euphonic reasons, we will use the expression "almost
finitely generated projective'' to denote an $A$-module
which is almost projective and almost finitely generated.
This convention does not give rise to ambiguities, since
we will never consider projective almost modules : indeed,
the following example \ref{ex_useless} explains why the
categorical notion of projectivity is useless in the setting
of almost ring theory.
\end{definition}

\begin{example}\label{ex_useless}
First of all we remark that the functor $M\mapsto M_!$ preserves
(categorical) projectivity, since it is left adjoint to
an exact functor. Moreover, if $P_!$ is a projective $V$-module
$P$ is a projective $V^a$-module, as one checks easily using
the fact the the functor $M\mapsto M_!$ is right exact.
Hence one has an equivalence from the full subcategory of
projective $V^a$-modules, to the full subcategory of projective
$V$-modules $P$ such that $P=\tilde\fm\otimes_VP$. The latter
condition is equivalent to $P=\fm P$; indeed, as $P$ is flat,
we have $\fm P=\fm\otimes_VP$. 

As an example, suppose that $V$ is local; then every projective
$V$-module is free, so if $\fm\neq V$ the condition $P=\fm P$
implies $P=0$, {\em ergo}, there are no non-trivial projective
$V^a$-modules.
\end{example}

\begin{lemma}\label{lem_annihilate.and.forget} 
Let $M$ be an almost finitely generated
$A$-module and $B$ a flat $A$-algebra. Then 
$\Ann_B(B\otimes_AM)=B\otimes_A\Ann_A(M)$.
\end{lemma}
\begin{proof} Using lemma \ref{lem_ancora}(i),(v) we reduce easily
to the case of a finitely generated $A$-module $M$.
Then, let $x_1,...,x_k\in M_*$ be a finite set of generators
for $M$; we have: $\Ann_A(M)=\Ker(\phi:A\to M^k)$,
where $\phi$ is defined by the rule: 
$a\mapsto(a\cdot x_1,...,a\cdot x_k)$ for every $a\in A_*$.
Since $B$ is flat, we have 
$\Ker(\one_B\otimes_A\phi)\simeq B\otimes_A\Ker\,\phi$,
whence the claim.
\end{proof}

\begin{lemma}\label{lem_used.once} 
Let $\bP$ be one of the properties : ``flat'',
``almost projective'', ``almost finitely generated'', ``almost 
finitely presented''. If $B$ is a $\bP$ $A$-algebra, and
$M$ is a $\bP$ $B$-module, then $M$ is $\bP$ as an $A$-module. 
\end{lemma}
\begin{proof} Left to the reader.
\end{proof}

\sset\subsubsection{}
Let $R$ be a $V$-algebra and $M$ a flat (resp. faithfully flat)
$R$-module (in the usual sense, see \cite[p.45]{Mat}). Then $M^a$ 
is a flat (resp. faithfully flat) $R^a$-module. Indeed, the 
functor $M\otimes_R-$ preserves the Serre subcategory of almost 
zero modules, so by general facts it induces an exact functor 
on the localized categories (cp. \cite[p.369]{Ga}). For the 
faithfullness we have to show that an $R$-module $N$ is almost 
zero whenever $M\otimes_RN$ is almost zero. However, $M\otimes_RN$
is almost zero $\Leftrightarrow$ $M\otimes_R(\fm\otimes_VN)=0$
$\Leftrightarrow$ $\fm\otimes_VN=0$ $\Leftrightarrow$
$N$ is almost zero.
It is clear that $A\Mod$ has enough almost projective (resp. flat) 
objects. 

\sset\subsubsection{}\label{subsec_der.adjunction}
Let $R$ be a $V$-algebra. The localisation functor induces a 
functor $G:\sD(R\Mod)\to\sD(R^a\Mod)$ and, in view of corollary 
\ref{cor_left.exact}, $M\mapsto M_!$ induces a functor 
$F:\sD(R^a\Mod)\to\sD(R\Mod)$. We have a natural isomorphism 
$G\circ F\simeq\one_{\sD(R^a\Mod)}$ and a natural transformation
$F\circ G\to\one_{\sD(R\Mod)}$. These satisfy the triangular 
identities of \cite[p.83]{Ma}  so $F$ is a left adjoint to $G$. 
If $\Sigma$ denotes
the multiplicative set of morphisms in $\sD(R\Mod)$ which induce almost
isomorphisms on the cohomology modules, then the localised category
$\Sigma^{-1}\sD(R\Mod)$ exists (see {\em e.g.\/} \cite[Th.10.3.7]{We})
and by the same argument we get an equivalence of categories
$\Sigma^{-1}\sD(R\Mod)\simeq\sD(R^a\Mod)$.

\sset\subsubsection{}\label{subsec_alExt}
\index{$\Tor_i^A(M,N)$, $\AlExt_A^i(M,N)$|indref{subsec_alExt}}
Given an $A$-module $M$, we can derive the functors  
$M\otimes_A-$ (resp. $\Alhom_A(M,-)$, resp. $\Alhom_A(-,M)$) by taking 
flat (resp. injective, resp. almost projective) resolutions : 
one remarks that bounded above exact complexes of flat (resp. 
almost projective) $A$-modules are acyclic for the functor 
$M\otimes_A-$ (resp. $\Alhom_A(-,M)$) (recall the standard 
argument: if $F_\bullet$ is a bounded above exact complex 
of flat $A$-modules, let $\Phi_\bullet$ be a flat resolution of $M$; 
then $\mathrm{Tot}(\Phi_\bullet\otimes_AF_\bullet)\to 
M\otimes_AF_\bullet$
is a quasi-isomorphism since it is so on rows, and 
$\mathrm{Tot}(\Phi_\bullet\otimes_AF_\bullet)$ is acyclic 
since its columns are; similarly, if $P_\bullet$ is a complex of 
almost projective objects, one considers the double complex 
$\Alhom_A(P_\bullet,J^\bullet)$ where $J^\bullet$ is an 
injective resolution of $M$; cp. \cite[\S 2.7]{We}); then one
uses the construction detailed in \cite[Th.10.5.9]{We}. We 
denote by $\Tor_i^A(M,-)$ (resp. $\AlExt_A^i(M,-)$, resp. 
$\AlExt^i_A(-,M)$) the corresponding derived functors. If $A:=R^a$
for some $V$-algebra $R$, we obtain easily natural isomorphisms
\set\begin{equation}\label{eq_belated}
\Tor^R_i(M,N)^a\simeq\Tor_i^A(M^a,N^a)
\end{equation} 
for all $R$-modules $M,N$. A similar result holds for $\Ext^i_R(M,N)$. 

\begin{remark}\label{rem_flproj}
(i) Clearly, an $A$-module $M$ is flat (resp. almost projective)
if and only if $\Tor^A_i(M,N)=0$ (resp. $\AlExt^i_A(M,N)=0$) 
for all $A$-modules $N$ and all $i>0$. In particular, an
almost projective $A$-module is flat, because for every $\eps\in\fm$
the scalar multiplication by $\eps$ : $M\to M$ factors through
a free module.

(ii) Let $M,N$ be two flat (resp. almost projective) $A$-modules.
Then $M\otimes_AN$ is a flat (resp. almost projective) 
$A$-module and for any $A$-algebra $B$, the $B$-module 
$B\otimes_AM$ is flat (resp. almost projective).

(iii) Resume the notation of proposition \ref{prop_al.small}.
If $M$ is almost finitely presented, then one has also that the 
natural morphism 
$\colim{\Lambda}\AlExt_A^1(M,N_\lambda)\to
\AlExt_A^1(M,\colim{\Lambda}N_\lambda)$ is a monomorphism.
This is deduced from proposition \ref{prop_al.small}(ii), 
using the fact that $(N_\lambda)$ can be injected into an
inductive system $(J_\lambda)$ of injective $A$-modules 
({\em e.g.\/} $J_\lambda=E^{\Hom_A(N_\lambda,E)}$, where 
$E$ is an injective cogenerator for $A\Mod$), and by 
applying $\AlExt$ sequences.
\end{remark}

\begin{lemma}\label{lem_reduce.Tor.to.count} 
Resume the notation of \eqref{subsec_reduce.to.count}.
The following conditions are equivalent:
\begin{enumerate}
\item
$M^a$ is a flat (resp. almost projective) $R^a$-module.
\item
$M^a_\lambda$ is a flat (resp. almost projective) 
$R^a_\lambda$-module for every $\lambda\in\Lambda$.
\end{enumerate}
\end{lemma}
\begin{proof} It is a straightforward consequence of 
\eqref{eq_belated} and its analogue for $\Ext^i_R(M,N)$.
\end{proof}

\begin{lemma}\label{lem_Exts.are.same} Let $R$ be a $V$-algebra.
\begin{enumerate}
\item
There is a natural isomorphism: 
$\Ext^i_R(M_!,N)\simeq\Ext^i_{R^a}(M,N^a)$ for every $R^a$-module 
$M$, every $R$-module $N$ and every integer $i\in\N$.
\item
If $P$ is an almost projective $R^a$-module, then 
$\mathrm{hom.dim}_RP_!\leq\mathrm{hom.dim_V}\tilde\fm$.
\end{enumerate}
\end{lemma}
\begin{proof} (here $\mathrm{hom.dim}$ denotes homological
dimension). (i) is a straightforward consequence of the
existence of the adjunction $(F,G)$ of \eqref{subsec_der.adjunction}.
Next we consider, for arbitrary $R$-modules $M$ and $N$, 
the spectral sequence:
$$
E^{p,q}_2:=\Ext_V^p(\tilde\fm,\Ext^q_R(M,N))\Rightarrow
\Ext^{p+q}_R(\tilde\fm\otimes_VM,N).
$$
(this spectral sequence
is constructed {\em e.g.\/} from the double complex 
$\Hom_V(F_p,\Hom_R(F'_q,N))$ where
$F_\bullet$ (resp. $F'_\bullet$) is a projective resolution
of $\tilde\fm$ (resp. of $M$)). If now we let $M:=P_!$,
we deduce from (i) that $\Ext_V^p(\tilde\fm,\Ext^q_R(P_!,N))\simeq
\Ext_{V^a}^p(V^a,\Ext^q_R(P_!,N)^a)\simeq 0$ for every $p\in\N$ 
and every $q>0$. Since $\tilde\fm\otimes_VP_!\simeq P_!$,
assertion (ii) follows easily.
\end{proof}

\begin{lemma}\label{lem_thetwo}
Let $M$ be an almost finitely generated $A$-module. Then 
$M$ is almost projective if and only if, for arbitrary 
$\eps\in\fm$, there exist $n(\eps)\in\N$ and $A$-linear 
morphisms
\set\begin{equation}\label{eq_uv}
M \stackrel{u_\eps}{\longrightarrow}A^{n(\eps)}
\stackrel{v_\eps}{\longrightarrow}M 
\end{equation}
such that $v_\eps\circ u_\eps=\eps\cdot\one_M$.
\end{lemma}
\begin{proof} Let morphisms as in \eqref{eq_uv} be given.
Pick any $A$-module $N$ and apply the functor 
$\AlExt^i_A(-,N)$ to \eqref{eq_uv} to get morphisms
$$\AlExt^i_A(M,N)\to\AlExt^i_A(A^{n(\eps)},N)\to
\AlExt^i_A(M,N)$$
whose composition is again the scalar multiplication by
$\eps$; hence $\eps\cdot\AlExt^i_A(M,N)=0$ for all $i>0$.
Since $\eps$ is arbitrary, it follows from remark 
\ref{rem_flproj}(i) that $M$ is almost projective.
Conversely, suppose that $M$ is almost projective; by 
hypothesis, for arbitrary  $\eps\in\fm$ we can find 
$n:=n(\eps)$ and a morphism $\phi_\eps:A^n\to M$ such that
$\eps\cdot\Coker\,\phi_\eps=0$. Let $M_\eps$ be the image of
$\phi_\eps$, so that $\phi_\eps$ factors as 
$A^{n(\eps)}\stackrel{\psi_\eps}{\longrightarrow}M_\eps 
\stackrel{j_\eps}{\longrightarrow}M$.
Also $\eps\cdot\one_M:M\to M$ factors as 
$M\stackrel{\gamma_\eps}{\longrightarrow}M_\eps
\stackrel{j_\eps}{\longrightarrow}M.$
Since by hypothesis $M$ is almost projective, the natural 
morphism induced by $\psi_\eps$ :
$$
\Alhom_A(M,A^n)\stackrel{\psi_\eps^*}{\longrightarrow}
\Alhom_A(M,M_\eps)
$$
is an epimorphism. Then for arbitrary $\delta\in\fm$ the 
morphism $\delta\cdot\gamma_\eps$ is in the image of 
$\psi_\eps^*$, in other words, there exists an $A$-linear 
morphism $u_{\eps\delta}:M\to A^n$ such that 
$\psi_\eps\circ u_{\eps\delta}=\delta\cdot\gamma_\eps$. 
If now we take $v_{\eps\delta}:=\phi_\eps$, it is clear that 
$v_{\eps\delta}\circ u_{\eps\delta}=\eps\cdot\delta\cdot\one_M$. 
This proves the claim.
\end{proof}

\begin{lemma}\label{lem_Tor.Ext}
Let $R$ be any ring, $M$ any $R$-module and 
$C:=\Coker(\phi:R^n\to R^m)$ any finitely presented (left) 
$R$-module. Let $C':=\Coker(\phi^*:R^m\to R^n)$ be the cokernel 
of the transpose of the map $\phi$. Then there is a natural 
isomorphism
$$\Tor^R_1(C',M)\simeq\Hom_R(C,M)/\Img(\Hom_R(C,R)\otimes_RM).$$
\end{lemma}
\begin{proof} We have a spectral sequence :
$$E^2_{ij}:=\Tor^R_i(H_j(\Cone\,\phi^*),M)\Rightarrow
H_{i+j}(\Cone(\phi^*)\otimes_RM).$$ 
On the other hand we have also natural isomorphisms 
$$\Cone(\phi^*)\otimes_RM\simeq 
\Hom_R(\Cone\,\phi,R)[1]\otimes_RM\simeq 
\Hom_R(\Cone\,\phi,M)[1].$$
Hence :
$$\begin{array}{r@{\:\simeq\:}l}
E^2_{10}\simeq E^\infty_{10}\simeq 
H_1(\Cone(\phi^*)\otimes_RM)/E^\infty_{01} &
H^0(\Hom_R(\Cone\,\phi,M))/\Img(E^2_{01}) \\
& \Hom_R(C,M)/\Img(\Hom_R(C,R)\otimes_RM)
\end{array}$$
which is the claim.
\end{proof}

\begin{proposition}\label{prop_converse}
Let $A$ be a $V^a$-algebra.
\begin{enumerate}
\item
Every almost finitely generated projective
$A$-module is almost finitely presented.
\item
Every almost finitely presented flat $A$-module is almost 
projective.
\end{enumerate}
\end{proposition}
\begin{proof} (ii) : let $M$ be such an $A$-module. Let 
$\eps,\delta\in\fm$ and pick a three term complex
$$
A^m\stackrel{\psi}{\longrightarrow}A^n
\stackrel{\phi}{\longrightarrow}M
$$ 
such that
$\eps\cdot\Coker\,\phi=\delta\cdot\Ker(\phi)/\Img(\psi)=0$. 
Set $P:=\Coker\,\psi_*$; this is a finitely presented
$A_*$-module and $\phi_*$ factors through a morphism
$\bar\phi_*:P\to M_*$.  Let $\gamma\in\fm$; from lemma 
\ref{lem_Tor.Ext} we see that $\gamma\cdot\bar\phi$ is
the image of some element 
$\sum_{j=1}^n\phi_j\otimes m_j\in\Hom_{A_*}(P,A_*)\otimes_{A_*}M_*$.
If we define $L:=A_*^n$ and $v:P\to L$, $w:L\to M_*$ by
$v(x):=(\phi_1(x),...,\phi_n(x))$ and 
$w(y_1,...,y_n):=\sum^n_{j=1}y_j\cdot m_j$, then clearly 
$\gamma\cdot\bar\phi=w\circ v$. Let $K:=\Ker\,\bar\phi_*$. 
Then $\delta\cdot K^a=0$ and
the map $\delta\cdot\one_{P^a}$ factors through a morphism 
$\sigma:(P/K)^a\to P^a$. Similarly the map $\eps\cdot\one_M$
factors through a morphism $\lambda:M\to(P/K)^a$.
Let $\alpha:=v^a\circ\sigma\circ\lambda:M\to L^a$
and $\beta:=w^a:L^a\to M$. The reader can check that 
$\beta\circ\alpha=\eps\cdot\delta\cdot\gamma\cdot\one_M$.
By lemma \ref{lem_thetwo} the claim follows.

(i) : let $P$ be such an almost finitely generated 
projective $A$-module. For any finitely generated 
ideal $\fm_0\subset\fm$ pick a morphism $\phi:A^r\to P$
such that $\fm_0\cdot\Coker\,\phi=0$.
If $\eps_1,...,\eps_k$ is a set of generators for $\fm_0$,
a standard argument shows that, for any $i\le k$, 
$\eps_i\cdot\one_P$ lifts to a morphism 
$\psi_i:P\to A^r/\Ker\,\phi$; then, since $P$ is almost 
projective, $\eps_j\psi_i$ lifts to a morphism 
$\psi_{ij}:P\to A^r$. Now claim \ref{cl_fin.pres} applies 
with $F_1:=A^r$, $F_2:=M=P$, $p:=\phi$, $q:=\one_P$ and 
$\psi:=\psi_{ij}$ and shows that $\Ker\,\phi$ has a finitely 
generated submodule $M_{ij}$ containing 
$\eps_i\cdot\eps_j\cdot\Ker\,\phi$. Then the span of all 
such $M_{ij}$ is a finitely generated submodule 
of $\Ker\,\phi$ containing $\fm_0^2\cdot\Ker\,\phi$. By 
proposition \ref{prop_equiv.cond}(ii), the claim follows.
\end{proof}

In general a flat almost finitely generated $A$-module
is not necessarily almost finitely presented, but one can 
give the following criterion, which extends 
\cite[\S1, Exerc.13]{BouAH}.

\begin{proposition}\label{prop_crit.projectivity} 
If $A\to B$ is a monomorphism of $V^a$-algebras
and $M$ is an almost finitely generated flat $A$-module such that
$B\otimes_AM$ is almost finitely presented over $B$, then $M$
is almost finitely presented over $A$.
\end{proposition}
\begin{proof} Let $\fm_0$ be a finitely generated subideal of
$\fm$ and an $A$-linear morphism $\phi:A^n\to M$ such that 
$\fm_0\cdot\Coker\,\phi=0$.
By assumption and claim \ref{cl_fin.pres.was.lemma}, we can find
a finitely generated $B$-submodule $R$ of
$\Ker(\one_B\otimes_A\phi)$ such that 
\set\begin{equation}\label{eq_onepunch}
\fm_0^2\cdot\Ker(\one_B\otimes_A\phi)\subset R.
\end{equation} 
By a Tor sequence we have $\fm_0\cdot\Coker
(\one_B\otimes_A\Ker(\phi)\to\Ker(\one_B\otimes_A\phi))=0$, hence
$\fm\cdot\fm_0\cdot\Coker
(\one_{B_*}\otimes_{A_*}\Ker(\phi)_*\to\Ker(\one_B\otimes_A\phi)_*)=0$,
therefore there exists a finitely generated submodule $R_0$ of
$\Ker\,\phi$ such that 
\set\begin{equation}\label{eq_twopunch}
\fm_0^2 R\subset B\cdot\Img(R_0\to\Ker(\one_B\otimes_A\phi)).
\end{equation}
By lemma \ref{lem_Tor.Ext}, for every $\eps\in\fm$ the morphism
$\eps\cdot\bar\phi:A^n/R_0\to M$ factors through a morphism 
$\psi:A^n/R_0\to F$, where $F$ is a finitely generated
free $A$-module. Since $F\subset B\otimes_AF$, we deduce
easily $\Ker(A^n/R_0\to B^n/B\cdot R_0)\subset\Ker\,\psi$;
on the other hand, by \eqref{eq_onepunch} and \eqref{eq_twopunch}
we derive 
$\fm_0^4\cdot\Ker\,\phi\subset\Ker(A^n/R_0\to B^n/B\cdot R_0)$.
Thus $\psi$ factors through $M':=A^n/(R_0+\fm_0^4\cdot\Ker\,\phi)$.
Clearly $\fm_0^4\cdot\Ker(M'\to M)=\fm_0\cdot\Coker(M'\to M)=0$;
hence, for every $A$-module $N$, the kernel of the induced morphism
\set\begin{equation}\label{eq_zeroExt}
\AlExt^1_A(M,N)\to\AlExt^1_A(M',N)
\end{equation}
is annihilated by $\fm_0^5$; however \eqref{eq_zeroExt} factors 
through $\AlExt^1_A(F,N)=0$, therefore $\fm_0^5\cdot\AlExt^1_A(M,N)=0$
for every $A$-module $N$.
This shows that $M$ is almost projective, which is equivalent 
to the conclusion, in view of proposition \ref{prop_converse}(i).
\end{proof}

\begin{definition}\label{def_dual_module}
\index{Almost module(s)!$M^*$ : dual of an|indref{def_dual_module}}
\index{Almost module(s)!$\ev_{M/A}$ : evaluation morphism of an|indref{def_dual_module}}
\index{Almost module(s)!$\cE_{M/A}$ : evaluation ideal|indref{def_dual_module}}
\index{Almost module(s)!reflexive|indref{def_dual_module}}
\index{Almost module(s)!invertible|indref{def_dual_module}}
\index{Almost module(s)!$\cE_{M/A}(f)$ : evaluation ideal of an almost
  element in an|indref{def_dual_module}}
Let $M$ be an $A$-module, $f:A\to M$ an almost element of $M$. 
\begin{enumerate}
\item
    The {\em dual $A$-module\/} of $M$ is the $A$-module 
    $M^*:=\Alhom_A(M,A)$. 
\item
    The {\em evaluation morphism\/} is the morphism
    $\ev_{M/A}:M\otimes_AM^*\to A~:~m\otimes\phi\mapsto\phi(m)$.
\item
    The {\em evaluation ideal\/} of $M$ is the ideal
    $\cE_{M/A}:=\Img\,\ev_{M/A}$.
\item
    The {\em evaluation ideal\/} of $f$ is the ideal
    $\cE_{M/A}(f):=\Img(\ev_{M/A}\circ(f\otimes_A\one_{M^*}))$.
\item
    We say that $M$ is {\em reflexive\/} if the natural morphism 
\set\begin{equation}\label{eq_reflex}
M\to (M^*)^*\qquad m\mapsto (f\mapsto f(m))
\end{equation}
is an isomorphism of $A$-modules. 
\item
   We say that $M$ is {\em invertible\/} if 
   $M\otimes_AM^*\simeq A$.
\end{enumerate}
\end{definition}
\begin{remark}\label{rem_B.struct} Notice that if $B$ is an 
$A$-algebra and $M$ any $B$-module, then by ``restriction of 
scalars'' $M$ is also an $A$-module and the dual $A$-module 
$M^*$ has a natural structure of $B$-module. 
This is defined by the rule $(b\cdot f)(m):=f(b\cdot m)$ 
($b\in B_*$, $m\in M_*$ and $f\in M^*_*$). 
With respect to this structure \eqref{eq_reflex} becomes a 
$B$-linear morphism. Incidentally, notice that the two 
meanings of ``$M^*_*$'' coincide, {\em i.e.\/} 
$(M_*)^*\simeq (M^*)_*$.
\end{remark}

\sset\subsubsection{}\label{subsec_define_omega}
\index{$\omega_{P/A}$|indref{subsec_define_omega}}
If $E$, $F$ and $N$ are $A$-modules, there is a natural 
morphism :
\set\begin{equation}\label{eq_three.mods}
E\otimes_A\Alhom_A(F,N)\to\Alhom_A(F,E\otimes_AN).
\end{equation}
Let $P$ be an $A$-module. As a special case of \eqref{eq_three.mods} 
we have the morphism: 
$$\omega_{P/A}:P\otimes_AP^*\to\End_A(P)^a:=\Alhom_A(P,P)$$
such that $\omega_{P/A}(p\otimes\phi)(q):=p\cdot\phi(q)$ for every
$p,q\in P_*$ and $\phi:P\to A$.

\begin{proposition}\label{prop_eval.ideal}
Let $P$ be an almost projective $A$-module. 
\begin{enumerate}
\item
For every morphism of algebras $A\to B$ we have
$\cE_{B\otimes_AP/B}=\cE_{P/A}\cdot B$.
\item
$\cE_{P/A}=\cE_{P/A}^{~2}$.
\item
$P=0$ if and only if $\cE_{P/A}=0$.
\item
$P$ is faithfully flat if and only if $\cE_{P/A}=A$. 
\item
$\cE_{P/A}(f)$ is the smallest of the ideals $J\subset A$
such that $f\in(JP)_*$.
\end{enumerate}
\end{proposition}
\begin{proof} Pick an indexing set $I$ large enough, and an 
epimorphism $\phi:F:=A^{(I)}\to P$. For every $i\in I$ we have 
the standard morphisms 
$A\stackrel{e_i}{\to}F\stackrel{\pi_i}{\to}A$
such that $\pi_i\circ e_j=\delta_{ij}\cdot\one_A$ and 
$\sum_{i\in I}e_i\circ\pi_i=\one_F$. For every $x\in\fm$ choose 
$\psi_x\in\Hom_A(P,F)$ such that $\phi\circ\psi_x=x\cdot\one_P$.
It is easy to check that $\cE_{P/A}$ is generated by the almost
elements $\pi_i\circ\psi_x\circ\phi\circ e_j$ 
($i,j\in I$, $x\in\fm$). (i) follows already. For (iii), the 
``only if'' is clear; if $\cE_{P/A}=0$, then $\psi_x\circ\phi=0$ 
for all $x\in\fm$, hence $\psi_x=0$ and therefore $x\cdot\one_P=0$, 
{\em i.e.} $P=0$. Next, notice that (i) and (iii) imply 
$P/\cE_{P/A}P=0$, {\em i.e.} $P=\cE_{P/A} P$,
so (ii) follows directly from the definition of $\cE_{P/A}$.
Since $P$ is flat, to show (iv) we have only to verify that
the functor $M\mapsto P\otimes_AM$ is faithful. To this purpose, 
it suffices to check that $P\otimes_AA/J\neq 0$ for every 
proper ideal $J$ of $A$. This follows easily from (i) and (iii).
Finally, it clear that $\cE_{P/A}(f)\subset J$ for every ideal
$J\subset A$ with $f\in(JP)_*$. Conversely, for given 
$\eps\in\fm$ pick a sequence 
$P\stackrel{u}{\to} A^n\stackrel{v}{\to} P$ as in \eqref{eq_uv};
we have $u(f)=(a_1,...,a_n)$ with $a_1,...,a_n\in\cE_{P/A}(f)_*$.
(v) follows from the identity 
$\eps f=\sum_{i=1}^na_i\cdot v(e_i)$ (where $e_1,...,e_n$
is the standard basis of the free $A_*$-module $A_*^n$).
\end{proof}

\begin{lemma}\label{lem_three.mods}
Let $E$, $F$, $N$ be three $A$-modules.
\begin{enumerate}
\item
The morphism \eqref{eq_three.mods} is an isomorphism in the 
following cases :
\begin{enumerate}
\item
when $E$ is flat and $F$ is almost finitely presented;
\item
when either $E$ or $F$ is almost finitely generated projective;
\item
when $F$ is almost projective and $E$ is almost finitely
presented;
\item
when $E$ is almost projective and $F$ is almost finitely
generated.
\end{enumerate}
\item
The morphism \eqref{eq_three.mods} is a monomorphism in the 
following cases :
\begin{enumerate}
\item
when $E$ is flat and $F$ is almost finitely generated;
\item
when $E$ is almost projective.
\end{enumerate}
\item
The morphism \eqref{eq_three.mods} is an epimorphism when
$F$ is almost projective and $E$ is almost finitely generated.
\end{enumerate}
\end{lemma}
\begin{proof} If $F\simeq A^{(I)}$ for some finite set $I$, then 
$\Alhom_A(F,N)\simeq N^{(I)}$ and the claims are obvious. More 
generally, if $F$ is almost finitely generated projective, 
for any $\eps\in\fm$ there exists a finite set $I:=I(\eps)$ and 
morphisms 
\set\begin{equation}\label{eq_CD1}
F\stackrel{u_\eps}{\longrightarrow}A^{(I)}
\stackrel{v_\eps}{\longrightarrow}F
\end{equation}
such that $v_\eps\circ u_\eps=\eps\cdot\one_F$. We apply the natural 
transformation 
$$E\otimes_A\Alhom_A(-,N)\to\Alhom_A(-,E\otimes_AN)$$ 
to \eqref{eq_CD1} : an easy diagram chase allows then to conclude 
that the kernel and cokernel of \eqref{eq_three.mods} are killed 
by $\eps$. As $\eps$ is arbitrary, it follows that 
\eqref{eq_three.mods} is an isomorphism in this case. An analogous 
argument works when $E$ is almost finitely generated projective, 
so we get (i.b).
If $F$ is only almost projective, then we still have morphisms
of the type \eqref{eq_CD1}, but now $I(\eps)$ is no longer necessarily
finite. However, the cokernels of the induced morphisms 
$\one_E\otimes u_\eps$ and $\Alhom_A(v_\eps,E\otimes_AN)$ are
still annihilated by $\eps$. Hence, to show (iii) (resp. (i.c)) 
it suffices to consider the case when $F$ is free and $E$ is almost 
finitely generated (resp. presented). By passing to almost elements, 
we can further reduce to the analogous question for usual rings 
and modules, and by the usual juggling we can even replace $E$ 
by a finitely generated (resp. presented) $A_*$-module and $F$ 
by a free $A_*$-module. This case is easily dealt with, and 
(iii) and (i.c) follow. Case (i.d) (resp. (ii.b)) is similar : 
one considers almost elements and replaces $E_*$ by a free 
$A_*$-module (resp. and $F_*$ by a finitely generated 
$A_*$-module). In case (ii.a) (resp. (i.a)), for every finitely 
generated submodule $\fm_0$ of $\fm$ we can find, by proposition
\ref{prop_equiv.cond}, a finitely generated (resp. presented) 
$A$-module $F_0$ and a morphism $F_0\to F$ whose kernel and 
cokernel are annihilated by $\fm_0$. It follows easily that 
we can replace $F$ by $F_0$ and suppose that $F$ is finitely 
generated (resp. presented). Then the argument in 
\cite[Ch.I \S 2 Prop.10]{BouAC} can be taken over 
{\em verbatim\/} to show (ii.a) (resp. (i.a)).
\end{proof}

\begin{lemma}\label{lem_alhom} Let $B$ be an $A$-algebra.
\begin{enumerate}
\item
Let $P$ be an $A$-module. If either $P$ or $B$ is almost 
finitely generated projective as an $A$-module, 
the natural morphism 
\set\begin{equation}\label{eq_claim.tens}
B\otimes_A\Alhom_A(P,N)\to\Alhom_B(B\otimes_AP,B\otimes_AN)
\end{equation}
is an isomorphism for all $A$-modules $N$.
\item
Every almost finitely generated projective $A$-module
is reflexive.
\item
If $P$ is an almost finitely generated projective
$B$-module, the natural morphism
\set\begin{equation}\label{eq_natural.transf}
\Alhom_B(P,B)\otimes_B\Alhom_A(B,A)\to\Alhom_A(P,A)
\quad:\quad\phi\otimes\psi\mapsto\psi\circ\phi
\end{equation}
is an isomorphism of $B$-modules.
\end{enumerate}
\end{lemma}
\begin{proof} (i) is an easy consequence of lemma 
\ref{lem_three.mods}(i.b). To prove (ii), we apply the 
natural transformation \eqref{eq_reflex} to \eqref{eq_CD1} : 
by diagram chase one sees that the kernel and cokernel of 
the morphism $F\to(F^*)^*$ are killed by $\eps$. (iii) is 
analogous : one applies the natural transformation
\eqref{eq_natural.transf} to \eqref{eq_CD1}.
\end{proof}

\begin{lemma}\label{lem_limproj} Let 
$(M_n~;~\phi_n:M_n\to M_{n+1}~|~n\in\N)$ be a direct 
system of $A$-modules and suppose there exist sequences 
$(\eps_n~|~n\in\N)$ and $(\delta_n~|~n\in\N)$ of 
ideals of $V$ such that
\begin{enumerate}
\item
$\limdir{n\to\infty}\eps_n=V$ and 
$\prod^\infty_{j=0}\delta_j$ is a Cauchy product 
(see \eqref{subsec_Cauchy.prod});
\item
$\eps_n\cdot\AlExt^i_A(M_n,N)=
\delta_n\cdot\AlExt^i_A(\Coker\,\phi_n,N)=0$ for all
$A$-modules $N$, all $i>0$ and all $n\in\N$;
\item
$\delta_n\cdot\Ker\,\phi_n=0$ for all $n\in\N$.
\end{enumerate}
Then $\colim{n\in\N}M_n$ is an almost projective $A$-module.
\end{lemma}
\begin{proof} Let $M=\colim{n\in\N}M_n$. By remark 
\ref{rem_flproj}(i) it suffices to show that 
$\AlExt^i_A(M,N)$ vanishes for all $i>0$ 
and all $A$-modules $N$. The maps $\phi_n$ define
a map $\phi:\oplus_nM_n\to\oplus_nM_n$ such that we have a short 
exact sequence 
$0\to\oplus_nM_n\stackrel{\one-\phi}{\longrightarrow}\oplus_nM_n
\longrightarrow M\to 0$. Applying the long exact $\AlExt$ sequence
one obtains a short exact sequence (cp. \cite[3.5.10]{We})
$$0\to\liminv{n\in\N}^1\AlExt_A^{i-1}(M_n,N)\to
\AlExt^i_A(M,N)\to\liminv{n\in\N}\AlExt^i_A(M_n,N)\to 0.$$
Then lemma \ref{lem_invlim}(ii) implies that $\AlExt^i_A(M,N)\simeq 0$ 
for all $i>1$ and moreover  $\AlExt^1_A(M,N)$ is isomorphic to 
$\liminv{n\in\N}^1\Alhom_A(M_n,N)$. Let 
$$\phi^*_n:\Alhom_A(M_{n+1},N)\to\Alhom_A(M_n,N)\qquad 
f\mapsto f\circ\phi_n$$
be the transpose of $\phi_n$ and write $\phi_n$ as a composition 
$M_n\stackrel{p_n}{\longrightarrow}\Img(\phi_n)
\stackrel{q_n}{\hookrightarrow}M_{n+1}$,
so that $\phi^*_n=q_n^*\circ p_n^*$, the composition of the 
respective transposed morphims. We have monomorphisms
$$\begin{array}{rl}
&\Coker\,p^*_n\hookrightarrow\Alhom_A(\Ker\,\phi_n,N) \\
&\Coker\,q^*_n\hookrightarrow\AlExt^1_A(\Coker\,\phi_n,N)
\end{array}$$
for all $n\in\N$. Hence $\delta^2_n\cdot\Coker\,\phi^*_n=0$
for all $n\in\N$. Since $\prod^\infty_{n=0}\delta^2_n$ is a 
Cauchy product, lemma \ref{lem_invlim}(iii) shows that 
$\liminv{n\in\N}^1\Alhom_A(M_n,N)\simeq 0$ and the assertion 
follows.
\end{proof}

\begin{proposition}\label{prop_comp.supp}
Suppose that $\tilde\fm$ is a flat $V$-module. Then for 
any $V$-algebra $R$ the functor $M\mapsto M_!$ commutes with 
tensor products and takes flat $R^a$-modules to 
flat $R$-modules.
\end{proposition}
\begin{proof} Let $M$ be a flat $R^a$-module and  
$N\hookrightarrow N'$ an injective map of $R$-modules.
Denote by $K$ the kernel of the induced map 
$M_!\otimes_RN\to M_!\otimes_RN'$; we have $K^a\simeq 0$.
We obtain an exact sequence
$0\to\tilde\fm\otimes_VK\to
\tilde\fm\otimes_VM_!\otimes_RN\to
\tilde\fm\otimes_VM_!\otimes_RN'$.
But one sees easily that $\tilde\fm\otimes_VK=0$ and 
$\tilde\fm\otimes_VM_!\simeq M_!$, which shows that $M_!$ is a flat 
$R$-module. Similarly, let $M,N$ be two $R^a$-modules. 
Then the natural map $M_*\otimes_RN_*\to(M\otimes_{R^a}N)_*$
is an almost isomorphism and the assertion follows from
remark \ref{rem_almost.zero}(i).
\end{proof}

\subsection{Almost homotopical algebra}\label{sec_hot}

The formalism of abelian tensor categories provides a minimal 
framework wherein the rudiments of deformation theory can be
developed. 

\sset\subsubsection{}
Let $(\cC,\otimes,U)$ be an abelian tensor category; we assume
henceforth that $\otimes$ is a right exact functor. Let $A$ be a 
given $\cC$-monoid. Then, for any two-sided ideal $I$ of $A$, the 
quotient $A/I$ in the underlying abelian category $\cC$\ \ has a 
unique $\cC$-monoid structure such that $A\to A/I$ is a morphism
of monoids. $A/I$ is unitary if $A$ is. 
If $I$ is a two-sided ideal of $A$ such that $I^2=0$, then,
using the right exactness of $\otimes$ one checks that $I$ 
has a natural structure of an $A/I$-bimodule, unitary when $A$ is.
\begin{definition}\label{def_C-extensions}
\index{$\cC$ : Category(ies)!tensor!$\bExmon_\cC(B,I)$, $\bExun_\cC$, 
$\bExal_\cC$ : extensions of a monoid (or an algebra) in 
a |indref{def_C-extensions}{}, \indref{subsec_variant.ext}}
A {\em $\cC$-extension\/} of a $\cC$-monoid 
$B$ by a $B$-bimodule $I$ is a short exact sequence of objects 
of\/ $\cC$
\set\begin{equation}\label{eq_def.ext}
X:\qquad 0\to I\to C\stackrel{p}{\to} B\to 0
\end{equation}
such that $C$ is a $\cC$-monoid, $p$ is a morphism of\/
$\cC$-monoids, $I$ is a square zero two-sided ideal in $C$
and the $E/I$-bimodule structure on $I$ coincides with
the given bimodule structure on $I$. The $\cC$-extensions 
form a category $\bExmon_\cC$. The morphisms are commutative 
diagrams with exact rows
$$\diagram
~X:\dto & 
0 \rto & I \rto \dto^f & E \rto^p \dto^g & B \rto \dto^h & 0 \\
~X': & 0 \rto & I' \rto & E' \rto^{p'} & B' \rto & 0
\enddiagram$$
such that $g$ and $h$ are morphisms of\/ $\cC$-monoids. We 
let $\bExmon_\cC(B,I)$ be the subcategory of\/ $\bExmon_\cC$ 
consisting of all $\cC$-extensions of $B$ by $I$, where the 
morphisms are all short exact sequences as above such that 
$f:=\one_I$ and $h:=\one_B$.
\end{definition}

\sset\subsubsection{}\label{subsec_variant.ext}
We have also the variant in which all the $\cC$-monoids 
in \eqref{eq_def.ext} are required to be unitary (resp. 
to be algebras) and $I$ is a unitary $B$-bimodule (resp. 
whose left and right $B$-module actions coincide, {\em i.e.\/} 
are switched by composition with the ``commutativity 
constraints'' $\theta_{B|I}$ and $\theta_{I|B}$, see 
\eqref{subsec_comm.constr}); we will call $\bExun_\cC$ (resp. 
$\bExal_\cC$) the corresponding category.

\sset\subsubsection{}\label{subsec_more.ext.stuff}
\index{$\cC$ : Category(ies)!tensor!$\Exmon_\cC(B,I)$,
  $\Exal_\cC(B,I)$ : isomorphism classes of extensions in
a|indref{subsec_more.ext.stuff}{}, \indref{subsec_again.exts}}
\index{$\cC$ : Category(ies)!tensor!$X*\phi$, $\psi*X$ : push-out, 
pull-back of extensions in a|indref{subsec_more.ext.stuff}}
For a morphism $\phi:C\to B$ of $\cC$-monoids, and a $\cC$-extension 
$X$ in $\bExmon_\cC(B,I)$, we can pullback $X$ via $\phi$ to obtain
an exact sequence $X*\phi$ with a morphism $\phi^*:X*\phi\to X$;
one checks easily that there exists a unique structure of
$\cC$-extension on $X*\phi$ such that $\phi^*$ is a morphism
of $\cC$-extension; then $X*\phi$ is an object in $\bExmon_\cC(C,I)$. 
Similarly, given a $B$-linear morphism $\psi:I\to J$, we 
can push out $X$ and obtain a well defined object $\psi*X$ 
in $\bExmon_\cC(B,J)$ with a morphism $X\to\psi*X$ of $\bExmon_\cC$.
In particular, if $I_1$ and $I_2$ are two $B$-bimodules,
the functors $p_i*$ ($i=1,2$) associated to the natural 
projections $p_i:I_1\oplus I_2\to I_i$ establish an equivalence
of categories
\set\begin{equation}\label{eq_exal.p_i}
\bExmon_\cC(B,I_1\oplus I_2)\stackrel{\sim}{\to}
\bExmon_\cC(B,I_1)\times\bExmon_\cC(B,I_2)
\end{equation}
whose quasi-inverse is given by 
$(E_1,E_2)\mapsto (E_1\oplus E_2)*\delta$, where 
$\delta:B\to B\oplus B$ is the diagonal morphism.
A similar statement holds for $\bExal$ and $\bExun$.
These operations can be used to induce an abelian group 
structure on the set $\Exmon_\cC(B,I)$ of isomorphism 
classes of objects of $\bExmon_\cC(B,I)$ as follows. For any 
two objects $X,Y$ of $\bExmon_\cC(B,I)$ we can form $X\oplus Y$ 
which is an object of $\bExmon_\cC(B\oplus B,I\oplus I)$. Let  
$\alpha:I\oplus I\to I$ be the addition morphism of $I$. Then we 
set $X+Y:=\alpha*(X\oplus Y)*\delta$. One can check that $X+Y\simeq Y+X$ 
for any $X,Y$ and that the trivial split $\cC$-extension $B\oplus I$ 
is a neutral element for $+$. Moreover every isomorphism class 
has an inverse $-X$. The functors $X\mapsto X*\phi$ and 
$X\mapsto\psi*X$ commute with the operation thus defined, and 
induce group homomorphisms
$$\begin{array}{l}
*\phi:\Exmon_\cC(B,I)\to\Exmon_\cC(C,I) \\
\psi*:\Exmon_\cC(B,I)\to\Exmon_\cC(B,J).
\end{array}$$

\sset\subsubsection{}\label{subsec_again.exts}
\index{$\bExal_A$, $A^\dagger\UniMod$|indref{subsec_again.exts}{},
\indref{subsec_A-dagger}}
We will need the variant $\Exal_\cC(B,I)$ defined in the same 
way, starting from the category $\bExal_\cC(B,I)$. For instance, 
if $A$ is an almost algebra (resp. a commutative ring), we can 
consider the abelian tensor category $\cC=A\Mod$. 
In this case the $\cC$-extensions will be called 
simply $A$-extensions, and we will write $\bExal_A$ rather 
than $\bExal_\cC$. In fact the commutative unitary case will 
soon become prominent in our work, and the more general setup 
is only required for technical reasons, in the proof of 
proposition \ref{prop_lift.idemp} below, which is the abstract 
version of a well-known result on the lifting of idempotents
over nilpotent ring extensions. 

\sset\subsubsection{}\label{subsec_A-dagger}
Let $A$ be a $\cC$-monoid. We form the biproduct 
$A^\dagger:=U\oplus A$ in $\cC$. We denote by $p_1$, $p_2$ 
the associated projections from $A^\dagger$ onto $U$ and
respectively $A$. Also, let $i_1$, $i_2$ be the natural
monomorphisms from $U$, resp. $A$ to $A^\dagger$. $A^\dagger$
is equipped with a unitary monoid structure 
$$\mu^\dagger:=i_2\circ\mu\circ(p_2\otimes p_2)+
i_2\circ\ell^{-1}_A\circ(p_1\otimes p_2)+
i_2\circ r^{-1}_A\circ(p_2\otimes p_1)+
i_1\circ u^{-1}\circ(p_1\otimes p_1)$$
where $\ell_A$, $r_A$ are the natural isomorphisms 
provided by \cite[Prop. 1.3]{DeM} and $u:U\to U\otimes U$ 
is as in \cite[\S 1]{DeM}. In terms of the ring
$A^\dagger_*\simeq U_*\oplus A_*$ this is the multiplication
$(u_1,b_1)\cdot(u_2,b_2):=
(u_1\cdot u_2,b_1\cdot b_2+b_1\cdot u_2+u_1\cdot b_2)$.
Then $i_2$ is a morphism of monoids and one verifies that
the ``restriction of scalars'' functor $i_2^*$ defines an 
equivalence from the category $A^\dagger\UniMod$ of unitary 
$A^\dagger$-modules to the category $A\Mod$ of all 
$A$-modules; let $j$ denote the inverse functor. A similar 
discussion applies to bimodules. 

\sset\subsubsection{}
Similarly, we derive equivalences of categories
$$\xymatrix{
\bExun_\cC(A^\dagger,j(M))\ar@<.5ex>[r]^-{*i_2} &
\bExmon_\cC(A,M)\ar@<.5ex>[l]^-{(-)^\dagger}
}$$
for all $A$-bimodules $M$.

\sset\subsubsection{}
Next we specialise to $A:=U$ : for a given $U$-module $M$ 
let $e_M:=\sigma_{M/U}\circ\ell_M:M\to M$; working out the 
definitions one finds that the condition that $(M,\sigma_{M/U})$ 
is a module structure is equivalent to $e_M^2=e_M$. Let 
$U\times U$ be the product of $U$ by itself in the category 
of $\cC$-monoids. There is an isomorphism of unitary 
$\cC$-monoids $\zeta:U^\dagger\to U\times U$
given by $\zeta:=i_1\circ p_1+i_2\circ p_1+i_2\circ p_2$.
Another isomorphism is $\tau\circ\zeta$, where $\tau$
is the flip $i_1\circ p_2+i_2\circ p_1$. Hence we get
equivalences of categories
$$\xymatrix{
U\Mod \ar@<.5ex>[r]^-{j} & U^\dagger\UniMod 
\ar@<.5ex>[rr]^-{(\zeta^{-1})^*} 
\ar@<.5ex>[l]^-{i_2^*}&&
(U\times U)\UniMod
\ar@<.5ex>[ll]^-{(\tau\circ\zeta)^*}.
}$$
The composition 
$i_2^*\circ(\zeta^{-1}\circ\tau\circ\zeta)^*\circ j$
defines a self-equivalence of $U\Mod$ which associates
to a given $U$-module $M$ the new $U$-module 
$M^\flip$ whose underlying object in $\cC$ is $M$ 
and such that $e_{M^\flip}=\one_M-e_M$.
The same construction applies to $U$-bimodules and finally
we get equivalences 
\set\begin{equation}\label{eq_flip}
\bExmon_\cC(U,M)\stackrel{\sim}{\to}\bExmon_\cC(U,M^\flip)
\qquad X\mapsto X^\flip
\end{equation}
for all $U$-bimodules $M$. If 
$X:=(0\to M\to E\stackrel{\pi}{\to}U\to 0)$
is an extension and $X^\flip:=(0\to M^\flip\to E^\flip\to U\to 0)$,
then one verifies that there is a natural isomorphism $X^\flip\to X$ 
of complexes in $\cC$ inducing $-\one_M$ on $M$, the identity on 
$U$ and carrying the multiplication morphism on $E^\flip$ to
$$-\mu_E+\ell_E^{-1}\circ(\pi\otimes\one_E)+
r_E^{-1}\circ(\one_E\otimes\pi):E\otimes E\to E.$$
In terms of the associated rings, this corresponds to 
replacing the given multiplication $(x,y)\mapsto x\cdot y$ 
of $E_*$ by the new operation
$(x,y)\mapsto\pi_*(x)\cdot y+\pi_*(y)\cdot x-x\cdot y$.
\begin{lemma}\label{lem_split.ext}
If $M$ is a $U$-bimodule whose left and right 
actions coincide, then every extension of $U$ by $M$ splits
uniquely.
\end{lemma}
\begin{proof} Using the idempotent $e_M$ we get a $U$-linear
decomposition $M\simeq M_1\oplus M_2$ where the bimodule 
structure on $M_1$ is given by the zero morphisms and the
bimodule structure on $M_2$ is given by $\ell_M^{-1}$
and $r_M^{-1}$. We have to prove that $\bExmon_\cC(U,M)$
is equivalent to a one-point category. By \eqref{eq_exal.p_i} 
we can assume that $M=M_1$ or $M=M_2$. By \eqref{eq_flip} we
have $\bExmon_\cC(U,M_2)\simeq\bExmon_\cC(U,M_2^\flip)$
and on $M_2^\flip$ the bimodule actions are the zero morphisms.
So it is enough to consider $M=M_1$. In this case,
if $X:=(0\to M\to E\to U\to 0)$ is any extension, 
$\mu_E:E\otimes E\to E$ factors through a morphism
$U\otimes U\to E$ and composing with $u:U\to U\otimes U$
we get a right inverse of $E\to U$, which shows that
$X$ is the split extension. Then it is easy to see that
$X$ does not have any non-trivial automorphisms, which
proves the assertion.
\end{proof}
\begin{proposition}\label{prop_lift.idemp}
Let $X:=(0\to I\to A\stackrel{p}{\to} A'\to 0)$ be a 
$\cC$-extension. 
\begin{enumerate}
\item
Let $e'\in A'_*$ be an idempotent element 
whose left action on the $A'$-bimodule $I$ coincides 
with its right action. Then there exists a unique 
idempotent $e\in A_*$ such that $p_*(e)=e'$.
\item
Especially, if $A'$ is unitary and $I$ is a unitary
$A'$-bimodule, then every extension of $A'$ by $I$ is unitary.
\end{enumerate}
\end{proposition}
\begin{proof} (i) : the hypothesis $e^{\prime 2}=e'$ implies that 
$e':U\to A'$ is a morphism of (non-unitary) $\cC$-monoids. 
We can then replace $X$ by $X*e'$ and thereby assume that
$A'=U$, $p:A\to U$ and $I$ is a (non-unitary) $U$-bimodule 
and the right and left actions on $I$ coincide. 
The assertion to prove is that $\underline 1_U$ lifts to
a unique idempotent $e\in A_*$. However, this follows easily 
from lemma \ref{lem_split.ext}. To show (ii), we observe that,
by (i), the unit $\underline 1_{A'}$ of $A'_*$ lifts uniquely
to an idempotent $e\in A_*$. We have to show that $e$ is a unit
for $A_*$. Let us show the left unit property. Via $e:U\to A$
we can view the extension $X$ as an exact sequence of left 
$U$-modules. We can then split $X$ as the direct sum 
$X_1\oplus X_2$ where $X_1$ is a sequence of unitary $U$-modules
and $X_2$ is a sequence of $U$-modules with trivial actions.
But by hypothesis, on $I$ and on $A$ the $U$-module structure
is unitary, so $X=X_1$ and this is the left unit property.
\end{proof}

\sset\subsubsection{}\label{subsec_so-much}
So much for the general nonsense; we now return to almost algebras.
As already announced, {\em from here on, we assume throughout that
$\tilde\fm$ is a flat $V$-module}.
As an immediate consequence of proposition \ref{prop_lift.idemp} 
we get natural equivalences of categories 
\set\begin{equation}\label{eq_equi.prod}
\bExal_{A_1}(B_1,M_1)\times\bExal_{A_2}(B_2,M_2)\stackrel{\sim}{\to}
\bExal_{A_1\times A_2}(B_1\times B_2,M_1\oplus M_2)
\end{equation}
whenever $A_1$, $A_2$ are $V^a$-algebras, $B_i$ is 
a $A_i$-algebra and $M_i$ is a (unitary) $B_i$-module, $i=1,2$.

\sset\subsubsection{}
Notice that, if $A:=R^a$ for some $V$-algebra $R$, $S$ (resp. 
$J$) is a $R$-algebra (resp. an $S$-module) and $X$ is any object 
of $\bExal_R(S,J)$, then by applying termwise the localisation 
functor we get an object $X^a$ of $\bExal_A(S^a,J^a)$.
With this notation we have the following lemma.
\begin{lemma}\label{lem_eq.Ext} 
Let $B$ be any $A$-algebra and $I$ a $B$-module. 
\begin{enumerate}
\item
The natural functor
\set\begin{equation}\label{eq_funct.ext}
\bExal_{A_{!!}}(B_{!!},I_*)\to\bExal_A(B,I)\qquad X\mapsto X^a
\end{equation}
is an equivalence of categories.
\item
The equivalence \eqref{eq_funct.ext} induces a group isomorphism
$\Exal_{A_{!!}}(B_{!!},I_*)\stackrel{\sim}{\to}\Exal_A(B,I)$
functorial in all arguments.
\end{enumerate}
\end{lemma}
\begin{proof} Of course (ii) is an immediate consequence of (i).
To show (i), let $X:=(0\to I\to E\to B\to 0)$ be any object of
$\bExal_A(B,I)$. Using corollary \ref{cor_left.exact} one sees 
easily that the sequence $X_!:=(0\to I_!\to E_{!!}\to B_{!!}\to 0)$ 
is right exact; $X_!$ won't be exact in general, unless $B$
(and therefore $E$) is an exact algebra. In any case, the 
kernel of $I_!\to E_{!!}$ is almost zero, so we get an extension
of $B_{!!}$ by a quotient of $I_!$ which maps to $I_*$. In
particular we get by pushout an extension $X_{!*}$ by $I_*$,
{\em i.e.\/} an object of $\bExal_{A_{!!}}(B_{!!},I_*)$ and 
in fact the assignment $X\mapsto X_{!*}$ is a quasi-inverse 
for the functor \eqref{eq_funct.ext}.
\end{proof}
\begin{remark}\label{rem_subseq} 
\index{$\Iperv$|indref{rem_subseq}}
By inspecting the proof, 
we see that one can replace $I_*$ by $I_{!*}:=\Img(I_!\to I_*)$ 
in (i) and (ii) above. When $B$ is exact, also $I_!$ will do.
\end{remark}

In \cite[II.1.2]{Il} it is shown how to associate to any ring
homomorphism $R\to S$ a natural simplicial complex of $S$-modules
denoted $\L_{S/R}$ and called the cotangent complex of $S$ over
$R$.

\begin{definition}\label{def_cotangent_cpx}
\index{Almost algebra(s)!$\L_{B/A}$ : cotangent complex of an|indref{def_cotangent_cpx}}
Let $A\to B$ be a morphism of almost 
$V$-algebras. The {\em almost cotangent complex\/} of $B$ over
$A$ is the simplicial $B_{!!}$-module 
$$\L_{B/A}:=B_{!!}\otimes_{(V^a\times B)_{!!}}
\L_{(V^a\times B)_{!!}/(V^a\times A)_{!!}}.$$
\end{definition}

\sset\subsubsection{}\label{subsec_usually}
\index{$\hExt^p_R(E,F)$|indref{subsec_usually}}
Usually we will want to view $\L_{B/A}$ as an object of the derived 
category $\sD(s.B_{!!}\Mod)$ of simplicial $B_{!!}$-modules. 
Indeed, the hyperext functors computed in this category relate 
the cotangent complex to a number of important invariants. Recall
that, for any simplicial ring $R$ and any two $R$-modules $E,F$ 
the hyperext of $E$ and $F$ is the abelian group defined as
$$\hExt^p_R(E,F):=
\colim{n\ge -p}\Hom_{\sD(R\Mod)}(\sigma^nE,\sigma^{n+p}F)$$
(where $\sigma$ is the suspension functor of \cite[I.3.2.1.4]{Il}).

Let us fix an almost algebra $A$. First we want to establish 
the relationship with differentials.

\begin{definition}\label{def_derivation}
\index{Almost algebra(s)!$\Der_A(B,M)$ : derivations of an|indref{def_derivation}}
\index{Almost algebra(s)!$\Omega_{B/A}$ : relative differentials of an
|indref{def_derivation}}
\index{$I_{B/A}$|indref{def_derivation}}
\index{$V^a\AlgMorph$, $V^a\AlgMod$|indref{def_derivation}}
Let $B$ be any $A$-algebra, $M$ any $B$-module. 
\begin{enumerate}
\item
An {\em $A$-derivation\/} of $B$ with values in $M$ 
is an $A$-linear morphism $\partial:B\to M$ such that 
$\partial(b_1\cdot b_2)=b_1\cdot\partial(b_2)+
b_2\cdot\partial(b_1)$
for $b_1,b_2\in B_*$. The set of all $M$-valued
$A$-derivations of $B$ forms a $V$-module $\Der_A(B,M)$ and 
the almost $V$-module $\Der_A(B,M)^a$ has a natural structure 
of $B$-module.
\item
We reserve the notation $I_{B/A}$ for the ideal 
$\Ker(\mu_{B/A}:B\otimes_AB\to B)$. The {\em module 
of relative differentials\/} of $\phi$ is defined as the 
(left) $B$-module $\Omega_{B/A}:=I_{B/A}/I^2_{B/A}$. 
It is endowed with a natural $A$-derivation 
$\delta:B\to\Omega_{B/A}$ defined by 
$b\mapsto\underline 1\otimes b-b\otimes\underline 1$ for all 
$b\in B_*$. The assignment $(A\to B)\mapsto\Omega_{B/A}$ 
defines a functor
$$\Omega:V^a\AlgMorph\to V^a\AlgMod$$
from the category of morphisms $A\to B$ of almost $V$-algebras
to the category denoted $V^a\AlgMod$, consisting of all pairs $(B,M)$ 
where $B$ is an almost $V$-algebra and $M$ is a $B$-module.
The morphisms in $V^a\AlgMorph$ are the commutative squares; the 
morphisms $(B,M)\to(B',M')$ in $V^a\AlgMod$ are all pairs $(\phi,f)$
where $\phi:B\to B'$ is a morphism of almost $V$-algebras and
$f:B'\otimes_BM\to M'$ is a morphism of $B'$-modules.
\end{enumerate}
\end{definition}

\sset\subsubsection{}
The module of relative differentials enjoys the familiar 
universal properties that one expects. In particular
$\Omega_{B/A}$ represents the functor $\Der_A(B,-)$, {\em i.e.}
for any (left) $B$-module $M$ the morphism
\set\begin{equation}\label{eq_repres}
\Hom_B(\Omega_{B/A},M)\to\Der_A(B,M)\qquad f\mapsto f\circ\delta
\end{equation}
is an isomorphism. As an exercise, the reader can supply the proof
for this claim and for the following standard proposition.
\begin{proposition}\label{prop_standard}
Let $B$ and $C$ be two $A$-algebras. 
\begin{enumerate}
\item
There is a natural isomorphism: 
$$\Omega_{C\otimes_AB/C}\simeq C\otimes_A\Omega_{B/A}.$$
\item
Suppose that $C$ is a $B$-algebra. Then there is a natural exact 
sequence of $C$-modules:
$$C\otimes_B\Omega_{B/A}\to\Omega_{C/A}\to\Omega_{C/B}\to 0.$$
\item
Let $I$ be an ideal of $B$ and let $C:=B/I$ be the quotient 
$A$-algebra. Then there is a natural exact sequence:
$I/I^2\to C\otimes_B\Omega_{B/A}\to\Omega_{C/A}\to 0$.
\item
The functor $\Omega:V^a\AlgMorph\to V^a\AlgMod$ commutes with all
colimits.\qed
\end{enumerate}
\end{proposition}

We supplement these generalities with one more statement 
which is in the same vein as lemma \ref{lem_uniform.estimate} 
and which will be useful in section \ref{sec_val.rings} to 
calculate the Fitting ideals of modules of differentials.

\begin{lemma}\label{lem_Omega.converge} 
Let $\phi:B\to B'$ be a morphism of $A$-algebras such 
that $I\cdot\Ker(\phi)=I\cdot\Coker(\phi)=0$ for an 
ideal $I\subset A$. Let 
$d\phi:\Omega_{B/A}\otimes_BB'\to\Omega_{B'/A}$ be the
natural morphism. Then $I\cdot\Coker\,d\phi=0$ and
$I^4\cdot\Ker\,d\phi=0$. 
\end{lemma}
\begin{proof} We will use the standard presentation
\set\begin{equation}\label{eq_present_Omega}
H(B/A)~:~B\otimes_AB\otimes_AB\stackrel{\partial}{\to} 
B\otimes_AB\stackrel{d}{\to}\Omega_{B/A}\to 0
\end{equation}
where $d$ is defined by : $b_1\otimes b_2\mapsto b_1\cdot db_2$
and $\partial$ is the differential of the Hochschild
complex :
$$
b_1\otimes b_2\otimes b_3\mapsto 
b_1b_2\otimes b_3-b_1\otimes b_2b_3+b_1b_3\otimes b_2.
$$
By naturality of $H(B)$, we deduce a morphism
of complexes : $B'\otimes_BH(B/A)\to H(B'/A)$. Then,
by snake lemma, we derive an exact sequence :
$\Ker(\one_{B'}\otimes_A\phi)\to\Ker\,d\phi\to X$, where 
$X$ is a quotient of $B'\otimes_A\Coker(\phi\otimes_A\phi)$.
Using the Tor exact sequences we see that 
$\Ker(\one_{B'}\otimes_A\phi)$ is annihilated by $I^2$.
It follows easily that $I^4$ annihilates $\Ker\,d\phi$.
Similarly, $\Coker\,d\phi$ is a quotient of 
$\Coker(\one_{B'}\otimes_A\phi)$, so $I\cdot\Coker\,d\phi=0$.
\end{proof}

\begin{lemma}\label{lem_differ} For any $A$-algebra $B$ 
there is a natural isomorphism of $B_{!!}$-modules
$$(\Omega_{B/A})_!\simeq\Omega_{B_{!!}/A_{!!}}.$$
\end{lemma}
\begin{proof} Using the adjunction \eqref{eq_repres}
we are reduced to showing that the natural map
$$\phi_M:\Der_{A_{!!}}(B_{!!},M)\to\Der_A(B,M^a)$$ 
is a bijection for all $B_{!!}$-modules $M$. Given 
$\partial:B\to M^a$ we construct $\partial_!:B_!\to M^a_!\to M$. 
We extend $\partial_!$ to $V\oplus B_!$ by setting it equal 
to zero on $V$. Then it is easy to check that the 
resulting map descends to $B_{!!}$, hence giving an  
$A$-derivation $B_{!!}\to M$.
This procedure yields a right inverse $\psi_M$ to 
$\phi_M$. To show that $\phi_M$ is injective, suppose
that $\partial:B_{!!}\to M$ is an almost zero $A$-derivation.
Composing with the natural $A$-linear map $B_!\to B_{!!}$
we obtain an almost zero map $\partial':B_!\to M$. But 
$\fm\cdot B_!=B_!$, hence $\partial'=0$. This implies
that in fact $\partial=0$, and the assertion follows.
\end{proof}

\begin{proposition}\label{prop_H_0}
Let $M$ be a $B$-module. There exists a natural 
isomorphism of $B_{!!}$-modules
$$\hExt^0_{B_{!!}}(\L_{B/A},M_!)\simeq\Der_A(B,M).$$
\end{proposition}
\begin{proof} To ease notation, set $\tilde A:=V^a\times A$
and $\tilde B:=V^a\times B$. We have natural isomorphisms :
$$\begin{array}{l@{\:\simeq\:}ll}
\hExt^0_{B_{!!}}(\L_{B/A},M_!) &
\hExt^0_{\tilde B_{!!}}(\L_{\tilde B_{!!}/\tilde A_{!!}},M_!)
& \quad\text{by \cite[I.3.3.4.4]{Il}} \\
& \Der_{\tilde A_{!!}}(\tilde B_{!!},M_!)
& \quad\text{by \cite[II.1.2.4.2]{Il}} \\
& \Der_{\tilde A}(\tilde B,M) 
& \quad\text{by lemma \ref{lem_differ}.} 
\end{array}$$
But it is easy to see that the natural map
$\Der_A(B,M)\to\Der_{\tilde A}(\tilde B,M)$ is an 
isomorphism.
\end{proof}
\begin{theorem}\label{th_main}
There is a natural isomorphism
\set\begin{equation}\label{eq_th-main}
\Exal_A(B,M)\stackrel{\sim}{\to}\hExt^1_{B_{!!}}(\L_{B/A},M_!).
\end{equation}
\end{theorem}
\begin{proof} With the notation of the proof of proposition
\ref{prop_H_0} we have natural isomorphisms 
$$\begin{array}{l@{\:\simeq\:}ll}
\hExt^1_{B_{!!}}(\L_{B/A},M_!) 
& \hExt^1_{\tilde B_{!!}}(\L_{\tilde B_{!!}/\tilde A_{!!}},M_!) 
& \quad\text{by \cite[I.3.3.4.4]{Il}} \\
& \Exal_{\tilde A_{!!}}(\tilde B_{!!},M_!) 
& \quad\text{by \cite[III.1.2.3]{Il}} \\
& \Exal_{\tilde A}(\tilde B,M) 
\end{array}$$
where the last isomorphism follows directly from lemma 
\ref{lem_eq.Ext}(ii) and the subsequent remark \ref{rem_subseq}.
Finally, \eqref{eq_equi.prod} shows that 
$\Exal_{\tilde A}(\tilde B,M)\simeq\Exal_A(B,M)$, as required.
\end{proof}
Moreover we have the following transitivity theorem as 
in \cite[II.2.1.2]{Il}.
\begin{theorem}\label{th_transit}
\index{Almost algebra(s)!$\L_{B/A}$ : cotangent complex of an!transitivity of the
|indref{th_transit}}
Let $A\to B\to C$ be a sequence of morphisms of almost $V$-algebras. 
There exists a natural distinguished triangle of\/ 
$\sD(s.C_{!!}\Mod)$ 
$$C_{!!}\otimes_{B_{!!}}\L_{B/A}\stackrel{u}{\to}\L_{C/A}
\stackrel{v}{\to}\L_{C/B}\to C_{!!}\otimes_{B_{!!}}\sigma\L_{B/A}$$
where the morphisms $u$ and $v$ are obtained by functoriality of\/ $\L$.
\end{theorem}
\begin{proof} It  follows directly from {\em loc. cit\/}. 
\end{proof}

\begin{proposition} \label{prop_L-and-colimits}
Let $(A_\lambda\to B_\lambda)_{\lambda\in I}$ be
a system of almost $V$-algebra morphisms indexed by a 
small filtered category $I$. Then there is a natural 
isomorphism in 
$\sD(s.\colim{\lambda\in I}B_{\lambda!!}\Mod)$
$$
\colim{\lambda\in I}\L_{B_\lambda/A_\lambda}\simeq
\L_{\colim{\lambda\in I}B_\lambda/
\colim{\lambda\in I}A_\lambda}.
$$
\end{proposition}
\begin{proof} Remark \ref{rem_adjoints}(i) gives an 
isomorphism : $\colim{\lambda\in I}A_{\lambda!!}
\stackrel{\sim}{\to}(\colim{\lambda\in I}A_\lambda)_{!!}$ 
(and likewise for $\colim{\lambda\in I}B_\lambda$). Then the 
claim follows from \cite[II.1.2.3.4]{Il}.
\end{proof}

Next we want to prove the almost version of the flat base 
change theorem \cite[II.2.2.1]{Il}. To this purpose we need 
some preparation.

\begin{proposition}\label{prop_Torindep} Let $B$ and $C$ be two 
$A$-algebras and set $T_i:=\Tor_i^{A_{!!}}(B_{!!},C_{!!})$. 
\begin{enumerate}
\item
If $A$, $B$, $C$ and $B\otimes_AC$ are all exact, 
then for every $i>0$ the natural morphism 
$\tilde\fm\otimes_VT_i\to T_i$ is an isomorphism.
\item
If, furthermore, $\Tor_i^A(B,C)\simeq 0$ for some $i>0$, then
the corresponding $T_i$ vanishes.
\end{enumerate}
\end{proposition}
\begin{proof} (i): for any almost $V$-algebra $D$ we let $k_D$ denote
the complex of $D_{!!}$-modules 
$[\tilde\fm\otimes_VD_{!!}\to D_{!!}]$ placed in degrees $-1,0$;
we have a distiguished triangle
$$
\cT(D)~:~\tilde\fm\otimes_VD_{!!}\to D_{!!}\to 
k_D\to \tilde\fm\otimes_VD_{!!}[1].
$$
By assumption,  the natural map $k_A\to k_B$ is a 
quasi-isomorphism and $\tilde\fm\otimes_VB_{!!}\simeq B_!$.
On the other hand, for all $i\in\N$ we have 
$$\Tor ^{A_{!!}}_i(k_B,C_{!!})\simeq
\Tor ^{A_{!!}}_i(k_A,C_{!!})\simeq H^{-i}(k_A\otimes_{A_{!!}}C_{!!})=
H^{-i}(k_C).$$
In particular $\Tor ^{A_{!!}}_i(k_B,C_{!!})=0$ for all $i>1$.
As $\tilde\fm$ is flat over $V$, we have
$\tilde\fm\otimes_VT_i\simeq
\Tor^{A_{!!}}_i(\tilde\fm\otimes_VB_{!!},C_{!!})$.
Then by the long exact Tor sequence associated to 
$\cT(B)\derotimes_{A_{!!}}C_{!!}$ we get the assertion 
for all $i>1$. Next we consider the natural 
map of distinguished triangles 
$\cT(A)\derotimes_{A_{!!}}A_{!!}\to
\cT(B)\derotimes_{A_{!!}}C_{!!}$;
writing down the associated morphism of long exact Tor
sequences, we obtain a diagram with exact rows :
$$\xymatrix{
0 \ar[r]  & 
\Tor ^{A_{!!}}_1(k_A,A_{!!}) \ar[r]^-\partial \ar[d] &
(\tilde\fm\otimes_VA_{!!})\otimes_{A_{!!}}A_{!!} \ar[r]^-i \ar[d] &
A_{!!}\otimes_{A_{!!}}A_{!!} \ar[d] \\ 
& \Tor ^{A_{!!}}_1(k_B,C_{!!}) \ar[r]^-{\partial'} &
(\tilde\fm\otimes_VB_{!!})\otimes_{A_{!!}}C_{!!} \ar[r]^-{i'} &
B_{!!}\otimes_{A_{!!}}C_{!!}.
}$$
By the above, the leftmost vertical map is an isomorphism;
moreover, the assumption gives 
$\Ker\,i\simeq\Ker(\tilde\fm\to V)\simeq\Ker\,i'$. Then,
since $\partial$ is injective, also $\partial'$ must be
injective, which implies our assertion for the remaining
case $i=1$. (ii): follows directly from (i).
\end{proof}

\begin{theorem}\label{th_flat.base.ch} Let $B$, $A'$ be 
two $A$-algebras. Suppose that the natural morphism 
$B\derotimes_AA'\to B':=B\otimes_AA'$ is an isomorphism in 
$\sD(s.A\Mod)$. Then the natural morphisms
$$\begin{array}{l}
B'_{!!}\otimes_{B_{!!}}\L_{B/A}\to\L_{B'/A'} \\
(B'_{!!}\otimes_{B_{!!}}\L_{B/A})\oplus
(B'_{!!}\otimes_{A'_{!!}}\L_{A'/A})\to\L_{B'/A}
\end{array}$$
are quasi-isomorphisms.
\end{theorem}
\begin{proof} Let us remark that the functor 
$D\mapsto V^a\times D$ : $A\Alg\to(V^a\times A)\Alg$ 
commutes with tensor products; hence the same holds for the 
functor $D\mapsto(V^a\times D)_{!!}$ (see remark 
\ref{rem_adjoints}(i)). Then, in view of proposition 
\ref{prop_Torindep}(ii), the theorem is reduced immediately 
to \cite[II.2.2.1]{Il}.
\end{proof}

As an application we obtain the vanishing of the almost cotangent 
complex for a certain class of morphisms.

\begin{theorem}\label{th_vanish.L} Let $R\to S$ be a morphism 
of almost algebras such that 
$$\Tor^R_i(S,S)\simeq 0\simeq\Tor^{S\otimes_RS}_i(S,S)
\qquad\text{for all $i>0$}$$
(for the natural $S\otimes_RS$-module structure induced by 
$\mu_{S/R}$). Then $\L_{S/R}\simeq 0$ in $\sD(S_{!!}\Mod)$.
\end{theorem}
\begin{proof} Since $\Tor^R_i(S,S)=0$ for all $i>0$, theorem
\ref{th_flat.base.ch} applies (with $A:=R$ and $B:=A':=S$), giving 
the natural isomorphisms
\set\begin{equation}\label{eq_two.equa}
\begin{array}{l}
(S\otimes_RS)_{!!}\otimes_{S_{!!}}\L_{S/R}\simeq\L_{S\otimes_RS/S} \\
((S\otimes_RS)_{!!}\otimes_{S_{!!}}\L_{S/R})\oplus
((S\otimes_RS)_{!!}\otimes_{S_{!!}}\L_{S/R})\simeq\L_{S\otimes_RS/R}.
\end{array}
\end{equation}
Since $\Tor^{S\otimes_RS}_i(S,S)=0$, the same theorem also applies
with $A:=S\otimes_RS$, $B:=S$, $A':=S$, and we notice that in this
case $B'\simeq S$; hence we have 
\set\begin{equation}\label{eq_thanksSR}
\L_{S/S\otimes_RS}\simeq 
S_{!!}\otimes_{S_{!!}}\L_{S/S\otimes_RS}\simeq\L_{S/S}\simeq 0.
\end{equation}
Next we apply transitivity to the sequence $R\to S\otimes_RS\to S$,
to obtain (thanks to \eqref{eq_thanksSR}) 
\set\begin{equation}\label{eq_one.down}
S_{!!}\otimes_{S\otimes_RS_{!!}}\L_{S\otimes_RS/R}\simeq\L_{S/R}.
\end{equation}
Applying $S_{!!}\otimes_{S\otimes_RS_{!!}}-$ to the second isomorphism
\eqref{eq_two.equa} we obtain 
\set\begin{equation}\label{eq_two.to.go}
\L_{S/R}\oplus\L_{S/R}\simeq
S_{!!}\otimes_{S\otimes_RS_{!!}}\L_{S\otimes_RS/R}.
\end{equation}
Finally, composing \eqref{eq_one.down} and \eqref{eq_two.to.go}
we derive 
\set\begin{equation}\label{eq_over}
\L_{S/R}\oplus\L_{S/R}\stackrel{\sim}{\to}\L_{S/R}.
\end{equation}
However, by inspection, the isomorphism \eqref{eq_over} is the
sum map. Consequently $\L_{S/R}\simeq 0$, as claimed.
\end{proof}

The following proposition shows that $\L_{B/A}$ is already
determined by $\L_{B/A}^a$.
\begin{proposition}\label{prop_quasi.exact}
Let $A\to B$ be a morphism of exact almost $V$-algebras. 
Then the natural map 
$\tilde\fm\otimes_V\L_{B_{!!}/A_{!!}}\to\L_{B_{!!}/A_{!!}}$ 
is a quasi-isomorphism.
\end{proposition}
\begin{proof} By transitivity we may assume $A=V^a$. Let 
$P_\bullet:=P_V(B_{!!})$ be the standard resolution of $B_{!!}$
(see \cite[II.1.2.1]{Il}). Each $P[n]^a$ contains $V$ as a 
direct summand, hence it is exact, so that we have an exact 
sequence of simplicial $V$-modules 
$0\to s.\tilde\fm\to s.V\oplus(P^a_\bullet)_!\to
(P^a_\bullet)_{!!}\to 0$.
The augmentation $(P^a_\bullet)_!\to(B^a_{!!})_!\simeq B_!$
is a quasi-isomorphism and we deduce that 
$(P^a_\bullet)_{!!}\to B_{!!}$ is a quasi-isomorphism; hence 
$(P^a_\bullet)_{!!}\to P_\bullet$ is a quasi-isomorphism as 
well. We have $P[n]\simeq\Sym(F_n)$ for a free $V$-module 
$F_n$ and the map $(P[n]^a)_{!!}\to P[n]$ is identified with 
$\Sym(\tilde\fm\otimes_VF_n)\to\Sym(F_n)$, whence 
$\Omega_{P[n]^a_{!!}/V}\otimes_{P[n]^a_{!!}}P[n]\to
\Omega_{P[n]/V}$ is identified with 
$\tilde\fm\otimes_V\Omega_{P[n]/V}\to\Omega_{P[n]/V}$.
By \cite[II.1.2.6.2]{Il} the map 
$\L^\Delta_{(P^a_\bullet)_{!!}/V}\to\L^\Delta_{P_\bullet/V}$
is a quasi-isomorphism. In view of \cite[II.1.2.4.4]{Il} we 
derive that $\Omega_{(P_\bullet^a)_{!!}/V}\to\Omega_{P_\bullet/V}$
is a quasi-isomorphism, {\em i.e.\/} 
$\tilde\fm\otimes_V\Omega_{P_\bullet/V}\to\Omega_{P_\bullet/V}$
is a quasi-isomorphism. Since $\tilde\fm$ is flat and 
$\Omega_{P_\bullet/V}\to
\Omega_{P_\bullet/V}\otimes_{P_\bullet}B_{!!}=\L_{B_{!!}/V}$
is a quasi-isomorphism, we get the desired conclusion.
\end{proof}

Finally we have a fundamental spectral sequence as in \cite[III.3.3.2]{Il}. 
\begin{theorem}\label{th_main.spec.seq} 
Let $\phi:A\to B$ be a morphism of almost algebras such that 
$B\otimes_AB\simeq B$ ({\em e.g.} such that $B$ is a quotient of $A$). 
Then there is a first quadrant homology spectral sequence of bigraded 
almost algebras
$$E^2_{pq}:=H_{p+q}(\mbox{\rm Sym}^q_B(\L^a_{B/A}))\Rightarrow
\Tor^{A}_{p+q}(B,B).$$
\end{theorem}
\begin{proof} We replace $\phi$ by $\one_{V^a}\times\phi$ and
apply the functor $B\mapsto B_{!!}$ (which commutes with tensor 
products by remark \ref{rem_adjoints}(i)) thereby reducing 
the assertion to \cite[III.3.3.2]{Il}.
\end{proof}

\newpage

\section{Almost ring theory}\label{ch_alm.ring.th}

With this chapter we begin in earnest the study 
of almost commutative algebra: in section \ref{sec_earnest}
the classes of flat, unramified and {\'e}tale morphisms are
defined, together with some variants.
In section \ref{sec_lift} we derive the infinitesimal lifting 
theorems for {\'e}tale algebras (theorem \ref{th_liftetale})
and for almost projective modules (theorem \ref{th_liftmod}).

In section \ref{sec_descent} we turn to study some cases of
non-flat descent; when we specialize to usual rings, we recover
known theorems (of course, standard commutative algebra is a
particular case of almost ring theory). But if the result is not
new, the argument is : indeed, we believe that our treatment,
even when specialized to usual rings, improves upon the method
found in the literature.

The last section of chapter \ref{ch_alm.ring.th} calls on
stage the Frobenius endomorphism of an almost algebra of
positive characteristic. The main results are
invariance of {\'e}tale morphisms under pull-back by Frobenius
maps (theorem \ref{th_Frobenius}) and theorem \ref{th_easypurity}, 
which can be interpreted as a purity theorem. Perhaps the
most remarkable aspect of the latter result is how cheap
it is : in positive characteristic, the availability of the 
Frobenius map allows for a quick and easy proof. Philosophically,
this proof is not too far removed from the method devised by
Faltings for his more recent proof of purity in mixed
characteristic.

Taken together, the above-mentioned four sections leave us
with a decent understanding of the morphisms ``of relative
dimension zero'' (this expression should be taken with a
grain of salt, since we do not try to define the dimension
of an almost algebra). On one hand, a good hold on the
case of relative dimension zero is all that is required for
the applications currently in sight (especially for the proof
of the almost purity theorem, but also for the needs of our
chapters \ref{ch_val.theory} and \ref{ch_analytic}); on the
other hand, having reached this stage, one cannot help wondering
what lies ahead, for instance whether there is a good notion of
smooth morphism of almost algebras. The full answer to this question
shall be delayed until chapter \ref{ch_hensel} : there we will
introduce a class of morphisms that generalize ``in higher dimension''
the class of weakly \'etale morphisms, and that specialize to
formally smooth morphisms in the ``classical limit'' $V=\fm$.
We will present evidence that our notion of smoothness is well
behaved and worthwhile; however we shall also see that smoothness
``in higher dimension'' exhibits some extra twists that have no
analogue in standard commutative algebra, and cannot
be easily guessed just by extrapolating from the case of \'etale
morphisms of $V^a$-algebras (which, after all, reproduce very
faithfully the behaviour of the classical notion defined in EGA). 

Such extra twists are already foreshadowed by our results on the
nilpotent deformation of almost projective modules : we show that
such deformations exist, but are not unique; however, any two such
deformations are ``very close'' in a precise sense (proposition
\ref{prop_best}). 

In (usual) algebraic geometry one can also
regard projective modules of finite rank $n$ as $\mathrm{GL}_n$-torsors
(say for the Zariski topology); in almost ring theory this description
carries through at least for the class of almost projective modules
{\em of finite rank\/} (to be defined in section \ref{sec_alm.fin.rk}).
From this vantage, one is naturally led to ask how much of the standard
deformation theory for torsors over arbitrary group schemes generalizes
to the almost world. We take up this question in section \ref{sec_def.tors},
focusing especially on the case of {\em smooth affine almost group schemes},
as a warm up to the later study of general smooth morphisms. Having
committed seriously to almost group schemes and almost torsors, it is
only a short while before one grows impatient at the limited expressive
range afforded by the purely algebraic terminology introduced thus far,
which by this point starts feeling a little like a linguistic straightjacket.
That is why we find ourselves compelled to introduce a more geometric
language : for our purpose an {\em affine almost scheme\/} is just an
object of the opposite of the category of $V^a$-algebras; similarly
we define quasi-coherent modules on affine almost schemes, as well as
some suggestive notation to go with it, mimicking the standard usage
in algebraic geometry.

Once these preliminaries are in place, the theory proceeds as in
\cite{Il2} : the techniques are rather sophisticated, but all the
hard work has already been done in the previous sections, and we can
just adapt Illusie's treatise without much difficulty.

\subsection{Flat, unramified and {\'e}tale morphisms}\label{sec_earnest}
Let $A\to B$ be a morphism of almost $V$-algebras. Using the natural
``multiplication'' morphism of $A$-algebras $\mu_{B/A}:B\otimes_AB\to B$ 
we can view $B$ as a $B\otimes_AB$-algebra.

\begin{definition}\label{def_morph}
\index{Almost algebra(s)!(faithfully) flat|indref{def_morph}}
\index{Almost algebra(s)!almost projective|indref{def_morph}}
\index{Almost algebra(s)!(uniformly) almost finite|indref{def_morph}}
\index{Almost algebra(s)!(weakly) unramified|indref{def_morph}}
\index{Almost algebra(s)!(weakly) {\'e}tale|indref{def_morph}}
\index{Almost algebra(s)!almost finite projective|indref{def_morph}}
Let $\phi:A\to B$ be a morphism of almost $V$-algebras.
\begin{enumerate}
\item
We say that $\phi$ is a {\em flat\/} (resp. {\em faithfully flat\/},
resp. {\em almost projective\/}) {\em morphism\/} if 
$B$ is a flat (resp. faithfully flat, resp. almost projective) 
$A$-module.
\item
We say that $\phi$ is {\em (uniformly) almost finite\/} 
(resp. {\em finite\/}) if $B$ is a (uniformly)  almost finitely 
generated (resp. finitely generated) $A$-module.
\item
We say that $\phi$ is {\em weakly unramified\/} (resp. 
{\em unramified\/}) if $B$ is a flat (resp. almost projective) 
$B\otimes_AB$-module (via the morphism $\mu_{B/A}$ defined above).
\item
$\phi$ is {\em weakly {\'e}tale\/} (resp. {\em {\'e}tale\/}) 
if it is flat and weakly unramified (resp. unramified).
\end{enumerate}
Furthermore, in analogy with definition \ref{def_flat_modules}, 
we shall write ``(uniformly) almost finite projective" to denote 
a morphism $\phi$ which is both (uniformly) almost finite and 
almost projective.
\end{definition}

\begin{lemma}\label{lem_itoiv}
Let $\phi:A\to B$ and $\psi:B\to C$ be morphisms of almost 
$V$-algebras. 
\begin{enumerate}
\item
Let $A\to A'$ be any morphism of $V^a$-algebras; if $\phi$ is 
flat (resp. almost projective, resp. faithfully flat, resp. 
almost finite, resp. weakly unramified, resp. unramified, resp. 
weakly {\'e}tale, resp. {\'e}tale)  then the same holds for 
$\phi\otimes_A\one_{A'}$.
\item
If both $\phi$ and $\psi$ are flat (resp. almost projective,
resp. faithfully flat, resp. almost finite, resp. weakly unramified,
resp. unramified, resp. weakly {\'e}tale, resp. {\'e}tale), then so is 
$\psi\circ\phi$.
\item
If $\phi$ is flat and $\psi\circ\phi$ is faithfully flat, then 
$\phi$ is faithfully flat.
\item
If $\phi$ is weakly unramified and $\psi\circ\phi$ is flat 
(resp. weakly {\'e}tale), then $\psi$ is flat (resp. weakly 
{\'e}tale).
\item
If $\phi$ is unramified and $\psi\circ\phi$ is {\'e}tale, then
$\psi$ is {\'e}tale.
\item
$\phi$ is faithfully flat if and only if it is a monomorphism
and $B/A$ is a flat $A$-module.
\item
If $\phi$ is almost finite and weakly unramified, then $\phi$
is unramified.
\item
If $\psi$ is faithfully flat and $\psi\circ\phi$ is flat (resp.
weakly unramified), then $\phi$ is flat (resp. weakly unramified). 
\end{enumerate}
\end{lemma}
\begin{proof} For (vi) use the Tor sequences. In view of proposition 
\ref{prop_converse}(ii), to show (vii) it suffices to know that $B$ 
is an almost finitely presented $B\otimes_AB$-module; but 
this follows from the existence of an epimorphism of 
$B\otimes_AB$-modules $(B\otimes_AB)\otimes_AB\to\Ker\,\mu_{B/A}$ 
defined by
$x\otimes b\mapsto x\cdot(\underline 1\otimes b-b\otimes\underline 1)$.
Of the remaining assertions, only (iv) and (v) are not obvious, but 
the proof is just the ``almost version'' of a well-known argument. 
Let us show (v); the same argument applies to (iv). We remark that 
$\mu_{B/A}$ is an {\'e}tale morphism, since $\phi$ is unramified. 
Define $\Gamma_\psi:=\one_C\otimes_B\mu_{B/A}$. By (i), 
$\Gamma_\psi$ is {\'e}tale. Define also 
$p:=(\psi\circ\phi)\otimes_A\one_B$. By (i), $p$ is flat (resp. 
{\'e}tale). The claim follows by remarking that 
$\psi=\Gamma_\psi\circ p$ and applying (ii).
\end{proof}

\begin{remark}\label{rem_naive} 
(i) Suppose we work in the classical limit case, that is, $\fm:=V$
(cp. example \ref{ex_rings}(ii)). Then we caution the reader that
our notion of ``{\'e}tale morphism'' is more general than the usual 
one, as defined in \cite{SGA1}. The relationship between the usual 
notion and ours is discussed in the digression \eqref{subsec_digression}.

(ii) The naive hope that the functor $A\mapsto A_{!!}$ 
might preserve flatness is crushed by the following counterexample.
Let $(V,\fm)$ be as in example \ref{ex_rings}(i) and let $k$ 
be the residue field of $V$. Consider the flat map 
$V\times V\to V$ defined as $(x,y)\mapsto x$. We get a flat morphism 
$V^a\times V^a\to V^a$ in $V^a\Alg$; applying the left adjoint 
to localisation yields a map $V\times_kV\to V$ that is not flat.
On the other hand, faithful flatness {\em is\/} preserved.
Indeed, let $\phi:A\to B$ be a morphism of almost algebras.
Then $\phi\/$ is a monomorphism if and only if $\phi_{!!}\/$
is injective; moreover, $B_{!!}/\Img(A_{!!})\simeq B_!/A_!$,
which is flat over $A_{!!}$ if and only if $B/A$ is flat over
$A$, by proposition \ref{prop_comp.supp}.
\end{remark}

We will find useful to study certain ``almost idempotents'', 
as in the following proposition.

\begin{proposition}\label{prop_idemp}
\index{Almost algebra(s)!(weakly) unramified!$e_{B/A}$ : idempotent
associated to an|indref{prop_idemp}}
A morphism $\phi:A\to B$ is unramified if and only if there exists 
an almost element $e_{B/A}\in B\otimes_AB_*$ such that
\begin{enumerate}
\item
$e_{B/A}^2=e_{B/A}$;
\item
$\mu_{B/A}(e_{B/A})=\underline 1$;
\item
$x\cdot e_{B/A}=0$ for all $x\in I_{B/A*}$.
\end{enumerate}
\end{proposition}
\begin{proof} Suppose that $\phi$ is unramified. We start by 
showing that for every $\eps\in\fm$ there exist almost elements 
$e_\eps$ of $B\otimes_AB$ such that 
\set\begin{equation}\label{eq_esta}
e_\eps^2=\eps\cdot e_\eps\qquad
\mu_{B/A}(e_\eps)=\eps\cdot\underline 1\qquad
I_{B/A*}\cdot e_\eps=0.
\end{equation}
Since $B$ is an almost projective $B\otimes_AB$-module, for 
every $\eps\in\fm$ there exists an ``approximate splitting'' for the
epimorphism $\mu_{B/A}:B\otimes_AB\to B$, {\em i.e.\/} a 
$B\otimes_AB$-linear morphism $u_\eps:B\to B\otimes_AB$ such that
$\mu_{B/A}\circ u_\eps=\eps\cdot\one_B$. Set 
$e_\eps:=u_\eps\circ\underline 1:A\to B\otimes_AB$.
We see that $\mu_{B/A}(e_\eps)=\eps\cdot\underline 1$. To show that 
$e^2_\eps=\eps\cdot e_\eps$ we use the $B\otimes_AB$-linearity 
of $u_\eps$ to compute
$$e^2_\eps=e_\eps\cdot u_\eps(\underline 1)=
u_\eps(\mu_{B/A}(e_\eps)\cdot\underline 1)=
u_\eps(\mu_{B/A}(e_\eps))=\eps\cdot e_\eps.$$
Next take any almost element $x$ of $I_{B/A}$ and compute
$$x\cdot e_\eps=x\cdot u_\eps(\underline 1)=
u_\eps(\mu_{B/A}(x)\cdot\underline 1)=0.$$
This establishes \eqref{eq_esta}. Next let us take any other
$\delta\in\fm$ and a corresponding almost element $e_\delta$. 
Both $\eps\cdot\underline 1-e_\eps$ and 
$\delta\cdot\underline 1-e_\delta$ are elements of $I_{B/A*}$, 
hence we have $(\delta\cdot\underline 1-e_\delta)\cdot e_\eps=0
=(\eps\cdot\underline 1-e_\eps)\cdot e_\delta$
which implies 
\set\begin{equation}\label{eq_fund.e}
\delta\cdot e_\eps=\eps\cdot e_\delta\qquad
\text{for all $\eps,\delta\in\fm$.}
\end{equation} 
Let us define a map $e_{B/A}:\fm\otimes_V\fm\to B\otimes_AB_*$ 
by the rule
\set\begin{equation}\label{eq_define.e}
\eps\otimes\delta\mapsto\delta\cdot e_\eps\qquad
\text{for all $\eps,\delta\in\fm$}.
\end{equation}
To show that \eqref{eq_define.e} does indeed determine a well defined
morphism, we need to check that 
$\delta\cdot v\cdot e_\eps=\delta\cdot e_{v\cdot\eps}$ and
$\delta\cdot e_{\eps+\eps'}=\delta\cdot(e_\eps+e_{\eps'})$
for all $\eps,\eps',\delta\in\fm$ and all $v\in V$. However,
both identities follow easily by a repeated application of 
\eqref{eq_fund.e}. It is easy to see that $e_{B/A}$ defines 
an almost element with the required properties.

Conversely, suppose an almost element $e_{B/A}$ of $B\otimes_AB$ 
is given with the stated properties. We define $u:B\to B\otimes_AB$ 
by $b\mapsto e_{B/A}\cdot(1\otimes b)$ ($b\in B_*$) and 
$v:=\mu_{B/A}$. 
Then (iii) says that $u$ is a $B\otimes_AB$-linear morphism and (ii) 
shows that $v\circ u=\one_B$. Hence, by lemma \ref{lem_thetwo}, 
$\phi$ is unramified.
\end{proof}

\begin{remark}\label{rem_idemp} The proof of proposition \ref{prop_idemp} 
shows that if $I$ is an ideal in an almost $V$-algebra $A$, then $A/I$
is almost projective over $A$ if and only if $I$ is generated by 
an idempotent of $A_*$. This idempotent is uniquely determined.
\end{remark}

\begin{corollary}\label{cor_unram} Under the hypotheses and notation 
of the proposition, the ideal $I_{B/A}$ has a natural structure 
of $B\otimes_AB$-algebra, with unit morphism given
by $\underline 1:=\underline 1_{B\otimes_AB/A}-e_{B/A}$ 
and whose multiplication is the restriction of $\mu_{B\otimes_AB/A}$ 
to $I_{B/A}$. Moreover the natural morphism
$$B\otimes_AB\to I_{B/A}\oplus B \qquad x\mapsto 
(x\cdot\underline 1\oplus \mu_{B/A}(x))$$
is an isomorphism of $B\otimes_AB$-algebras.
\end{corollary}
\begin{proof} Left to the reader as an exercise.
\end{proof}

\subsection{Nilpotent deformations of almost algebras and modules}
\label{sec_lift}
Throughout the following, the terminology ``epimorphism of 
$V^a$-alge\-bras'' will refer to a morphism of $V^a$-algebras
that induces an epimorphism on the underlying $V^a$-modules. 
\begin{lemma}\label{lem_epim}
Let $A\to B$ be an epimorphism of almost $V$-algebras with 
kernel $I$. Let $U$ be the $A$-extension $0\to I/I^2\to A/I^2\to B\to 0$.
Then the assignment $f\mapsto f*U$ defines a natural isomorphism
\set\begin{equation}\label{eq_Hom.Exal}
\Hom_B(I/I^2,M)\stackrel{\sim}{\to}\Exal_A(B,M).
\end{equation}
\end{lemma}
\begin{proof} Let $X:=(0\to M\to E\stackrel{p}{\to}B\to 0)$ be any 
$A$-extension of $B$ by $M$. The composition $g:A\to E\stackrel{p}{\to}B$
of the structural morphism for $E$ followed by $p$ coincides with
the projection $A\to B$. Therefore $g(I)\subset M$ and $g(I^2)=0$.
Hence $g$ factors through $A/I^2$; the restriction of $g$ to $I/I^2$
defines a morphism $f\in\Hom_B(I/I^2,M)$ and a morphism of $A$-extensions
$f*U\to X$. In this way we obtain an inverse for \eqref{eq_Hom.Exal}.
\end{proof}

\sset\subsubsection{}
Now consider any morphism of $A$-extensions 
\set\begin{equation}\label{eq_A.ext}{
\diagram 
~\tilde B: \dto^{\tilde f} & \qquad 0 \rto & I \rto \dto^u & 
B \rto \dto^f & B_0 \rto \dto^{f_0} & 0 \\
~\tilde C: & \qquad 0 \rto & J \rto & C \rto & C_0 \rto & 0.
\enddiagram
}\end{equation}
The morphism $u$ induces by adjunction a morphism of $C_0$-modules
\set\begin{equation}\label{eq_adj} 
C_0\otimes_{B_0}I\to J
\end{equation}
whose image is the ideal $I\cdot C$, so that the square diagram of 
almost algebras defined by $\tilde f$ is cofibred ({\em i.e.\/} 
$C_0\simeq C\otimes_BB_0$) if and only if \eqref{eq_adj} is an 
epimorphism.
\begin{lemma}\label{lem_flat}
Let $\tilde f:\tilde B\to\tilde C$ be a morphism of $A$-extensions as 
above, such that the corresponding square diagram of almost algebras 
is cofibred. Then the morphism $f:B\to C$ is flat if and only if 
$f_0:B_0\to C_0$ is flat and \eqref{eq_adj} is an isomorphism.
\end{lemma}
\begin{proof} It follows directly from the (almost version of the) local
flatness criterion (see \cite[Th. 22.3]{Mat}).
\end{proof}

We are now ready to put together all the work done so far and
begin the study of deformations of almost algebras. 

\sset\subsubsection{}\label{subsec_obstruction}
\index{$\omega(\tilde B,f_0,u)$, $\omega(\tilde B,f_0)$ :
obstruction classes|indref{subsec_obstruction}{}, \indref{subsec_obstr_again}}
The morphism $u:I\to J$ is an element in $\Hom_{B_0}(I,J)$; by lemma 
\ref{lem_epim} the latter group is naturally isomorphic to 
$\Exal_B(B_0,J)$. On the other hand, in view of proposition
\ref{prop_quasi.exact} and lemma \ref{lem_Exts.are.same}(i)
we have natural isomorphisms:
\set\begin{equation}\label{eq_back-to-almost}
\hExt^i_{C_{0!!}}(\L_{C_0/B_0},M_!)\simeq\hExt^i_{C_0}(\L^a_{C_0/B_0},M)
\end{equation}
for every $i\in\N$ and every $C_0$-module $M$.
By applying transitivity (theorem \ref{th_transit}) to the sequence of 
morphisms $B\to B_0\stackrel{f_0}{\to} C_0$ we deduce an exact sequence 
of abelian groups
$$\Exal_{B_0}(C_0,J)\to\Exal_B(C_0,J)\to\Hom_{B_0}(I,J)
\stackrel{\partial}{\to}\hExt^2_{C_0}(\L^a_{C_0/B_0},J).$$
Hence we can form the element
$\omega(\tilde B,f_0,u):=\partial(u)\in
\hExt^2_{C_0}(\L^a_{C_0/B_0},J)$.
The proof of the next result goes exactly as in \cite[III.2.1.2.3]{Il}.
\begin{proposition}\label{prop_obstr}
Let the $A$-extension $\tilde B$, the $B_0$-linear morphism $u:I\to J$ 
and the morphism of $A$-algebras $f_0:B_0\to C_0$ be given as above.
\begin{enumerate}
\item
There exists an $A$-extension $\tilde C$ and a morphism 
$\tilde f:\tilde B\to\tilde C$ completing diagram \eqref{eq_A.ext} 
if and only if $\omega(\tilde B,f_0,u)=0$. (\/{\em i.e.\/} 
$\omega(\tilde B,f_0,u)$ is the obstruction to the lifting 
of $\tilde B$ over $f_0$.)
\item
Assume that the obstruction $\omega(\tilde B,f_0,u)$ vanishes. 
Then the set of isomorphism classes of $A$-extensions $\tilde C$ as 
in {\em(i)} forms a torsor under the group:
$$\Exal_{B_0}(C_0,J)\simeq\hExt^1_{C_0}(\L^a_{C_0/B_0},J).$$
\item
The group of automorphisms of an $A$-extension $\tilde C$ as 
in {\em (i)} is naturally isomorphic to $\Der_{B_0}(C_0,J)$ 
(\/$\simeq\hExt^0_{C_0}(\L^a_{C_0/B_0},J)$\/).
\qed\end{enumerate}
\end{proposition} 

\sset\subsubsection{}\label{subsec_obstr_again}
The obstruction $\omega(\tilde B,f_0,u)$ depends functorially on $u$.
More exactly, if we denote by
$$\omega(\tilde B,f_0)\in
\hExt^2_{C_0}(\L^a_{C_0/B_0},C_0\otimes_{B_0}I)$$
the obstruction corresponding to the natural morphism 
$I\to C_0\otimes_{B_0}I$, then for any other morphism $u:I\to J$ we have 
$$\omega(\tilde B,f_0,u)=v_!\circ\omega(\tilde B,f_0)$$
where $v$ is the morphism \eqref{eq_adj}. Taking lemma \ref{lem_flat} 
into account we deduce 
\begin{corollary}\label{cor_flat.obstr}
Suppose that $B_0\to C_0$ is flat. Then
\begin{enumerate}
\item
The class $\omega(\tilde B,f_0)$ is the obstruction to the existence
of a flat deformation of\/ $C_0$ over $B$, {\em i.e.\/} of a 
$B$-extension $\tilde C$ as in \eqref{eq_A.ext} such that $C$ 
is flat over $B$ and $C\otimes_BB_0\to C_0$ is an isomorphism.
\item
If the obstruction $\omega(\tilde B,f_0)$ vanishes, then the set of 
isomorphism classes of flat deformations of\/ $C_0$ over $B$ forms 
a torsor under the group $\Exal_{B_0}(C_0,C_0\otimes_{B_0}I)$.
\item
The group of automorphisms of a given flat deformation of\/ $C_0$
over $B$ is naturally isomorphic to $\Der_{B_0}(C_0,C_0\otimes_{B_0}I)$.
\qed\end{enumerate}
\end{corollary}

\sset\subsubsection{}
Now, suppose we are given two $A$-extensions $\tilde C^1,\tilde C^2$ 
with morphisms of $A$-extensions 
$$\diagram 
~\tilde B: \dto^{\tilde f^i} & \qquad 0 \rto & I \rto \dto^{u^i} & 
B \rto \dto^{f^i} & B_0 \rto \dto^{f^i_0} & 0 \\
~\tilde C^i: & \qquad 0 \rto & J^i \rto & C^i \rto & C^i_0 \rto & 0
\enddiagram$$
and morphisms $v:J^1\to J^2$, $g_0:C_0^1\to C_0^2$ such that 
\set\begin{equation}\label{eq_problem}
u^2=v\circ u^1 \qquad \text{and} \qquad f_0^2=g_0\circ f_0^1. 
\end{equation}
We consider the problem of finding a morphism of $A$-extensions 
\set\begin{equation}\label{eq_A.ext.C}{
\diagram
~\tilde C^1: \dto^{\tilde g} & \qquad 0 \rto & J^1 \rto \dto^v & 
C^1 \rto \dto^g & C^1_0 \rto \dto^{g_0}& 0 \\
~\tilde C^2: & \qquad 0 \rto & J^2 \rto & C^2 \rto & C^2_0 \rto & 0
\enddiagram
}\end{equation}
such that $\tilde f^2=\tilde g\circ\tilde f^1$. Let us denote 
by $e(\tilde C^i)\in\hExt^1_{C^i_0}(\L^a_{C^i_0/B},J^i)$ 
the classes defined by the $B$-extensions $\tilde C^1,\tilde C^2$ 
via the isomorphisms \eqref{eq_th-main} and \eqref{eq_back-to-almost},
and by
$$
\begin{array}{r@{~:~}l}
v* & \hExt^1_{C^1_0}(\L^a_{C^1_0/B},J^1)\to
\hExt^1_{C^1_0}(\L^a_{C^1_0/B},J^2) \\
*g_0 & \hExt^1_{C^2_0}(\L^a_{C^2_0/B},J^2)\to
\hExt^1_{C^2_0}(C^2_0\otimes_{C^1_0}\L^a_{C^1_0/B},J^2)
\end{array}
$$
the canonical morphisms defined by $v$ and $g_0$. Using the natural
isomorphism
$$\hExt^1_{C^1_0}(\L^a_{C^1_0/B},J^2)\simeq
\hExt^1_{C^2_0}(C^2_0\otimes_{C^1_0}\L_{C^1_0/B},J^2)$$
we can identify the target of both $v*$ and $*g$ with 
$\hExt^1_{C^1_0}(\L^a_{C^1_0/B},J^2)$.
It is clear that the problem
admits a solution if and only if the $A$-extensions $v*\tilde C^1$ 
and $\tilde C^2*g_0$ coincide, {\em i.e.\/} if and only if 
$v*e(\tilde C^1)-e(\tilde C^2)*g_0=0$. By applying transitivity 
to the sequence of morphisms $B\to B_0\to C^1_0$ we obtain an 
exact sequence
$$
\hExt^1_{C^1_0}(\L^a_{C^1_0/B_0},J^2)\hookrightarrow
\hExt^1_{C^1_0}(\L^a_{C^1_0/B},J^2)
\to\Hom_{C^1_0}(C^1_0\otimes_{B_0}I,J^2).
$$
It follows from \eqref{eq_problem} that the image of 
$v*e(\tilde C^1)-e(\tilde C^2)*g_0$ in the group 
$\Hom_{C^1_0}(C^1_0\otimes_{B_0}I,J^2)$ vanishes, therefore 
\set\begin{equation}\label{eq_obstr}
v*e(\tilde C^1)-e(\tilde C^2)*g_0\in
\hExt^1_{C^1_0}(\L^a_{C^1_0/B_0},J^2).
\end{equation}
In conclusion, we derive the following result as in \cite[III.2.2.2]{Il}.
\begin{proposition}\label{prop_obstr2}
With the above notations, the class \eqref{eq_obstr} is 
the obstruction to the existence of a morphism of $A$-extensions 
$\tilde g:\tilde C^1\to\tilde C^2$ as in \eqref{eq_A.ext.C} such that
$\tilde f^2=\tilde g\circ\tilde f^1$. When the obstruction vanishes,
the set of such morphisms forms a torsor under the group 
$\Der_{B_0}(C^1_0,J^2)$ (the latter being identified with 
$\hExt^0_{C^2_0}(C^2_0\otimes_{C^1_0}\L^a_{C^1_0/B_0},J^2)$).
\qed\end{proposition}

\sset\subsubsection{}\label{subsec_etale_category}
\index{Almost algebra(s)!(weakly) {\'e}tale!$A\wEt$, $A\Et$ : category of|indref{subsec_etale_category}}
For a given almost $V$-algebra $A$, we define the category $A\wEt$ 
(resp. $A\Et$) as the full subcategory of $A\Alg$ consisting of all 
weakly {\'e}tale (resp. {\'e}tale) $A$-algebras. Notice that, by lemma 
\ref{lem_itoiv}(iv) all morphisms in $A\wEt$ are weakly {\'e}tale.

\begin{theorem}\label{th_liftetale} Let $A$ be a $V^a$-algebra.
\begin{enumerate}
\item
Let $B$ be a weakly {\'e}tale $A$-algebra, $C$ any $A$-algebra 
and $I\subset C$ a nilpotent ideal. Then the natural morphism
$$\Hom_{A\Alg}(B,C)\to\Hom_{A\Alg}(B,C/I)$$
is bijective.
\item
Let $I\subset A$ a nilpotent ideal and $A':=A/I$. Then the 
natural functor
$$A\wEt\to A'\wEt\qquad
(\phi:A\to B)\mapsto(\one_{A'}\otimes_A\phi:A'\to A'\otimes_AB)$$
is an equivalence of categories.
\item
The equivalence of {\em (ii)} restricts to an equivalence 
$A\Et\to A'\Et$.
\end{enumerate}
\end{theorem}
\begin{proof}By induction we can assume $I^2=0$. Then (i) follows 
directly from proposition \ref{prop_obstr2} and theorem 
\ref{th_vanish.L}. We show (ii) : by corollary \ref{cor_flat.obstr} 
(and again theorem \ref{th_vanish.L}) a given weakly {\'e}tale 
morphism $\phi':A'\to B'$ can be lifted to a {\em unique\/} 
flat morphism $\phi:A\to B$. We need to prove that $\phi$ is 
weakly {\'e}tale, {\em i.e.\/} that $B$ is $B\otimes_AB$-flat. 
However, it is clear that $\mu_{B'/A'}:B'\otimes_{A'}B'\to B'$ 
is weakly {\'e}tale, hence it has a flat lifting 
$\tilde\mu:B\otimes_AB\to C$. Then the composition 
$A\to B\otimes_AB\to C$ is flat and it is a lifting 
of $\phi'$. We deduce that there is an isomorphism of 
$A$-algebras $\alpha:B\to C$ lifting $\one_{B'}$ and moreover
the morphisms $b\mapsto\tilde\mu(b\otimes\underline 1)$ and
$b\mapsto\tilde\mu(\underline 1\otimes b)$ coincide with $\alpha$. 
Claim (ii) follows. To show (iii), suppose that $A'\to B'$ is {\'e}tale
and let $I_{B'/A'}$ denote as usual the kernel of $\mu_{B'/A'}$.  
By corollary \ref{cor_unram} there is a natural morphism of 
almost algebras $B'\otimes_{A'}B'\to I_{B'/A'}$ which is clearly 
{\'e}tale. Hence $I_{B'/A'}$ lifts to a weakly {\'e}tale
$B\otimes_AB$-algebra $C$, and the isomorphism
$B'\otimes_{A'}B'\simeq I_{B'/A'}\oplus B'$ lifts to an isomorphism
$B\otimes_AB\simeq C\oplus B$ of $B\otimes_AB$-algebras.
It follows that $B$ is an almost projective 
$B\otimes_AB$-module, {\em i.e.\/} $A\to B$ is {\'e}tale, as claimed.
\end{proof}

We conclude with some results on deformations of almost modules.
These can be established independently of the theory of the 
cotangent complex, along the lines of \cite[IV.3.1.12]{Il}.

\sset\subsubsection{}\label{subsec_webeginby}
We begin by recalling some notation from {\em loc. cit.\/}
Let $R$ be a ring and $J\subset R$ an ideal with
$J^2=0$. Set $R':=R/J$; an extension of $R$-modules
$\underline M:=(0\to K\to M\stackrel{p}{\to} M'\to 0)$ 
where $K$ and $M'$ are killed by $J$, defines a natural morphism
of $R'$-modules $u(\underline M):J\otimes_{R'}M'\to K$
such that $u(\underline M)(x\otimes m')=xm$ for $x\in J$, $m\in M$
and $p(m)=m'$. By the local flatness criterion (\cite[Th. 22.3]{Mat})
$M$ is flat over $R$ if and only if $M'$ is flat over $R'$ and
$u(\underline M)$ is an isomorphism. One can then show the following.
\begin{proposition}\label{prop_liftmodules}(cp. \cite[IV.3.1.5]{Il}) 
\index{$\omega(R,u')$ : obstruction class|indref{prop_liftmodules}}
With the notation of \eqref{subsec_webeginby} we have:
\begin{enumerate}
\item 
Given $R'$-modules $M'$ and $K$ and a morphism 
$u':J\otimes_{R'}M'\to K$ there exists an obstruction 
$\omega(R,u')\in\Ext_{R'}^2(M',K)$ whose vanishing is necessary 
and sufficient for the existence of an extension of $R$-modules 
$\underline M$ of $M'$ by $K$ such that $u(\underline M)=u'$. 
\item
When $\omega(R,u')=0$, the set of isomorphism classes of 
such extensions $\underline M$ forms a torsor under 
$\Ext^1_{R'}(M',K)$; the group of automorphisms of such
an extension is isomorphic to $\Hom_{R'}(M',K)$.
\qed\end{enumerate}
\end{proposition}

\begin{lemma}\label{lem_Gruson} 
Let $R$ be a ring, $M$ a finitely generated $R$-module
such that $\Ann_R M$ is a nilpotent ideal. Then $R$ admits a 
filtration $0=J_m\subset...\subset J_1\subset J_0=R$ such that
each $J_i/J_{i+1}$ is a quotient of a direct sum of copies of $M$. 
\end{lemma}
\begin{proof} This is \cite[1.1.5]{Gru}; for the convenience of
the reader we reproduce the proof. Let $I:=F_0(M)\subset R$;
if $M$ is generated by $k$ elements, we have 
$(\Ann_RM)^k\subset I\subset\Ann_RM$, especially $I$ is nilpotent.
\begin{claim}\label{cl_reduce-to-R/I}
It suffices to show that $\bar R:=R/I$ admits a filtration as above.
\end{claim}
\begin{pfclaim}
Indeed, if $0=J'_0\subset J'_1\subset...\subset J'_{n-1}\subset J'_n=R/I$
is such a filtration, we deduce filtrations 
$0\subset J'_1(I^t/I^{t+1})\subset...
\subset J'_{n-1}(I^t/I^{t+1})\subset (I^t/I^{t+1})$ for every
$t\in\N$; the graded module associated to this filtration is a direct
sum of quotients of modules of the form 
$(J'_k/J'_{k+1})\otimes_R(I^t/I^{t+1})$, so the claim follows easily.
\end{pfclaim}

Let $F_0(M/IM)$ be the $0$-th Fitting ideal of the $\bar R$-module
$M/IM$; we have $F_0(M/IM)=F_0(M)\cdot\bar R=0$, and 
$(\Ann_{\bar R}M/IM)^k\subset F_0(M/IM)$, {\em i.e.} $\Ann_{\bar R}M/IM$
is a nilpotent ideal. Thus, thanks to claim \ref{cl_reduce-to-R/I}
we can replace $R$ and $M$ by $\bar R$ and $M/IM$, and thereby
reduce to the case where $F_0(M)=0$. We claim that the filtration
$0=F_0(M)\subset F_1(M)\subset...\subset F_{k+1}(M)=R$ will do in
this case. Indeed, let $L_1\stackrel{\phi}{\to}L_0\to M\to 0$
be a presentation of $M$, where $L_0$ and $L_1$ are free $R$-modules
and the rank of $L_0$ equals $k$; by definition $F_j(M)$ is the image
of the map 
$\Lambda_R^{k-j}L_1\otimes_R\Lambda_R^jL_0\to\Lambda_R^kL_0\simeq R$
defined by the rule $x\otimes y\mapsto\Lambda_R^{k-j}\phi(x)\wedge y$.
We deduce easily that the induced surjection
$\Lambda_R^{k-j}L_1\otimes_R\Lambda_R^jL_0\to F_j(M)/F_{j-1}(M)$
factors through the module
$\Lambda_R^{k-j}L_1\otimes_R\Lambda_R^jM$; however the latter
is a quotient of sums of $M$, at least when $j\geq 1$, so the
claim follows.
\end{proof}

\begin{lemma}\label{lem_frequent}
Let $A\to B$ be a finite morphism of almost algebras with nilpotent
kernel. There exists $m\geq 0$ such that the following holds.
For every $A$-linear morphism $\phi:M\to N$, set 
$\phi_B:=\phi\otimes_A\one_B:M\otimes_AB\to N\otimes_AB$; then :
\begin{enumerate}
\item
$\Ann_A(\Coker\,\phi_B)^m\subset\Ann_A(\Coker\,\phi)$.
\item
$(\Ann_V(\Ker\,\phi_B)\cdot\Ann_V(\Tor^A_1(B,N))\cdot
\Ann_V(\Coker\,\phi))^m\subset\Ann_A(\Ker\,\phi)$.
\end{enumerate}
\noindent If $B=A/I$ for some nilpotent ideal $I$, and $I^n=0$, then
we can take $m=n$ in {\em (i)} and {\em (ii)}.
\end{lemma}
\begin{proof} Under the assumptions, we can find a finitely 
generated $A_*$-module $Q$ such that $\fm\cdot B_*\subset Q\subset B_*$.
By lemma \ref{lem_Gruson} there exists a finite filtration 
$0=J_m\subset ...\subset J_1\subset J_0=A_*$ such that each $J_i/J_{i+1}$
is a quotient of a direct sum of copies of $Q$. This implies that 
\set\begin{equation}\label{eq_Annihil}
\Ann_A(M\otimes_AB)^m\subset\Ann_A(M)
\end{equation}
for every $A$-module $M$; 
(i) follows easily. Notice that if $B=A/I$ and $I^n=0$, then we can 
take $m=n$ in \eqref{eq_Annihil}. For (ii) let $C^\bullet:=\Cone\,\phi$. 
We estimate $H:=H^{-1}(C^\bullet\derotimes_AB)$ in two ways. 
By the first spectral sequence of hyperhomology we have an
exact sequence $\Tor_1^A(N,B)\to H\to\Ker\,\phi_B$.
By the second spectral sequence for hyperhomology 
we have an exact sequence 
$\Tor^A_2(\Coker\,\phi,B)\to\Ker(\phi)\otimes_AB\to H$.
Hence $\Ker(\phi)\otimes_AB$ is annihilated by the product 
of the three annihilators in (ii) and the result follows by 
applying \eqref{eq_Annihil} with $M:=\Ker\,\phi$.
\end{proof}

\begin{lemma}\label{lem_flat.proj} Keep the assumptions of lemma 
{\em\ref{lem_frequent}} and let $M$ be an $A$-module.
\begin{enumerate}
\item
If $A\to B$ is an epimorphism, $M$ is flat and $M_B:=B\otimes_AM$ 
is almost projective over $B$, then $M$ is almost projective over 
$A$.
\item
If $M_B$ is an almost finitely generated $B$-module
then $M$ is an almost finitely generated $A$-module.
\item
If $\Tor_1^A(B,M)=0$ and $M_B$ is almost finitely presented
over $B$, then $M$ is almost finitely presented over $A$.
\end{enumerate}
\end{lemma}
\begin{proof} (i) : we have to show that $\Ext^1_A(M,N)$ is almost
zero for every $A$-module $N$. Let $I:=\Ker(A\to B)$; by 
assumption $I$ is nilpotent, so by the usual devissage we may 
assume that $IN=0$. If $\chi\in\Ext^1_A(M,N)$ 
is represented by an extension $0\to N\to Q\to M\to 0$ then 
after tensoring by $B$ and using the flatness of $M$ we get 
an exact sequence of $B$-modules 
$0\to N\to B\otimes_AQ\to M_B\to 0$. Thus $\chi$ comes from an 
element of $\Ext^1_{B}(M_B,N)$ which is almost zero by assumption. 

(ii) : for a given finitely generated subideal $\fm_0\subset\fm$,
let $N\subset M_B$ be a finitely generated $B$-submodule such that 
$\fm_0M_B\subset N$. Since the induced map 
$M_*\otimes_{A_*}B_*\to(M_B)_*$ is almost surjective, 
we can find a finitely generated $A$-submodule $N_0\subset M$
such that $\fm_0N\subset\Img((N_0)_B\to M_B)$; by
lemma \ref{lem_frequent}(i) it follows that 
$\fm_0^{2n}(M/N_0)=0$ for some $n\geq 0$ depending 
only on $B$, whence the claim.

(iii) : Let $\fm_0$ be as above. By (ii), $M$ is almost 
finitely generated over $A$, so we can choose a morphism
$\phi:A^r\to M$ such that $\fm_0\cdot\Coker\,\phi=0$.
Consider $\phi_B:=\phi\otimes_A\one_B:B^r\to M_B$.
By claim \ref{cl_fin.pres.was.lemma}, there is a finitely 
generated submodule $N$ of $\Ker\,\phi_B$ containing
$\fm_0^2\cdot\Ker\,\phi_B$. Notice that $\Ker(\phi)\otimes_AB$
maps onto $\Ker(B^r\to\Img(\phi)\otimes_AB)$ and
$\Ker(\Img(\phi)\otimes_AB\to M_B)\simeq
\Tor^A_1(B,\Coker\,\phi)$ is annihilated by $\fm_0$.
Hence $\fm_0\cdot\Ker\,\phi_B$ is contained in the image
of $\Ker\,\phi$ and therefore we can lift a finite generating 
set $\{x_1',...,x_n'\}$ for $\fm_0^2 N$ to almost 
elements $\{x_1,...,x_n\}$ of $\Ker\,\phi$. If we quotient 
$A^r$ by the span of these $x_i$, we get a finitely presented 
$A$-module $F$ with a morphism $\bar\phi:F\to M$
such that $\Ker(\bar\phi\otimes_AB)$ is annihilated by 
$\fm_0^4$ and $\Coker\,\bar\phi$ is annihilated by $\fm_0$.
By lemma \ref{lem_frequent}(ii) we derive 
$\fm_0^{5m}\cdot\Ker\:\bar\phi=0$ for some $m\geq 0$. 
Since $\fm_0$ is arbitrary, this proves the result.
\end{proof}

\begin{remark}\label{rem_descent} 
(i) Inspecting the proof, one sees that parts (ii)
and (iii) of lemma \ref{lem_flat.proj} hold whenever 
\eqref{eq_Annihil} holds. For instance, if $A\to B$ is any 
faithfully flat morphism, then \eqref{eq_Annihil} holds 
with $m:=1$. 

(ii) Consequently, if $A\to B$ is faithfully flat
and $M$ is an $A$-module such that $M_B$ is flat 
(resp. almost finitely generated, resp. almost finitely 
presented) over $B$, then $M$ is flat (resp. almost finitely
generated, resp. almost finitely presented) over $A$.

(iii) On the other hand, we do not know whether a general
faithfully flat morphism $A\to B$ descends almost projectivity.
However, using (ii) and proposition \ref{prop_converse}
we see that if the $B$-module $M_B$ is almost finitely 
generated projective, then $M$ has the same property.

(iv) Furthermore, if $B$ is faithfully flat and almost finitely
presented as an $A$-module, then $A\to B$ does descend 
almost projectivity, as can be easily deduced from lemma 
\ref{lem_alhom}(i) and proposition \ref{prop_converse}(ii).
\end{remark}

\sset\subsubsection{}\label{subsec_defafp}
\index{Almost algebra(s)!$A\Et_\mathrm{afp}$ : {\'e}tale almost
finitely presented|indref{subsec_defafp}}
We denote by $A\Et_\mathrm{afp}$ the full subcategory of $A\Et$
consisting of all {\'e}tale $A$-algebras $B$ such that $B$
is almost finitely presented as an $A$-module.

\begin{theorem}\label{th_liftmod} Let $I\subset A$ be a nilpotent
ideal, and set $A':=A/I$. 
\begin{enumerate}
\item
Suppose that\/ $\tilde\fm$ is a (flat) $V$-module of homological 
dimension $\leq 1$. Let $P'$ be an almost projective $A'$-module.
\begin{enumerate}
\item
There is an almost projective $A$-module $P$ 
with $A'\otimes_AP\simeq P'$. 
\item
If $P'$ is almost finitely presented, then $P$ is almost 
finitely presented.
\end{enumerate}
\item
The equivalence of theorem {\em\ref{th_liftetale}(ii)} restricts 
to an equivalence $A\Et_\mathrm{afp}\to A'\Et_\mathrm{afp}$.
\end{enumerate}
\end{theorem}
\begin{proof} (i.a): as usual we reduce to $I^2=0$. Then proposition 
\ref{prop_liftmodules}(i) applies with $R:=A_*$, $J:=I_*$, 
$R':=A_*/I_*$, $M':=P'_!$, $K:=I_*\otimes_{R'}P'_!$ and 
$u':=\one_K$. We obtain a class 
$\omega(A_*,u')\in\Ext^2_{R'}(P'_!,I_*\otimes_{R'}P'_!)$ which 
gives the obstruction to the existence of a flat $A_*$-module 
$F$ lifting $P'_!$. Since $P'_!$ is almost projective, lemma
\ref{lem_Exts.are.same}(ii) says that
$\Ext^2_{R'}(P'_!,I_*\otimes_{R'}P'_!)=0$, so such $F$ can 
always be found, and then the $A$-module $P=F^a$ is a flat 
lifting of $P'$; by lemma \ref{lem_flat.proj}(i) we see that 
$P$ is almost projective. Now (i.b) follows from (i.a), lemma 
\ref{lem_flat.proj}(ii) and proposition \ref{prop_converse}(i).

(ii): in view of theorem \ref{th_liftetale}(iii), we only have 
to show that an {\'e}tale $A$-algebra $B$ is almost finitely presented 
as an $A$-module whenever $B\otimes_AA'$ is almost finitely 
presented as an $A'$-module. However, the assertion is a direct
consequence of lemma \ref{lem_flat.proj}(iii).
\end{proof}

\begin{remark}\label{rem_counter} (i) According to proposition
\ref{th_countably.pres}(ii.b), theorem \ref{th_liftmod}(i) applies
especially when $\fm$ is countably generated as a $V$-module.

(ii) For $P$ and $P'$ as in theorem 
\ref{th_liftmod}(i.b) let $\sigma_P:P\to P'$ be the projection. 
It is natural to ask whether the pair $(P,\sigma_P)$ is uniquely 
determined up to isomorphism, {\em i.e.\/} whether, for any 
other pair $(Q,\sigma_Q:Q\to P')$ for which theorem 
\ref{th_liftmod}(i.b) holds, there exists an $A$-linear isomorphism 
$\phi:P\to Q$ such that $\sigma_Q\circ\phi=\sigma_P$. 
The answer is negative in general. Consider the case $P':=A'$. 
Take $P:=Q:=A$ and let $\sigma_P$ be the natural projection, while 
$\sigma_Q:=(u'\cdot\one_{A'})\circ\sigma_P$, where $u'$ is a unit 
in $A'_*$. Then the uniqueness question amounts to whether 
every unit in $A'_*$ lifts to a unit of $A_*$. The 
following counterexample is related to the fact that the 
completion of the algebraic closure $\bar\Q_p$ of $\Q_p$ 
is not maximally complete. Let $V:=\bar\Z_p$, the integral 
closure of $\Z_p$ in $\bar\Q_p$. Then $V$ is a non-discrete 
valuation ring of rank one, and we take for $\fm$ the maximal
ideal of $V$, $A:=(V/p^2V)^a$ and $A':=A/pA$. Choose a 
compatible system of roots of $p$. An almost element of $A'$ 
is just a $V$-linear morphism 
$\phi:\colim{n>0}p^{1/n!}V\to V/pV$.
Such a $\phi$ can be represented (in a non-unique way) by an 
infinite series of the form $\sum^\infty_{n=1}a_np^{1-1/n!}$
($a_n\in V$). The meaning of this expression is as follows. 
For every $m>0$, scalar multiplication by the element 
$\sum^{m}_{n=1}a_np^{1-1/n!}\in V$ defines a morphism 
$\phi_m:p^{1/m!}V\to V/pV$. For $m'>m$, let 
$j_{m,m'}:p^{1/m!}V\to p^{1/m'!}V$ be the imbedding. Then we
have $\phi_{m'}\circ j_{m,m'}=\phi_m$, so that we can define 
$\phi:=\colim{m>0}\phi_m$. Similarly, every almost element of 
$A$ can be represented by an expression of the form 
$a_0+\sum^\infty_{n=1}a_np^{2-1/n!}$. Now, if $\sigma:A\to A'$
is the natural projection, the induced map 
$\sigma_*:A_*\to A'_*$ is given by:
$a_0+\sum^\infty_{n=1}a_np^{2-1/n!}\mapsto a_0$.
In particular, its image is the subring $V/p\subset(V/p)_*=A'_*$.
For instance, the unit $\sum^\infty_{n=1}p^{1-1/n!}$ of $A'_*$
does not lie in the image of this map.
\end{remark}
In the light of remark \ref{rem_counter}, the best one can achieve 
in general is the following result.
\begin{proposition}\label{prop_best} Assume ({\bf A}) (see 
\eqref{subsec_conditions.A.B}) and keep the notation of theorem 
{\em\ref{th_liftmod}(i)}. Suppose that $(Q,\sigma_Q:Q\to P')$ 
and $(P,\sigma_P:P\to P')$ are two pairs as in remark 
{\em\ref{rem_counter}(ii)}. Then for every $\eps\in\fm$ there exist 
$A$-linear morphisms $t_\eps:P\to Q$ and $s_\eps:Q\to P$ 
such that 
\smallskip

\noindent{\bf PQ}\rm{(}$\eps$\rm{)}$\qquad\qquad\qquad
\begin{array}{r@{\:=\:}lr@{\:=\:}l}
\sigma_Q\circ t_\eps & \eps\cdot\sigma_P & \qquad
\sigma_P\circ s_\eps & \eps\cdot\sigma_Q \\
s_\eps\circ t_\eps & \eps^2\cdot\one_P & \qquad
t_\eps\circ s_\eps & \eps^2\cdot\one_Q.
\end{array}$
\end{proposition}
\begin{proof}
Since both $Q$ and $P$ are almost projective and $\sigma_P,\sigma_Q$ 
are epimorphisms, there exist morphisms $\bar t_\eps:P\to Q$ and 
$\bar s_\eps:Q\to P$ such that 
$\sigma_Q\circ\bar t_\eps=\eps\cdot\sigma_P$ and 
$\sigma_P\circ\bar s_\eps=\eps\cdot\sigma_Q$. Then we have 
$\sigma_P\circ(\bar s_\eps\circ\bar t_\eps-\eps^2\cdot\one_P)=0$
and
$\sigma_Q\circ(\bar t_\eps\circ\bar s_\eps-\eps^2\circ\one_Q)=0$,
{\em i.e.\/} the morphism 
$u_\eps:=\eps^2\cdot\one_P-\bar s_\eps\circ\bar t_\eps$ (resp.
$v_\eps:=\eps^2\cdot\one_Q-\bar t_\eps\circ\bar s_\eps$) has image 
contained in the almost submodule $IP$ (resp. $IQ$). Since $I^m=0$ 
this implies $u^m_\eps=0$ and $v^m_\eps=0$. Hence
$$\eps^{2m}\cdot\one_P=(\eps^2\one_P)^m-u^m_\eps=
(\sum_{a=0}^{m-1}\eps^{2a}u^{m-1-a}_\eps)\circ
\bar s_\eps\circ\bar t_\eps.$$
Define $\bar s_{(2m-1)\eps}:=
(\sum_{a=0}^{m-1}\eps^{2a}u^{m-1-a}_\eps)\circ\bar s_\eps$. 
Notice that  
$\bar s_{(2m-1)\eps}=
\bar s_\eps\circ(\sum_{a=0}^{m-1}\eps^{2a}v^{m-1-a}_\eps).$
This implies the equalities
$\bar s_{(2m-1)\eps}\circ\bar t_\eps=\eps^{2m}\cdot\one_P$ and
$\bar t_\eps\circ\bar s_{(2m-1)\eps}=\eps^{2m}\cdot\one_Q$.
Then the pair $(\bar s_{(2m-1)\eps},\eps^{2(m-1)}\cdot\bar t_\eps)$
satisfies {\bf PQ}($\eps^{2m-1}$). Under ({\bf A}), every element
of $\fm$ is a multiple of an element of the form $\eps^{2m-1}$, 
therefore the claim follows for arbitrary $\eps\in\fm$.
\end{proof}

\subsection{Nilpotent deformations of torsors}\label{sec_def.tors}
For the considerations that follow it will be convenient
to extend yet further our basic setup. Namely, suppose
that $T$ is any topos; we can define a {\em basic setup
relative to $T$} as a pair $(V,\fm)$ consisting of a
$T$-ring $V$ and an ideal $\fm\subset V$ satisfying the
usual assumptions \eqref{subsec_basic.setup}. Then most
of the discussion of chapter \ref{ch_homol-th} extends
{\em verbatim\/} to this relative setting. Accordingly,
we generally continue to use the same notation as in 
{\em loc.cit.}; however, if it is desirable for clarity's
sake, we may sometimes stress the dependance on $T$ by mentioning
it explicitly. For instance, an almost $T$-ring is an object
of the category $V^a\Alg$ of associative, commutative and unitary
monoids of the abelian tensor category of $V^a$-modules, and
sometimes the same category is denoted $(T,V,\fm)^a\Alg$.

\sset\subsubsection{}\label{subsec_topos.basic}
Let $T$ be a topos and $(V,\fm)$ a basic setup relative to $T$.
For every object $U$ of $T$, we can consider the restriction
$(V_{/U},\fm_{/U})$ of $(V,\fm)$ to $T_{/U}$, which is a basic
setup relative to the latter topos. As usual, for every object 
$X\to U$ of $T_{/U}$, one defines 
$V_{/U}(X):=\Hom_{T_{/U}}(X,V)$ and similarly for $\fm_{/U}(X)$.
The restriction functor $V\Mod\to V_{/U}\Mod$ clearly preserves
almost isomorphisms, whence a restriction functor
$$
(T,V,\fm)^a\Mod\to(T_{/U},V_{/U},\fm_{/U})^a\Mod\qquad 
M\mapsto M_{/U}.
$$
Similar functors exist for the categories of $V^a$-algebras,
and more generally, for $A$-algebras, where $A$ is any 
$(T,V,\fm)^a$-algebra.

\sset\subsubsection{}
For every almost $T$-module $M$, we define a functor 
$M_*:T^o\to\Z\Mod$ by the rule: 
$U\mapsto\Hom_{V^a_{/U}}(V^a_{/U},M_{/U})$. Let $N$ be a
$V$-module; using the natural isomorphism
$$
\Hom_{V_{/U}}(\tilde\fm_{/U},N_{/U})\stackrel{\sim}{\to}N^a_*(U)
$$
(cp. \eqref{eq_alm.morph}) one sees that $M_*$ is a sheaf for
the canonical topology of $T$, hence it representable by an
abelian group object of $T$, which we denote by the same name.
It is then easy to check that $M_*$ is a $V$-module, and that
the functor $M\mapsto M_*$ is right adjoint to the localization
functor $V\Mod\to V^a\Mod$.
We can then generalize to the case of the localization functors
$A\Mod\to A^a\Mod$ and $A\Alg\to A^a\Alg$, for an arbitrary 
$V$-algebra $A$. Likewise, the left adjoint functors to localization
$M\mapsto M_!$ and $B\mapsto B_{!!}$ are obtained in the same
way as in the earlier treatment of the one-point topos.

\sset\subsubsection{}\label{subsec_general-schemes}
\index{Almost scheme(s)!quasi-coherent module on
an|indref{subsec_general-schemes}}
\index{Almost scheme(s)!affine, morphism of
affine|indref{subsec_general-schemes}}
\index{$\phi^\sharp$ (for a morphism $\phi$ of almost
schemes)|indref{subsec_general-schemes}}
\index{Affine almost group scheme(s)|indref{subsec_general-schemes}}
\index{Almost algebra(s)!$\Spec\,A$ : spectrum of
an|indref{subsec_general-schemes}}
Let $T$ and $(V,\fm)$ be as in \eqref{subsec_topos.basic}
and let $R$ be any $V$-algebra. An {\em affine $R$-scheme\/} is
an object of the category $R\Alg^o$. An {\em affine almost
$R$-scheme\/} (or an {\em affine $R^a$-scheme}) is an object
of the category $R^a\Alg^o$. If $X$ is an affine $R^a$-scheme,
then we may write $\cO_X$ in place of the $R^a$-algebra $X^o$,
and an $X^o$-module will also be called a {\em quasi-coherent
$\cO_X$-module}. A morphism $\phi:X\to Y$ of affine almost schemes
is the same as the morphism of almost algebras 
$\phi^\sharp:=\phi^o:\cO_Y\to\cO_X$; moreover, $\phi$
induces pullback ({\em i.e.} tensor product $-\otimes_{\cO_Y}\cO_X$)
and direct image functors ({\em i.e.} restriction of scalars), 
which we will denote
\set\begin{equation}\label{eq_guess-what}
\phi^*:\cO_Y\Mod\to\cO_X\Mod\qquad\text{and}\qquad
\phi_*:\cO_X\Mod\to\cO_Y\Mod
\end{equation}
respectively. In the same vein, if $A$ is any $R^a$-algebra,
we may write $\Spec\,A$ to denote the object $A^o$ represented
by $A$ in the opposite category $R^a\Alg^o$.
We say that $X$ is {\em flat} over $R^a$ if $\cO_X$ is a flat
$R^a$-algebra. Clearly the category of affine $R^a$-schemes
admits arbitrary products. Hence, we can define an 
{\em affine $R^a$-group scheme\/} as a group object in the 
category of affine $R^a$-schemes (and likewise for the notion
of affine $R$-group scheme).

\sset\subsubsection{}\label{subsec_adj-for-schemes}
The functors $B\mapsto B_*$ and $B\mapsto B_{!!}$ induce functors
on almost schemes, which we denote in the same way. Notice that
the functor $X\mapsto X_*$ (resp. $X\mapsto X_{!!}$) from affine
$R^a$-schemes to $R$-schemes (resp. to $R^a_{!!}$-schemes) commutes
with all colimits (resp. with all limits). Especially, if $G$ is an
affine $R^a$-group scheme, then $G_{!!}$ is an affine $R^a_{!!}$-group
scheme.

\sset\subsubsection{}\label{subsec_def.G-action}
\index{Affine almost group scheme(s)!(right, left) action
on an almost scheme by a|indref{subsec_def.G-action}}
\index{Affine almost group scheme(s)!(right, left) action
on an almost scheme of a!$G^\bullet_X$ : nerve of an|indref{subsec_def.G-action}}
\index{Affine almost group scheme(s)!action
on a quasi-coherent module, equivariant quasi-coherent module|indref{subsec_def.G-action}}
Throughout the rest of this section we fix a $V^a$-algebra $A$
and let $S:=\Spec A$.  

Let $X$ an affine $A$-scheme, $G$ an affine $A$-group scheme.
A {\em right action\/} of $G$ on $X$ is a morphism of $S$-schemes
$\rho:G\times_S X\to X$ fulfilling the usual conditions. 
To the datum $(X,G,\rho)$ one assigns its {\em nerve\/}
$G^\bullet_X$ which is a simplicial affine $S$-scheme whose
component in degree $n$ is $G^n_X:=G^n\times_S X$, and whose
face and degeneracy morphisms
$$
\partial_i:G^{n+1}_X\to G^n_X\qquad 
\sigma_j:G^n_X\to G^{n+1}_X\qquad
\text{$i=0,...,n+1$; $j=0,...,n$}
$$
are defined for every $n\in\N$ as in \cite[Ch.VI, \S2.5]{Il2}.

\sset\subsubsection{}\label{subsec_def-G-mods}
Let $(X,G,\rho)$ be as in \eqref{subsec_def.G-action}, and
let $M$ be a quasi-coherent $\cO_X$-module. A {\em $G$-action\/}
on $M$ is a morphism of quasi-coherent $\cO_{G^1_X}$-modules
\set\begin{equation}\label{eq_for-beta}
\beta:\partial^*_0 M\to\partial^*_1 M
\end{equation}
such that the following diagram commutes:
\set\begin{equation}\label{eq_cocycle-action}
{\diagram
\partial_0^*\partial_0^*M \ar[r]^{\partial_0^*\beta} \ddouble &
\partial_0^*\partial_1^*M \rdouble & 
\partial_2^*\partial_0^*M \ar[d]^{\partial_2^*\beta} \\
\partial_1^*\partial_0^*M \ar[r]^{\partial_1^*\beta} &
\partial_1^*\partial_1^*M \rdouble & \partial_2^*\partial_1^*M
\enddiagram}
\end{equation}
and such that 
\set\begin{equation}\label{eq_rigidify}
\sigma_0^*\beta=\one_M.
\end{equation}
One also says that $(M,\beta)$ is a  {\em $G$-equivariant $\cO_X$-module}.
One defines in the obvious way the morphisms of $G$-equivariant
$\cO_X$-modules, and we denote by $\cO_X\Mod_G$ the category of
$G$-equivariant $\cO_X$-modules.

\begin{lemma} For every $G$-equivariant $\cO_X$-module $M$,
the morphism \eqref{eq_for-beta} is an isomorphism.
\end{lemma}
\begin{proof} We let $\tau:G^1_X\to G^2_X$ be the
morphism given on $T$-points by the rule: $(g,x)\mapsto(g,g^{-1},x)$.
Working out the identifications, one checks easily that 
$\tau^*\eqref{eq_cocycle-action}$ boils down to the diagram
$$
\xymatrix{ \partial_0^*M \ar[rr]^-\beta && \partial_1^*M \\
& \partial_1^*M \ar[ul]^{\tau^*\partial_0^*\beta} 
\ar[ru]_-{\partial_1^*\sigma_0^*\beta}.
}$$
However, $\partial_1^*\sigma_0^*\beta=\one_{\partial_1^*M}$, in view
of \eqref{eq_rigidify}, hence $\beta$ is an epimorphism and
$\tau^*\partial_0^*\beta$ is an monomorphism. Since 
$\partial_0\circ\tau$ is an isomorphism, we deduce that
$\beta$ is already a monomorphism, whence the claim.
\end{proof}

\begin{lemma}\label{lem_when.Gmod.abel}
If $G$ is a flat affine $S$-group scheme acting
on an affine $S$-scheme $X$, then the category of $G$-equivariant
quasi-coherent $\cO_X$-modules is abelian.
\end{lemma}
\begin{proof} We need to verify that every ($G$-equivariant)
morphism $\phi:M\to N$ of $G$-equivariant $\cO_X$-modules admits
a kernel and a cokernel (and then it will be clear that
the kernel of $N\to\Coker\,\phi$ equals the cokernel of $\Ker\,\phi\to M$,
since the same holds for $\cO_X$-modules). The obvious candidates are
the kernel $K$ and cokernel $C$ taken in the category of $\cO_X$-modules,
and one has only to show that the $G$-actions of $M$ and $N$
induce $G$-actions on $K$ and $C$. This is always the case for
$C$ (even when $G$ is not flat). To deal with $K$, one remarks
that both morphisms $\partial_0,\partial_1:G^1_X\to X$
are flat; indeed, this is clear for $\partial_1$.
Then the same holds for $\partial_0:=\partial_1\circ\omega$,
where $\omega:=(\rho,\one_X):G^1_X\to G^1_X$ is the
isomorphism deduced from the action $\rho:G^1_X\to X$.
It follows that $\Ker\,(\partial^*_i\phi)\simeq\partial^*_i(\Ker\,\phi)$
for $i=0,1$, whence the claim.
\end{proof}

\sset\subsubsection{}\label{subsec_define-iota}
\index{Almost scheme(s)!square zero deformation of an|indref{subsec_define-iota}}
\index{Almost scheme(s)!$\L_{X/Y}$ : cotangent complex of an|indref{subsec_define-iota}}
Let $\phi:X\to Y$ be a morphism of affine $S$-schemes. A {\em square 
zero deformation\/} of $X$ over $Y$ is a datum of the form
$(j:X\to X',\cI,\beta)$, consisting of:
\begin{enumerate}
\renewcommand{\labelenumi}{(\alph{enumi})}
\item
a morphism of $Y$-schemes $j:X\to X'$ such that the induced morphism
of $\cO_Y$-algebras $j^\sharp:\cO_{X'}\to\cO_X$ is an epimorphism and 
$\cJ:=\Ker\,j^\sharp$ is a  square zero ideal, and
\item
a quasi-coherent $\cO_X$-module $\cI$ with an isomorphism
of $\cO_X$-modules $\beta:j^*\cJ\stackrel{\sim}{\to}\cI$.
\renewcommand{\labelenumi}{(\roman{enumi})}
\end{enumerate}
The square zero deformations form a category $\bExal_Y(X,\cI)$, with
morphisms defined in the obvious way. As in the case of the one-point
topos, we can compute the isomorphism classes of square zero
deformations of $X$ in terms of an appropriate cotangent complex.
And, just as in the earlier treatment, we have to make sure that
we are dealing with exact algebras, hence the right definition
of the cotangent complex of $\phi$ is:
$$
\L_{X/Y}:=\iota^*\L_{(X\amalg\Spec V^a)_{!!}/(Y\amalg\Spec V^a)_{!!}}
$$
where $\iota:X_{!!}\to(X\amalg\Spec V^a)_{!!}$ is the natural morphism
of schemes. This is an object of $\sD(\cO_{X_{!!}}\Mod)$.

\sset\subsubsection{}\label{subsec_equiv-deforms}
\index{Almost scheme(s)!$\bExal_Y(X/G,\cI)$ : equivariant square zero
deformation of an|indref{subsec_equiv-deforms}}
Next, let $G$ be a flat affine $S$-group scheme acting on $X$ and $Y$, 
in such a way that $\phi$ is $G$-equivariant; it follows from
lemma \ref{lem_when.Gmod.abel} that $\cJ$ has a natural $G$-action.
A {\em $G$-equivariant square zero deformations\/} of $X$ over $Y$,
is a datum $(j:X\to X',\cI,\beta)$ as above, such that, additionally,
$X'$ and $\cI$ are endowed with a $G$-action and both $j$ and $\beta$
are $G$-equivariant. Let us denote by $\bExal_Y(X/G,\cI)$ the
category of such $G$-equivariant deformations. We aim to classify
the isomorphism classes of objects of $\bExal_Y(X/G,\cI)$, and more
generally, study the $G$-equivariant deformation theory of $X$ by
means of an appropriate cotangent complex cohomology.
This is achieved by the following device.

\sset\subsubsection{}\label{subsec_fibred-fctr}
\index{Topos!fibred|indref{subsec_fibred-fctr}}
\index{Topos!fibred!$\mathrm{Top}(X)$ : topos associated to a|indref{subsec_fibred-fctr}}
Let $I$ be a (small) category; recall (\cite[Ch.VI, \S5.1]{Il2})
that a {\em fibred topos over $I$\/} is a pseudo-functor $X$
of $I$ in the $2$-category of of topoi, {\em i.e.\/} the datum
consisting of:
\begin{enumerate}
\renewcommand{\labelenumi}{(\alph{enumi})}
\item
for every object $i$ of $I$, a topos $X_i$
\item
for every arrow $f:i\to j$ in $I$, a morphism of
topoi $X_f:X_i\to X_j$ (sometimes denoted $f$)
\item
for every composition $i\stackrel{f}{\to}j\stackrel{g}{\to}k$,
a {\em transitivity isomorphism\/} 
$X_{f,g}:X_g\circ X_f\stackrel{\sim}{\to}X_{gf}$, submitted
to certain compatibility conditions (in practice, we will
omit from the notation the transitivity isomorphisms).
\renewcommand{\labelenumi}{(\alph{enumi})}
\end{enumerate}
Given a fibred topos $X$ over $I$, one denotes by $\mathrm{Top}(X)$
the following category, which is easily seen to be a topos
(\cite[Ch.VI, \S5.2]{Il2}). An object of $E$ of $\mathrm{Top}(X)$
is the datum of
\begin{enumerate}
\renewcommand{\labelenumi}{(\alph{enumi})}
\item
for every $i\in I$, an object $E_i$ of $X_i$
\item
for every arrow $f:i\to j$, a morphism $E_f:f^*E_j\to E_i$
such that, for every composition
$i\stackrel{f}{\to}j\stackrel{g}{\to}k$, one has
$E_{gf}=E_f\circ f^*E_j$, provided one identifies 
$(gf)^*E_k$ with $f^*g^*E_k$ via $X_{f,g}$.
\renewcommand{\labelenumi}{(\alph{enumi})}
\end{enumerate}
As an example we have the topos $s.T$ whose objects are the
cosimplicial objects of $T$; indeed this is the topos $\mathrm{Top}(F)$
associated to a fibred topos $F$, whose indexing category is the 
category $\Delta^o$ defined in \eqref{sec_alm.cat} : to every object
$[n]$ of $\Delta^o$ one assigns $F_{[n]}:=T$ and for every
morphism $f:[n]\to[m]$ of $\Delta^o$, one takes $F_f:=\one_T$, 
the identity functor of the topos $T$.

There is an obvious functor $T\to s.T$ that assigns to any
object $Z$ of $T$ the constant cosimplicial $T$-object $s.Z$ 
associated to $Z$. Especially, if $(T,V,\fm)$ is a basic setup
for $T$, then $(s.T,s.V,s.\fm)$ is a basic setup relative to
$s.T$.

\sset\subsubsection{}\label{subsec_obvious-equiv}
Suppose now that $M_\bullet$ is a cosimplicial $V^a$-module.
By applying termwise the functor $N\mapsto N_*$ we deduce
a cosimplicial $V$-module $(M_\bullet)_*$, whence an object of
$s.T$ which we denote by the same name. Clearly $(M_\bullet)_*$
is a $s.V$-module and we can therefore take its image in
the localized category $(s.V)^a\Mod$. This defines a 
functor
\set\begin{equation}\label{eq_equiv.simpl}
s.(V^a\Mod)\to(s.V)^a\Mod
\end{equation}
and it is not difficult to see that \eqref{eq_equiv.simpl}
is an equivalence of categories. Similar equivalences then
follow for categories of cosimplicial $V^a$-algebras (a.k.a.
simplicial $V^a$-schemes) and the like. For instance, for 
$(X,G,\rho)$ as in \eqref{subsec_def.G-action} we can regard
the nerve $G^\bullet_X$ as an affine $s.S$-scheme. 

\sset\subsubsection{}\label{subsec_more-general}
More generally, let $X_\bullet$ be any simplicial $S$-scheme.
A quasi-coherent $\cO_{X_\bullet}$-module is the same as
a cosimplicial $\cO_S$-module $M_\bullet$, such that 
$M_n$ is an $\cO_{X_n}$-module for every $n\in\N$,
and the coface (resp. codegeneracy) morphisms 
$\partial^i:M_n\to M_{n+1}$ (resp. $\sigma^j:M_n\to M_{n-1}$)
are $\partial_i$-linear (resp. $\sigma_j$-linear), {\em i.e.} they 
induce $M_n$-linear morphisms $M_n\to\partial_{i*}M_{n+1}$
(resp. $M_n\to\sigma_{j*}M_{n-1}$) (notation of
\eqref{eq_guess-what}). It is convenient to introduce the following
notation: for every $i,n\in\N$ we let
$$
\bar\partial{}^i:\partial_i^*M_n\to M_{n+1}
$$
the $M_{n+1}$-linear morphism deduced from $\partial^i$
by extension of scalars (and likewise for $\bar\sigma{}^j$).

\sset\subsubsection{}\label{subsec_define.pi}
For every $S$-scheme $X$ on which $G$ acts, and for every $n\in\N$,
set $\pi_{X,n}:=\partial_1\circ\partial_2\circ...\circ\partial_n:
G^n_X\to X$. For any $G$-equivariant quasi-coherent
$\cO_X$-module $M$, we define a quasi-coherent
$\cO_{G^\bullet_X}$-module $\pi_X^*M$ as follows. According
to \eqref{subsec_obvious-equiv}, this is the same as defining a
module $\pi_X^*M$ over the cosimplicial almost algebra
$\cO_{G^\bullet_X}$; then we set 
$\pi_X^*M_n:=\pi_{X,n}^*M$ for every object $[n]$ of 
$\Delta^o$. Next, we remark that 
$\pi_{X,n-1}\circ\partial_i=\pi_{X,n}$ for every $i,n>0$, hence
we have natural isomorphisms 
$\pi_X^*M_n\stackrel{\sim}{\to}\partial_i^*\pi_X^*M_{n-1}$,
from which we deduce the coface morphisms
$\partial^i:\pi_X^*M_{n-1}\to\pi_X^*M_n$ for every $i,n>0$.
Finally we use the morphism \eqref{eq_for-beta} and the cartesian
diagram:
$$
\xymatrix{ G^n_X \ar[r]^{\partial_0} 
\ar[d]_{\tau:=\partial_2\circ...\circ\partial_n} &
G^{n-1}_X \ar[d]^{\pi_{X,n-1}} \\
G^1_X \ar[r]^{\partial_0} & X
}$$
to define $\partial^0$ as the composition:
$$
\xymatrix{
M_{n-1} \ar[r] & 
\partial^*_0M_{n-1}=\partial_0^*\pi_{X,n-1}^*M=\tau^*\partial_0^*M 
\ar[r]^-{\tau^*\beta} & \pi_{X,n}^*M=M_n.
}$$
We leave to the reader the task of defining the codegeneracy
morphisms; using the cocycle relation encoded in 
\eqref{eq_cocycle-action} one can then verify the required
cosimplicial identities. In this way we obtain a functor
\set\begin{equation}\label{eq_lift.with.pi}
\cO_X\Mod_G\to\cO_{G^\bullet_X}\Mod\qquad M\mapsto\pi_X^*M.
\end{equation}

\begin{proposition} The functor \eqref{eq_lift.with.pi} is fully
faithful, and its essential image is the full subcategory of all
quasi-coherent $\cO_{G^\bullet_X}$-modules $(M_n~|~n\in\N)$ such
that 
$\bar\partial{}^n:\partial_n^*M_{n-1}\to M_n$
is an isomorphism for every $n\in\N$ (notation of
\eqref{subsec_more-general}).
\end{proposition}
\begin{proof} Let $\phi_\bullet:\pi_X^*M\to\phi_X^*N$ be a
morphism of $\cO_{G^\bullet_X}$-modules. Since $M_n=\pi_{X,n}^*M$
for every $n\in\N$, we see that $\phi_\bullet$ is already
determined by its component $\phi_0:M\to N$, whence the full
faithfulness of \eqref{eq_lift.with.pi}. On the other hand,
let $M_\bullet$ be a quasi-coherent $\cO_{G^\bullet_X}$-module
satisfying the condition of the lemma, and set $M:=M_0$.
We define a morphism 
$\beta:\partial_0^*M\to\partial_1^*M$ as
the composition $\partial^*_0M\stackrel{\bar\partial{}^0}{\to}M_1
\stackrel{(\bar\partial{}^1)^{-1}}{\to}\partial^*_1M$. 
\begin{claim} $\beta$ defines an action of $G$ on $M$.
\end{claim}
\begin{pfclaim} The identity $\sigma_0^*\beta=\one_M$ is immediate,
hence it suffices to show that \eqref{eq_cocycle-action} commutes. 
This is the same as showing the commutativity of the following diagram:
$$
\xymatrix{ & \partial^*_0\partial^*_1M=\partial^*_2\partial^*_0M
\ar[ld]_-{\partial^*_0(\bar\partial{}^1)} 
\ar[rd]^-{\partial^*_2(\bar\partial{}^0)} \\
\partial^*_0M_1 \ar[rd]^-{\bar\partial{}^0} & & 
\partial^*_2M_1 \ar[ld]_-{\bar\partial{}^2} \\
\partial_0^*\partial_0^*M \ar[u]^{\partial_0^*(\bar\partial{}^0)} \ddouble
& M_2 & \partial_2^*\partial_1^*M \ddouble 
\ar[u]_{\partial_2^*(\bar\partial{}^1)} \\
\partial_1^*\partial_0^*M \ar[r]^-{\partial^*_1(\bar\partial{}^0)} & 
\partial_1^*M_1 \ar[u]_-{\bar\partial{}^1} &
\partial_1^*\partial_1^*M \ar[l]_-{\partial_1^*(\bar\partial{}^1)}
}$$
which can be checked separately on each of its three quadrangular
subdiagrams. Let us verify for instance the commutativity of the 
bottom left diagram. By linearity, we can simplify down to the 
diagram:
$$
\xymatrix{ M_1 \ar[r]^{\partial^0} & M_2 \\
M \ar[u]^{\partial^0} \ar[r]^{\partial_0} & M_1 
\ar[u]_{\partial^1}
}$$
whose commutativity expresses one of the identities defining the
cosimplicial module $M_\bullet$.
\end{pfclaim}

To conclude the proof, it suffices to exhibit an isomorphism
$\gamma_\bullet:\pi_X^*M\stackrel{\sim}{\to}M_\bullet$. For
every $n\in\N$, let $\gamma_n:\pi_{X,n}^*M\to M_n$ be the
morphism induced by the composition
$\partial^n\circ\partial^{n-1}\circ...\circ\partial^1$; under 
our assumption, $\gamma_n$ is an isomorphism, and we leave to
the reader the verification that $(\gamma_n~|~n\in\N)$ defines
a morphism of cosimplicial modules.
\end{proof}

\sset\subsubsection{}\label{subsec_change-topos}
In case $G$ is flat over $S$, $\cI$ is $G$-equivariant,
and $\phi:X\to Y$ is an equivariant morphism of affine $S$-schemes
on which $G$ acts, we deduce a natural functor:
\set\begin{equation}\label{eq_change-topos}
\bExal_Y(X/G,\cI)\to\bExal_{G^\bullet_Y}(G^\bullet_X,\pi_X^*\cI).
\end{equation}
Namely, to any $G$-equivariant square zero deformation 
$D:=(j:X\to X',\cI,\beta)$ one assigns the datum
$G_D:=(G_j:G^\bullet_X\to G^\bullet_{X'},\pi_X^*\cI,\pi_X^*\beta)$.
The flatness of $G$ ensures that 
$\Ker(\cO_{G^n_{X'}}\to\cO_{G^n_X})\simeq\pi_{X,n}^*\cI$
for every $n\in\N$, which means that $G_D$ is indeed a
deformation of $G^\bullet_X$ over $G^\bullet_Y$.

The basic observation is contained in the following:

\begin{lemma}\label{lem_change-topos} 
Under the assumptions of \eqref{subsec_change-topos}, the functor
\eqref{eq_change-topos} is an equivalence of categories.
\end{lemma}
\begin{proof} Let $(j:G^\bullet_X\to X'_\bullet,\pi_X^*\cI,\beta)$
be a deformation of $G^\bullet_X$. For every $n\in\N$ we have a
commutative diagram of affine $S$-schemes:
$$
\xymatrix{ X'_n \ar[r]^-\partial \ar[d] & X'_0 \ar[d] \\
G^n_Y \ar[r]^-{\pi_Y} & Y
}$$
where $\partial:=\partial_1\circ\partial_2\circ...\circ\partial_n$
and $\pi_{Y,n}$ is the natural projection as in \eqref{subsec_define.pi}.
We deduce a unique morphism 
$\alpha_n:X'_n\to G^n_Y\times_YX'_0\simeq G^n_{X'_0}$.
By construction, $\alpha_n^\sharp$ fits into a commutative diagram
$$
\xymatrix{ 0 \ar[r] & \cO_{G^n}\otimes_{\cO_S}\cI \ar[r] \ddouble
& \cO_{G^n}\otimes_{\cO_S}\cO_{X'_0} \ar[d]^{\alpha_n^\sharp} \ar[r]
& \cO_{G^n}\otimes_{\cO_S}\cO_X \ar[r] \ddouble & 0
\\
0 \ar[r] & \cO_{G^n}\otimes_{\cO_S}\cI \ar[r] & \cO_{X'_n} \ar[r] &
\cO_{G^n}\otimes_{\cO_S}\cO_X \ar[r] & 0
}$$
hence $\alpha_n$ is an isomorphism for every $n\in\N$. To conclude,
it suffices to verify that the system of morphisms $(\alpha_n~|~n\in\N)$
defines a morphism of simplicial $S$-schemes: 
$\alpha_\bullet:X'_\bullet\to G^\bullet_{X'_0}$. This amounts to
showing that the $\alpha_n$ commute with the face and degeneracy
morphisms, which however is easily checked from the definition.
\end{proof}

\sset\subsubsection{}\label{subsec_comb-verbatim}
Combining lemmata \ref{lem_change-topos} and \ref{lem_eq.Ext}
(which holds {\em verbatim\/} in the present context) we derive
a natural equivalence of categories:
$$
\bExal_Y(X/G,\cI)\to
\bExal_{G^\bullet_{Y!!}}(G^\bullet_{X!!},(\pi_X^*\cI)_*).
$$
This enables us to use the usual theory of the cotangent complex to
classify the $G$-equivariant deformations of $X$. 

\sset\subsubsection{}
According to \cite[Ch.VI, \S5.3]{Il2}, for every $n\in\N$ we have 
a morphism of topoi 
$$
[n]_T:T\to s.T
$$
called {\em restriction to the $n$-th level}. It is given by
a pair of adjoint functors $([n]_T^*,[n]_{T*})$ such that
$[n]_T^*:s.T\to T$ is the functor that assigns to any cosimplicial
object $(X_{[k]}~|~k\in\N)$ the object $X_{[n]}$ of $T$.
For every $k\in\N$ set $N_{X,k}:=G^k_X\amalg\Spec\,V^a$; 
the system $(N_{X,k}~|~k\in\N)$ defines a simplicial $\Spec\,V^a$-scheme
$N_X$ (and likewise we define $N_Y$). In view of
\cite[Ch.II, (1.2.3.5)]{Il}, we deduce natural isomorphisms of
simplicial complexes of flat $\cO_{N_{X,n!!}}$-modules:
$$
[n]^*_T\L_{N_{X!!}/N_{Y!!}}\stackrel{\sim}{\to}
\L_{N_{X,n!!}/N_{Y,n!!}}\qquad\text{for every $n\in\N$}
$$
whence natural isomorphisms of simplicial complexes of flat
$\cO_{G^n_{X!!}}$-modules:
$$
[n]^*_T\L_{G^\bullet_X/G^\bullet_Y}:=
[n]_T^*\iota_\bullet^*\L_{N_{X!!}/N_{Y!!}}
\stackrel{\sim}{\to}\iota_{[n]}^*[n]_T^*\L_{N_{X!!}/N_{Y!!}}
\stackrel{\sim}{\to}\iota_{[n]}^*\L_{N_{X,n!!}/N_{Y,n!!}}=:
\L_{G^n_X/G^n_Y}
$$
where $\iota_\bullet:G^\bullet_{X!!}\to N_{X!!}$ is the morphism of
simplicial schemes defined as in \eqref{subsec_define-iota}. In other
words, $\L_{G^\bullet_X/G^\bullet_Y}$ is a mixed simplicial-cosimplicial
module $\L_{\bullet\bullet}$ whose rows $\L_{\bullet n}$ are the
cotangent complexes of the morphisms $G_\phi:G^n_X\to G^n_Y$.
Furthermore, for every $n\in\N$ and every $i\leq n+1$ we have a cartesian
diagram
$$
\xymatrix{ G^{n+1}_X \ar[r]^{\partial_i} 
\ar[d]_{G^{n+1}_\phi} & G^n_X
\ar[d]^{G^n_\phi} \\
G^{n+1}_Y \ar[r]^{\partial_i} & G^n_Y
}$$
whose horizontal arrows are flat morphisms (since $G$ is $S$-flat
by assumption). Whence, taking into account theorem
\ref{th_flat.base.ch}, a natural isomorphism in 
$\sD(\cO_{N_{X,n+1!!}}\Mod)$ :
$$
\bar\partial{}^i_\L:\partial^*_i\L_{G^n_X/G^n_Y}\stackrel{\sim}{\to}
\L_{G_X^{n+1}/G^{n+1}_Y}
$$
and, by unwinding the definitions, one sees that $\bar\partial{}^i_\L$
is induced by the $i$-th coface morphism 
$\L_{\bullet n}\to\L_{\bullet n+1}$ of the double complex
$\L_{\bullet\bullet}$.

\begin{definition}\label{def_G-torsors}
\index{Affine almost group scheme(s)!torsor over an|indref{def_G-torsors}}
Let $G$ be an affine $S$-group scheme, 
$\phi:X\to Y$ a $G$-equivariant morphism of affine $S$-schemes
on which $G$ acts. We say that $\phi$ is a 
{\em $G$-torsor over $Y$\/} if the action of $G$ on $Y$ is
trivial ({\em i.e.} $\rho:G\times_SY\to Y$ is the natural projection)
and there exists a cartesian diagram of affine $S$-schemes 
\set\begin{equation}\label{eq_trivial}{
\diagram G\times_SZ \ar[r]^-g \ar[d]_-{p_Z} & X \ar[d]^-\phi \\
Z \ar[r]^-f & Y
\enddiagram}
\end{equation}
such that $f$ is faithfully flat, $g$ is $G$-equivariant for the
action on $G\times_SZ$ deduced from $G$, and $p_Z$ is the natural
projection.
\end{definition}

\sset\subsubsection{}\label{subsec_equiv-cot-cpx}
\index{$\L^G_{X/Y,k}$}
Suppose that $G$ is a flat affine $S$-group scheme and
$X\to Y$ is a $G$-torsor over $Y$. Then the groupoid
$\xymatrix{G\times_SX\ar@<.5ex>[r] \ar@<-.5ex>[r] & X}$
defined by the $G$-action on $X$ is an effective equivalence
relation (cp. \eqref{subsec_interpret}) and $Y\simeq X/G$.
Furthermore, let $X^{\times_Yn}:=X\times_YX\times_Y....\times_YX$
(the $n$-th fold cartesian power of $X$ over $Y$);
there are natural identifications $G^n_X\simeq X^{\times_Yn+1}$,
amounting to an isomorphism of (semi-)simplicial $S$-schemes augmented
over $Y$:
$$
G^\bullet_X\stackrel{\sim}{\to}[Y|X].
$$
(cp. the notation of \cite[Exp.V$^\mathrm{bis}$, (1.2.7)]{SGA4-2}).
We can regard the cosimplicial $T$-ring $\cO_{G^\bullet_X}$
as a ring of the topos $\Gamma(\Delta\times T)$ deduced from
the simplicial topos $\Delta\times T$ (notation of 
\cite[Exp.V$^\mathrm{bis}$, (1.2.8)]{SGA4-2}), whence an
augmentation of fibred topoi
\set\begin{equation}\label{eq_augment-topos}
\theta:(\Delta\times T,\cO_{G^\bullet_X})\to (\Delta\times T,\cO_Y)
\end{equation}
and it follows from the foregoing and from remark \ref{rem_naive}
that \eqref{eq_augment-topos} is an 
{\em augmentation of $2$-cohomological descent}
(see \cite[Exp.V$^\mathrm{bis}$, D{\'e}f.2.2.6]{SGA4-2}). Denote by
$$
\bar\theta{}^*:(T,\cO_Y)\Mod\to(\Gamma(\Delta\times T),\cO_{G^\bullet_X})\Mod
$$
the morphism obtained as the composition of the constant functor
$\eps^*:(T,\cO_Y)\Mod\to(\Gamma(\Delta\times T),\cO_Y)\Mod$ and
the functor
$\Gamma(\theta^*):(\Gamma(\Delta\times T),\cO_Y)\Mod\to\cO_{G^\bullet_X}\Mod$.
Since the morphisms $\bar\partial{}^i_\L$ are isomorphisms, it then follows
from general cohomological descent 
(\cite[Exp.V$^\mathrm{bis}$, Prop.2.2.7]{SGA4-2})
that, for every $k\in\N$, the truncated system
$(\tau_{[-k}\L_{G^n_X/G^n_Y}~|~n\in\N)$ is in the essential
image of the functor $L^+\bar\theta{}^*$. In the following we will 
be only interested in the case where the cotangent complex is
concentrated in degree zero, in which case one can avoid the 
recourse to cohomological descent, and rather appeal to more
down-to-earth faithfully flat descent. In any case, the foregoing
shows that there exists a uniquely determined pro-object 
$(\L^G_{X/Y,k}~|~k\in\N)$ of $\sD(\cO_Y\Mod)$ such that
$$
L^+\bar\theta{}^*(\L^G_{X/Y,k})\simeq
\tau_{[-k}\L^a_{G^\bullet_X/G^\bullet_Y}
\qquad\text{for every $k\in\N$}.
$$
By trivial duality there follow natural isomorphisms for every
$k\in\N$ and every complex $K^\bullet\in\sD^+(\cO_{G^\bullet_X}\Mod)$
such that $K=\tau_{[-k}K$ :
\set\begin{equation}\label{eq_equate-Exts}
\hExt^k_{\cO_{G^\bullet_X}}(\L^a_{G^\bullet_X/G^\bullet_Y},K^\bullet)
\simeq\hExt^k_{\cO_Y}(\L^G_{X/Y,k},R^+\bar\theta_*K^\bullet).
\end{equation}

\begin{definition}\label{def_colie}
\index{Affine almost group scheme(s)!$\ell_{G/S}$ : colie complex
of an|indref{def_colie}}
Let $G$ be an affine $S$-scheme and
$e:S\to G$ its unit section. The {\em colie complex\/} of
$G$ is the complex of $\cO_S$-modules
$\ell_{G/S}:=e^*\L^a_{G/S}$.
\end{definition}

\begin{proposition}\label{prop_colie} 
Let $G$ be a flat $S$-group scheme, $\phi:X\to Y$ be a $G$-torsor
over $Y$ and $\pi_Y:Y\to S$ the structure morphism. Then, for every
$k\in\N$, the complex $\L^G_{X/Y,k}$ is locally isomorphic to 
$\pi_Y^*\tau_{[-k}\ell_{G/S}$ in the fpqc topology of $Y$.
\end{proposition}
\begin{proof} We will exhibit more precisely a faithfully flat 
morphism $f:Z\to Y$ such that, for every $k\in\N$ there exist
isomorphisms 
$R^+f^*\L^G_{X/Y,k}\simeq R^+f^*\pi_Y^*\tau_{[-k}\ell_{G/S}$
in $\sD(\cO_Z\Mod)$. Let $\pi_G:G\to S$ be
the structure morphism; first of all we remark:
\begin{claim}\label{cl_about-colie}
$\L^a_{G/S}\simeq\pi_G^*\ell_{G/S}$.
\end{claim}
\begin{pfclaim} Indeed, $\pi_G$ is trivially
a $G$-torsor over $S$, hence we have a compatible system of
isomorphisms $\tau_{[-k}\L^a_{G/S}\simeq\pi_Y^*\L^G_{G/S,k}$.
If $e:S\to G$ is the unit section, we deduce:
$\tau_{[-k}\ell_{G/S}\simeq e^*\tau_{[-k}\L^a_{G/S}\simeq\L^G_{G/S,k}$
for every $k\in\N$. After taking $\pi_G^*$ of the two sides, the
claim follows.
\end{pfclaim}

Let now $f:Z\to Y$ be a faithfully flat morphism that trivializes
the given $G$-torsor $\phi:X\to Y$, so we have a cartesian diagram
\eqref{eq_trivial}. Denoting by $p_G:G\times_SZ\to G$
the natural projection, we deduce (since $G$ is $S$-flat) an
isomorphism $g^*\L_{X/Y}\simeq p_G^*\L_{G/S}$. Whence, in view of
claim \ref{cl_about-colie} :
\set\begin{equation}\label{eq_whaddoyoknow}
p_Z^*f^*\L^G_{X/Y,k}\simeq p_G^*\pi_G^*\tau_{[-k}\ell_{G/S}\simeq
p_Z^*\pi_Z^*\tau_{[-k}\ell_{G/S}\simeq 
p_Z^*f^*\pi_Y^*\tau_{[-k}\ell_{G/S}.
\end{equation}
Let $e_Z:=e\times_S\one_Z:Z\to G\times_SZ$; the claim follows
after applying $e_Z^*$ to \eqref{eq_whaddoyoknow}.
\end{proof}

\sset\subsubsection{}\label{subsec_wrapup}
Finally we can wrap up this section with a discussion of equivariant
deformations of torsors. Hence, let $\phi:X\to Y$ be a $G$-torsor
over $Y$, and $j_Y:Y\to Y'$ a morphism of affine $S$-schemes such that
$\cI:=\Ker(j_Y^\sharp:\cO_{Y'}\to\cO_Y)$ is a square zero ideal of
$\cO_{Y'}$. We wish to classify the {\em square zero deformations
of the torsor $\phi$ over $Y'$}, that is, the isomorphism classes of
cartesian diagrams
$$
\xymatrix{ X \ar[d]_-\phi \ar[r]^{j_X} & X' \ar[d]^-{\phi'} \\
Y \ar[r]^{j_Y} & Y'
}$$
such that $\phi'$ is a $G$-torsor over $Y'$ and $j_X$ is
$G$-equivariant.

\begin{theorem}\label{th_wrapup}
Suppose that $G$ is flat over $S$, that $H_i(\ell_{G/S})=0$
for $i>0$ and that $H_0(\ell_{G/S})$ is an almost finitely
generated projective $\cO_S$-module. Furthermore, suppose
that the homological dimension of $\tilde\fm$ is $\leq 1$.
Then, in the situation of \eqref{subsec_wrapup}, we have:
\begin{enumerate}
\item
The pro-object $(\L^G_{X/Y,k}~|~k\in\N)$ is constant, isomorphic
to a complex of\/ $\sD(\cO_Y\Mod)$ concentrated
in degree zero, that we shall denote by $\L^G_{X/Y}$.
\item
The set of isomorphism classes of square zero deformations of
the torsor $\phi:X\to Y$ over $Y'$ is a torsor under the group
$\Ext^1_{\cO_Y}(\L^G_{X/Y},\cI)$.
\item
The group of automorphisms of a square zero deformation
$\phi':X'\to Y'$ as in \eqref{subsec_wrapup} is naturally
isomorphic to $\Ext^0_{\cO_Y}(\L^G_{X/Y},\cI)$.
\end{enumerate}
\end{theorem}
\begin{proof} (i) follows easily from proposition \ref{prop_colie};
moreover it follows from remark \ref{rem_descent}(iii) that 
$H_0(\L^G_{X/Y})^a$ is an almost finitely generated projective
$\cO_Y$-module. Let $\cJ:=\Ker(j^\sharp:\cO_{X'}\to\cO_X)$;
by flatness, the natural morphism $\phi^*\cI\to\cJ$
is a $G$-equivariant isomorphism. We notice that, for every
quasi-coherent $\cO_Y$-module $M$, there is a natural isomorphism
$$
\bar\theta{}^*M_!\stackrel{\sim}{\to}(\pi_X^*\phi^*M)_!.
$$ 
By cohomological descent (or else, by plain old-fashioned 
faithfully flat descent) it follows that the counit of the adjunction:
$$
\cI_!\to R^+\bar\theta_*\pi_X^*\cJ_!
$$
is an isomorphism, whence, in light of \eqref{eq_equate-Exts},
natural isomorphisms:
\set\begin{equation}\label{eq_inlightofup}
\hExt^k_{\cO_{G^\bullet_X}}(\L^a_{G^\bullet_X/G^\bullet_Y},\pi_X^*\cJ)
\simeq\Ext^k_{\cO_Y}(H_0(\L^G_{X/Y}),\cI)\qquad\text{for every $k\in\N$}.
\end{equation}
\begin{claim}\label{cl_diagr.of-leftists} 
The diagram of functors:
$$
\xymatrix{\cO_Y\Mod \ar[r]^-{(-)_!} \ar[d]_{\phi^*} & 
          \cO_{Y!!}\Mod \ar[d]^{\phi_{!!}^*} \\
          \cO_X\Mod \ar[r]^{(-)_!} & \cO_{X!!}\Mod
}$$
is essentially commutative.
\end{claim}
\begin{pfclaim} Indeed, both of the two possible composition
of arrows represent a functor $\cO_Y\Mod\to\cO_{X!!}\Mod$
which is left adjoint to the functor $M\mapsto\phi_*M^a$.
\end{pfclaim}
\begin{claim}\label{cl_L-G.is.good} The counit of the adjunction:
$\eps:H_0(\L^G_{X/Y})\to H_0(\L^G_{X/Y})^a_!$ is an isomorphism.
\end{claim}
\begin{pfclaim} 
It follows easily from claim \ref{cl_diagr.of-leftists} and
lemma \ref{lem_differ} that $\phi^*(\eps)$ is an isomorphism.
Since $\phi$ is faithfully flat, the assertion follows.
\end{pfclaim}

Combining \eqref{eq_equate-Exts}, claim \ref{cl_L-G.is.good}
and lemma \ref{lem_Exts.are.same}(i),(ii), we deduce that
$$
\hExt^2_{\cO_{G^\bullet_X}}
(\L_{G^\bullet_X/G^\bullet_Y},\pi_X^*\cJ)=0.
$$
However, by \eqref{subsec_comb-verbatim} and the usual arguments
(cp. section \ref{sec_lift}) one knows that the obstruction to deforming
the torsor $\phi$ is a class in the latter group; since the
obstruction vanishes, one deduces (ii). (iii) follows in the same way.
\end{proof}

\subsection{Descent}\label{sec_descent}
Faithfully flat descent in the almost setting presents no particular
surprises: since the functor $A\mapsto A_{!!}$ preserves faithful
flatness of morphisms (see remark \ref{rem_naive}) many well-known 
results for usual rings and modules extend {\em verbatim\/}
to almost algebras. 

\sset\subsubsection{}
So for instance, faithfully flat morphisms are
of universal effective descent for the fibred categories 
$F:V^a\AlgMod^o\to V^a\Alg^o$ and $G:V^a\AlgMorph^o\to V^a\Alg^o$
(see definition \ref{def_derivation}: for an almost $V$-algebra $B$, 
the fibre $F_B$ (resp. $G_B$) is the opposite of the category of 
$B$-modules (resp. $B$-algebras)).
Then, using remark \ref{rem_descent}, we deduce also universal 
effective descent for the fibred subcategories of flat (resp. almost 
finitely generated, resp. almost finitely presented, resp. almost 
finitely generated projective) modules.
Likewise, a faithfully flat morphism is of universal effective
descent for the fibred subcategories $\fibEt{}^o\to V^a\Alg^o$ of {\'e}tale 
(resp. $\fibwEt{}^o\to V^a\Alg^o$ of weakly {\'e}tale) algebras.

\sset\subsubsection{}
More generally, since the functor $A\mapsto A_{!!}$ preserves pure
morphisms in the sense of \cite{Oli}, and since, by a theorem of
Olivier ({\em loc. cit.\/}), pure morphisms are of universal effective 
descent for modules, the same holds for pure morphisms of almost 
algebras.

\sset\subsubsection{}
Non-flat descent is more delicate. Our results are not as complete
here as it could be wished, but nevertheless, they suffice for
the further study of {\'e}tale and unramified morphisms that shall
be taken up in chapter \ref{ch_hensel}.
(and they also cover the cases needed in \cite{Fa2}).
Our first statement is the almost version of a theorem of Gruson
and Raynaud (\cite[Part II, Th. 1.2.4]{Gr-Ra}).

\begin{proposition}\label{prop_descent}
A finite monomorphism of almost algebras descends flatness.
\end{proposition}
\begin{proof} Let $\phi:A\to B$ be such a morphism. Under the 
assumption, we can find a finite $A_*$-module $Q$ such that 
$\fm B_*\subset Q\subset B_*$. One sees easily that
$Q$ is a faithful $A_*$-module, so by 
\cite[Part II, Th. 1.2.4 and lemma 1.2.2]{Gr-Ra}, $Q$ satisfies 
the following condition :
\set\begin{equation}\label{eq_cond.Q}{
\parbox{14cm}{If\/ $(0\to N\to L\to P\to 0)$ is an exact sequence of 
$A_*$-modules with $L$ flat, such that $\Img(N\otimes_{A_*}Q)$ is 
a pure submodule of $L\otimes_{A_*}Q$, then $P$ is flat.}
}\end{equation}
Now let $M$ be an $A$-module such that $M\otimes_AB$ is flat.
Pick an epimorphism $p:F\to M$ with $F$ free over $A$.
Then ${\underline Y}:=
(0\to\Ker(p\otimes_A\one_B)\to F\otimes_AB\to M\otimes_AB\to 0)$
is universally exact over $B$, hence over $A$. Consider the 
sequence 
${\underline X}:=(0\to\Img(\Ker(p)_!\otimes_{A_*}Q)\to F_!\otimes_{A_*}Q
\to M_!\otimes_{A_*}Q\to 0)$.
Clearly ${\underline X}^a\simeq{\underline Y}$. However, it is easy to
check that a sequence $\underline E$ of $A$-modules is universally
exact if and only if the sequence ${\underline E}_!$ is universally exact
over $A_*$. We conclude that $\underline{X}=({\underline X}^a)_!$ is 
a universally exact sequence of $A_*$-modules, hence, by condition 
\eqref{eq_cond.Q}, $M_!$ is flat over $A_*$, {\em i.e.\/} $M$ is flat 
over $A$ as required. 
\end{proof}

\begin{corollary} Let $A\to B$ be a finite morphism of almost algebras,
with nilpotent kernel. If $C$ is a flat $A$-algebras such that 
$C\otimes_AB$ is weakly {\'e}tale (resp. {\'e}tale) over $B$, then $C$ is
weakly {\'e}tale (resp. {\'e}tale) over $A$.
\end{corollary}
\begin{proof} In the weakly {\'e}tale case, we have to show that 
the multiplication morphism $\mu:C\otimes_AC\to C$ is flat. As 
$N:=\Ker(A\to B)$ is nilpotent, the local flatness criterion reduces 
the question to the situation over $A/N$. So we may assume that 
$A\to B$ is a monomorphism. Then 
$C\otimes_AC\to(C\otimes_AC)\otimes_AB$ is a monomorphism, but 
$\mu\otimes_{C\otimes_AC}\one_{(C\otimes_AC)\otimes_AB}$
is the multiplication morphism of $C\otimes_AB$, which is flat
by assumption. Therefore, by proposition \ref{prop_descent}, $\mu$ 
is flat.

For the {\'e}tale case, we have to show that $C$ is almost finitely 
presented as a $C\otimes_AC$-module. By hypothesis $C\otimes_AB$
is almost finitely presented as a $C\otimes_AC\otimes_AB$-module
and we know already that $C$ is flat as a $C\otimes_AC$-module, so 
by lemma \ref{lem_flat.proj}(iii) (applied to the finite morphism 
$C\otimes_AC\to C\otimes_AC\otimes_AB$) the claim follows. 
\end{proof}

\sset\subsubsection{}\label{subsec_cartes.diagr}
\index{$\cD\Mod$ : modules over a diagram $\cD$ of almost 
algebras|indref{subsec_cartes.diagr}}
Suppose now that we are given a cartesian diagram $\cD$ of almost 
algebras
\set\begin{equation}\label{eq_cartesian}{
\diagram
A_0 \ar[r]^{f_2} \ar[d]_{f_1} & A_2 \ar[d]^{g_2} \\
A_1 \ar[r]^{g_1} & A_3
\enddiagram}\end{equation}
such that one of the morphisms $A_i\to A_3$ ($i=1,2$) is an 
epimorphism. Diagram $\cD$ induces an essentially commutative 
diagram for the corresponding categories $A_i\Mod$, where the 
arrows are given by the ``extension of scalars'' functors. We 
define the category of {\em $\cD$-modules\/} as the $2$-fibre product 
$\cD\Mod:=A_1\Mod\times_{A_3\Mod}A_2\Mod$. Recall (see 
\cite[Ch.VII \S 3]{Bass} or \cite[Ch.I]{Hak} for generalities 
on $2$-categories and $2$-fibre products) that $\cD\Mod$ is the 
category whose objects are the triples $(M_1,M_2,\xi)$, where 
$M_i$ is an $A_i$-module ($i=1,2$) and 
$\xi:A_3\otimes_{A_1}M_1\stackrel{\sim}{\to}A_3\otimes_{A_2}M_2$
is an $A_3$-linear isomorphism. There follows a natural functor 
$$
\pi:A_0\Mod\to\cD\Mod.
$$
Given an object $(M_1,M_2,\xi)$ of $\cD\Mod$,
let us denote $M_3:=A_3\otimes_{A_2}M_2$; we have a natural morphism
$M_2\to M_3$, and $\xi$ gives a morphism $M_1\to M_3$, so we can
form the fibre product $T(M_1,M_2,\xi):=M_1\times_{M_3}M_2$.
In this way we obtain a functor $T:\cD\Mod\to A_0\Mod$,
and we leave to the reader the verification that $T$ is right
adjoint to $\pi$. Let us denote by $\eps:\one_{\cM_0}\to T\circ\pi$
and $\eta:\pi\circ T\to\one_{\cD\Mod}$ the unit and counit of 
the adjunction.

\begin{lemma}\label{lem_adjunct} The functor $\pi$ induces 
an equivalence of full subcategories :
$$
\{X\in\mathrm{Ob}(A_0\Mod)~|~\eps_X\ \text{\rm{is an isomorphism}}\}
\stackrel{\pi}{\rightarrow}
\{Y\in\mathrm{Ob}(\cD\Mod)~|~\eta_Y\ 
\text{\rm{is an isomorphism}}\}
$$
having $T$ as quasi-inverse.
\end{lemma}
\begin{proof} General nonsense.
\end{proof}

\begin{lemma}\label{lem_eps.epi} Let $M$ be any $A_0$-module. 
Then $\eps_M$ is an epimorphism.  If $\Tor_1^{A_0}(M,A_3)=0$,
then $\eps_M$ is an isomorphism.
\end{lemma}
\begin{proof} Indeed, 
$\eps_M:M\to
(A_1\otimes_{A_0}M)\times_{A_3\otimes_{A_0}M}(A_2\otimes_{A_0}M)$
is the natural morphism. So, the assertions follow by applying 
$-\otimes_{A_0}M$ to the short exact sequence of $A_0$-modules
\set\begin{equation}\label{eq_short.ex.A_0}
0\to A_0\stackrel{f}{\to} A_1\oplus A_2\stackrel{g}{\to} A_3\to 0
\end{equation}
where $f(a):=(f_1(a),f_2(a))$ and $g(a,b):=g_1(a)-g_2(b)$.
\end{proof}

\sset\subsubsection{}\label{subsec_case.of.int}
There is another case of interest, in which $\eps_M$ is an isomorphism.
Namely, suppose that one of the morphisms $A_i\to A_3$ ($i=1,2$),
say $A_1\to A_3$, has a section. Then also the morphism 
$A_0\to A_2$ gains a section $s:A_2\to A_0$ and we have the
following :
\begin{lemma}\label{lem_special.case} In the situation
of \eqref{subsec_case.of.int}, suppose that the $A_0$-module 
$M$ arises by extension of scalars from an $A_2$-module $M'$, 
via the section $s:A_2\to A_0$. Then $\eps_M$ is an isomorphism.
\end{lemma}
\begin{proof} Indeed, in this case \eqref{eq_short.ex.A_0} 
is split exact as a sequence of $A_2$-modules, and it
remains such after tensoring by $M'$.
\end{proof}
 
\begin{lemma}\label{lem_key} $\eta_{(M_1,M_2,\xi)}$ is an 
isomorphism for all objects $(M_1,M_2,\xi)$.
\end{lemma}
\begin{proof} To fix ideas, suppose that $A_1\to A_3$ is an 
epimorphism. Consider any $\cD$-module $(M_1,M_2,\xi)$. Let 
$M:=T(M_1,M_2,\xi)$; we deduce a natural morphism 
$$\phi:(M\otimes_{A_0}A_1)\times_{M\otimes_{A_0}A_3}(M\otimes_{A_0}A_2)
\to M_1\times_{M_3}M_2$$ 
such that $\phi\circ\eps_M=\one_M$. It follows that $\eps_M$ is 
injective, hence it is an isomorphism, by lemma \ref{lem_eps.epi}.
We derive a commutative diagram with exact rows :
$$\xymatrix{
0 \ar[r] & M \ar[r] \ddouble & 
(M\otimes_{A_0}A_1)\oplus(M\otimes_{A_0}A_2) \ar[d]^{\phi_1\oplus\phi_2}
\ar[r] & M\otimes_{A_0}A_3 \ar[r] \ar[d]^{\phi_3} & 0 \\
0 \ar[r] & M \ar[r] & M_1\oplus M_2 \ar[r] & M_3 \ar[r] & 0.
}$$ 
From the snake lemma we deduce 
$$\begin{array}{ll}
(*)\qquad & \Ker(\phi_1)\oplus\Ker(\phi_2)\simeq\Ker(\phi_3) \\
(**)\qquad & \Coker(\phi_1)\oplus\Coker(\phi_2)\simeq\Coker(\phi_3).
\end{array}$$
Since $M_3\simeq M_1\otimes_{A_1}A_3$ we have 
$A_3\otimes_{A_1}\Coker\,\phi_1\simeq\Coker\,\phi_3$. But by 
assumption $A_1\to A_3$ is an epimorphism, so 
also $\Coker\,\phi_1\to\Coker\,\phi_3$ is an epimorphism. Then 
$(**)$ implies that $\Coker\,\phi_2=0$. But 
$\phi_3=\one_{A_3}\otimes_{A_2}\phi_2$, thus $\Coker\,\phi_3=0$
as well. We look at the exact sequence 
$0\to\Ker\,\phi_1\to M\otimes_{A_0}A_1\stackrel{\phi_1}{\to} M_1\to 0$ :
applying $A_3\otimes_{A_1}-$ we obtain an epimorphism 
$A_3\otimes_{A_1}\Ker\,\phi_1\to\Ker\,\phi_3$. From $(*)$ it follows
that $\Ker\,\phi_2=0$. In conclusion, $\phi_2$ is an isomorphism.
Hence the same is true for $\phi_3=\one_{A_3}\otimes_{A_2}\phi_2$,
and again $(*)$, $(**)$ show that $\phi_1$ is an isomorphism as well, 
which implies the claim.
\end{proof}

\begin{lemma}\label{lem_superFerrand} 
In the situation of \eqref{subsec_cartes.diagr}, let $M$ be any
$A_0$-module and $n\geq 1$ an integer. The following conditions 
are equivalent:
\begin{enumerate}
\renewcommand{\labelenumi}{(\alph{enumi})}
\item
$\Tor^{A_0}_j(M,A_i)=0$ for every $j\leq n$ and $i=1,2,3$.
\item
$\Tor^{A_i}_j(A_i\otimes_{A_0}M,A_3)=0$ for every $j\leq n$ and 
$i=1,2$.
\renewcommand{\labelenumi}{(\alph{enumi})}
\end{enumerate}
\end{lemma}
\begin{proof} There is a base change spectral sequence
\set\begin{equation}\label{eq_spectralseq}
E_{pq}^2:=\Tor_p^{A_i}(\Tor_q^{A_0}(M,A_i),A_3)\Rightarrow
\Tor_{p+q}^{A_0}(M,A_3).
\end{equation}
If now (a) holds, we deduce that 
$\Tor_p(A_i\otimes_{A_0}M,A_3)\simeq\Tor_p(M,A_3)$ whenever $p\leq n$,
and then (b) follows. Conversely, suppose that (b) holds; we show
(a) by induction on $n$. Say that the morphism $A_1\to A_3$ is
an epimorphism, so that the same holds for the morphism $A_0\to A_2$,
and denote by $I$ the common kernel of these morphisms. For $n=1$ 
and $i=1$, the assumption is equivalent to saying that the natural 
morphism $I\otimes_{A_1}(A_1\otimes_{A_0}M)\to(A_1\otimes_{A_0}M)$
is injective. It follows that the same holds for the morphism
$I\otimes_{A_0}M\to M$, which already shows that
$\Tor_1^{A_0}(M,A_2)=0$. Next, the assumption for $i=2$ means that
the term $E_{10}^2$ of the spectral sequence \eqref{eq_spectralseq}
vanishes, whence an isomorphism 
$\Tor_1^{A_0}(M,A_3)\simeq\Tor_1^{A_0}(M,A_2)\otimes_{A_2}A_3\simeq 0$.
Finally, we use the long exact Tor sequence arising from the short
exact sequence \eqref{eq_short.ex.A_0} to deduce that 
also $\Tor_1^{A_0}(M,A_1)$ vanishes. Let now $n>1$ and suppose that 
assertion (a) is already known for all $j<n$. We choose a 
presentation 
\set\begin{equation}\label{eq_choose-pres}
0\to R\to F\to M\to 0
\end{equation}
with $F$ flat over $A_0$; using the long exact Tor sequence we
deduce that $\Tor^{A_0}_j(M,A_i)\simeq\Tor^{A_0}_{j-1}(R,A_i)$ for
every $j>1$ and every $i\leq 3$. Moreover, assertion (a) taken with
$j=1$ shows that the sequence $A_i\otimes_{A_0}$\eqref{eq_choose-pres}
is again exact for $i\leq 3$, therefore $\Tor_1^{A_0}(R,A_i)$ 
vanishes as well for $i\leq 3$. In other words, the $A_0$-module
$R$ fulfills condition (b) for every $j<n$, hence the inductive
assumption shows that $\Tor_j^{A_0}(R,A_i)=0$ for every $j<n$
and $i=1,2,3$; in turns this implies the sought vanishing for $M$.
\end{proof}

The following lemma \ref{lem_Ferrand}(iii) generalizes a result
of D.Ferrand (\cite[lemma]{Fer}).

\begin{lemma}\label{lem_Ferrand} Let $M$ be any $A_0$-module. We have:
\begin{enumerate}
\item$
\Ann_{A_0}(M\otimes_{A_0}A_1)\cdot\Ann_{A_0}(M\otimes_{A_0}A_2)\subset
\Ann_{A_0}(M)$.
\item
$M$ admits a three-step filtration 
$0\subset\mathrm{Fil}_0M\subset\mathrm{Fil}_1M\subset\mathrm{Fil}_2M=M$
such that $\mathrm{Fil}_0M$ and $\gr_2M$ are $A_2$-modules and
$\gr_1M$ is an $A_1$-module.
\item
If $(A_1\times A_2)\otimes_{A_0}M$ is flat over $A_1\times A_2$, 
then $M$ is flat over $A_0$.
\end{enumerate}
\end{lemma}
\begin{proof} To fix ideas, suppose that $A_1\to A_3$ is an
epimorphism, and let $I$ be its kernel; let also 
$M_i:=A_i\otimes_{A_0}M$ for $i=1,2,3$.

(i): clearly $I\simeq\Ker(A_0\to A_2)$, therefore $aM\subset IM$ 
for every $a\in\Ann_{A_0}(M_2)$. On the other hand, 
the natural morphism 
$I\otimes_{A_1}M_1\stackrel{\sim}{\to}I\otimes_{A_0}M\to IM$ 
is obviously an epimorphism, whence the assertion.

(ii): we set $\mathrm{Fil}_0M:=\Ker(\eps_M:M\to M_1\times_{M_3}M_2)$; 
using the short exact sequence \eqref{eq_spectralseq} we see that
$\mathrm{Fil}_0M\simeq
\Tor_1^{A_0}(M,A_3)/\Tor_1^{A_0}(A_1\oplus A_2,M)$. Obviously $I$
annihilates $\mathrm{Fil}_0M$ (since it annihilates already
$\Tor_1^{A_0}(M,A_3)$), hence the latter is an $A_2$-module.
By lemma \ref{lem_eps.epi}, we have
$M':=M/\mathrm{Fil}_0M\simeq M_1\times_{M_3}M_2$. We can then filter
the latter module by defining 
$\mathrm{Fil}_0M':=0$, 
$\mathrm{Fil}_1M':=\Ker(M'\to M_2)\simeq\Ker(M_1\to M_3)$, which
is a $A_1$-module, and $\mathrm{Fil}_2M':=M'$. Since 
$\gr_2M'\simeq M_2$, the assertion follows.

(iii): in view of lemma \ref{lem_superFerrand}, it suffices to show 
the following:
\begin{claim} $M$ is flat over $A_0$ if and only if $M_i$ is $A_i$-flat 
and $\Tor_1^{A_0}(M,A_i)=0$ for $i\leq 2$.
\end{claim}
\begin{pfclaim} It suffices to prove the non-obvious implication,
and in view of (ii) we are reduced to showing that 
$\Tor_1^{A_0}(M,L)=0$ whenever $L$ is an $A_i$-module, for $i=1,2$.
However, for any $A_i$-module $L$ we have a base change spectral
sequence $E_{pq}^2:=\Tor_p^{A_i}(\Tor_q^{A_0}(M,A_i),L)\Rightarrow
\Tor_{p+q}^{A_0}(M,L)$. If $\Tor_1^{A_0}(M,A_i)=0$, this yields
$\Tor_1^{A_0}(M,L)\simeq\Tor_1^{A_i}(M_i,L)$, which vanishes when
$M_i$ is $A_i$-flat.
\end{pfclaim}
\end{proof}

\sset\subsubsection{}
For any $V^a$-algebra $A$, let $A\Mod_\mathrm{fl}$ (resp.
$A\Mod_\mathrm{proj}$, resp. $A\Mod_\mathrm{afpr}$) denote 
the full subcategory of $A\Mod$ consisting of all flat (resp. almost
projective, resp. almost finitely generated projective) $A$-modules.
For any integer $n\geq 1$, let $A_0\Mod_n$ be the full subcategory
of all $A_0$-modules satisfying condition (a) of lemma 
\ref{lem_superFerrand}; let also $A_i\Mod_n$ (for $i=1,2$) be
the full subcategory of $A_i\Mod$ consisting of all
$A_i$-modules $M$ such that $\Tor_j^{A_i}(M,A_3)=0$ for
every $j\leq n$. Finally, Let $A\Alg_\mathrm{fl}$ be the 
full subcategory of $A\Alg$ consisting of all flat $A$-algebras.

\begin{proposition}\label{prop_equiv.2-prod} In the situation
of \eqref{subsec_cartes.diagr}, the natural essentially commutative
diagram:
$$
\xymatrix{ A_0\Mod_\mathrm{?} \ar[r] \ar[d] & 
           A_2\Mod_\mathrm{?} \ar[d] \\
           A_1\Mod_\mathrm{?} \ar[r] & A_3\Mod_\mathrm{?}
}$$
is $2$-cartesian ({\em i.e.} cartesian in the $2$-category 
of categories, cp. \cite[Ch.I]{Hak}) whenever $?=``\mathrm{fl}"$ 
or $?=``\mathrm{proj}"$ or $?=``\mathrm{afpr}"$, or $?=n$,
for any integer  $n\geq 1$.
\end{proposition}
\begin{proof} The assertion for flat almost modules follows 
directly from lemmata \ref{lem_adjunct}, \ref{lem_eps.epi}, 
\ref{lem_key} and \ref{lem_Ferrand}(iii). Similarly the 
assertion for the categories $A_i\Mod_n$ follows from the
same lemmata and from lemma \ref{lem_superFerrand}. Set 
$B:=A_1\times A_2$. To establish the assertion for projective
modules, it suffices to show that, if $P$ is an $A_0$-module
such that $B\otimes_{A_0}P$ is almost projective over $B$,
then $P$ is almost projective over $A_0$, or which is the same,
that $\AlExt^i_{A_0}(P,N)\simeq 0$ for all $i>0$ and any 
$A_0$-module $N$. We know already that $P$ is flat. Let $M$ be
any $A_0$-module and $N$ any $B$-module. The standard isomorphism 
$R\Hom_B(B\derotimes_{A_0}M,N)\simeq R\Hom_{A_0}(M,N)$ yields 
a natural isomorphism 
$\AlExt^i_B(B\otimes_{A_0}M,N)\simeq\AlExt^i_{A_0}(M,N)$,
whenever $\Tor^{A_0}_j(B,M)=0$ for every $j>0$. In particular,
we have $\AlExt^i_{A_0}(P,N)\simeq 0$ whenever $N$ comes from 
either an $A_1$-module, or an $A_2$-module. In view of lemma
\ref{lem_Ferrand}(ii), we deduce that the sought vanishing
holds in fact for every $A_0$-module $N$. Finally, suppose that
$P\otimes_AB$ is almost finitely generated projective over $B$;
we have to show that $P$ is almost finitely generated.
To this aim, notice that 
$\Ann_{A_0}(P\otimes_AB)^2\subset\Ann_{A_0}(P)$, in view 
of lemma \ref{lem_Ferrand}(i); then the claim follows from 
remark \ref{rem_descent}(i).
\end{proof}
\begin{corollary}\label{cor_equiv.2-prod} In the situation 
of \eqref{subsec_cartes.diagr}, the natural essentially commutative
diagram:
$$
\xymatrix{ A_0\text{-}\cC \ar[r] \ar[d] & 
           A_2\text{-}\cC \ar[d] \\
           A_1\text{-}\cC \ar[r] & A_3\text{-}\cC
}$$
is $2$-cartesian whenever $\cC$ is one of the categories
$\mathbf{Alg}_\mathrm{fl}$, $\mathbf{\acute{E}t}$, 
$\mathbf{w.\acute{E}t}$, $\mathbf{\acute{E}t}_\mathrm{afp}$
(notation of \eqref{subsec_defafp}).
\qed\end{corollary}

\sset\subsubsection{}\label{subsec_desc_category}
\index{$\Desc(\cC,Y/X)$|indref{subsec_desc_category}}
Next we want to reinterpret the equivalences of proposition
\ref{prop_equiv.2-prod} in terms of descent data. 
If $F:\cC\to V^a\Alg^o$ is a fibred category over the opposite
of the category of almost algebras, and if $X\to Y$ is a
given morphism of almost algebras, we shall denote by 
$\Desc(\cC,Y/X)$ the category of objects of the fibre category
$F_Y$, endowed with a descent datum relative to the morphism
$X\to Y$ (cp. \cite[Ch.II \S1]{Gi}). In the arguments hereafter, 
we consider morphisms of almost algebras and modules, and one 
has to reverse the direction of the arrows to pass to morphisms 
in the relevant fibred category. 
Denote by $p_i:Y\to Y\otimes_XY$ ($i=1,2$), resp. 
$p_{ij}:Y\otimes_XY\to Y\otimes_XY\otimes_XY$ ($1\leq i<j\leq 3$) 
the usual morphisms. 

\sset\subsubsection{}
As an example, $\Desc(V^a\AlgMod^o,Y/X)$ 
consists of the pairs $(M,\beta)$ where $M$ is a $Y$-module 
and $\beta$ is a $Y\otimes_XY$-linear isomorphism 
$\beta:p_2^*(M)\stackrel{\sim}{\to}p_1^*(M)$ such that 
\set\begin{equation}\label{eq_cocycle}
p_{12}^*(\beta)\circ p_{23}^*(\beta)=p_{13}^*(\beta).
\end{equation} 

\sset\subsubsection{}\label{subsec_alligator}
Let now $I\subset X$ be an ideal, and set $\bar X:=X/I$, 
$\bar Y:=Y/IY$. For any $F:\cC\to V^a\Alg^o$ as in
\eqref{subsec_desc_category}, one has an essentially commutative 
diagram:
\set\begin{equation}\label{eq_Desce}{
\diagram
\Desc(\cC,Y/X) \ar[r] \ar[d] & 
\Desc(\cC,\bar Y/\bar X) \ar[d] \\
F_Y \ar[r] & F_{\bar Y}.
\enddiagram}\end{equation}
This induces a functor :
\set\begin{equation}\label{eq_Desc}
\Desc(\cC,Y/X)\to
\Desc(\cC,\bar Y/\bar X)\times_{F_{\bar Y}}F_Y.
\end{equation}
\begin{lemma}\label{lem_Desc} In the situation of
\eqref{subsec_alligator}, suppose moreover that the natural 
morphism $I\to IY$ is an isomorphism. 
Then diagram \eqref{eq_Desce} is $2$-cartesian whenever
$\cC$ is one of the fibred categories $V^a\AlgMod^o$, 
$V^a\AlgMorph^o$, $\fibEt{}^o$, $\fibwEt{}^o$.
\end{lemma}
\begin{proof} For any $n>0$, denote by $Y^{\otimes n}$ (resp.
$\bar Y{}^{\otimes n}$) the $n$-fold tensor product of 
$Y$ (resp. $\bar Y$) with itself over $X$ (resp. $\bar X$), 
and by $\rho_n:Y^{\otimes n}\to\bar Y{}^{\otimes n}$ the natural
morphism. First of all we claim that, for every $n>0$, the 
natural diagram of almost algebras
\set\begin{equation}\label{eq_Y-cart}
{\diagram
Y^{\otimes n} \ar[r]^-{\rho_n} \ar[d]_{\mu_n} & 
\bar Y{}^{\otimes n} \ar[d]^{\bar\mu_n} \\
Y \ar[r]^-{\rho_1} & \bar Y
\enddiagram}\end{equation}
is cartesian (where $\mu_n$ and $\bar\mu_n$  are $n$-fold 
multiplication morphisms). For this, we need to verify that, 
for every $n>0$, the induced morphism $\Ker\,\rho_n\to\Ker\,\rho_1$
(defined by multiplication of the first two factors)
is an isomorphism. It then suffices to check that the
natural morphism $\Ker\,\rho_n\to\Ker\,\rho_{n-1}$
is an isomorphism for all $n>1$. Indeed, consider the 
commutative diagram
$$\xymatrix{
I\otimes_XY^{\otimes n-1} \ddouble \ar[r]^-p & 
IY^{\otimes n-1} \ar[r]^-i \ar[d]_-\psi & 
Y^{\otimes n-1} \ar[d]_{\phi\otimes\one_{Y^{\otimes n-1}}} 
\ar[rrd]^{\one_{Y^{\otimes n-1}}}  \\
I\otimes_XY^{\otimes n-1} \ar[r]_-{p'} & 
\Ker\,\rho_n \ar[r]_-{i'} & 
Y^{\otimes n} \ar[rr]_-{\mu_{Y/X}\otimes\one_{Y^{\otimes n-2}}} 
& & Y^{\otimes n-1}.
}$$
From $IY=\phi(Y)$, it follows that $p'$ is an epimorphism.
Hence also $\psi$ is an epimorphism. Since $i$ is a monomorphism,
it follows that $\psi$ is also a monomorphism, hence $\psi$
is an isomorphism and the claim follows easily.

We consider first the case $\cC:=V^a\AlgMod^o$; we see that 
\eqref{eq_Y-cart} is a diagram of the kind considered in 
\eqref{eq_cartesian}, hence, for every $n>0$, we have the 
associated functor 
$\pi_n:Y^{\otimes n}\Mod\to
\bar Y{}^{\otimes n}\Mod\times_{\bar Y\Mod}Y\Mod$
and also its right adjoint $T_n$.
Denote by 
$\bar p_i:\bar Y\to\bar Y{}^{\otimes 2}$ ($i=1,2$) the usual 
morphisms, and similarly define  
$\bar p_{ij}:\bar Y{}^{\otimes 2}\to\bar Y{}^{\otimes 3}$.
Suppose there is given a descent datum $(\bar M,\bar\beta)$ 
for $\bar M$, relative to $\bar X\to\bar Y$. The 
cocycle condition \eqref{eq_cocycle} implies easily that 
$\bar\mu_2^*(\bar\beta)$ is the identity on 
$\bar\mu_2^*(\bar p_i^*\bar M)=\bar M$. It follows
that the pair $(\bar\beta,\one_M)$ defines an isomorphism
$\pi_2(p_1^*M)\stackrel{\sim}{\to}\pi_2(p_2^*M)$ in 
the category
$\bar Y{}^{\otimes 2}\Mod\times_{\bar Y\Mod}Y\Mod$.
Hence 
$T_2(\bar\beta,\one_M):T_2\circ\pi_2(p_1^*M)\to
T_2\circ\pi_2(p_2^*M)$ is an isomorphism. However, 
we remark that either morphism $\bar p_i$ yields a section 
for $\mu_2$, hence we are in the situation 
contemplated in lemma \ref{lem_special.case}, and we derive
an isomorphism $\beta:p_2^*(M)\stackrel{\sim}{\to}p_1^*(M)$.
We claim that $(M,\beta)$ is an object of $\Desc(\cC,Y/X)$, 
{\em i.e.\/} that $\beta$ verifies the cocycle condition 
\eqref{eq_cocycle}. Indeed, we can compute: 
$\pi_3(p_{ij}^*\beta)=
(\rho_3^*(p_{ij}^*\beta),\mu_3^*(p_{ij}^*\beta))$ and by
construction we have 
$\rho_3^*(p_{ij}^*\beta)=\bar p_{ij}^*(\bar\beta)$ and
$\mu_3^*(p_{ij}^*\beta)=\mu_2^*(\beta)=\one_M$. Therefore,
the cocycle identity for $\bar\beta$ implies the equality
$\pi_3(p_{12}^*\beta)\circ\pi_3(p_{23}^*\beta)=
\pi_3(p_{13}^*\beta)$. If we now apply the functor $T_3$
to this equality, and then invoke again lemma 
\ref{lem_special.case}, the required cocycle identity
for $\beta$ will ensue. Clearly $\beta$ is the only 
descent datum on $M$ lifting $\bar\beta$. This proves
that \eqref{eq_Desc} is essentially surjective. The same 
sort of argument also shows that the functor \eqref{eq_Desc} 
induces bijections on morphisms, so the lemma follows
in this case. Next, the case $\cC:=V^a\AlgMorph^o$ can
be deduced formally from the previous case, by applying
repeatedly natural isomorphisms of the kind
$p_i^*(M\otimes_YN)\simeq 
p_i^*(M)\otimes_{Y\otimes_XY}p_i^*(N)$ ($i=1,2$).
Finally, the ``{\'e}taleness'' of an object of 
$\Desc(V^a\AlgMorph^o,Y/X)$ can be checked on its
projection onto $Y\Alg^o$, hence also the cases
$\cC:=\fibwEt{}^o$ and $\cC:=\fibEt{}^o$ follow directly.
\end{proof}

\sset\subsubsection{}
Now, let $B:=A_1\times A_2$; to an object $(M,\beta)$ in 
$\Desc(V^a\AlgMod^o,B/A)$ we assign a $\cD$-module $(M_1,M_2,\xi)$ 
(notation of \eqref{subsec_cartes.diagr}) as follows. Set 
$M_i:=A_i\otimes_BM$ ($i=1,2$) and $A_{ij}:=A_i\otimes_{A_0}A_j$. 
We can write $B\otimes_{A_0}B=\prod_{i,j=1}^2A_{ij}$ and $\beta$
gives rise to the $A_{ij}$-linear isomorphisms
$\beta_{ij}:A_{ij}\otimes_{B\otimes_{A_0}B}p_2^*(M)
\stackrel{\sim}{\to}A_{ij}\otimes_{B\otimes_{A_0}B}p_1^*(M)$.
In other words, we obtain isomorphisms 
$\beta_{ij}:A_i\otimes_{A_0}M_j\to M_i\otimes_{A_0}A_j$.
However, we have a natural isomorphism $A_{12}\simeq A_3$
(indeed, suppose that $A_1\to A_3$ is an epimorphism with 
kernel $I$; then $I$ is also an ideal of $A_0$ and 
$A_0/I\simeq A_2$; now the claim follows by 
remarking that $IA_1=I$). Hence we can choose 
$\xi=\beta_{12}$. In this way we obtain a functor :
\set\begin{equation}\label{eq_funct.desc}
\Desc(V^a\AlgMod^o,B/A_0)\to\cD\Mod^o.
\end{equation}
\begin{proposition}\label{prop_equiv.desc}
The functor \eqref{eq_funct.desc} is an equivalence of 
categories.
\end{proposition}
\begin{proof} Let us say that $A_1\to A_3$ is an epimorphism
with kernel $I$. Then $I$ is also an ideal of $B$ and we have
$B/I\simeq A_3\times A_2$ and $A_0/I\simeq A_2$. We intend to
apply lemma \ref{lem_Desc} to the morphism $A_0\to B$.
However, the induced morphism 
$\bar B:=B/I\to\bar A_0:=A_0/I$ in $V^a\Alg^o$ has a section, 
and hence it is of universal effective descent for every 
fibred category. Thus, we can replace in \eqref{eq_Desc} 
the category 
$\Desc(V^a\AlgMod^o,\bar B/\bar A_0)$ by $\bar A_0\Mod^o$,
and thereby, identify (up to equivalence) the target of 
\eqref{eq_Desc} with the 2-fibred product
$(A_1\Mod\times A_2\Mod)^o\times_{(A_3\Mod\times A_2\Mod)^o}
A_2\Mod^o$. The latter is equivalent to the category 
$\cD\Mod^o$ and the resulting functor 
$\Desc(V^a\AlgMod^o,B/A_0)\to\cD\Mod^o$ is canonically 
isomorphic to \eqref{eq_funct.desc}, whence the claim.
\end{proof}

Putting together propositions \ref{prop_equiv.2-prod} and
\ref{prop_equiv.desc} we obtain the following :
\begin{corollary}\label{cor_eff-desc} In the situation of 
\eqref{eq_cartesian}, the morphism $A_0\to A_1\times A_2$ is 
of effective descent for the fibred categories of flat
modules and of almost projective modules.
\qed\end{corollary}

\sset\subsubsection{}
Next we would like to give sufficient conditions to ensure
that a morphism of almost algebras is of effective descent for
the fibred category $\fibwEt{}^o\to V^a\Alg^o$ of weakly {\'e}tale
algebras (resp. for {\'e}tale algebras). To this aim we are led 
to the following :
\begin{definition}\label{def_strictly_finite}
\index{Almost algebra(s)!strictly finite|indref{def_morph}}
A morphism $\phi:A\to B$ of almost algebras is
said to be {\em strictly finite\/} if $\Ker\,\phi$ is nilpotent
and $B\simeq R^a$, where $R$ is a finite $A_*$-algebra.
\end{definition}

\begin{theorem}\label{th_strictly} Let $\phi:A\to B$ be a strictly 
finite morphism of almost algebras. Then :
\begin{enumerate}
\item
For every $A$-algebra $C$, the induced morphism 
$C\to C\otimes_AB$ is again strictly finite.
\item
If $M$ is a flat $A$-module and $B\otimes_AM$ is
almost projective over $B$, then $M$ is almost projective 
over $A$.
\item
$A\to B$ is of universal effective descent for the fibred 
categories of weakly {\'e}tale (resp. {\'e}tale) almost algebras.
\end{enumerate}
\end{theorem}
\begin{proof} (i): suppose that $B=R^a$ for a finite $A_*$-algebra
$R$; then $S:=C_*\otimes_{A_*}R$ is a finite $C_*$-algebra and
$S^a\simeq C$. It remains to show that $\Ker(C\to C\otimes_AB)$
is nilpotent. Suppose that $R$ is generated by $n$ elements as 
an $A_*$-module and let $F_{A_*}(R)$ (resp. $F_{C_*}(S)$) be 
the Fitting ideal of $R$ (resp.of $S$); we have 
$\Ann_{C_*}(S)^n\subset F_{C_*}(S)\subset\Ann_{C_*}(S)$ (see
\cite[Ch.XIX Prop.2.5]{Lan}); on the other hand
$F_{C_*}(S)=F_{A_*}(R)\cdot C_*$, so the claim is clear.

(iii): we shall consider the fibred category  
$F:\fibwEt{}^o\to V^a\Alg^o$; the same argument applies also to 
{\'e}tale almost algebras. We begin by establishing a very 
special case :
\begin{claim}\label{cl_special} Assertion (iii) holds when 
$B=(A/I_1)\times(A/I_2)$, where $I_1$ and $I_2$ are ideals 
in $A$ such that $I_1\cap I_2$ is nilpotent.
\end{claim}
\begin{pfclaim} First of all we remark that the situation 
considered in the claim is stable under arbitrary base change, 
therefore it suffices to show that $\phi$ is of $F$-2-descent 
in this case. Then we factor $\phi$ as a composition
$A\to A/\Ker\,\phi\to B$ and we remark that $A\to A/\Ker\,\phi$
is of $F$-2-descent by theorem \ref{th_liftetale}; since a composition
of morphisms of $F$-2-descent is again of $F$-2-descent, we are reduced
to show that $A/\Ker\,\phi\to B$ is of $F$-2-descent, {\em i.e.\/} 
we can assume that $\Ker\,\phi\simeq 0$. However, in this case the
claim follows easily from corollary \ref{cor_eff-desc}.
\end{pfclaim}

\begin{claim}\label{cl_more-ideals} More generally, assertion (iii) 
holds when $B=\prod_{i=1}^nA/I_i$, where $I_1,...,I_n$ 
are ideals of $A$, such that $\bigcap^n_{i=1}I_i$ is nilpotent.
\end{claim}
\begin{pfclaim} We prove this by induction on $n$, the case $n=2$ 
being covered by claim \ref{cl_special}. Therefore, suppose that 
$n>2$, and set $B':=A/(\bigcap_{i=1}^{n-1}I_j)$. By induction,
the morphism $B'\to\prod_{i=1}^{n-1}A/I_i$ is of universal 
$F$-2-descent. However, according to \cite[Ch.II Prop.1.1.3]{Gi}, 
the sieves of universal $F$-2-descent form a topology on $V^a\Alg^o$;
for this topology, $\{A,B\}$ is a covering family of $A\times B$ and
$(A\to B'\times(A/I_n))^o$ is a covering morphism, hence $\{B',A/I_n\}$
is a covering family of $A$, and then, by composition of covering 
families, $\{\prod_{i=1}^{n-1}A/I_i,A/I_n\}$ is a covering family
of $A$, which is equivalent to the claim.
\end{pfclaim}

Now, let $A\to B$ be a general strictly finite morphism, so
that $B=R^a$ for some finite $A_*$-algebra $R$. Pick generators 
$f_1,...,f_m$ of the $A_*$-module $R$, and monic
polynomials $p_1(X),...,p_m(X)$ such that $p_i(f_i)=0$ 
for $i=1,...,m$. 
\begin{claim}\label{cl_splitter} There exists a finite and faithfully 
flat extension $C$ of $A_*$ such that the images  in $C[X]$ of 
$p_1(X)$,...,$p_m(X)$ split as products of monic linear factors. 
\end{claim}
\begin{pfclaim} This extension $C$ can be obtained as follows. 
It suffices to find, for each $i=1,...,m$, an extension $C_i$ 
that splits $p_i(X)$, because then 
$C:=C_1\otimes_{A_*}...\otimes_{A_*}C_m$ will split 
them all, so we can assume that $m=1$ and $p_1(X)=p(X)$; moreover, 
by induction on the degree of $p(X)$, it suffices to find an extension 
$C'$ such that $p(X)$ factors in $C'[X]$ as a product of the form
$p(X)=(X-\alpha)\cdot q(X)$, where $q(X)$ is a monic polynomial
of degree $\deg(p)-1$. Clearly we can take $C':=A_*[T]/(p(T))$.
\end{pfclaim}

Given a $C$ as in claim \ref{cl_splitter}, we remark that the morphism 
$A\to C^a$ is of universal $F$-2-descent. Considering again the 
topology of universal $F$-2-descent, it follows that $A\to B$ 
is of universal $F$-2-descent if and only if the same holds for the 
induced morphism $C^a\to C^a\otimes_AB$. Therefore, in proving 
assertion (iii) we can replace $\phi$ by $\one_C\otimes_A\phi$ 
and assume from start that the polynomials $p_i(X)$ factor in $A_*[X]$ 
as product of linear factors. Now, let $d_i:=\deg(p_i)$
and $p_i(X):=\prod^{d_i}_j(X-\alpha_{ij})$ (for $i=1,...,m$). 
We get a surjective homomorphism of $A_*$-algebras 
$D:=A_*[X_1,...,X_m]/(p_1(X_1),...,p_m(X_m))\to R$ by
the rule $X_i\mapsto f_i$ ($i=1,...,m$). Moreover, any 
sequence $\underline\alpha:=
(\alpha_{1,j_1},\alpha_{2,j_2},...,\alpha_{m,j_m})$
yields a homomorphism $\psi_{\underline\alpha}:D\to A_*$, 
determined by the assignment $X_i\mapsto\alpha_{i,j_i}$.
A simple combinatorial argument shows that 
$\prod_{\underline\alpha}\Ker\,\psi_{\underline\alpha}=0$, where 
$\underline\alpha$ runs over all the sequences as above.
Hence the product map
$\prod_{\underline\alpha}\psi_{\underline\alpha}:
D\to\prod_{\underline\alpha}A_*$
has nilpotent kernel. We notice that the $A_*$-algebra 
$(\prod_{\underline\alpha}A_*)\otimes_DR$ is a quotient of 
$\prod_{\underline\alpha}A_*$, hence it can be written as a 
product of rings of the form $A_*/I_{\underline\alpha}$, 
for various ideals $I_{\underline\alpha}$. By (i), the
kernel of the induced homomorphism 
$R\to\prod_{\underline\alpha}A_*/I_{\underline\alpha}$ is nilpotent,
hence the same holds for the kernel of the 
composition $A\to\prod_{\underline\alpha}A/I^a_{\underline\alpha}$, 
which is therefore of the kind considered in claim 
\ref{cl_more-ideals}. Hence  
$A\to\prod_{\underline\alpha}A/I^a_{\underline\alpha}$
is of universal $F$-2-descent. Since such 
morphisms form a topology, it follows that 
also $A\to B$ is of universal $F$-2-descent, which 
concludes the proof of (iii).

Finally, let $M$ be as in (ii) and pick again $C$ as in 
the proof of claim \ref{cl_splitter}. By remark 
\ref{rem_descent}(iv), $M$ is almost 
projective over $A$ if and only if $C^a\otimes_AM$ 
is almost projective over $C^a$; hence we can replace 
$\phi$ by $\one_{C^a}\otimes_A\phi$, and by arguing 
as in the proof of (iii), we can assume from start 
that $B=\prod_{j=1}^nA/I_j$ for ideals
$I_j\subset A$, $j=1,...,n$ such that 
$I:=\bigcap^n_{j=1}I_j$ is nilpotent. By an easy induction, 
we can furthermore reduce to the case $n=2$. We factor
$\phi$ as $A\to A/I\to B$; by proposition 
\ref{prop_equiv.2-prod} it follows that $(A/I)\otimes_AM$
is almost projective over $A/I$, and then lemma 
\ref{lem_flat.proj}(i) says that $M$ itself is almost 
projective.
\end{proof}

\begin{remark} It is natural to ask whether theorem 
\ref{th_strictly} holds if we replace everywhere ``strictly 
finite'' by ``finite with nilpotent kernel'' (or even by 
``almost finite with nilpotent kernel''). 
We do not know the answer to this question.
\end{remark}

\sset\subsubsection{}\label{subsec_finer.top}
\index{$\tau_e$, $\tau_w$ : topologies on $V^a\Alg$|indref{subsec_finer.top}}
\index{$\Idemp(R)$ : idempotent elements of a ring|indref{subsec_finer.top}}
On the category $V^a\Alg$ (taken in some universe) consider
the topologies $\tau_e$ (resp. $\tau_w$) of universal effective 
descent for the fibred category $\fibEt{}^o$ (resp. $\fibwEt{}^o$).
For a ring $R$ denote by $\Idemp(R)$ the set of idempotents
of $R$.

\begin{proposition} With the notation of \eqref{subsec_finer.top}
we have:
\begin{enumerate}
\item
The presheaf $A\mapsto\Idemp(A_*)$ 
is a sheaf for both $\tau_e$ and $\tau_w$.
\item
If $f:A\to B$ is an {\'e}tale (resp. weakly {\'e}tale) morphism 
of almost $V$-algebras and there is a covering family
$\{(A\to A_\alpha)^o\}$ for $\tau_e$ (resp. $\tau_w$) such that
$A_\alpha\to A_\alpha\otimes_AB$ is an almost projective
epimorphism for all $\alpha$, then $f$ is an almost projective 
epimorphism.
\item
$\tau_e$ is finer than $\tau_w$.
\end{enumerate}
\end{proposition}
\begin{proof} (i): use descent of morphisms and the bijection 
$$\Hom_{A\Alg}(A\times A,A)\stackrel{\sim}{\to}\Idemp(A_*)\qquad
\phi\mapsto\phi_*(1,0).$$
(ii): by remark \ref{rem_idemp}, 
$\Ker(A_\alpha\to A_\alpha\otimes_AB)$ is generated by 
$e_\alpha\in\Idemp(A_{\alpha *})$. $e_\alpha$ and $e_\beta$ agree 
in $(A_\alpha\otimes_AA_\beta)_*$, so by (i) there is an idempotent
$e\in A_*$ that restricts to $e_\alpha$ in $\Idemp(A_{\alpha *})$,
for each $\alpha$. The $A$-algebras $B$ and $A/eA$ become
isomorphic after applying $-\otimes_AA_\alpha$; these isomorphisms
are unique and are compatible on $A_\alpha\otimes_AA_\beta$, hence
they patch to an isomorphism $B\simeq A/eA$.

(iii): we have to show that if $A$ is an almost $V$-algebra, $R$
a sieve of $(V^a\Alg)^o/A$ and $R$ is of universal $\fibwEt{}^o$-2-descent,
then $R$ is of universal $\fibEt{}^o$-2-descent. Since the assumption
is stable under base change, it suffices to show that $R$ is of
$\fibEt{}^o$-2-descent. Descent of morphisms is clear. Let $R$ be
the sieve generated by a family of morphisms $\{(A\to A_\alpha)^o\}$. Any
descent datum consisting of {\'e}tale $A_\alpha$-algebras $B_\alpha$
and isomorphisms $A_\alpha\otimes_AB_\beta\simeq B_\alpha\otimes_AA_\beta$
satisfying the cocycle condition, becomes effective when we pass
to $\fibwEt{}^o$. So one has to verify that if $B$ is a weakly {\'e}tale
$A$-algebra such that $B\otimes_AA_\alpha$ is {\'e}tale over $A_\alpha$
for all $\alpha$, then $B$ is {\'e}tale over $A$. Indeed, an application
of (ii) gives that $B\otimes_AB\to B$ is almost projective. 
\end{proof}

\sset\subsubsection{}\label{subsec_digression}
\index{$\fibuEt^o$|indref{subsec_digression}}
We conclude with a digression to explain the relationship
between our results and known facts that can be extracted
from the literature. So, we now place ourselves in the 
``classical limit'' $\fm:=V$ (cp. example \ref{ex_rings}(ii)). 
In this case, weakly {\'e}tale morphisms had already been 
considered in some earlier work, and they were called 
``absolutely flat'' morphisms. A ring homomorphism $A\to B$
is {\'e}tale in the usual sense of \cite{SGA1} if and only if 
it is absolutely flat and of finite presentation. Let us denote
by $\fibuEt{}^o$ the fibred category over $V\Alg^o$, whose fibre
over a $V$-algebra $A$ is the opposite of the category of 
{\'e}tale $A$-algebras in the usual sense. We claim that, if 
a morphism $A\to B$ of $V$-algebras is of universal effective 
descent for the fibred category $\fibwEt{}^o$ (resp. $\fibEt{}^o$), 
then it is a morphism of universal effective descent for 
$\fibuEt{}^o$. Indeed, let $C$ be an {\'e}tale $A$-algebra (in the 
sense of definition \ref{def_morph}) and such that $C\otimes_AB$ 
is {\'e}tale over $B$ in the usual sense. We have to show that $C$
is {\'e}tale in the usual sense, {\em i.e.\/} that it is of
finite presentation over $A$. This amounts to showing that,
for every filtered inductive system 
$(A_\lambda)_{\lambda\in\Lambda}$ of $A$-algebras, we have
$\colim{\lambda\in\Lambda}\Hom_{A\Alg}(C,A_\lambda)\simeq
\Hom_{A\Alg}(C,\colim{\lambda\in\Lambda}A_\lambda)$. Since, 
by assumption, this is known after extending scalars to $B$
and to $B\otimes_AB$, it suffices to show that, for any 
$A$-algebra $D$, the natural sequence
$$\xymatrix{\Hom_{A\Alg}(C,D)\ar[r] & \Hom_{B\Alg}(C_B,D_B)
\ar@<.5ex>[r] \ar@<-.5ex>[r] &
\Hom_{B\otimes_AB\Alg}(C_{B\otimes_AB},D_{B\otimes_AB})
}$$
is exact. For this, note that 
$\Hom_{A\Alg}(C,D)=\Hom_{D\Alg}(C_D,D)$ (and similarly for 
the other terms) and by hypothesis $(D\to D\otimes_AB)^o$ is
a morphism of 1-descent for the fibred category $\fibwEt{}^o$ 
(resp. $\fibEt{}^o$). 

As a consequence of these observations and 
of theorem \ref{th_strictly}, we see that any finite ring 
homomorphism $\phi:A\to B$ with nilpotent kernel is of universal 
effective descent for the fibred category of {\'e}tale algebras. 
This fact was known as follows. By \cite[Exp.IX, 4.7]{SGA1}, 
$\Spec(\phi)$ is of universal effective descent for the fibred 
category of separated {\'e}tale morphisms of finite type. One 
has to show that if $X$ is such a scheme over $A$, such that 
$X\otimes_AB$ is affine, then $X$ is affine. This follows by 
reduction to the noetherian case and \cite[Ch.II, 6.7.1]{EGA}.

\subsection{Behaviour of {\'e}tale morphisms under Frobenius}
\label{subsec_behaviour}
\index{$\cB$ : category of basic setups|indref{subsec_behaviour}}
We consider the following category $\cB$ of basic setups. The 
objects of $\cB$ are the pairs $(V,\fm)$, where $V$ is a ring
and $\fm$ is an ideal of $V$ with $\fm=\fm^2$ and $\tilde\fm$ 
is flat. The morphisms $(V,\fm_V)\to(W,\fm_W)$ between two 
objects of $\cB$ are the ring homomorphisms $f:V\to W$ such 
that $\fm_W=f(\fm_V)\cdot W$.

\sset\subsubsection{}\label{subsec_cB-modules}
\index{$\cB\Mod$ : category of $\cB$-modules|indref{subsec_cB-modules}}
We have a fibred and cofibred category $\cB\Mod\to\cB$ (see 
\cite[Exp.VI \S5,6,10]{SGA1} for generalities on fibred
categories). An object of $\cB\Mod$ (which we may call a 
``$\cB$-module'') consists of a pair $((V,\fm),M)$, where $(V,\fm)$
is an object of $\cB$ and $M$ is a $V$-module. Given two objects
$X:=((V,\fm_V),M)$ and $Y:=((W,\fm_W),N)$, the morphisms $X\to Y$
are the pairs $(f,g)$, where $f:(V,\fm_V)\to(W,\fm_W)$ is a 
morphism in $\cB$ and $g:M\to N$ is an $f$-linear map.

\sset\subsubsection{}\label{subsec_cB-algebra}
\index{$\cB$-monoid(s)|indref{subsec_cB-algebra}}
\index{$\cB\Alg$ : category of $\cB$-algebras|indref{subsec_cB-algebra}}
\index{$\cB$-monoid(s)!$\cB\Mon$ : category of|indref{subsec_cB-algebra}}
\index{$V\Mon$ : category of $V$-monoids|indref{subsec_cB-algebra}}
Similarly one has a fibred and cofibred category $\cB\Alg\to\cB$
of $\cB$-algebras. We will also need to consider the fibred and
cofibred category $\cB\Mon\to\cB$ of non-unitary commutative
$\cB$-monoids: an object of $\cB\Mon$ is a pair $((V,\fm),A)$ where 
$A$ is a $V$-module endowed with a morphism $A\otimes_VA\to A$
subject to associativity and commutativity conditions, as 
discussed in section \ref{sec_alm.cat}. The fibre over
an object $(V,\fm)$ of $\cB$, is the category of $V$-monoids
denoted $(V,\fm)\Mon$ or simply $V\Mon$. 

\sset\subsubsection{}\label{subsec_locally-small}
\index{$\cB^a\Mod$, $\cB^a\Alg$, $\cB^a\Mon$|indref{subsec_locally-small}}
The almost isomorphisms in the fibres of $\cB\Mod\to\cB$ give a 
multiplicative system $\Sigma$ in $\cB\Mod$, admitting a calculus
of both left and right fractions. The ``locally small''
conditions are satisfied (see \cite[p.381]{We}), so that one 
can form the localised category $\cB^a\Mod:=\Sigma^{-1}(\cB\Mod)$. 
The fibres of the localised category over the objects of 
$\cB$ are the previously considered categories of almost 
modules. Similar considerations hold for $\cB\Alg$ and $\cB\Mon$, 
and we get the fibred and cofibred categories $\cB^a\Mod\to\cB$, 
$\cB^a\Alg\to\cB$ and $\cB^a\Mon\to\cB$. In particular, for 
every object $(V,\fm)$ of $\cB$, we have an obvious notion of
almost $V$-monoid and the category consisting of these is 
denoted $V^a\Mon$. 

\sset\subsubsection{}
The localisation functors 
$$
\text{$\cB\Mod\to\cB^a\Mod\ :\ M\mapsto M^a \qquad 
       \cB\Alg\to\cB^a\Alg\ :\ B\mapsto B^a$}
$$
have left and right adjoints. These adjoints can be chosen
as functors of categories over $\cB$ such that the adjunction
units and counits are morphisms over identity arrows in $\cB$.
On the fibres these induce the previously considered left and 
right adjoints $M\mapsto M_!$, $M\mapsto M_*$, $B\mapsto B_{!!}$, 
$B\mapsto B_*$.
We will use the same notation for the corresponding functors 
on the larger categories. Then it is easy to check that the 
functor $M\mapsto M_!$ is cartesian and cocartesian ({\em i.e.\/}
it sends cartesian arrows to cartesian arrows and cocartesian
arrows to cocartesian arrows), $M\mapsto M_*$ and $B\mapsto B_*$
are cartesian, and $B\mapsto B_{!!}$ is cocartesian.

\sset\subsubsection{}\label{subsec_full-subcat}
\index{$\cB/\F_p$, $\cB\Alg/\F_p$, $\cB\Mon/\F_p$, $\cB^a\Alg/\F_p$
|indref{subsec_full-subcat}}
Let $\cB/\F_p$ be the full subcategory of $\cB$ consisting of
all objects $(V,\fm)$ where $V$ is an $\F_p$-algebra. Define 
similarly $\cB\Alg/\F_p$, $\cB\Mon/\F_p$ and $\cB^a\Alg/\F_p$, 
$\cB^a\Mon/\F_p$, so that we have again fibred and cofibred categories 
$\cB^a\Alg/\F_p\to\cB/\F_p$ and $\cB^a\Alg/\F_p\to\cB/\F_p$ (resp.
the same for non-unitary monoids). We remark that the categories
$\cB^a\Alg/\F_p$ and $\cB^a\Mon/\F_p$ have small limits and colimits, 
and these are preserved by the projection to $\cB/\F_p$. Especially,
if $A\to B$ and $A\to C$ are two morphisms in $\cB^a\Alg/\F_p$
or $\cB^a\Mon/\F_p$, we can define $B\otimes_AC$ as such a colimit.

\sset\subsubsection{}\label{subsec_Frobenius}
\index{$\cB$-monoid(s)!$\Phi_A$ : Frobenius endomorphism of|indref{subsec_Frobenius}}
If $A$ is a (unitary or non-unitary) $\cB$-monoid over $\F_p$, 
we denote by $\Phi_A:A\to A$ the {\em Frobenius endomorphism\/}: 
$x\mapsto x^p$. If $(V,\fm)$ is an object of $\cB/\F_p$, it 
follows from proposition \ref{prop_less.obv}(ii) that 
$\Phi_V:(V,\fm)\to(V,\fm)$ is a morphism in $\cB$. 
If $B$ is an object of $\cB\Alg/\F_p$ (resp. $\cB\Mon/\F_p$) 
over $V$, then the Frobenius map induces a morphism 
$\Phi_B:B\to B$ in $\cB\Alg/\F_p$ (resp. $\cB\Mon/\F_p$) over 
$\Phi_V$. In this way we get a natural transformation from the 
identity functor of $\cB\Alg/\F_p$ (resp. $\cB\Mon/\F_p$) to 
itself that induces a natural transformation on the identity 
functor of $\cB^a\Alg/\F_p$ (resp. $\cB^a\Mon/F_p$).

\sset\subsubsection{}\label{subsec_pull_by_Frob}
\index{$B_{(m)}$|indref{subsec_pull_by_Frob}}
Using the pull-back functors, any object $B$ of $\cB\Alg$ over 
$V$ defines new objects $B_{(m)}$ of $\cB\Alg$ ($m\in\N$) over 
$V$, where $B_{(m)}:=(\Phi^m_V)^*(B)$, which is just $B$
considered as a $V$-algebra via the homomorphism 
$V\stackrel{\Phi^m}{\longrightarrow} V\to B$. These 
operations also induce functors $B\mapsto B_{(m)}$ on 
almost $\cB$-algebras.

\begin{definition}\label{def_invert.phi} 
\index{$\cB$-monoid(s)!invertible morphism up to $\Phi^m$ of
|indref{def_invert.phi}}
\index{$\cB$-monoid(s)!Frobenius nilpotent|indref{def_invert.phi}}
Let $(V,\fm)$ be an object of $\cB/\F_p$.
\begin{enumerate}
\item
We say that a morphism $f:A\to B$ of almost $V$-algebras 
(resp. almost $V$-monoids) is {\em invertible up to\/} $\Phi^m$ 
if there exists a morphism $f':B\to A$ in $\cB^a\Alg$ (resp. 
$\cB^a\Mon$) over $\Phi^m_V$, such that $f'\circ f=\Phi^m_A$ 
and $f\circ f'=\Phi^m_B$.
\item
We say that an almost $V$-monoid $I$ ({\em e.g.\/} an
ideal in a $V^a$-algebra) is {\em Frobenius nilpotent\/}
if $\Phi_I$ is nilpotent.
\end{enumerate}
\end{definition}

\sset\subsubsection{}
Notice that a morphism $f$ of $V^a\Alg$ (or $V^a\Mon$) is 
invertible up to $\Phi^m$ if and only if $f_*:A_*\to B_*$ 
is so as a morphism of $\F_p$-algebras.

\begin{lemma}\label{lem_invert} Let $(V,\fm)$ be an 
object of $\cB/\F_p$ and let $f:A\to B$, $g:B\to C$ 
be  morphisms of almost $V$-algebras or almost $V$-monoids. 
\begin{enumerate}
\item
If $f$ (resp. $g$) is invertible up to $\Phi^n$ (resp. $\Phi^m$), 
then $g\circ f$ is invertible up to $\Phi^{m+n}$.
\item
If $f$ (resp. $g\circ f$) is invertible up to $\Phi^n$ 
(resp. $\Phi^m$), then $g$ is invertible up to $\Phi^{m+n}$.
\item
If $g$ (resp.$g\circ f$) is invertible up to $\Phi^n$ (resp. 
$\Phi^m$), then $f$ is invertible up to $\Phi^{m+n}$.
\item
The Frobenius morphisms induce $\Phi_V$-linear morphisms
({\em i.e.\/} morphisms in $\cB^a\Mod$ over $\Phi_V$)
$\Phi':\Ker\,f\to\Ker\,f$ and 
$\Phi'':\Coker\,f\to\Coker\,f$, and $f$ is
invertible up to some power of $\Phi$ if and only
if both $\Phi'$ and $\Phi''$ are nilpotent.
\item
Consider a map of short exact sequences of almost 
$V$-monoids :
$$\xymatrix{0 \ar[r] & A' \ar[r] \ar[d]_{f'} & 
                       A \ar[r] \ar[d]_f &
                       A'' \ar[r] \ar[d]_{f''} & 0 \\
            0 \ar[r] & B' \ar[r] & B \ar[r] & B'' \ar[r] & 0
}$$
and suppose that two of the morphisms $f',f,f''$ are
invertible up to a power of $\Phi$. Then also the third
morphism has this property.
\end{enumerate}
\end{lemma}
\begin{proof} (i): if $f'$ is an inverse of $f$ up to
$\Phi^n$ and $g'$ is an inverse of $g$ up to $\Phi^m$,
then $f'\circ g'$ is an inverse of $g\circ f$ up to
$\Phi^{m+n}$. (ii): given an inverse $f'$ of $f$ up 
to $\Phi^n$ and an inverse $h'$ of $h:=g\circ f$ up to 
$\Phi^m$, let $g':=\Phi^n_B\circ f\circ h'$. We compute :
$$\begin{array}{l@{\: =\:}l}
g\circ g' & g\circ\Phi^n_B\circ f\circ h'=
\Phi^n_C\circ g\circ f\circ h=\Phi^n_C\circ\Phi^m_C \\
g'\circ g & \Phi^n_B\circ f\circ h'\circ g=
f\circ h'\circ g\circ\Phi^n_B=
f\circ h'\circ g\circ f\circ f'\\
& f\circ\Phi^m_A\circ f'=\Phi^m_B\circ f\circ f'=
\Phi^m_B\circ\Phi^n_B.
\end{array}$$
(iii) is similar and (iv) is an easy diagram chasing
left to the reader. (v) follows from (iv) and the snake 
lemma.
\end{proof}

\begin{lemma}\label{lem_pull.and.push}
Let $(V,\fm)$ be an object of $\cB/\F_p$.
\begin{enumerate}
\item
If $f:A\to B$ is a morphism of almost $V$-algebras, 
invertible up to $\Phi^n$, then so is $A'\to A'\otimes_AB$ 
for every morphism $A\to A'$ of almost algebras.
\item
If $f:(V,\fm_V)\to(W,\fm_W)$ is a morphism in 
$\cB/\F_p$, the functors 
$$f_*:(V,\fm_V)^a\Alg\to(W,\fm_W)^a\Alg\quad\text{and}\quad
f^*:(W,\fm_W)^a\Alg\to(V,\fm_V)^a\Alg$$ 
preserve the class of morphisms invertible up to $\Phi^n$.
\end{enumerate}
\end{lemma}
\begin{proof} (i): given $f':B\to A_{(m)}$, construct a 
morphism $A'\otimes_AB\to A'_{(m)}$ using the morphism 
$A'\to A'_{(m)}$ coming from $\Phi^m_{A'}$ and $f'$.
(ii): the assertion for $f^*$ is clear, and the assertion
for $f_*$ follows from (i).
\end{proof}

\begin{remark}\label{rem_pull.push} Statements like those 
of lemma \ref{lem_pull.and.push} hold for the classes of 
flat, (weakly) unramified, (weakly) {\'e}tale morphisms.
\end{remark}

\begin{theorem}\label{th_Frobenius} 
Let $(V,\fm)$ be an object of $\cB/\F_p$ and 
$f:A\to B$ a weakly {\'e}tale morphism of almost 
$V$-algebras.
\begin{enumerate}
\item
If $f$ is invertible up to $\Phi^n$ ($n\geq 0$), 
then it is an isomorphism.
\item
For every integer $m\geq 0$ the natural square 
diagram
\set\begin{equation}\label{eq_Frobenius}{
\diagram A \ar[r]^f \ar[d]_{\Phi^m_A} & 
         B \ar[d]^{\Phi^m_B} \\
         A_{(m)} \ar[r]^{f_{(m)}} & B_{(m)}
\enddiagram
}\end{equation}
is cocartesian.
\end{enumerate}
\end{theorem}
\begin{proof} (i): we first show that $f$ is faithfully 
flat. Since $f$ is flat, it remains to show that if $M$ 
is an $A$-module such that $M\otimes_AB=0$, then 
$M=0$. It suffice to do this for $M:=A/I$, for an arbitrary 
ideal $I$ of $A$. After base change by $A\to A/I$, we 
reduce to show that $B=0$ implies $A=0$. However, 
$A_*\to B_*$ is invertible up to $\Phi^n$, so 
$\Phi^n_{A_*}=0$ which means $A_*=0$. In particular, 
$f$ is a monomorphism, hence the proof is complete in 
case that $f$ is an epimorphism. In general, consider 
the composition
$B\stackrel{\one_B\otimes f}{\longrightarrow}
B\otimes_AB\stackrel{\mu_{B/A}}{\longrightarrow}B$.
From lemma \ref{lem_pull.and.push}(i) it follows that 
$\one_B\otimes f$ is invertible up to $\Phi^n$; then 
lemma \ref{lem_invert}(ii) says that $\mu_{B/A}$ is 
invertible up to $\Phi^n$. The latter is also weakly
{\'e}tale; by the foregoing we derive that it is an
isomorphism. Consequently $\one_B\otimes f$ is an
isomorphism, and finally, by faithful flatness, $f$
itself is an isomorphism.

(ii): the morphisms $\Phi^m_A$ and $\Phi^m_B$ are 
invertible up to $\Phi^m$. By lemma 
\ref{lem_pull.and.push}(i) it follows that 
$\one_B\otimes\Phi^m_A:B\to B\otimes_AA_{(m)}$ is
invertible up to $\Phi^m$; hence, by lemma 
\ref{lem_invert}(ii), the morphism 
$h:B\otimes_AA_{(m)}\to B_{(m)}$ induced by $\Phi^m_B$
and $f_{(m)}$ is invertible up to $\Phi^{2m}$ (in fact 
one verifies that it is invertible up to $\Phi^m$). 
But $h$ is a morphism of weakly {\'e}tale $A_{(m)}$-algebras, 
so it is weakly {\'e}tale, so it is an isomorphism by (i).
\end{proof}

\begin{remark} Theorem \ref{th_Frobenius}(ii) extends
a statement of Faltings (\cite[p.10]{Fa2}) for his notion
of almost {\'e}tale extensions. 
\end{remark} 

\sset\subsubsection{}\label{subsec_bicovering}
\index{bicovering morphism in a site|indref{subsec_bicovering}}
We recall (cp. \cite[Ch.0, 3.5]{Gi}) that a morphism 
$f:X\to Y$ of objects in a site is called {\em bicovering\/} 
if the induced map of associated sheaves of sets is an 
isomorphism; if $f$ is squarable (``quarrable'' in French), 
this is equivalent to the condition that both $f$ and the 
diagonal morphism $X\to X\times_YX$ are covering morphisms.

\sset\subsubsection{}
Let $F\to E$ be a fibered category and $f:P\to Q$ a
squarable morphism of $E$. Consider the following
condition:
\set\begin{equation}\label{eq_bicover}{
\parbox{14cm}
{for every base change $P\times_QQ'\to Q'$ of $f$, the 
inverse image functor $F_{Q'}\to F_{P\times_QQ'}$ is 
an equivalence of categories.}
}\end{equation}
Inspecting the arguments in \cite[Ch.II, \S 1.1]{Gi}
one can show:
\begin{lemma}\label{lem_bicover} With the above notation, 
let $\tau$ be the topology of  universal effective descent 
relative to $F\to E$. Then we have :
\begin{enumerate}
\item
if \eqref{eq_bicover} holds, then $f$ is a covering morphism 
for the topology $\tau$.
\item
$f$ is bicovering for $\tau$ if and only if \eqref{eq_bicover} 
holds both for $f$ and for the diagonal morphism $P\to P\times_QP$. 
\qed\end{enumerate}
\end{lemma}
\begin{remark} In \cite[Ch.II, 1.1.3(iv)]{Gi} it
is stated that ``la r{\'e}ciproque est vraie si $i=2$'', 
meaning that \eqref{eq_bicover} is equivalent to the condition 
that $f$ is bicovering for $\tau$. (Actually the cited statement
is given in terms of presheaves, but one can show that
\eqref{eq_bicover} is equivalent to the corresponding condition 
for the fibered category $F^+\to\hat E_U$ considered in 
{\em loc.cit.}) However, this fails in general : as a counterexample 
we can give the following. Let $E$ be the category of
schemes of finite type over a field $k$; set $P=\A^1_k$,
$Q=\Spec\,k$. Finally let $F\to E$ be the discretely
fibered category defined by the presheaf 
$X\mapsto H^0(X,\Z)$. Then it is easy to show that $f$
satisfies \eqref{eq_bicover} but the diagonal map does not, 
so $f$ is not bicovering. The mistake in the proof is in 
\cite[Ch.II, 1.1.3.5]{Gi}, where one knows that 
$F^+(d)$ is an equivalence of categories (notation 
of {\em loc.cit.}) but one needs it also after base 
changes of $d$. 
\end{remark}

\begin{lemma}\label{lem_invert.Et} (i)
Let $f:A\to B$ be a morphism of $V^a$-algebras.
\begin{enumerate}
\item
If $f$ is invertible up to $\Phi^m$, then the induced functors 
$A\Et\to B\Et$ and $A\wEt\to B\wEt$ are equivalences 
of categories.
\item
If $f$ is weakly {\'e}tale and $C\to D$ is a morphism 
of $A$-algebras invertible up to $\Phi^m$, then the induced
map: $\Hom_{A\Alg}(B,C)\to\Hom_{A\Alg}(B,D)$ is bijective.
\end{enumerate}
\end{lemma} 
\begin{proof} We first consider (i) for the special case 
where $f:=\Phi^m_A:A\to A_{(m)}$. The functor 
$(\Phi^m_V)^*:V^a\Alg\to V^a\Alg$ induces a functor
$(-)_{(m)}:A\Alg\to A_{(m)}\Alg$, and by restriction 
(see remark \ref{lem_pull.and.push}) we obtain a functor 
$(-)_{(m)}:A\Et\to A_{(m)}\Et$; by theorem 
\ref{th_Frobenius}(ii), the latter is isomorphic to
the functor $(\Phi^m)_*:A\Et\to A_{(m)}\Et$ of the 
lemma. Furthermore, from remark \ref{rem_almost.zero}(ii) and
\eqref{eq_alm.morph} we derive a natural ring isomorphism 
$\omega:A_{(m)*}\simeq A_*$, hence an essentially
commutative diagram
$$\xymatrix{
A\Et \ar[r] \ar[d]_{(\Phi^m)_*} & 
A\Alg \ar[r]^-\alpha \ar[d]_{(-)_{(m)}} & 
(A_*,\fm\cdot A_*)^a\Alg \ar[d]_{\omega^*} \\
A_{(m)}\Et \ar[r] & A_{(m)}\Alg \ar[r]^-\beta &
(A_{(m)*},\fm\cdot A_{(m)*})^a\Alg 
}$$
where $\alpha$ and $\beta$ are the equivalences 
of remark \ref{rem_base.comp}. Clearly $\alpha$
and $\beta$ restrict to equivalences on the 
corresponding categories of {\'e}tale algebras, 
hence the lemma follows in this case. 

For the general case of (i), let $f':B\to A_{(m)}$ be
a morphism as in definition \ref{def_invert.phi}.
Diagram \eqref{eq_Frobenius} induces an essentially
commutative diagram of the corresponding categories
of algebras, so by the previous case, the functor
$(f')_*:B\Et\to A_{(m)}\Et$ has both a left
quasi-inverse and a right quasi-inverse;
these quasi-inverses must be isomorphic, so
$f_*$ has a quasi-inverse as desired.
Finally, we remark that the map in (ii) is the same
as the map 
$\Hom_{C\Alg}(B\otimes_AC,C)\to\Hom_{D\Alg}(B\otimes_AD,D)$,
and the latter is a bijection in view of (i).
\end{proof}
\begin{remark} Notice that lemma \ref{lem_invert.Et}(ii)
generalises the lifting theorem \ref{th_liftetale}(i) (in 
case $V$ is an $\F_p$-algebra). Similarly, it follows 
from lemmata \ref{lem_invert.Et}(i) and \ref{lem_invert}(iv) 
that, in case $V$ is an $\F_p$-algebra, one can replace 
``nilpotent'' in theorem \ref{th_liftetale}(ii),(iii) 
by ``Frobenius nilpotent''. 
\end{remark}

\sset\subsubsection{}
In the following, $\tau$ will denote indifferently the topology 
of universal effective descent defined by either of the fibered 
categories $\fibwEt{}^o\to V^a\Alg^o$ or $\fibEt{}^o\to V^a\Alg^o$.
\begin{proposition}\label{prop_invert.bicov}
If $f:A\to B$ is a morphism of almost 
$V$-algebras which is invertible up to $\Phi^m$, then $f^o$\/
is bicovering for the topology $\tau$.
\end{proposition}
\begin{proof} In light of lemmata \ref{lem_bicover}(ii) and 
\ref{lem_invert.Et}(i), it suffices to show that $\mu_{B/A}$ 
is invertible up to a power of $\Phi$. For this, factor the 
identity morphism of $B$ as 
$B\stackrel{\one_B\otimes f}{\longrightarrow}
B\otimes_AB\stackrel{\mu_{B/A}}{\longrightarrow}B$ and argue 
as in the proof of theorem \ref{th_Frobenius}.
\end{proof}
\begin{proposition} Let $A\to B$ be a morphism of almost 
$V$-algebras and $I\subset A$ an ideal. Set $\bar A:=A/I$ 
and $\bar B:=B/IB$. Suppose that either
\renewcommand{\labelenumi}{(\alph{enumi})}
\begin{enumerate}
\item
$I\to IB$ is an epimorphism with nilpotent kernel, or
\item
$V$ is an $\F_p$-algebra and $I\to IB$ is invertible
up to a power of $\Phi$.
\end{enumerate}
Then we have :
\renewcommand{\labelenumi}{(\roman{enumi})}
\begin{enumerate}
\item
conditions {\em (a)} and {\em (b)} are stable under any base change 
$A\to C$.
\item
$(A\to B)^o$ is covering (resp. bicovering) for $\tau$
if and only if $(\bar A\to\bar B)^o$ is.
\end{enumerate}
\end{proposition}
\begin{proof} Suppose first that $I\to IB$ is an isomorphism;
in this case we claim that $IC\to I(C\otimes_AB)$ is
an epimorphism and  $\Ker(IC\to I(C\otimes_AB))^2=0$ for 
any $A$-algebra $C$. Indeed, since by assumption $I\simeq IB$, 
$C\otimes_AB$ acts on $C\otimes_AI$, hence $\Ker(C\to C\otimes_AB)$
annihilates $C\otimes_AI$, hence annihilates its image $IC$, 
whence the claim. If, moreover, $V$ is an $\F_p$-algebra, lemma 
\ref{lem_invert}(iv) implies that $IC\to I(C\otimes_AB)$
is invertible up to a power of $\Phi$.

In the general case, consider the intermediate almost $V$-algebra 
$A_1:=\bar A\times_{\bar B}B$ equipped with the ideal 
$I_1:=0\times_{\bar B}(IB)$. In case (a), $I_1=IA_1$ 
and $A\to A_1$ is an epimorphism with nilpotent kernel, hence 
it remains such after any base change $A\to C$. To prove (i)
in case (a), it suffices then to consider the morphism $A_1\to B$, 
hence we can assume from start that $I\to IB$ is an isomorphism,
which is the case already dealt with. To prove (i) in case (b),
it suffices to consider the cases of $(A,I)\to(A_1,I_1)$
and $(A_1,I_1)\to(B,I B)$. The second case is treated above.
In the first case, we do not necessarily have $I_1=IA_1$
and the assertion to be checked is that, for every $A$-algebra 
$C$, the morphism $IC\to I_1(A_1\otimes_AC)$ is
invertible up to a power of $\Phi$. We apply lemma 
\ref{lem_invert}(v) to the commutative diagram with exact rows:
$$\xymatrix{
0 \ar[r] & I \ar[r] \ar[d] & A \ar[r] \ar[d] & 
A/I \ar[r] \ddouble & 0 \\
0 \ar[r] & IB \ar[r] & A_1 \ar[r] & A/I \ar[r] & 0 \\
}$$
to deduce that $A\to A_1$ is invertible up to some power of 
$\Phi$, hence so is $C\to A_1\otimes_AC$, which implies the
assertion.

As for (ii), we remark that  the ``only if'' part is trivial;
and we assume therefore that $(\bar A\to\bar B)^o$ is 
$\tau$-covering (resp. $\tau$-bicovering). Consider first the 
assertion for ``covering''. We need to show that $(A\to B)^o$
is of universal effective descent for $F$, where $F$ is either
one of our two fibered categories. In light of (i), this is 
reduced to the assertion that $(A\to B)^o$ is of effective 
descent for $F$. We notice that $(A\to A_1)^o$ is bicovering 
for $\tau$ (in case (a) by theorem \ref{th_liftetale} and lemma 
\ref{lem_bicover}(ii), in case (b) by proposition 
\ref{prop_invert.bicov}). As $(\bar A\to A_1/I_1)^o$ is 
an isomorphism, the assertion is reduced to the case where 
$I\to IB$ is an isomorphism. In this case, by lemma 
\ref{lem_Desc}, there is a natural equivalence: 
$\Desc(F,B/A)\stackrel{\sim}{\to}
\Desc(F,\bar B/\bar A)\times_{F_{\bar B}}F_B$.
Then the assertion follows easily from corollary 
\ref{cor_equiv.2-prod}. Finally suppose that 
$(\bar A\to\bar B)^o$ is bicovering. The foregoing 
already says that $(A\to B)^o$ is covering, so
it remains to show that $(B\otimes_AB\to B)^o$ is
also covering. The above argument again reduces to 
the case where $I\to IB$ is an isomorphism.
Then, as in the proof of lemma \ref{lem_Desc}, the induced 
morphism $I(B\otimes_AB)\to IB$ is an isomorphism
as well. Thus the assertion for ``bicovering'' is reduced to 
the assertion for ``covering''.
\end{proof}

We conclude this section with a result of a more special
nature, which can be interpreted as an easy case of almost
purity in positive characteristic.

\sset\subsubsection{}\label{subsec_easypurity}
We suppose now that the basic setup $(V,\fm)$ consists of
a {\em perfect\/} $\F_p$-algebra $V$, {\em i.e.} such that
the Frobenius endomorphism $\Phi_V:V\to V$ is bijective;
moreover we assume that there exists a non-zero-divisor 
$\eps\in\fm$ such that $\fm=\bigcup_{n>0}\eps^{1/p^m}V$.
Let us denote by $V^a\Et_\mathrm{uafp}$ (resp. 
$V[\eps^{-1}]\uEt_\mathrm{fp}$) the category of uniformly 
almost finite projective {\'e}tale $V^a$-algebras (resp. of finite
{\'e}tale $V[\eps^{-1}]$-algebra in the usual sense of \cite{SGA1}).
We will be concerned with the natural functor:
\set\begin{equation}\label{eq_easypurity}
V^a\Et_\mathrm{uafp}\to V[\eps^{-1}]\uEt_\mathrm{fp}
\quad:\quad A\mapsto A_*[\eps^{-1}].
\end{equation}

\begin{theorem}\label{th_easypurity} 
Under the assumptions of \eqref{subsec_easypurity}, the functor 
\eqref{eq_easypurity} is an equivalence of categories.
\end{theorem}
\begin{proof} Let $R$ be a finite {\'e}tale $V[\eps^{-1}]$-algebra.
Since $V[\eps^{-1}]$ is perfect, the same holds for $R$, in
view of theorem \ref{th_Frobenius}(ii) (applied in the classical 
limit case of example \ref{ex_rings}(ii)).
Let us choose a finite $V$-algebra $R_0\subset R$ such that
$R_0[\eps^{-1}]=R$ and define 
$R_1:=\bigcup_{n\in\N}\Phi^{-n}_R(R_0)$.
\begin{claim}\label{cl_independent} The $V^a$-algebra $R_1^a$ 
does not depend on the choice of $R_0$.
\end{claim}
\begin{pfclaim} Let $R_0'\subset R$ be another finite $V$-algebra
such that $R_0'[\eps^{-1}]=R$; clearly we have 
$\eps^m R_0\subset R_0'\subset\eps^{-m} R_0$ for
$m\in\N$ sufficiently large. It follows that 
$\eps^{m/p^n}\Phi^{-n}_R(R_0)\subset\Phi^{-n}_R(R'_0)
\subset\eps^{-m/p^n}\Phi^{-n}_R(R_0)$ for every $n\in\N$.
The claim readily follows.
\end{pfclaim}

\begin{claim}\label{cl_unramified}
$R_1^a$ is an unramified $V^a$-algebra.
\end{claim}
\begin{pfclaim} Let $e\in R\otimes_{V[\eps^{-1}]}R$ be the 
idempotent provided by proposition \ref{prop_idemp}; for $m\in\N$
large enough, $\eps^m\cdot e$ is contained in the subring 
$R_0\otimes_VR_0$. Hence, for every $n\in\N$, 
$\eps^{m/p^n}\cdot e\in\Phi_R^{-n}(R_0)\otimes_V\Phi_R^{-n}(R_0)$, so
$e$ defines an almost element in $(R_1\otimes_VR_1)^a$ which
fulfills the conditions (i)-(iii) of proposition \ref{prop_idemp},
and the claim follows.
\end{pfclaim}

\begin{claim}\label{cl_uniformly} 
$R_1^a$ is a uniformly almost finite $V^a$-algebra.
\end{claim}
\begin{pfclaim} For large enough $m\in\N$ we have: 
$R_0\subset\Phi_R^{-1}(R_0)\subset\eps^{-m}R_0$, therefore
$\Phi_R^{-n}(R_0)\subset\Phi_R^{-(n+1)}(R_0)\subset
\eps^{-m/p^n}\Phi_R^{-n}(R_0)$ for every $n\in\N$.
By an easy induction we deduce: 
$\Phi_R^{-(n+k)}(R_0)\subset\prod^k_{j=0}\eps^{-m/p^{n+j}}\cdot
\Phi_R^{-n}(R_0)\subset\eps^{-m/p^{n-1}}\Phi_R^{-n}(R_0)$
for every $n,k\in\N$. Finally, this implies that 
$R_1\subset\eps^{-m/p^{n-1}}\Phi_R^{-n}(R_0)$ for every $n\in\N$
and the claim follows.
\end{pfclaim}

\begin{claim}\label{cl_integr.clos} 
Let $S$ be the integral closure of $V$ in $R$; then $R_1^a=S^a$.
\end{claim}
\begin{pfclaim} Let us endow $R$ with the unique ring topology $\tau$
such that the induced subspace topology on $R_0$ is the $\eps$-adic
topology and $R_0$ is open in $R$. It is easy to check that $S$ 
consists of power-bounded elements of $R$ relative to the topology
$\tau$. Since clearly $R_1\subset S$, it suffices therefore to
show that $(R_1^a)_*\subset R$ is the subring of all power-bounded
elements of $R$. However, $(R_1^a)_*$ can be characterized as the
subring of all $x\in R$ such that $\fm\cdot x\subset R_1$; this already
implies that $(R_1^a)_*$ consists of power-bounded elements.
On the other hand, if $x\in R$ is power-bounded, it follows that
$\delta\cdot x$ is topologically nilpotent for every $\delta\in\fm$;
since $R_0$ is open in $R$, it follows that, for every $\delta\in\fm$ 
there exists $n_0\in\N$ such that $(\delta\cdot x)^n\in R_0$ for
every $n>n_0$. By taking $n:=p^k$ for sufficiently large $k\in\N$,
we deduce that $\Phi_R^k(\delta\cdot x)\in R_0$, that is 
$\delta\cdot x\in R_1$, and the claim follows.
\end{pfclaim}

\begin{claim}\label{cl_projective} 
$R_1^a$ is an almost projective $V^a$-algebra.
\end{claim}
\begin{pfclaim}
As a special case of claim \ref{cl_integr.clos}, let $W$
be the integral closure of $V$ in $V[\eps^{-1}]$; then: 
\set\begin{equation}\label{eq_almost.the.same}
W^a=V^a.
\end{equation}
Next, let $\Tr_{R/V[\eps^{-1}]}:R\to V[\eps^{-1}]$ be the
trace map of the finite {\'e}tale extension $V[\eps^{-1}]\to R$;
recall that $\Tr_{R/V[\eps^{-1}]}$ sends elements integral
over $V$ to elements integral over $V$ (to see this, we can
assume that $R$ has constant rank $n$ over $V[\eps^{-1}]$;
then the assertion can be checked after a faithfully flat base 
change $V[\eps^{-1}]\to S$, so we can further suppose that 
$R\simeq V[\eps^{-1}]^n$, in which case everything is clear);
it then follows from claim \ref{cl_integr.clos} and
\eqref{eq_almost.the.same} that $\Tr_{R/V[\eps^{-1}]}^a$
restricts to a morphism $T:R_1^a\to V^a$. Furthermore, let
$e\in R\otimes_{V[\eps^{-1}]}R$ be the idempotent defining
the diagonal imbedding; by claim \ref{cl_unramified}, for every
$\delta\in\fm$ we can write $\delta\cdot e=\sum^n_ix_i\otimes y_i$ 
for certain $x_i,y_i\in R_1$. By remark \ref{rem_zeta.identity}
(whose proof does not use theorem \ref{th_easypurity}) we deduce 
the identity: $\delta\cdot b=\sum^n_i x_i\cdot T(b\cdot y_i)$ for
every $b\in(R_1^a)_*$. This allows us to define morphisms
$\alpha:R_1^a\to(V^a)^n$, $\beta:(V^a)^n\to R_1^a$ with
$\beta\circ\alpha=\delta\cdot\one_{R^a_1}$, namely 
$\alpha(b)=(T(b\cdot y_1),...,T(b\cdot y_n))$ and
$\beta(v_1,...,v_n)=\sum^n_ix_i\cdot v_i$ for every 
$b\in R_1^a$ and $v_1,...,v_n\in V^a_*$. By lemma \ref{lem_thetwo},
the claim follows.
\end{pfclaim}

\begin{claim}\label{cl_faithful} 
The functor \eqref{eq_easypurity} is fully faithful.
\end{claim}
\begin{pfclaim} First of all, it is clear that, for every flat
$V^a$-algebras $A$, $B$, the natural map 
\set\begin{equation}\label{eq_now.and.ag}
\Hom_{V^a\Alg}(A,B)\to
\Hom_{V[\eps^{-1}]\Alg}(A_*[\eps^{-1}],B_*[\eps^{-1}])
\end{equation} 
is injective, since $A_*\subset A_*[\eps^{-1}]$ and similarly for $B$.
Suppose now that $A$ and $B$ are {\'e}tale and almost finite over $V^a$; 
then $\Phi_A$ and $\Phi_B$ are automorphisms, due to theorem 
\ref{th_Frobenius}(ii) and the assumption that $V$ is perfect.
Let $\psi:A_*[\eps^{-1}]\to B_*[\eps^{-1}]$ be any map of
$V[\eps^{-1}]$-algebras; since $A$ is almost finite, we have
$\psi(A_*)\subset\eps^{-m}B_*$ for $m\in\N$ large enough.
Since Frobenius commutes with every ring homomorphism, we deduce
$\psi(A_*)=\psi(\Phi_{A_*}^{-n}(A_*))\subset
\eps^{-m/p^n}\Phi_{B_*}^{-n}(B_*)=\eps^{-m/p^n}B_*$
for every $n\in\N$, so $\psi$ induces a morphism $\psi^a:A\to B$,
which shows that \eqref{eq_now.and.ag} is surjective.
\end{pfclaim}

It now follows from claims \ref{cl_independent}, \ref{cl_unramified},
\ref{cl_uniformly}, \ref{cl_projective} that the assignment 
$R\to R_1^a$ defines a quasi-inverse to \ref{eq_easypurity};
together with claim \ref{cl_faithful}, this concludes the proof 
of the theorem.
\end{proof}


\newpage

\section{Fine study of almost projective modules}\label{ch_traces}

An alternative title for this chapter could have been
``Everything you can do with traces". Right at the outset
we find the definition of the trace map of an almost projective 
almost finitely generated $A$-module. The whole purpose of the 
chapter is to showcase the versatility of this construction, 
a real swiss-knife of almost linear algebra. For instance, we apply
it to characterize {\'e}tale morphisms (theorem \ref{th_proj.etale});
more generally, it is used to define the {\em different ideal\/}
of an almost finite $A$-algebra. In section \ref{sec_alm.fin.rk} 
it is employed in an essential way to study the important class 
of $A$-modules of {\em finite rank}, {\em i.e.} those almost 
projective $A$-modules $P$ such that $\Lambda^i_AP=0$ for 
sufficiently large $i\in\N$. A rather complete and satisfactory
description is achieved for such $A$-modules (proposition 
\ref{prop_decomp.fin.rank}). This is further generalized in
theorem \ref{th_structure.alm.fin.rk}, to arbitrary $A$-modules 
so called of {\em almost finite rank} (see definition 
\ref{def_formal_rank}(ii)). The interest of the latter class is 
that it contains basically all the almost projective modules 
found in nature; indeed, we cannot produce a single example of
an almost projective module that is almost finitely generated 
but has not almost finite rank (but we suspect that they do exist). 
In any case, almost finitely generated modules whose rank is not 
almost finite are certainly rather weird beasts : some clue about 
their looks can be gained by analyzing the structure of 
{\em invertible\/} modules : we do this at the end of section 
\ref{sec_local.flat.site}.

The other main construction of chapter \ref{ch_traces} is
the {\em splitting algebra\/} of an almost projective module,
introduced in section \ref{sec_local.flat.site} : with its
aid we show that $A$-modules of finite rank are locally
free in the flat topology of $A$. It should be clear that 
this is a very pleasant and useful culmination for our study
of almost projective modules; we put it to use right away in 
the following section \ref{sec_groupoids}, where we show that 
an {\'e}tale groupoid of finite rank over the category of affine
almost schemes (more prosaically, the opposite of the category
of almost algebras) is universally effective, that is, it admits a good
quotient, as in the classical algebro-geometric setting.

\subsection{Almost traces}
Let $A$ be a $V^a$-algebra. 

\begin{definition}\label{def_almost_trace}
\index{Almost module(s)!almost projective!$\tr_{P/A}$ : trace morphism of an
|indref{def_almost_trace}}
\index{$\zeta_P$|indref{def_almost_trace}}
Let $P$ be an almost finitely generated  projective $A$-module. 
Then $\omega_{P/A}$ is an isomorphism by lemma
\ref{lem_three.mods}(b). The 
{\em trace morphism} of $P$ is the $A$-linear morphism
$$\tr_{P/A}:=\ev_{P/A}\circ\omega_{P/A}^{-1}:
\End_A(P)^a\to A.$$
We let $\zeta_P$ be the unique almost element of 
$P\otimes_AP^*$ such that $\omega_{P/A}(\zeta_P)=\one_P$.
\end{definition}

\begin{lemma}\label{lem_swap.maps} Let $M$, $N$ be 
almost finitely generated projective $A$-modules, 
and $\phi:M\to N$, $\psi:N\to M$ two $A$-linear morphisms. 
Then :
\begin{enumerate}
\item
$\tr_{M/A}(\psi\circ\phi)=\tr_{N/A}(\phi\circ\psi)$.
\item
If $\psi\circ\phi=a\cdot\one_M$ and 
$\phi\circ\psi=a\cdot\one_N$ for some $a\in A_*$, and
if, furthermore, there exist $u\in\End_A(M)$, $v\in\End_A(N)$
such that $v\circ\phi=\phi\circ u$, then 
$a\cdot(\tr_{M/A}(u)-\tr_{N/A}(v))=0$.
\end{enumerate}
\end{lemma}
\begin{proof} (i) : by lemma \ref{lem_three.mods}(i),
the natural morphism $N\otimes_A\Alhom_A(M,A)\to\Alhom_A(M,N)$
is an isomorphism (and similarly when we exchange the roles
of $M$ and $N$). By $A$-linearity, we can therefore assume
that $\phi$ (resp. $\psi$) is of the form 
$x\mapsto n\cdot\alpha(x)$ for some $n\in N_*$, 
$\alpha:M\to A$ (resp. of the form 
$x\mapsto m\cdot\beta(x)$ for some $m\in M_*$, 
$\beta:N\to A$). Then a simple computation yields:
$$\phi\circ\psi=\omega_{N/A}(n\cdot\alpha(m)\otimes\beta)
\qquad
\psi\circ\psi=\omega_{M/A}(m\cdot\beta(n)\otimes\alpha)$$
and the claim follows directly from the definition of the 
trace morphism.
For (ii) we compute using (i) : 
$a\cdot\tr_{M/A}(u)=\tr_{M/A}(\psi\circ\phi\circ u)=
\tr_{M/A}(\psi\circ v\circ\phi)=\tr_{N/A}(v\circ\phi\circ\psi)=
a\cdot\tr_{N/A}(v)$.
\end{proof}

\begin{lemma}\label{lem_trace.tensors} 
Let $M$, $N$ be two almost finitely generated projective 
$A$-modules, $\phi\in\End_A(M)$ and $\psi\in\End_A(N)$. Then 
$\tr_{M\otimes_AN/A}(\phi\otimes\psi)=
\tr_{M/A}(\phi)\cdot\tr_{N/A}(\psi)$.
\end{lemma}
\begin{proof} As usual we can suppose that 
$\phi=\omega_{M/A}(m\otimes\alpha)$, 
$\psi=\omega_{N/A}(n\otimes\beta)$ for some $\alpha\in M^*$,
$\beta\in N^*$. Then $\phi\otimes\psi=\omega_{M\otimes_AN/A}
((m\otimes n)\otimes(\alpha\otimes\beta))$ and the sought
identity follows by explicit calculation.
\end{proof}

\begin{proposition}\label{prop_short.ex.seq.trace} Let 
$\underline M=
(0\to M_1\stackrel{i}{\to}M_2\stackrel{p}{\to}M_3\to 0)$ 
be an exact sequence of almost finitely generated 
projective $A$-modules, and let 
$\underline u=(u_1,u_2,u_3):\underline M\to\underline M$ be 
an endomorphism of $\underline M$, given by endomorphisms 
$u_i:M_i\to M_i$ ($i=1,2,3$). Then we have 
$\tr_{M_2/A}(u_2)=\tr_{M_1/A}(u_1)+\tr_{M_3/A}(u_3)$.
\end{proposition}
\begin{proof} Suppose first that there exists a splitting
$s:M_3\to M_2$ for $p$, so that we can view $u_2$ as a matrix 
$\left(\begin{array}{ll} u_1 & v \\ 0 & u_3\end{array}\right)$,
where $v\in\Hom_A(M_3,M_1)$. By additivity of the trace, we are 
then reduced to show that $\tr_{M_2/A}(i\circ v\circ p)=0$.
By lemma \ref{lem_swap.maps}(i), this is the same as 
$\tr_{M_3/A}(p\circ i\circ v)$, which obviously vanishes.
In general, for any $a\in\fm$ we consider the morphism
$\mu_a=a\cdot\one_{M_3}$ and the pull back morphism 
$\underline{M}*\mu_a\to\underline M$ : 
$$\xymatrix{
0 \ar[r] & M_1 \ar[r]^i & M_2 \ar[r]^p & M_3 \ar[r] & 0 \\
0 \ar[r] & M_1 \ar[r] \udouble & P \ar[r]^-{p'} 
\ar[u]^{\phi} & M_3 \ar[r] \ar[u]_{\mu_a} & 0.
}$$
Pick a morphism $j:M_3\to M_2$ such that 
$p\circ j=a\cdot\one_{M_3}$; the pair $(j,\one_{M_3})$
determines a morphism $\sigma:M_3\to P$ such that 
$\sigma\circ p'=\one_{M_3}$, {\em i.e.\/} the sequence
$\underline{M}*\mu_a$ is split exact; this sequence also
inherits the endomorphism $\underline{u}*\mu_a=(u_1,v,u_3)$,
for a certain $v\in\End_A(P)$. The pair of morphisms
$(a\cdot\one_{M_2},p)$ determines a morphism 
$\psi:M_2\to P$, and it is easy to check that 
$\phi\circ\psi=a\cdot\one_{M_2}$ and 
$\psi\circ\phi=a\cdot\one_P$. We can therefore apply 
lemma \ref{lem_swap.maps} to deduce that 
$a\cdot(\tr_{P/A}(v)-\tr_{M/A}(u_2))=0$. By the foregoing 
we know that $\tr_{P/A}(v)=\tr_{M_1/A}(u_1)+\tr_{M_3/A}(u_3)$,
so the claim follows.
\end{proof}

\begin{lemma}\label{lem_explain.invert} 
Let $A$ be a $V^a$-algebra.
\begin{enumerate}
\item
If $P:=M\otimes_AN$ is an almost projective and faithful 
(resp. and almost finitely generated) $A$-module, then so 
are $M$ and $N$.
\item
If $M\otimes_AN\simeq A$, then the evaluation
map $\ev_M:M\otimes_AM^*\to A$ is an isomorphism.
\item
An invertible $A$-module is faithful and almost finitely 
generated  projective.
\item
An epimorphism $\phi:M\to N$ of invertible $A$-modules is an 
isomorphism.
\end{enumerate}
\end{lemma}
\begin{proof} Clearly (iii) is just a special case
of (i). We show (i): by proposition \ref{prop_eval.ideal}(iv)
we know that $\cE_{P/A}=A$; however, one checks easily
that $\cE_{P/A}\subset\cE_{N/A}$, whence
\set\begin{equation}\label{eq_show.faithful}
\cE_{N/A}=A.
\end{equation}
Therefore $N$ will be faithful, as soon as it is shown
to be almost projective, again by virtue of proposition
\ref{prop_eval.ideal}(iv). In any case, 
\eqref{eq_show.faithful} means that, for every $\eps\in\fm$, 
we can find an almost element of the form 
$\sum_{i=1}^nx_i\otimes\phi_i\in N\otimes_AN^*$, such
that $\sum_{i=1}^n\phi_i(x_i)=\eps$. We use such an
element to define morphisms $A\to N^n\to A$ whose
composition equals $\eps\cdot\one_N$. After tensoring
by $M$, we obtain morphisms $M\to P\to M$ whose
composition is $\eps\cdot\one_M$. Then, since $P$
is almost projective, it follows easily that so must
be $M$; similarly, if $P$ is almost finitely generated,
the same follows for $M$. By symmetry, the same holds
for $N$.

(ii): notice that, by (i), we know already that $M$
and $N$ are almost finitely generated  projective. 
By lemma \ref{lem_trace.tensors} we deduce that
$\tr_{M/A}(\one_M)\cdot\tr_{N/A}(\one_N)=1$, so both
factors are invertible in $A_*$. It follows that the
morphism $A\to\End_A(M)$ given by $a\mapsto a\cdot\one_M$
provides a splitting for the trace morphism (and similarly 
for $N$ in place of $M$). Thus we can write 
$\End_A(M)\simeq A\oplus X$, $\End_A(N)\simeq A\oplus Y$
for some $A$-modules $X$, $Y$. However, on one
hand we have a natural isomorphism 
$\End_A(M)\otimes_A\End_A(N)\simeq A$; on the other
hand, we have a decomposition 
$\End_A(M)\otimes_A\End_A(N)\simeq 
A\oplus X\oplus Y\oplus(X\otimes_AY)$; working out
the identifications, one sees that the induced isomorphism
$A\oplus X\oplus Y\oplus(X\otimes_AY)\simeq A$ restricts
to the identity morphism on the direct summand $A$; it
follows that $X=Y=0$, which readily implies the claim.

(iv): in view of (ii) we can replace $\phi$ by 
$\phi\otimes_A\one_{M^*}$, and thereby assume that $M=A$.
Then $N\simeq A/\Ker(\phi)$; it is clear that such a module
is faithful if and only if $\Ker(\phi)=0$. By (iii), the
claim follows.
\end{proof}

Lemma \ref{lem_explain.invert} explains why we do not
insist, in the definition of an invertible $A$-module, 
that it should be almost projective or almost finitely
generated: both conditions can be deduced.

\sset\subsubsection{}\label{subsec_trace.of.algebra}
\index{Almost algebra(s)!almost projective!$\Tr_{B/A}$ : trace morphism of an
|indref{subsec_trace.of.algebra}}
Suppose now that $B$ is an almost finite projective $A$-algebra. 
For any $b\in B_*$, denote by $\mu_b:B\to B$ the $B$-linear morphism 
$b'\mapsto b\cdot b'$. The map $b\mapsto\mu_b$ defines a $B$-linear 
monomorphism $\mu:B\to\End_A(B)^a$. The composition 
$$\Tr_{B/A}:=\tr_{B/A}\circ\mu:B\to A$$
will also be called the almost trace morphism of the 
$A$-algebra $B$.

\begin{proposition}\label{prop_traces}
Let $A$ and $B$ be as in \eqref{subsec_trace.of.algebra}.
\begin{enumerate}
\item
If $\phi:A\to B$ is an isomorphism, then $\Tr_{B/A}=\phi^{-1}$.
\item
If $C$ any other $A$-algebra, then 
$\Tr_{C\otimes_AB/C}=\one_C\otimes_A\Tr_{B/A}$.
\item
If $C$ is an almost finite projective $B$-algebra, then 
$\Tr_{C/A}=\Tr_{B/A}\circ\Tr_{C/B}$.
\end{enumerate}
\end{proposition}
\begin{proof} (i) and (ii) are left as exercises for 
the reader. We verify (iii). It comes down to checking 
that the following diagram commutes:
$$\xymatrix{
C\otimes_B\Alhom_B(C,B) \ar[r] \ar[d]^{\ev_{C/B}} &
C\otimes_B\Alhom_A(C,B) \ar[r]^-\sim &
C\otimes_A\Alhom_A(C,A) \ar[d]^{\ev_{C/A}} \\
B \ar[rr]^-{\Tr_{B/A}} & & A.
}$$
Therefore, pick $c\in C_*$ and $\phi\in\Hom_B(C,B)$.
For every $\eps\in\fm$ we can find elements 
$b_1,...,b_k\in B_*$ and $\phi_1,...,\phi_k\in\Hom_A(C,A)$
such that $\eps\cdot\phi(x)=\sum_ib_i\cdot\phi_i(x)$
for every $x\in C_*$. The $B$-linearity of $\phi$
translates into the identity:
\set\begin{equation}\label{eq_translate.B-lin}
\sum_ib_i\cdot\phi_i(b\cdot x)=
\sum_ib\cdot b_i\cdot\phi_i(x)\qquad
\text{for all $b\in B_*$, $x\in C_*$}.
\end{equation}
Then $\eps\cdot\ev_{C/B}(c\otimes\phi)=
\sum_ib_i\cdot\phi_i(c)$ and we need to show that
\set\begin{equation}\label{eq_Traces.A.B}
\Tr_{B/A}(\sum_ib_i\cdot\phi_i(c))=
\sum_i\phi_i(c\cdot b_i).
\end{equation}
For every $i\leq k$, let $\mu_i:A\to B$ be the morphism
$a\mapsto b_i\cdot a$ (for all $a\in A_*$); furthermore, 
let $\mu_c:B\to C$ be the morphism $b\mapsto c\cdot b$
(for all $b\in B_*$). In view of \eqref{eq_translate.B-lin},
the left-hand side of \eqref{eq_Traces.A.B} is
equal to $\tr_{B/A}(\sum_i\mu_i\circ\phi_i\circ\mu_c)$.
By lemma \ref{lem_swap.maps}(i), we have
$\tr_{B/A}(\mu_i\circ\phi_i\circ\mu_c)=
\tr_{A/A}(\phi_i\circ\mu_c\circ\mu_i)=\phi_i(c\cdot b_i)$
for every $i\leq k$. The claim follows.
\end{proof}
\begin{corollary} Let $A\to B$ be a faithfully flat almost 
finitely presented  and {\'e}tale morphism of almost 
$V$-algebras. Then $\Tr_{B/A}:B\to A$ is an epimorphism.
\end{corollary}
\begin{proof} Under the stated hypotheses, $B$ is an almost 
projective $A$-module (by proposition \ref{prop_converse}). 
Let $C=\Coker(\Tr_{B/A})$ and $\Tr_{B/B\otimes_AB}$ the trace 
morphism for the morphism of almost $V$-algebras $\mu_{B/A}$. 
By faithful flatness, the natural morphism 
$C\to C\otimes_AB=\Coker(\Tr_{B\otimes_AB/B})$ is a 
monomorphism, hence it suffices to show that 
$\Tr_{B\otimes_AB/B}$ is an 
epimorphism (here $B\otimes_AB$ is considered as a
$B$-algebra via the second factor). However, from proposition 
\ref{prop_traces}(i) and (iii) we see that $\Tr_{B/B\otimes_AB}$ 
is a right inverse for $\Tr_{B\otimes_AB/B}$. The claim follows.
\end{proof}

\sset\subsubsection{}\label{subsec_define.tau}
\index{Almost algebra(s)!almost projective!$t_{B/A}$,
$\tau_{B/A}$ : trace form of an|indref{subsec_define.tau}}
It is useful to introduce the $A$-linear morphism
$$t_{B/A}:=\Tr_{B/A}\circ\mu_{B/A}:B\otimes_AB\to A.$$
We can view $t_{B/A}$ as a bilinear form; it
induces an $A$-linear morphism
$$\tau_{B/A}:B\to B^*=\Alhom_A(B,A)$$
characterized by the equality 
$t_{B/A}(b_1\otimes b_2)=\tau_{B/A}(b_1)(b_2)$ for all
$b_1,b_2\in B_*$. We say that $t_{B/A}$ is
{\em a perfect pairing\/} if $\tau_{B/A}$ is an isomorphism.
\begin{lemma}\label{lem_eta.BC}
\index{$\eta_{B,C}$|indref{lem_eta.BC}}
Let $A\to B$ be an almost finite projective morphism of 
$V^a$-algebras and $C$ any $A$-algebra. Denote 
by $\eta_{B,C}:C\otimes_A\Alhom_A(B,A)\to
\Alhom_C(C\otimes_AB,C)$ the natural isomorphism
provided by lemma {\em\ref{lem_alhom}(i)}. Then :
\begin{enumerate}
\item 
$\tau_{B/A}$ is $B$-linear (for the natural $B$-module
structure of $B^*$ defined in remark {\em\ref{rem_B.struct}});
\item
$\eta_{B,C}$ is $C\otimes_AB$-linear;
\item
$\eta_{B,C}\circ(\one_C\otimes\tau_{B/A})=
\tau_{C\otimes_AB/C}$.
\end{enumerate}
\end{lemma}
\begin{proof} For any $b\in B_*$, let 
$\xi_b:B\to A$ the $A$-linear morphism defined by
the rule $b'\mapsto\Tr_{B/A}(b'\cdot b)$ for all 
$b'\in B_*$. Then, directly from the definition we
can compute: 
$(\eta_{B,C}\circ(\one_C\otimes\tau_{B/A}))(c\otimes b)
(c'\otimes b')=\eta_{B,C}(c\otimes\xi_b)(c'\otimes b')=
c\cdot c'\cdot\Tr_{B/A}(b'\cdot b)$ for all 
$b,b'\in B_*$, $c,c'\in C_*$. But by
proposition \ref{prop_traces}(ii), the latter expression 
can be rewritten as 
$\tau_{C\otimes_AB/C}(c\otimes b)(c'\otimes b')$,
which shows (iii). The proofs of (i) and (ii) are 
similar direct verifications: we show (i) and leave
(ii) to the reader. Let us pick any $b,b',b''\in B_*$; 
then $(b\cdot\tau_{B/A}(b'))(b'')=\tau_{B/A}(b')(bb'')=
\Tr_{B/A}(bb'b'')=(\tau_{B/A}(bb'))(b'')$.
\end{proof}

\begin{theorem}\label{th_proj.etale} An almost finite projective
morphism $\phi:A\to B$ of almost $V$-algebras is {\'e}tale if and only 
if the trace form $t_{B/A}$ is a perfect pairing.
\end{theorem}
\begin{proof} By lemma \ref{lem_eta.BC}, we have a 
commutative diagram:
\set\begin{equation}\label{eq_big.diagram.eta}{
{\diagram
(B\otimes_AB)\otimes_BB \ar[r]^-\sim 
\ar[d]_{\one_{B\otimes_AB}\otimes_B\tau_B} &
B\otimes_AB \ar[d]_{\one_B\otimes_A\tau_B}
\ar[dr]^{\tau_{B\otimes_AB/B}} \\
(B\otimes_AB)\otimes_BB^* 
\ar[r]^-\sim & B\otimes_AB^* \ar[r]^-{\eta_{B,B}} &
\Alhom_B(B\otimes_AB,B)
\enddiagram}}\end{equation}
in which all the morphisms are $B\otimes_AB$-linear
(here we take the $B$-module structure on $B\otimes_AB$
given by multiplication on the right factor).
Suppose now that $\phi$ is {\'e}tale; then, by corollary 
\ref{cor_unram}, there is an isomorphism of 
$B$-algebras: $B\otimes_AB\simeq I_{B/A}\oplus B$.
It follows that $\tau_{B\otimes_AB/B}=
\tau_{B/B}\oplus\tau_{I_{B/A}/B}$. Especially,
$\one_B\otimes_{B\otimes_AB}\tau_{B\otimes_AB/B}$
is the identity morphism of $B$ (by proposition 
\ref{prop_traces}(i)). This means that in the diagram
$B\otimes_{B\otimes_AB}\eqref{eq_big.diagram.eta}$
all the arrows are isomorphisms. In particular,
$\tau_{B/A}$ is an isomorphism, as claimed.

To prove the converse, we consider the almost element
$\zeta_B$ of the $B\otimes_AB$-module $B\otimes_AB^*$.
Viewing $B^*$ as a $B$-module in the natural way
(cp. remark \ref{rem_B.struct}), we also get a scalar
multiplication morphism $\sigma_{B^*/B}:B\otimes_AB^*\to B^*$
(see \eqref{subsec_comm.constr}).
\begin{claim}\label{claim_zeta.props} 
With the above notation we have: $I_{B/A}\cdot\zeta_B=0$ and 
$\sigma_{B^*/B}(\zeta_B)=\Tr_{B/A}$.
\end{claim}
\begin{pfclaim} Notice that $\omega_{B/A}$ is also
$B\otimes_AB$-linear for the $B\otimes_AB$-module
structure on $\End_A(B)$ such that 
$((b\otimes b')\cdot\phi)(b'')=b'\cdot\phi(b\cdot b")$
for every $b,b',b"\in B_*$ and every $\phi\in\End_A(B)$.
We compute $\omega_{B/A}((b\otimes b')\cdot\zeta_B)(b")=
((b\otimes b')\cdot\omega_{B/A}(\zeta_B))(b")=
b\cdot b'\cdot b"$. Whence $x\cdot\zeta_B=
\mu_{B/A}(x)\cdot\zeta_B$ for every $x\in B\otimes_AB_*$
which implies the first claimed identity. Next we compute:
$\sigma_{B^*/B}(\zeta_B)(b)=
\ev_{B}((1\otimes b)\cdot\zeta_B)=
(\tr_{B/A}\circ\omega_{B/A})((1\otimes b)\cdot\zeta_B)=
\tr_{B/A}((1\otimes b)\cdot\omega_{B/A}(\zeta_B))=
\tr_{B/A}((1\otimes b)\cdot\one_B)=\Tr_{B/A}(b)$ for every
$b\in B_*$. The claim follows.

\end{pfclaim}
Suppose now that $\tau_{B/A}$ is an isomorphism. Then
we can define $e:=(\one_B\otimes\tau_{B/A}^{-1})(\zeta_B)$.
From claim \ref{claim_zeta.props} and lemma 
\ref{lem_eta.BC}(i) we derive that $I_{B/A}\cdot e=0$
and $\tau_{B/A}(\sigma_{B/B}(e))=\Tr_{B/A}$. The latter
equality implies that $\sigma_{B/B}(e)=1$, in other
words $\mu_{B/A}(e)=1$. We see therefore that $e$
satisfies conditions (ii) and (iii) of proposition
\ref{prop_idemp} and therefore also condition (i),
since the latter is an easy consequence of the other
two. Thus $A\to B$ is an {\'e}tale morphism, as claimed.
\end{proof}

\begin{remark}\label{rem_zeta.identity} By inspection 
of the proof of theorem \ref{th_proj.etale}, we see that 
the following has been shown. Let $A\to B$ be an {\'e}tale 
morphism of $V^a$-algebras. Then 
$(\one_B\otimes\tau_{B/A})(e_{B/A})=\zeta_B$.
\end{remark}

\begin{definition}\label{def_nilradical}
\index{Almost algebra(s)!$\nil(A)$ : nilradical of an|indref{def_nilradical}}
\index{Almost algebra(s)!reduced|indref{def_nilradical}}
The {\em nilradical\/} of an almost algebra $A$ is the ideal 
$\nil(A)=\nil(A_*)^a$ (where, for a ring $R$, we denote
by $\nil(R)$ the ideal of nilpotent elements in $R$).
We say that $A$ is {\em reduced\/} if $\nil(A)\simeq 0$.
\end{definition}

\sset\subsubsection{}
Notice that, if $R$ is a $V$-algebra, then every nilpotent 
ideal in $R^a$ is of the form $I^a$, where $I$ is a nilpotent
ideal in $R$ (indeed, it is of the form $I^a$ where $I$ is
an ideal, and $\fm\cdot I$ is seen to be nilpotent). It
follows easily that $\nil(A)$ is the colimit of the nilpotent
ideals in $A$; moreover $\nil(R)^a=\nil(R^a)$. Using this
one sees that $A/\nil(A)$ is reduced.

\begin{proposition}\label{prop_reduced.etale} Let $A\to B$ 
be an {\'e}tale almost finitely presented morphism of almost 
algebras. If $A$ is reduced then $B$ is reduced as well.
\end{proposition}
\begin{proof} Under the stated hypothesis, $B$ is an almost 
projective $A$-module (by virtue of proposition 
\ref{prop_converse}(ii)). Hence, 
for given $\eps\in\fm$, pick a sequence of morphisms 
$B\stackrel{u_\eps}{\to}A^n\stackrel{v_\eps}{\to}B$ such that
$v_\eps\circ u_\eps=\eps\cdot\one_B$; let $\mu_b:B\to B$
be multiplication by $b\in B_*$ and define $\nu_b:A^n\to A^n$
by $\nu_b=v_\eps\circ\mu_b\circ u_\eps$. One verifies easily 
that $\nu^m_b=\eps^{m-1}\cdot\nu_{b^m}$ for all integers $m>0$.
Now, suppose that $b\in\nil(B_*)$. It follows that $b^m=0$ 
for $m$ sufficiently large, hence $\nu^m_b=0$ for $m$ 
sufficiently large. Let $\fp$ be any prime ideal of $A_*$; 
let $\pi:A_*\to A_*/\fp$ be the natural projection and 
$F$ the fraction field of $A_*/\fp$. 
The $F$-linear morphism $\nu_{b*}\otimes_{A_*}\one_F$ 
is nilpotent on the vector space $F^n$, hence 
$\pi\circ\tr_{A^n_*/A_*}(\nu_{b*})=
\tr_{F^n/F}(\nu_{b*}\otimes_{A_*}\one_F)=0$.
This shows that $\tr_{A^n_*/A_*}(\nu_{b*})$ lies in the 
intersection of all prime ideals of $A_*$, hence it 
is nilpotent. Since by hypothesis $A$ is reduced, we get 
$\tr_{A^n_*/A_*}(\nu_{b*})=0$, whence $\tr_{A^n/A}(\nu_b)=0$. 
Using lemma \ref{lem_swap.maps}(i) we deduce 
$\eps\cdot\tr_{B/A}(\mu_b)=0$, and finally, $\tr_{B/A}(b)=0$.  
Now, for any $b'\in B_*$, the almost element $bb'$ will 
be nilpotent as well, so the same conclusion applies to 
it. This shows that $\tau_{B/A}(b)=0$. But by hypothesis 
$B$ is {\'e}tale over $A$, hence theorem \ref{th_proj.etale} 
yields $b=0$, as required.
\end{proof}
\begin{remark}\label{rem_regular}
\index{Almost module(s)!$M_*$ : almost element(s) of an!$M$-regular|indref{rem_regular}}
Let $M$ be an $A$-module. We say that an almost 
element $a$ of $A$ is $M$-{\em regular\/} if the multiplication 
morphism $m\mapsto am~:~M\to M$ is a monomorphism. Assume 
({\bf A}) (see \eqref{subsec_conditions.A.B}) and suppose 
furthermore that $\fm$ is generated by a multiplicative system 
$\cS$ which is a cofiltered semigroup under the preorder 
structure $(\cS,\succ)$ induced by the divisibility relation 
in $V$. We say that $\cS$ is {\em archimedean\/} if, for all 
$s,t\in\cS$ there exists $n>0$ such that $s^n\succ t$. 
Suppose that $\cS$ is  archimedean and that $A$ is a reduced 
almost algebra. Then $\cS$ consists of $A$-regular elements. 
Indeed, by hypothesis $\nil(A_*)^a=0$; since the annihilator 
of $\cS$ in $A_*$ is $0$ we get $\nil(A_*)=0$. Suppose 
that $s\cdot a=0$ for some $s\in\cS$ and $a\in A_*$. Let 
$t\in\cS$ be arbitrary and pick $n>0$ such that $t^n\succ s$. 
Then $(ta)^n=0$ hence $ta=0$ for all $t\in\cS$, hence $a=0$.
\end{remark}

\begin{definition}\label{def_different}
\index{Almost algebra(s)!almost projective!$\cD_{B/A}$ : different ideal of an
|indref{def_different}}
Let $\phi:A\to B$ be an almost finite projective 
morphism of $V^a$-algebras. By \eqref{subsec_define.tau}, we 
can assign to $\phi$ a $B$-linear trace morphism $\tau_{B/A}:B\to B^*$.
The {\em different ideal\/} of the morphism $\phi$ is the ideal 
$\cD_{B/A}:=\Ann_B(\Coker\,\tau_{B/A})\subset B$.
\end{definition}

\begin{lemma}\label{lem_Anns.multiply} 
Let $M_1\stackrel{\phi}{\to}M_2\stackrel{\psi}{\to}M_3$
be two $A$-linear morphisms of invertible $A$-modules
$M_i$ ($i\leq 3$) and $C$ an $A$-algebra. Then:
\begin{enumerate}
\item 
$\Ann_A(\Coker(\psi\circ\phi))=
\Ann_A(\Coker\,\phi)\cdot\Ann_A(\Coker\,\psi)$.
\item
$\Ann_C(\Coker(\one_C\otimes_A\psi))=C\cdot\Ann_A(\Coker\,\psi)$.
\end{enumerate}
\end{lemma}
\begin{proof} (i): Since, by lemma \ref{lem_explain.invert}(iii), 
$M_3$ is faithfully flat, lemma \ref{lem_annihilate.and.forget}
yields
\set\begin{equation}\label{eq_base.change.inv}
\Ann_A(\Coker\,\psi)=\Ann_A(\Coker(\psi\otimes_A\one_{M^*_3}))
\end{equation}
and likewise for $\phi$; hence we can replace $M_i$
by $M_i\otimes_AM_3^*$ and suppose that $M_3=A$.
Moreover, since $M_2$ is invertible, $\ev_{M_2/A}$ is
an isomorphism, by lemma \ref{lem_explain.invert}(ii). Let 
$\tilde\ev_{M_2/A}:M_2^*\otimes_AM_2\to A$ be the map 
given by the rule: 
$\phi\otimes x\mapsto\ev_{M_2/A}(x\otimes\phi)$, for
every $\phi\in(M^*_2)_*$ and $x\in M_2$. Set 
$\lambda:=\ev_{M_2/A}\circ(\phi\otimes_A\one_{M_2^*}):
M_1\otimes_AM_2^*\to A$; then 
$\phi\circ(\one_{M_1}\otimes_A\tilde\ev_{M_2/A})=
\lambda\otimes_A\one_{M_2}$, so that 
$\Ann_A(\phi)=\Ann_A(\lambda\otimes_A\one_{M_2})=
\Ann_A(\lambda)$.
Thus, we can replace $\phi$ by $\lambda\otimes_A\one_{M_2}$
and then we have to show that 
$$
\Ann_A(\Coker\,\psi\circ(\lambda\otimes_A\one_{M_2}))=
\Ann_A(\Coker\,\psi)\cdot\Ann_A(\Coker\,\lambda).
$$ 
However, quite generally we have: 
\set\begin{equation}\label{eq_quite.generally}
\Ann_A(\Coker(M\to A))=\Img(M\to A)
\end{equation}
for any $A$-linear morphism $M\to A$. Hence we compute:
$\Ann_A(\Coker\,\psi\circ(\lambda\otimes_A\one_{M_2}))=
\Img(\psi\circ(\lambda\otimes_A\one_{M_2}))=
\psi(\Img(\lambda\otimes_A\one_{M_2}))=
\psi(\Img(\lambda)\cdot M_2)=\Img(\lambda)\cdot\Img(\psi)=
\Ann_A(\Coker\,\lambda)\cdot\Ann_A(\Coker\,\psi)$.

(ii): again, using \eqref{eq_base.change.inv} we reduce
to the case where $M_3=A$; then the claim follows easily
from \eqref{eq_quite.generally}.
\end{proof}

\begin{lemma}\label{lem_complicated} Let $\phi:A\to B$ 
be a morphism of\/ $V^a$-algebras as in definition 
{\em\ref{def_different}}. Let $C$ be an $A$-algebra. Suppose
that either $C$ is flat over $A$, or $B^*$ is an invertible
$B$-module for its natural $B$-module structure. Then  
$\cD_{C\otimes_AB/C}=\cD_{B/A}\cdot(C\otimes_AB)$.
\end{lemma}
\begin{proof} Under the stated assumptions, 
$\Alhom_A(B,A)$ is an almost finitely generated
projective $A$-module. In particular, 
$\Coker\,\tau_{B/A}$ is almost finitely generated; 
If $C$ is flat over $A$, it follows by lemma 
\ref{lem_annihilate.and.forget} that
$\Ann_{C\otimes_AB}(C\otimes_A\Coker\,\tau_{B/A})=
\cD_{B/A}\cdot(C\otimes_AB)$; if $B^*$ is an invertible
$B$-module, the same holds by virtue of lemma 
\ref{lem_Anns.multiply}(ii). 
However, by lemma 
\ref{lem_eta.BC}(iii), the trace pairing is preserved 
under arbitrary base changes, so: 
$C\otimes_A\Coker\,\tau_{B/A}\simeq
\Coker(\one_C\otimes_A\tau_{B/A})\simeq
\Coker\,\tau_{C\otimes_A/B}$, which shows the claim.
\end{proof}

\begin{proposition}\label{prop_diff.in.towers} 
Let $B\to C$ be a morphism of $A$-algebras, and suppose 
that $B$ (resp. $C$) is an almost finite projective $A$-algebra 
(resp. $B$-algebra). 
Suppose moreover that $B^*:=\Alhom_A(B,A)$ (resp. 
$C^*:=\Alhom_B(C,B)$) is an invertible $B$-module (resp.
$C$-module) for its natural $B$-module (resp. $C$-module) 
structure. Then
$$\cD_{C/A}=\cD_{C/B}\cdot\cD_{B/A}.$$
\end{proposition}
\begin{proof} Let $C^*_{/A}:=\Alhom_A(C,A)$ and
define a $C$-linear morphism $\xi:\Alhom_B(C,B^*)\to C^*_{/A}$
by the rule: $\phi\mapsto(c\mapsto\phi(c)(1))$ for every
$\phi\in\Hom_B(C,B^*)$ and $c\in C_*$.
\begin{claim}\label{cl_theta.iso} $C^*_{/A}$ is an invertible 
$C$-module and $\xi$ is an isomorphism.
\end{claim}
\begin{pfclaim} By lemma \ref{lem_three.mods}(i), the natural 
morphism $\lambda:C^*\otimes_BB^*\to\Alhom_B(C,B^*)$ is a
$C$-linear isomorphism. It suffices therefore to show that 
$\xi\circ\lambda^{-1}:C^*\otimes_BB^*\to C^*_{/A}$ is an 
isomorphism. One verifies easily that $\xi\circ\lambda^{-1}$ 
is defined by the rule: $\phi\otimes\psi\mapsto\psi\circ\phi$, 
and then the claim follows from lemma \ref{lem_alhom}(iii).
\end{pfclaim}

Unwinding the definitions, one verifies that the following 
diagram commutes:
$$
\xymatrix{
C \ar[rrr]^-{\tau_{C/A}} \ar[d]_{\tau_{C/B}} & & & C^*_{/A} \\
C^* \ar[rrr]^-{\Alhom_B(C,\tau_{B/A})} & & &
\Alhom_B(C,B^*) \ar[u]_\xi.
}
$$
Thus, taking into account claim \ref{cl_theta.iso}, and
lemma \ref{lem_Anns.multiply}(i), we have 
$\Ann_C(\Coker\,\tau_{C/A})=\Ann_C(\Coker\,\tau_{C/B})\cdot
\Ann_C(\Coker(\Alhom_B(C,\tau_{B/A})))$. However,
$\Coker(\Alhom_B(C,\tau_{B/A}))\simeq 
C^*\otimes_B\Coker\,\tau_{B/A}$ by lemma \ref{lem_three.mods}(b).
By lemma \ref{lem_explain.invert}(iii), $C^*$ is faithfully 
flat; consequently:
$$
\Ann_C(\Coker(\Alhom_B(C,\tau_{B/A})))=
\Ann_C(\Coker\,\tau_{B/A})
$$ 
which implies the assertion.
\end{proof}

\begin{lemma}\label{lem_etale.different}
Let $\phi:A\to B$ be a morphism of $V^a$-algebras
as in definition {\em\ref{def_different}}. Suppose moreover that 
$B^*$ is an invertible $B$-module for its natural $B$-module
structure. Then $\phi$ is {\'e}tale if and only if $\cD_{B/A}=B$.
\end{lemma}
\begin{proof} By theorem \ref{th_proj.etale} it follows easily 
that $\cD_{B/A}=B$ whenever $\phi$ is {\'e}tale. Conversely, 
suppose that $\cD_{B/A}=B$; it then follows that $\tau_{B/A}$
is an epimorphism. Again by theorem \ref{th_proj.etale}, we need
only show that $\tau_{B/A}$ is an isomorphism. This follows from
lemma \ref{lem_explain.invert}(iv).
\end{proof}

The following lemma will be useful when we will compute
the different ideal in situations such as those contemplated
in proposition \ref{prop_firstblood}.

\begin{lemma}\label{lem_net.of.algs}
Let $A$ be a $V^a$-algebra, $B$ an almost finite almost projective 
$A$-algebra, and let $\{B_\alpha~|~\alpha\in J\}$ be a net of 
$A$-subalgebras of $B$, with $B_\alpha$ almost finite projective 
over $A$ for every $\alpha\in J$, such that 
$\liminv{\alpha\in J}\,B_\alpha=B$ in $\cI_A(B)$. 
Then $\liminv{\alpha\in J}\,\cD_{B_\alpha/A}=\cD_{B/A}$.
\end{lemma}
\begin{proof} For given $\alpha\in J$, let $\eps\in V$ such that
$\eps B\subset B_\alpha$; lemma \ref{lem_swap.maps}(ii) 
implies that $\eps\cdot\Tr_{B_\alpha/A}(b)=
\eps\cdot\Tr_{B/A}(b)$ for every $b\in B_{\alpha*}$.
Hence the diagrams:
$$
\xymatrix{
B \ar[r]^{\mu_\eps} \ar[d]_{\eps\cdot\tau_{B/A}} &
B_\alpha \ar[d]^{\tau_{B_\alpha/A}} & & 
B_\alpha \ar[d]_{\tau_{B_\alpha/A}} \ar[r] & 
B \ar[d]^{\eps\cdot\tau_{B/A}} \\
B^* \ar[r] & B_\alpha^* & & 
B_\alpha^* \ar[r]^{\mu_\eps^*} & B^*
}
$$
commute. The rightmost diagram implies that 
$\cD_{B_\alpha/A}\cdot\Img\,\mu_\eps^*\subset
\Img(\eps\cdot\tau_{B/A})\subset\Img\,\tau_{B/A}$.
Hence $\eps\cdot\cD_{B_\alpha/A}\subset\cD_{B/A}$, so
finally $\Ann_V(B/B_\alpha)\cdot\cD_{B_\alpha/A}\subset\cD_{B/A}$.
From the leftmost diagram we deduce that $\eps\cdot\cD_{B/A}$
(which is an ideal in $B_\alpha$) annihilates 
$\Coker(\eps\cdot\tau_{B/A}:B\to B^*)$ and on 
the other hand $\Ann_V(B/B_\alpha)$ obviously annihilates 
$\Coker(B^*\to B_\alpha^*)$; we deduce that 
$\Ann_V(B/B_\alpha)^2\cdot\cD_{B/A}\subset\cD_{B_\alpha/A}$,
whence the claim.
\end{proof}

\subsection{Endomorphisms of $\hat\G_m$.}
This section is dedicated to a discussion of the universal 
ring that classifies endomorphisms of the formal group
$\hat\G_m$. The results of this section will be used in 
sections \ref{sec_alm.fin.rk} and \ref{sec_local.flat.site}.

\sset\subsubsection{}\label{subsec_truncated_G_m}
\index{${\bf Bud}(n,d,R)$ : $n$-buds of a formal 
group!$\G_{m,R}(n)$ : $n$-bud of $\hat\G_m$|indref{subsec_truncated_G_m}}
For every ring $R$ and every integer $n\geq 0$ we introduce 
the "$n$-truncated" version of $\hat\G_{m,R}$. This is 
the scheme $\G_{m,R}(n):=\Spec\,R[T]/(T^{n+1})$, endowed 
with the multiplication morphism which is associated
to the co-multiplication map 
$$R[T]/(T^{n+1})\to R[T,S]/(T,S)^{n+1}\qquad
T\mapsto T+S+T\cdot S.$$
Then in the category of formal schemes we have a natural 
identification
$\hat\G_{m,R}\simeq\colim{n\in\N}\G_{m,R}(n)$.

\sset\subsubsection{}\label{subsec_buds}
\index{${\bf Bud}(n,d,R)$ : $n$-buds of a formal group|indref{subsec_buds}}
In the terminology of \cite[\S II.4]{MLazard}, $\G_m(n)$
is the $n$-bud of $\hat\G_m$. We will be mainly interested
in the endomorphisms of $\G_m(n)$, but before we can get
to that, we will need some complements on buds over
artinian ring. Therefore, suppose we have a cartesian
diagram of artinian rings
\set\begin{equation}\label{eq_again.diag.algs}
{\diagram
R_3 \ar[r] \ar[d] & R_1 \ar[d] \\
R_2 \ar[r] & R_0
\enddiagram}
\end{equation}
such that one of the two maps $R_i\to R_0$ ($i=1,2$)
is surjective. For any ring $R$, we define the category 
${\bf Bud}(n,d,R)$ of $n$-buds over $R$ whose underlying 
$R$-algebra is isomorphic to 
$R[T_1,...,T_d]/(T_1,...,T_d)^{n+1}$.

\begin{lemma}\label{lem_augmented.alg} 
Let $S$ be a finite flat augmented algebra 
over a local noetherian ring $R$; let $I$ be the augmentation 
ideal, and suppose that $I^{n+1}=0$. Let $\kappa$ be the 
residue field of $R$, and suppose that 
$S\otimes_R\kappa\simeq
\kappa[t_1,...,t_d]/(t_1,...,t_d)^{n+1}$.
Then $S\simeq R[T_1,...,T_d]/(T_1,...,T_d)^{n+1}$.
\end{lemma}
\begin{proof} Let $\eps:S\to R$ be the augmentation
map. For every $i=1,...,d$, pick a lifting $T'_i\in S$
of $t_i$; set $T_i=T'_i-\eps(T'_i)$. By Nakayama's lemma, 
the monomials $T_1^{a_1}\cdot...\cdot T_d^{a_d}$ with 
$\sum_{i=1}^da_i\leq n$ generate the $R$-module $S$.
Furthermore, under the stated hypothesis, $S$ is
a free $R$-module, and its rank is equal to 
$\dim_\kappa S\otimes_R\kappa$; hence the above
monomials form an $R$-basis of $S$. Clearly the elements
$T_i$ lie in the augmentation ideal of $S$, therefore
every product of $n+1$ of them equals zero; in other
words, the natural morphism $R[X_1,...,X_d]\to S$ given
by $X_i\mapsto T_i$ is surjective, with kernel containing
$J:=(X_1,...,X_d)^{n+1}$; but by comparing the ranks over 
$R$ we see that this kernel cannot be larger than $J$. 
The assertion follows.
\end{proof}

\begin{proposition}\label{prop_two.fib.prod.buds}
In the situation of \eqref{eq_again.diag.algs}, the 
natural functor
$${\bf Bud}(n,d,R_3)\to{\bf Bud}(n,d,R_1)
\times_{{\bf Bud}(n,d,R_0)}{\bf Bud}(n,d,R_2)$$
is an equivalence of categories.
\end{proposition}
\begin{proof} Let $\cM_{i,\mathrm{proj}}$ ($i=0,...,3$)
be the category of projective $R_i$-modules. By our 
previous discussion on descent, we already know that 
\eqref{eq_again.diag.algs} induces a natural equivalence 
between $\cM_{3,\mathrm{proj}}$ and the $2$-fibered product 
$\cM_{1,\mathrm{proj}}\times_{\cM_{0,\mathrm{proj}}}
\cM_{2,\mathrm{proj}}$. It is easy to see that this 
equivalence respects the rank of $R_i$-modules, hence 
induces a similar equivalence for the categories 
$\cM_{i,\mathrm{f.f.}}$ of
free $R_i$-modules of finite rank. Given two objects 
$M:=(M_1,M_2,\alpha:M_1\otimes_{R_1}R_0
\stackrel{\sim}{\to}M_2\otimes_{R_2}R_0)$ and
$N:=(N_1,N_2,\beta:N_1\otimes_{R_1}R_0
\stackrel{\sim}{\to}N_2\otimes_{R_2}R_0)$, define
the tensor product $M\otimes N:=(M_1\otimes_{R_1}N_1,
M_2\otimes_{R_2}N_2,\alpha\otimes_{R_0}\beta)$.
Then one checks easily that the above equivalences
respect tensor products. It follows formally that
one has analogous equivalences for the categories
of finite flat $R_i$-algebras. From there, one further
obtains equivalences on the categories of such 
$R_i$-algebras that are augmented over $R_i$, and
even on the subcategories $R_i\Alg_{\mathrm{aug.fl.}}^{(n)}$
of those augmented $R_i$-algebras such that the 
$(n+1)$-th power of the augmentation ideal
vanishes. These categories admit finite coproducts,
that are constructed as follows. For augmented
$R_i$-algebras $\eps_A:A\to R_i$ and $\eps_B:B\to R_i$,
set $(A\to R_i)\otimes(B\to R_i):=
(A\otimes_{R_i}B/\Ker(\eps_A\otimes_{R_i}\eps_B)^{n+1}
\to R_i)$; 
this is a coproduct of $A$ and $B$. By formal reasons,
the foregoing equivalences of categories respect these
coproducts. Finally, an object of ${\bf Bud}(n,d,R_i)$ 
can be defined as a commutative group object in 
$(R_i\Alg_{\mathrm{aug.fl.}}^{(n)})^o$, such that its
underlying $R_i$-algebra is isomorphic to 
$R_i[T_1,...,T_d]/(T_1,...,T_d)^{n+1}$.
By formal categorical considerations we see that 
the foregoing equivalence induces equivalences on
the commutative group objects in the respective categories.
It remains to check that an $R_3$-algebra $S$ such
that $S\otimes_{R_3}R_i\simeq 
R_i[T_1,...,T_d]/(T_1,...,T_d)^{n+1}$, (for $i=1,2$)
is itself of the form 
$R_i[T_1,...,T_d]/(T_1,...,T_d)^{n+1}$. However, this 
follows readily from lemma \ref{lem_augmented.alg} and 
the fact that one of the maps $R_3\to R_i$ ($i=1,2$) is 
surjective.
\end{proof}

\sset\subsubsection{}\label{subsec_universal.ring}
\index{$\cG_n$, $\cG_\infty$|indref{subsec_universal.ring}{}, 
\indref{def_formal_rank}}
For a given ring $R$, the endomorphisms of $\G_{m,R}(n)$ 
are all the polynomials
$f(T):=a_0+a_1\cdot T+...+a_n\cdot T^n$ such that
$f(T)+f(S)+f(T)\cdot f(S)\equiv
f(T+S+T\cdot S)$ (mod $(T,S)^{n+1}$).
This relationship translates into a finite set of 
polynomial identities for the coefficients $a_0,...,a_n$, 
and using these identities we can therefore define a 
quotient $\cG_n$ of the ring in $n$ indeterminates
$\Z[X_1,...,X_n]$ which will be the ``universal ring
of endomorphisms'' of $\G_m(R)$, {\em i.e.}, such that
$X_1\cdot T+X_2\cdot T^2+...+X_n\cdot T^n$ is
an endomorphism of $\G_{m,\cG_n}(n)$ and such that, for
every ring $R$, and every $f(T)$ as above, the
map $\Z[X_1,...,X_n]\to R$ given by $X_i\mapsto a_i$
($i=1,...,n$) factors through a (necessarily unique) map 
$\cG_n\to R$. One of the main results of this section will
be a simple and explicit description of the ring $\cG_n$.

\begin{proposition}\label{prop_cG_n.smooth} 
$\cG_n$ is a smooth $\Z$-algebra.
\end{proposition}
\begin{proof} We know already that $\cG_n$ is of finite 
type over $\Z$, therefore it suffices to show that, for 
every prime ideal $\fp$ of $\cG_n$, the local ring 
$\cG_{n,\fp}$ is formally smooth for the $\fp$-adic 
topology (see \cite[Ch.IV, Prop.17.5.3]{EGA4}). Therefore, 
let $R_1\to R_0$ be a surjective homomorphism of local 
artinian rings; we need to show that the natural map 
$\End(\G_{m,R_1}(n))\to\End(\G_{m,R_0}(n))$
is surjective. Let $f\in\End(\G_{m,R_0}(n))$; we define
an automorphism $\chi$ of 
$\G_{m,R_0}\times_{R_0}\G_{m,R_0}:=R_0[T,S]/(T,S)^{n+1}$, 
the $n$-bud of $\hat\G_m\times\hat\G_m$, by setting
$(T,S)\mapsto(T,f(T)+S+f(T)\cdot S)$. Then, thanks 
to proposition \ref{prop_two.fib.prod.buds}, we obtain 
an $n$-bud $X_n$ over $R_2:=R_1\times_{R_0}R_1$, by gluing 
two copies of $\G_{m,R_1}\times\G_{m,R_1}$ along the 
automorphism $\chi$.
\begin{claim}\label{cl_best.left.reader} 
The $n$-bud $X_n$ is isomorphic to 
$\G_{m,R_2}\times\G_{m,R_2}$ if and only if $\chi$
lifts to an automorphism of $\G_{m,R_1}\times\G_{m,R_1}$.
\end{claim}
\begin{pfclaim} Taking into account the decription
of $B(n,d,R_2)$ as $2$-fibered product of categories,
the proof amounts to a simple formal verification,
which is best left to the reader.
\end{pfclaim}

\begin{claim} There exists a compatible system of 
$k$-buds $X_k$ over $R_2$ for every $k>n$, such that $X_k$ 
reduces to $X_{k-1}$ over $R_2$, and specializes to 
$\G^2_{m,R_1}(k)$ over the quotient $R_1$ of $R_2$.
\end{claim}
\begin{pfclaim}
In case $R_2$ is a torsion-free $\Z$-algebra, this follows
from \cite[Ch.II, \S 4.10]{MLazard} and an easy induction.
If $R_2$ is a general artinian ring, choose a torsion-free
$\Z$-algebra $R_3$ with a surjective homomorphism $R_3\to R_2$.
By {\em loc.~cit.} (and an easy induction) we can find an 
$n$-bud $Y_n$ over $R_3$ such that $Y_n$ specialises to $X_n$
on the quotient $R_2$, and $Y_n$ reduces to $\G^2_{m,R_3}(1)$
over $R_3$. Then, again by {\em loc.~cit}, we can find a
compatible system of $k$-buds $Y_k$ on $R_3$ for every $k>n$, 
such that $Y_k$ reduces to $Y_{k-1}$ over $R_3$ and specializes
to $\G^2_{m,R_1}(k)$ over the quotient $R_1$ of $R_3$.
The claim holds if we take $X_k$ equal to the specialization
of $Y_k$ over $R_2$.
\end{pfclaim}

The direct limit (in the category of
formal schemes) of the system $(X_k)_{k\geq n}$
is a formal group $\hat X$ over $R_2$, such 
that $\hat X\otimes_{R_2}R_1\simeq
\hat\G_{m,R_1}\times\hat\G_{m,R_1}$. This formal group
gives rise to a $p$-divisible group $(\hat X(n))_{n\geq 0}$,
where $\hat X(n)$ is the kernel of multiplication by $p^n$
in $\hat X$. For every $m\in\N$, $\hat X(m)$ is a finite
flat group scheme over $R_2$, such that 
$\hat X(m)\times_{R_2}R_1\simeq
\mu_{p^m,R_1}\times\mu_{p^m,R_1}$.
Denote by $\hat X(m)^*$ the Cartier dual of $\hat X(m)$ 
(cp. \cite[\S III.14]{Mumford}). 
Then $\hat X(m)^*\times_{R_2}R_1\simeq(\Z/p^m\Z)_{R_1}^2$, 
in particular it has $p^{2m}$ connected components.
Since the pair $(R_2,R_1)$ is henselian, it follows that
$\hat X(m)^*$ must have $p^{2m}$ connected components as
well, and consequently
$\hat X(m)^*\simeq(\Z/p^m\Z)_{R_2}^2$. Finally, this shows
that $\hat X(n)\simeq\mu_{p^m,R_2}\times\mu_{p^m,R_2}$,
whence $\hat X\simeq\G_{m,R_2}\times\G_{m,R_2}$.
From claim \ref{cl_best.left.reader}, we deduce that 
$\chi$ lifts to an automorphism
$\chi'$ of $\G_{m,R_1}(n)\times\G_{m,R_1}(n)$. Let 
$i:\G_{m,R_1}(n)\to\G_{m,R_1}(n)\times\G_{m,R_1}(n)$,
$\pi:\G_{m,R_1}(n)\times\G_{m,R_1}(n)\to\G_{m,R_1}(n)$
be respectively the imbedding of the first factor,
and the projection onto the second factor; clearly
$\pi\circ\chi'\circ i$ yields a lifting of $f(T)$, as
required.
\end{proof}

\sset\subsubsection{}
Next, let us remark that, for every $n\in\N$, the 
polynomial 
$(1+T)^X-1:=X\cdot T+\binom{X}{2}\cdot T^2+...+
\binom{X}{n}\cdot T^n\in\Q[X,T]$ is an endomorphism
of $\G_{m,\Q[X]}(n)$. As a consequence, there is
a unique ring homomorphism 
$\cG_n\to\Z[X,\binom{X}{2},...,\binom{X}{n}]$ 
representing this endomorphism. The following theorem
will show that this homomorphism is an isomorphism. 

\begin{theorem} The functor 
$$\Z\Alg\to\Set\qquad R\mapsto\End_R(\G_{m,R}(n))$$
is represented by the ring 
$\Z[X,\binom{X}{2},...,\binom{X}{n}]$.
\end{theorem}
\begin{proof} The above discussion has already furnished
us with a natural surjective map 
$\rho:\cG_n\to\Z[X,\binom{X}{2},...,\binom{X}{n}]$. Therefore,
it suffices to show that this map is injective. 
\begin{claim}\label{cl_isom.over.Q}
$\rho\otimes_\Z\one_\Q$ is an isomorphism.
\end{claim}
\begin{pfclaim} First of all, the map $\rho$ can be
characterized in the following way. The identity map
$\cG_n\to\cG_n$ determines an endomorphism 
$f(T):=a_0+a_1\cdot T+...+a_n\cdot T^n$ of $\G_{m,\cG_n}(n)$;
then $\rho$ is the unique ring homomorphism such that
$\rho(f):=f(a_0)+f(a_1)\cdot T+...+f(a_n)\cdot T^n=(1+T)^X-1$.
On the other hand, the ring $\cG_n\otimes_\Z\Q$ represents 
endomorphisms of $\G_m(n)$ in the category of $\Q$-algebras. 
However, for every $n\in\N$ and for every $\Q$-algebra $R$,
there is an isomorphism
$$\log:
\G_{m,R}(n)\stackrel{\sim}{\longrightarrow}\G_{a,R}(n)$$
to the $n$-bud of the additive formal group $\hat\G_{a,R}$.
The endomorphism group of $\G_{a,R}(n)$ is easily computed,
and found to be isomorphic to $R$. In other words, the
universal ring representing endomorphisms of $\G_a(n)$
over $\Q$-algebras is just $\Q[X]$, and the bijection 
$\Hom_{\Q\Alg}(\Q[X],R)\simeq\End(\G_{a,R}(n))$ assigns
to a homomorphism $\phi:\Q[X]\to R$, the endomorphism
$g_\phi(T):=\phi(X)\cdot T$. It follows that, for any 
$\Q$-algebra $R$ there is a natural bijection 
$\Hom_{\Q\Alg}(\Q[X],R)\simeq\End(\G_{m,R}(n))$ given by:
$(\phi:\Q[X]\to R)\mapsto
\exp(\phi(X)\cdot\log(1+T))-1=(1+T)^{\phi(X)}-1$.
Especially, $f(T)$ can be written in the form 
$(1+T)^{\psi(X)}-1$ for a unique ring homomorphism 
$\psi:\Q[X]\to\cG_n\otimes_\Z\Q$. Clearly $\psi$ is
inverse to $\rho\otimes_\Z\one_\Q$.
\end{pfclaim}

In view of claim \ref{cl_isom.over.Q}, we are thus
reduced to show that $\cG_n$ is a flat $\Z$-algebra,
which follows readily from proposition \ref{prop_cG_n.smooth}.
\end{proof}

\sset\subsubsection{}
Furthermore, $\cG_n$ is endowed with a co-addition,
{\em i.e.} a ring homomomorphism 
$\cG_n\to\cG_n\otimes_\Z\cG_n$ satisfying the usual
co-associativity and co-commutativity conditions.
The co-addition is given by the rule:
$$\mathrm{coadd}:\cG_n\to\cG_n\otimes_\Z\cG_n\qquad
\binom{X}{k}\mapsto
\sum_{i+j=k}\binom{X}{i}\otimes\binom{X}{j}.$$
Moreover, for every $k\in\Z$, we have a ring homomorphism
$\pi_k:\cG_n\to\Z$, which corresponds to the endomorphism
of $\G_{m,\Z}(n)$ given by the rule: 
$T\mapsto(1+T)^k-1$ (raising to the $k$-th power in
$\G_{m,\Z}(n)$). Hence we derive, for every $k\in\Z$,
a ring homomorphism 
\set\begin{equation}\label{eq_comp.coadd}
\xymatrix{
\cG_n \ar[rr]^-{\mathrm{coadd}} & &
\cG_n\otimes_\Z\cG_n \ar[rr]^-{\one_{\cG_n}\otimes\pi_k} & &
\cG_n.
}\end{equation}
\begin{remark}\label{rem_choose.plus.i}
(i) On $\cG_n\otimes_\Z\Q=\Q[X]$, \eqref{eq_comp.coadd} is the 
unique map such that $\binom{X}{i}\mapsto\binom{X+k}{i}$ for 
all $i\leq n$, therefore we see that 
$\binom{X+k}{i}\in\cG_n$ for all $k\in\Z$, $n\geq 0$ and
$0\leq i\leq n$. Moreover, \eqref{eq_comp.coadd} is clearly
an automorphism for every $k\in\Z$.

(ii) It is also interesting (though it will not be needed
in this work) to remark that $\cG_n$ is endowed additionally
with a co-composition structure, so that $\cG_n$ is actually
a co-ring, and it represents the functor 
$R\mapsto\End(\G_{m,R}(n))$ from $\Z$-algebras to unitary 
rings. One can check that the co-composition map is given 
by the rule: 
$$\binom{X}{k}\mapsto
\sum_\phi\binom{X}{\phi}\otimes
\prod_{j\in\N^*}\binom{Y}{j}^{\phi(j)}$$
where $\phi$ ranges over all the functions 
$\phi:\N^*:=\N\setminus\{0\}\to\N$ subject to the condition
that $\sum_{j\in\N^*}j\cdot\phi(j)=k$, and 
$\binom{X}{\phi}:= 
\frac{X(X-1)\cdot...\cdot(X-\sum_{j\in\N^*}\phi(j)+1)}
{\prod_{j\in\N^*}\phi(j)!}$. To show that 
$\binom{X}{\phi}\in\cG_n$, one notices that 
$\binom{X}{\phi}=\prod_{j\in\N^*}
\binom{X-\sum^{j-1}_{i=1}\phi(i)}{\phi(j)}$ and then
uses (i).
\end{remark}

\sset\subsubsection{}
For the rest of this section we fix a prime number $p$
and we let $v_p:\Q\to\Z\cup\{\infty\}$ be the $p$-adic valuation.
\begin{lemma}\label{lem_reduce.to.pows}
The ring $\cG_{n,(p)}:=\cG_n\otimes_\Z\Z_{(p)}$ 
is the $\Z_{(p)}$-algebra generated by the polynomials 
$X$, $\binom{X}{p}$,$\binom{X}{p^2}$,...,$\binom{X}{p^k}$, 
where $k$ is the unique integer such that $p^k\leq n<p^{k+1}$.
\end{lemma}
\begin{proof} We proceed by induction on $n$. It suffices to 
prove that $\binom{X}{n}$ is contained in the $\Z_{(p)}$-algebra
$R:=\Z_{(p)}[\binom{X}{p},\binom{X}{p^2},...,\binom{X}{p^k}]$.  
We will use the following (easily verified) identity which holds 
in $\Q[X]$ for every $i,j\in\N$ : 
\set\begin{equation}\label{eq_choose.i.plus.j}
\binom{X}{i+j}=\binom{X}{i}\cdot
               \binom{X-i}{j}\cdot\binom{i+j}{j}^{-1}.
\end{equation} 
Suppose first that $n$ is a multiple of $p^k$, and write
$n=(b+1)p^k$ for some $b<p-1$. If $b=0$, there is nothing
to prove, so we can even assume that $b>0$. We apply 
\eqref{eq_choose.i.plus.j} with $i=b\cdot p^k$ and $j=p^k$.
By remark \ref{rem_choose.plus.i}(i), $\binom{X-b\cdot p^k}{p^k}$
is in $R$, and so is $\binom{X}{b\cdot p^k}$, by induction.
The claim will therefore follow in this case, if we show 
that $\binom{(b+1)p^k}{p^k}$ is invertible in $\Z_{(p)}$.
However, this is clear, since $v_p(i)=v_p(i+b\cdot p^k)$
for every $i=1,...,p^k$. Finally, it remains consider the 
case where $n=b\cdot p^k+a$ for some $b>0$ and $0<a<p^k$.
This is dealt with in the same way: apply 
\eqref{eq_choose.i.plus.j} with $i=b\cdot p^k$ and $j=a$
and use the previous case.
\end{proof}

\begin{lemma}\label{lem_combinatorics} 
Let $k\in\N$. If $R$ is a flat $\Z_{(p)}$-algebra and 
$f\in R$, then the following two conditions are equivalent:
\begin{enumerate}
\item
$\binom{f}{p^i}\in R$ for every $i=1,...,k$;
\item
locally on $\Spec R$ there exists $j\in\Z$
such that $f\equiv j\pmod{p^k}$.
\end{enumerate}
\end{lemma}
\begin{proof} We may assume that $R$ is local. For $k=0$
there is nothing to prove. For $k=1$ we have 
$\binom{f}{p}=u\cdot p^{-1}\cdot\prod^{p-1}_{i=0}(f-i)$
for a unit $u$ of $R$. Then the assertion holds since
all but possibly one of the $f-i$ are invertible.
For $k>1$, by induction we can write $f=i+p\cdot g$ for 
some $g\in R$ and $0\leq i<p$. Since 
$v_p(p^k!)=1+p+p^2+...+p^{k-1}$, we have
$$\binom{f}{p^k}=
u\cdot p^{-1-p-p^2-...-p^{k-1}}\cdot\prod_
{\buildrel j\equiv i\pmod{p}\over
{{}_{0\leq j<p^k}}}(f-j)=u'\cdot\binom{g}{p^{k-1}}
$$
for some units $u,u'\in R$. The claim follows.
\end{proof}

\sset\subsubsection{}\label{subsec_glue.X}
For every integer $k\geq 0$, we construct a scheme $X_k$
by gluing the affine schemes 
$\Spec\,\Z_{(p)}[\frac{X-i}{p^k}]$ ($0\leq i<p^k$) 
along their general fibres. For every $k\in\N$
and every $i\in\N$ with $0\leq i<p^{k+1}$ there 
is an obvious imbedding 
$\Z_{(p)}[\frac{X-i}{p^k}]\subset
\Z_{(p)}[\frac{X-i}{p^{k+1}}]$. By gluing the
duals of these imbeddings, we obtain, for every 
$k\in\N$, a morphism of schemes 
$\rho_k:X_{k+1}\to X_k$. Let also 
$\xi_k:\Spec\,\cG_{p^{k+1}}\to\Spec\,\cG_{p^k}$ 
be the morphism which is dual to the imbedding 
$\cG_{p^k}\subset\cG_{p^{k+1}}$.

\begin{proposition}\label{prop_scheme.descript}
With the notation of \eqref{subsec_glue.X} we have:
\begin{enumerate}
\item
For given $n>0$, let $k$ be the unique
integer such that $p^k\leq n<p^{k+1}$. Then there 
is a natural isomorphism of schemes: 
$\pi_k:X_k\stackrel{\sim}{\longrightarrow}
\Spec\,\cG_n\otimes_\Z\Z_{(p)}$.
\item
For every $k\in\N$ the diagram of schemes:
$$\xymatrix{
X_{k+1} \ar[r]^-{\pi_{k+1}} \ar[d]_-{\rho_k} &
\Spec\,\cG_{p^{k+1}}\otimes_\Z\Z_{(p)} 
\ar[d]^{\xi_k\otimes_\Z\Z_{(p)}} \\
X_k \ar[r]^-{\pi_k} &
\Spec\,\cG_{p^k}\otimes_\Z\Z_{(p)}
}$$
commutes.
\end{enumerate}
\end{proposition}
\begin{proof} By lemma \ref{lem_reduce.to.pows} we
may assume that $n=p^k$. By lemma \ref{lem_combinatorics}, 
we see that both $X_k$ and 
$\Spec\,\cG_{p^k}\otimes_\Z\Z_{(p)}$ represent the
same functor from the category of flat $\Z_{(p)}$-schemes
to the category of sets.
Since both schemes are flat over $\Spec\,\Z_{(p)}$,
(i) follows. It is similarly clear that 
$\xi_k\otimes_\Z\Z_{(p)}$ and $\rho_k$ represent
the same natural transformation of functors, so (ii)
follows.
\end{proof}

\begin{corollary}\label{cor_scheme.descrip}
{\em (i)} For given $n\in\N$, let $k$ be the unique integer
such that $p^{k-1}\leq n<p^k$. Then there is a natural
ring isomorphism
$$\Z/p^k\Z\stackrel{\sim}{\longrightarrow}
\End\,(\G_{m,\F_p}(n))\quad
i\mapsto(1+T)^i-1$$
\begin{enumerate}
\addtocounter{enumi}{1}
\item
Let $R$ be a ring such that
$\F_p\subset R$. Then there is a natural ring
isomorphism
$$\cC^0(\Spec\,R,\Z_p)\stackrel{\sim}{\longrightarrow}
\End_R\,(\hat\G_{m,R})\ \qquad\beta\mapsto(1+T)^\beta-1.$$
\end{enumerate}
\end{corollary}
\begin{proof} (i): by lemma \ref{lem_reduce.to.pows}
we can assume $n=p^k-1$. In this case, it is clear
that the polynomials $(1+T)^i-1$ are all distinct for
$i=0,...,p^k-1$ and they form a subring of 
$\End(\G_{m,\F_p}(n))$. However, an endomorphism
of $\G_{m,\F_p}(n)$ corresponds to a unique point
in $\Spec\,\cG_n(\F_p)$. From proposition 
\ref{prop_scheme.descript}(i) we derive that 
$\Spec\,\cG_n\otimes_\Z\F_p$ is the union
of the special fibres of the affine schemes 
$\Spec\,\Z[\frac{X-i}{p^{k-1}}]$, for $i=0,...,p^{k-1}-1$.
Each of those contribute an affine line $\A^1_{\F_p}$,
so $\Spec\,\cG_n\otimes_\Z\F_p$ consists of
exactly $p^{k-1}$ connected components. In total,
we have therefore exactly $p^k$ points in
$\Spec\,\cG_n(\F_p)$, so (i) follows.

(ii): to give an endomorphism of $\hat\G_m$
is the same as giving a compatible system of endomorphisms
of $\G_m(n)$, one for each $n\in\N$. In case 
$\F_p\subset R$, lemma \ref{lem_reduce.to.pows} shows 
that this is also equivalent to the datum of a compatible 
system of morphisms 
$\phi_k:\Spec R\to\Spec\,\cG_{p^k}\otimes_\Z\F_p$, for 
every $k\geq 0$.  From proposition
\ref{prop_scheme.descript}(ii) we can further
deduce that, under the morphism $\xi_k$, each
of the $p^{k+1}$ connected components of
$\Spec\,\cG_{p^{k+1}}\otimes_\Z\F_p$ gets mapped
onto one of the $p^{k+1}$ rational points of 
$\Spec\,\cG_{p^k}\otimes_\Z\F_p$.
Since $\phi_{k-1}=\xi_k\circ\phi_k$, we see that the 
image of $\phi_{k-1}$ is contained in 
$\Spec\,\cG_{p^{k-1}}(\F_p)$, for every $k>0$. 
Taking (i) into account, we see that an endomorphism of 
$\hat\G_{m,R}$ is the same as the datum of a compatible 
system of continuous maps 
$\Spec R\to\Z/p^k\Z$. Since the $p$-adic topology of $\Z_p$
is the inverse limit of the discrete topologies on the
$\Z/p^k\Z$, the claim follows.
\end{proof}

\subsection{Modules of almost finite rank}
\label{sec_alm.fin.rk} 
Let $A$ be a $V^a$-algebra, $P$ an almost finitely generated 
projective $A$-module and $\phi\in\End_A(P)$. 

\sset\subsubsection{}\label{subsec_lambda_nilpotent}
\index{Almost module(s)!$\Lambda$-nilpotent endomorphism of an
|indref{subsec_lambda_nilpotent}}
\index{$\det(\one_P+\phi)$|indref{subsec_lambda_nilpotent}}
We say that $\phi$ is {\em $\Lambda$-nilpotent\/} if there 
exists an integer $i>0$ such that $\Lambda^i_A\phi=0$. 
Notice that the $\Lambda$-nilpotent endomorphisms of $P$
form a bilateral ideal of the unitary ring $\End_A(P)$. 
Notice also that $\Lambda^i_AP$ is an almost projective 
$A$-module for every $i\geq 0$; indeed, this is easily 
shown by means of lemma \ref{lem_thetwo}. For a 
$\Lambda$-nilpotent endomorphism $\phi$ we introduce the 
notation
$$\det(\one_P+\phi):=
\sum_{i\geq 0}\tr_{\Lambda^i_AP/A}(\Lambda_A^i\phi).$$
Notice that the above sum consists of only finitely
many non-zero terms, so that $\det(\one_P+\phi)$ is a well
defined element of $A_*$.
\begin{lemma}\label{lem_commute.det} Let $P$ be an 
almost finitely generated projective $A$-module.
\begin{enumerate}
\item
If $\phi$ is a $\Lambda$-nilpotent endomorphism of $P$ and
$\alpha:A\to A'$ is any morphism of $V^a$-algebras, set 
$P':=P\otimes_AA'$. Then: 
$\det(\one_{P'}+\phi\otimes_A\one_{A'})=
\alpha(\det(\one_P+\phi))$.
\item
Let $\phi,\psi\in\End_A(P)$ such that $\phi\circ\psi$ 
and $\psi\circ\phi$ are $\Lambda$-nilpotent. Then:
$\det(\one_P+\phi\circ\psi)=\det(\one_P+\psi\circ\phi)$.
\end{enumerate}
\end{lemma}
\begin{proof} (i) is a straightforward consequence of the 
definitions. As for (ii), it is clear that $\psi\circ\phi$ 
is $\Lambda$-nilpotent and the stated identity follows 
directly from lemma \ref{lem_swap.maps}(i).
\end{proof}

\sset\subsubsection{}
Now, let $\phi,\psi\in\End_A(P)$ be two endomorphisms.
Set $B:=A[X,Y]/(X^n,Y^n)$ and $P_B:=P\otimes_AB$; $\phi$ and 
$\psi$ induce endomorphisms of $P_B$ that we denote again
by the same letters. Clearly $X\cdot\phi$ and $Y\cdot\psi$ 
are $\Lambda$-nilpotent; hence we get elements 
$\det(\one_{P_B}+X\cdot\phi)$,
$\det(\one_{P_B}+Y\cdot\psi)$ and 
$\det(\one_{P_B}+X\cdot\phi+Y\cdot\psi+XY\cdot\psi\circ\phi)$ 
in $B_*$. Notice that any element of $B_*$ can be written 
uniquely as an $A_*$-linear combination of the monomials 
$X^iY^j$ with $0\leq i,j<n$. Moreover, it is clear that 
$\det(\one_{P_B}+X\cdot\phi)=\sum_{0\leq i<n}
\tr_{\Lambda^i_AP/A}(\Lambda_A^i\phi)\cdot X^i$,
and similarly for $\psi$.
\begin{proposition}\label{prop_determ.prods}
With the above notation, the following identity holds :
\set\begin{equation}\label{eq_determ}
\det(\one_{P_B}+X\cdot\phi)\cdot\det(\one_{P_B}+Y\cdot\psi)=
\det(\one_{P_B}+X\cdot\phi+Y\cdot\psi+XY\cdot\psi\circ\phi).
\end{equation}
\end{proposition}
\begin{proof} First of all we remark that, when $P$ is
a free $A$-module of finite rank, the above identity is
well-known, and easily verified by working with matrices
with entries in $A_*$. Suppose next that $P$
is arbitrary, but $\phi=\eps\cdot\phi'$, $\psi=\eps\cdot\psi'$
for some $\phi',\psi'\in\End_A(P)$ and $\eps\in\fm$. 
Pick a free $A$-module $F$ of finite rank, and morphisms 
$u:P\to F$, $v:F\to P$ such that $v\circ u=\eps\cdot\one_P$. Set 
$\phi_\eps:=u\circ\phi'\circ v:F\to F$ and define similarly 
$\psi_\eps$. Clearly 
$\det(\one_{P_B}+X\cdot\phi)=
\det(\one_{P_B}+X\cdot\eps\cdot\phi')=
\det(\one_{P_B}+X\cdot v\circ u\circ\phi')
=\det(\one_{P_B}+X\cdot\phi_\eps)$ (by lemma 
\ref{lem_commute.det}(ii)) and similarly for the other terms 
appearing in \eqref{eq_determ}. Thus we have reduced this 
case to the case of a free $A$-module.
Finally, we deal with the general case. The foregoing shows that
the sought identity is known at least when $\phi$ and $\psi$
are replaced by $\eps\cdot\phi$, resp. $\eps\cdot\psi$, for
any $\eps\in\fm$. Equivalently, consider the $A$-algebra
endomorphism $\alpha:B\to B$ defined by $X\mapsto\eps\cdot X$, 
$Y\mapsto\eps\cdot Y$ and let $C$ be the $B$-algebra
structure on $B$ determined by $\alpha$; by lemma 
\ref{lem_commute.det}(i) we have 
$$\det(\one_{P_B}+X\cdot\eps\cdot\phi)=
\det(\one_{P_C}+(X\cdot\phi)\otimes_B\one_C)=
\alpha(\det(\one_{P_B}+X\cdot\phi))$$
and similarly for the other terms appearing in 
\eqref{eq_determ}. Thus, the images under $\alpha$
of the two members of \eqref{eq_determ} coincide. But
applying $\alpha$ to a monomial of the form $a\cdot X^iY^j$
has the effect of multiplying it by $\eps^{i+j}$;
by ({\bf B}), the $(i+j)$-powers of 
elements of $\fm$ generate $\fm$, hence the claim
follows easily.
\end{proof}
\begin{corollary} If $\phi,\psi\in\End_A(P)$ are two
$\Lambda$-nilpotent endomorphisms, then
$$\det(\one_{P}+\phi)\cdot\det(\one_P+\psi)=
\det(\one_P+\phi+\psi+\phi\circ\psi).$$
\end{corollary}
\begin{proof} For an arbitrary $\alpha\in\End_A(P)$, 
one can define $P':=P\otimes_AA[[X]]$ and
$\det(\one_{P'}+X\cdot\alpha):=\sum_{i\geq 0}
\tr_{\Lambda^i_AP/A}(\Lambda_A^i\alpha)\cdot X^i
\in A_*[[X]]$. Then proposition \ref{prop_determ.prods}
implies that the analogue of \eqref{eq_determ} holds 
in $A_*[[X,Y]]$. But if $\alpha$ is $\Lambda$-nilpotent, 
the power series $\det(\one_{P'}+X\cdot\alpha)$ is
actually a polynomial in $A_*[X]$; the claim then follows
by evaluating the polynomials 
$\det(\one_{P'}+X\cdot\phi)$, $\det(\one_{P'}+Y\cdot\psi)$
and $\det(\one_{P'}+
X\cdot\phi+Y\cdot\psi+XY\cdot\phi\circ\psi)$ for $X=Y=1$.
\end{proof}

\sset\subsubsection{}\label{subsec_ch.pow-srs}
\index{$\chi_P(X)$, $\psi_P(X)$|indref{subsec_ch.pow-srs}}
Next, for $P$ as above, set 
$$\begin{array}{r@{~:=~}l}
\chi_P(X) &
\sum_{i\geq 0}\tr_{\Lambda^i_AP/A}(\Lambda_A^i\one_P)\cdot X^i
\in A_*[[X]]\\
\psi_P(X) &
\sum_{i\geq 0}\tr_{\Sym^i_AP/A}(\one_{\Sym^i_AP})\cdot X^i
\in A_*[[X]].
\end{array}$$

\begin{corollary}\label{cor_psi.Sym} Let $P$ be an almost 
finitely generated projective $A$-module. Then:
\begin{enumerate}
\item
the power series $\chi_P(X)$ defines an endomorphism of 
the formal group $\hat\G_{m,A_*}$.
\item
$\chi_P(X)\cdot\psi_P(-X)=1$.
\item
$\chi_P(X)\in 1+\cE_{P/A*}[[X]]$.
\end{enumerate}
\end{corollary}
\begin{proof} (i) is immediate. For (ii), recall that, 
for every $n>0$ there is an acyclic Koszul complex
(cp. \cite[Ch.X, \S9, n.3, Prop.3]{BouAH})
$$0\to\Lambda^n_AQ\to(\Lambda^{n-1}_AQ)\otimes_A(\Sym^1_AQ)
\to...\to(\Lambda^1_AQ)\otimes_A(\Sym^{n-1}_AQ)\to
\Sym^n_AQ\to 0.$$
From proposition \ref{prop_short.ex.seq.trace} we derive,
by a standard argument, that the trace is an additive 
function on arbitrary bounded acyclic complexes. Then,
taking into account lemma \ref{lem_trace.tensors} we
obtain: $\sum_{i=0}^n(-1)^i\cdot
\tr_{\Lambda^{n-i}_AP/A}(\one_{\Lambda^{n-i}_AP})\cdot
\tr_{\Sym^i_AP/A}(\one_{\Sym^i_AP})=0$
for every $n>0$. This is equivalent to the sought 
identity. To show (iii) we remark more precisely
that $\cE_{\Lambda^r_AP/A}\subset\cE_{P/A}$ for every $r>0$.
Indeed, set $B:=A/\cE_{P/A}$. Then, by proposition 
\ref{prop_eval.ideal}(i),(iii):
$\cE_{\Lambda^r_AP/A}\cdot B=
\cE_{\Lambda^r_B(P\otimes_AB)/B}=0$, whence the claim.
\end{proof}

\begin{definition}\label{def_formal_rank}
\index{Almost module(s)!almost projective!$\frk_A(P)$ : formal rank of an
|indref{def_formal_rank}}
\index{Almost module(s)!almost projective!of almost finite rank
|indref{def_formal_rank}}
\index{Almost module(s)!almost projective!of finite rank
|indref{def_formal_rank}}
\index{Almost module(s)!almost projective!of constant rank
|indref{def_formal_rank}}
Let $P$ be an almost finitely generated projective $A$-module. 
\begin{enumerate}
\item
The {\em formal rank\/} of $P$ is the ring homomorphism 
$\frk_A(P):\cG_\infty:=\Z[\alpha,\binom{\alpha}{2},...]\to A_*$ 
associated to $\chi_P(X)$.
\item
We say that $P$ is {\em of almost finite rank\/} if, 
for every $\eps\in\fm$, there exists an integer $i\geq 0$ 
such that $\eps\cdot\Lambda^i_AP=0$.
\item
We say that $P$ is {\em of finite rank} if there
exists an integer $i\geq 0$ such that $\Lambda^i_AP=0$.
\item
Let $r\in\N$; we say that $P$ has {\em constant rank equal to $r$\/}
if $\Lambda^{r+1}_AP=0$ and $\Lambda^r_AP$ is an invertible
$A$-module.
\end{enumerate}
\end{definition}

\begin{remark}\label{rem_continuous.fct} 
(i): It follows easily from lemma \ref{lem_ancora}(vi) that 
every uniformly almost finitely generated projective $A$-module 
is of finite rank.

(ii): Notice that if $P$ is of finite rank,
then $\chi_P(X)$ is a polynomial, whence it defines
an endomorphism of the algebraic group $\G_{m,A_*}$. 
In this case, it follows that $\chi_P(X)$ is of the
form $(1+X)^\alpha$, where $\alpha:\Spec\,A_*\to\Z$
is a continuous function (where $\Z$ is seen as a discrete
topological space). More precisely, there is an obvious
injective ring homomorphism
\set\begin{equation}\label{eq_continuous.fct}
\cC^0(\Spec\,A_*,\Z)\to\End_A(\G_{m,A_*}) \qquad
\beta\mapsto(1+X)^\beta
\end{equation}
which allows to identify the continuous function
$\alpha$ with the formal rank of $P$. Moreover, if 
$\Lambda^i_AP=0$, it is clear that 
$\alpha(\Spec\,A_*)\subset\{0,...,i-1\}$. 
\end{remark}

The main result of this section is theorem 
\ref{th_structure.alm.fin.rk}, which describes  general
modules of almost finite rank as infinite products of
modules of finite rank. The first step is lemma 
\ref{lem_lambda.two}, concerned with the case of an 
$A$-module of rank one. 

\begin{lemma}\label{lem_lambda.two} 
Let $P$ be an almost finitely generated
projective $A$-module such that $\Lambda^2_AP=0$. 
There exists $V^a$-algebras $A_0$, $A_1$ and an 
isomorphism of $V^a$-algebras $A\simeq A_0\times A_1$ 
such that $P\otimes_AA_0=0$ and $P\otimes_AA_1$ is an 
invertible $A_1$-module.
\end{lemma}
\begin{proof} Since the natural map
$P\times P\to\Lambda^2_AP$ is universal for alternating
$A$-bilinear maps on $P\times P$, we have 
\set\begin{equation}\label{eq_alternate}
f(p)\cdot q=f(q)\cdot p\qquad
\text{for every $f\in (P^*)_*$ and $p,q\in P_*$.}
\end{equation}
Using \eqref{eq_alternate} we derive
$\omega_{P/A}(p\otimes f)(q)=
\tr_{P/A}(\omega_{P/A}(p\otimes f))\cdot q$
for every $f\in (P^*)_*$ and $p,q\in P_*$. In other
words, $\omega_{P/A}(p\otimes f)=
\tr_{P/A}(\omega_{P/A}(p\otimes f))\cdot\one_P$, for every
$f\in (P^*)_*$, $p\in P_*$. By linearity we finally
deduce
\set\begin{equation}\label{eq_trace.multiples}
\phi=\tr_{P/A}(\phi)\cdot\one_P \qquad 
\text{for all $\phi\in\End_A(P)$}.
\end{equation}
Now, by remark \ref{rem_continuous.fct}(ii), the hypothesis 
$\Lambda^2_AP=0$ also implies that $\chi_P(X)=(1+X)^\alpha$, 
for a continuous function $\alpha:\Spec\,A_*\to\{0,1\}$.
We can decompose accordingly $A=A_0\times A_1$, so
that $\alpha(\Spec\,A_{i*})=i$, which gives the sought
decomposition. We can now treat separately the two cases 
$A=A_0$ and $A=A_1$. In case $\frk_A(P)=0$, then 
$\tr_{P/A}(\one_P)=0$, and then \eqref{eq_trace.multiples} 
implies that $P=0$. In case $\frk_A(P)=1$, then 
$\tr_{P/A}(\one_P)=1$ and \eqref{eq_trace.multiples} 
implies that the natural map $
A\to\End_A(P)~:~a\mapsto a\cdot\one_P$ is an inverse
for $\tr_{P/A}$, thus $A\simeq P\otimes_AP^*$.
\end{proof}

The next step consists in analyzing the structure
of $A$-modules of finite rank. To this purpose we
need some preliminaries of multi-linear algebra.

\sset\subsubsection{}
For every $n\geq 0$ let $\mathbf n:=\{1,...,n\}$; for a 
subset $I\subset{\mathbf n}$ let $|I|$ be the cardinality
of $I$; for a given partition ${\mathbf n}=I\cup J$, let 
$\prec$ denote the total ordering on $\mathbf n$ that 
restricts to the usual ordering on $I$ and on $J$, and 
such that $i\prec j$ for every $i\in I$, $j\in J$. Finally 
let $\eps_{IJ}$ be the sign of the unique order-preserving 
bijection $({\mathbf n},<)\to({\mathbf n},\prec)$.

Let $M$ be any $A$-module. Given elements $m_1,m_2,...,m_n$
in $M_*$, and $I\subset{\mathbf n}$ a subset of elements
$i_1<i_2<...<i_{|I|}$, let
$m_I:=m_{i_1}\wedge...\wedge m_{i_{|I|}}
\in\Lambda^{|I|}_AM_*$ (with the convention that 
$m_\emptyset=1\in A_*=\Lambda_A^0M_*$).

\sset\subsubsection{}
Let $M,N$ be any two $A$-modules. 
For every $i,j\geq 0$ there is a natural morphism
\set\begin{equation}\label{eq_map.wedge}
\Lambda^i_AM\otimes_A\Lambda^j_AN\to
\Lambda^{i+j}_A(M\oplus N)
\end{equation}
determined by the rule:
$$m_1\wedge...\wedge m_i\otimes n_1\wedge...\wedge n_j
\mapsto(m_1,0)\wedge(m_2,0)\wedge...\wedge
(0,n_1)\wedge...\wedge(0,n_j)
$$
for all $m_1,...,m_i\in M_*$ and $n_1,...,n_j\in N$.
The morphisms \eqref{eq_map.wedge} assemble to
an isomorphism of $A$-modules
\set\begin{equation}\label{eq_iso.wedge}
\Lambda^\bullet_AM\otimes_A\Lambda^\bullet_AN\to
\Lambda_A^\bullet(M\oplus N).
\end{equation}
Clearly, there is a unique graded $A$-algebra structure
on $\Lambda^\bullet_AM\otimes_A\Lambda^\bullet_AN$
such that \eqref{eq_iso.wedge} is an isomorphism
of (graded-commutative) $A$-algebras. Explicitly,
given $x_i\in\Lambda^{a_i}_AM$, 
$y_i\in\Lambda ^{b_i}_AN$ ($i=1,2$) one verifies
easily that the product on 
$\Lambda^\bullet_AM\otimes_A\Lambda^\bullet_AN$
is fixed by the rule
\set\begin{equation}\label{eq_rules.product}
(x_1\otimes y_1)\cdot(x_2\otimes y_2)=
(-1)^{a_2b_1}\cdot(x_1\wedge x_2)\otimes(y_1\wedge y_2).
\end{equation}
Then $\Lambda^\bullet_AM\otimes_A\Lambda^\bullet_AN$
is even a bigraded $A$-algebra, if we let 
$\Lambda^i_AM\otimes_A\Lambda^j_AN$ be the graded
component of bidegree $(i,j)$.

\sset\subsubsection{}
Next, let $\delta:M\to M\oplus M$ be the diagonal
morphism $m\mapsto(m,m)$ (for all $m\in M_*$).
It induces a morphism 
$\Lambda^\bullet_A\delta:
\Lambda^\bullet_AM\to\Lambda^\bullet_A(M\oplus M)$
of $A$-algebras.
We let $\Delta:\Lambda^\bullet_AM\to
\Lambda^\bullet_AM\otimes_A\Lambda^\bullet_AN$
be the composition of the morphism 
$\Lambda^\bullet_A\delta$ and the inverse of the
isomorphism \eqref{eq_iso.wedge}. 
For every $a,b\geq 0$ we also let 
$\Delta_{a,b}:\Lambda^\bullet_AM\to
\Lambda^a_AM\otimes_A\Lambda^b_AN$ be the composition
of $\Delta$ and the projection onto the graded
component of bidegree $(a,b)$. The morphisms 
$\Delta_{a,b}$ are usually called "co-multiplication 
morphisms". An easy calculation shows that:
\set\begin{equation}\label{eq_Delta.a.b}
\Delta_{a,b}(x_1\wedge x_2\wedge...\wedge x_{a+b})
=\sum_{I,J}\eps_{IJ}\cdot x_I\otimes x_J
\end{equation}
where the sum ranges over all the partitions 
$\mathbf{a+b}=I\cup J$ such that $|I|=a$.
Let now $x_1,...,x_a,y_1,...,y_b\in M_*$.
Since $\Delta$ is a morphism of $A$-algebras,
we have $\Delta(x_{\mathbf a}\wedge y_{\mathbf b})=
\Delta(x_{\mathbf a})\cdot\Delta(y_{\mathbf b})$.
Hence, using \eqref{eq_rules.product} and 
\eqref{eq_Delta.a.b} one deduces easily:
\set\begin{equation}\label{eq_Delta.explode}
\Delta_{a,b}(x_{\mathbf a}\wedge y_{\mathbf b})=
\sum_{I,J,K,L}\eps_{IJ}\cdot\eps_{KL}\cdot(-1)^{|J|}
\cdot(x_I\wedge y_K)\otimes(x_J\wedge y_L).
\end{equation}
where the sum runs over all partitions 
$I\cup J={\mathbf a}$, $K\cup L={\mathbf b}$ such that
$|J|=|K|$.
\begin{lemma}\label{lem_identity.Lambda}
Suppose that $\Lambda^{a+1}_AM=0$ for some integer
$a\geq 0$. Let $0<b\leq a$ and 
$x_1,...,x_a,y_1,...,y_b\in M_*$. 
Then the following identity holds in 
$\Lambda^a_AM\otimes_A\Lambda^b_AM$ :
$$x_{\mathbf a}\otimes y_{\mathbf b}=
\sum_{I,J}\eps_{JI}\cdot
(x_J\wedge y_{\mathbf b})\otimes x_I$$
where the sum ranges over all the partitions 
${\mathbf a}=I\cup J$ such that $|I|=b$.
\end{lemma}
\begin{proof} For a given subset $B\subset{\mathbf b}$ 
we let 
$$\gamma(y_B):=
\sum_{I,J}\eps_{IJ}\cdot
(x_I\wedge y_B)\otimes x_J-
x_{\mathbf a}\otimes y_B$$
where the sum is taken over all the partitions 
$I\cup J={\mathbf a}$ such that $|J|=|B|$. Notice
that $\gamma(y_\emptyset)=0$.
We have to show that $\gamma(y_{\mathbf b})=0$.
To this purpose we show the following:
\begin{claim}\label{cl_gamma.equal} If $|B|>0$, then
\set\begin{equation}\label{eq_claim.gamma}
\Delta_{a,|B|}(x_{\mathbf a}\wedge y_B)=
\sum_{K,L}\eps_{KL}\cdot(-1)^{|K|}\cdot
\gamma(y_K)\wedge y_L
\end{equation}
where the sum ranges over all the partitions $K\cup L=B$.
\end{claim}
\begin{pfclaim} Using \eqref{eq_Delta.explode}, the 
difference between the two sides of \eqref{eq_claim.gamma} 
is seen to be equal to 
$\sum_{K,L}\eps_{KL}\cdot(-1)^{|K|}\cdot 
x_{\mathbf a}\otimes(y_K\wedge y_L)=
\sum_{K,L}\cdot(-1)^{|K|}\cdot 
x_{\mathbf a}\otimes y_B$, where the sum runs over
all partitions $K\cup L=B$. A standard combinatorial
argument shows that this expression can be rewritten
as $x_{\mathbf a}\otimes y_B\cdot
\sum_{k=0}^{|B|}(-1)^k\cdot\binom{|B|}{k}$, which vanishes
if $|B|>0$.
\end{pfclaim}

To conclude the proof of the lemma, we remark that 
$\Delta_{a,b}$ vanishes if $b>0$ because by assumption
$\Lambda^{a+1}_AM=0$; then the claim follows by induction
on $|B|$, using claim \ref{cl_gamma.equal}.
\end{proof}
\begin{lemma}\label{lem_lambda.of.lambda} 
Let $P$ be an $A$-module such that
$\Lambda^{n+1}_AP=0$ and assume that either $P$
is flat or $2$ is invertible in $A_*$. Then
$\Lambda^2_A(\Lambda^n_AP)=0$.
\end{lemma}
\begin{proof} For any $A$-module $M$ and $r\geq 0$ 
there exists an antisymmetrizer operator (cp. 
\cite[Ch. III, \S 7.4, Remarque]{Bourbaki})
$$a_r:M^{\otimes r}\to M^{\otimes r}\qquad
m_1\otimes...\otimes m_r\mapsto\sum_{\sigma\in S_r}
\mathrm{sgn}(\sigma)\cdot m_{\sigma(1)}\otimes...
\otimes m_{\sigma(r)}.$$
Clearly $a_r$ factors thorugh $\Lambda^r_AM$, and
in case $M$ is free of finite rank, it is easy
to check (just by arguing with basis elements) that 
the induced map 
$\bar a_r:\Lambda^r_AM\to\Img(a_r)$ is an isomorphism.
This is still true also in case $r!$ is invertible
in $A_*$, since in that case one checks that
$a_r/r!$ is idempotent (see {\em loc.~cit.\/}). More 
generally, if $M$ is flat then, by \cite[Ch.I, Th.1.2]{La}, 
$M_!$ is the filtered colimit of a direct system of 
free $A_*$-modules of finite rank, so also in this 
case $\bar a_r$ is an isomorphism. Notice that, again 
by \cite[Ch.I, Th.1.2]{La}, if $P$ is flat, then 
$\Lambda^k_AP$ is also flat, for every $k\geq 0$. 
Hence, to prove the lemma, it suffices to verify that 
$\Img(a_2:(\Lambda^n_AP)^{\otimes 2}\to
(\Lambda^n_AP)^{\otimes 2})=0$ when $\Lambda^{n+1}_AP=0$.
However, this follows easily from lemma 
\ref{lem_identity.Lambda}.
\end{proof}

We are now ready to return to $A$-modules of finite rank.

\begin{proposition}\label{prop_decomp.fin.rank}
Let $P$ be an almost projective $A$-module of finite rank;
say that $\Lambda^r_AP=0$.
There exists a natural decomposition
$A\simeq A_0\times A_1\times...\times A_{r-1}$ such that
$P_i:=P\otimes_AA_i$ is an $A_i$-module of constant rank 
equal to $i$ for every $i=0,...,r-1$.
\end{proposition}
\begin{proof} We proceed by induction on $r$; the
case $r=2$ is covered by lemma \ref{lem_lambda.two}.
By lemma \ref{lem_lambda.of.lambda} we have 
$\Lambda^2_A(\Lambda^{r-1}_AP)=0$, so by lemma 
\ref{lem_lambda.two}, there is a decomposition
$A\simeq A_{r-1}'\times A'_{r-2}$, such that for 
$P_i:=P\otimes_AA'_i$ ($i=r-2,r-1$) the following holds. 
$\Lambda^{r-1}_{A'_{r-2}}(P_{r-2})=0$ and 
$\Lambda^{r-1}_{A'_{r-1}}(P_{r-1})$ is an invertible
$A_{r-1}$-module.
It follows in particular that $\chi_{P_{r-1}}(X)$ is 
a polynomial of degree $r-1$, and its leading coefficient 
is invertible in $A_{r-1}$. Hence 
$\chi_{P_{r-1}}(X)=(1+X)^{r-1}$. By induction,
$A'_{r-2}$ admits a decomposition 
$A'_{r-2}\simeq A_{r-2}\times...\times A_0$ with the
stated properties; it suffices then to take 
$A_{r-1}:=A'_{r-1}$.
\end{proof}

\begin{theorem}\label{th_structure.alm.fin.rk} 
Let $P$ be an almost projective $A$-module
of almost finite rank. Then there exists a natural
decomposition $A\simeq\prod_{i=0}^\infty A_i$ such that:
\begin{enumerate}
\item
$\liminv{i\to\infty}\Ann_{V^a}(A_i)=V^a$\ \ (for the uniform
structure of definition {\em\ref{def_unif.on.A}});
\item
for $i\in\N$, let $P_i:=P\otimes_AA_i$; then 
$P\simeq\prod_{i=0}^\infty P_i$ and every $P_i$ is
an $A_i$-module of finite constant rank equal to $i$.
\end{enumerate}
\end{theorem}
\begin{proof} Let $\{\fm_\lambda\}_{\lambda\in I}$ be 
the filtered family of finitely generated subideals of 
$\fm$. For every $\lambda\in I$, let 
$A_\lambda:=A/\Ann_A(\fm_\lambda)$. By hypothesis,
$P_\lambda:=P\otimes_AA_\lambda$ is an $A_\lambda$-module
of finite rank; say that the rank is $r(\lambda)$.
By proposition \ref{prop_decomp.fin.rank} we have natural 
decompositions 
$A_\lambda\simeq 
A_{\lambda,0}\times...\times A_{\lambda,r(\lambda)}$
such that $P\otimes_AA_{\lambda,i}$ is an $A_{\lambda,i}$-module
of constant rank equal to $i$ for every $i\leq r(\lambda)$. 
The naturality of the decomposition
means that for every $\lambda,\mu\in I$ such that 
$\fm_\lambda\subset\fm_\mu$, we have 
$A_{\mu,i}\otimes_AA_{\lambda}\simeq A_{\lambda,i}$.
for every $i\leq r(\mu)$. In particular, 
$\fm_\lambda A_{\mu,i}=0$ for every $i>r(\lambda)$.
By considering the short exact sequence of cofiltered
systems of $A$-modules:
$$
0\to(\Ann_A(\fm_\lambda))_{\lambda\in I}\to 
(A)_{\lambda\in I}\to(A_\lambda)_{\lambda\in I}\to 0
$$
we deduce easily that 
$A\simeq\liminv{\lambda\in I}A_\lambda$ and therefore
we obtain a decomposition $A\simeq\prod_{i=0}^\infty A_i$,
with $A_i:=\liminv{\lambda\in I}A_{\lambda,i}$ for 
every $i\in\N$. Notice that, for every $i\in\N$
and every $\lambda\in I$, the natural morphism
$A_i\to A_{\lambda,i}$ is surjective with kernel
killed by $\fm_\lambda$. It follows easily that
$\Lambda^{i+1}_{A_i}(P\otimes_AA_i)=0$ and
$\Lambda^i_{A_i}(P\otimes_AA_i)$ is an invertible
$A_i$-module. Furthermore, for every $\lambda\in I$,
$\fm_\lambda\cdot A_i=0$ for all $i>r(\lambda)$,
which implies (i). Finally, for every $\lambda\in I$,
$\fm_\lambda$ kills the kernel of the
projection $P\to\prod^{r(\lambda)}_{i=0}P_i$, so 
$P$ is isomorphic to the infinite product of the $P_i$. 
\end{proof}
\subsection{Localisation in the flat site}
\label{sec_local.flat.site}
Throughout this section $P$ denotes an almost finitely 
generated projective $A$-module. The following 
definition introduces the main tool used in this section. 
\begin{definition}\label{def_splitter}
\index{Almost module(s)!almost projective!$\Split(A,P)$ : splitting
  algebra of an|indref{def_splitter}}
The {\em splitting algebra\/} of $P$ is the $A$-algebra:
$$\Split(A,P):=\Sym^\bullet_A(P\oplus P^*)/(1-\zeta_P).$$
We endow $\Sym^\bullet_A(P\oplus P^*)$ with the structure
of graded algebra such that $P$ is placed in degree one
and $P^*$ in degree $-1$. Then $\zeta_P$ is a
homogeneous element of degree zero, and consequently 
$\Split(A,P)$ is also a graded $A$-algebra.
\end{definition} 

\sset\subsubsection{}
We define a functor $S:A\Alg\to\Set$ by assigning
to every $A$-algebra $B$ the set $S(B)$ of all pairs 
$(x,\phi)$ where $x\in(P\otimes_AB)_*$, 
$\phi:P\otimes_AB\to B$ such that $\phi(x)=1$.

\begin{lemma}\label{lem_Split.generals}
{\em (i)} $\Split(A,P)$ is a flat $A$-algebra.
\begin{enumerate}
\addtocounter{enumi}{1}
\item
$\Split(A',P\otimes_AA')\simeq\Split(A,P)\otimes_AA'$
for every $A$-algebra $A'$.
\item
$\Split(A,P)$ represents the functor $S$.
\end{enumerate}
\end{lemma}
\begin{proof}
For every $k\geq 0$ we have 
$\zeta_P^k\in(\Sym^k_AP)\otimes_A(\Sym^k_AP^*)\subset
\gr^0(\Sym^{2k}_A(P\oplus P^*))$. It is easy to verify the 
formula:
\set\begin{equation}\label{eq_gr.of.splitter}
\gr^k\Split(A,P)\simeq\colim{j\in\Z}
(\Sym^{k+j}_AP)\otimes_A(\Sym^j_AP^*)
\end{equation}
where the transition maps in the direct system are 
given by multiplication by $\zeta_P$. In particular,
it is clear that $\gr^k\Split(A,P)$ is a flat $A$-module,
so (i) holds. (ii) is immediate. To show (iii), let
us introduce the functor $T:A\Alg\to\Set$ that
assigns to every $A$-algebra $B$ the set of all
pairs $(x,\phi)$  where $x\in(P\otimes_AB)_*$ and 
$\phi:P\otimes_AB\to B$. So $S$ is a subfunctor of
$T$. 
\begin{claim}\label{cl_T.functor} 
The functor $T$ is represented by the $A$-algebra 
$\Sym_A^\bullet(P\oplus P^*)$.
\end{claim}
\begin{pfclaim} Indeed, there are natural bijections:
$$(P\otimes_AB)_*\stackrel{\sim}{\longrightarrow}
\Hom_B(P^*\otimes_AB,B)\stackrel{\sim}{\longrightarrow}
\Hom_A(P^*,B)\stackrel{\sim}{\longrightarrow}
\Hom_{A\Alg}(\Sym_A^\bullet P^*,B)$$
which show that the functor $B\mapsto(P\otimes_AB)_*$
is represented by the $A$-algebra $\Sym_AP^*$.
Working out the definitions, one finds that the
composition of these bijections assigns to an element
$x\in(P\otimes_AB)_*$ the unique $A$-algebra morphism
$f_x:\Sym_A^\bullet P^*\to B$ such that $f_x(\psi)=\psi(x)$ 
for every $\psi\in\Hom_A(P,A)$. 
Similarly, the functor $B\mapsto\Hom_B(P\otimes_AB,B)$
is represented by $\Sym_A^\bullet P$, and again, one checks
that the bijection assigns to $\phi:P\otimes_AB\to B$
the unique $A$-algebra morphism 
$g_\phi:\Sym^\bullet_AP\to B$ such that $g_\phi(p)=\phi(p)$
for every $p\in P_*$. It follows that $T$ is represented
by $(\Sym_A^\bullet P)\otimes_A(\Sym_A^\bullet P^*)\simeq
\Sym_A^\bullet(P\oplus P^*)$.
\end{pfclaim}

For every $A$-algebra $B$ we have a natural map:
$\alpha_B:T(B)\to B_*$ given by $(x,\phi)\mapsto\phi(x)$.
This defines a natural transformation of functors
$\alpha:T\to (-)_*$. Moreover, let us consider the
trivial map $\beta_B:T(B)\to B_*$ that sends everything 
onto the element $1\in B_*$. $\beta$ is another
natural transformation from $T$ to the almost elements
functor. Clearly :
$$
S(B)=\mathrm{Equal}(
\xymatrix{T(B) \ar@<.5ex>[r]^-{\alpha_B} 
\ar@<-.5ex>[r]_-{\beta_B} & B_*}).
$$
We remark that the functor $B\mapsto B_*$ on $A$-algebras,
is represented by $A[X]:=\Sym_A^\bullet A$. Therefore
there are morphisms $\alpha^*,\beta^*:A[X]\to
\Sym_A^\bullet(P\oplus P^*)$ that represent these
natural transformations. It follows that $S$ is represented
by the $A$-algebra 
$$
\mathrm{Coequal}(
\xymatrix{A[X] \ar@<.5ex>[r]^-{\alpha^*} 
\ar@<-.5ex>[r]_-{\beta^*} & 
\Sym_A^\bullet(P\oplus P^*)}).
$$
To determine $\alpha^*$ and $\beta^*$ it suffices
to calculate them on the element $X\in A[X]_*$.
It is easy to see that $\alpha^*(X)=1$. To conclude 
the proof it suffices therefore to show:
\begin{claim}\label{cl_beta.is.zeta} 
$\beta^*(X)=\zeta_P$.
\end{claim}
\begin{pfclaim} In view of the definitions, and using
the notation of the proof of claim \ref{cl_T.functor},
the claim amounts to the identity: 
$(g_\phi\otimes f_x)(\zeta_P)=\phi(x)$ for every
$(x,\phi)\in T(B)$. By naturality, it suffices to
show this for $B=A$. Now, for every $\eps\in\fm$,
we can write $\eps\cdot\zeta_P=\sum_iq_i\otimes\psi_i$ 
for some $q_i\in P$, $\psi_i\in P^*$ and we have
$\sum_i q_i\cdot\psi_i(b)=\eps\cdot b$
for all $b\in P_*$. Hence 
$(g_\phi\otimes f_x)(\sum_i q_i\otimes\psi_i)=
\sum_i\phi(q_i)\cdot\psi_i(x)=
\phi(\sum_i q_i\cdot\psi_i(x))=\phi(\eps x)$ and the 
claim follows.
\end{pfclaim}
\end{proof}

\begin{remark} The construction of the splitting algebra
occurs already, in a tannakian context, in Deligne's paper
\cite{De} : see the proof of lemma 7.15 in {\em loc.cit.}
\end{remark}

\sset\subsubsection{}\label{subsec_divlittle}
We recall that for every $k\geq 0$ there are natural 
morphisms
$$\Sym^k_AP\stackrel{\alpha_P}{\longrightarrow}
\Gamma^k_AP\stackrel{\beta_P}{\longrightarrow}\Sym^k_AP$$
such that $\beta_P\circ\alpha_P=k!\cdot\one_{\Sym^k_AP}$
and $\alpha_P\circ\beta_P=k!\cdot\one_{\Gamma^k_AP}$.
(to obtain the morphisms, one can consider the flat
$A_*$-module $P_!$, thus one can assume that $P$ is
a module over a usual ring; then $\alpha_P$ is obtained 
by extending multiplicatively the identity morphism 
$\Sym^1_AP=P\to P=\Gamma^1_AP$, and $\beta_P$ is deduced 
from the homogeneous degree $k$ polynomial law 
$P\otimes_AB\to(P\otimes_AB)^{\otimes k}$ defined by
$x\mapsto x^{\otimes k}$).
Moreover $(\Gamma^k_AP)^*\simeq\Sym^k_AP^*$.
\begin{lemma}\label{lem_Casimir} 
With the notation of \eqref{subsec_divlittle} we have:
$(\alpha_P\otimes_A\one_{\Sym^k_AP^*})
(\zeta^k_P)=k!\cdot\zeta_{\Gamma^k_AP}$.
\end{lemma}
\begin{proof} Suppose first that $P$ is a free 
$A$-module, let $e_1,...,e_n$ be a base of $P_*$
and $e^*_1,...,e^*_n$ the dual base of $P^*$. 
Then $\Gamma^k_AP$ is the free $A$-module generated 
by the basis $e_1^{[n_1]}\cdot...\cdot e^{[n_k]}_k$
where $0\leq n_i\leq k$ for $i=1,...,k$ and 
$\sum_j n_j=k$. The dual of this basis is the
basis of $\Sym^k_AP^*$ consisting of the elements
$e^{*n_1}_1\cdot...\cdot e^{*n_k}_k$.
Furthermore, $\zeta_P=\sum_ie_i\otimes e^*_i$
and therefore $\zeta^k_P=\sum_{\underline n}
\binom{k}{\underline{n}}
(e_1^{n_1}\cdot...\cdot e^{n_k}_k)\otimes
(e^{*n_1}_1\cdot...\cdot e^{*n_k}_k)$, where 
$\underline n:=(n_1,...,n_k)$ ranges over the 
multi-indices submitted to the above conditions and 
$\binom{k}{\underline{n}}:=
\frac{k!}{n_1!\cdot...\cdot n_k!}$.
Then the claim follows straightforwardly from the
identity: $\alpha_P(e_1^{n_1}\cdot...\cdot e^{n_k}_k)=
n_1!\cdot...\cdot n_k!\cdot
e_1^{[n_1]}\cdot...\cdot e^{[n_k]}_k$.
For the general case we shall use the following
\begin{claim}\label{claim_simple.comput}
Let $M$ be an almost finitely generated projective 
$A$-module and pick, for a given $\eps\in\fm$, morphisms 
$u:M\to F$ and $v:F\to M$ with $v\circ u=\eps\cdot\one_M$. 
Then $v\otimes u^*(\zeta_F)=\eps\cdot\zeta_M$.
\end{claim}
\begin{pfclaim} We have a commutative diagram
$$\xymatrix{
F\otimes_AF^* \ar[rr]^-{v\otimes u^*} 
\ar[d] && M\otimes_AM^* \ar[d] \\
\End_A(F) \ar[rr] && \End_A(M)
}$$
where the vertical morphisms are the natural ones, 
and where the bottom morphism is given by 
$\phi\mapsto v\circ\phi\circ u$. Then the claim follows
by an easy diagram chase.
\end{pfclaim}

Pick morphisms $u:P\to F$ and $v:F\to P$ with 
$v\circ u=\eps\cdot\one_P$. We consider the commutative 
diagram
$$\xymatrix{
\Sym^k_AF\otimes_A\Sym^k_AF^* 
\ar[rrr]^-{\Sym_A^kv\otimes\Sym_A^ku^*}
\ar[d]_{\alpha_F\otimes\one_{\Sym^k_AF^*}}
& & &\Sym^k_AP\otimes_A\Sym^k_AP^* 
\ar[d]^{\alpha_P\otimes\one_{\Sym^k_AP^*}}
\\
\Gamma^k_AF\otimes_A\Sym^k_AF^*
\ar[rrr]^-{\Gamma_A^kv\otimes\Sym_A^ku^*} & & &
\Gamma^k_AP\otimes_A\Sym^k_AP^*.
}$$
By claim \ref{claim_simple.comput}, we have 
$v\otimes u^*(\zeta_F)=\eps\cdot\zeta_P$; whence
 $\Sym_A^kv\otimes\Sym_A^ku^*(\zeta^k_F)=
\eps^k\cdot\zeta^k_P$. Moreover, we remark that
$(\Gamma^k_Av)\circ(\Sym_A^ku^*)^*=
(\Gamma^k_Av)\circ(\Gamma^k_Au)=
\eps^k\cdot\one_{\Gamma_A^kP}$,
therefore claim \ref{claim_simple.comput} (applied
for $M=\Gamma^k_AP$) yields
$\Gamma_A^kv\otimes\Sym_A^ku^*(\zeta_{\Gamma_A^kF})=
\eps^k\cdot\zeta_{\Gamma_A^kP}$. Since we already know
the lemma for $F$, a simple diagram chase shows that
$(\alpha\otimes_A\one_{\Sym^k_AP^*})(\eps^k\cdot\zeta^k_P)=
\eps^k\cdot k!\cdot\zeta_{\Gamma_A^kP}$. Since the 
$k$-powers of elements of $\fm$ generate $\fm$, the
claim follows.
\end{proof}
\begin{lemma}\label{lem_zeta.nilpotent} 
Let $P$ be as above and suppose that  $\zeta_P$ is 
nilpotent in $\Sym^\bullet_A(P\oplus P^*)$.
\begin{enumerate}
\item
If\/ $\Q\subset A_*$, then $\chi_P(X)=(1+X)^{-\alpha}$ 
for some continuous function $\alpha:\Spec\,A_*\to\N$. 
\item
If\/ $\F_p\subset A_*$, then $\cE_{P/A}$ is a 
Frobenius-nilpotent ideal.
\end{enumerate}
\end{lemma}
\begin{proof} 
(i): Since $\Q\subset A_*$, 
then $k!$ is invertible in $A_*$ for every $k\geq 0$; 
by lemma \ref{lem_Casimir} it follows that $\Sym^k_AP=0$, 
thus $\psi_P(X)$ is a polynomial. By corollary 
\ref{cor_psi.Sym},
$\psi_P(-X)$ defines an endomorphism of $\G_{m,A_*}$,
therefore $\psi_P(X)=(1+X)^\alpha$ for some continuous
function $\alpha:\Spec\,A_*\to\N$; then the claim
follows by corollary \ref{cor_psi.Sym}(ii).

(ii): Let $(f,p)\in(P^*\oplus P)_*$; in the notation
of the proof of lemma \ref{lem_Split.generals}, we 
can write $(f,p)\in T(A)$. It follows that $(f,p)$
corresponds to a morphism of $A$-algebras
$\Sym_A^\bullet(f,p)^*:\Sym^\bullet_A(P\oplus P^*)\to A$. 
In particular
$\Sym_A^\bullet(f,p)^*(\zeta^k_P)=
\Sym_A^\bullet(f,p)^*(\zeta_P)^k$ for every $k\geq 0$.
By inspecting the proof of claim \ref{cl_beta.is.zeta}, 
we deduce $\Sym^\bullet_A(f,p)^*(\zeta_P)=f(p)$ for
every $f\in P^*$ and $p\in P$. By hypothesis,
$\zeta^{p^n}_P=0$ for every sufficiently large $n$.
It follows that $\cE_{P/A}$ is Frobenius nilpotent.
\end{proof}

\begin{lemma}\label{lem_reduce.to.fields}
Let $R_0$ be a noetherian commutative ring,
$R$ an $R_0$-algebra and $M$ a flat $R$-module. Then $M=0$
if and only if $M\otimes_{R_0}\kappa=0$ for every residue
field $\kappa$ of $R_0$ ({\em i.e.}, for every field $\kappa$
of the form ${\rm Frac}(R_0/\fp)$, where $\fp$ is some
prime ideal of $R_0$).
\end{lemma}
\begin{proof} Clearly we have only to show the direction
$\Leftarrow$. It suffices to show that $M_\fp=0$ for every 
prime ideal of $R_0$. Hence we can assume that $R_0$ is
local, in particular of finite Krull dimension. We proceed
by induction on the dimension of $R_0$. If $\dim R_0=0$,
then $R_0$ is a local artinian ring, hence a power of
its maximal ideal $\fm$ is equal to $0$. By assumption,
$M/\fm M=0$, {\em i.e.\/} $M=\fm M$. Then
$M=\fm^kM$ for every $k\geq 0$, so $M=0$.
Next, suppose that $\dim R_0=d$ and the lemma already
known for all rings of dimension strictly less than
$d$. Assume first that $R_0$ is an integral domain and
pick $f\in S:=R_0\setminus\{0\}$. Then $R_0/fR_0$ has
dimension strictly less than $d$, so by induction we have
$M/fM=0$, {\em i.e.\/} $M=fM$. Due to the 
flatness of $M$, we have: $\Ann_M(f)=\Ann_R(f)\cdot M=
\Ann_R(f)\cdot fM=0$. This implies that the kernel
of the natural map $M\to S^{-1}M$ is trivial. On the
other hand, by hypothesis $S^{-1}M=0$, whence $M=0$ in
this case. For a general $R_0$ of dimension $d$, notice
that the above argument implies that, for every minimal
prime ideal $\fp$ of $R_0$, we have $\fp M=M$.
But the product of all (finitely many) minimal prime
ideals is contained in the nilpotent radical ${\mathfrak R}$
of $R_0$, whence ${\mathfrak R}M=M$, and finally $M=0$
as claimed.
\end{proof}

\sset\subsubsection{}\label{subsec_formal-rk.at.p}
\index{Almost module(s)!almost projective!$\frk_A(P,\fp)$ : 
localization at a prime ideal of the formal rank of 
an|indref{subsec_formal-rk.at.p}}
Let now $\fp\in\Spec\,A_*$. By composing $\frk_A(P)$ with 
the map $A_*\to A_*/\fp$, we obtain a ring homomorphism
$$\frk_A(P,\fp):\cG_\infty\to A_*/\fp.$$
In case $(A_*/\fp)^a\neq 0$, we can interpret 
$\frk_A(P,\fp)$ as the formal rank of $P\otimes_A(A/\fp^a)$.
More precisely, let $\pi:A\to A/\fp^a$ be the natural
projection; then $\pi_*$ factors through a map 
$\pi':A_*/\fp\to(A/\fp^a)_*$ and we have:
$\frk_{A/\fp^a}(P\otimes_A(A/\fp^a))=
\pi'\circ\frk_A(P,\fp)$.

Even if $(A_*/\fp)^a=0$, the morphism $\frk_A(P,\fp)$
can still be interpreted as the map associated to an
endomorphism of $\hat\G_{m,A_*/\fp}$, so it still makes
sense to ask whether $\frk_A(P,\fp)$ is an integer,
as indicated in remark \ref{rem_continuous.fct}(ii).

\begin{lemma}\label{cl_always.same.rk} 
Let $P$ be an almost finitely generated projective 
$A$-module and $\fp\in\Spec\,A_*$. If $B$ is any 
$A_{*\fp}^a$-algebra and $\fq\in\Spec\,B_*$ such that 
$r(\fq):=\frk_B(P\otimes_AB,\fq)$ is an integer, 
then $r(\fp):=\frk_A(P,\fp)$ is also an integer 
and $r(\fp)=r(\fq)$.
\end{lemma}
\begin{proof} Indeed, let us consider the natural
maps $\cG_\infty\to A_*\to B_*$; under the assumptions,
the contraction of $\fq$ in $A_*$ is contained in $\fp$.
Since the image of $\binom{\alpha}{i}$ in $B_*/\fq$ is
$\binom{r(\fq)}{i}$, it follows that the same holds in 
$A_*/\fp$.
\end{proof}
\begin{definition}\label{def_admits_splittings}
\index{Almost module(s)!that admits infinite splittings
|indref{def_admits_splittings}}
We say that an $A$-module $P$ {\em admits infinite splittings\/}
if there is an infinite chains of decompositions of the form:
$P\simeq A\oplus P_1$, $P_1\simeq A\oplus P_2$, 
$P_2\simeq A\oplus P_3$, ...
\end{definition}

\begin{theorem}\label{th_charact.almost.fin.rk}
Let $P$ be an almost finitely generated projective
$A$-module. The following conditions are equivalent:
\begin{enumerate}
\item
$P$ is of almost finite rank.
\item
For all $A$-algebras $B\neq 0$, we have: 
$\bigcap_{r>0}\cE_{\Lambda^r_B(P\otimes_AB)/B}=0$.
\item
For all $A$-algebras $B\neq 0$, $P_B:=P\otimes_AB$ 
does not admit infinite splittings, and moreover if 
$P_B\simeq B^n\oplus Q$ for some $B$-module $Q$ and 
$\chi_Q(X)=(1+X)^{-\alpha}$ for some continuous function 
$\alpha:\Spec\,B_*\to\N$, then $Q=0$.
\item
For all $A$-algebras $B\neq 0$, $P_B$ does not admit 
infinite splittings, and moreover if $P_B\simeq B^n\oplus Q$ 
for some $B$-module $Q$, then :
\begin{enumerate}
\item
If\/ $\F_p\subset B_*$ and $Q=IQ$ for
a Frobenius-nilpotent ideal $I\subset B$, then $Q=0$;
\item
If\/ $\Q\subset B_*$ and $\Sym^r_BQ=0$ for some 
$r\geq 1$, then $Q=0$.
\end{enumerate}
\end{enumerate}
\end{theorem}
\begin{proof} (i) $\Rightarrow$ (ii) : indeed, from proposition
\ref{prop_eval.ideal}(iii) one sees that, for every $A$-module
of almost finite rank, every $A$-algebra $B$ and every 
$\eps\in\fm$, there exists $r\geq 0$ such that 
$\eps\cdot\cE_{\Lambda^r_B(P\otimes_AB)/B}=0$.

(ii) $\Rightarrow$ (iii) : let $B\neq 0$ be an $A$-algebra; 
by hypothesis, there exists $r\geq 0$ such that 
$J_r:=\cE_{\Lambda^r_B(P\otimes_AB)/B}\neq B$.
Suppose that $P\otimes_AB$ admits infinite splittings.
The $B/J_r$-module $P_B/J_rP_B$ has rank $<r$,
and at the same time it admits infinite splittings, a
contradiction.

Suppose next, that there is a decomposition
$P\otimes_AB\simeq B^n\oplus Q$; then $Q$ is obviously
of almost finite rank. Suppose that $\chi_Q(X)=(1+X)^\alpha$ 
has the shape described in (ii). We reduce easily to the
case where $\alpha$ is a constant function.
However, $\chi_{Q/J_rQ}(X)$ is a polynomial of
degree $<r$, thus $\alpha=0$, and then $Q=0$ by theorem
\ref{th_structure.alm.fin.rk}(ii). 

(iii) $\Rightarrow$ (iv) : suppose that 
$\F_p\subset B_*$ and $Q=IQ$ for some 
Frobenius-nilpotent ideal $I$. Then 
$\chi_Q(X)\in 1+I_*[[X]]$, which means that
the image of $\chi_Q(X)$ in $\End(\hat\G_{m,B_*/I_*})$
is the trivial endomorphism. But we have a commutative
diagram 
$$\xymatrix{\cC^0(\Spec\,B_*,\Z_p) \ar[r] \ar[d] &
\End(\hat\G_{m,B_*}) \ar[d] \\
\cC^0(\Spec\,B_*/I_*,\Z_p) \ar[r] &
\End(\hat\G_{m,B_*/I_*})
}$$
where the horizontal maps are those defined in remark
\ref{rem_continuous.fct}(ii), and are bijective by corollary
\ref{cor_scheme.descrip}. The left vertical map is induced 
by restriction to the closed subset $\Spec\,B_*/I_*$, and 
since $I$ is Frobenius-nilpotent, it is a bijection as well.
It follows that the right vertical map is bijective,
whence $\frk_B(Q)=0$, and finally $Q=0$ by (iii).

Next, consider the case when $\Q\subset B_*$ and 
$\Sym^r_BQ=0$ for some $r\geq 1$. It follows that 
$\zeta_Q$ is nilpotent in $\Sym_B^\bullet(Q\oplus Q^*)$.
By lemma \ref{lem_zeta.nilpotent}(i), 
$\chi_Q(X)=(1+X)^{-\alpha}$ for some continuous
$\alpha:\Spec\,B_*\to\N$. Then (iii) implies that $Q=0$.  

To show that (iv) $\Rightarrow$ (i),
we will use the following:
\begin{claim}\label{cl_assume.iv} 
Assume (iv). Then:~ $\Split(B,Q)=0\Rightarrow Q=0$.
\end{claim}
\begin{pfclaim} Suppose $Q\neq 0$ and $\Split(B,Q)=0$; 
then \eqref{eq_gr.of.splitter} implies that for every
$\eps\in\fm$ there exists $j\geq 0$ such that 
$\eps\cdot\zeta_Q^j=0$. We have $\eps Q\neq 0$ for 
some $\eps\in\fm$. From the flatness of $Q$, we derive 
$\Ann_{\Sym^\bullet_B(Q\oplus Q^*)}(\eps)=
\Ann_B(\eps)\cdot\Sym^\bullet_B(Q\oplus Q^*)$, hence 
we can replace $B$ by $B/\Ann_B(\eps)$, $Q$ by 
$Q/\Ann_B(\eps)\cdot Q$, thereby achieving that 
$\zeta_Q$ is nilpotent in $\Sym^\bullet_B(Q\oplus Q^*)$
and still $Q\neq 0$. Using lemma \ref{lem_reduce.to.fields}
(and the functoriality of $\Split(B,Q)$ for base extensions
$B\to B'$) we can further assume that $B_*$ contains
either $\Q$ or one of the finite fields $\F_p$.
If $\Q\subset B_*$, then $k!$ is invertible in 
$B_*$ for every $k\geq 0$; by lemma \ref{lem_Casimir}
it follows that $\Sym^k_BQ=0$, whence $Q=0$ by (iv),
a contradiction. If $\F_p\subset B_*$, then by lemma 
\ref{lem_zeta.nilpotent}(ii), $\cE_{Q/B}$ is 
Frobenius-nilpotent. However, from proposition 
\ref{prop_eval.ideal}(i) it follows easily that 
$Q=\cE_{Q/B}\cdot Q$, whence $Q=0$, again by (iv),
and again a contradiction. In either case, this shows 
that $\Split(B,Q)\neq 0$, as claimed.
\end{pfclaim}

\begin{claim}\label{cl_A.over.I_P}
Assume (iv). Then $A/\cE_{P/A}$ is a flat $A$-algebra.
\end{claim}
\begin{pfclaim} It suffices to show that 
$A_\fp/(\cE_{P/A})_\fp$ is a flat $A_\fp$-algebra
for every prime ideal $\fp\subset A_*$.
If $\cE_{P/A*}$ is not contained in $\fp$, then 
$(\cE_{P/A})_\fp=A_\fp$, so there is nothing to 
prove in this case. We assume therefore that 
\set\begin{equation}\label{eq_assume}
\cE_{P/A*}\subset\fp.
\end{equation}
We will show that $P_\fp=0$ in such case, whence 
$A_\fp/(\cE_{P/A})_\fp=A_\fp$, so the
claim will follow. From \eqref{eq_assume} and corollary 
\ref{cor_psi.Sym}(iii) we know already that
\set\begin{equation}\label{eq_residual.rk}
\frk_A(P,\fp)=0.
\end{equation}
Suppose that $P_\fp\neq 0$; then there exists $\eps\in\fm$ 
such that $\eps P_\fp\neq 0$. Define inductively 
$A_0:=A_\fp$, $Q_0:=P_\fp$, $A_{i+1}:=\Split(A_i,Q_i)$ and 
$Q_{i+1}$ as an $A_{i+1}$-module such that 
$Q_i\otimes_{A_i}A_{i+1}\simeq A_{i+1}\oplus Q_{i+1}$, 
for every $i\geq 0$ (the existence of $Q_{i+1}$ is
assured by lemma \ref{lem_Split.generals}(iii)). Then 
$\colim{n\in\N}A_n\simeq 0$, since, after base change 
to this $V^a$-algebra, $P$ admits infinite splittings. 
This implies that there exists $n\in\N$ such that 
$\eps A_{n+1}=0$ and $\eps A_n\neq 0$. However, 
since $A_{n+1}$ is flat over $A_n$, we have: 
$A_{n+1}=\Ann_{A_{n+1}}(\eps)=\Ann_{A_n}(\eps)\cdot A_{n+1}$. 
Set $A':=A_n/\Ann_{A_n}(\eps)$; then
$\Split(A',Q_n\otimes_{A_n}A')=0$, so 
$Q_n\otimes_{A_n}A'=0$ by claim \ref{cl_assume.iv}. 
By flatness of $Q_n$, this means that $\eps Q_n=0$;
in particular, $n>0$. By definition, 
$Q_{n-1}\otimes_{A_{n-1}}A_n\simeq A_n\oplus Q_n$;
it follows that 
$Q':=Q_{n-1}\otimes_{A_{n-1}}A'\simeq A'$, in 
particular $\frk_{A'}(Q')=1$ and consequently
$\frk_{A'}(P\otimes_AA')=n>0$; in view of lemma
\ref{cl_always.same.rk}, this contradicts 
\eqref{eq_residual.rk}, therefore $P_\fp=0$, as required.
\end{pfclaim}

\begin{claim}\label{cl_faith.flat}
Assuming (iv), the natural morphism 
$\phi:A\to(A/\cE_{P/A})\times \Split(A,P)$ is faithfully flat.
\end{claim}
\begin{pfclaim} The flatness is clear from claim
\ref{cl_A.over.I_P}. Hence, to prove the claim, it
suffices to show that
$((A/\cE_{P/A})\times\Split(A,P))\otimes_A(A/J)\neq 0$
for every proper ideal $J\subset A$,.
But the construction of $\phi$ commutes with arbitrary
base changes $A\to A'$, therefore we are reduced to
verify that $(A/\cE_{P/A})\times\Split(A,P)\neq 0$ when
$A\neq 0$.
By claim \ref{cl_assume.iv}, this can fail only
when $P=0$; but in this case $\cE_{P/A}=0$, so
the claim follows.
\end{pfclaim}

We can now conclude the proof of the theorem: define 
inductively as in the proof of claim \ref{cl_A.over.I_P}:
$A_0:=A$, $Q_0:=P$, $A_{i+1}:=\Split(A_i,Q_i)$ 
and $Q_{i+1}$ as an $A_{i+1}$-module such that 
$Q_i\otimes_{A_i}A_{i+1}=A_{i+1}\oplus Q_{i+1}$, for
every $i\geq 0$. The same argument as in {\em loc.~cit.} 
shows that, for every $\eps\in\fm$, there exists $n\in\N$ 
such that $\eps A_n=0$. We may assume that 
$\eps A_{n-1}\neq 0$. Moreover, by claim 
\ref{cl_faith.flat} (and an easy induction),
$B:=A_0/\cE_{Q_0/A_0}\times A_1/\cE_{Q_1/A_1}\times...
\times A_{n-1}/\cE_{Q_{n-1}/A_{n-1}}\times A_n$ is 
a faithfully flat $A$-algebra. However, one checks 
easily by induction that $P\otimes_A(A_i/\cE_{Q_i/A_i})$
is a free $A_i/\cE_{Q_i/A_i}$-module of rank $i$, for 
every $i<n$. Hence, $\Lambda^n_B(P\otimes_AB)\simeq
\Lambda^n_{A_n}(P\otimes_AA_n)$, which is therefore
killed by $\eps$. By faithful flatness, so is 
$\Lambda^n_AP$. The proof is concluded.
\end{proof}

\begin{proposition}\label{prop_finally.faith.flat}
If $P$ is a faithfully flat almost projective $A$-module
of almost finite rank, then $\Split(A,P)$ is 
faithfully flat over $A$. 
\end{proposition}
\begin{proof}
If $A=0$ there is nothing to prove, so we assume that
$A\neq 0$. In this case, it suffices to show that 
$\Split(A,P)\otimes_AA/I\neq 0$ for every proper ideal 
$I$ of $A$. However, 
$\Split(A,P)\otimes_AA/I\simeq\Split(A/I,P/IP)$, 
and since $P$ is faithfully flat, $P/IP\neq 0$; 
hence we are reduced to showing that $\Split(A,P)\neq 0$ 
when $P$ is faithfully flat. Suppose that $\Split(A,P)=0$;
then \eqref{eq_gr.of.splitter} implies that for every
$\eps\in\fm$ there exists $j\geq 0$ such that 
$\eps\cdot\zeta_P^j=0$. Since $P\neq 0$, we have 
$\eps P\neq 0$ for some $\eps\in\fm$. From
the flatness of $P$, we derive 
$\Ann_{\Sym^\bullet_A(P\oplus P^*)}(\eps)=
\Ann_A(\eps)\cdot\Sym^\bullet_A(P\oplus P^*)$, hence 
we can replace $A$ by $A/\Ann_A(\eps)$, $P$ by 
$P/\Ann_A(\eps)\cdot P$, which allows us to assume 
that $\zeta_P$ is nilpotent in $\Split(A,P)$.
Using lemmata \ref{lem_reduce.to.fields} and 
\ref{lem_Split.generals}(ii), we can further assume that 
$A_*$ contains either $\Q$ or one of the finite fields 
$\F_p$. If $\Q\subset A_*$, then $k!$ is invertible in 
$A_*$ for every $k\geq 0$; by lemma \ref{lem_Casimir}
it follows that $\Sym^k_AP=0$, whence $P=0$ by
theorem \ref{th_charact.almost.fin.rk}(iv),
which contradicts our assumptions, so the proposition
is proved in this case. Finally, suppose that
$\F_p\subset A_*$, then by lemma 
\ref{lem_zeta.nilpotent}(ii), $\cE_{P/A}$ is 
Frobenius-nilpotent. However, since $P$ is faithfully
flat, proposition \ref{prop_eval.ideal}(iv) says
that $\cE_{P/A}=A$, so $A=0$, which again contradicts
our assumptions.
\end{proof}

\sset\subsubsection{}\label{subsec_fpqc}
For any $V^a$-algebra $A$ we have a (large) fpqc site 
on the category $(A\Alg)^o$ (in some fixed universe!); 
as usual, this site is defined by the pretopology whose 
covering families are the finite families 
$\{\Spec\,C_i\to\Spec\,B~|~i=1,...,n\}$ such that the induced 
morphism $B\to C_1\times...\times C_n$ is faithfully flat
(notation of \eqref{subsec_general-schemes}).

\begin{theorem}\label{th_free.in.fpqc} Every almost projective 
$A$-module of finite rank is locally free of finite rank in the 
fpqc topology of $(A\Alg)^o$.
\end{theorem}
\begin{proof} We iterate the construction of $\Split(A,P)$
to split off successive free submodules of rank one.
We use the previous characterization of modules of finite 
rank (proposition \ref{prop_decomp.fin.rank}) to show that 
this procedure stops after finitely many iterations. By
proposition \ref{prop_finally.faith.flat}, the output of 
this procedure is a faithfully flat $A$-algebra. 
\end{proof}

Theorem \ref{th_free.in.fpqc} allows to prove easily results
on almost projective modules of finite rank, by reduction
to the case of free modules. Here are a few examples of this method.

\begin{lemma}\label{lem_reduce.to.free} 
Let $P$ be an almost projective $A$-module of constant 
rank equal to $r\in\N$. Then, for every integer $0\leq k\leq r$, 
the natural morphism
\set\begin{equation}\label{eq_wedging}
\Lambda_A^kP\otimes_A\Lambda_A^{r-k}P\to\Lambda^r_AP
\qquad x\otimes y\mapsto x\wedge y
\end{equation}
is a perfect pairing.
\end{lemma}
\begin{proof} By theorem \ref{th_free.in.fpqc}, there exists
a faithfully flat $A$-algebra $B$ such that $P_B:=P\otimes_AB$ 
is a free $B$-module of rank $r$. It suffices to prove the 
assertion for the $B$-module $P_B$, in which case the claim
is well known.
\end{proof}

\sset\subsubsection{}\label{subsec_Cramer}
Keep the assumptions of lemma \ref{lem_reduce.to.free}. Taking
$k=1$ in \eqref{eq_wedging}, we derive a natural isomorphism
$$
\beta_P:(\Lambda^{r-1}_AP)^*\stackrel{\sim}{\to} 
P\otimes_A(\Lambda^r_AP)^*.
$$
Now, let us consider an $A$-linear morphism $\phi:P\to Q$
of $A$-modules of constant rank equal to $r$. We set
$$
\psi:=\beta_P\circ(\Lambda^{r-1}_A\phi)^*\circ\beta_Q^{-1}:
Q\otimes_A(\Lambda_A^rQ)^*\to P\otimes_A(\Lambda_A^rP)^*.
$$

\begin{proposition}\label{prop_Cramer} 
With the notation of \eqref{subsec_Cramer}, we have:
$$
\psi\circ(\phi\otimes_A\one_{(\Lambda^r_AQ)^*})=
\one_P\otimes_A(\Lambda^r_A\phi)^*\quad\text{and}\quad
(\phi\otimes_A\one_{(\Lambda^r_AP)^*})\circ\psi=
\one_Q\otimes_A(\Lambda^r_A\phi)^*.
$$
Especially, $\phi$ is an isomorphism if and only if the
same holds for $\Lambda^r_A\phi$. 
\end{proposition}
\begin{proof} After faithfully flat base change, we can
assume that $P$ and $Q$ are free modules of rank $r$. Then
we recognize Cramer's rule in the above identities.
\end{proof}

\sset\subsubsection{}\label{subsec_Cayley}
Keep the assumption of lemma \ref{lem_reduce.to.free}
and let $\phi$ be an $A$-linear endomorphism of $P$.
Set $B:=A[X,X^{-1}]$ and $P_B:=P\otimes_AB$. Obviously $\phi$
induces a $\Lambda$-nilpotent endomorphism of $P_B$,
which we denote by the same letter. Hence we can define 
$$
\chi_\phi(X):=
X^r\cdot\det(\one_{P_B}-X^{-1}\cdot\phi)\in A_*[X]
$$
(notation of \eqref{subsec_lambda_nilpotent}).
\begin{proposition}\label{prop_Cayley-Ham} 
With the notation of \eqref{subsec_Cayley},
we have $\chi_\phi(\phi)=0$ in $\End_A(P)$.
\end{proposition}
\begin{proof} Again, we can reduce to the case of
a free $A$-module of constant finite rank, in which
case we conclude by Cayley-Hamilton.
\end{proof}

\begin{corollary}\label{cor_Cayley} Keep the assumptions
of \eqref{subsec_Cayley}, and suppose that $\phi$ is integral
over a subring $S\subset A_*$. Then the coefficients of
$\chi_\phi$ are integral over $S$.
\end{corollary}
\begin{proof} The assumption means that $\phi$ satisfies an
equation of the kind $\phi^n+\sum_{i=0}^{n-1}a_i\phi^i=0$, where $a_i\in S$
for every $i=0,...,n-1$. We can assume that $P$ is free, in which
case we reduce to the case of an endomorphism of a free $R$-module of
finite rank, where $R$ is a usual commutative ring containing $S$;
we can further suppose that $R$ is of finite type over $\Z$,
and it is easily seen that we can even replace $R$ by its associated
reduced ring $R/\nil(R)$. Then $R$ injects into a finite product
of fields $\prod_iK_i$, and we can replace $S$ by its integral
closure in $\prod_iK_i$, which allows us to reduce to the case
where $R$ is a field. In this case the coefficients of $\chi_\phi$
are elementary symmetric polynomials in the eigenvalues of $\phi$,
so it suffices to show that these eigenvalues $e_1,...,e_r$ are
integral over $S$. But this is clear, since we have more precisely
$e_j^n+\sum^{n-1}_{i=0}a_ie^i=0$ for every $j\leq r$.
\end{proof}

To conclude this section, we want to apply the previous
results to analyze in some detail the structure of invertible 
modules : it turns out that the notion of invertibility is  
rather more subtle than for usual modules over rings.

\begin{definition}\label{def_strictly_invertible}
\index{Almost module(s)!strictly invertible
|indref{def_strictly_invertible}}
\index{$u_M$|indref{def_strictly_invertible}}
Let $M$ be an invertible $A$-module.
Clearly $M\otimes_AM$ is invertible as well, consequently
the map $A\to\End_A(M\otimes_AM)\to A~:
~a\mapsto a\cdot\one_{M\otimes_AM}$ is an isomorphism 
(by the proof of lemma \ref{lem_explain.invert}(iii)). 
Especially, for the transposition endomorphism $\theta_{M|M}$ 
of $M\otimes_AM~:~x\otimes y\mapsto y\otimes x$, there 
exists a unique element $u_M\in A_*$ such that 
$\theta_{M|M}=u_M\cdot\one_{M\otimes_AM}$. Clearly $u^2_M=1$. 
We say that $M$ is {\em strictly invertible\/} if $u_M=1$.
\end{definition}

\begin{lemma}\label{lem_charact.strictly.inv} 
For an invertible $A$-module the following are equivalent:
\begin{enumerate}
\item
$M$ is strictly invertible;
\item
$\Lambda^2_AM=0$;
\item
$M$ is of almost finite rank;
\item
there exists a faithfully flat $A$-algebra $B$
such that $M\otimes_AB\simeq B$.
\end{enumerate}
\end{lemma}
\begin{proof} (i) $\Rightarrow$ (ii): indeed, the
condition $u_M=1$ says that the antisymmetrizer 
operator $a_2:M^{\otimes 2}\to M^{\otimes 2}$ 
vanishes (cp. the proof of lemma 
\ref{lem_lambda.of.lambda}); since $M$ is flat, (ii)
follows. 

(ii) $\Rightarrow$ (iii) and (iv) $\Rightarrow$ (i)
are obvious. To show that (iii) $\Rightarrow$ (iv)
let us set $B:=\Split(A,M)$; by proposition 
\ref{prop_finally.faith.flat} $B$ is faithfully flat
over $A$, and $B\otimes_AM\simeq B\oplus X$ for
some $B$-module $X$. Clearly $B\otimes_AM$ is an
invertible $B$-module, therefore, by lemma 
\ref{lem_explain.invert}(ii), the evaluation
morphism gives an isomorphism 
$(B\oplus X)\otimes_A(B\oplus X^*)\simeq 
B\oplus X\oplus X^*\oplus(X\otimes_AX^*)\simeq B$. 
By inspection, the restriction of the latter morphism
to the direct summand $B$ equals the identity of $B$;
hence $X=0$ and (iv) follows.
\end{proof}

\begin{lemma}\label{lem_figure.traces} 
If $M$ is invertible, then $\tr_{M/A}(\one_M)=u_M$.
\end{lemma}
\begin{proof} Pick arbitrary $f\in M^*_*$, $m,n\in M_*$.
Then, directly from the definition of $u_M$ we
deduce that $f(m)\cdot n=u_M\cdot f(n)\cdot m$. In other 
words, $\omega_{M/A}(n\otimes f)=
u_M\cdot\ev_{M/A}(n\otimes f)\cdot\one_M$. By linearity
we deduce that $\phi=u_M\cdot\tr_{M/A}(\phi)\cdot\one_M$
for every $\phi\in\End_A(M)$. By letting $\phi:=\one_M$,
and taking traces on both sides, we obtain:
$\tr_{M/A}(\one_M)=u_M\cdot\tr_{M/A}(\one_M)^2$.
But since $M$ is invertible, $\tr_{M/A}(\one_M)$ 
is invertible in $A_*$, whence 
$u_M\cdot\tr_{M/A}(\one_M)=1$,
which is equivalent to the sought identity.
\end{proof}

\begin{proposition} Let $M$ be an invertible $A$-module.
Then:
\begin{enumerate}
\item
$M\otimes_AM$ is strictly invertible.
\item
There exists a natural decomposition 
$A\simeq A_1\times A_{-1}$ where $M\otimes_AA_1$ is
strictly invertible, $A_{-1*}$ is a $\Q$-algebra and
$\Sym_{A_{-1}}^2(M\otimes_AA_{-1})=0$.
\end{enumerate}
\end{proposition}
\begin{proof} (i): it is clear that $M^{\otimes n}$
is invertible for every $n$. Let $\sigma\in S_n$ be
any permutation; it is easy to verify that the morphism 
$\sigma_M:M^{\otimes n}\to M^{\otimes n}~:~
x_1\otimes x_2\otimes...\otimes x_n\mapsto
x_{\sigma(1)}\otimes x_{\sigma(2)}\otimes...\otimes
x_{\sigma(n)}$ equals 
$u_M^{\mathrm{sgn}(\sigma)}\cdot\one_{M^{\otimes n}}$.
Especially, the transposition operator on
$(M\otimes_AM)^{\otimes 2}$ acts via the permutation:
$x\otimes y\otimes z\otimes w\mapsto 
z\otimes w\otimes x\otimes y$ whose sign is even.
Therefore $u_{M\otimes_AM}=1$, which is (i).    

(ii): It follows from (i) that the antisymmetrizer operator 
$a_2$ on $(M\otimes_AM)^{\otimes 2}$ vanishes; 
{\em a fortiori\/} it vanishes on the quotient 
$(\Lambda^2_AM)^{\otimes 2}$, therefore 
$\Lambda^2_A(\Lambda^2_AM)\simeq
\Img(a_2:\Lambda^2_AM\to\Lambda^2_AM)=0$.
Then lemma \ref{lem_lambda.two} says that there exists a 
natural decomposition $A\simeq A_1\times A_{-1}$ such that 
$(\Lambda^2_AM)\otimes_AA_1=0$ and 
$(\Lambda^2_AM)\otimes_AA_{-1}$ is invertible.
To show that $A_{-1*}$ is a $\Q$-algebra, it is enough
to show that $A_{-1}/pA_{-1}=0$ for every prime
$p$. Up to replacing $A$ by $A/pA$, we reduce to
verifying that, if $\F_p\subset A_*$ and $M$ is invertible,
then $M$ is of almost finite rank. To this aim, it
suffices to verify that the equivalent condition (iv)
of theorem \ref{th_charact.almost.fin.rk} is satisfied.
If $B\neq 0$ and $M_B:=M\otimes_AB\simeq B\oplus X$, then
the argument in the proof of lemma 
\ref{lem_charact.strictly.inv} shows that $X=0$ and
therefore $M_B$ does not admit infinite splittings.
Finally, it remains only to verify condition (a) of 
{\em loc.~cit.\/} So suppose that $M_B\simeq B^n\oplus Q$.
If $n>0$, we have just seen that $Q=0$; if $n=0$,
and $Q/IQ=0$ for some ideal $I$, then by the 
faithfulness of $M$ (lemma \ref{lem_explain.invert}(iii)) 
we must have $I=B$; if $I$ is Frobenius nilpotent it
follows that $B=0$. Finally, set $M_{-1}:=M\otimes_AA_{-1}$; 
notice that, since $A_{-1*}$ is a $\Q$-algebra, the 
endomorphism group of $\hat\G_{m,A_{-1*}}$ is isomorphic 
to $A_{-1*}$, and therefore $\chi_{M_{-1}}(X)=(1+X)^\alpha$, 
where $\alpha$ is an element of $A_{-1*}$ which can be 
determined by looking at the coefficient of 
$\chi_{M_{-1}}(X)$ in degree $1$. One finds 
$\alpha=\tr_{M_{-1}/A_{-1}}(\one_{M_{-1}})$.
In view of lemma \ref{lem_figure.traces}, we can rewrite
$\alpha=u_{M_{-1}}$; therefore 
\set\begin{equation}\label{eq_too.many.minus}
\tr_{\Lambda^2_{A_{-1}}M_{-1}/A_{-1}}
(\one_{\Lambda^2_{A_{-1}}M_{-1}})=\binom{u_{M_{-1}}}{2}.
\end{equation}
On the other hand, since $\Lambda^2_{A_{-1}}M_{-1}$
is an invertible $A_{-1}$-module of finite rank, we
know that the left-hand side of \eqref{eq_too.many.minus}
equals $1$; consequently $u_{M_{-1}}=-1$. This means that,
in $M^{\otimes 2}_{-1}$, the identity 
$x\otimes y=-y\otimes x$ holds for every $x,y\in M_{-1*}$;
therefore, the kernel of the projection 
$M^{\otimes 2}_{-1}\to\Sym^2_{A_{-1}}M_{-1}$ contains
all the elements of the form $2\cdot x\otimes y$; in
other words, multiplication by $2$ is the zero morphism
in $\Sym^2_{A_{-1}}M_{-1}$; since $A_{-1*}$ is a
$\Q$-algebra, this at last shows that 
$\Sym^2_{A_{-1}}M_{-1}$ vanishes, and concludes the proof
of the proposition.
\end{proof}

\subsection{Construction of quotients by flat equivalence
relations}\label{sec_groupoids}
We will need to recall some generalities on groupoids, which 
we borrow from \cite[Exp. V]{SGA3}. 
\sset\subsubsection{}\label{subsec_interpret}
\index{Groupoid(s)|indref{subsec_interpret}}
\index{Groupoid(s)!with trivial automorphisms|indref{subsec_interpret}}
If $\cC$ is any category admitting fibred products and a final object,
a {\em $\cC$-groupoid\/} is the datum of two objects $X_0$, $X_1$ of 
$\cC$, together with "source" and "target" morphisms 
$s,t:X_1\to X_0$, an "identity" morphism $\iota:X_0\to X_1$ 
and a further "composition" morphism 
$c:X_2\to X_1$, where $X_2$ is the fibre product in the cartesian 
diagram:
$$
\xymatrix{ X_2 \ar[r]^{t'} \ar[d]_{s'} & X_1 \ar[d]^s \\
           X_1 \ar[r]^t & X_0. }
$$
The datum $(X_0,X_1,s,t,c,\iota)$ is subject to the following
condition. For every object $S$ of $\cC$, the set 
$X_0(S):=\Hom_{\cC}(S,X_0)$ is a groupoid, with set of morphisms
given by $X_1(S)$, and for every $\phi\in X_1(S)$, the source
and target of $\phi$ are respectively $s(\phi):=\phi\circ s$
and $t(s):=\phi\circ s$; furthermore the composition law in 
$X_0(S)$ is given by $c(S):X_1(S)\times_{X_0(S)}X_1(S)\to X_1(S)$.
The above conditions amount to saying that 
$s\circ\iota=t\circ\iota=\one_{X_0}$ and the commutative diagrams 
\set\begin{equation}\label{eq_upperlower}{
\diagram X_2 \ar@<.5ex>[r]^c \ar@<-.5ex>[r]_{t'}
               \ar[d]_{s'} & X_1 \ar[d]^s & 
         X_2 \ar[r]^c \ar[d]_{t'} & X_1 \ar[d]^t \\
           X_1 \ar@<.5ex>[r]^s \ar@<-.5ex>[r]_t & X_0 &
         X_1 \ar[r]^t & X_0
\enddiagram
}\end{equation}
are cartesian in $\cC$ both for the square made up from the 
upper arrows and for the square made up from the lower arrows 
(cp. \cite[Exp. V \S1]{SGA3}).

One says that the groupoid $G:=(X_0,X_1,s,t,c,\iota)$ {\em has trivial
automorphisms}, if the morphism $(s,t):X_1\to X_0\times X_0$ 
is a (categorical) monomorphism. (This translates in categorical
terms the requirement that for every object $S$ of $\cC$, and
every $x\in X_0(S)$, the automorphism group of $x$ in $X_0(S)$
is trivial).

It is sometimes convenient to denote by $X\times_{(\alpha,\beta)}Z$
the fibre product of two morphisms $\alpha:X\to Y$ and $\beta:Z\to Y$.

\sset\subsubsection{}
Given a groupoid $G$, and a morphism $X_0\to X_0'$, we
obtain a new groupoid $G\times_{X_0}X_0'$ by taking the
datum $(X_0',X_1\times_{X_0}X_0',s\times_{X_0}\one_{X_0'},
t\times_{X_0}\one_{X_0'},c\times_{X_0}\one_{X_0'},
\iota\times_{X_0}\one_{X_0'})$.

Moreover, suppose that $\cC$ admits finite coproducts and 
that all such coproducts are disjoint universal (cp. 
\cite[Exp.II, Def.4.5]{SGA4-1}). Denote by $Y\amalg Z$ the 
coproduct of two objects $Y$ and $Z$ of $\cC$. Let 
$G':=(X'_0,X'_1,s',t',c',\iota')$ be another groupoid;
one can define a groupoid $G\amalg G'$ by taking the
datum $(X_0\amalg X_0',X_1\amalg X_1',s\amalg s',t\amalg t',
c\amalg c',\iota\amalg\iota')$.

\sset\subsubsection{}
In the following we will be concerned with groupoids in the
category $A\Alg^o$ of $A$-schemes, where $A$ is any $V^a$-algebra.
We will use the general terminology for almost schemes introduced
in \eqref{subsec_general-schemes}, and complemented by the following: 

\begin{definition}\label{def_gen.groupoids}
\index{Groupoid(s)!$X_0/G$ : quotient 
by an action of a|indref{def_gen.groupoids}}
\index{Groupoid(s)!closed equivalence relation, flat, \'etale,
almost finite, almost projective, of finite rank|indref{def_gen.groupoids}}
\index{Almost scheme(s)!closed imbedding of|indref{def_gen.groupoids}}
\index{Almost scheme(s)!open and closed imbedding of|indref{def_gen.groupoids}}
\index{Almost scheme(s)!almost finite, \'etale, flat,
almost projective morphism of|indref{def_gen.groupoids}}
Let $\phi:X\to Y$ be a morphism of $A$-schemes, 
$G:=(X_0,X_1,s,t,c,\iota)$ a groupoid in the category of $A$-schemes.
\begin{enumerate}
\item
We say that $\phi$ is a {\em closed imbedding\/} 
(resp. is {\em almost finite}, resp. is {\em {\'e}tale}, resp. is 
{\em flat}, resp. is {\em almost projective}) if the corresponding 
morphism $\phi^\sharp:\cO_Y\to\cO_X$ is an epimorphism of $\cO_Y$-modules 
(resp. enjoys the same property).
We say that $\phi$ is an {\em open and closed imbedding\/}
if it induces an isomorphism $X\stackrel{\sim}{\to}Y_1$
onto one of the factors of a decomposition $Y=Y_1\amalg Y_2$.
\item
We say that $G$ is a {\em closed equivalence relation\/} if 
the morphism $(s,t):X_1\to X_0\times X_0$ is a a closed 
imbedding. We say that $G$ is {\em flat} (resp. {\em {\'e}tale},
resp. {\em almost finite}, resp. {\em almost projective}) if
the morphism $s:X_1\to X_0$ enjoys the same property.
We say that $G$ is {\em of finite rank\/} if $\cO_{X_1}$ is an 
almost projective $\cO_{X_0}$-module of finite rank. 
Furthermore, we set $X_0/G:=\Spec\,\cO_{X_0}^G$, where 
$\cO_{X_0}^G\subset\cO_{X_0}$ is the equalizer of the morphisms
$s^\sharp$ and $t^\sharp$.
\end{enumerate}
\end{definition}

\sset\subsubsection{}\label{subsec_decompose.gpd}
Let $G:=(X_0,X_1,s,t,c,\iota)$ be a groupoid of finite rank in 
$A\Alg^o$. by assumption
$\cO_{X_1}$ is an almost projective $\cO_{X_0}$-module of finite rank,
hence, by proposition \ref{prop_decomp.fin.rank}, there is a 
decomposition $\cO_{X_0}\simeq\prod^r_{i=0}B_i$ such that
$C_i:=\cO_{X_1}\otimes_{\cO_{X_0}}B_i$ is of constant rank equal to $i$,
for $i=0,...,r$. Set $X_{0,i}:=\Spec\,B_i$.

\begin{lemma}\label{lem_decompose.gpd}
In the situation of \eqref{subsec_decompose.gpd}, there
is a natural isomorphism of groupoids:
$$
G\simeq(G\times_{X_0}X_{0,1})\amalg...\amalg(G\times_{X_0}X_{0,r}).
$$
\end{lemma}
\begin{proof} For every $i\leq r$, let $\alpha_i:X_{0,i}\to X_0$
be the open and closed imbedding defined by \eqref{subsec_decompose.gpd}.
Set $X_{1,i}:=X_{0,i}\times_{(\alpha_i,s)}X_1$ (so $X_{1,i}=\Spec\,C_i)$.
Moreover, let $X'_{1,i}:=X_{0,i}\times_{(\alpha_i,t)}X_1$,
$\beta_i:X'_{1,i}\to X_1$ the open and closed imbedding
(obtained by pulling back $\alpha_i$), 
$X_{2,i}:=X_{1,i}\times_{(\beta_i,s')}X_2$ and
$X_{2,i}:=X'_{1,i}\times_{(\beta_i,s')}X_2$. There follow
natural decompositions $X_1\simeq X'_{1,1}\amalg...\amalg X'_{1,r}$
and $X_2\simeq X'_{2,1}\amalg...\amalg X'_{2,r}$, such that $s'$
decomposes as a coproduct of morphisms $X'_{2,i}\to X_{1,i}'$.
By the construction of $X_{0,i}$, it is clear that $X^{\prime o}_{2,i}$
has rank equal to $i$ as an $X^{\prime o}_{1,i}$-module, for every $i\leq r$.
In other words, the above decompositions fulfill the conditions
of proposition \ref{prop_decomp.fin.rank}. Similarly, we
obtain decompositions $X_1\simeq X_{1,1}\amalg...\amalg X_{1,r}$
and $X_2\simeq X_{2,1}\amalg...\amalg X_{2,r}$ which fulfill
the same conditions. However, these conditions characterize
uniquely the factors occuring in it, thus $X_{1,i}=X_{1,i}'$
for $i\leq r$. The claim follows easily.
\end{proof}

\begin{lemma}\label{lem_groupoid} Let $G:=(X_0,X_1,s,t,c,\iota)$
be a groupoid of finite rank in $A\Alg^o$. If $G$ has trivial 
automorphisms, then it is a closed equivalence relation. 
\end{lemma}
\begin{proof} Using lemma 
\ref{lem_decompose.gpd} we reduce easily to the case where 
$\cO_{X_1}$ is of constant rank, say equal to $r\in\N$.
Let $Y:=X_0\times X_0$; since $G$ has trivial automorphisms,
the morphism $(s,t):X_1\to Y$ is a monomorphism; equivalently,
the natural projections $\pr_1,\pr_2:X_1\times_YX_1\to X_1$ are 
isomorphisms. Let $D:=\Img((s,t)^\sharp:\cO_Y\to\cO_{X_1})$;
it follows that the natural morphisms
$\pr_1^\sharp,\pr_2^\sharp:\cO_{X_1}\to\cO_{X_1}\otimes_D\cO_{X_1}$
are isomorphisms and consequently, 
\set\begin{equation}\label{eq_knowalready}
(\cO_{X_1}/D)\otimes_D\cO_{X_1}=0.
\end{equation}
We need to show that $\cO_{X_1}=D$, or equivalently, that
$\cO_{X_1}/D=0$. However, by theorem \ref{th_free.in.fpqc}, we can
find a faithfully flat $\cO_{X_0}$-algebra $B$ such that 
$C:=B\otimes_{\cO_{X_0}}\cO_{X_1}\simeq B^r$.
Let $D':=B\otimes_BD$; it follows that $C$ is a faithful finitely
generated $D'$-module. It suffices to show that $C/D'=0$, and we
know already from \eqref{eq_knowalready} that $(C/D')\otimes_{D'}C=0$.
By proposition \ref{prop_descent} it follows that $C/D'$ is
a flat $D'$-module; consequently $C/D'\subset(C/D')\otimes_{D'}C$, 
and the claim follows.
\end{proof}

\sset\subsubsection{}
Let $B$ be an $A$-algebra, $P$ an almost finitely 
generated  projective $B$-module.
For every integer $i\geq 0$, we define a $B$-linear morphism
\set\begin{equation}\label{eq_noname}
\Gamma^i_B(\End_B(P)^a)\to B
\end{equation}
as follows (see \eqref{subsec_app.gammas.etc} for the definition 
of the functor $\Gamma^i_B:B\Mod\to B\Mod$). Let $R$ be any
$B_*$-algebra; we remark that the natural map
$\beta_P:\End_B(P)^a_!\otimes_{B_*}R\to
\End_{R^a}(R^a\otimes_AP)^a_!$ is an isomorphism. Hence,
for every $i\geq 0$, we can define a map of sets
$$
\lambda^i_R:
\End_B(P)^a_!\otimes_{B_*}R\to\End_B(\Lambda^i_BP)^a_!\otimes_{B_*}R
$$
by letting 
$\phi\mapsto\beta^{-1}_{\Lambda_B^iP}(\Lambda_{R^a}^i\beta_P(\phi)^a)$.
In the terminology of \cite{Ro}, the system of maps $\lambda^i_R$ 
forms a homogeneous polynomial law of degree $i$ from 
$\End_B(P)^a_!$ to $\End_B(\Lambda^i_BP)^a_!$, so it induces a
$B_*$-linear map 
$\bar{\lambda^i}:\Gamma_{B_*}^i(\End_B(P)^a_!)\to\End_B(\Lambda^i_BP)^a_!$.
After passing to almost modules, we obtain a $B$-linear morphism
\set\begin{equation}\label{eq_pass.to.almost}
\Gamma_B^i(\End_B(P)^a)\to\End_B(\Lambda^i_BP)^a.
\end{equation}
Then \eqref{eq_noname} is defined as the composition of 
\eqref{eq_pass.to.almost} and the trace morphism $\tr_{\Lambda^i_BP/B}$.

\sset\subsubsection{}\label{subsec_symm.gammas}
Let $C$ be an almost finite projective $B$-algebra.
Define $\mu:C\to\End_B(C)^a$ as in \eqref{subsec_trace.of.algebra}.
By composition of $\Gamma^i_B\mu$ and \eqref{eq_noname} we
obtain a $B$-linear morphism
$$
\Gamma^i_BC\to B
$$
characterized by the condition: 
$c^{[i]}\mapsto\sigma_i(c):=\tr_{C/B}(\Lambda^i_C\mu(c))$.

\sset\subsubsection{}\label{subsec_symm.groups}
The construction of \eqref{subsec_symm.gammas} applies
especially to an almost finite projective groupoid
$G:=(X_0,X_1,s,t,c,\iota)$. In such case one verifies, using the
cartesian diagrams \eqref{eq_upperlower}, that 
$\sigma_i(t^\sharp(f))\in\cO_{X_0}^G$ for every $f\in\cO_{X_{0*}}$
and every $i\leq r$: the argument is the same as in the proof of
\cite[Exp.V, Th.4.1]{SGA3}. In this way one obtains 
$\cO_{X_0}^G$-linear morphisms
\set\begin{equation}\label{eq_exhibit}
T_{G,i}:\Gamma^i_{\cO_{X_0}^G}\cO_{X_0}
\stackrel{\Gamma^it^\sharp}{\longrightarrow}
\Gamma^i_{\cO_{X_0}}\cO_{X_1}\to\cO_{X_0}^G\quad
f^{[i]}\mapsto\sigma_i(t^\sharp(f)).
\end{equation}
\begin{theorem}\label{th_case.etaleclosed} 
Let $G:=(X_0,X_1,s,t,c,\iota)$ be an {\'e}tale almost finite and closed 
equivalence relation in $A\Alg^o$. Then $G$ is effective 
and the natural morphism $X_0\to X_0/G$ is {\'e}tale and almost finite
projective.
\end{theorem}
\begin{proof} See \cite[Exp.IV, \S 3.3]{SGA3} for the definition
of effective equivalence relation. By \eqref{eq_upperlower}, we have an 
identification $X_2\simeq X_1\times_{(s,s)}X_1$.; therefore, the 
natural diagonal morphism $X_1\to X_1\times_{(s,s)}X_1$ gives a section
$\delta:X_1\to X_2$ of the morphism $s':X_2\to X_1$. Furthermore,
since $X_2=X_1\times_{(s,t)}X_1$, the pair of morphisms 
$(\one_{X_1},\iota\circ s):X_1\to X_1$ induces another morphism
$\psi_0:X_1\to X_2$; similarly, let $\psi_1:X_1\to X_2$
be the morphism induced by the pair $(\iota\circ t,\one_{X_1})$
(these are the degeneracy maps of the simplicial complex associated
to $G$: cp. \cite[Exp. V, \S 1]{SGA3}). By arguing
with $T$-points (and exploiting the interpretation 
\eqref{subsec_interpret} of $X_0(T)$, $X_1(T)$, etc.) one checks
easily, first that $\psi_1=\delta$, and second,  that the 
two commutative diagrams
\set\begin{equation}\label{eq_yet.more.diags}
{\diagram
           X_0 \ar[r]^\iota \ar[d]_\iota & X_1 \ar[d]^{\psi_1} &
           X_1 \ar[r]^t \ar[d]_{\psi_1} & X_0 \ar[d]^\iota \\
           X_1 \ar[r]^{\psi_0} & X_2 & X_2 \ar[r]^{t'} & X_1 
\enddiagram}
\end{equation}
are cartesian. Since by assumption $s$ is {\'e}tale, corollary 
\ref{cor_unram} implies that $\delta$ is an open and closed imbedding;
consequently the same holds for $\iota$.
Let $e_0\in\cO_{X_{1*}}$ (resp. $e_1\in\cO_{X_{2*}}$) be the idempotent 
corresponding to the open and closed imbedding $\iota$
(resp. $\delta$); since $G$ is a closed equivalence relation, 
for every $\eps\in\fm$ we can write 
$\eps\cdot e_0=\sum^n_is^\sharp(b_i)\cdot t^\sharp(b'_i)$ for some 
$b_i,b_i'\in\cO_{X_{1*}}$. In view of \eqref{eq_yet.more.diags}
we deduce that $\eps\cdot e_1=
\sum^n_i(t^{\prime\sharp}\circ s^\sharp(b_i))\cdot
(t^{\prime\sharp}\circ t^\sharp(b'_i))$.
However, $s\circ t'=t\circ s'$ and $t\circ t'=t\circ c$, consequently
$$
\eps\cdot e_1=\sum^n_i(s^{\prime\sharp}\circ t^\sharp(b_i))\cdot
(c^\sharp\circ t^\sharp(b'_i)).
$$
Finally, thanks to remark \ref{rem_zeta.identity}, and again
\eqref{eq_upperlower}, we can write:
\set\begin{equation}\label{eq_go.home}
\eps\cdot f=\sum^n_is^\sharp\circ
\Tr_{X_1/X_0}(f\cdot t^\sharp(b_i))\cdot t^\sharp(b_i')
\quad\text{for every $f\in\cO_{X_{1*}}$.}
\end{equation}
If we now let $f:=t^\sharp(g)$ in \eqref{eq_go.home} we deduce:
$\eps\cdot t^\sharp(g)=
\sum_is^\sharp(T_{G,1}(g\cdot b_i))\cdot t^\sharp(b_i')
=\sum_it^\sharp(T_{G,1}(g\cdot b_i))\cdot t^\sharp(b_i')$
for every $g\in\cO_{X_{0*}}$. Since $t^\sharp$ is injective, this means that:
\set\begin{equation}\label{eq_dropping.t}
\eps\cdot g=\sum^n_iT_{G,1}(g\cdot b_i)\cdot b_i'\qquad
\text{for every $g\in\cO_{X_{0*}}$.}
\end{equation}
It follows easily that $\cO_{X_0}$ is an almost finitely generated projective
$\cO_{X_0}^G$-module. Furthermore, let us introduce the bilinear pairing
$t_G:=T_G\circ\mu_{\cO_{X_0}/\cO_{X_0}^G}:
\cO_{X_0}\otimes_{\cO_{X_0}^G}\cO_{X_0}\to\cO_{X_0}^G$.

\begin{claim}\label{cl_T.is.a.trace} $t_G$ is a perfect pairing.
\end{claim}
\begin{pfclaim} We have to show that the associated $\cO_{X_0}$-linear
morphism 
$$
\tau_G:\cO_{X_0}\to\cO_{X_0}^*:=
\Alhom_{\cO_{X_0}^G}(\cO_{X_0},\cO_{X_0}^G)
$$
is an isomorphism. From \eqref{eq_dropping.t} it follows easily that
$\tau_G$ is a monomorphism. Let $\phi:\cO_{X_0}\to\cO_{X_0}^G$ be a
$\cO_{X_0}^G$-linear morphism; it remains only to show that, for every
$\eps\in\fm$, there exists $b\in\cO_{X_{0*}}$ such that 
$\tau_G(b)=\eps\cdot\phi$.  
Let $\alpha:\cO_{X_1}\to\cO_{X_0}^n$, $\beta:\cO_{X_0}^n\to\cO_{X_1}$
be defined by the rules: 
$c\mapsto(\Tr_{X_1/X_0}(c\cdot b_1),...,\Tr_{X_1/X_0}(c\cdot b_n))$
and $(x_1,...,x_n)\mapsto\sum^n_ix_i\cdot b_i'$ for all $c\in\cO_{X_{1*}}$,
$x_1,...,x_n\in\cO_{X_{0*}}$. We remark that 
$\Img(\alpha\circ t^\sharp)\subset(\cO_{X_0}^G)^n$ and 
$\Img(\beta\circ(t^\sharp)^n)\subset\cO_{X_0}$, so that 
we deduce, by restriction, morphisms 
$\alpha_0:\cO_{X_0}\to(\cO_{X_0}^G)^n$ and 
$\beta_0:(\cO_{X_0}^G)^n\to\cO_{X_0}$.
Let $\psi:=((\phi\circ\beta_0)\otimes_{\cO_{X_0}^G}\one_{\cO_{X_0}})
\circ\alpha:\cO_{X_1}\to\cO_{X_0}$.
By theorem \ref{th_proj.etale} we can find, for every $\eps\in\fm$,
an element $c\in\cO_{X_{1*}}$ such that 
$\eps\cdot\psi=\tau_{X_1/X_0}(c)$.
Using \eqref{eq_go.home} we derive easily 
$\eps\cdot c=\eps\cdot\sum^n_i\psi(t^\sharp(b_i))\cdot t^\sharp(b'_i)=
\eps\cdot\sum^n_i\phi\circ\beta_0\circ\alpha_0(b_i)\cdot t^\sharp(b_i')=
\eps\cdot\sum^n_i\phi(\eps\cdot b_i)\cdot t^\sharp(b_i)$. In particular,
$\eps\cdot c=t^\sharp(b)$ for some $b\in B$, so the claim follows.
\end{pfclaim}

By assumption, the morphism
$\pi:C:=\cO_{X_0}\otimes_{\cO_{X_0}^G}\cO_{X_0}\to\cO_{X_1}$ induced
by the pair $(s^\sharp,t^\sharp)$ is an epimorphism. Moreover, by
construction, we have the identity:
$$
\Tr_{X_1/X_0}\circ\pi=\one_{\cO_{X_0}}\otimes_{\cO_{X_0}^G}T_{G,1}.
$$
By claim \ref{cl_T.is.a.trace} we see that
$\one_{\cO_{X_0}}\otimes_{\cO_{X_0}^G}T_{G,1}$ induces a perfect pairing
$C\otimes_{\cO_{X_0}}C\to\cO_{X_0}$; on the other hand,
$\Tr_{X_1/X_0}$ is already
a perfect pairing, by theorem \ref{th_proj.etale}. It then follows
that $\pi$ must be a monomorphism, hence $C\simeq\cO_{X_1}$, which
shows that $G$ is effective; then it is easy to verify that $T_{G,1}$ 
is actually the trace of the $\cO_{X_0}^G$-algebra $\cO_{X_0}$, which
is consequently {\'e}tale over $\cO_{X_0}^G$.
\end{proof}

\begin{proposition}\label{prop_groupoids} 
Let $G=(X_0,X_1,s,t,c,\iota)$ be a groupoid of finite rank. Then 
$\cO_{X_{0*}}$ is integral over $\cO^G_{X_{0*}}$.
\end{proposition}
\begin{proof} By lemma \ref{lem_decompose.gpd} we can reduce to 
the case where the rank of $\cO_{X_1}$ is constant, say equal to $r$.
The assertion is then a direct consequence of the following: 
\begin{claim}
Let $f\in\cO_{X_{0*}}$. With the notation of \eqref{subsec_symm.groups}
we have:
$$
(t^\sharp(f))^r+T_{G,1}(f)\cdot (t^\sharp(f))^{r-1}+T_{G,2}(f)\cdot
(t^\sharp(f))^{r-2}+...+T_{G,r}(f)=0.
$$
\end{claim}
\begin{pfclaim} We apply proposition \ref{prop_Cayley-Ham}
({\em i.e.} Cayley-Hamilton's theorem) to the endomorphism 
$t^\sharp(f)\cdot\one_{\cO_{X_1}}:\cO_{X_1}\to\cO_{X_1}$.
\end{pfclaim}
\end{proof}

\begin{proposition}\label{prop_all.hyps} 
Let $G$ be an {\'e}tale closed equivalence relation of finite rank. 
Then $G$ is universally effective and the morphism $X_0\to X_0/G$ 
is {\'e}tale, faithfully flat and almost finite projective.
\end{proposition}
\begin{proof} Everything is known by theorem \ref{th_case.etaleclosed}, 
except for the faithfulness, which follows from the following:
\begin{claim} Under the assumptions of proposition \ref{prop_groupoids},
let $I\subset\cO_{X_0}^G$ be an ideal such that $I\cdot\cO_{X_0}=\cO_{X_0}$.
Then $I=\cO_{X_0}^G$.
\end{claim}
\begin{pfclaim} First of all, let $\cO_{X_0}\simeq\prod_{i=0}^rB_i$ be
the decomposition as in \eqref{subsec_decompose.gpd}; one derives easily a 
corresponding decomposition $\cO_{X_0}^G\simeq\prod_{i=0}^rB^G_i$, 
so we can assume that the rank of $\cO_{X_1}$ is constant, 
equal to $r$. Let $J\subset\cO_{X_0}^G$ be any ideal, and set
$C:=\cO_{X_0}^G/J$. We have a natural isomorphism
$$
\Gamma^i_C(\cO_{X_0}\otimes_{\cO_{X_0}^G}C)\simeq
\Gamma^i_{\cO_{X_0}^G}(\cO_{X_0})\otimes_{\cO_{X_0}^G}C.
$$
Composing with \eqref{eq_exhibit}$\:\otimes_{\cO_{X_0}^G}C$, we derive
a $C$-linear morphism: 
$\psi_i:\Gamma^i_C(\cO_{X_0}\otimes_{\cO_{X_0}^G}C)\to C$.
By inspecting the construction, one shows easily that 
$\psi_i(1^{[i]})=\binom{r}{i}$ (indeed, by flat base change one reduces 
easily to the case where $\cO_{X_1}$ is a free $B$-module of rank $r$,
in which case the result is obvious). Let us now take $J=I$.
Then $\Gamma^i_C(\cO_{X_0}\otimes_{\cO_{X_0}^G}C)=0$ for every $i>0$,
whence $1=\psi_r(1^{[r]})=0$ in $C$, and the claim follows.
\end{pfclaim}
\end{proof}

\begin{proposition} Keep the assumptions of proposition
{\em\ref{prop_all.hyps}}. If $X_0$ is an almost finite (resp. 
almost finitely presented, resp. flat, resp. almost projective,
resp. weakly unramified, resp. unramified, resp. weakly {\'e}tale, 
resp. {\'e}tale) $\Spec\,A$-scheme, then the same holds for the
$\Spec\,A$-scheme $X_0/G$.
\end{proposition}
\begin{proof} By proposition \ref{prop_all.hyps}, $\cO_{X_0}$
is a faithful almost finitely generated projective 
$\cO^G_{X_0}$-module, hence $\cE_{\cO_{X_0}/\cO_{X_0}^G}=\cO^G_{X_0}$
by proposition \ref{prop_eval.ideal}(iv). 
It follows easily that, for every $\eps\in\fm$ there exists $n\in\N$ 
such that $\eps\cdot\one_{\cO_{X_0}^G}$ factors as a composition of 
$\cO_{X_0}^G$-linear morphisms:
\set\begin{equation}\label{eq_as-usual}
\cO_{X_0}^G\to\cO_{X_0}^n\to\cO_{X_0}^G.
\end{equation}
The assertions for ``almost finite" and for ``almost projective" 
are immediate consequences.
To prove the assertion for ``almost finitely presented" we use the
criterion of proposition \ref{prop_al.small}(ii). Indeed, let 
$(N_\lambda,\phi_{\lambda\mu}~|~\lambda)$ be a filtered system 
of $A$-modules; we apply the natural transformation \eqref{eq_was.nu} 
to the sequence of morphisms \eqref{eq_as-usual} : since $\cO_{X_0}$ is 
almost finitely presented, so is $\cO_{X_0}^n$, hence the claim follows 
by a little diagram chase. The assertions for ``flat" and 
``weakly unramified" are easy and shall be left to the reader.
To conclude, it suffices to consider the assertion for 
``unramified". Now, by proposition \ref{prop_all.hyps} it
follows that $\cO_{X_0}\otimes_A\cO_{X_0}$ is an almost finitely
generated projective $\cO_{X_0}^G\otimes_A\cO_{X_0}^G$-module; since
by assumption $\cO_{X_0}$ is an almost projective
$\cO_{X_0}\otimes_A\cO_{X_0}$-module, we deduce from
lemma \ref{lem_used.once} that $\cO_{X_0}$ is an almost projective 
$\cO_{X_0}^G\otimes_A\cO_{X_0}^G$-module. Using \eqref{eq_as-usual}
we deduce that $\cO_{X_0}^G$ is almost projective over 
$\cO_{X_0}^G\otimes_A\cO_{X_0}^G$ as well.
\end{proof}


\newpage


\def\fSet{\mathbf{f.Set}}
\renewcommand\emptyset{\varnothing}
\def\rad{\mathrm{rad}}
\def\sH{\text{\sf H}}
\def\sm{\mathrm{sm}}
\def\fpqc{\mathrm{fpqc}}
\def\sj{\text{\sf j}}
\def\Proj{\mathrm{Proj}}
\def\Iso{\underline{\mathrm{Iso}}}
\def\GL{\mathrm{GL}}


\section{Henselization and completion of almost algebras}\label{ch_hensel}
This chapter deals with more advanced aspects of almost commutative
algebras : we begin with the definitions of {\em Jacobson radical}
$\rad(A)$ of an almost algebra $A$, of {\em henselian pair} and
{\em henselization} of a pair $(A,I)$, where $I\subset A$ is an
ideal contained in $\rad(A)$. These notions are especially well
behaved when $I$ is a {\em tight} ideal (definition \ref{def_tight}),
in which case we can also prove a version of Nakayama's lemma
(lemma \ref{lem_Naka}).

In section \ref{sec_lin-topol} we explain what is a {\em linear topology\/}
on an $A$-module and an $A$-algebra; as usual one is most interested
in the case of $I$-adic topologies. In case $A$ is $I$-adically complete
and $I$ is tight, we show that the functor $B\mapsto B/IB$ from almost
finitely presented \'etale $A$-algebras to almost finitely presented
\'etale $A/I$-algebras is an equivalence (theorem \ref{th_lift.etale.cpte}).
For the proof we need some criteria to ensure that an $A$-algebra is
unramified under various conditions : such results are collected in
section \ref{sec_criteria}, especially in theorem \ref{th_unramif.critter}
and its corollary \ref{cor_oftheabove}.

In section \ref{sec_lift.hens}, theorem \ref{th_lift.etale.cpte} is
further generalized to the case where the pair $(A,I)$ is tight henselian
(see theorem \ref{th_lifts.proj.hensel}, that also contains an analogous
statement concerning almost finitely generated projective $A$-modules).
The proof is a formal patching argument, which can be outlined as follows.
First one reduces to the case where $I$ is principal, say generated
by $f$, and since $I$ is tight, one can assume that $f\in\fm$; hence,
we want to show that a given \'etale almost finitely presented
$A/fA$-algebra $B_0$ lifts uniquely to an $A$-algebra $B$ of the same
type; in view of section \ref{sec_lin-topol} one can lift $B_0$ to
an \'etale algebra $B^\wedge$ over the $f$-adic completion $A^\wedge$
of $A$; on the other hand, the almost spectrum $\Spec\,A$ is a usual
scheme away from the closed subscheme defined by $I$, so we can use
standard algebraic geometry to lift $B_0[f^{-1}]$ to an \'etale algebra
$B'$ over $A[f^{-1}]$. Finally we need to show that $B^\wedge$ and $B'$
can be patched in a unique way over $\Spec\,A$; this amounts to showing
that certain commutative diagrams of functors are $2$-cartesian
(proposition \ref{prop_descent.hensel}).

The techniques needed to construct $B'$ are borrowed from Elkik's
paper \cite{Elk}; for our purpose we need to extend and refine slightly
Elkik's results, to deal with non-noetherian rings. This material is
presented in section \ref{sec_elkik}; its usefulness transcends the
modest applications to almost ring theory presented here.

The second main thread of the chapter is the study of the {\em smooth
locus\/} of an almost scheme; first we consider the affine case: as
usual, an affine scheme $X$ over $S:=\Spec\,A$ can be identified with the
fpqc sheaf that it represents; then the smooth locus $X_\sm$ of $X$
is a certain natural subsheaf, defined in terms of the cotangent
complex $\L_{X/S}$. To proceed beyond simple generalities one
needs to impose some finiteness conditions on $X$, whence the definition
of {\em almost finitely presented\/} scheme over $S$. For such affine
$S$-schemes we can show that the smooth locus enjoys a property which
we could call {\em almost formal smoothness}. Namely, suppose that
$I\subset A$ is a tight ideal such that the pair $(A,I)$ is henselian;
suppose furthermore that $\sigma_0:S_0:=\Spec\,A/I\to X$ is a section
lying in the smooth locus of $X$; in this situation it does not
necessarily follow that $\sigma_0$ extends to a full section
$\sigma:S\to X$, however $\sigma$ always exists if $\sigma_0$
extends to a section over some closed subscheme of the form
$\Spec\,A/\fm_0I$ (for a finitely generated subideal $\fm_0\subset\fm$).

Next we consider quasi-projective almost schemes; if $X$ is
such a scheme, the invertible sheaf $\cO_X(1)$ defines a
quasi-affine $\G_m$-torsor $Y\to X$, and we define the smooth
locus $X_\sm$ just as the projection of the smooth locus of $Y$.
This is presumably not the best way to define $X_\sm$, but anyway
it suffices for the applications of section \ref{sec_lift-tors}.
In the latter we consider again a tight henselian pair $(A,I)$,
and we study some deformation problems for $G$-torsors, where
$G$ is a closed subgroup scheme of $\mathrm{GL}_n$ defined over
$\Spec\,A$ and fulfilling certain general assumptions (see
\eqref{subsec_G-torsors}).
For instance, theorem \ref{th_last-one} says that every $G$-torsor
over the closed subscheme $\Spec\,A/I$ extends to a $G$-torsor
over the whole of $\Spec\,A$; the extension is however not
unique, but any two such extensions are close in a precise
sense (theorem \ref{th_lift-iso-tors}) : here the almost formal
smoothness of the quasi-projective almost scheme $(\GL_n/G)^a$
comes into play and accounts for the special quirks of the
situation.

\subsection{Henselian pairs}

\begin{definition}\label{def_Jacobson}
\index{Almost algebra(s)!$\rad(A)$ : Jacobson radical|indref{def_Jacobson}}
Let $A$ be a $V^a$-algebra. The {\em Jacobson radical\/}
of $A$ is the ideal $\rad(A):=\rad(A_*)^a$ (where, for a
ring $R$, we have denoted by $\rad(R)$ the usual Jacobson
ideal of $R$).
\end{definition}

\begin{lemma}\label{lem_Jacobson} 
Let $R$ be a $V$-algebra, $I\subset R$ an
ideal. Then $I^a\subset\rad(R^a)$ if and only if 
$\fm I\subset\rad(R)$.
\end{lemma}
\begin{proof} Let us remark the following:
\begin{claim}\label{cl_Jacobson} 
If $S$ is any ring, $J\subset S$ an ideal, then 
$J\subset\rad(S)$ if and only if, for every $x\in J$
there exists $y\in J$ such that $(1+x)\cdot(1+y)=1$.
\end{claim}
\begin{pfclaim} Suppose that $J\subset\rad(S)$ and let
$x\in J$; then $1+x$ is not contained in any maximal
ideal of $S$, so it is invertible. Find some $u\in S$
with $u\cdot(1+x)=1$; setting $y:=u-1$, we derive 
$y=-x-xy\in J$. Conversely, suppose that the condition 
of the claim holds for all $x\in J$. Let $x\in J$;
we have to show that $x\in\rad(S)$. If this were not
the case, there would be a maximal ideal $\fq\subset S$
such that $x\notin\fq$; then we could find $a\in S$ such
that $x\cdot a\equiv -1\pmod\fq$, so $1+x\cdot a\in\fq$,
especially $1+a\cdot x$ is not invertible, which contradicts
the assumption.
\end{pfclaim}

Let $\phi:\fm I\to R\to R^a_*$ be the natural composed
map.
\begin{claim}\label{cl_equiv.Jack} 
$\Img\,\phi$ is an ideal of $R^a_*$ and 
$\fm I\subset\rad(R)$ if and only if
$\Img\,\phi\subset\rad(R^a_*)$.
\end{claim}
\begin{pfclaim} The first assertion is easy to check, and
clearly we only have to verify the "if" direction of the
second assertion, so suppose that $\Img\,\phi\subset\rad(R^a_*)$. 
Notice that $\Ker\,\phi$ is a square-zero ideal of $R$. Then, 
using claim \ref{cl_Jacobson}, we deduce that
for every $x\in\fm\cdot I$ there exists $z\in\fm\cdot I$
such that $(1+x)\cdot(1+z)=1+a$, where $a\in\Ker\,\phi$, so
$a^2=0$. Consequently $(1+x)\cdot(1+z)\cdot(1-a)=1$, so
the element $y:=z-a-z\cdot a$ fulfills the condition of
claim \ref{cl_Jacobson}.
\end{pfclaim}

On the other hand, it is clear that 
$I^a\subset\rad(R^a):=\rad(R^a_*)^a$ if
and only if $\Img\,\phi\subset\rad(R^a_*)$, so the
lemma follows from claim \ref{cl_equiv.Jack}.
\end{proof}

Given a $V^a$-algebra $A$ and an ideal $I\subset\rad(A)$,
one can ask whether the obvious analogue of Nakayama's
lemma holds for almost finitely generated $A$-modules.
As stated in lemma \ref{lem_Naka}, this is indeed the case, 
at least if the ideal $I$ has the property singled out
by the following definition, which will play a constant
role in the sequel.

\begin{definition}\label{def_tight}
\index{Almost algebra(s)!ideal in an!tight|indref{def_tight}}
Let $I$ be an ideal of a $V^a$-algebra $A$. We say that $I$
is {\em tight\/} if there exists a finitely generated subideal
$\fm_0\subset\fm$ and an integer $n\in\N$ such that
$I^n\subset\fm_0A$.
\end{definition}

\begin{lemma}\label{lem_Naka} Let $A$ be a $V^a$-algebra, 
$I\subset\rad(A)$ a tight ideal. If $M$ is an almost 
finitely generated $A$-module with $I M=M$, we have $M=0$.
\end{lemma}
\begin{proof} Under the assumptions of the lemma we can find
a finitely generated $A$-module $Q$ such that 
$\fm_0 M\subset Q\subset M$. It follows that 
$M=I^n M\subset\fm_0 M\subset Q$, so $M=Q$
and actually $M$ is finitely generated. Let 
$M_0\subset M_*$ be a finitely generated $A_*$-submodule
with $M_0^a=M$; clearly $\fm M_0\subset I_* M_0$ and
$\fm(I_*)^n\subset\fm_0A_*$. Therefore: 
$$
\fm_0 M_0\subset\fm^{n+2} M_0\subset 
I_*\fm^{n+1} M_0\subset...\subset I_*^{n+1}\fm M_0
\subset I_*\fm\cdot\fm_0 M_0\subset\fm_0 M_0
$$
and consequently 
$\fm M_0=\fm_0 M_0=\fm I_*\cdot\fm_0 M_0$.
By lemma \ref{lem_Jacobson} we have $\fm I_*\subset\rad(A_*)$,
hence $\fm_0 M_0=0$ by Nakayama's lemma, thus finally
$\fm M_0=0$, {\em i.e.} $M=0$, as claimed.
\end{proof}

\begin{corollary}\label{cor_Naka} Let $A$, $I$ be as in lemma 
{\em\ref{lem_Naka}}; suppose that $\phi:N\to M$ is an $A$-linear 
morphism of almost finitely generated projective $A$-modules 
such that $\phi\otimes_A\one_{A/I}$ is an isomorphism. Then 
$\phi$ is an isomorphism.
\end{corollary}
\begin{proof} From the assumptions we derive that
$\Coker(\phi)\otimes_A(A/I)=0$, hence $\Coker\,\phi=0$, in
view of lemmata \ref{lem_long.forgotten}(i) and \ref{lem_Naka}. 
By lemma \ref{lem_long.forgotten}(ii) it follows that
$\Ker\,\phi$ is almost finitely generated; moreover, since 
$M$ is flat, 
$\Ker(\phi)\otimes_A(A/I)=\Ker(\phi\otimes_A\one_{A/I})=0$,
whence $\Ker\,\phi=0$, again by lemma \ref{lem_Naka}.
\end{proof}

\begin{definition}\label{def_hens-pais}
\index{Almost algebra(s)!(tight) henselian pair|indref{def_hens-pais}}
Let $A$ be a $V^a$-algebra, $I\subset\rad(A)$ 
an ideal. We say that $(A,I)$ is a {\em henselian pair\/}
if $(A_*,\fm\cdot I_*)$ is a henselian pair. If in addition,
$I$ is tight, we say that $(A,I)$ is a {\em tight henselian pair}.
\end{definition}

\begin{remark}\label{rem_hensel}
For the convenience of the reader, we recall without proofs
a few facts about henselian pairs. 

(i) For a ring $R$, let 
$R_\mathrm{red}:=R/\nil(R)$ and denote by $\sqrt{I}$ the 
radical of the ideal $I$. Then the pair $(R,I)$ is henselian
if and only if the same holds for the pair 
$(R_\mathrm{red},\sqrt{I\cdot R_\mathrm{red}})$.

(ii) Suppose that $I\subset\rad(R)$; then $(R,I)$ is a henselian 
pair if and only if the same holds for the pair $(\Z\oplus I,I)$, 
where $\Z\oplus I$ is endowed with the ring structure such that 
$(a,x)\cdot(b,y):=(ab,ay+bx+xy)$ for every $a,b\in\Z$, $x,y\in I$. 
Indeed, this follows easily from the following criterion (iii), which
is shown in \cite[Ch.XI, \S2, Prop.1]{Ray}. 

(iii) Assume $I\subset\rad(R)$; then the pair $(R,I)$ is henselian
if and only if every monic polynomial $p(X)\in R[X]$ such that 
$p(X)\equiv(X^2-X)^m\pmod{I[X]}$ decomposes as a product 
$p(X)=q(X)\cdot r(X)$ where $q(X),r(X)$ are monic polynomials in 
$R[X]$ and $q(X)\equiv X^m\pmod{I[X]}$, $r(X)\equiv(X-1)^m\pmod{I[X]}$.

(iv) Let $J\subset I$ be a subideal. If the pair $(R,I)$
is henselian, the same holds for $(R,J)$.

(v) If $(R,I)$ is a henselian pair, and $R\to S$ is an integral
ring homomorphism, then the pair $(S,IS)$ is henselian 
(\cite[Ch.XI, \S2, Prop. 2]{Ray}).
\end{remark}

\begin{lemma}\label{lem_hensel} 
Let $R$ be a $V$-algebra, $I\subset\rad(R)$ an ideal.
Then the pair $(R^a,I^a)$ is henselian if and only the same
holds for the pair $(R,\fm I)$.
\end{lemma}
\begin{proof} It comes down to checking that $(R,\fm I)$ is
henselian if and only if $(R^a_*,\fm I^a_*)$ is. To this 
aim, let $\phi:R\to R^a_*$ be the natural map. Let $S:=\Img\,\phi$ 
and $J:=\phi(\fm I)\subset S$. Since $\Ker\,\phi$ is a
square-zero ideal in $R$, it follows from remark \ref{rem_hensel}(i)
that $(R,\fm I)$ is henselian if and only if $(S,J)$
is. However, it is clear that the induced map $J\to\fm I^a_*$
is bijective, hence remark \ref{rem_hensel}(ii) and lemma 
\ref{lem_Jacobson} say that $(S,J)$ is henselian if and only if 
$(R^a_*,\fm I^a_*)$ is.
\end{proof}

\sset\subsubsection{}\label{subsec_plenty.hensel}
Suppose that $(R,I)$ is a henselian pair, with $R$ a $V$-algebra.
Then, in view of lemma \ref{lem_hensel} and remark 
\ref{rem_hensel}(iv) we see that $(R^a,I^a)$ is also henselian. 
This gives a way of producing plenty of henselian pairs.

\begin{lemma}\label{lem_plenty} Let $B$ be an almost finite 
$A$-algebra, $I\subset A$ an ideal.
\begin{enumerate}
\item
The induced ring homomorphism $A_*\to\phi(A_*)+\fm B_*$
is integral.
\item
If $I\subset\rad(A)$, then $I B\subset\rad(B)$.
\item
If $(A,I)$ is a henselian pair, the same holds for the pair
$(B,IB)$.
\end{enumerate}
\end{lemma}
\begin{proof} (i): for a given finitely generated subideal 
$\fm_0\subset\fm$, pick a finitely generated submodule 
$Q\subset B_*$ with $\fm_0B_*\subset Q$; notice that 
$(\fm_0 Q)^2\subset\fm_0 Q$, hence 
$\phi(A_*)+\fm_0 Q$ is a subring of $B_*$, finite over
$A_*$. As $\phi(A_*)+\fm B_*$ is a filtered union
of such subrings, the assertion follows. (ii) follows from
(i) and from lemma \ref{lem_Jacobson}. (iii) is a direct
consequence of (i), of remark \ref{rem_hensel}(v) and of 
\eqref{subsec_plenty.hensel}.
\end{proof}

\sset\subsubsection{}\label{subsec_hensel.is.fflat}
\index{Almost algebra(s)!henselization of|indref{subsec_hensel.is.fflat}}
Given a $V^a$-algebra $A$ and an ideal $I\subset A$, we define
the {\em henselization\/} of the pair $(A,I)$ as the unique
pair (up to unique isomorphism) $(A^\mathrm{h},I^\mathrm{h})$ 
which satisfies the (almost version of the) usual universal 
property (cp. \cite[Ch.XI, \S2, D{\'e}f.4]{Ray}). Suppose 
$A=R^a$ and $I=J^a$ for a $V$-algebra $R$ and an ideal $J\subset R$, 
and let $(R',J')$ be a henselization of the pair $(R,\fm J)$; 
one can easily check that $(R^{\prime a},J^{\prime a})$ is a
henselization of $(A,I)$. It follows in particular that, if
$I\subset\rad(A)$ and $(A^\mathrm{h},I^\mathrm{h})$ is a 
henselization of $(A,I)$, then the morphism $A\to A^\mathrm{h}$ 
is faithfully flat (\cite[Ch.XI, \S2]{Ray}).

\begin{lemma}\label{lem_betise.deux} 
Let $\fm_0\subset\fm$ be a finitely generated subideal. 
Then there exists an integer $n:=n(\fm_0)>0$ such that
$((\fm_0 A)_*)^n\subset\fm_0 A_*$ for any 
$V^a$-algebra $A$.
\end{lemma}
\begin{proof} We proceed by induction on the number $k$ of generators 
of $\fm_0$. To start out, let $B$ be any $V^a$-algebra and $f\in B_*$ 
any almost element. The endomorphism $B\to B~:~x\mapsto f\cdot x$ 
induces an isomorphism 
$\alpha:\bar B:=B/\Ann_B(f)\stackrel{\sim}{\to}fB$, 
and a commutative diagram:
$$
\xymatrix{ \bar B_*\otimes_V\bar B_* \ar[r]^-{\mu_{\bar B*}}
\ar[d]_{\alpha_*\otimes\alpha_*} & \bar B_* \ar[d]^{f\cdot\alpha_*} \\
(fB)_*\otimes_V(fB)_* \ar[r]^-{\mu_{B*}} & (fB)_*. 
}$$
It follows that 
\set\begin{equation}\label{eq_betise}
f\cdot(fB)_*=((fB)_*)^2
\end{equation} 
as ideals of $B_*$, which takes care of the case $k=1$. Suppose 
now that $k>1$; let us write $\fm_0=x_1\cdot V+\fm_1$, where 
$\fm_1$ is generated by $k-1$ elements. We apply \eqref{eq_betise}
with $B:=A/\fm_1A$ and $f=x_1$ to deduce: 
$$
(\fm_0A/\fm_1A)_*\cdot(\fm_0A/\fm_1A)_*=
x_1\cdot(\fm_0A/\fm_1A)_*\subset\Img((x_1A)_*\to(A/\fm_1A)_*).
$$
Therefore $((\fm_0 A)_*)^2\subset(x_1A)_*+(\fm_1A)_*$.
After raising the latter inclusion to some high power, the
inductive assumption on $\fm_1$ allows to conclude.
\end{proof}

\begin{corollary}\label{cor_hens.pairs} 
Let $A$ be a $V^a$-algebra, $I\subset\rad(A)$ a tight ideal. Then:
\begin{enumerate}
\item
$I_*\subset\rad(A_*)$.
\item
If $(A,I)$ is a henselian pair, then the same holds for the
pair $(A_*,I_*)$.
\end{enumerate}
\end{corollary}
\begin{proof} For integers $n,m$ large enough we can write:
$(I_*)^{nm}\subset(I^n_*)^m\subset(\fm_0A)_*^m\subset\fm_0A_*$,
thanks to lemma \ref{lem_betise.deux}. Then by lemma 
\ref{lem_Jacobson} we deduce $(I_*)^{nm+1}\subset\rad(A_*)$,
which implies (i). Similarly, we deduce from lemma
\ref{lem_betise.deux} that $\sqrt{\fm I}=\sqrt{I}$,
so (ii) follows from remark \ref{rem_hensel}(i).
\end{proof}

\begin{proposition}\label{prop_bij.tight.hens} 
Let $(A,I)$ be a tight henselian pair.
The natural morphism $A\to A/I$ induces a bijection from
the set of idempotents of $A_*$ to the set of idempotents 
of $(A/I)_*$.
\end{proposition}
\begin{proof} Pick an integer $m>0$ and a finitely generated
subideal $\fm_0\subset\fm$ such that $I^m\subset\fm_0\cdot A$.
We suppose first that $\tilde\fm$ has homological 
dimension $\leq 1$. By corollary \ref{cor_hens.pairs}(ii), the
quotient map $A_*\to A_*/I_*$ induces a bijection on 
idempotents. So we are reduced to showing that the
natural injective map $A_*/I_*\to(A/I)_*$ induces a surjection
on idempotents. To this aim, it suffices to show that
an idempotent almost element $\bar e:V^a\to A/I$ always lifts
to an almost element $e:V^a\to A$; indeed, the image
of $e$ inside $A_*/I_*$ will then necessarily agree with
$\bar e$. Now, the obstruction to the existence of $e$ is a 
class $\omega_1\in\Ext^1_{V^a}(V^a,I)$. On the other hand,
proposition \ref{prop_lift.idemp}(i)  ensures that, for
every integer $n>0$, $\bar e$ admits a unique idempotent
lifting $\bar e_n:V^a\to A/I^n$. Let more generally
$\omega_n\in\Ext^1_{V^a}(V^a,I^n)$ be the obstruction
to the existence of a lifting of $\bar e_n$ to an almost
element of $A$; the imbedding $I^n\subset I$ induces 
a map 
\set\begin{equation}\label{eq_map.of.Exts}
\Ext^1_{V^a}(V^a,I^n)\to\Ext^1_{V^a}(V^a,I)
\end{equation}
and clearly the image of $\omega_n$ under \eqref{eq_map.of.Exts}
agrees with $\omega_1$. Thus, the proposition will follow
in this case from the following:
\begin{claim}\label{cl_vanishing.trick} 
For $n$ sufficiently large, the map \eqref{eq_map.of.Exts} 
vanishes identically.
\end{claim} 
\begin{pfclaim} 
We will prove more precisely that the natural map
\set\begin{equation}\label{eq_natural.Exts2}
\Ext^1_{V^a}(V^a,\fm_0 I)\to\Ext^1_{V^a}(V^a,I)
\end{equation} 
vanishes; the claim will follow for $n:=m+1$. A choice
of generators $\eps_1,...,\eps_k$ for $\fm_0$ determines
an epimorphism $\phi:I^{\oplus k}\to\fm_0I$; notice that
$\Ext^2_{V^a}(V^a,\Ker\,\phi)=0$ due to lemma 
\ref{lem_Exts.are.same}(i),(ii)
and the assumption on the homological dimension of $\tilde\fm$.
Thus, the induced map
\set\begin{equation}\label{eq_natural.Exts}
\Ext^1_{V^a}(V^a,I^{\oplus k})\to
\Ext^1_{V^a}(V^a,\fm_0 I)
\end{equation}
is surjective. To prove the claim, it suffices therefore 
to show that the composition of \eqref{eq_natural.Exts} and
\eqref{eq_natural.Exts2} vanishes, which is obvious, since
the modules in question are almost zero.
\end{pfclaim}

Finally suppose that $\fm$ is arbitrary; by theorem
\ref{th_countably.pres}(i.b),(ii.b), $\fm$ is the colimit of a 
filtered family of subideals 
$(\fm_\alpha\subset\fm~|~\lambda\in \Lambda)$, such that 
$\fm_0\subset\fm_\lambda=\fm_\lambda^2$ and 
$\tilde\fm_\lambda:=\fm_\lambda\otimes_V\fm_\lambda$
is $V$-flat for every $\lambda\in\Lambda$, and moreover each
$\tilde\fm_\lambda$ has homological dimension $\leq 1$.
We may suppose that $(A,I)=(R,J)^a$ for a henselian pair
$(R,I)$, where $R$ is a $V$-algebra and $J^m\subset\fm_0 R$. 
Each pair $(V,\fm_\lambda)$ is a basic setup, and we denote
by $(A_\lambda,I_\lambda)$ the tight henselian pair
corresponding to $(R,J)$ in the almost category associated
to $(V,\fm_\lambda)$ (so $A_\lambda$ is the image of $R$ under 
the localization functor $V\Alg\to(V,\fm_\lambda)^a\Alg$).
Let $\bar e$ be an idempotent almost element of $A/I$ (which is
an object of $(V,\fm)^a\Alg$); $\bar e$ is represented by
a unique $V$-linear map $f:\tilde\fm\to R/J$ and by the 
foregoing, for every $\lambda\in\Lambda$ the restriction 
$f_\lambda:\tilde\fm_\lambda\to R/J$ lifts to a unique
idempotent map $g_\lambda:\tilde\fm_\lambda\to R$. By 
uniqueness, the maps $g_\lambda$ glue to a map 
$\colim{\lambda\in\Lambda}g_\lambda:\tilde\fm\to R$
which is the sought lifting of $\bar e$.
\end{proof}

\subsection{Criteria for unramified morphisms}\label{sec_criteria}
The following lemma generalizes a case of
\cite[Partie II, lemma 1.4.2.1]{Gr-Ra}.

\begin{lemma}\label{lem_Tor.orgy} 
Let $A\to C$ be a morphism of $V^a$-algebras, $f\in A_*$ 
any almost element and $M$ a $C$-module. Suppose that $M[f^{-1}]$ 
is a flat $C[f^{-1}]$-module, $M/fM$ is a flat $C/fC$-module,
$\Tor^A_1(C,A/fA)=0$ and $\Tor^A_i(M,A/fA)=0$  for $i=1,2$. 
Then $M$ is a flat $C$-module.
\end{lemma}
\begin{proof}  Let $N$ be any $C$-module; We need to show that 
$\Tor_1^C(M,N)=0$. To this aim we consider the short
exact sequence $0\to K\to N\stackrel{j}{\to} N[f^{-1}]\to L\to 0$
where $j$ is the natural morphism. We have 
$\Tor_1^C(M,N[f^-1])=\Tor_1^{C[f^{-1}]}(M[f^{-1}],N[f^{-1}])=0$,
therefore we are reduced to showing : 
\begin{claim} $\Tor^C_1(M,K)=\Tor_2^C(M,L)=0$.
\end{claim}
\begin{pfclaim} Notice that $K=\bigcup_{k>0}\Ann_K(f^k)$, and
similarly for $L$, so it suffices to show that $\Tor_i^C(M,Q)=0$
for $i=1,2$ and every $C$-module $Q$ such that 
$f^kQ=0$ for some integer $k>0$. By considering 
the short exact sequence
$0\to\Ann_Q(f^{k-1})\to Q\to Q'\to 0$, an easy induction on
$k$ further reduces to showing that $\Tor_i^C(M,Q)=0$ for
$i=1,2$, in case $fQ=0$. However, the morphisms $C\to C/fC$
and $A\to C$ determine base change spectral sequences 
(cp. \cite[Th.5.6.6]{We})
$$\begin{array}{r@{\:{:=}}l}
E^2_{pq} & 
\Tor^{C/fC}_p(\Tor^C_q(M,C/fC),Q)\Rightarrow\Tor^C_{p+q}(M,Q) \\
F^2_{pq} &
\Tor^C_p(\Tor^A_q(C,A/fA),M)\Rightarrow\Tor^A_{p+q}(M,A/fA).
\end{array}
$$
The spectral sequence $F$ yields an exact sequence:
$$\xymatrix{
\Tor^A_2(M,A/fA) \ar[r] & 
\Tor^C_2(C/fC,M) \ar[r] & 
\Tor^A_1(C,A/fA)\otimes_CM \ar `r[d] `[l] `[lld] `[dr] [dl] \\
& \Tor^A_1(M,A/fA)\ar[r] &\Tor^C_1(C/fC,M) \ar[r] & 0 
}$$
which implies $\Tor^C_i(M,C/fC)=0$  for $i=1,2$.
Thus $E^2_{pq}=0$ for $q=1,2$ and we get
$$
\Tor^C_i(M,Q)\simeq\Tor^{C/fC}_i(M/fM,Q)\qquad \text{for $i=1,2$.}
$$ 
Since $M/fM$ is a flat $C/fC$-module, the claim follows.
\end{pfclaim}
\end{proof}

\begin{lemma}\label{lem_vanishing.critter} 
Let $M$ be an $A$-module, $f\in A_*$ an almost element,
and denote by $A ^\wedge:=\liminv{n\in\N}(A/f^nA)$ the $f$-adic
completion of $A$. Then we have:
\begin{enumerate}
\item
$\Ann_A(M[f^{-1}])\cdot\Ann_A(M\otimes_AA^\wedge)\subset
\Ann_A(M)$.
\item
If $M$ is almost finitely generated and $M/fM=M[f^{-1}]=0$ 
then $M=0$.
\end{enumerate}
\end{lemma}
\begin{proof} (i): let $a\in\Ann_A(M[f^{-1}])_*$, 
$b\in\Ann_A(M\otimes_AA^\wedge)_*$ and denote by $\mu_a,\mu_b:M\to M$
the scalar multiplication morphisms. From $a\cdot M[f^{-1}]=0$
we deduce that $aM\subset\bigcup_{n>0}\Ann_M(f^n)$; it follow
that the natural morphism $aM\to(aM)\otimes_AA^\wedge$ is an
isomorphism. Now the claim follows by inspecting the commutative
diagram:
$$
\xymatrix{ M \ar[r] \ar[d]_{\mu_a} & 
M\otimes_AA^\wedge \ar[rr]^-{0=\mu_b\otimes\one_{A^\wedge}} & &
M\otimes_AA^\wedge \ar[d]^{\mu_a\otimes\one_{A^\wedge}} \\
aM \ar[r]^-{\mu_b} & (aM)\otimes_AA^\wedge \ar[rr]^-\sim & &
(aM)\otimes_AA^\wedge.
}$$

(ii): for a given finitely generated subideal
$\fm_0\subset\fm$, pick a finitely generated $A$-module
$Q\subset M$ such that $\fm_0 M\subset Q$. By 
assumption, $Q[f^{-1}]=0$; hence there exists an integer 
$n\geq 0$ such that $\fm_0\cdot f^n Q=0$,
whence $\fm_0^2\cdot f^n M=0$. However $M=fM$ by
assumption, thus $\fm_0^2 M=0$, and finally $M=0$,
since $\fm_0$ is arbitrary.
\end{proof}

\begin{theorem}\label{th_gen.by.idemps} 
Let $A$ be a $V^a$-algebra, $f\in A_*$ any
almost element and $I\subset A$ an almost finitely generated
ideal with $I^2=I$. Suppose that both $I[f^{-1}]\subset A[f^{-1}]$
and $(I+fA)/fA\subset A/fA$ are generated by idempotents
in $A[f^{-1}]_*$ and respectively $(A/fA)_*$. Then $I$
is generated by an idempotent of $A_*$.
\end{theorem}
\begin{proof} We start out with the following:
\begin{claim}\label{cl_invocation} Let $I,J\subset A$ be two 
ideals such that $J$ is nilpotent and $I^2=I$. If 
$\bar I:=\Img(I\to A/J)$ is generated by an idempotent of 
$(A/J)_*$, then $I$ is generated by an idempotent of $A_*$. 
\end{claim}
\begin{pfclaim} We apply proposition \ref{prop_lift.idemp}
to the non-unitary extension $0\to J\cap I\to I\to\bar I\to 0$
to derive that the idempotent $\bar e$ that generates 
$\bar I$ lifts uniquely to an idempotent $e\in I_*$.
Then $e$ induces decompositions $A\simeq A_0\times A_1$
and $I=(IA_0)\times(IA_1)$ such that $e=0$
(resp. $e=1$) on $A_0$ (resp. on $A_1$). We can therefore
reduce to the cases $A=A_0$ or $A=A_1$. If $A=A_0$, then 
$I\subset J$, hence $I=I^n=0$ for $n$ sufficiently large.
If $A=A_1$ then $1\in I_*$, so $I=A$; in either case,
the claim holds.
\end{pfclaim}

\begin{claim}\label{cl_crocodile} 
In the situation of \eqref{subsec_cartes.diagr},
suppose that $I\subset A_0$ is an ideal such that 
$I A_1$ and $I A_2$ are generated by idempotents.
Then $I$ is generated by an idempotent of $A_{0*}$.
\end{claim}
\begin{pfclaim} By applying termwise the functor $M\mapsto M_*$
to \eqref{eq_cartesian}, we obtain a commutative diagram
\eqref{eq_cartesian}$_*$ which is still cartesian.
By assumption, we can find idempotents $e_1\in A_{1*}$, 
$e_2\in A_{2*}$ with $e_i A_i=I A_i$ ($i=1,2$).
It is then easy to see that the images of $e_1$ and $e_2$
agree in $A_{3*}$; consequently there exists a unique
element $e_0\in A_{0*}$ whose image in $A_i$ agrees
with $e_i$ for $i=1,2$. Such element is necessarily an
idempotent; moreover, since the natural $A_0$-linear
morphism $I\to(I A_1)\times(I A_2)$ is a monomorphism,
we deduce easily that $(1-e_0)\cdot I=0$, whence 
$I\subset e_0 A_0$. Let $M:=(e_0 A_0)/I$; clearly 
$M\otimes_{A_0}A_i=0$ for $i=1,2$. Using lemma \ref{lem_Ferrand}(i) 
we deduce that $M=0$, whence $e_0\in I$, as stated.
\end{pfclaim}

Suppose next that $f$ is a regular element of $A_*$. Denote
by $A^\wedge$ the $f$-adic completion of $A$. One verifies
easily that the natural commutative diagram
\set\begin{equation}\label{eq_situation}{
\diagram
A \ar[r] \ar[d] & A[f^{-1}] \ar[d] \\
A^\wedge \ar[r] & A^\wedge[f^{-1}]
\enddiagram}\end{equation}
is cartesian.
\begin{claim}\label{cl_situation} 
In the situation of \eqref{eq_situation}, let $I\subset A$ be 
an ideal such that both $I[f^{-1}]\subset A[f^{-1}]$ and 
$I\cdot A^\wedge\subset A^\wedge$ are generated by idempotents.
Then $I$ is generated by an idempotent.
\end{claim}
\begin{pfclaim} As in the proof of claim \ref{cl_crocodile},
we find in $A_*$ an idempotent $e$ such that 
$e A^\wedge=I A^\wedge$ and $e A[f^{-1}]=I[f^{-1}]$,
and deduce that $I\subset e A$. Let $M:=(e A)/I$;
on one hand we have $M\otimes_AA^\wedge=0$. On the other hand,
$M[f^{-1}]=0$, {\em i.e.} $M=\bigcup_{n\in\N}\Ann_M(f^n)$,
which implies that $M\otimes_AA^\wedge=M$. Hence $M=0$, and
the claim follows.
\end{pfclaim}

\begin{claim}\label{cl_complete.is.OK} 
The theorem holds in case $A$ is $f$-adically complete.
\end{claim}
\begin{pfclaim} Say that $\Img(I\to A/fA)=\bar e_1\cdot(A/fA)$,
for an idempotent $\bar e_1\in(A/fA)_*$. For every $n>0$ set
$\bar I_n:=\Img(I\to A/f^nA)$; applying proposition 
\ref{prop_lift.idemp} to the non-unitary extensions 
$$
0\to (f^nA/f^{n+1}A)\cap\bar I_{n+1}\to
\bar I_{n+1}\to\bar I_n\to 0
$$
we construct recursively a compatible system of idempotents 
$\bar e_n\in\bar I_{n*}$, and by claim \ref{cl_invocation},
$\bar e_n$ generates $\bar I_n$ for every $n>0$; whence an 
element $e\in\liminv{n\in\N}\,(A/f^nA)_*\simeq
(\liminv{n\in\N}\,A/f^nA)_*\simeq A_*$. Clearly $e$ is
again an idempotent, and it induces decompositions
$A\simeq A_0\times A_1$, $I\simeq(I A_0)\times(I A_1)$,
such that $e=0$ (resp. $e=1$) in $A_0$ (resp. in $A_1$).
We can thus assume that either $A=A_0$ or $A=A_1$. If
$A=A_0$, then $\bar e_n=0$ for every $n>0$, so 
$I\subset\bigcap_{n>0}f^nA=0$. If $A=A_1$, then 
$\bar e_n=1$ for every $n>0$, {\em i.e.} $\bar I_n=A/f^nA$
for every $n>0$, that is, $I$ is dense in $A$. It follows 
easily that $\Ann_A(I)=0$. Using lemma \ref{lem_annihilate.and.forget}
we deduce $\Ann_{A[f^{-1}]}(I[f^{-1}])=0$; since $I[f^{-1}]$
is generated by an idempotent, this means that 
$I[f^{-1}]=A[f^{-1}]$. Set $M:=A/I$; the foregoing shows
that $M$ fulfills the hypotheses of lemma 
\ref{lem_vanishing.critter}(ii), whence $M=0$, which implies 
the claim.
\end{pfclaim}

\begin{claim}\label{cl_non.0.div.OK} 
The theorem holds in case $f$ is a non-zero-divisor
of $A_*$.
\end{claim}
\begin{pfclaim} According to claim \ref{cl_situation}, it suffices 
to prove that the ideal $I A^\wedge\subset A^\wedge$ is generated
by an idempotent. However, since 
$A^\wedge/fA^\wedge\simeq A/fA$, this follows directly from
our assumptions and from claim \ref{cl_complete.is.OK}.
\end{pfclaim}

After these preparations, we are ready to prove the theorem.
Let $T:=\bigcup_{n>0}\Ann_A(f^n)$; by claim \ref{cl_non.0.div.OK},
we know already that the ideal $I\cdot(A/T)$ is generated by an 
idempotent, and the same holds for the ideal $I A^\wedge$.
Since $T\otimes_AA^\wedge\simeq T$, it follows that the
natural morphism $\phi:T\to T A^\wedge$ is an epimorphism;
one verifies easily that $J:=\Ker\,\phi=T\cap(\bigcap_{n>0}f^nA)$, 
and then it is clear that $J^2=0$. Furthermore, it follows that 
the natural commutative diagram
$$
\xymatrix{ A/J \ar[r] \ar[d] & A/T \ar[d] \\
A^\wedge \ar[r] & A^\wedge/T A^\wedge
}$$
is cartesian. Taking claim \ref{cl_crocodile} into account, we 
deduce that the image of $I$ in $A/J$ is generated by an idempotent.
Lastly, we invoke claim \ref{cl_invocation} to show that $I$ is indeed
generated by an idempotent.
\end{proof}

\begin{remark} One may wonder whether every almost finitely
generated ideal $I\subset A$ such that $I=I^2$, is generated by
an idempotent. We do not know the answer to this question.
\end{remark}

\begin{theorem}\label{th_unramif.critter} 
Let $\phi:A\to B$ be a morphism of $V^a$-algebras,
$I,J\subset A$ any two ideals and $f\in A_*$ any almost element. 
\begin{enumerate}
\item
If $\phi\otimes_A\one_{A/I}$ and $\phi\otimes_A\one_{A/J}$ are
weakly unramified, then the same holds for $\phi\otimes_A\one_{A/IJ}$.
\item
If $\phi\otimes_A\one_{A/I}$ and $\phi\otimes_A\one_{A/J}$ are
unramified, then the same holds for $\phi\otimes_A\one_{A/IJ}$.
\item
If $\phi$ is flat and both $\phi\otimes_A\one_{A/fA}$ and
$\phi\otimes_A\one_{A[f^{-1}]}$ are weakly unramified, 
then $\phi$ is weakly unramified.
\item
If $\phi$ is almost finite  and both $\phi\otimes_A\one_{A/fA}$ and 
$\phi\otimes_A\one_{A[f^{-1}]}$ are unramified, then $\phi$ is
unramified.
\end{enumerate}
\end{theorem}
\begin{proof} To start out, we remark the following:
\begin{claim}\label{cl_base.ch.Tor} 
Let $A'$ be any $A$-algebra and set $B':=B\otimes_AA'$;
then for every $B'\otimes_{A'}B'$-module $M$ we have a natural
isomorphism:
$\Tor_1^{B\otimes_AB}(B,M)\simeq\Tor_1^{B'\otimes_{A'}B'}(B',M)$.
\end{claim}
\begin{pfclaim} Notice that the short exact sequence of
$B\otimes_AB$-modules
\set\begin{equation}\label{eq_stays.split}
0\to I_{B/A}\to B\otimes_AB\to B\to 0.
\end{equation}
is split exact as a sequence of $B$-modules (and {\em a fortiori\/}
as a sequence of $A$-modules); hence 
\eqref{eq_stays.split}\,$\otimes_AN$ is again exact for every
$A$-module $N$. We deduce easily that 
$\Tor_1^{B\otimes_AB}(B,N\otimes_AB\otimes_AB)=0$ for every
$A$-module $N$. On the other hand, the morphism 
$B\otimes_AB\to B'\otimes_{A'}B'$ determines a base change
spectral sequence (cp. \cite[Th.5.6.6]{We})
$$
E^2_{pq}:=
\Tor_p^{B'\otimes_{A'}B'}(\Tor_q^{B\otimes_AB}(B,B'\otimes_{A'}B'),M)
\Rightarrow\Tor^{B\otimes_AB}_{p+q}(B,M)
$$
for every $B'\otimes_AB'$-module $M$. The foregoing yields
$E^2_{p1}=0$ for every $p\in\N$, so the claim follows.
\end{pfclaim}

In order to show (i) we may replace $A$ by $A/IJ$ and therefore
assume that $IJ=0$; then we have to prove that
$\Tor_1^{B\otimes_AB}(B,N)=0$ for every $B\otimes_AB$-module
$N$. By a simple devissage we can further reduce to
the case where either $IN=0$ or $JN=0$.
Say $IN=0$; we apply claim \ref{cl_base.ch.Tor} with 
$A':=A/I$, $B':=B/IB$ to deduce 
$\Tor_1^{B\otimes_AB}(B,N)\simeq\Tor_1^{B'\otimes_{A'}B'}(B',N)$; 
however the latter module vanishes by our assumption on 
$\phi\otimes_A\one_{B'}$.

For (ii) we can again assume that $IJ=0$. We need to show
that $B$ is an almost projective $B\otimes_AB$-module, and
we know already that it is flat by (i). Set 
$I':=I(B\otimes_AB)$, $J':=J(B\otimes_AB)$; by assumption
$B/IB$ (resp. $B/JB$) is almost projective over $B\otimes_AB/I'$
(resp. $B\otimes_AB/J'$); in light of proposition 
\ref{prop_equiv.2-prod} we deduce that $B/(IB\cap JB)$ is almost 
projective over $B\otimes_AB/(I'\cap J')$. Then, since 
$(I'\cap J')^2\subset I'J'=0$, the assertion follows from lemma
\ref{lem_flat.proj}(i).

Next we remark that, under the assumptions of (iii), the 
hypotheses of lemma \ref{lem_Tor.orgy} are fulfilled by 
$C:=B\otimes_AB$ and $M:=B$, so assertion (iii) holds. 
Finally, suppose that $\phi$ is almost finite; hence
$\Omega_{B/A}$ is an almost finitely generated $A$-module, 
and by assumption we have 
$\Omega_{B/A}\otimes_A(A/fA)=\Omega_{B/A}[f^{-1}]=0$, therefore
$\Omega_{B/A}=0$, in view of lemma \ref{lem_vanishing.critter}(ii).
Set $I:=\Ker(B\otimes_AB\to B)$; it follows that $I^2=I$,
and $I$ is almost finitely generated as an $A$-module, so
{\em a fortiori\/} as a $B\otimes_AB$-module. Moreover, 
from proposition \ref{prop_idemp} and our assumption on 
$\phi\otimes_A\one_{A/fA}$ and $\phi[f^{-1}]$, it follows
that both $I[f^{-1}]\subset B\otimes_AB[f^{-1}]$ and
$\Img(I\to B\otimes_A(B/fB))$ are generated by idempotents.
Then theorem \ref{th_gen.by.idemps} says that also $I$ 
is generated by an idempotent of $(B\otimes_AB)_*$, so that
$B$ is almost projective as a $B\otimes_AB$-module, by
remark \ref{rem_idemp}, which proves (iv).
\end{proof}

\begin{corollary}\label{cor_oftheabove}
Let $I\subset\rad(A)$ be a tight ideal. If $B$ is an almost
finite $A$-algebra and $B/IB$ is unramified over $A/I$, then
$B$ is unramified over $A$.
\end{corollary}
\begin{proof} Under the stated assumptions, $\Omega_{B/A}$
is an almost finitely generated $A$-module such that 
$\Omega_{B/A}\otimes_AA/I=0$; by lemma \ref{lem_Naka} it
follows that 
\set\begin{equation}\label{eq_omega.vanish}
\Omega_{B/A}=0.
\end{equation} 
Set 
$I_{B/A}:=\Ker(\mu_{B/A}:B\otimes_AB\to B)$; we derive that 
$I_{B/A}=I^2_{B/A}$.
Pick $n>0$ and a finitely generated subideal 
$\fm_0\subset\fm$ with $I^n\subset\fm_0 A$; let
$\eps_1,...,\eps_k$ be a set of generators for $\fm_0$.
\begin{claim}\label{cl_when.invert.eps} $I_{B/A}[\eps_i^{-1}]$ 
is generated by an idempotent of $B\otimes_AB[\eps_i^{-1}]$, 
for every $i\leq k$.
\end{claim}
\begin{pfclaim} Notice that $A[\eps_i^{-1}]$ is a (usual) 
$V$-algebra (that is, the localization functor 
$A[\eps_i^{-1}]_*\Mod\to A[\eps_i^{-1}]\Mod$
is an equivalence) and $B[\eps_i^{-1}]$ is a finite 
$A[\eps_i^{-1}]$-algebra; from \eqref{eq_omega.vanish}
and \cite[Exp.I, Prop.3.1]{SGA1} we deduce that 
$B[\eps_i^{-1}]$ is unramified over $A[\eps_i^{-1}]$,
which implies the claim. 
\end{pfclaim}

\begin{claim}\label{cl_when.modout.eps} 
$I_{B/A}/\fm_0I_{B/A}$ is generated by an idempotent as well.
\end{claim}
\begin{pfclaim} We apply theorem \ref{th_unramif.critter}(ii) with 
$J=I^m$ to deduce, by induction on $m$, that $B/I^{m+1}B$
is an unramified $A/I^{m+1}$-algebra for every $m\in\N$;
it follows that $B/\fm_0B$ is unramified over $A/\fm_0A$,
whence the claim.
\end{pfclaim}

Using theorem \ref{th_unramif.critter}(iv) and claims
\ref{cl_when.invert.eps}, \ref{cl_when.modout.eps}, the 
assertion is now easily verified by induction on the number 
$k$ of generators of $\fm_0$,
\end{proof}

The following lemma shows that morphisms of unramified $A$-algebras
can be read off from their "graphs". This result will be useful in 
section \ref{sec_lift.hens}.

\sset\subsubsection{}\label{subsec_graph.rule}
Let $A\to B$ be an unramified morphism of $V^a$-algebras, 
and $C$ any $A$-algebra. Suppose that $e\in(B\otimes_AC)_*$
is an idempotent; to $e$ we can attach a morphism of $A$-algebras:
$$
\Gamma(e):C\to B\otimes_AC\to e\cdot(B\otimes_AC)\qquad 
c\mapsto 1\otimes c\mapsto e\cdot(1\otimes c)\quad (c\in C_*)
$$
and, in case $\Gamma(e)$ is an isomorphism, we obtain a morphism
$f_e:B\to C$ after composing $\Gamma(e)^{-1}$ and the natural 
morphism 
$$
\Delta(e):B\to B\otimes_AC\to e\cdot(B\otimes_AC)
\qquad 
b\mapsto b\otimes 1\mapsto e\cdot(b\otimes 1)\quad (b\in B_*)
$$ 
Conversely, To a given a morphism 
$\phi:B\to C$ of $A$-algebras, we can associate an idempotent 
$e_\phi\in(B\otimes_AC)_*$, by setting 
$e_\phi:=(\one_B\otimes_A\phi)(e_{B/A})$, where
$e_{B/A}$ is the idempotent provided by proposition 
\ref{prop_idemp}. Then, since the commutative diagram
$$
\xymatrix{ 
B\otimes_AB \ar[rr]^-{\one_B\otimes\phi} \ar[d]_{\mu_{B/A}} &&
B\otimes_AC \ar[d]^{\mu_{C/A}\circ(\phi\otimes\one_C)} \\
B \ar[rr]^-\phi && C
}$$
is cocartesian, we deduce from proposition \ref{prop_idemp}
that the corresponding morphism $\Gamma(e_\phi)$ is an isomorphism.
In fact, the following holds:

\begin{lemma}\label{lem_graph.rule} The rule $\phi\mapsto e_\phi$ 
of \eqref{subsec_graph.rule} induces a natural bijection from 
$\Hom_{A\Alg}(B~,C)$ onto the set of idempotents 
$e\in(B\otimes_AC)_*$ such that $\Gamma(e)$ is an isomorphism.
Its inverse is the rule $e\mapsto\phi_e$.
\end{lemma}
\begin{proof} Let $\phi:B\to C$ be given and set $e:=e_\phi$.
We need to show that $f_e=\phi$, or equivalently 
$\Delta(e)=\Gamma(e)\circ\phi$. This translates into the equality: 
\set\begin{equation}\label{eq_identity.e}
e\cdot(b\otimes 1)=e\cdot(1\otimes\phi(b))\quad 
\text{for every $b\in B_*$.}
\end{equation}
The latter can be rewritten as 
$(1\otimes\phi)(e_{B/A}\cdot(1\otimes b-b\otimes 1))=0$ which
follows from identity (iii) in proposition \ref{prop_idemp}. 
Conversely, for $e$ given such that
$\Gamma(e)$ is an isomorphism, set $\phi:=\phi_e$; we have 
to check that $e_\phi=e$. By construction of $\phi$ we have
$\Delta(e)=\Gamma(e)\circ\phi$, which means that \eqref{eq_identity.e}
holds for the pair $(e,\phi)$. Let $J\subset B\otimes_AC$
be the ideal generated by the elements of the form
$b\otimes 1-1\otimes\phi(b)$, for all $b\in B_*$; we know
that $eJ=0$, which means that the natural morphism
$B\otimes_AC\to e(B\otimes_AC)$ factors through the
natural morphism $B\otimes_AC\to e_\phi(B\otimes_AC)$.
Since both $\Gamma(e)$ and $\Gamma(e_\phi)$ are isomorphisms,
it follows that $e(B\otimes_AC)=e_\phi(B\otimes_AC)$,
so $e=e_\phi$.
\end{proof}

\subsection{Topological algebras and modules}\label{sec_lin-topol}
\begin{definition}\label{def_a-topologies}
\index{Almost module(s)!linear topology on an|indref{def_a-topologies}}
\index{Almost module(s)!linear topology on an!generated by a
family of submodules|indref{def_a-topologies}}
\index{Almost module(s)!linear topology on an!topological
closure of a submodule in a|indref{def_a-topologies}}
\index{Almost module(s)!linear topology on
an!complete|indref{def_a-topologies}}
\index{Almost module(s)!linear topology on
an!$I$-preadic, $I$-c-preadic, $I$-adic|indref{def_a-topologies}}
\index{Almost module(s)!linear topology on an!ideal of definition
of a|indref{def_a-topologies}}
Let $A$ be a $V^a$-algebra, $M$ an $A$-module.
\begin{enumerate}
\item
A {\em linear topology\/} on $M$ is a non-empty family $\cL$ of 
submodules of $M$ that satisfies the following conditions. If 
$I,J\in\cL$, then $I\cap L\in\cL$; if $I\in\cL$ and $I\subset J$, then
$J\in\cL$. We say that a submodule $I$ is {\em open\/} if $I\in\cL$.
\item
Let $\cF$ be a family of submodules of $M$. The {\em topology
generated by $\cF$\/} is the smallest linear topology containing $\cF$.
\item
Let $\cL$ be a linear topology on $M$ and $I\subset M$ a submodule. The 
{\em closure\/} of $I$ is the submodule $\bar I:=\bigcap_{J\in\cL}(I+J)$.
\item
We say that the topology $\cL$ is {\em complete} if the 
natural morphism $M\to\liminv{I\in\cL}\,M/I$ is an isomorphism.
\item
Let $\cL_M$, $\cL_A$ be topologies on $M$, resp. on $A$, and $I$ 
an open ideal of $A$. We say that the topology $\cL_M$ is 
{\em $I$-preadic\/} (resp. {\em $I$-c-preadic\/})
if the family of submodules $(I^n M~|~n\in\N)$ (resp.
$(\bar{\bar{I^n}\cdot M}~|~n\in\N)$) is a cofinal subfamily in $\cL_M$.
In either case, we say that $I$ is an {\em ideal of definition\/}
for $\cL_M$. Furthermore, we say that the topology $\cL_M$ is 
{\em $I$-adic\/} if it is $I$-preadic and complete. We similarly 
define an {\em $I$-c-adic\/} topology.
\end{enumerate}
\end{definition}

\sset\subsubsection{}\label{subsec_continuous-etc}
\index{Almost module(s)!linear topology on an!continuous morphism
of|indref{subsec_continuous-etc}}
\index{Almost module(s)!linear topology on an!induced topology
on a submodule|indref{subsec_continuous-etc}}
\index{Almost module(s)!linear topology on
an!completion of a|indref{subsec_continuous-etc}}
One defines as usual the notions of {\em continuous\/}, resp. 
{\em open\/} morphism of $A$-modules, of {\em induced
topology\/} on submodules of a linearly topologized module
and of {\em completion\/} of a topological module. 
Notice that, if $A$ is an $I$-c-preadic topological 
$V^a$-algebra and $M$ an $I$-c-preadic $A$-module 
(for some open ideal $I\subset A$) then we have 
$\bar{I^n M}=\bar{\bar{I^n}\cdot M}$ for every $n\in\N$.

\begin{lemma}\label{lem_from.gr.to.fil} Let $M$, $N$ be 
two $A$-modules endowed with descending filtrations 
$\cF_M:=(M_n~|~n\in\N)$, $\cF_N:=(N_n~|~n\in\N)$ 
so that $M_0=M$, $N_0=N$ and $M$ (resp. $N$) is complete
for the topology generated by $\cF_M$ (resp. by $\cF_N$).
Suppose furthermore that $\phi:M\to N$ is a morphism of
$A$-modules which respect the filtrations, and such that
the induced morphisms $\gr^i\phi:\gr^iM\to\gr^iN$ are epimorphisms
for every $i\in\N$. Then $\phi$ is an epimorphism.
\end{lemma}
\begin{proof} Under the assumptions, the induced morphism
$\phi_n:M/M_n\to N/N_n$ is an epimorphism for every $n\in\N$.
Clearly it suffices to show then that $\liminv{n\in\N}^1\,\Ker(\phi_n)$
vanishes. However, from $\Coker(\gr^i\phi)=0$ one deduces by snake 
lemma that the induced morphism $\Ker(\phi_{i+1})\to\Ker(\phi_i)$
is an epimorphism, so the claim follows.
\end{proof}

\sset\subsubsection{}\label{subsec_A_infty}
\index{Almost algebra(s)!projective topology on a
complete linearly topologized|indref{subsec_A_infty}}
Let $(A_n;~\phi_n:A_{n+1}\to A_n~|~n\in\N)$ be a projective 
system of\/ $V^a$-algebras. Then $A_\infty:=\liminv{n\in\N}\,A_n$ 
is naturally endowed with a linear topology, namely the topology 
$\cL$ generated by the family $(\Ker(A_\infty\to A_n)~|~n\in\N)$. 
We call $\cL$  the {\em projective topology\/} defining $A_\infty$.

\begin{lemma}\label{lem_proj.lim.tops}  
Let $(A_n;~\phi_n:A_{n+1}\to A_n~|~n\in\N)$ be a projective system 
of\/ $V^a$-algebras, $A_\infty$ its projective limit, and suppose that:
\begin{enumerate}
\renewcommand{\labelenumi}{(\alph{enumi})}
\item
$\phi_n$ is an epimorphism for every $n\in\N$.
\item
$\Ker\,\phi_n=\Ker(\phi_0\circ...\circ\phi_n:A_{n+1}\to A_0)^{n+1}$
for every $n\in\N$.
\renewcommand{\labelenumi}{(\roman{enumi})}
\end{enumerate}
Then we have:
\begin{enumerate}
\item
Let $I:=\Ker(A_\infty\to A_0)$. The projective topology
of $A_\infty$ is $I$-c-preadic. This topology
is complete and moreover $\bar{I^{n+1}}=\Ker(A_\infty\to A_n)$ 
for every $n\in\N$.
\item
Conversely, suppose that $A$ is any complete linearly topologized 
$V^a$-algebra and $I\subset A$ an open ideal such that the topology 
of $A$ is $I$-c-preadic. 
Set $A_n:=A/\bar{I^{n+1}}$ for every $n\in\N$. Then the projective 
system $(A_n~|~n\in\N)$ satisfies conditions {\em (a)} and {\em (b)}.
\end{enumerate}
\end{lemma}
\begin{proof} To show (i) it suffices to show that 
$\bar{I^{n+1}}=\Ker(A_\infty\to A_n)$ for every $n\in\N$.
We endow $\bar{I^{n+1}}$ (resp. $\Ker(A_\infty\to A_n)$) with the 
descending filtration $F_i:=\bar{I^{n+1+j}}$  
(resp. $G_i:=\Ker(A_\infty\to A_{n+j})$) ($i\in\N$).
Notice first that the natural morphism $A_\infty\to A_{k+1}$
induces an isomorphism 
$I/G_{k+1}\stackrel{\sim}{\to}\Ker(A_{k+1}\to A_0)$.
By (b) we deduce isomorphisms 
\set\begin{equation}\label{eq_graded.ideals}
(I^{k+1}+G_{k+1})/G_{k+1}\stackrel{\sim}{\to}\Ker\,\phi_k
\qquad\text{for every $k\in\N$.}
\end{equation}
Especially, $I^{k+1}\subset\Ker(A_\infty\to A_k)$. Since 
the latter is a closed ideal, we deduce 
$\bar{I^{k+1}}\subset\Ker(A_\infty\to A_k)$ for every $k\in\N$.
Therefore, we obtain an imbedding 
$\bar{I^{n+1}}\subset\Ker(A_\infty\to A_n)$ that respects the
filtrations $F_\bullet$ and $G_\bullet$. Since 
$\gr^iG_\bullet\simeq\Ker\:\phi_{i+n}$, we derive
easily from \eqref{eq_graded.ideals} that the induced morphism
$\gr^iF_\bullet\to\gr^iG_\bullet$ is an epimorphism for every
$i\in\N$, so (i) follows by lemma \ref{lem_from.gr.to.fil}.
Under the assumptions of (ii) it is obvious that condition
(a) holds. To show (b) comes down to verifying the identity
$I^{n+1}+\bar{I^{n+2}}=\bar{I^{n+1}}$, which is obvious
since $\bar{I^{n+2}}$ is an open ideal.
\end{proof}

\begin{remark} Notice that, in the situation of lemma 
\ref{lem_proj.lim.tops}, the pair $(A_\infty,I)$ is henselian.
One sees this as follows, using the criterion of remark 
\ref{rem_hensel}(iii). Let $R_n:=A_*/\bar{I^{n+1}}_*$,
$J_n:=\bar{I^n}_*/\bar{I^{n+1}}_*$ for every $n\in\N$;
let also $R:=(\liminv{n\in\N}\,R_n)$ and denote by $J$ the 
image of $\liminv{n\in\N}\,J_n$ in $R$. Since 
$J_n$ is a nilpotent ideal of $R_n$, the pair $(R_n,J_n)$
is henselian. Let now $p(X)\in R[X]$ and denote by $p_n(X)$
the image of $p(X)$ in $R_n(X)$, for every $n\in\N$; suppose 
that $p_0(X)$  decomposes in $R_0[X]$ as 
$p_0(X)=q_0(X)\cdot r_0(X)$; by the cited criterion we derive 
a compatible family of decompositions
$p_n(X)=q_n(X)\cdot r_n(X)$ in $R_n[X]$ for every $n\in\N$; 
hence $p(X)$ admits a similar decomposition in $R[X]$,
so $(R,J)$ is a henselian pair, which implies the contention,
in view of lemma \ref{lem_hensel}.
\end{remark}

\begin{lemma}\label{lem_when.J.open} Let $A_\infty$ be the 
projective limit of a system of $V^a$-algebras satisfying 
conditions {\em (a)} and {\em (b)} of lemma {\em\ref{lem_proj.lim.tops}}.
Let $J\subset A_\infty$ be a finitely generated ideal
such that $\bar J$ is open. Then:
\begin{enumerate}
\item
 $J=\bar J$. 
\item
Every finite system of generators $x_1,...,x_r$ of $J$ defines 
an open morphism $A^r_\infty\to J$ (for the topologies on 
$A_\infty^r$ and $J$ induced from $A_\infty$).
\item
$J^k$ is open for every $k\in\N$.
\end{enumerate}
\end{lemma}
\begin{proof} By lemma \ref{lem_proj.lim.tops} the topology
on $A_\infty$ is c-preadic for some defining ideal $I$. The 
assumption means that $\bar{I^n}\subset J$ for a sufficiently 
large integer $n$. 
\begin{claim}\label{cl_proj.lim.tops} The natural morphism:
$(J\cap\bar{I^m})/(J\cap\bar{I^{m+1}})\to\bar{I^m}/\bar{I^{m+1}}$
\ is an isomorphism for every $m\geq n$.
\end{claim}
\begin{pfclaim} It is an easy verification that shall be left to
the reader.
\end{pfclaim}

Let $x_1,...,x_r\in J_*$ be a system of generators for $J$, and 
$f:A^r_\infty\to J$ the corresponding epimorphism. Set 
$M_k:=f^{-1}(J\cap\bar{I^k})$. To start out, 
$M_k$ contains $(\bar{I^k})^{\oplus r}$, whence $M_k$ is open in
$A_\infty^r$ for every $k\in\N$. This already shows that $f$
is continuous. We define a descending filtration on $M_n$ 
by setting $F_k:=\bar{I^k M_n}$ for every $k\in\N$. 
In view of claim \ref{cl_proj.lim.tops}, the morphism $f$
induces an epimorphism $M_n\to\bar{I^n}/\bar{I^{n+1}}$, hence
for every $k\in\N$, the morphisms:
\set\begin{equation}\label{eq_epimaps}
\bar{I^k}\cdot M_n\to
(\bar{I^k}\cdot\bar{I^n})/(\bar{I^k}\cdot\bar{I^{n+1}})\to
\bar{I^{n+k}}/\bar{I^{n+k+1}}
\end{equation}
are epimorphisms as well. The composition of the morphisms
\eqref{eq_epimaps} extends to an epimorphism 
$F_k\to\bar{I^{n+k}}/\bar{I^{n+k+1}}$, for every $k\in\N$.
In other words, if we endow $\bar{I^n}$ with its $I$-c-preadic 
filtration, then $f$ induces a morphism $\phi:M_n\to\bar{I^n}$ 
of filtered modules, such that $\gr^\bullet\phi$ is an epimorphism
of graded modules. By assumption $\bar{I^n}$ is complete for its 
filtration. Similarly, $M_n$ is complete for its filtration 
$F_\bullet$; indeed, this follows by remarking that 
$(\bar{I^{n+k}})^{\oplus r}\subset F_k\subset(\bar{I^k})^{\oplus r}$ 
for every $k\in\N$. Then lemma \ref{lem_from.gr.to.fil} says
that $\phi$ is an epimorphism, {\em i.e.} $f(M_n)=\bar{I^n}$;
it follows that $J$ is an open ideal, hence equal to its closure.
More generally, the same argument proves that for every $k\in\N$, 
$f(F_k)=\bar{I^{n+k}}$. Since by the foregoing, $(F_k~|~k\in\N)$
is a cofinal system of open submodules of $A_\infty^r$, we deduce
that $f$ is an open morphism. Finally the identity 
$$
f((J^n)^{\oplus r})=J^{n+1}\qquad \text{for every $n\in\N$}
$$
together with (ii) and an easy induction, yields (iii).
\end{proof}

\sset\subsubsection{}
Let $I$ be an ideal of definition for the projective topology 
on a ring $A_\infty$ as in lemma \ref{lem_proj.lim.tops}(i).
We wish now to give a criterion to ensure that the projective
topology on $A_\infty$ is actually $I$-adic. This can be achieved
in case $I$ is tight, as shown by proposition
\ref{prop_fin.gen.subideal}, which generalizes
\cite[Ch.0, Prop.7.2.7]{EGAI}.

\begin{proposition}\label{prop_fin.gen.subideal} 
Let $A_\infty$ be as in lemma {\em\ref{lem_when.J.open}}
and $I$ an ideal of definition for the projective topology of
$A_\infty$. Assume that $I$ is tight and moreover that $I/\bar{I^2}$ 
is an almost finitely generated $A_\infty$-module. Then :
\begin{enumerate}
\item
There exists $n\in\N$ and a finitely generated ideal 
$J\subset A_\infty$ such that $I^n\subset J\subset I$.
\item
The topology of $A_\infty$ is $I$-adic.
\end{enumerate}
\end{proposition}
\begin{proof} Using the natural epimorphism 
$(I/\bar{I^2})^{\otimes k}\to\bar{I^k}/\bar{I^{k+1}}$ we deduce that 
$\bar{I^k}/\bar{I^{k+1}}$ is almost finitely generated for every 
$k\geq 0$; then the same holds for $I/\bar{I^{k+1}}$. Let 
$\fm_0\subset\fm$ be a finitely generated subideal and $n\in\N$
such that $I^n\subset\fm_0$; pick a finitely generated ideal 
$\fm_1\subset\fm$ with $\fm_0\subset\fm_1^2$; we can find a 
finitely generated $A_\infty$-module $Q\subset I/\bar{I^{n+2}}$
such that $\fm_1\cdot(I/\bar{I^{n+2}})\subset Q$; up to replacing
$Q$ by $\fm_1\cdot Q$, we can then achieve that $Q$ is generated
by the images of finitely many almost elements $x_1,...,x_t$ of 
$I$ and moreover $\fm_0\cdot(I/\bar{I^{n+2}})\subset Q$. Then,
since $\bar{I^{n+2}}$ is open we deduce 
\set\begin{equation}\label{eq_include.in.Q}
\bar{I^{n+1}}/\bar{I^{n+2}}\subset Q.
\end{equation} 
Let $J\subset A_\infty$ be the ideal generated by $x_1,...,x_t$. 
From \eqref{eq_include.in.Q} we deduce that the natural morphism
$(J\cap\bar{I^{n+1}})/(J\cap\bar{I^{n+2}})\to
\bar{I^{n+1}}/\bar{I^{n+2}}$ is an isomorphism. Since
$I^k(J\cap\bar{I^{n+1}})\subset J\cap\bar{I^{n+1+k}}$,
the same holds more generally for the morphisms 
$(J\cap\bar{I^{n+1+k}})/(J\cap\bar{I^{n+2+k}})\to
\bar{I^{n+1+k}}/\bar{I^{n+2+k}}$, for every $k\in\N$.
This easily implies that $J':=J\cap\bar{I^{n+1}}$ is {\em dense\/}
in $\bar{I^{n+1}}$, {\em i.e.\/} $\bar{J'}=\bar{I^{n+1}}$;
especially, $\bar J$ is open in $A_\infty$ and therefore
$J=\bar J$ by lemma \ref{lem_when.J.open}(i), which proves 
assertion (i). Assertion (ii) follows easily from (i) and lemma 
\ref{lem_when.J.open}(iii).
\end{proof}

\sset\subsubsection{}\label{subsec_2-limit}
\index{$A_\infty\Mod_\mathrm{top}$ : topological modules on a
complete almost algebra|indref{subsec_2-limit}}
Let now $(A_n~|~n\in\N)$ be an inverse system of $V^a$-algebras
satisfying conditions (a) and (b) of lemma \ref{lem_proj.lim.tops}
and let $I:=\Ker(A_\infty\to A_0)$.
The induced functors $A_{n+1}\Mod\to A_n\Mod$ : 
$M\mapsto A_n\otimes_{A_{n+1}}M$ define an inverse system of
categories $(A_n\Mod~|~n\in\N)$. In such situation one can
define a natural functor
\set\begin{equation}\label{eq_2-limit}
A_\infty\Mod_\mathrm{top}\to \text{2-}\liminv{n\in\N}~A_n\Mod
\end{equation}
from the category of topological $A_\infty$-modules whose topology
is $I$-c-adic (and of continuous $A_\infty$-linear morphisms) to 
the 2-limit of the foregoing inverse system of categories
(see \cite[Ch.I]{Hak} for generalities on 2-categories).
Namely, let $I:=\Ker(A_\infty\to A_0)$; then to an $A_\infty$-module
$M$ one associates the compatible system 
$(M/\bar{I^{n+1} M}~|~n\in\N)$.

\begin{lemma}\label{lem_2-limit} 
The functor \eqref{eq_2-limit} is an equivalence of categories.
\end{lemma}
\begin{proof} We claim that a quasi-inverse to \eqref{eq_2-limit}
can be given by associating to any compatible system $(M_n~|~n\in\N)$
the $A_\infty$-module $M_\infty:=\liminv{n\in\N}M_n$, with the
linear topology generated by the submodules 
$K_n:=\Ker(M_\infty\to M_n)$, for all $n\in\N$. In order
to show that this topology on $M_\infty$ is $I$-c-adic
(where $I:=\Ker(A_\infty\to A_0)$), it suffices to verify
that $K_n=\bar{I^{n+1} M_\infty}$ for every $n\in\N$.
Clearly $I^{n+1} M_\infty\subset K_n$, hence we are
reduced to showing that $K_n\subset K_m+I^{n+1} M_\infty$
for every $m>n$ or equivalently, that $K_n/K_m$ equals the
image of $I^{n+1} M_\infty$ inside $M_m$, which is obvious.
\end{proof}

In the following we will seek conditions under which
\eqref{eq_2-limit} can be refined to equivalences between
interesting subcategories, for instance to almost finitely
generated, or almost projective modules.

\begin{lemma}\label{lem_about.N} 
Let $M_\infty$ be a topological $A_\infty$-module 
whose topology is $I$-c-adic (for some ideal of definition 
$I\subset A_\infty$). Suppose that $N\subset M_\infty$ is
a finitely generated submodule such that $\bar N$ is open
in $M_\infty$. Then :
\begin{enumerate}
\item
$N=\bar N$.
\item
Every finite set of generators $x_1,...,x_r$ of $N$ determines
an open morphism $A_\infty^r\to N$.
\item
For every open ideal $J\subset A_\infty$, the submodule $JN$
is open in $M_\infty$.
\end{enumerate}
\end{lemma}
\begin{proof} {\em Mutatis mutandis}, this is the same as the
proof of lemma \ref{lem_when.J.open}; the details can thus be 
safely entrusted to the reader.
\end{proof}

\begin{lemma}\label{lem_fin.gen.at.infty} 
Keep the notation of \eqref{subsec_2-limit} and
assume that $I$ is a tight ideal. 
Let $M$ be a topological $A_\infty$-module whose topology 
is $I$-c-adic and let $(M_n~|~n\in\N)$ be the image of $M$ 
under \eqref{eq_2-limit}. Suppose also that $M_0$ is an 
almost finitely generated $A_0$-module. Then:
\begin{enumerate}
\item
$M$ admits an open finitely generated submodule.
\item
$M$ is almost finitely generated.
\item
$M\otimes_{A_\infty}A_n\simeq M_n$ for every $n\in\N$.
\end{enumerate}
\end{lemma}
\begin{proof} Using lemma \ref{lem_flat.proj}(ii) and the
hypothesis on $M_0$ we deduce that $M_n$ is almost finitely
generated for every $n\in\N$. Choose finitely generated 
$\fm_0\subset\fm$ and $n\geq 0$ such that $I^n\subset\fm_0\cdot A$; 
it follows easily that there are
almost elements $x_1,...,x_t$ of $M$ whose images in $M_{n+1}$
generate a submodule $Q$ such that $\fm_0\cdot M_{n+1}\subset Q$.
This implies that 
$\bar{I^{n+1} M}/\bar{I^{n+2} M}\subset Q$. Let
$N\subset M$ be the submodule generated by $x_1,...,x_t$. Then, 
arguing as in the proof of proposition \ref{prop_fin.gen.subideal} 
we see that $\bar{I^{n+1} M}\subset\bar N$. Hence $N=\bar N$
in view of lemma \ref{lem_about.N}(i), so (i) holds.
Furthermore, by lemma \ref{lem_about.N}(iii), $\bar{I^n}\cdot M$
is open for every $n\in\N$ (since it contains $\bar{I^n}\cdot N$),
so (iii) follows easily. Finally, $M$ is almost finitely generated
because the natural morphism $M_n\to M/N$ is an epimorphism.
\end{proof}

\sset\subsubsection{}\label{subsec_infinite.proj}
In the situation of \eqref{subsec_2-limit}, let $M$ be 
an $A_\infty$-module with an $I$-c-adic topology, and let 
$(M_n~|~n\in\N)$ be its image under \eqref{eq_2-limit}. 
We suppose now that $M_n$ is an almost projective $A_n$-module 
for every $n\in\N$. Let us choose an epimorphism 
$f:A^{(S)}_\infty\to M$. We endow the free 
$A_\infty$-module $A^{(S)}_\infty$ with the topology generated
by the family of submodules  
$(\bar{I^n}\cdot A^{(S)}_\infty~|~n\in\N)$. Notice that $f$ is 
continuous for this topology, hence it extends to a continuous
morphism $f^\wedge:F_\infty\to M$ on the completion $F_\infty$
of $A^{(S)}_\infty$.

\begin{lemma}\label{lem_top.split} 
With the notation of \eqref{subsec_infinite.proj}, we have:
\begin{enumerate}
\item
$f^\wedge$ is {\em topologically almost split}, {\em i.e.},
for every $\eps\in\fm$ there exists a continuous morphism 
$g:M\to F_\infty$ such that $f^\wedge\circ g=\eps\cdot\one_M$.
\item
If $M$ is almost finitely generated, then $M$ is almost 
projective.
\end{enumerate}
\end{lemma}
\begin{proof} By lemma \ref{lem_2-limit} the set of continuous
morphisms $M\to F_\infty$ is in natural bijection with 
$\liminv{n\in\N}~\Hom_{A_n}(M_n,F_n)$, where 
$F_n:=F_\infty/\bar{I^{n+1} F_\infty}$ for every $n\in\N$.
On the other hand, under the assumptions of the lemma, $f$
induces an epimorphism:
$$
\liminv{n\in\N}~\Alhom_{A_n}(M_n,F_n)\to
\liminv{n\in\N}~\Alhom_{A_n}(M_n,M_n).
$$ 
Assertion (i) follows easily. Suppose now that $M$ is almost
finitely generated; then for every finitely generated subideal
$\fm_0\subset\fm$ we can find $m\geq 0$ and a morphism 
$h:A^m_\infty\to M$ such that $\fm_0\cdot\Coker\,h=0$. Let 
$\eps_1,...,\eps_t$ be a system of generators for $\fm_0$;
then for every $i,j\leq t$ we can pick a morphism 
$\phi:A^{(S)}_\infty\to A^m_\infty$ such that 
$h\circ\phi=\eps_i\cdot\eps_j\cdot f$; after taking completions,
this relation becomes 
$h\circ\phi^\wedge=\eps_i\cdot\eps_j\cdot f^\wedge$.
Choose $g:M\to F_\infty$ as in (i); we deduce 
$\eps_i\cdot\eps_j\cdot\eps\cdot\one_M=
h\circ\phi^\wedge\circ g$, which as usual shows that 
$\eps_i\cdot\eps_j\cdot\eps$ annihilates 
$\AlExt^1_{A_\infty}(M,N)$ for every $A_\infty$-module $N$,
so (ii) holds.
\end{proof}

\sset\subsubsection{}\label{subsec_in.the.situation}
For a $V^a$-algebra $A$, let us denote by $A\Mod_\mathrm{afpr}$
the full subcategory of $A\Mod$ consisting of all almost finitely
generated projective $A$-modules. Let now $(A_n~|~n\in\N)$ and
$I\subset A_\infty$ be as in \eqref{subsec_2-limit}. We define 
natural functors
\set\begin{equation}\label{eq_yet.another.fct}
A_\infty\Mod_\mathrm{afpr}\to
\text{2-}\liminv{n\in\N}~A_n\Mod_\mathrm{afpr}
\end{equation}
respectively (see \eqref{subsec_defafp}),
\set\begin{equation}\label{eq_yet.one.fct.etale}
A_\infty\Et_\mathrm{afp}\to
A_0\Et_\mathrm{afp}
\end{equation}
by assigning to a given object $M$ of $A_\infty\Mod_\mathrm{afpr}$,
the compatible system $(M\otimes_{A_\infty}A_n~|~n\in\N)$
(resp. to an object $B$ of $A_\infty\Et_\mathrm{afp}$, the
$A_0$-algebra $B\otimes_{A_\infty}A_0$).

\begin{theorem}\label{th_for.A-mods} In the situation of 
\eqref{subsec_in.the.situation}, suppose that $I$ is tight.
Then \eqref{eq_yet.another.fct} is an equivalence of categories.
\end{theorem}
\begin{proof} We claim that a quasi-inverse to
\eqref{eq_yet.another.fct} is obtained by assigning to
any compatible system $(M_n~|~n\in\N)$ the $A_\infty$-module
$M:=\liminv{n\in\N}~M_n$. Indeed, by the proof of lemma 
\ref{lem_2-limit} there results that $M$ is endowed with
a natural $I$-c-adic topology; then by lemma 
\ref{lem_fin.gen.at.infty}(ii) we deduce that $M$ is almost 
finitely generated, and finally lemma \ref{lem_top.split}(ii) 
says that $M$ is almost projective. From lemma 
\ref{lem_fin.gen.at.infty}(iii) it follows that the functor
thus defined is a right quasi-inverse to \eqref{eq_yet.another.fct},
so the latter is essentially surjective. Full faithfulness
is a consequence of the following:
\begin{claim}\label{cl_proj.on.P} 
Let $P$ be an almost finitely generated projective
$A_\infty$-module. Then the natural morphism 
\set\begin{equation}\label{eq_natural.tr.lim}
P\to\liminv{n\in\N}\,(P\otimes_{A_\infty}A_n)
\end{equation} 
is an isomorphism. 
\end{claim}
\begin{pfclaim} This is clear if $P=A_\infty^r$ for some $r\geq 0$.
For a general $P$ one chooses, for every $\eps\in\fm$, a sequence 
$P\to A_\infty^n\to P$ as in lemma \ref{lem_thetwo} and applies 
to it the natural transformation \eqref{eq_natural.tr.lim}; the 
claim follows by the usual diagram chase.
\end{pfclaim}
\end{proof}

\begin{theorem}\label{th_lift.etale.cpte} 
In the situation of \eqref{subsec_in.the.situation}, suppose that 
$I$ is tight. Then \eqref{eq_yet.one.fct.etale} is an equivalence 
of categories.
\end{theorem}
\begin{proof} It follows easily from theorem \ref{th_liftmod}(ii)
that the natural functor
$\text{2-}\liminv{n\in\N}\,A_n\Et_\mathrm{afp}\to A_0\Et_\mathrm{afp}$
is an equivalence. Hence we are reduced to showing that the
natural functor 
\set\begin{equation}\label{eq_yawn}
A_\infty\Et_\mathrm{afp}\to
\text{2-}\liminv{n\in\N}\,A_n\Et_\mathrm{afp}
\end{equation} 
is an equivalence. To this aim, we claim that the rule
$(B_n~|~n\in\N)\mapsto\liminv{n\in\N}B_n$ defines a quasi-inverse
to \eqref{eq_yawn}. Taking into account theorem \ref{th_for.A-mods},
this will follow once we have shown:
\begin{claim} Let $B$ be an almost finitely generated projective 
$A_\infty$-algebra such that $B_0:=B\otimes_{A_\infty}A_0$ is
unramified. Then $B$ is unramified.
\end{claim}
\begin{pfclaim} To start out, we apply theorem 
\ref{th_unramif.critter}(ii) with $J=I^n$ to deduce, by 
induction on $n$, that $B/I^{n+1}B$ is unramified 
over $A_\infty/I^{n+1}$ for every $n\in\N$. Set 
$B_n:=B\otimes_{A_\infty}A_n$; it follows that $B_n$ is 
unramified over $A_n$ for every $n\in\N$. By proposition
\ref{prop_idemp} we deduce that there exists a compatible 
system of idempotent almost elements $e_n$ of 
$B_n\otimes_{A_n}B_n$, such that $e_n\cdot I_{B_n/A_n}=0$ 
and $\mu_{B_n/A_n}(e_n)=1$ for every $n\in\N$. Hence we obtain
an idempotent in $(\liminv{n\in\N}\,B_n\otimes_{A_n}B_n)_*
\simeq\liminv{n\in\N}\,(B_n\otimes_{A_n}B_n)_*$; however, the
latter is isomorphic to $(B\otimes_{A_\infty}B)_*$, in view of
claim \ref{cl_proj.on.P}.
Then conditions (ii) and (iii) of proposition \ref{prop_idemp}
follow easily by remarking that 
$\mu_{B/A_\infty}=\liminv{n\in\N}\,\mu_{B_n/A_n}$ and
$I_{B/A_\infty}=\liminv{n\in\N}\,I_{B_n/A_n}$.
\end{pfclaim}
\end{proof}

\subsection{Henselian approximation of structures over adically
complete rings}\label{sec_elkik}
This section reviews and complements some results of Elkik's
article \cite{Elk}. Especially, we wish to show how the main
theorems of {\em loc.cit.} generalize to the case of not necessarily
noetherian rings. In principle, this is known 
(\cite[Ch.III, \S4, Rem.2, p.587]{Elk} explains briefly how to adapt 
the proofs to make them work in some non-noetherian situations), 
but we feel that it is worthwhile to give more details.

\begin{definition}\label{def_H}
\index{$\sH_R(F,J)$|indref{def_H}}
Let $R$ be a ring, $F:=R[X_1,...,X_N]$
a free $R$-algebra of finite type, $J\subset F$ a finitely
generated ideal, and let $S:=F/J$. We define an ideal of $F$
by setting: 
$$
\sH_R(F,J):=\Ann_F\,\hExt^1_S(\L_{S/R},J/J^2).
$$
\end{definition}

\begin{lemma}\label{lem_H}
In the situation of definition {\em\ref{def_H}}, we have:
\begin{enumerate}
\item
For any $R$-algebra $R'$ let $F':=R'\otimes_RF$. Then
$\sH_R(F,J)\cdot F'\subset\sH_{R'}(F',JF')$.
\item
The open subset $\Spec\,S\setminus V(\sH_R(F,J)\cdot S)$ is the smooth
locus of $S$ over $R$.
\item
 $\sH_R(F,J)$ annihilates $\hExt_S^1(\L_{S/R},N)$ for every
$S$-module $N$.
\end{enumerate}
\end{lemma}
\begin{proof} According to \cite[Ch.III, Cor.1.2.9.1]{Il}
there is a natural isomorphism in $\sD(S\Mod)$:
\set\begin{equation}\label{eq_refer-out}
\tau_{[-1}\L_{S/R}\simeq(0\to J/J^2\stackrel{\partial}{\to}
S\otimes_F\Omega_{F/R}\to 0)
\end{equation}
where $\partial$ is induced by the universal derivation
$d:F\to\Omega_{F/R}$. Since $\Omega_{F/R}$ is a free $F$-module,
we derive a natural isomorphism:
\set\begin{equation}\label{eq_calcul-Ext}
\hExt^1_S(\L_{S/R},J/J^2)\simeq
\End_S(J/J^2)/\partial^*\Hom_F(\Omega_{F/R},J/J^2).
\end{equation}
Let now $R'$ be an $R$-algebra and $h\in F$; in view of
\eqref{eq_calcul-Ext}, the condition $h\in\sH_R(F,J)$ means
precisely that the scalar multiplication $h:J/J^2\to J/J^2$
factors through $\partial$. It follows that 
$h:(J/J^2)\otimes_RR'\to(J/J^2)\otimes_RR'$ factors through
$\partial\otimes_S\one_{S'}$. However, the latter map admits
a factorization:
$$
\partial\otimes_S\one_{S'}:(J/J^2)\otimes_RR'\stackrel{\alpha}{\to}
J'/J^{\prime 2}\stackrel{\partial'}{\to}
S'\otimes_{F'}\Omega_{F'/R'}
$$
where $\alpha$ is a surjective map. It follows easily that
the scalar multiplication $h:J'/J^{\prime 2}\to J'/J^{\prime 2}$
factors through $\partial'$, which shows (i).

To show (ii), let us first pick any $h\in\sH_R(F,J)$; we have to
prove that $S_h$ is smooth over $R$. However, the foregoing
shows that multiplication by $h$ on $J/J^2$ factors through
$S\otimes_F\Omega_{F/R}$; if now $h$ is invertible, this means
that $\partial$ is a split imbedding, therefore $S_h$ is formally
smooth over $R$ by the Jacobian criterion \cite[Ch.0, Th.22.6.1]{EGAIV}.
Since $S$ is of finite presentation over $R$, the assertion follows.
Conversely, let $h\in F$ be an element such that $S_h$ is smooth
over $R$; we wish to show that $h^n\in\sH_R(F,J)$ for $n$ large enough.
One can either prove this directly using the definition of $\sH_R(F,J)$,
or else by using the following lemma \ref{lem_old-elkik} and the well
known fact that $h^n\in H_S$ for $n$ large enough.

Finally we recall the natural isomorphism:
$$
\Hom_{\sD(S\Mod)}(\L_{S/R},N[1])\simeq
\Hom_{\sD(S\Mod)}(\tau_{[-1}\L_{S/R},N[1])
$$
which, in view of \eqref{eq_refer-out}, shows that every
morphism $\tau_{[-1}\L_{S/R}\to N[1]$ factors through a map
$\tau_{[-1}\L_{S/R}\to J/J^2[1]$, whence (iii).
\end{proof}

\sset\subsubsection{}\label{subsec_H-of-Elkik}
In the situation of definition \ref{def_H}, choose a finite system
of generators $f_1,...,f_q$ for $J$; it is shown in
\cite[Ch.0, \S2]{Elk} how to construct an ideal, called $H_S$
in {\em loc.cit.} with the same properties as in lemma
\ref{lem_H}(i),(ii). However, the definition of $H_S$ depends
explicitly on the choice of the generators $f_1,...,f_q$.
This goes as follows. For every integer $p$ and every multi-index
$(\alpha)=(\alpha_1,...,\alpha_p)\in\N^p$ with 
$1\leq\alpha_1<\alpha_2...<\alpha_p\leq q$, set $|\alpha|:=p$
and let $J_\alpha\subset J$ be the subideal generated by
$f_{\alpha_1},...,f_{\alpha_p}$ and $\Delta_\alpha\subset F$
the ideal generated by the determinants of the minors of order
$p$ of the Jacobian matrix 
$(\partial f_{\alpha_i}/\partial X_j~|~1\leq i\leq p, 1\leq j\leq N)$.
We set 
$$
H_S:=\sum_{p\geq 0}\sum_{|\alpha|=p}\Delta_\alpha\cdot(J_\alpha:J).
$$
Though we won't be using the ideal $H_S$, we want to explain
how it relates to the more intrinsic $\sH_R(F,J)$. This is
the purpose of the following:

\begin{lemma}\label{lem_old-elkik} With the notation of
\eqref{subsec_H-of-Elkik}, we have: $H_S\subset\sH_R(F,J)$.
\end{lemma}
\begin{proof} Let $p\in\N$, $\alpha$ a multi-index with $|\alpha|=p$,
$\delta\in\Delta_p$ and $x\in F$ such that $xJ\subset J_\alpha$.
We can suppose that $\delta=\det(M)$, where 
$M=(\partial f_{\alpha_j}/\partial X_{\beta_i}~|~1\leq i,j\leq p)$
for a certain multi-index $\beta$ with $|\beta|=p$. We consider the
maps $\phi_\alpha:R^p\to J/J^2$ : $e_i\mapsto f_{\alpha_i}\pmod{J^2}$
and $\pi_\beta:S\otimes_F\Omega_{F/R}\to R^p$ : 
$dX_{\beta_j}\mapsto e_j$ for $j=1,...,p$ and $dX_k\mapsto 0$ if
$k\notin\{\beta_1,...,\beta_p\}$. Then the matrix of the composed
$S$-linear map $\pi_\beta\circ\partial\circ\phi_\alpha:S^p\to S^p$
is none else than $M$. Let $M'$ be the adjoint matrix of $M$, so that
$M\cdot M'=\delta\cdot\one_{S^p}$, and we can compute:
$$
\phi_\alpha\circ M'\circ\pi_\beta\circ\partial\circ\phi_\alpha=
\phi_\alpha\circ(\delta\cdot\one_{S^p})=
(\delta\cdot\one_{J/J^2})\circ\phi_\alpha.
$$
In other words, the maps
$\phi_\alpha\circ M'\circ\pi_\beta\circ\partial$ and
$\delta\cdot\one_{J/J^2}$ agree on $\Img\,\phi_\alpha=(J_\alpha+J^2)/J^2$.
Therefore: 
$(x\cdot\phi_\alpha)\circ
M'\circ\pi_\beta\circ\partial=x\cdot\delta\cdot\one_{J/J^2}$,
{\em i.e} the scalar multiplication by $x\cdot\delta$ on
$J/J^2$ factors through $\partial$; as in the proof of lemma
\ref{lem_H} this implies that $x\cdot\delta\in\sH_R(F,J)$, as
claimed.
\end{proof}

\sset\subsubsection{}\label{subsec_new-Elkik-setup} 
Let $R$ be a (not necessarily noetherian) ring, $t\in R$ 
a non-zero-divisor, $I\subset R$ an ideal,  $S$ an $R$-algebra
of finite presentation, which we write in the form $S=F/J$ for
some finitely generated ideal $J\subset F:=R[X_1,...,X_N]$.
For given $a:=(a_1,...,a_N)\in R^N$, let $\fp_a\subset F$
be the ideal generated by $(X_1-a_1,...,X_N-a_N)$. 

\begin{lemma}\label{lem_improve}
In the situation of \eqref{subsec_new-Elkik-setup},
let $n,h\geq 0$ be two integers and $a\in R^N$ such that 
\set\begin{equation}\label{eq_time-0}
t^h\in\sH_R(F,J)+\fp_a\qquad J\subset\fp_a+t^nIF\qquad n>2h.
\end{equation}
Then there exists $b\in R^N$ such that 
$$
b-a\in t^{n-h}IR^N\quad\text{and}\quad
J\subset\fp_b+(t^{n-h}I)^2F.
$$
\end{lemma}
\begin{proof} We are given a morphism $\sigma:\Spec\,R/t^nI\to\Spec\,S$
whose restriction to the closed subscheme $\Spec\,R/t^{n-h}I$
we denote by $\sigma_0$.
The claim amounts to saying that there exists a lifting
$\tilde\sigma:\Spec\, R/(t^{n-h}I)^2\to\Spec\,S$ of $\sigma_0$.
By proposition \ref{prop_obstr2}, the obstruction to the existence
of an extension of $\sigma$ to a morphism
$\Spec\,R/(t^nI)^2\to\Spec\,S$ is a class
$\omega\in\hExt_S^1(\L_{S/R},t^nI/(t^nI)^2)$. We have a commutative
diagram:
$$
\xymatrix{ t^nI/(t^nI)^2 \ar[r]^-\alpha \ar[d]_\beta &
t^{n-h}I/(t^{n-h}I)^2 \ar[d]^\gamma \\
t^nI/(t^nI)^2 \ar[r]^-\delta & t^nI/t^{2n-h}I^2
}$$
where $\alpha$ is induced by the inclusion $t^nI\subset t^{n-h}I$,
$\beta$ is the scalar multiplication by $t^h$, $\delta$ is the
restriction to $t^nI$ of the natural projection 
$R/(t^nI)^2\to R/t^{2n-h}I^2$, and $\gamma$ is an isomorphism
induced by scalar multiplication by 
$t^h:t^{n-h}I\stackrel{\sim}{\to}t^nI$. Since the $S$-module
structure on $t^nI/(t^nI)^2$ is induced by extension of scalars
via $\sigma$, it is clear that $\fp_a\cdot\omega=0$; on the
other hand, by lemma \ref{lem_H}(iii) we know that
$\sH_R(F,J)\cdot\omega=0$, hence 
$t^h\cdot\omega=\hExt^1_S(\L_{S/R},\beta)(\omega)=0$,
and consequently $\hExt^1_S(\L_{S/R},\alpha)(\omega)=0$. Since
the latter class is the obstruction to the existence of
$\tilde\sigma$, the assertion follows.
\end{proof}

Lemma \ref{lem_improve} is the basis of an inductive procedure
that allows to construct actual sections of $X:=\Spec\,S$, starting
from approximate solutions of the system of equations defined by
the ideal $J$. The section thus obtained live in $X(R^\wedge)$,
where $R^\wedge$ is the $tI$-adic completion of $R$. In case
$I$ is finitely generated, $R^\wedge$ is (separated) complete
for the $(tI)^\wedge$-adic topology, where $(tI)^\wedge$ is the topological
closure of $tI$ in $R^\wedge$; however, in later sections we will
find situations where the relevant ideal $I$ is not finitely
generated; in such case, $R^\wedge$ is complete only for the
$(tI)^\wedge$-c-adic topology (see definition \ref{def_a-topologies}(v)).
In this section we carry out a preliminary study of some topologies
on $R$-modules and on sets of sections of $R$-schemes; one of
the main themes is to compare the topologies of, say $Z(R)$
where $Z$ is an $R$-scheme, and of $Z(R^\wedge)$, where $R^\wedge$
is an adic completion of $R$. For this reason, it is somewhat
annoying that, in passing from $R$ to its completion, one is forced
to replace a preadic topology by a c-preadic one.
That is why we prefer to use a slightly coarser topology
on $R$, described in the following:

\begin{definition}\label{def_thereyougo}
\index{$(t,I)$-adic topology on a module|indref{def_thereyougo}}
Let $R$ be any ring, $t\in R$ a non-zero-divisor, $I\subset R$ an
ideal, $M_0$ (resp. $N$) a finitely generated $R$-module
(resp. $R[t^{-1}]$-module).
\begin{enumerate}
\item
The {\em $(t,I)$-preadic topology} of $M_0$ is the linear topology
that admits the family of submodules $(t^nIM_0~|~n\in\N_0)$ as a 
cofinal system of open neighborhoods. 
\item
The {\em $(t,I)$-preadic topology} of $N$ is the linear topology
$\cL_N$ defined as follows. Choose a finitely generated
$R$-module $N_0$ with an $R$-linear
map $\phi:N_0\to N$ such that $\phi[t^{-1}]:N_0[t^{-1}]\to N$
is onto; then $\cL_N$ is the finest topology such that $\phi_0$
becomes an open map when we endow $N_0$ with the $(t,I)$-preadic
topology. Therefore, the family of $R$-submodules 
$(\phi(t^kIN_0)~|~k\in\N)$ is a cofinal system of open neighborhoods
of $0\in N$ for the topology $\cL_N$.
\end{enumerate}
\end{definition}

The advantage of the $(t,I)$-adic topology is that the $(t,I)$-adic
completion $R^\wedge$ of $R$ is complete for a topology of the same
type, namely for the $(t,I^\wedge)$-topology (where 
$I^\wedge\subset R^\wedge$ is the topological closure of $I$). Hence,
for every $R$-scheme $Z$, the topologies on $Z(R)$ and $Z(R^\wedge)$
admit a uniform description. After this caveat, let us stress
nevertheless that all of the results of this section hold as well
for the c-adic topologies.

\sset\subsubsection{}
In the situation of definition \ref{def_thereyougo},
it is easy to verify that the $(t,I)$-preadic topology on $N$
does not depend on the choice of $N_0$ and $\phi$.
Indeed, if $\phi':N_1\to N$ is another choice, then
one checks that $\phi(t^kN_0)\subset\phi'(N_1)$
for $k\in\N$ large enough, and symmetrically 
$\phi_1(t^kN_1)\subset\phi_0(N_0)$, which implies the
assertion.

Moreover, it is easy to check that every $R[t^{-1}]$-linear
map $M\to N$ of finitely generated $R[t^{-1}]$-modules is
continuous for the respective $(t,I)$-preadic topologies.

\begin{lemma}\label{lem_Elkik-Noether} 
Let $R$ be a noetherian ring, $\cI\subset R$ an ideal, and
suppose that the pair $(R,\cI)$ is henselian. Let $\bar R$ be the
$\cI$-adic completion of $R$ and set $\bar S:=\bar R\otimes_RS$. 
Let $U\subset\Spec\,S$ be an open subset smooth over $\Spec\,R$, 
and $\bar U\subset\Spec\,\bar S$ the preimage of $U$. Then, for
every integer $n$ and every $\bar R$-section 
$\bar\sigma:\Spec\,\bar R\to\Spec\,\bar S$ whose restriction 
to $\Spec\,\bar R\setminus V(\cI\bar R)$ factors through 
$\bar U$, there exists an $R$-section $\sigma:\Spec\,R\to\Spec\,S$
congruent to $\bar\sigma$ modulo $\cI^n$, and whose restriction
to $\Spec\,R\setminus V(\cI)$ factors through $U$.
\end{lemma}
\begin{proof} It is \cite[Ch.II, Th.2 bis]{Elk}.
\end{proof}

\begin{proposition}\label{prop_new-Elkik} 
Keep the notation of \eqref{subsec_new-Elkik-setup} and let 
$I\subset R$ be an ideal such that $(R,tI)$ is a henselian pair.
Let $a\in R^N$ and $n,h\geq 0$ such that \eqref{eq_time-0} holds.
Then there exists $b\in R^N$ such that $b-a\in t^{n-h}IR^N$
and $J\subset\fp_b$.
\end{proposition}
\begin{proof} We consider first the following special case:
\begin{claim}\label{cl_ok.for.cplte.noeth} 
The proposition holds if $R$ is complete for the $tI$-adic topology.
\end{claim}
\begin{pfclaim} We apply repeatedly lemma \ref{lem_improve}
to obtain a $tI$-adically convergent sequence of elements
$(a_m\in R^N~|~m\in\N)$, with $a_0:=a$ and such that
$a_m\equiv a\pmod{t^{n-h}IR^N}$ for every $m\in\N$. The
limit $b$ of the sequence $(a_m~|~m\in\N)$ will do.
\end{pfclaim}

Let next $R$ be a general ring; let $H\subset\sH_R(F,J)$
be a finitely generated subideal such that $t^h\in H$.
We can write $R=\bigcup_{\lambda\in\Lambda}R_\lambda$ for a filtered
family of noetherian subrings $(R_\lambda~|~\lambda\in\Lambda)$
such that $t\in R_\lambda$ and $a\in R^N_\lambda$ for every 
$\lambda\in\Lambda$, and we set $I_\lambda:=I\cap R_\lambda$ for every
$\lambda\in\Lambda$. We can also assume that the pair
$(R_\lambda, tI_\lambda)$ is henselian for every 
$\lambda\in\Lambda$. Let $f_1,...,f_q\in F$
be a finite set of generators for $J$ and 
$g_1,...,g_r\in F$ a finite set of generators for 
$H$; up to restricting to a cofinal family, we can then assume that 
$f_i,g_j\in F_\lambda:=R_\lambda[X_1,...,X_N]$ for every
$\lambda\in\Lambda$ and every $i\leq q$, $j\leq r$. 
Let $J_\lambda\subset F_\lambda$ be the ideal generated by
$f_1,...,f_q$ and set $S_\lambda:=F_\lambda/J_\lambda$;
let also $H_\lambda\subset F_\lambda$ be the ideal generated
by $g_1,...,g_r$. Again after replacing $\Lambda$ by a cofinal
subfamily, we can achieve that 
$H_\lambda\subset\sH_{R_\lambda}(F_\lambda,J_\lambda)$ for every
$\lambda\in\Lambda$. Finally, let 
$\fp_{\lambda,a}\subset F_\lambda$ be the
ideal generated by $X_1-a_1,...,X_N-a_N$; we can assume
that $t^h\in H_\lambda+\fp_{\lambda,a}$ and 
$J_\lambda\subset\fp_{\lambda,a}+t^nI_\lambda F_\lambda$
for every $\lambda\in\Lambda$. With this setup, let 
$\bar R_\lambda$ be the $tI_\lambda$-adic completion of $R_\lambda$;
we can apply claim \ref{cl_ok.for.cplte.noeth} to deduce that 
there exists $c\in\bar R{}^N_\lambda$ such that 
$c-a\in t^{n-h}I_\lambda\bar R{}^N_\lambda$ and such that 
$J_\lambda\subset\bar\fp_{\lambda,c}$ (where $\bar\fp_{\lambda,c}$
denotes the ideal generated by $X_1-c_1,...,X-c_n$ in
$\bar R_\lambda[X_1,...,X_N]$). Let 
$U_\lambda:=\Spec\,S_\lambda\setminus V(H_\lambda\cdot S_\lambda)$
and $\bar S_\lambda:=\bar R_\lambda\otimes_{R_\lambda}S$;
from lemma \ref{lem_H}(ii) we know that $U_\lambda$
is smooth over $\Spec\,R_\lambda$, and by construction $c$
determines a $\bar R_\lambda$-section of $\Spec\,\bar S_\lambda$
whose restriction to $\Spec\,\bar R_\lambda\setminus V(tI_\lambda)$
factors through the preimage of $U_\lambda$. Therefore lemma 
\ref{lem_Elkik-Noether} ensures that there exists an 
$R_\lambda$-section $\sigma:\Spec\,R_\lambda\to\Spec\,S_\lambda$
that agrees with $c$ modulo $(tI_\lambda)^{n-h}$ and whose restriction
to $\Spec\,R_\lambda\setminus V(tI_\lambda)$ factors through $U$. 
Let $\pi_\lambda:F_\lambda\to S_\lambda$ be the natural projection; 
there is a unique $b\in R^N_\lambda$ such that
$\fp_{\lambda,b}=
\pi_\lambda^{-1}(\Ker\,\sigma^\sharp:S_\lambda\to R_\lambda)$; 
this point $b$ has the sought properties.
\end{proof}

\sset\subsubsection{}\label{subsec_adic-on.schs}
\index{$(t,I)$-adic topology on a scheme|indref{subsec_adic-on.schs}}
Let now $X$ be an affine scheme of finite type over $\Spec\,R[t^{-1}]$;
for every $n\in\N$ let $\A^n_{R[t^{-1}]}$ be the $n$-dimensional
affine space over $\Spec\,R[t^{-1}]$; for $n$ large enough we
can find a closed imbedding $j:X\to\A^n_{R[t^{-1}]}$
of $\Spec\,R[t^{-1}]$-schemes. The choice of coordinates on
$\A^n_{R[t^{-1}]}$, yields a bijection 
$\A^n_{R[t^{-1}]}(R[t^{-1}])\simeq R[t^{-1}]^n$, and then 
$j$ induces an injective map $j_*:X(R[t^{-1}])\subset R[t^{-1}]^n$.
The {\em $(t,I)$-adic topology\/} of $X(R[t^{-1}])$ is defined as the
subspace topology induced by $j_*$ (where $R[t^{-1}]^n$ is endowed
with its $(t,I)$-preadic topology).

\begin{lemma}\label{lem_t-adic.on.cshs} 
Let $X$ be as in \eqref{subsec_adic-on.schs}; we have:
\begin{enumerate}
\item
The $(t,I)$-adic topology of $X(R[t^{-1}])$ is independent
of the choice of closed imbedding $j$ and of coordinates on 
$\A^n_{R[t^{-1}]}$.
\item
If $Y$ is another affine $R[t^{-1}]$-scheme, then the natural
map $(X\times_{R[t^{-1}]}Y)(R[t^{-1}])\to 
X(R[t^{-1}])\times Y(R[t^{-1}])$ is a homeomorphism for
the $(t,I)$-adic topologies of the corresponding schemes.
\item
If $U\subset X$ is any open subset, and $\sigma\in U(R[t^{-1}])$,
there exists $f\in\cO_X(X)$ such that 
$D(f):=X\setminus V(f)\subset U$ and such that $\sigma$ factors
through $D(f)$.
\item
If $U\subset X$ is an affine open subset, then the $(t,I)$-adic
topology on $U(R[t^{-1}])$ agrees with the topology induced 
from the $(t,I)$-adic topology on $X(R[t^{-1}])$, and $U(R[t^{-1}])$
is an open subset of $X(R[t^{-1}])$.
\end{enumerate}
\end{lemma}
\begin{proof}(i): suppose $j_1:X\to\A^n_{R[t^{-1}]}$ and
$j_2:X\to\A^m_{R[t^{-1}]}$ are two closed imbeddings, $\tau_1$
and $\tau_2$ the respective topologies on $X(R[t^{-1}])$;
by symmetry it suffices to show that the identity map 
$(X(R[t^{-1}]),\tau_1)\to(X(R[t^{-1}]),\tau_2)$ is continuous.
However, $j_2$ can be extended to some morphism
$\phi:\A^n_{R[t^{-1}]}\to\A^m_{R[t^{-1}]}$, and we come down
to showing that the induced map $\phi_*:R[t^{-1}]^n\to R[t^{-1}]^m$
is continuous for the $(t,I)$-adic topology. We can further reduce
to the case of $m=1$, in which case $\phi$ is given by a
single polynomial $f\in R[t^{-1},T_1,...,T_n]$, and $\phi_*$
is the map $(a_1,...,a_n)\mapsto f(a_1,...,a_n)$.
Let $x_0\in R[t^{-1}]^n$ and set $y_0:=f(x_0)$; we have to show
that, for every $k\in\N$ there exists $h\in\N$ such that
$f(x_0+t^hI)\subset y+t^kI$. The Taylor formula gives an
identity of the form:
$f(T)=y_0+\sum_{r\in\N^n\setminus\{0\}}a_r\cdot(T-x_0)^r$,
where $a_r\in R[t^{-1}]$ for every $r\in\N^n$ and $a_r=0$ for all
but finitely many $r$. Let $s\in\N$ be an integer
large enough so that $t^sa_r\in R$ for every $r\in\N^n$; obviously
$h:=k+s$ will do. (ii) is easily reduced to the corresponding
statement for $(t,I)$-adic topologies on direct sums of 
$R[t^{-1}]$-modules; we leave the details to the reader.

(iii): we can write $U=X\setminus V(J)$, where $J\subset S:=\cO_X(X)$,
and a section $\sigma\in U(R[t^{-1}])$ induces a map
$\phi:S\to R[t^{-1}]$ such that $\phi(J)R=R$.
Let $a_i\in R$ and $f_i\in J$, $i=1,...,n$ with 
$\sum_ia_i\cdot\phi(f_i)=1$ and set $f:=\sum_if_ia_i$; it 
follows that the image of $\sigma$ is contained in the affine 
open subset $\Spec\,S[f^{-1}]\subset U$.
To show (iv), we can suppose, thanks to (iii), that $U=X\setminus V(f)$
for some $f\in\cO_X(X)$. Let $\phi_f:X\to\A_{R[t^{-1}]}^1$
be the morphism defined by $f$, and 
$\Gamma_f:X\to X\times_{R[t^{-1}]}\A_{R[t^{-1}]}^1$ its
graph; choose a closed imbedding $X\subset\A_{R[t^{-1}]}^n$.
The composition 
$j:X\to X\times_{R[t^{-1}]}\A_{R[t^{-1}]}^1\to\A_{R[t^{-1}]}^{n+1}$
is another closed imbedding, which induces the same
topology on $X(R[t^{-1}])$ in view of (i). We have 
$j^{-1}(\A_{R[t^{-1}]}^n\times_{R[t^{-1}]}\G_{m,R[t^{-1}]})=U$,
consequently we are reduced to showing the assertion for
the open imbedding 
$\A_{R[t^{-1}]}^n\times_{R[t^{-1}]}\G_{m,R[t^{-1}]}\subset
\A_{R[t^{-1}]}^{n+1}$. Using (ii) we further reduce to
considering the imbedding $\G_{m,R[t^{-1}]}\subset\A_{R[t^{-1}]}^1$.
This case can be dealt with by explicit calculations.
\end{proof}
\begin{lemma}\label{lem_unionofsects} 
Let $R$ be any ring, $X$ a quasi-projective $R$-scheme. Then every
$R$-section $\Spec\,R\to X$ of $X$ factors through an open imbedding
$U\subset X$, where $U$ is an affine $R$-scheme.
\end{lemma}
\begin{proof} By assumption $X$ is a locally closed subset
of $\P^n_R$, for some $n\in N$. 
\begin{claim}\label{cl_closed.subsch} 
The lemma holds if $X$ is a closed subscheme of $\P^n_R$.
\end{claim}
\begin{pfclaim} Indeed, in this case we can assume $X=\P^n_R$.
Therefore, let $\sigma\in\P^n_R(R)$; by \cite[Ch.II, Th.4.2.4]{EGA}
$\sigma$ corresponds to a rank one locally free quotient $L$ of 
$R^{n+1}$. We choose a section of the projection $\pi_1:R^{n+1}\to L$,
hence a decomposition $R^{n+1}\simeq L\oplus\Ker\,\pi_1$;
the induced projection $\pi_2:R^{n+1}\to\Ker\,\pi_1$ determines
a closed imbedding $\P(\Ker\,\pi_1)\to\P_R^n$ representing
the transformation of functors that assigns to every rank
one locally free quotient of $\Ker\,\pi_1$ the same module, 
seen as a quotient of $R^{n+1}$ via the projection $\pi_2$.
It is clear that the image of $\sigma$ does not intersect the 
image of $\P(\Ker\,\pi_1)$, and the complement 
$\P^n_R\setminus\P(\Ker\,\pi_1)$ is affine ($\pi_1$ is a section
in $\Gamma(\P^n_R,\cO_{\P^n_R}(1))$).
\end{pfclaim}

Thanks to claim \ref{cl_closed.subsch} we can assume that
$X$ is an open subscheme of an affine scheme 
$Y$ of finite type over $\Spec R$. In this case, the claim
reduces to lemma \ref{lem_t-adic.on.cshs}(iii) (applied with 
$t=1$).
\end{proof}

\sset\subsubsection{}
Let now $X$ be a quasi-projective $R[t^{-1}]$-scheme; each affine
open subset of $X$ comes with a natural $(t,I)$-adic 
topology. By lemma \ref{lem_t-adic.on.cshs}(iv) these topologies
agree on the intersections of any two such affine open subsets, 
and according to lemma \ref{lem_unionofsects} we have 
$X(R[t^{-1}])=\bigcup_UU(R[t^{-1}])$, where $U$ ranges on the
family of all affine open subsets $U\subset X$, so $X(R[t^{-1}])$ 
can be endowed with a well defined $(t,I)$-adic topology, independent
of all choices.

\begin{lemma}\label{lem_matrix.locus} 
The closed subscheme $U_n$ of $\Spec\,\Z[x_{11},...,x_{nn}]$ that 
classifies the $n\times n$ idempotent matrices is smooth 
over $\Spec\,\Z$.
\end{lemma}
\begin{proof} Clearly $U_n$ is of finite type over $\Spec\,\Z$,
hence it suffices to show that $U_n$ is formally smooth. Therefore, 
let $R_0$ be a ring and $I\subset R_0$ an ideal with $I^2=0$;
we need to show that the induced map $U_n(R_0)\to U_n(R_0/I)$ is 
surjective, {\em i.e.} that every $n\times n$ idempotent matrix 
$\bar M$ with entries in $R_0/I$ lifts to an idempotent matrix 
with entries in $R_0$. Pick an arbitrary matrix $M\in M_n(R_0)$ 
that lifts $\bar M$; let $E:=R_0[M]\subset M_n(R_0)$ be the 
commutative $R_0$-algebra generated by $M$, 
$\bar E\subset M_n(R_0/I)$ be the image of $E$,
$J\subset E$ the kernel of the induced map $E\to\bar E$.
We have $J^2=0$, so we can apply proposition 
\ref{prop_lift.idemp}(i) to lift $\bar M$ to some idempotent
matrix in $E$.
\end{proof}

\begin{proposition}\label{prop_denseimage} 
Let $t\in R$ be a non-zero-divisor, $I\subset R$ an ideal,
$R^\wedge:=\liminv{n\in\N}\,R/t^nI$ the $(t,I)$-adic completion
of $R$, $I^\wedge$ the topological closure of $I$ in $R^\wedge$,
and suppose that the pair $(R,tI)$ is henselian. Let $X$ be a
smooth quasi-projective $R[t^{-1}]$-scheme, and endow $X(R[t^{-1}])$
(resp. $X(R^\wedge[t^{-1}])$) with its $(t,I)$-adic (resp.
$(t,I^\wedge)$-adic) topology. Then the natural map 
$X(R[t^{-1}])\to X(R^\wedge[t^{-1}])$ has dense image.
\end{proposition}
\begin{proof} We begin with the following special case:
\begin{claim}\label{cl_denseaffine} 
The proposition holds if $X$ is affine.
\end{claim}
\begin{pfclaim} Say that $X=\Spec\,S$, where $S$ is some
finitely presented smooth $R[t^{-1}]$-algebra, and let 
$\sigma:S\to R^\wedge[t^{-1}]$ be any element of 
$X(R^\wedge[t^{-1}])$. We have to show that there are 
elements of $X(R[t^{-1}])$ in every $(t,I)$-adic neighborhood
of $\sigma$. To this aim, we pick a finitely presented 
$R$-algebra $S_0$ such that $S_0[t^{-1}]\simeq S$; after
clearing some denominators we can assume that $\sigma$
extends to a map $S_0\to R^\wedge$. Let
\set\begin{equation}\label{eq_introducing}
S_0=R[X_1,...,X_N]/J
\end{equation}
be a finite presentation of $S_0$, and set 
$H:=\sH_R(R[X_1,...,X_N],J)$ (notation of definition \ref{def_H}). By
assumption we have
\set\begin{equation}\label{eq_assmuption}
t^h\in H+J
\end{equation}
for $h\in\N$ large enough. The presentation 
\eqref{eq_introducing} defines a closed imbedding 
$X\subset\A^N_R$, and we can then find a section
$\sigma_0:\Spec\,R\to\A^N_R$ that is $(t,I)$-adically close
to $\sigma$, so that the restrictions of $\sigma$
and $\sigma_0$ agree on $\Spec\,R/t^nI$. Let 
$\fp:=\Ker\,\sigma_0^\sharp\subset R[X_1,...,X_N]$
be the ideal corresponding to $\sigma_0$; it then follows
that
\set\begin{equation}\label{eq_assmuption2}
J\subset\fp+t^nIR[X_1,...,X_N].
\end{equation}
From \eqref{eq_assmuption} and \eqref{eq_assmuption2}
we deduce that $t^h\in H+t^nIR[X_1,...,X_N]+\fp$; up to
enlarging $n$ (which is harmless) we can assume that 
$n>2h$. Let $\bar H:=(H+\fp)/\fp\subset R$; we derive:
$t^h\in\bar H+t^nI$, whence $t^h(1+at^k)\in\bar H$ for some
$a\in I$ and $k>0$; since $tI\subset\rad(R)$, it follows:
\set\begin{equation}\label{eq_assumption3}
t^h\in H+\fp.
\end{equation}
By assumption we have $\Ann_R(t)=0$, hence from 
\eqref{eq_assmuption}, \eqref{eq_assumption3} and
proposition \ref{prop_new-Elkik} we deduce the
contention.
\end{pfclaim}

\begin{claim}\label{cl_sec.through.aff}
Suppose that $X=\P^r_{R[t^{-1}]}$. Then
every section $\sigma:\Spec\,R^\wedge[t^{-1}]\to X$ 
factors through an open imbedding $U\subset X$, where
$U$ is an affine $R$-scheme. 
\end{claim}
\begin{pfclaim} By \cite[Ch.II, Th.4.2.4]{EGA}, $\sigma$ 
corresponds to a rank one projective quotient $L$ of 
$(R^\wedge[t^{-1}])^{r+1}$; we can then find an idempotent
$(r+1)\times(r+1)$ matrix $e$ such that $\Coker(e)=L$.
By lemma \ref{lem_matrix.locus}, the scheme $U$ that represents
the $(r+1)\times(r+1)$ idempotent matrices is smooth;
clearly $U$ decomposes as a disjoint union of open and closed
subschemes $U=\bigcup_{n=0}^{r+1}U_n$, where $U_n$ represents
the subfunctor that classifies all $(r+1)\times(r+1)$ idempotent
matrices of rank $n$, for every $n=0,...,r+1$. Therefore
each $U_n$ is smooth and affine over $\Spec\,\Z$ and then by claim 
\ref{cl_denseaffine} it follows that $e$ can be approximated
closely by an idempotent matrix $e_0\in M_{r+1}(R[t^{-1}])$
whose rank equals the rank of $e$. We have 
$e_0\cdot e=e_0\cdot(I_{r+1}-e_0+e)$ and if $e_0$ is sufficiently
$(t,I)$-adically close to $e$, the matrix $I_{r+1}-e_0+e$ is
invertible in $M_{r+1}(R)$, hence we can assume that:
\set\begin{equation}\label{eq_sameimage}
\Img(e_0\cdot e)=\Img(e_0)\subset R^\wedge[t^{-1}]^{r+1}.
\end{equation}
The projection 
$$
\pi:R[t^{-1}]^{r+1}\to\Img\,e_0\quad:\quad x\mapsto e_0(x)
$$
determines a closed imbedding $\P(\Img\,e_0)\to\P^r_{R[t^{-1}]}$
representing the transformation of functors that, to every
$R[t^{-1}]$-algebra $S$ and every rank one projective quotient of 
$S\otimes_{R[t^{-1}]}\Img(e_0)$ assigns the same module seen as
a quotient of $S\otimes_{R[t^{-1}]}R[t^{-1}]^{r+1}$ via $\pi$. 
By \eqref{eq_sameimage}, the restriction of 
$\one_{R^\wedge}\otimes_R\pi$ to $R^\wedge\otimes_R\Img(e)$
is an isomorphism, hence it remains such after every base change
$R^\wedge[t^{-1}]\to S$; this means that the image of $\sigma$ lands
in the complement of $\P(\Img\,e_0)$, which implies the contention.
\end{pfclaim} 
\begin{claim}\label{cl_Xgeneral} 
Let $X$ be any quasi-projective $R[t^{-1}]$-scheme,
$\sigma:\Spec\,R^\wedge[t^{-1}]\to X$ any section.
Then $\sigma$ factors through an open imbedding $U\subset X$, where
$U$ is an affine $R$-scheme. 
\end{claim}
\begin{pfclaim} Due to claim \ref{cl_sec.through.aff} we can
assume that $X$ is an open subscheme of an affine $R[t^{-1}]$-scheme
of finite type. Thus, we can write $X=Y\setminus V(J)$, where
$J\subset S:=\cO_Y(Y)$, and a section $\sigma\in X(R^\wedge[t^{-1}])$
induces a map $\phi:S\to R^\wedge[t^{-1}]$ such that 
$\phi(J)R^\wedge[t^{-1}]=R^\wedge[t^{-1}]$. Let 
$a_i\in R^\wedge[t^{-1}]$, $f_i\in J$ such that 
$\sum_ia_i\cdot\phi(f_i)=1$; we choose $b_i\in R[t^{-1}]$,  
$(t,I)$-adically close to $a_i$ ($i=1,...,n$), so that 
$\sum_i\phi(f_i)\cdot(a_i-b_i)\in tIR^\wedge$.
Since $tIR^\wedge\subset\rad(R^\wedge)$, we deduce easily 
that $\sum_i\phi(f_i)b_i\in R^\wedge$
is invertible in $R^\wedge$. Set $f:=\sum_if_ib_i$; it follows that the 
image of $\sigma$ is contained in the affine open subset 
$\Spec\,S[f^{-1}]\subset X$.
\end{pfclaim}

The proposition follows easily from claims \eqref{cl_denseaffine}
and \eqref{cl_Xgeneral}.
\end{proof}

\begin{proposition}\label{prop_openwhensmooth} 
Resume the assumptions of proposition {\em\ref{prop_denseimage}} 
and let $\phi:X\to Y$ be a morphism of quasi-projective 
$R[t^{-1}]$-schemes; then we have:
\begin{enumerate}
\item
The map $\phi_*:X(R[t^{-1}])\to Y(R[t^{-1}])$ is 
continuous for the $(t,I)$-adic topologies.
\item
If $\phi$ is smooth, $\phi_*$ is an open map.
\end{enumerate}
\end{proposition}
\begin{proof} (i): we can factor $\phi$ as a closed imbedding
$X\to X\times Y$ followed by a projection $X\times Y\to Y$, so
it suffices to prove the claim for the latter maps. The case
of a projection is immediate, and the case of a closed imbedding
follows straightforwardly from lemma \ref{lem_t-adic.on.cshs}(i).

(ii): due to lemmata \ref{lem_t-adic.on.cshs}(iv) and 
\ref{lem_unionofsects} we can assume that both $X$ and 
$Y$ are affine schemes, say $X=\Spec\,S$, $Y=\Spec\,T$ where 
$S$ and $T$ are finitely generated $R[t^{-1}]$-algebras and $S$
is smooth over $T$, especially $S$ is a finitely presented
$T$-algebra. Let $\sigma:\Spec R[t^{-1}]\to X$ be an element 
of $X(R[t^{-1}])$; after choosing presentations for $S$ and 
$T$ and clearing some denominators, we can find morphisms of 
finitely presented $R$-algebras $f:T_0\to S_0$ and $g:S_0\to R$ 
such that $f[t^{-1}]=\phi^\sharp:T\to S$ and 
$g[t^{-1}]=\sigma^\sharp:S\to R[t^{-1}]$. Let $X_0:=\Spec\,S_0$, 
$Y_0:=\Spec\,T_0$, $\sigma_0:=\Spec(g)$. We have to show
that, for every section $\xi:\Spec\,R[t^{-1}]\to Y$
sufficiently close to $\sigma$, there is $\xi'\in X(R[t^{-1}])$
such that $\phi_*(\xi')=\xi$.
\begin{claim}\label{cl_sect-extends} 
There is an open neighborhood $U\subset Y(R[t^{-1}])$
of $\phi_*(\sigma)$ such that every $\xi\in U$ extends to a
section $\xi_0:\Spec\,R\to Y_0$.
\end{claim}
\begin{pfclaim} By construction, $\phi_*(\sigma)$ is induced by the
map $g\circ f:S_0\to R$; any other section $\xi\in Y(R[t^{-1}])$
is determined by a map $\xi^\sharp:S\to R[t^{-1}]$; we have $S=S_0[t^{-1}]$
and $S_0$ is of finite type over $R$, say $S_0=R[x_1,...,x_N]/J$.
Therefore $a_i:=g\circ f(x_i)\in R$ for every $i\leq N$. The claim
follows by observing that the set of sections $\xi\in Y(R[t^{-1}])$ 
such that $\xi^\sharp(x_i)-a_i\in R$ for every $i\leq N$ forms an open
neighborhood of $\phi_*(\sigma)$. 
\end{pfclaim}

Next, choose a presentation 
\set\begin{equation}\label{eq_presentation}
S_0\simeq T_0[x_1,...,x_N]/J
\end{equation}
and let $H:=\sH_{T_0}(T_0[x_1,...,x_N],J)$ (notation of definition
\ref{def_H}). Since $S$ is smooth over $T$, we have
\set\begin{equation}\label{eq_uniformity}
t^h\in H+J
\end{equation}
for $h\in\N$ large enough. Let $\xi\in Y(R[t^{-1}])$; by claim
\ref{cl_sect-extends} we can suppose that $\xi$ extends to a
section $\xi_0:\Spec\,R\to Y_0$. Define $X_0(\xi_0)$ as the fibre 
product in the cartesian diagram
$$
\xymatrix{ X_0(\xi_0) \ar[r] \ar[d] & X_0 \ar[d]^{\phi_0} \\
           \Spec\,R \ar[r]^-{\xi_0} & Y_0.
}$$
The presentation \eqref{eq_presentation} induces a closed
imbedding $X_0(\xi_0)\subset\A_R^N$ whose defining ideal 
$J(\xi_0)\subset R[x_1,...,x_N]$ is the image of $J$;
define similarly $H(\xi_0)\subset R[x_1,...,x_N]$ as the
image of $H$. From \eqref{eq_uniformity} we deduce that
\set\begin{equation}\label{eq_specialize}
t^h\in H(\xi_0)+J(\xi_0).
\end{equation} 
Suppose now that $\xi$ is sufficiently close to $\sigma$; this means 
that the restrictions of $\xi_0$ and $\phi_{0*}(\sigma_0)$ agree 
on some closed subset $\Spec\,R/t^nI\subset\Spec\,R$. Hence, let 
$\bar\sigma_0:\Spec\,R/t^nI\to X_0(\xi)$ be the restriction
of $\sigma_0$, and choose any extension of $\bar\sigma_0$
to a morphism $\sigma_1:\Spec\,R\to\A ^N_R$. Finally,
let $\fp:=\Ker\,\sigma_1^\sharp\subset R[x_1,...,x_N]$ 
be the ideal corresponding to $\sigma_1$. By construction
we have
\set\begin{equation}\label{eq_sameoldstory}
J(\xi_0)\subset\fp+t^nIR[x_1,...,x_N].
\end{equation}
Now, arguing as in the proof of claim \ref{cl_denseaffine} we 
see that \eqref{eq_specialize} and \eqref{eq_sameoldstory} 
imply $t^h\in H(\xi_0)+\fp$, at least if $n>2h$, which can 
always be arranged. Finally, proposition 
\ref{prop_new-Elkik} shows that $\bar\sigma_0$ can 
be extended to a section $\xi':\Spec\,R\to X_0(\xi_0)$, as
required.
\end{proof}

\sset\subsubsection{}\label{subsec_group.schemes}
Resume the assumptions of proposition \ref{prop_denseimage} and
let $(X_0,X_1,s,t,c,\iota)$ be a groupoid of quasi-projective
$R[t^{-1}]$-schemes (see \eqref{subsec_interpret} for our general
notations concerning groupoids in a category). We have a natural
map
\set\begin{equation}\label{eq_onpizero}
\pi_0(X_0(R[t^{-1}]))\to\pi_0(X_0(R^\wedge[t^{-1}]))
\end{equation}
where, for any groupoid of sets 
$G:=(G_0,G_1,s_G,t_G,c_G,\iota_G)$, we denote by $\pi_0(G)$ the set 
of isomorphism classes of elements of $G_0$; this is the same as
the set of connected components of the geometric realization of
the simplicial set associated to the groupoid $G$.

\begin{theorem}\label{th_groupoid} 
Keep the notation of \eqref{subsec_group.schemes},
and suppose that $X_0$ is smooth over $R[t^{-1}]$ and that the
morphism $(s,t):X_1\to X_0\times_{R[t^{-1}]}X_0$ is smooth.
Then \eqref{eq_onpizero} is a bijection.
\end{theorem}
\begin{proof} Let $\sigma\in X_0(R^\wedge[t^{-1}])$; since
$X_0$ is smooth over $R[t^{-1}]$ we can find sections 
$\sigma_0\in X_0(R^\wedge[t^{-1}])$ arbitrarily $(t,I)$-adically
close to $\sigma$ (proposition \ref{prop_denseimage}). 
Then $(\sigma,\sigma_0)$ can be made arbitrarily close to 
$(\sigma,\sigma)\in X_0\times_{R[t^{-1}]}X_0(R^\wedge[t^{-1}])$; 
since the latter lies in the image of $X_1(R^\wedge[t^{-1}])$
under the morphism $(s,t)$, it follows that the same holds
for $(\sigma,\sigma_0)$, provided $\sigma_0$ is sufficiently 
close to $\sigma$ (proposition \ref{prop_openwhensmooth}(ii)).
This shows that \eqref{eq_onpizero} is onto.
Next, suppose that $\sigma,\tau\in X_0(R[t^{-1}])$ and
that their images $\sigma^\wedge,\tau^\wedge$ in 
$X_0(R^\wedge[t^{-1}])$ lie in the same homotopy class.
By definition, this means that 
$(\sigma^\wedge,\tau^\wedge)=(s,t)(\alpha)$ for some 
$\alpha\in X_1(R^\wedge[t^{-1}])$; since $X_1$ is smooth
over $R[t^{-1}]$ it follows that there exist sections
$\alpha_0\in X_1(R[t^{-1}])$ arbitrarily close to $\alpha$
(again proposition \ref{prop_denseimage}). Since $(s,t)_*$
is continuous (proposition \ref{prop_openwhensmooth}(i)),
it follows that $(\sigma_0,\tau_0):=(s,t)_*(\alpha_0)$
can be made arbitrarily close to $(\sigma,\tau)$ in
$X_0\times_{R[t^{-1}]}X_0(R[t^{-1}])$. This means that the
pair $(\sigma,\sigma_0)$ can be made arbitrarily close to
$(\sigma,\sigma)$, hence $(\sigma,\sigma_0)$ is in the image
of $(s,t)_*$ provided $\sigma$ and $\sigma_0$ are sufficiently
close. So $\sigma$ is in the same homotopy class as $\sigma_0$.
Likewise we can arrange that $\tau$ and $\tau_0$ are in
the same homotopy class, which shows that \eqref{eq_onpizero}
is injective as well.
\end{proof}

\sset\subsubsection{}\label{subsec_rigidify}
Let $n\in\N$ and define $U_n$ as in lemma \ref{lem_matrix.locus};
sometimes we identify $U_n$ to the functor which it represents. 
We define a functor from $F_n:\Z\Alg\to\Set$ as follows. Given a 
ring $R$, we let $F_n(R)$ be the set of all data of the form 
$(S,T,\phi,\psi)$, where $S,T\in M_n(R)$ are two idempotent 
matrices and $\phi,\psi\in M_n(R)$ are two other matrices 
submitted to the following conditions:
\begin{enumerate}
\renewcommand{\labelenumi}{(\alph{enumi})}
\item
$\phi\cdot S=0=T\cdot\phi$.
\item
$\psi\cdot T=0=S\cdot\psi$.
\item
$(I_n-S)\cdot(\psi\cdot\phi-I_n)=0=(I_n-T)\cdot(\phi\cdot\psi-I_n)$.
\renewcommand{\labelenumi}{(\roman{enumi})}
\end{enumerate}
The meaning of (a) is that $\phi$ induces a
map $\bar\phi:\Coker(S)\to\Coker(T)$; likewise, (b) means that
$\psi$ induces a map $\bar\psi:\Coker(T)\to\Coker(S)$.
Finally (c) means that $\bar\phi\circ\bar\psi$ is the
identity of $\Coker(T)$, and likewise for $\bar\psi\circ\bar\phi$.
It is easy to see from this description that the functor
$F_n$ is representable by an affine $\Z$-scheme of finite type, 
which we shall denote by the same name.
Moreover, the rule $(S,T,\phi,\psi)\mapsto(S,T)$ defines
a transformation of functors $F_n\to U_n\times U_n$, whence
a natural morphism of schemes:
\set\begin{equation}\label{eq_GL_nU_n}
F_n\to U_n\times_{\Z}U_n.
\end{equation}

\begin{lemma}\label{lem_GL_nU_n}
The morphism \eqref{eq_GL_nU_n} is smooth.
\end{lemma}
\begin{proof} Since \eqref{eq_GL_nU_n} is clearly of finite
presentation, it suffices to verify that it is formally smooth.
Hence, let $R\to R_0$ be a surjective ring homomorphism with
nilpotent kernel $J$, let $(S,T)\in U_n\times_\Z U_n(R)$ and
$(S_0,T_0,\phi_0,\psi_0)\in F_n(R_0)$ such that $(S_0,T_0)$
coincides with the image of $(S,T)$ in $U_n\times_\Z U_n(R_0)$.
We need to show that there exist $\phi,\psi\in M_n(R)$ lifting
$\phi_0$ and $\psi_0$, such that $(S,T,\phi,\psi)\in F_n(R)$.
However, since $S$ and $T$ are idempotent, $P:=\Coker(S)$
and $Q:=\Coker(T)$ are finitely generated projective $R$-modules.
According to \eqref{subsec_rigidify}, the induced map
$\bar\phi_0:P_0=:P\otimes_RR_0\to Q_0:=Q\otimes_RR_0$ is an
isomorphism with inverse $\bar\psi_0$. Let $\pi_P:P\to P_0$,
$\pi_Q:Q\to Q_0$ be the projections; since 
$\bar\psi_0\circ\pi_Q:Q\to P_0$ is
surjective, we can find a map $\bar\phi:P\to Q$ such that
$\bar\psi_0\circ\pi_Q\circ\bar\phi=\pi_P$. Using Nakayama's
lemma one checks easily that $\bar\phi$ is an isomorphism
that lifts $\bar\phi_0$. Let $\alpha:R^n\to Q$ and 
$\beta:R^n\to P$ be the natural projections.
We set $\phi:=(I_n-T)\circ\bar\phi\circ\beta$ and
$\psi:=(I_n-S)\circ\bar\phi{}^{-1}\circ\alpha$ and
leave to the reader the verification of the identities
(a)-(c) of \eqref{subsec_rigidify}.
\end{proof}

\begin{corollary}\label{cor_elkik}
Resume the assumptions of proposition {\em\ref{prop_denseimage}}.
Then the base change functor
$R[t^{-1}]\Mod\to R^\wedge[t^{-1}]\Mod~:~
M\mapsto M\otimes_RR^\wedge$ induces a bijection from the
set of isomorphism classes of finitely generated projective 
$R[t^{-1}]$-modules to the set of isomorphism classes of
finitely generated projective $R^\wedge[t^{-1}]$-modules.
\end{corollary}
\begin{proof} Resume the notation of \eqref{subsec_rigidify}
and let $s,t:F_n\to U_n$ be the morphisms obtained by
composing \eqref{eq_GL_nU_n} with the two projections onto
$U_n$. The datum $(U_n,F_n,s,t)$ can be completed to
a groupoid of schemes, by letting 
$\iota:U_n\to F_n$ be the morphism representing the
transformation of functors: $S\mapsto(S,S,I_n,I_n)$,
and $c:F_n\times_{U_n}F_n\to F_n$ the morphism representing
the transformation:
$((S,T,\phi_1,\psi_1),(S,T,\phi_2,\psi_2))\mapsto
(S,T,\phi_2\cdot\phi_1,\psi_2\cdot\psi_1)$.
Since every finitely generated projective
module can be realized as the cokernel of an idempotent
endomorphism of a free module of finite rank, the assertion 
is a straightforward consequence of theorem \ref{th_groupoid} 
and lemmata \ref{lem_matrix.locus}, \ref{lem_GL_nU_n}.
\end{proof}

\sset\subsubsection{}\label{subsec_alg.structureP}
Let $S$ be an $R[t^{-1}]$-scheme, $\cP$ a coherent $\cO_S$-module. 
An {\em $S$-algebra structure\/} on $\cP$ is a datum 
$(\mu,\underline 1)$ consisting of a map 
$\mu:\cP\otimes_{\cO_S}\cP\to\cP$ and a global section 
$\underline 1\in\cP(S)$, such that $(\cP,\mu,\underline 1)$ 
is an $\cO_S$-algebra, {\em i.e.}
such that $\mu$ and $\underline 1$ satisfy the following 
conditions for every open subset $U\subset S$ and every 
local sections $x,y,z\in\cP(U)$ :
\begin{enumerate}
\renewcommand{\labelenumi}{(\alph{enumi})}
\item
$\mu(\mu(x\otimes y)\otimes z)=\mu(x\otimes\mu(y\otimes z))$. 
\item
$\mu(x\otimes y)=\mu(y\otimes x)$.
\item
$\mu(\underline 1\otimes x)=x$.
\renewcommand{\labelenumi}{(\roman{enumi})}
\end{enumerate}
We say that an $S$-algebra structure on $\cP$ is {\'e}tale
if the datum $(\cP,\mu,\underline 1)$ is an {\'e}tale $\cO_S$-algebra.
We denote by $\mathrm{Alg}_S(\cP)$ (resp. $\mathrm{Et}_S(\cP)$)
the set of all $S$-algebra structures (resp.
{\'e}tale $S$-algebra structures) on $\cP$. 
If now $P$ is a finitely presented $R[t^{-1}]$-module, 
we obtain a functor: 
\set\begin{equation}\label{eq_fctor-Alg}
R[t^{-1}]\text{-}\mathbf{Scheme}\to\Set\qquad
S\mapsto\mathrm{Alg}_S(\cO_S\otimes_{R[t^{-1}]}P).
\end{equation}

\begin{lemma}\label{lem_alg.structureP} 
Suppose that $P$ is a finitely generated
projective $R[t^{-1}]$-module. Then the functor 
\eqref{eq_fctor-Alg} is representable by a finitely 
presented $R[t^{-1}]$-algebra.
\end{lemma}
\begin{proof} Clearly \eqref{eq_fctor-Alg} is a sheaf on
the fppf topology of $\Spec\,R[t^{-1}]$, hence it suffices
to show that we can cover $\Spec\,R[t^{-1}]$ by finitely
many affine Zariski open subsets $U_i$, such that the
restriction of \eqref{eq_fctor-Alg} to $U_i$ is representable
and finitely presented.
However, $P$ is locally free of finite rank on the Zariski 
topology of $\Spec\,R[t^{-1}]$, hence we can assume that $P$
is a free $R[t^{-1}]$-module. Let $e_1,...,e_n$ be a basis
of $P$; then a multiplication law $\mu$ is determined by its
values $a_{ij}\in P$ on $e_i\otimes e_j$; by writing 
$a_{ij}=\sum_ka_{ijk}e_k$ we obtain $n^3$ elements of
$R[t^{-1}]$; likewise, $\underline 1$ is represented by
elements $b_1,...,b_n\in R[t^{-1}]$, and conditions (a),(b) 
and (c) of \eqref{subsec_alg.structureP} translate as a finite 
system of polynomial identities for the $a_{ijk}$ and the $b_l$;
in other words, our functor is represented by a quotient of 
the free polynomial algebra 
$R[t^{-1},X_{ijk}, Y_l~|~i,j,k,l=1,...,n]$ by
a finitely generated ideal, which is the contention.
\end{proof}

\begin{lemma}\label{lem_callitU_P} 
Keep the assumptions of lemma {\em\ref{lem_alg.structureP}}.
Let $X_P$ be an affine $R[t^{-1}]$-scheme representing
the functor \eqref{eq_fctor-Alg}. Then  the functor
\set\begin{equation}\label{eq_fctor.etale}
R[t^{-1}]\text{-}\mathbf{Scheme}\to\Set\qquad
S\mapsto\mathrm{Et}_S(\cO_S\otimes_{R[t^{-1}]}P)
\end{equation}
is represented by an affine open subset $U_P\subset X_P$.
Moreover, $U_P$ is smooth over $\Spec\,R[t^{-1}]$.
\end{lemma}
\begin{proof} Let us show first that the functor
\eqref{eq_fctor.etale} is formally smooth. Indeed,
suppose that $Z$ is an affine $R[t^{-1}]$-scheme and 
$Z_0\subset Z$ is a closed subset defined by a nilpotent 
ideal. Let $(\cO_Z\otimes_{R[t^{-1}]}P,\mu_0)$ be an 
{\'e}tale $Z_0$-algebra structure; we can lift it to some
{\'e}tale $\cO_Z$-algebra $(\cQ,\mu)$ and then $\cQ$
is necessarily a locally free sheaf of finite rank 
(for instance by lemma \ref{lem_flat.proj}(i),(iii));
from proposition \ref{prop_best} we deduce that
$\cQ\simeq\cO_Z\otimes_{R[t^{-1}]}P$, whence the claim.
Next, let $(\cP,\mu,\underline 1)$ be the universal 
$\cO_{X_P}$-algebra structure on 
$P\otimes_{R[t^{-1}]}\cO_{X_P}$; let also
$\delta\in\cO_{X_P}(X_P)$ be the discriminant of the
trace form of $\cP$. By theorem \ref{th_proj.etale},
a point $x\in X_P$ is in the support of $\delta$ if
and only if $\cP_x$ is not {\'e}tale over $\cO_{X_P,x}$;
therefore, the subset $U_P$ over which $\cP$ is
{\'e}tale is indeed open and affine. To conclude, it suffices
to show that $U_P$ represents the functor 
\eqref{eq_fctor.etale}. This amounts to showing that,
for every morphism $f:S\to X_P$ of $R[t^{-1}]$-schemes,
the algebra $f^*\cP$ is {\'e}tale over $\cO_S$ if and
only if the image of $f$ lands in $U_P$. However,
the latter statement follows easily from 
\cite[Ch.IV, Cor.17.6.2]{EGA4}.
\end{proof}

\sset\subsubsection{}
Let $S$ be an $R[t^{-1}]$-scheme, $\cP$ a coherent $\cO_S$-module.
We denote by $\Aut_{\cO_S}(\cP)$ the group of $\cO_S$-linear
automorphisms of $\cP$. Then, for a given finitely presented
$R[t^{-1}]$-module $P$ we obtain a group-valued functor
\set\begin{equation}\label{eq_auto.of.P}
R[t^{-1}]\text{-}\mathbf{Scheme}\to\mathbf{Grp}\qquad
S\mapsto \Aut_{\cO_S}(\cO_S\otimes_{R[t^{-1}]}P).
\end{equation}
\begin{lemma} Keep the assumptions of lemma 
{\em\ref{lem_alg.structureP}}. Then the functor \eqref{eq_auto.of.P}
is representable by a finitely presented $R[t^{-1}]$-group scheme.
\end{lemma}
\begin{proof} It is analogous to the proof of lemma 
\ref{lem_alg.structureP} : up to restricting to a Zariski
open subset, we can assume that $P$ is free of some rank 
$n$. Then the group scheme representing our functor is
just $\GL_{n,R[t^{-1}]}$.
\end{proof}

\sset\subsubsection{}
Let now $P$ be a finitely generated projective $R[t^{-1}]$-module.
For any $R[t^{-1}]$-scheme, let $\cP_S:=\cO_S\otimes_{R[t^{-1}]}P$; 
we wish to define an action of $\Aut_{\cO_S}(\cP_S)$ on the set
$\mathrm{Et}_S(\cP_S)$. Indeed, if
$\mu:\cP_S\otimes_{\cO_S}\cP_S\to\cP_S$ is any {\'e}tale
$S$-algebra structure and $g\in\Aut_{\cO_S}(\cP_S)$, let
$\mu^g$ be the unique $S$-algebra structure on $\cP_S$
such that $g$ is an isomorphism of {\'e}tale $\cO_S$-algebras:
$g:(\cP_S,\mu)\stackrel{\sim}{\to}(\cP_S,\mu^g)$. It is
obvious that the rule $(g,\mu)\mapsto\mu^g$ is a functorial
group action. Let $U_P$ be as in lemma \ref{lem_callitU_P},
and let $\Aut_P$ be a $R[t^{-1}]$-scheme representing the
functor \ref{eq_auto.of.P}; the functorial map 
$(g,\mu)\mapsto(\mu^g,\mu)$ is represented by a morphism
of schemes:
\set\begin{equation}\label{eq_morphing}
\Aut_P\times_{R[t^{-1}]}U_P\to U_P\times_{R[t^{-1}]}U_P.
\end{equation}

\begin{lemma}\label{lem_grpetale} 
The morphism \eqref{eq_morphing} is {\'e}tale.
\end{lemma}
\begin{proof} The map is clearly of finite presentation, hence
it suffices to show that it is formally {\'e}tale. Therefore,
let $\phi:Z\to U_P\times_{R[t^{-1}]}U_P$ be a morphism of
$R[t^{-1}]$-schemes, $Z_0\subset Z$ a closed subscheme
defined by a nilpotent ideal, and suppose that the 
restriction of $\phi$ to $Z_0$ lifts to a morphism
$\psi_0:Z_0\to \Aut_P\times_{R[t^{-1}]}U_P$. We need to show
that $\phi$ lifts uniquely to a morphism $\psi$ that extends 
$\psi_0$. 
However, the datum of $\phi$ is equivalent to the datum
consisting of a pair of $Z$-algebra structures $(\cP_Z,\mu_1)$
and $(\cP_Z,\mu_2)$. The datum of $\psi_0$ is equivalent to the
datum of a $Z_0$-algebra structure $\mu_0$ on $\cP_{Z_0}$, and
of an automorphism $g_0$ of $\cP_{Z_0}$. Finally, the fact
that $\psi_0$ lifts the restriction of $\phi$ means that
$\mu_0=\mu_2\otimes_Z\one_{Z_0}$, and 
$g_0:(\cP_{Z_0},\mu_1\otimes_Z\one_{Z_0})\stackrel{\sim}{\to}
(\cP_{Z_0},\mu_2\otimes_Z\one_{Z_0})$ is an isomorphism
of {\'e}tale $\cO_{Z_0}$-algebras. By theorem \ref{th_liftetale}(iii), 
such an isomorphism extends uniquely to an isomorphism of {\'e}tale 
$\cO_Z$-algebras $g:(\cP_Z,\mu_1)\stackrel{\sim}{\to}(\cP_Z,\mu_2)$. 
The datum $(g,\mu_2)$ is equivalent to the sought map $\psi$.
\end{proof}

\begin{proposition}\label{prop_elkik.etale}
Resume the assumptions of proposition {\em\ref{prop_denseimage}}.
Then the base change functor
$R[t^{-1}]\Alg\to R^\wedge[t^{-1}]\Alg$ induces an equivalence
of categories from the category of finite {\'e}tale 
$R[t^{-1}]$-algebras to the category of finite {\'e}tale
$R^\wedge[t^{-1}]$-algebras.
\end{proposition}
\begin{proof} Let $(P^\wedge,\mu^\wedge)$ be a finite {\'e}tale 
$R^\wedge[t^{-1}]$-algebra; in particular,  $P^\wedge$ is
a finitely generated projective $R^\wedge[t^{-1}]$-module,
hence by corollary \ref{cor_elkik} we can find a finitely
generated projective $R[t^{-1}]$-module $P$ such that 
$R^\wedge[t^{-1}]\otimes_{R[t^{-1}]}P\simeq P^\wedge$.
The functorial action \eqref{eq_morphing} of $\Aut_P$ on the 
set of {\'e}tale algebra structures on $P$ defines a groupoid
of quasi-projective schemes. It then follows from theorem
\ref{th_groupoid} and lemmata \ref{lem_callitU_P},
\ref{lem_grpetale} that the base change map from the set of
isomorphism classes of $R[t^{-1}]$-algebra structures on $P$ 
to the set of isomorphism classes of $R^\wedge[t^{-1}]$-algebra 
structures on $P^\wedge$ is bijective. This shows 
that the base change functor is essentially surjective on 
finite {\'e}tale $R^\wedge[t^{-1}]$-algebras. 
To prove full faithfulness, let $Y_1$, $Y_2$ be any two finite
{\'e}tale schemes over $X:=R[t^{-1}]$; we let $H$ be the functor
that assigns to every $R[t^{-1}]$-scheme $Z$ the set 
$\Hom_Z(Z\times_XY_1,Z\times_XY_2)$. Since the functors represented 
by $Y_1$ and $Y_2$ are locally constant sheaves, say in the fppf 
topology of $X:=\Spec\,R[t^{-1}]$, the same holds for the functor
$H$, hence the latter is represented by a finite {\'e}tale $X$-scheme
(proposition \ref{prop_covers}), which we denote by the same name.
We can view $H$ as a trivial groupoid ({\em i.e.} such that
for every $X$-scheme $Z$, the only morphisms of the groupoid
$H(Z)$ are the identity morphisms of its objects). In this
case, the associated morphism $(s,t):H\to H\times_XH$ is none
other than the diagonal morphism; especially, the latter is
an open imbedding, hence theorem \ref{th_groupoid} applies and 
yields the sought assertion.
\end{proof}

\subsection{Lifting theorems for henselian pairs}\label{sec_lift.hens}
For the considerations that follow, it will be useful to
generalize a little our usual setup : we wish to work with
sheaves of almost modules (or almost algebras) on a 
scheme. This is just a matter of introducing the relevant
language, so we will proceed somewhat briskily.

\sset\subsubsection{}\label{subsec_globalize}
\index{$\cO^a_X\Mod$, $\cO_X^a\Alg$ : sheaves of $\cO^a_X$-modules,
resp. $\cO^a_X$-algebras on a scheme $X$|indref{subsec_globalize}}
\index{$\cO^a_X\Mod$, $\cO_X^a\Alg$ : sheaves of $\cO^a_X$-modules,
resp. $\cO^a_X$-algebras on a scheme $X$!$\cO_X^a\Mod_\mathrm{qcoh}$,
$\cO_X^a\Alg_\mathrm{qcoh}$ : quasi-coherent|indref{subsec_globalize}}
Let $X$ be a scheme over $\Spec(V)$. For every open
subset $U\subset X$, $\Gamma(U,\cO_X)$ is a $V$-algebra,
hence we obtain a sheaf of $V^a$-algebras $\cO^a_X$ on
$X$ by setting $\Gamma(U,\cO^a_X):=\Gamma(U,\cO_X)^a$
for every open $U\subset X$. We refer to \cite[Ch.0, \S3.1]{EGAI}
for generalities on sheaves with values in arbitrary
categories; in particular the sheaves of $V^a$-modules 
on $X$ form an abelian tensor category, and hence we
can define a sheaf $\cF$ of $\cO_X^a$-modules on $X$ (briefly: a
$\cO_X^a$-module) as a sheaf of almost modules endowed
with a scalar multiplication $\cO_X^a\otimes_{V^a}\cF\to\cF$.
Those gadgets form a category that we denote by
$\cO_X^a\Mod$. There is a functor
\set\begin{equation}\label{eq_sheafify.alm}
\Gamma(X,\cO_U^a)\Mod\to\cO_X^a\Mod\qquad M\mapsto M^\sim
\end{equation}
defined as one expects. We say that $\cF$ is {\em quasi-coherent\/} 
if we can cover $X$ by affine open subsets $U_i$, such
that $\cF_{|U_i}$ is in the essential image of a functor
\eqref{eq_sheafify.alm}. We denote by $\cO_X^a\Mod_\mathrm{qcoh}$
the full subcategory of quasi-coherent $\cO_X^a$-modules.
Similarly, we denote by $\cO_X^a\Alg$
(resp. $\cO_X^a\Alg_\mathrm{qcoh}$) the category of 
$\cO_X^a$-algebras (resp. quasi-coherent $\cO_X^a$-algebras)
defined as one expects.

\sset\subsubsection{}\label{subsec_usual-gens}
\def\pervcF{\cF_!}
\index{$\cO^a_X\Mod$, $\cO_X^a\Alg$ : sheaves of $\cO^a_X$-modules,
resp. $\cO^a_X$-algebras on a scheme $X$!$\pervcF$, $\cF_*$ (for an
$\cO^a_X$-module $\cF$)|indref{subsec_usual-gens}}
Since the functors $M\to M_!$ and $M\mapsto M_*$ from $A$-modules 
to $A_*$-modules are right exact, we can globalize them to the
situation of \eqref{subsec_globalize}. Thus, for every $V$-scheme
$X$ there are functors
$$
\cO_X^a\Mod\to\cO_X\Mod\qquad \cF\mapsto\cF_!\quad
\text{(resp. $\cF\mapsto\cF_*$)}
$$
which are left (resp. right) adjoint to the localization
functor $\cO_X\Mod\to\cO_X^a\Mod$. The functor $\cF\mapsto\cF_!$
is exact and preserves quasi-coherence (as can be easily deduced 
from proposition \ref{eq_natural.transf}), hence it provides
a left adjoint to the localization functor 
$\cO_X\Mod_\mathrm{qcoh}\to\cO_X^a\Mod_\mathrm{qcoh}$.
The functor $\cF\mapsto\cF_*$ does not preserve quasi-coherence,
in general.

\sset\subsubsection{}\label{subsec_affine.alm.schemes}
Let $R$ be a $V$-algebra and set $X:=\Spec(R)$.
Using the full faithfulness of the functor $\cF\mapsto\cF_!$
one can easily verify that the functor $M\mapsto M^\sim$
from $R^a$-modules to quasi-coherent $\cO_X^a$-modules is 
an equivalence, whose quasi-inverse is given by the global 
section functor.

After these preliminaries, we are ready to state the
following descent result which will be crucial for the
proof of theorem \ref{th_lifts.proj.hensel}.

\begin{proposition}\label{prop_descent.hensel} 
{\em (i)}\ \ Let $R$ be a $V$-algebra, $J\subset R$ a finitely
generated ideal and $R\to R'$ a flat morphism inducing an 
isomorphism $R/J\to R'/JR'$. 
Let $X:=\Spec(R)$, $X':=\Spec(R')$, $U:=X\setminus V(J)$ 
and $U':=U\times_XX'$. 
Then the natural commutative diagrams of functors:
$$
\xymatrix{ \cO_X^a\Mod_\mathrm{qcoh} \ar[r] \ar[d] & 
\cO_{X'}^a\Mod_\mathrm{qcoh} \ar[d] \\
\cO_U^a\Mod_\mathrm{qcoh} \ar[r] & 
\cO_{U'}^a\Mod_\mathrm{qcoh}
}\quad
\xymatrix{ \cO_X^a\Alg_\mathrm{qcoh} \ar[r] \ar[d] & 
\cO_{X'}^a\Alg_\mathrm{qcoh} \ar[d] \\
\cO_U^a\Alg_\mathrm{qcoh} \ar[r] & 
\cO_{U'}^a\Alg_\mathrm{qcoh}
}$$
are $2$-cartesian (that is, cartesian in the category
of\/ $2$-categories).
\begin{enumerate}
\addtocounter{enumi}{1}
\item Let $A$ be a $V^a$-algebra, $f\in A_*$ a non-zero-divisor,
$A^\wedge$ the $f$-adic completion of $A$. Denote by 
$A\Mod_f$ (resp. $A\Alg_f$) the full subcategory of $f$-torsion 
free $A$-modules (resp. $A$-algebras), and similarly define 
$A^\wedge\Mod_f$ (resp. $A^\wedge\Alg_f$).
Then the natural commutative diagrams of functors
$$
\xymatrix{ A\Mod_f \ar[r] \ar[d] & A^\wedge\Mod_f \ar[d] \\
A[f^{-1}]\Mod \ar[r] & A^\wedge[f^{-1}]\Mod
}\quad
\xymatrix{ A\Alg_f \ar[r] \ar[d] & A^\wedge\Alg_f \ar[d] \\
A[f^{-1}]\Alg \ar[r] & A^\wedge[f^{-1}]\Alg
}$$
are $2$-cartesian.
\end{enumerate}
\end{proposition}
\begin{proof} (i): for the functors on $\cO_X^a$-modules, one applies 
the functor $\cF\mapsto\cF_!$, thereby reducing to the corresponding 
assertion for quasi-coherent $\cO_X$-modules. Under assumption
(a), the latter is proved in \cite[Prop.4.2]{Fe-Ra} (actually, in 
{\em loc.cit.} one assumes that $X'$ is faithfully flat 
over $X$, but one can reduce to such case after replacing 
$X'$ by $X'\amalg(U_1\amalg ...\amalg U_n)$, where $(U_i~|~i\leq n)$
is a finite cover of $U$ by affine open subsets; notice also that 
{\em loc.cit.} omits the assumption that $J$ is finitely generated, 
but the proof works only under such assumption). Since all the 
functors involved commute with tensor products, the assertion 
about $\cO_X^a$-algebras follows formally.

(ii): for modules one argues as in the proof of (i), except that
instead of invoking \cite{Fe-Ra}, one uses \cite[Theorem]{Be-La}.
For algebras, one has to proceed a little more carefully, since
the tensor product of two $f$-torsion free modules may
fail to be $f$-torsion free. Hence, let $(B_1,B_2,\beta)$ the
datum consisting of an $A[f^{-1}]$-algebra $B_1$, an 
$A^\wedge$-algebra $B_2$ and an isomorphism 
$\beta:B_1\otimes_AA^\wedge\stackrel{\sim}{\to}B_2[f^{-1}]$ of
$A^\wedge[f^{-1}]$-algebras. Let 
$I:=\bigcup_{n>0}\Ann_{B_2\otimes_{A^\wedge}B_2}(f^n)$ and
set $C:=B_2\otimes_{A^\wedge}B_2/I$; $\beta$ induces an isomorphism
$\gamma:B_1\otimes_AB_1\otimes_AA^\wedge\stackrel{\sim}{\to}C$ of
$A^\wedge$-modules, so by the foregoing there exists an $A$-module
$D$ such that $D[f^{-1}]\simeq B_1\otimes_AB_1$ and 
$D\otimes_AA^\wedge\simeq C$. Furthermore, the multiplication
morphism $\mu_{B_2/A^\wedge}$ factors through a morphism 
$\tilde\mu:D\to B_2$, and consequently the datum 
$(\mu_{B_1/A[f^{-1}]},\tilde\mu,\gamma,\beta)$ determines a unique
morphism $D\to B$. Let $I':=\bigcup_{n>0}\Ann_{B\otimes_AB}(f^n)$;
one verifies easily that $(B\otimes_AB/J)\otimes_AA^\wedge\simeq C$,
so again the same sort of arguments show that 
$D\simeq B\otimes_AB/J$, hence we obtain a morphism
$\mu_{B/A}:B\otimes_AB\to B$ that lifts $\mu_{B_1/A[f^{-1}]}$
and $\mu_{B_2/A^\wedge}$. Arguing along the same lines one can 
now verify easily that $(B,\mu_{B/A})$ is really a $B$-algebra :
we leave the details to the reader.
\end{proof}

\begin{theorem}\label{th_lifts.proj.hensel}
Let $(A,I)$ be a tight henselian pair, $\bar P$
an almost finitely generated projective $A/I$-module, 
$\fm_1\subset\fm$ a finitely generated subideal. We have:
\begin{enumerate}
\item
If $\tilde\fm$ has homological dimension $\leq 1$, then there
exists an almost finitely generated projective $A$-module such
that $P\otimes_A(A/I)\simeq\bar P$.
\item
If $P_1$ and $P_2$ are two liftings of $P$ as in (i) and if 
there exists an isomorphism 
$\bar\beta:P_1\otimes_A(A/\fm_1I)\stackrel{\sim}{\to}
P_2\otimes_A(A/\fm_1I)$, then there exists an isomorphism 
$\beta:P_1\to P_2$ such that 
$\beta\otimes_A\one_{A/I}=\bar\beta\otimes_A\one_{A/I}$.
\item
With the notation of \eqref{subsec_in.the.situation}, the 
natural functor $A\Et_\mathrm{afp}\to(A/I)\Et_\mathrm{afp}$
is an equivalence of categories.
\end{enumerate}
\end{theorem}
\begin{proof} We begin by showing (ii): indeed, the obstruction 
to the existence of a morphism $\alpha:P_1\to P_2$ such that 
$\alpha\otimes_A\one_{(A/\fm_1I)}=\beta$ is a class 
$\omega\in\Ext^1_A(P_1,\fm_1IP_2)$. The same argument
used in the proof of claim \ref{cl_vanishing.trick} shows
that the natural map 
$\Ext^1_A(P_1,\fm_1IP_2)\to\Ext^1_A(P_1,IP_2)$
vanishes identically, and proves the assertion.

Next we wish to show that the functor of (iii) is
fully faithful. Therefore, let $B,C$ be two almost finitely
presented {\'e}tale $A$-algebras, and $\bar\phi:B/IB\to C/IC$
a morphism. According to lemma \ref{lem_graph.rule}, $\bar\phi$
is characterized by its associated idempotent, call it 
$\bar e\in(B\otimes_AC/IC)_*$. Set $D:=B\otimes_AC$;
according to lemma \ref{lem_plenty}(iii), the pair $(D,I D)$
is tight henselian. Then proposition \ref{prop_bij.tight.hens}
says that $\bar e$ lifts uniquely to an idempotent $e\in D_*$.

\begin{claim}\label{cl_about.Gamma}
The associated morphism $\Gamma(e)$ is an isomorphism (notation 
of \eqref{subsec_graph.rule}).
\end{claim}
\begin{pfclaim} Indeed, by naturality of $\Gamma$, we have 
$\Gamma(e)\otimes_A\one_{A/I}=\Gamma(\bar e)$, so the assertion
follows from corollary \ref{cor_Naka}.
\end{pfclaim}

By claim \ref{cl_about.Gamma} and lemma \ref{lem_graph.rule},
$e$ corresponds to a unique morphism $\phi:B\to C$ which
is the sought lifting of $\bar\phi$. The remaining steps to
complete the proof of (iii) will apply as well to the proof
of (i). Pick an integer $n>0$ and a finitely generated subideal 
$\fm_0\subset\fm$ such that $I^n\subset\fm_0A$; we notice
that assertions (i) and (iii) also hold when $I$ is nilpotent, 
since in this case they reduce to theorem \ref{th_liftmod}(i.b),(ii).
It follows easily that it suffices to prove the assertions 
for the pair $(A,I^n)$, hence we can and do assume throughout
that $I\subset\fm_0A$.

\begin{claim}\label{cl_principle} 
Assertions (i) and (iii) hold if $I$ is generated by a 
non-zero-divisor of $A_*$.
\end{claim}
\begin{pfclaim} Say that $I=fA$, for some non-zero-divisor
$f\in\fm_0A_*$, and let $A^\wedge$ be the $f$-adic completion of 
$A$. By theorem \ref{th_liftmod}(i), $\bar P$ lifts to a compatible 
system $(P_n~|~n\in\N)$ of almost finitely generated projective 
$A/I^{n+1}$-modules; by theorem \ref{th_for.A-mods}, the latter 
compatible system gives rise to a unique almost finitely generated 
projective $A^\wedge$-module $P^\wedge$.
Notice that $f$ is regular on $A^\wedge$, hence also on 
$P^\wedge$. Since $A^\wedge[f^{-1}]$ is a (usual) $V$-algebra, the 
$A^\wedge[f^{-1}]$-module $P^\wedge[f^{-1}]$ is finitely generated 
projective; it follows from corollary \ref{cor_elkik} that there 
exists a finitely generated projective $A[f^{-1}]$-module $Q$ with 
an isomorphism $\beta:Q\otimes_AA^\wedge\simeq P^\wedge[f^{-1}]$. 
By proposition \ref{prop_descent.hensel}(ii) the datum 
$(P^\wedge,Q,\beta)$ determines a unique $f$-torsion-free
$A$-module $P$ which lifts $\bar P$. Since $f$ is regular 
on both $P$ and $A$, we have $\Tor_i^A(P,A/fA)=0$ for $i=1,2$,
hence $P$ is $A$-flat, by virtue of lemma \ref{lem_Tor.orgy}.
Next, set $C:=A[f^{-1}]\times A^\wedge$; from lemma 
\ref{lem_vanishing.critter}(i) we deduce that 
$\Ann_C(P\otimes_AC)^2\subset\Ann_A(P)$, and since 
$\Tor^A_1(C,P)=0$, remark \ref{rem_descent}(i) implies
that $P$ is almost finitely presented, therefore almost
projective over $A$, which shows that (i) holds. 
Likewise, let $\bar B$ be an almost finitely presented {\'e}tale
$A/I$-algebra; by theorem \ref{th_lift.etale.cpte}, $\bar B$ 
admits a unique lifting to an almost finitely presented {\'e}tale
$A^\wedge$-algebra $B^\wedge$. Then $B^\wedge[f^{-1}]$ is
a finite {\'e}tale $A^\wedge[f^{-1}]$-algebra, hence by 
proposition \ref{prop_elkik.etale} there exists a unique 
finite {\'e}tale $A[f^{-1}]$-algebra $B_0$ with an isomorphism
$\beta:B_0\otimes_BB^\wedge\stackrel{\sim}{\to}B^\wedge[f^{-1}]$. 
By proposition \ref{prop_descent.hensel}(ii), the datum 
$(B^\wedge,B_0,\beta)$ determines a unique $f$-torsion-free 
$A$-algebra $B$; the foregoing proof of assertion (i)
applies to the $A$-module underlying $B$ and shows that
$B$ is an almost finitely generated projective $A$-algebra.
By construction $B/fB$ is unramified over $A/fA$, so
theorem \ref{th_unramif.critter}(iv) applies and shows 
that $B$ is unramified over $A$. Thus, we have shown that 
the functor of (iii) is essentially surjective under the
present assumptions; since it is already known in general
that this functor is fully faithful, assertion (iii) is
completely proved in this case.
\end{pfclaim}

\begin{claim}\label{cl_can.lift.princ}
Assertions (i) and (iii) hold if $I$ is a principal ideal. 
\end{claim}
\begin{pfclaim} Say that $I=fA$ for some $f\in\fm_0A_*$.
Let $J:=\bigcup_{n>0}\Ann_A(f^n)$; we have a cartesian
diagram
$$
\xymatrix{ A/(J\cap I) \ar[r] \ar[d] & A/I \ar[d] \\
A/J \ar[r] & A/(I+J).
}$$
Let $\bar P$ be as in (i) and let $\bar B$ be a
finitely presented {\'e}tale $A/I$-algebra; by claim 
\ref{cl_principle} we can find an almost finitely generated 
projective $A/J$-module $\bar P_1$ with an isomorphism 
$\beta:\bar P_1/I\bar P_1\stackrel{\sim}{\to}\bar P/J\bar P$
and a finitely presented {\'e}tale $A/J$-algebra $\bar B_1$
that lifts $\bar B/J\bar B$. By proposition \ref{prop_equiv.2-prod}, 
the datum $(\bar P,\bar P_1,\beta)$ determines a unique 
almost finitely generated projective $A/(I\cap J)$-module 
$\bar P_2$; likewise, using corollary \ref{cor_equiv.2-prod} 
we obtain an {\'e}tale almost finitely presented
$A/(I\cap J)$-algebra $\bar B_2$ that lifts $\bar B$.
Next, let $K:=\bigcap_{n>0}f^nA$ and set $N:=K\cap J\cap I$; 
we have a cartesian diagram
$$
\xymatrix{ A/N \ar[r] \ar[d] & A/(I\cap J) \ar[d] \\
A^\wedge \ar[r]& A^\wedge/(I\cap J)A^\wedge.
}$$
Due to theorem \ref{th_liftmod}(i.b) we can lift 
$\bar P_2\otimes_AA^\wedge$ to an almost finitely generated 
projective $A^\wedge/(I\cap J)^2A^\wedge$-module 
$\bar P{}^\wedge_2$. For the same reason, 
$\bar P{}^\wedge_2\otimes_AA/f^2A$
can be lifted to a compatible family $(\bar Q_n~|~n\in\N)$,
where $\bar Q_n$ is almost finitely generated projective
over $A/f^{n+2}A$ for every $n\in\N$. Finally, by theorem
\ref{th_for.A-mods}, the projective limit $Q$ of the system
$(\bar Q_n~|~n\in\N)$ is an almost finitely generated
projective $A^\wedge$-module. By construction, there
is an isomorphism 
$\bar\beta:Q/f^2Q\stackrel{\sim}{\to}
\bar P{}^\wedge_2\otimes_AA/f^2A$; by assertion (ii) it
follows that there exists an isomorphism 
$\beta:Q/(I\cap J)^2Q\stackrel{\sim}{\to}\bar P{}^\wedge_2$
that lifts $\bar\beta\otimes_A\one_{A/fA}$.
By proposition \ref{prop_equiv.2-prod}, the datum 
$(\bar P_2,Q,\beta\otimes_A\one_{A/(I\cap J)})$
determines a unique almost finitely generated projective 
$A/N$-module $P_1$; since $N^2=0$, $P_1$
can be further lifted to an almost finitely generated
projective $A$-module $P$, so assertion (i) holds in 
this case. The proof of assertion (iii) is analogous,
but easier : we need to show that $\bar B_2$ lifts to
an almost finitely presented {\'e}tale $A$-algebra $B$; to this aim,
it suffices to show that $\bar B_2$ lifts to an almost
finitely presented {\'e}tale $A/N$-algebra $B_1$, since 
in that case $B_1$ can be lifted to an almost finitely
presented {\'e}tale $A$-algebra $B$, by theorem \ref{th_liftmod}(ii). 
To obtain $B_1$ it suffices to find an almost finitely 
presented {\'e}tale $A^\wedge$-algebra $B^\wedge$ with an isomorphism 
$\beta:B^\wedge\otimes_AA/(I\cap J)\stackrel{\sim}{\to}
\bar B_2\otimes_AA^\wedge$; indeed, in this case the
datum $(B^\wedge,\bar B_2,\beta)$ determines a unique {\'e}tale
almost finitely presented $A/N$-algebra in view of corollary 
\ref{cor_equiv.2-prod}. Finally, we consider the natural functors
$$
A^\wedge\Et_\mathrm{afp}\to 
A^\wedge/(I\cap J)A^\wedge\Et_\mathrm{afp}\to
A/I\Et_\mathrm{afp}.
$$
By theorem \ref{th_lift.etale.cpte}, the composition of these 
two functors is an equivalences of categories and the rightmost 
functor is fully faithful by (ii), so the leftmost functor is an
equivalence, thus $B_1$ as sought can be found, which concludes 
the proof of (iii) in this case.
\end{pfclaim}

\begin{claim}\label{cl_can.lift.for.fg}
Assertions (i) and (iii) hold if $I$ is a finitely generated 
ideal. 
\end{claim}
\begin{pfclaim} We proceed by induction on the number $n$
of generators of $I$, the case $n=1$ being covered by claim
\ref{cl_principle}. So suppose $n>1$ and let 
$f_1,...,f_n\in I_*$ be a finite set of generators of $I$.
Let $\bar P$ be as in (i) and $\bar B$ any almost finitely
presented {\'e}tale $A/I$-algebra. We let $A':=A/f_1A$ and 
$J:=\Img(I\to A')$. By lemma \ref{lem_plenty}(iii) the pair 
$(A',J)$ is again tight henselian, so by inductive assumption 
we can find lift $\bar P$ (resp. $\bar B$) to an almost finitely 
generated projective $A'$-module (resp. to an almost finitely
presented {\'e}tale $A'$-algebra) $P'$ (resp. $B'$). Thence we 
apply claim \ref{cl_principle} to further lift $P'$ to a 
module $P$ (resp. to an algebra $B$) as stated. 
\end{pfclaim}

Let now $(A,I)$ be a general tight henselian pair; we can
find a henselian pair $(R,J)$ such that $R^a=A$, $J^a=I$
and $J\subset\fm_0$. Denote by 
$(R^\mathrm{h},\fm_0R^\mathrm{h})$ the henselization
of the pair $(R,\fm_0R)$, and let 
$\bar R{}^\mathrm{h}:= R^\mathrm{h}/JR^\mathrm{h}$.
Let us write $J=\colim{\lambda\in\Lambda}J_\lambda$, where 
$J_\lambda$ runs over the filtered family of all finitely 
generated subideals of $J$; set $R_\lambda:=R/J_\lambda$ and 
$R_\lambda^\mathrm{h}:=R^\mathrm{h}\otimes_RR_\lambda$ for 
every $\lambda\in\Lambda$. Furthermore, let $X:=\Spec\,R$,
$X_\lambda:=\Spec\,R_\lambda$, 
$X^\mathrm{h}:=\Spec\,R^\mathrm{h}$, 
$X^\mathrm{h}_\lambda:=X^\mathrm{h}\times_XX_\lambda$,
$U:=X\setminus V(\fm_0)$, 
$\phi^\mathrm{h}:U^\mathrm{h}:=U\times_XX^\mathrm{h}\to U$,
$\phi_\lambda^\mathrm{h}:
U_\lambda^\mathrm{h}:=U\times_XX^\mathrm{h}_\lambda\to 
U_\lambda:=U\times_XX_\lambda$
the natural morphisms of schemes. Now, let 
$\bar P$ be an almost finitely generated $A/I$-module
and $\bar B$ an almost finitely presented {\'e}tale 
$A/I$-algebra; $\bar P$ induces a quasi-coherent 
$\cO_X^a$-module $\bar P{}^\sim$, and by restriction we obtain 
a quasi-coherent $\cO_U$-module $\bar P{}^\sim_{|U}$. Furthermore, 
since $\bar P$ is almost finitely presented, we see that 
$\bar P{}^\sim_{|U}$ is finitely presented; by 
\cite[Ch.IV, Th.8.5.2(ii)]{EGAIV-3} it follows that for some 
$\lambda_0\in\Lambda$ there exists a quasi-coherent
finitely presented module $\cP$ on $U_{\lambda_0}$ 
whose restriction to the closed subset $U$ agrees with 
$\bar P{}^\sim_{|U}$. By restricting further, we can even
achieve that $\cP$ be locally free of finite rank
(\cite[Ch.IV, Prop.8.5.5]{EGAIV-3}). Similarly, we can
find a locally free $\cO_{\lambda_0}$-algebra $\cB$ 
such that $\cB_{|U}\simeq\bar B{}^\sim_{|U}$.

\begin{claim}\label{cl_liftalittle} 
For every almost finitely generated projective
$(\bar R{}^\mathrm{h})^a$-module $\bar Q$ there exists an almost 
finitely generated projective $(R^\mathrm{h})^a$-module $Q$ 
such that $Q\otimes_{(R^\mathrm{h})^a}(\bar R{}^\mathrm{h})^a
\simeq\bar Q$.
\end{claim}
\begin{pfclaim} By theorem \ref{th_liftmod}(i.b) we know that
$\bar Q$ lifts to an almost finitely generated projective
module $Q_1$ over $(R^\mathrm{h}/J^2R^\mathrm{h})^a$.
Then $Q_2:=Q_1\otimes_A(A/\fm_0^2A)$ is an almost finitely
generated projective $(R^\mathrm{h}/\fm_0^2R^\mathrm{h})^a$-module, 
therefore by claim \ref{cl_can.lift.for.fg} we can lift $Q_2$ to 
an almost finitely generated projective $R^\mathrm{h}$-module $Q$
(notice that $(R^\mathrm{h},\fm_0^2)$ is still a henselian pair).
It remains only to show that $Q$ is a lifting of $\bar Q$;
however, by construction we have 
$Q_1\otimes_{(R^\mathrm{h})^a}(R^\mathrm{h}/\fm_0^2R^\mathrm{h})^a\simeq
Q\otimes_{(R^\mathrm{h})^a}(R^\mathrm{h}/\fm_0^2R^\mathrm{h})^a$,
so it follows from (ii) that $Q_1\simeq Q/I^2Q$, whence the claim.
\end{pfclaim}

\begin{claim}\label{cl_liftalittlemore} The natural functor 
$(R^\mathrm{h})^a\Et_\mathrm{afp}\to
(\bar R{}^\mathrm{h})^a\Et_\mathrm{afp}$ is an
equivalence of categories.
\end{claim}
\begin{pfclaim} By claim \ref{cl_can.lift.for.fg}, the
natural functors 
$(\bar R{}^\mathrm{h})^a\Et_\mathrm{afp}\to
(R^\mathrm{h}/\fm_0R^\mathrm{h})^a\Et_\mathrm{afp}$
and $(R^\mathrm{h})^a\Et_\mathrm{afp}\to
(R^\mathrm{h}/\fm_0R^\mathrm{h})^a\Et_\mathrm{afp}$
are equivalences of categories. The claim is a formal consequence.
\end{pfclaim}

We apply claim \ref{cl_liftalittle} with
$\bar Q:=\bar P\otimes_A(R^\mathrm{h})^a$; let $Q^\sim$ be
the quasi-coherent $\cO_{X^\mathrm{h}}^a$-module associated to
$Q$. By construction, the restriction $Q^\sim_{|U^\mathrm{h}}$
is a quasi-coherent $\cO_{U^\mathrm{h}}$-module of finite
presentation, and we have an isomorphism 
$\bar\beta:Q^\sim_{|U^\mathrm{h}}\stackrel{\sim}{\to}
\phi^{\mathrm{h}*}(\bar P{}^\sim_{|U})$. It then follows by 
\cite[Ch.IV, Cor.8.5.2.5]{EGAIV-3} that there exists some 
$\mu\in\Lambda$ with $X_\mu\subset X_{\lambda_0}$, such
that the isomorphism $\bar\beta$ extends to an isomorphism
$\beta_\mu:Q^\sim_{|U^\mathrm{h}_\mu}\stackrel{\sim}{\to}
\phi_\mu^{\mathrm{h}*}(\cP_{|U_\mu})$. Similarly, by claim
\ref{cl_liftalittlemore}, we can find an almost finitely
presented {\'e}tale $(R^\mathrm{h})^a$-algebra $C$ with an
isomorphism $\gamma_\mu:C^\sim_{|U^\mathrm{h}_\mu}\simeq
\phi^{\mathrm{h}*}(\bar B{}^\sim_{|U^\mathrm{h}_\mu})$.
According to 
\eqref{subsec_affine.alm.schemes}, the global section functors: 
$$
\cO_{X_\lambda}^a\Mod_\mathrm{qcoh}\to R^a_\lambda\Mod\qquad
\cO_{X^\mathrm{h}_\lambda}^a\Mod_\mathrm{qcoh}
\to(R^\mathrm{h}_\lambda)^a\Mod
$$
are equivalences. Clearly the localization functors:
$$
\cO_{U_\lambda}\Mod_\mathrm{qcoh}\to
\cO_{U_\lambda}^a\Mod_\mathrm{qcoh}\qquad
\cO_{U^\mathrm{h}_\lambda}\Mod_\mathrm{qcoh}\to
\cO_{U^\mathrm{h}_\lambda}^a\Mod_\mathrm{qcoh}
$$
are equivalences as well, and similarly for the corresponding
categories of algebras. By proposition 
\ref{prop_descent.hensel}(i) it follows that the natural diagrams:
$$
\xymatrix{ R^a_\lambda\Mod \ar[r] \ar[d] &
(R^\mathrm{h}_\lambda)^a\Mod \ar[d] \\
\cO_{U_\lambda}\Mod_\mathrm{qcoh} 
\ar[r]^-{\phi_\lambda^{\mathrm{h}*}} &
\cO_{U^\mathrm{h}_\lambda}\Mod_\mathrm{qcoh}
}\qquad
\xymatrix{ R^a_\lambda\Alg \ar[r] \ar[d] &
(R^\mathrm{h}_\lambda)^a\Alg \ar[d] \\
\cO_{U_\lambda}\Alg_\mathrm{qcoh} 
\ar[r]^-{\phi_\lambda^{\mathrm{h}*}} &
\cO_{U^\mathrm{h}_\lambda}\Alg_\mathrm{qcoh}
}$$
are $2$-cartesian for every $\lambda\in\Lambda$.
Hence, the datum 
$(\cP_{|U_\mu},Q\otimes_{R^a}R^a_\mu,\beta_\mu)$
determines uniquely a $R^a_\mu$-module $P_\mu$ that
lifts $\bar P$, and the datum 
$(\cB_{|U_\mu},C\otimes_{R^a}R^a_\mu,\gamma_\mu)$
determines a $R^a_\mu$-algebra $B_\mu$ that lifts $\bar B$. 
Furthermore, since the natural morphism 
$U_\mu\amalg X^\mathrm{h}_\mu\to X_\mu$ is faithfully
flat, it follows from remark \ref{rem_descent}(ii) that 
$P_\mu$ and $B_\mu$ are almost finitely generated projective
over $R^a_\mu$. For the same reasons, $B_\mu$ is
unramified, hence {\'e}tale over $R^a_\mu$. Finally, we apply 
claim \ref{cl_can.lift.for.fg} to the henselian pair 
$(R^a,J^a_\mu)$ to lift $P_\mu$ and $B_\mu$ all the way to 
$A$, thereby concluding the proof of the theorem.
\end{proof}

\begin{lemma}\label{lem_lift.morphs} Suppose that $\tilde\fm$ 
has homological dimension $\leq 1$, and let $(A,I)$ be a tight 
henselian pair, $\bar A:=A/I$, $\bar Q$ an almost finitely 
generated projective $\bar A$-module, $M$ an $A$-module and 
$\bar\phi:\bar Q\to M/IM$ an $\bar A$-linear epimorphism. 
Then there exists an almost finitely generated projective 
$A$-module $Q$ and a morphism $\phi:Q\to M$ such that 
$\phi\otimes_A\one_{\bar A}=\bar\phi$.
\end{lemma}
\begin{proof} We begin with the following special case:
\begin{claim} The lemma holds if $I^2=0$. Furthermore,
in this case $\phi$ is an epimorphism.
\end{claim}
\begin{pfclaim} First of all, notice that $M/IM$ is almost finitely
generated, hence the same holds for $M$, in view of lemma 
\ref{lem_flat.proj}(ii). If now $Q$ is an almost finitely
generated projective $A$-module and $\phi:Q\to M$ is a morphism
that lifts $\bar\phi$, we have 
$\Coker(\phi)\otimes_AA/I\simeq\Coker\,\bar\phi=0$, whence
$\Coker\,\phi=0$ by lemm \ref{lem_Naka}. In other words,
the second assertion follows from the first.
Define the $A$-module $N$ as the fibre
product in the cartesian diagram of $A$-modules:
$$
\xymatrix{ N \ar[r]^\alpha \ar[d]_\beta & \bar Q \ar[d]^{\bar\phi} \\
M \ar[r]^-\pi & M/IM
}$$
(where $\pi$ is the natural projection). Notice that
$IN=\Ker\,\alpha$; indeed, clearly $\alpha(IN)=0$,
and on the other hand $\beta(IN)=IM=\Ker\,\pi\simeq\Ker\,\alpha$.
We derive an isomorphism $\bar\psi:\bar Q\stackrel{\sim}{\to}N/IN$,
and clearly it suffices to find a morphism $Q\to N$ that
lifts $\bar\psi$. Under our current assumptions, theorem 
\ref{th_lifts.proj.hensel}(i) provides an almost finitely 
generated projective $A$-module $Q_1$ such that 
$Q_1\otimes_A\bar A\simeq\bar Q$, which in turns determines
an extension of $A$-modules 
$\underline E:=(0\to I\otimes_A\bar Q\to Q_1\to\bar Q\to 0)$.
Furthermore, $\bar\psi$ induces an epimorphism 
$\chi:I\otimes_A\bar Q\to IN$, whence an extension 
$\chi *\underline E:=(0\to IN\to Q_2\to\bar Q\to 0)$.
Another such extension is defined by 
$\underline F:=(0\to IN\to N\stackrel{\alpha}{\to}\bar Q\to 0)$.
However, any extension $\underline X$ of $\bar Q$ by $IN$ induces 
a morphism $u(\underline X):I\otimes_{\bar A}\bar Q\to IN$, defined
as in \eqref{subsec_webeginby}. Directly on the definition
one can check that $u(\underline X)$ depends only on the
class of $\underline X$ in $\Ext^1_A(\bar Q,IN)$, and moreover,
if $\underline Y$ is any other such extension, then 
$u(\underline X+\underline Y)=u(\underline X)+u(\underline Y)$
(where $\underline X+\underline Y$ denotes the Baer sum
of the two extensions). We can therefore compute:
$u(\chi*\underline E-\underline F)=\chi\circ u(\underline E)-
u(\underline F)$; but the definition of $\underline E$ is such 
that $u(\underline E)=\one_{I\otimes_{\bar A}\bar Q}$ and
by inspecting the construction of $\underline F$ we get
$u(\underline F)=\chi$. So finally 
$u(\chi*\underline E-\underline F)=0$; this means that 
$\chi*\underline E-\underline F$ is an extension of 
$\bar A$-modules, that is, its class is contained in the
subgroup $\Ext^1_{\bar A}(\bar Q,IN)\subset\Ext^1_A(\bar Q,IN)$.
Notice now that $\Ext^2_{\bar A}(\bar Q,\Ker\,\chi)=0$,
due to lemma \ref{lem_Exts.are.same}(i),(ii); since $\chi$
is an epimorphism, it then follows that the induced map
$\Ext^1_{\bar A}(\bar Q,\chi)$ is surjective. Hence there
exists an extension 
$\underline X:=(0\to I\otimes_A\bar Q\to Q_3\to\bar Q\to 0)$
of $\bar A$-modules, such that 
$\chi*\underline X=\chi*\underline E-\underline F$, {\em i.e.}
$\underline F=\chi*(\underline E-\underline X)$. Say
$\underline E':=\underline E-\underline X=
(0\to I\otimes_A\bar Q\to Q\to\bar Q\to 0)$; by construction,
$u(\underline E')=u(\underline E)$, so $Q$ is a flat 
$A$-module (see \eqref{subsec_webeginby}) that lifts $\bar Q$. 
Finally, $Q$ is almost finitely generated projective by lemma 
\ref{lem_flat.proj}(i),(ii). The push-out 
$\underline E'\to\chi*\underline E'$ delivers the promised
morphism $Q\to N$. 
\end{pfclaim}

Next, since the pair $(A,I^{n+1})$ is still tight henselian for
every $n\in\N$, an easy induction shows that the lemma holds
when $I$ is a nilpotent ideal. For the general case, pick
$n>0$ and a finitely generated subideal $\fm_0\subset\fm$ with
$I^n\subset\fm_0A$; by the foregoing we can find an almost
finitely generated projective $A/I^{n+1}$-module $Q_{n+1}$
and an epimorphism $\phi_{n+1}:Q_{n+1}\to M/I^{n+1}M$. 
By theorem \ref{th_lifts.proj.hensel}(i), we can lift $Q_{n+1}$
to an almost finitely generated projective $A$-module $Q$;
the obstruction to lifting the induced morphism 
$Q\to M/I^{n+1}M$ to a morphism $Q\to M$ is a class
$\omega\in\Ext^1_A(Q,I^{n+1}M)$; by the argument of claim 
\ref{cl_vanishing.trick} we see that the image of $\omega$
in $\Ext^1_A(Q,IM)$ vanishes, whence the claim.
\end{proof}

\begin{corollary}\label{cor_crit.for.a-proj}
Let $A$ be a $V^a$-algebra, $I\subset\rad(A)$ a 
tight ideal and set $\bar A:=A/I$; let $P$ an almost finitely 
generated $A$-module, such that $\bar P:=P\otimes_A\bar A$ is 
an almost projective $\bar A$-module. Then the following conditions 
are equivalent:
\begin{enumerate}
\item
$P$ is an almost projective $A$-module.
\item
$P$ is a flat $A$-module.
\item
$P$ is an almost finitely presented $A$-module and 
$\Tor_1^A(P,\bar A)=0$.
\item
The natural morphism $P^*\to(\bar P)^*$ 
is an epimorphism.
\end{enumerate}
\end{corollary}
\begin{proof} Clearly (i) implies all the other assertions,
so it suffices to show that each of the assertions (ii)-(iv)
implies (i).
Let us first remark that, in view of lemmata
\ref{subsec_reduce.to.count}, \ref{lem_reduce.Tor.to.count} 
and theorem \ref{th_countably.pres}, we can assume that
the homological dimension of $\tilde\fm$ is $\leq 1$.
Furthermore, let $(A^\mathrm{h},I^\mathrm{h})$ be the
henselization of the pair $(A,I)$; according to 
\eqref{subsec_hensel.is.fflat}, the morphism 
$A\to A^\mathrm{h}$ is faithfully flat. In view of remark 
\ref{rem_descent}(ii) and lemma \ref{lem_alhom}(i) we deduce 
that each of the statements (i)-(iii) on $P$ and $\bar P$ is 
equivalent to the corresponding statement 
(i)$^\mathrm{h}$-(iii)$^\mathrm{h}$ made on the 
$A^\mathrm{h}$-module $P\otimes_AA^\mathrm{h}$ and the 
$(\bar A\otimes_AA^\mathrm{h})$-module
$\bar P\otimes_AA^\mathrm{h}$. Moreover, one checks easily
that (iv)$\Rightarrow$(iv)$^\mathrm{h}$. Thus, up to replacing 
$(A,I)$ by $(A^\mathrm{h},I^\mathrm{h})$ we can assume that 
$(A,I)$ is a tight henselian pair (notice as well that 
$I^\mathrm{h}=IA^\mathrm{h}$). By lemma 
\ref{lem_lift.morphs} we can find an almost finitely generated 
projective $A$-module $Q$ and a morphism $\phi:Q\to P$ such 
that $\phi\otimes_A\one_{A/I}$ is an isomorphism. By lemma 
\ref{lem_Naka} we deduce easily that $\phi$ is an epimorphism,
Suppose now that (ii) holds; then to deduce (i) it remains only 
to prove the following :
\begin{claim}
$\Ker\,\phi=0$.
\end{claim}
\begin{pfclaim} In view of proposition \ref{prop_eval.ideal}(v), 
it suffices to show that $\cE:=\cE_{Q/A}(x)=0$ for every 
$x\in\Ker\,\phi_*$ (see definition \ref{def_dual_module}(iv)). 
However, since $Q$ is $A$-flat, we can compute:
$$
0=\Tor^A_1(A/\cE,P)\simeq
\Ker(\Ker(\phi)\otimes_A(A/\cE)\to Q/\cE Q)\simeq
((\Ker\,\phi)\cap\cE Q)/(\cE\cdot\Ker\,\phi)
$$
that is, 
$(\Ker\,\phi)\cap\cE Q=\cE\cdot\Ker\,\phi$. By proposition
\ref{prop_eval.ideal}(v) we have $x\in(\cE Q)_*$; on the
other hand we also know that $\Ker\,\phi\subset IQ$, whence
$x\in(\cE I Q)_*$. We apply once again proposition 
\ref{prop_eval.ideal}(v) to derive $\cE=\cE I$, so finally
$\cE=0$ in view of lemma \ref{lem_Naka}. 
\end{pfclaim}

Next, assume (iii); we compute: 
$0=\Tor_1^A(P,\bar A)\simeq
\Ker(\Ker(\phi)\otimes_A\bar A\to\Ker(\phi\otimes_A\one_{\bar A}))=
\Ker(\phi)\otimes_A\bar A$. Since $P$ is almost finitely presented,
$\Ker\,\phi$ is almost finitely generated by lemma 
\ref{lem_long.forgotten}, whence $\Ker\,\phi=0$ by lemma
\ref{lem_Naka}, so (i) holds.
Finally, let $\phi^*:P^*\to Q^*$ be the transposed of the 
morphism $\phi$; by lemma \ref{lem_alhom}(i), the natural 
morphism $\psi:Q^*\to(Q/IQ)^*$ is an epimorphism. The 
composition $\psi\circ\phi^*$ factors through the transposed 
morphism 
$(\phi\otimes_A\one_{\bar A})^*:(\bar P)^*\to(Q/IQ)^*$, so it 
is an epimorphism when (iv) holds; then lemma \ref{lem_Naka}
implies easily that $\phi^*$ is an epimorphism; since it
is obviously a monomorphism, we deduce that
$P^*\stackrel{\sim}{\to} Q^*$. so 
$(\phi^*)^*:Q\simeq(Q^*)^*\to(P^*)^*$ is an isomorphism; 
since the latter factors through the natural morphism 
$P\to(P^*)^*$, we see that $\phi$ is a monomorphism
and (i) follows.
\end{proof}

\subsection{Smooth locus of an affine almost scheme}\label{sec_smoothloci}
Thoughout this section we fix a $V^a$-algebra $A$ and set
$S:=\Spec\,A$. Let $X$ be an affine $S$-scheme. We often identify
$X$ with the functor it represents:
$$
X:A\Alg\to\Set\qquad T\mapsto X(T^o):=\Hom_{A\Alg^o}(T^o,X).
$$
The usual argument from faithfully flat descent shows that
$X$ is a sheaf for the fpqc topology of $A\Alg^o$. In this
section we aim to study, for every such $X$, the {\em smooth
locus of $X$ over $S$}, which will be a certain natural subsheaf of
$X$. The starting point is the following:

\begin{definition}\label{def_smooth.locus}
\index{Almost scheme(s)!$X_\sm$ : smooth locus of
an|indref{def_smooth.locus}{}, \indref{def_smooth.for.proj}}
Let $S$ and $X$ be as in \eqref{sec_smoothloci}.
Given an affine $S$-scheme $T$ and $\sigma\in X(T)$, we
say that {\em $\sigma$ lies in the smooth locus of $X$ over $S$}
if the following two conditions hold:
\begin{enumerate}
\renewcommand{\labelenumi}{(\alph{enumi})}
\item
$H_1(L\sigma^*\L_{X/S}^a)=0$ and
\item
$H_0(L\sigma^*\L_{X/S}^a)$ is an almost finitely generated projective
$\cO_T$-module.
\renewcommand{\labelenumi}{(\roman{enumi})}
\end{enumerate}
We denote by $X_\sm(T)\subset X(T)$ the subset of all
the $T$-sections of $X$ that lie in the smooth locus of $X$ over $S$.
\end{definition}

\sset\subsubsection{}
Using remark \ref{rem_descent}(iii) one sees that $X_\sm$ is a
subsheaf of $X$. Just as for usual schemes, in order
to get a handle on the smooth locus $X_\sm$, one often needs to
assume that the almost scheme $X$ satisfies some finiteness
conditions. For our purposes, the following will do:

\begin{definition}\label{def_af-pres.scheme}
\index{Almost scheme(s)!almost finitely presented|indref{def_af-pres.scheme}}
We say that the affine almost $S$-scheme
$X$ is {\em almost finitely presented\/} if there exists an
almost finitely generated projective $\cO_S$-module, and an almost
finitely generated ideal $J$ of $S_P:=\Sym^\bullet_{\cO_S}(P)$, such
that $X\simeq\Spec\,S_P/J$.
\end{definition}

\begin{lemma}\label{lem_smooth-morphs}
Let $X=\Spec\,S_F/J$, where $S_F:=\Sym^\bullet_{\cO_S}(F)$ for
some flat $\cO_S$-module $F$, and $J$ is any ideal. Then there is a
natural isomorphism in $\sD(\cO_X\Mod)$:
$$
\tau_{[-1}\L^a_{X/S}\simeq(0\to J/J^2\to\cO_X\otimes_{\cO_S}P\to 0).
$$
\end{lemma}
\begin{proof} Let us remark the following:
\begin{claim}\label{cl_smooth-morphs}
With the notation of the lemma, there is a natural isomorphism
$\L^a_{\Spec\,S_F/S}\simeq S_F\otimes_{\cO_S}F[0]$ in $\sD(S_F\Mod)$.
\end{claim}
\begin{pfclaim} By a theorem of Lazard (\cite[Ch.I, Th.1.2]{La}), 
every flat $\cO_{S*}$-module is the filtered colimit of a family
of free modules of finite rank; in particular this holds for $F_!$;
since both functors $\Sym^\bullet$ and $\L$ commute with filtered
colimits (proposition \ref{prop_L-and-colimits}), we can then reduce
to the case where $S_F\simeq\cO_S[T_1,...,T_n]$ for some $n\in\N$.
In this case, we have natural isomorphisms
$\L_{\Spec\,S_F/S}^a\simeq\L^a_{\cO_{S*}[T_1,...,T_n]/\cO_{S*}}$
in light of proposition \ref{prop_Gabber}(ii). The claim follows.
\end{pfclaim}

Using claim \ref{cl_smooth-morphs}, the assertion can be shown
as in the proof of \cite[Ch.III, Cor.1.2.9.1]{Il}.
\end{proof}

\sset\subsubsection{}\label{subsec_invert-t}
Let $X$ be an almost finitely presented $S$-scheme and $t\in\fm$
any element. Then $\cO_S[t^{-1}]$ is a (usual) $V[t^{-1}]$-algebra,
and we let $S_t:=\Spec\,\cO_S[t^{-1}]$, $X_t:=X\times_SS_t$.
Both $S_t$ and $X_t$ are represented by (usual) affine schemes over
$\Spec\,V[t^{-1}]$, and obviously $X_t$ is finitely presented
over $S_t$. Using lemma \ref{lem_smooth-morphs} it is also easy
to see that the subfunctor $X_{\sm,t}:=X_\sm\cap X_t$ of the functor
$X$ is represented by the smooth locus of $X_t$ over $S_t$, which is an
open subscheme of the latter scheme.

\sset\subsubsection{}\label{subsec_add-I}
In the situation of \eqref{subsec_invert-t}, suppose moreover that
$I\subset\cO_S$ is a given ideal, and that $t$ is regular in $\cO_S$.
Let $R:=\cO_{S*}$; then $t$ is a non-zero-divisor in $R$, and
we have a well defined $(t,I_*)$-adic topology on 
$X_t(R[t^{-1}])=X_t(S_t)$. Furthermore, it is clear that
the restriction map:
\set\begin{equation}\label{eq_add-I}
X(S)\to X_t(S_t)\qquad \sigma\mapsto\sigma_*[t^{-1}]
\end{equation}
is injective. Consequently we can endow $X(S)$ with the
{\em $(t,I)$-adic topology}, defined as the topology induced
by the $(t,I_*)$-adic topology of $X_t(S_t)$.

\begin{lemma}\label{lem_add-I}
In the situation of \eqref{subsec_add-I}, the map
\eqref{eq_add-I} is an open imbedding for the respective
$(t,I)$-adic and $(t,I_*)$-adic topologies.
\end{lemma}
\begin{proof} Let us write $X=\Spec\,S_P/J$, where $S_P$ is the
symmetric algebra of an almost finitely generated projective 
$\cO_S$-module $P$, and set $\cP:=\Spec\,S_P$. Then $X$ is a
closed subscheme of $\cP$ and $X_t$ is a closed subscheme of
$\cP_t$, which is a vector bundle of finite rank over $S_t$.
The $(t,I_*)$-adic topology of $X_t(S_t)$ is induced by the
$(t,I_*)$-adic topology of $\cP_t(S_t)$, and consequently the
$(t,I)$-adic topology of $X(S)$ is induced by the $(t,I)$-adic
topology of $\cP(S)$. Since the commutative diagram of sets
\set\begin{equation}\label{eq_diagr-of-sets}
{\diagram X(S) \ar[r] \ar[d] & X_t(S_t) \ar[d] \\
           \cP(S) \ar[r] & \cP_t(S_t)
\enddiagram}
\end{equation}
is cartesian, we reduce to showing that the restriction map
$\cP(S)\to\cP_t(S_t)$ is open. To this aim, we set
$P[t^{-1}]^*:=\Hom_{\cO_S}(P,\cO_S[t^{-1}])$ and we consider the
diagram
$$
\xymatrix{
\Hom_{\cO_S}(P,\cO_S) \ar[r]^-\sim \ar[d] & \cP(S) \ar[d] \\
P[t^{-1}]^* \ar[r]^-\sim & \cP_t(S_t)
}$$
whose horizontal arrows are given by the rule:
$\phi\mapsto\Sym_{\cO_S}^\bullet\phi$. Since $P[t^{-1}]^*$ is a
finitely generated $\cO_S[t^{-1}]$-module, it is endowed with a
well defined $(t,I_*)$-preadic topology. We define a linear 
{\em $(t,I)$-adic topology} on $\Hom_{\cO_S}(P,\cO_S)$, by
declaring that the system of submodules 
$(\Hom_{\cO_S}(P,t^nI)~|~n\in\N)$ forms a cofinal family of
open neighborhoods of zero.
\begin{claim}\label{cl_horizontal} With these $(t,I)$-adic and
$(t,I_*)$-adic topologies, \eqref{eq_diagr-of-sets} is a diagram
of continuous maps, and the horizontal arrows are homeomorphisms.
\end{claim}
\begin{pfclaim} In case $P[t^{-1}]\simeq R[t^{-1}]^n$ for some $n\in\N$,
we have $\cP_t(S_t)\simeq R[t^{-1}]^n$, and then the bottom
arrow of \eqref{eq_diagr-of-sets} is a homeomorphism, essentially
by definition. The general case can be reduced to the case of
a free module, by writing $R[t^{-1}]^n=P[t^{-1}]\oplus Q$ for
some projective $R[t^{-1}]$-module $Q$, and remarking that the
$(t,I_*)$-adic topologies are compatible with cartesian products.
To prove that the top arrow is a homeomorphism, it suffices therefore
to show that the $(t,I)$-adic topology of $\Hom_{\cO_S}(P,\cO_S)$
is induced by the $(t,I_*)$-adic topology of $P[t^{-1}]^*$.
To this aim, pick a finitely generated $R$-module
$P_0\subset P_*$ with $tP_*\subset P_0$. Clearly
$\Hom_R(P_0,\cO_S[t^{-1}])=P[t^{-1}]^*$, and again by reducing
to the case where $P[t^{-1}]$ is free, one verifies that the
$(t,I_*)$-adic topology on $P[t^{-1}]^*$ admits the system
$(\Hom_R(P_0,t^nI_*)~|~n\in\N)$ as a cofinal family of open
neighborhoods of zero. For every $n\in\N$ set 
$U_n:=\{\phi:P\to\cO_S~|~\phi(P_0)\subset t^nI_*\}$; we have
\set\begin{equation}\label{eq_vertical}
\Hom_{\cO_S}(P,t^nI)\subset U_n\subset\Hom_{\cO_S}(P,t^{n-1}I)\qquad
\text{for every $n\in\N$}
\end{equation}
which implies the claim.
\end{pfclaim}

In view of claim \ref{cl_horizontal}, we are reduced to showing that
the map $\Hom_{\cO_S}(P,\cO_S)\to P[t^{-1}]^*$ is open. But this
is again a direct consequence of \eqref{eq_vertical}.
\end{proof}

\begin{proposition}\label{prop_smooth-is-topo}
Let $X$, $T$ be affine $S$-schemes and $\sigma\in X(T)$
a $T$-section, $I\subset\rad(\cO_T)$ an ideal, and set
$T_0:=\Spec\,\cO_T/I$; suppose that the restriction 
$\sigma_0\in X(T_0)$ of $\sigma$ lies in the smooth
locus of $X$. Suppose moreover that either:
\begin{enumerate}
\renewcommand{\labelenumi}{(\alph{enumi})}
\item
$I$ is nilpotent, or
\item
$I$ is tight and $X$ is almost finitely presented over $S$.
\renewcommand{\labelenumi}{(\roman{enumi})}
\end{enumerate}
Then $\sigma\in X_\sm(T)$.
\end{proposition}
\begin{proof} Suppose that (a) holds; for any quasi-coherent 
$\cO_T$-module $\cF$, let us denote by $\mathrm{Fil}^\bullet_I\cF$
the $I$-adic filtration on $\cF$. We can write
$\tau_{[-1}L\sigma^*\L^a_{X/S}\simeq(0\to N\stackrel{\phi}{\to} P\to 0)$
for two $\cO_T$-modules $N$ and $P$, and we can assume that $P$ is almost
projective over $\cO_T$, so that the natural morphism 
$$
\tau_{[-1}(\cO_{T_0}\derotimes_{\cO_T}L\sigma^*\L^a_{X/S})\to
(0\to N/IN\stackrel{\gr^0_I\phi}{\longrightarrow} P/IP\to 0)
$$
is an isomorphism in $\sD(\cO_{T_0}\Mod)$. Hence, the assumption on
$\sigma$ means that $\gr^0_I\phi$ is a monomorphism with
almost finitely generated projective cokernel over $\cO_{T_0}$. 
We consider, for every integer $i\in\N$ the commutative diagram:
$$
\xymatrix{ \gr^i_IN \ar[rrr]^{\gr_I^i\phi} & & & \gr^i_IP \\
\gr^0_IN\otimes_{\cO_T}\gr_I^i\cO_T \ar[u]^{\alpha_i}
\ar[rrr]^-{\gr^0_I\phi\otimes_{\cO_T}\one_{\gr_I^i\cO_T}} & & &
\gr^0_IP\otimes_{\cO_T}\gr_I^i\cO_T \ar[u]_{\beta_i}
}$$
Since $P$ is almost projective (especially, flat) $\beta_i$ is
an isomorphism. Moreover, since $\gr^0_I\phi$ is a monomorphism
with almost projective cokernel, the long exact Tor sequence
shows that $\gr^0_I\phi\otimes_{\cO_T}\one_{\gr_I^i\cO_T}$ is
a monomorphism. It follows that $\alpha_i$ is a monomorphism for every
$i\in\N$, and since it is obviously an epimorphism, we deduce that
$\alpha_i$ is an isomorphism and $\gr^\bullet_I\phi$ is a monomorphism,
therefore the same holds for $\phi$. Let $C:=\Coker(N\to P)$; we
deduce easily that $\Tor_1^{\cO_T}(\cO_{T_0},C)=0$, and then it 
follows from the local flatness criterion (see \cite[Ch.8,Th.22.3]{Mat})
that $C$ is a flat $\cO_T$-module. Finally lemma 
\ref{lem_flat.proj}(i),(ii) says that $C$ is almost finitely generated
projective, whence the claim, in case (a).

Next, suppose that assumption (b) holds; by lemma
\ref{lem_smooth-morphs} there is an isomorphism :
$$
\tau_{[-1}\L^a_{X/S}\simeq
(0\to N\stackrel{\phi}{\to}Q\to 0)
$$
where $N$ is almost finitely generated and $Q$ is almost finitely
generated projective over $\cO_X$.
Since $Q$ is almost projective, we have
$\tau_{[-1}(\cO_{T_0}\derotimes_{\cO_T}L\sigma^*\L^a_{X/S})\simeq
(0\to N/IN\stackrel{\phi_0}{\to}Q/IQ\to 0)$, and by assumption
$\Ker\,\phi_0=0$ and $\Coker\,\phi_0$ is an almost finitely generated
projective $\cO_{T_0}$-module. Using the long exact
Ext sequence we deduce that $N/IN$ is almost projective.
Thus, the bottom arrow of the natural commutative diagram
$$
\xymatrix{ Q^* \ar[r]^{\phi^*} \ar[d]_\alpha & N^* \ar[d]^\beta \\
(Q/IQ)^* \ar[r]^{\phi_0^*} & (N/IN)^*
}$$
is an epimorphism. Invoking twice corollary \ref{cor_crit.for.a-proj}
we find first that $\alpha$ is an epimorphism (whence so is $\beta$),
and then that $N$ is almost projective. It then follows that
$\phi^*\otimes_{\cO_S}\one_{\cO_{S_0}}$ is an epimorphism as
well, hence the same holds for $\phi^*$, by applying Nakayama's lemma
\ref{lem_Naka} to $\Coker\,\phi^*$. The long exact Ext sequence
then shows that $\Ker\,\phi^*$ is almost projective, and the
latter is almost finitely generated as well, since $Q$ is.
Dualizing, we see that $\phi$ is a monomorphism and 
$\Coker\,\phi\simeq(\Ker\,\phi^*)^*$ is almost finitely generated
projective, which is the claim.
\end{proof}

\begin{lemma}\label{lem_cartes-diagr-sects}
Let $X$ be an affine $S$-scheme, and suppose we are given
a cartesian diagram of $\cO_S$-algebras as \eqref{eq_cartesian}.
There follows a cartesian diagram of sets:
$$
\xymatrix{ X_\sm(\Spec\,A_0) \ar[r] \ar[d] &
X_\sm(\Spec\,A_2) \ar[d] \\
X_\sm(\Spec\,A_1) \ar[r] & X_\sm(\Spec\,A_3).
}$$
\end{lemma}
\begin{proof} We have to check that every section
$\sigma\in X(\Spec\,A_0)$ whose restrictions to
$\Spec\,A_1$ and $\Spec\,A_2$ lie in the smooth locus of
$X$ over $S$, lies itself in the smooth locus. Hence, let
$\tau_{[-1}L\sigma^*\L_{X/S}\simeq(0\to N\stackrel{\phi}{\to}P\to 0)$,
for some $A_0$-modules $N$ and $P$, chosen so that $P$ is
almost projective. Hence 
$A_i\derotimes_{A_0}\tau_{[-1}L\sigma^*\L_{X/S}
\simeq(0\to N_i\stackrel{\phi_i}{\to}P_i\to 0)$ (where
$N_i:=A_i\otimes_{A_0}N$ and likewise for $P_i$, $i=1,2$).
By assumption, $\Ker\,\phi_i=0$ and $\Coker\,\phi_i$ is almost
finitely generated projective over $A_i$ ($i=1,2$).
It follows that $N_i$ is almost projective for $i=1,2$.
Using proposition \ref{prop_equiv.2-prod} we deduce that
$\Coker\,\phi$ is almost finitely generated projective
over $A_0$ and $N$ is almost projective. In particular,
$N$ is flat and consequently $N\subset N_1\oplus N_2$,
so $\Ker\,\phi=0$, and the assertion follows.
\end{proof}

\begin{lemma}\label{lem_smooth-prod}
Let $X$ and $Y$ be two affine $S$-schemes, and suppose that
$\Tor^{\cO_S}_i(\cO_X,\cO_Y)=0$ for every $i>0$.
Then we have a natural isomorphism of sheaves:
$$
(X\times_SY)_\sm\simeq X_\sm\times Y_\sm.
$$
\end{lemma}
\begin{proof} Let $\pi_X:X\times_SY\to X$ be the natural projection,
and define likewise $\pi_Y$. By theorem \ref{th_flat.base.ch}, our
assumptions imply that the natural morphism
$\pi_X^*\L^a_{X/S}\oplus\pi_Y^*\L^a_{Y/S}\to\L^a_{X\times_SY/S}$
is a quasi-isomorphism. Let now $T$ be an affine $S$-scheme and
$(\sigma,\tau)\in X\times_SY(T)=X(T)\times Y(T)$; we derive
a natural isomorphism:
$$
L(\sigma,\tau)^*\L_{X\times_SY/S}\simeq 
L\sigma^*\L_{X/S}\oplus L\tau^*\L_{Y/S}
$$
from which the claim follows straightforwardly.
\end{proof}

\begin{theorem}\label{th_lift-smooth-sects}
Assume that $\mathrm{hom.dim}_V\tilde\fm\leq 1$. Let $X$ be an 
almost finitely presented affine $S$-scheme, $I\subset\cO_S$
an ideal such that the pair $(\cO_S,I)$ is henselian,
$\fm_0\subset\fm$ a finitely generated subideal, and set
$S_n:=\Spec\,\cO_S/\fm_0^nI$ for every $n\in\N$. Then we have:
\begin{enumerate}
\item
$\Img(X_\sm(S_1)\to X_\sm(S_0))=\Img(X_\sm(S)\to X_\sm(S_0))$.
\item
Set $X_\sm^\wedge(S):=\liminv{n\in\N}\,X_\sm(S_n)$, and endow
$X^\wedge_\sm(S)$ with the corresponding pro-discrete topology.
Then the natural map $X_\sm(S)\to X^\wedge_\sm(S)$ has dense image.
\end{enumerate}
\end{theorem}
\begin{proof} We begin with an easy reduction:
\begin{claim}\label{cl_reduce-to-tight}
In order to prove (i), we can assume that $I$ is a tight ideal.
\end{claim}
\begin{pfclaim} Indeed, let $\fm_0\subset\fm$ be any subideal,
and choose $\fm_1\subset\fm$ such that $\fm_0\subset\fm_1^2$.
Set $I':=\fm_1I$, $S'_0:=\Spec\,\cO_S/I'$, $S'_1:=\Spec\,\cO_S/\fm_1I'$.
Notice that $I'$ is tight, and suppose that the assertion
is known for this ideal. By a simple chase on the commutative
diagram:
$$
\xymatrix{ X_\sm(S) \ar[r] \ddouble & X_\sm(S_0) & 
X_\sm(S_1) \ar[d] \ar[l] \\
X_\sm(S) \ar[r] & X_\sm(S'_0) \ar[u] &
X_\sm(S'_1) \ar[l]
}$$
the assertion can then be deduced for $I$ as well.
\end{pfclaim}

Let $\sigma_0\in X_\sm(S_0)$ that admits an extension
$\tilde\sigma\in X_\sm(S_1)$, pick finitely
many generators $\eps_1,...,\eps_k$ for $\fm_0$,
and define a map $\phi:I^{\oplus k}\to\fm_0I$
by the rule: $(a_1,...,a_k)\mapsto\sum^k_i\eps_i\cdot a_i$.
\begin{claim}\label{cl_two-asserts} Suppose that $I^2=0$ (so
that $I$ is an $\cO_{S_0}$-module); then
$\hExt^1_{\cO_{S_0}}(L\sigma_0^*\L^a_{X/S},\phi)$ is
onto and $\hExt^1_{\cO_{S_0}}(L\sigma_0^*\L^a_{X/S},I)^a=0$.
\end{claim}
\begin{pfclaim} For the first assertion we use the spectral sequence 
$$
E_2^{pq}:=\Ext^p_{\cO_{S_0}}(H_q(L\sigma_0^*\L^a_{X/S}),I^{\oplus k})
\Rightarrow\hExt^{p+q}_{\cO_{S_0}}(L\sigma_0^*\L^a_{X/S},I^{\oplus k})
$$
and a similar one for $\fm_0I$.
Since $\sigma_0$ lies in the smooth locus of $X$, the $E^{01}_2$
terms vanish, hence we are reduced to verifying that the map
$\Ext^1_{\cO_{S_0}}(H_0(L\sigma_0^*\L^a_{X/S}),\phi)$
is surjective. However, the cokernel of this
map is a submodule of 
$\Ext^2_{\cO_{S_0}}(H_0(L\sigma_0^*\L^a_{X/S}),\Ker\,\phi)$;
again the assumption that $\sigma_0$ lies in the smooth locus of $X$,
and $\mathrm{hom.dim}_V\tilde\fm\leq 1$, imply that the latter Ext group
vanishes by lemma \ref{lem_Exts.are.same}(ii). The same spectral sequence
argument also proves the second assertion.
\end{pfclaim}

\begin{claim}\label{cl_case-when_I-sqzero}
Assertion (i) holds if $I^2=0$.
\end{claim}
\begin{pfclaim}
We need to show that there exists a morphism $\sigma:S\to X$
extending $\sigma_0$; then $\sigma\in X_\sm(S)$ in view of
proposition \ref{prop_smooth-is-topo}. In other words, we have 
to find $\sigma^\sharp$ that fits into a morphism of extensions 
of $\cO_S$-algebras:
$$
\xymatrix{ 0 \ar[r] & 0 \ar[r] \ar[d] &
\cO_X \rdouble \ar[d]_{\sigma^\sharp} &
\cO_X \ar[r] \ar[d]^{\sigma_0^\sharp} & 0 \\
0 \ar[r] & I \ar[r] & \cO_S \ar[r] & \cO_{S_0} \ar[r] & 0.
}$$
By proposition \ref{prop_obstr2}, the obstruction to the existence
of $\sigma^\sharp$ is a class 
$$
\omega\in\hExt^1_{\cO_X}(\L^a_{X/S},I)\simeq
\hExt^1_{\cO_{S_0}}(L\sigma_0^*\L^a_{X/S},I).
$$
Likewise, the obstruction to extending $\tilde\sigma$ is a class 
$\tilde\omega\in\hExt^1_{\cO_{S_0}}(L\sigma_0^*\L^a_{X/S},\fm_0I)$,
and $\omega$ is the image of $\tilde\omega$ under the map
$\hExt^1_{\cO_{S_0}}(L\sigma_0^*\L^a_{X/S},j)$
(where we have denoted by $j:\fm_0I\to I$ the inclusion).
From claim \ref{cl_two-asserts} one deduces easily first: that
$\hExt^1_{\cO_{S_0}}(L\sigma_0^*\L^a_{X/S},j\circ\phi)$ vanishes,
and therefore, second: that $\hExt^1_{\cO_{S_0}}(L\sigma_0^*\L^a_{X/S},j)$
must already vanish. Thus $\omega=0$, and the assertion holds. 
\end{pfclaim}

\begin{claim}\label{cl_inv-smooths} Choose a finitely generated
subideal $\fm_1\subset\fm$ such that $\fm_0\subset\fm_1^2$. The
section $\sigma_0$ can be lifted to an element of
$\liminv{n\in\N}\,X_\sm(\Spec\,\cO_S/\fm_1I^n)$.
\end{claim}
\begin{pfclaim} For every $n>0$, let $T_n:=\Spec\,\cO_S/\fm_1^2I^n$
and $j_n:\Spec\,\cO_S/\fm_1I^n\to T_n$ the natural morphism; we 
construct by induction on $n\in\N$ a sequence of section 
$\sigma_n\in X(T_n)$, such that the family $(\sigma_i\circ j_i~|~i>0)$
defines an element of $\liminv{n\in\N}\,X_\sm(\Spec\,\cO_S/\fm_1I^n)$.
To this aim, we take $\sigma_1$ equal to the restriction of $\tilde\sigma$;
suppose then that $n>1$ and that $\sigma_{n-1}$ is already given.
Notice that the image $J\subset \cO_{S_n}$ of $\fm_1I^{n-1}$ satisfies
$J^2=0$. By the claim \ref{cl_case-when_I-sqzero}, it follows that
$\sigma_{n-1}\circ j_{n-1}$ extends to a section in $X_\sm(T_n)$,
and this we call $\sigma_n$.
\end{pfclaim}

Let us now show how to deduce assertion (i) from (ii). By claim
\ref{cl_reduce-to-tight} we can suppose that $I^m\subset\fm_0\cO_S$
for some $m\geq 0$ and a finitely generated subideal
$\fm_0\subset\fm$; let $\fm_1\subset\fm$ be as in claim
\ref{cl_inv-smooths}; clearly for every $n\in\N$ there exists
$n\geq 0$ such that $\fm_1I^m\subset\fm_0^nI$, hence claim
\ref{cl_inv-smooths} shows that $\sigma_0$ can be lifted to
an element of $X^\wedge_\sm(S)$ and then (ii) yields (i) trivially.
Hence, it remains only to show (ii).

\begin{claim}\label{cl_suffices_m_0-princ}
In order to prove (ii) we can assume that $\fm_0$
is a principal ideal.
\end{claim}
\begin{pfclaim} We argue by induction on the number of generators
of $\fm_0$. Thus, let $\eps_1,...,\eps_k$ be a finite system of
generators for $\fm_0$, and suppose that the assertion is known
for all ideals generated by less than $k$ elements. Let $\fm_0^{(n)}$
be the ideal generated by $\eps_1^n,...,\eps_k^n$, and set
$S_{(n)}:=\Spec\,\cO_S/\fm_0^{(n)}I$; for every $n\in\N$ there exists
$N\in\N$ such that $\fm_0^N\subset\fm_0^{(n)}\subset\fm_0^n$,
whence an isomorphism of pro-discrete spaces:
$X^\wedge_\sm(S)\simeq\liminv{n\in\N}\,X_\sm(S_{(n)})$.
Thus, we can suppose that we are given a compatible system
of sections $\sigma_n\in X_\sm(S_{(n)})$, and wish to show
that, for every $n\in\N$, there exists $\sigma\in X_\sm(S)$ whose 
restriction to $S_{(n)}$ agrees with $\sigma_n$. Fix $N>0$,
set $T:=\Spec\,\cO_S/\eps_1^NI$ and let
$\fm_1\subset\fm$ be the ideal generated by $\eps_2,...,\eps_k$.
Let also $T_{(n)}:=T\times_SS_{(n)}$ for every $n\in\N$,
and denote by $\sigma_{n|T}\in X_\sm(T_{(n)})$ the restriction
of $\sigma_n$. Clearly $T_{(n)}\simeq\Spec\,\cO_T/\fm_1^{(n)}\cO_T$
for every $n\geq N$, and the pair 
$(\cO_T,I\cO_T)$ is henselian (cp. remark \ref{rem_hensel}(v)).
Hence, by inductive assumption, for every $n\geq N$ we can find
$\sigma_T\in X(T)$ whose restriction to $T_{(n)}$ agrees with
$\sigma_{n|T}$. We can then apply claim \ref{cl_inv-smooths},
in order to find a compatible system of sections 
$(\sigma'_n\in X(\Spec\,\cO_S/\eps_1^nI)~|~n\geq N)$,
whose restriction to $\Spec\,\cO_S/\eps^{N-1}_1I$ agrees with
the restriction of $\sigma_T$. Finally, if we assume that (ii) is known
whenever $\fm_0$ is principal, we can find a section
$\sigma\in X_\sm(S)$ whose restriction to
$\Spec\,\cO_S/\eps^{N-1}_1\cO_S$ agrees with $\sigma'_N$, whence
the claim.
\end{pfclaim}

For a given $\eps\in\fm$, set 
$K(\eps):=\bigcup_{n\in\N}\Ann_{\cO_S}(\eps^n)$.
\begin{claim}\label{cl_about-K(eps)}
$K(\eps)_*=\bigcup_{n\in\N}\Ann_{\cO_{S*}}(\eps^n)$.
\end{claim}
\begin{pfclaim} Clearly we have only to show the inclusion
$K(\eps)_*\subset\bigcup_{n\in\N}\Ann_{\cO_{S*}}(\eps^n)$.
By applying the left exact functor $M\mapsto M_*$ to
the left exact sequence 
$0\to\Ann_{\cO_S}(\eps^n)\to\cO_S\stackrel{\eps^n}{\to}\cO_S$
we deduce that $\Ann_{\cO_S}(\eps^n)_*=\Ann_{\cO_{S*}}(\eps^n)$.
However, 
$\eps\cdot K(\eps)_*\subset\bigcup_{n\in\N}\Ann_{\cO_S}(\eps^n)_*$,
so the claim follows easily.
\end{pfclaim}
\begin{claim}\label{cl_stuff-in-nil}
$(K(\eps)\cap\eps^nI)_*\subset\nil(\cO_{S*})$ for every $n>0$.
\end{claim}
\begin{pfclaim} We have
$(K(\eps)\cap\eps^nI)_*=K(\eps)_*\cap(\eps^nI)_*$. Let $x\in(\eps^nI)_*$;
according to lemma \ref{lem_betise.deux}, we have $x^k\in\eps^n\cO_{S*}$
for every sufficiently large $k\in\N$. On the other hand, it is easy
to check that 
$\Ann_{\cO_{S*}}(\eps^m)\cap\eps^n\cO_{S*}\subset\nil(\cO_{S*})$
for every $n,m>0$. In view of claim \ref{cl_about-K(eps)}, the
assertion follows.
\end{pfclaim}

\begin{claim}\label{cl_alm-nilpotent}
If $n>0$, every almost finitely generated subideal of $K(\eps)\cap\eps^nI$
is nilpotent.
\end{claim}
\begin{pfclaim} Let $\cI$ be such an ideal; we can find a finitely
generated ideal $\cI_0\subset(K(\eps)\cap\eps^nI)_*$ such that
$\eps\cI_*\subset\cI_0\subset\cI_*$. Since
$\cI\subset\eps^nI$, lemma \ref{lem_betise.deux} says
that there exists $N\in\N$ such that $(\cI_*)^N\subset\eps^n\cO_{S*}$,
hence $(\cI_*)^{N+1}\subset\cI_0$. From claim \ref{cl_stuff-in-nil}
we deduce that $\cI_0$ is a nilpotent ideal, so the same
holds for $\cI$.
\end{pfclaim}

\begin{claim} In order to prove (ii) we can assume that
$\fm_0=\eps V$, where $\eps\in\fm$ is an $\cO_S$-regular
element.
\end{claim}
\begin{pfclaim} Let us write $X=\Spec\,S_P/J$, where
$S_P:=\Sym^\bullet_{\cO_S}P$ for some almost finitely
projective $\cO_S$-module $P$, and $J\subset S_P$ is an almost
finitely generated ideal. By claim \ref{cl_suffices_m_0-princ} we
can assume that $\fm_0=\eps V$ for some $\eps\in\fm$. Set 
$\bar S:=\Spec\,\cO_S/K(\eps)$ and $\bar S_n:=\bar S\times_SS_n$
for every $n\in\N$. Let $\sigma^\wedge\in X^\wedge_\sm(S)$;
by definition $\sigma^\wedge$ is a compatible family of
sections $\sigma_n:S_n\to X$ lying in the smooth locus of $X$
over $S$. In turns, $\sigma_n$ can be viewed as the datum of
an $\cO_S$-linear morphism $\tau_n:P\to\cO_S/\eps^nI$, such
that the induced morphism of $\cO_S$-algebras
$\Sym^\bullet_{\cO_S}\tau:\Sym^\bullet_{\cO_S}P\to\cO_S/\eps^nI$
satisfies the condition: $\Sym^\bullet_{\cO_S}\tau(J)=0$.
For every $n\in\N$ we have a cartesian diagram of $\cO_S$-algebras:
$$
\xymatrix{ \cO_S/(K(\eps)\cap\eps^nI) \ar[r] \ar[d] &
\cO_{\bar S} \ar[d]^{\pi_n} \\
\cO_{S_n} \ar[r]^{p_n} & \cO_{\bar S_n}.
}$$
Notice that $\eps$ is regular in $\cO_{\bar S}$ and suppose that
assertion (ii) holds for all almost finitely presented $\bar
S$-schemes, especially for $\bar X:=X\times_S\bar S$. 
Let $\bar\sigma_n:\bar S_n\to\bar X$ be the restriction of
$\sigma_n$, for every $n\in\N$. It then follows that, for
every $n\in\N$, $\bar\sigma_n$ extends to a section 
$\bar\sigma\in\bar X_\sm(\bar S)$. The datum of $\bar\sigma$
is equivalent to the datum of an $\cO_{\bar S}$-linear morphism
$\bar\tau:P\to\cO_{\bar S}$ such that
$\Sym^\bullet_{\cO_S}\bar\tau(J)=0$ and such that 
$\pi_n\circ\bar\tau=p_n\circ\tau_n$. The pair $(\tau_n,\bar\tau)$
determines a morphism $\omega_n:P\to\cO_S/(K(\eps)\cap\eps^nI)$,
and by construction we have:
$\Sym^\bullet_{\cO_S}\omega_n(J)=0$, {\em i.e.} $\omega_n$
induces a section $\sigma':\Spec\,\cO_S/(K(\eps)\cap\eps^nI)\to X$,
and then $\sigma'$ must lie in the smooth locus of
$X$ over $S$, in view of lemma \ref{lem_cartes-diagr-sects}.
The obstruction to the existence of a lifting $\omega:P\to\cO_S$
of $\omega_n$, is a class
$\alpha_n\in\Ext^1_{\cO_S}(P,K(\eps)\cap\eps^nI)$.
A simple verification shows that
\set\begin{equation}\label{eq_simple-verif}
K(\eps)\cap\eps^nI=\eps\cdot(K(\eps)\cap\eps^{n-1}I)
\qquad\text{for all $n>0$.}
\end{equation}
From \eqref{eq_simple-verif}, an argument as in the proof
of claim \ref{cl_vanishing.trick} allows to conclude that
the image of $\alpha_n$ in $\Ext^1_{\cO_S}(P,K(\eps)\cap\eps^{n-1}I)$
vanishes; however, this image is none other than $\alpha_{n-1}$,
so actually $\alpha_n=0$, and the sought lifting can be found for
every $n>0$. By construction we have 
$\bar J:=\Sym^\bullet_{\cO_S}\omega(J)\subset K(\eps)\cap\eps^nI$.
Now, $\omega$ induces a section $\sigma'':\Spec\,\cO_S/\bar J\to X$,
which by construction extends $\sigma'$; moreover, the pair
$(\cO_S/\bar J,I/\bar J)$ is henselian (cp. remark
\ref{rem_hensel}(v)), hence $\eps^nI/\bar J$ is a tight radical ideal.
Thus, using proposition \ref{prop_smooth-is-topo}(ii) we derive that
$\sigma''$ lies as well in the smooth locus of $X$ over $S$. Since 
$J$ is almost finitely generated, $\bar J$ is nilpotent in view of
claim \ref{cl_alm-nilpotent}.
Moreover, by \eqref{eq_simple-verif} we can write $\bar J=\eps\cJ$
for some ideal $\cJ\subset K(\eps)\cap\eps^{n-1}I$; $\cJ$ is
nilpotent if $n\geq 2$, since in that case $\cJ^2\subset\bar J$.
Hence we can apply claim \ref{cl_inv-smooths}, to deduce that the
restriction of $\sigma''$ to $\Spec\,\cO_S/\cJ$ extends to a
section $\sigma\in X_\sm(S)$; by construction $\sigma$ extends
$\sigma_{n-1}$. Since $n$ can be taken to be arbitrarily large,
the claim follows.
\end{pfclaim}

So finally we suppose that $\fm_0=tV$, with $t$ an
$\cO_S$-regular element. Set $\cO^\wedge_S:=\liminv{n\in\N}\,\cO_{S_n}$
and $S^\wedge:=\Spec\,\cO^\wedge_S$. We have a natural bijection 
\set\begin{equation}\label{eq_identifies}
X(S^\wedge)\simeq X^\wedge(S):=\liminv{n\in\N}\,X(S_n).
\end{equation}
Let $I^\wedge_*\subset\cO^\wedge_{S*}\simeq\liminv{n\in\N}\,\cO_{S_n*}$
be the topological closure of $I_*$; since $t^nI^\wedge_*$ is the
topological closure of $t^nI_*$, one verifies easily that
\eqref{eq_identifies} identifies the pro-discrete topology of
$X^\wedge(S)$ with the $(t,I^\wedge_*)$-adic topology of
$X(S^\wedge)$ (see \eqref{subsec_add-I} and the proof of lemma 
\ref{lem_add-I}). 

\begin{claim}\label{cl_restricts-ide}
The homeomorphism \eqref{eq_identifies} induces a
bijection: $X_\sm(S^\wedge)\simeq X_\sm^\wedge(S)$.
\end{claim}
\begin{pfclaim} Let $(\sigma_n~|~n\in\N)$ be an element of
$X_\sm(S^\wedge)$, and $\sigma^\wedge\in X(S^\wedge)$ the
corresponding section. We have to show that $\sigma^\wedge$
lies in the smooth locus of $X$. Thus, let 
$I^\wedge:=\Ker(\cO^\wedge_S\to\cO_S/I)$; we have a natural
isomorphism: $\cO^\wedge_S/tI^\wedge\simeq\cO_S/tI$, and by
assumption, the restriction of $\sigma^\wedge$ to
$\Spec\,\cO^\wedge_S/tI^\wedge$ lies in the smooth locus of
$X$. The assertion then follows from proposition
\ref{prop_smooth-is-topo}.
\end{pfclaim}

Under the standing assumptions, we have a cartesian diagram
of almost algebras:
$$
\xymatrix{ \cO_S \ar[r] \ar[d] & \cO_S[t^{-1}] \ar[d] \\
           \cO^\wedge_S \ar[r] & \cO^\wedge_S[t^{-1}]
}$$
whence a cartesian diagram of sets:
$$
\xymatrix{ X(S) \ar[r] \ar[d] & X_t(S_t) \ar[d]^\alpha \\
           X(S^\wedge) \ar[r]^\beta & X_t(S^\wedge_t).
}$$
Now, let $(\sigma_n~|~n\in\N)\in X^\wedge_\sm(S)$; by claim
\ref{cl_restricts-ide}, the corresponding section 
$\sigma^\wedge:S^\wedge\to X$ lies in the smooth locus of $X$.
Let $\sigma^\wedge_t\in X_{\sm,t}(S^\wedge_t)$ be the restriction of
$\sigma^\wedge$.
By corollary \ref{cor_hens.pairs}(ii), the pair $(\cO_{S*},tI_*)$
is henselian; then by \eqref{subsec_invert-t} and proposition
\ref{prop_denseimage}, the restriction 
$\alpha_\sm:X_{\sm,t}(S_t)\to X_{\sm,t}(S^\wedge_t)$
has dense image for the $(t,I^\wedge_*)$-adic topology.
Hence, we can approximate $\sigma^\wedge_t$ arbitrarily
$(t,I^\wedge_*)$-adically close by a section of the form
$\tau^\wedge_t:=\alpha(\tau_t)$, where $\tau_t\in X_{\sm,t}(S_t)$.
Furthermore, $\beta$ is an open imbedding, by lemma \ref{lem_add-I}.
Hence, if $\tau^\wedge_t$ is close enough to $\sigma^\wedge$,
we can find $\tau^\wedge\in X(S^\wedge)$ such that 
$\beta(\tau^\wedge)=\tau^\wedge_t$. In view of proposition
\ref{prop_smooth-is-topo} we can also achieve that $\tau^\wedge$
lies in the smooth locus of $X$. The pair $(\tau_t,\tau^\wedge)$
determines a unique section $\tau\in X(S)$, and by construction
$\tau$ can be obtained as $(t,I)$-adically close to
$(\sigma_n~|~n\in\N)$ as desired. Especially, we can achieve that
$\tau$ lies in the smooth locus of $X$, which concludes the proof
of the theorem.
\end{proof}

\subsection{Quasi-projective almost schemes}\label{sec_qu-proj-schs}
\index{$\sj$ : morphism of fpqc-sites|indref{sec_qu-proj-schs}}
Let $R$ be a $V$-algebra; we denote by $(R\Alg)^o_\fpqc$
the large fpqc site of {\em affine\/} $R$-schemes, and similarly
for the site $(R^a\Alg)^o_\fpqc$. The localization functor
$R\Alg\to R^a\Alg$ defines a morphism of sites:
$$
\sj:(R^a\Alg)^o_\fpqc\to(R\Alg)^o_\fpqc.
$$
If $\cF$ is any sheaf on $(R\Alg)^o_\fpqc$, then $\cF^a:=\sj^*\cF$
can be described as the sheaf associated to the presheaf
$B\mapsto\cF(B_*)$. Especially, we can regard any $R$-scheme
$X$ as a sheaf on $(R\Alg)^o_\fpqc$, and hence we obtain the
{\em almost scheme\/} $X^a$ associated to $X$. If $X$ is affine,
this notation agrees with that of \eqref{subsec_general-schemes}.
For a general $X$, pick a Zariski hypercovering $Z_\bullet\to X$,
where each $Z_i$ is a disjoint union of affine $R$-schemes; then 
$X\simeq\colim{\Delta^o}Z_\bullet$ as fpqc-sheaves, and thus
$X^a\simeq\colim{\Delta^o}Z^a_\bullet$ in the topos of sheaves
on $(R^a\Alg)^o_\fpqc$. Furthermore, we have
$$
\sj_*\Spec\,B=\Spec\,B_{!!}\otimes_{R^a_{!!}}R\qquad
\text{for every $R^a$-algebra $B$}.
$$
If $G$ is an $R$-group scheme, then $G^a$ is clearly an
$R^a$-group scheme, and if $Y\to X$ is any $G$-torsor over an
$R$-scheme $X$ (for the fpqc topology of $X$), then $Y^a\to X^a$
is a $G^a$-torsor.

\sset\subsubsection{}\label{subsec_qu-affine-a-sch}
\index{Almost scheme(s)!quasi-affine|indref{subsec_qu-affine-a-sch}}
Let now $A$ be an $R^a$-algebra, $J\subset A_*$ a finitely
generated ideal. We define the {\em quasi-affine $R^a$-scheme\/}
$\Spec\,A\setminus V(J)$ as the almost scheme 
$\sj^*(\Spec\,A_*\setminus V(J))$. Let $f_1,...,f_k$ be a finite set
of generators of $J$; by \eqref{sec_qu-proj-schs} this can
be realized as the subsheaf
$\bigcup_{i=1}^k\Spec\,A[f^{-1}_i]\subset\Spec\, A$, where the
union takes place in the category of sheaves on $(R^a\Alg)^o_\fpqc$.

\begin{lemma}\label{lem_crit-f.flatnes}
Let $B$ be a $V^a$-algebra and $g_1,...,g_k\in B_*$.
Then $\sum_{i=1}^kg_iB=B$ if and only if the natural morphism
$B\to\prod_{i=1}^kB[g^{-1}_i]$ is faithfully flat.
\end{lemma}
\begin{proof} Supppose that $\sum_{i=1}^kg_iB=B$; we have to show
that $\prod_{i=1}^kC[g^{-1}_i]\neq 0$ for every non-zero quotient
$C$ of $B$. Replacing $B$ by $C$, we are reduced to showing that,
if $B[g^{-1}_i]=0$ for every $i\leq k$, then $B=0$. However, the
condition $B[g^{-1}_i]=0$ implies that for every $\eps\in\fm$ there
exists $n_i\in\N$ such that $\eps\cdot g^{n_i}_i=0$; set
$n:=\max(n_i~|~i\leq k)$. We have $\sum_{i=1}^kg^n_iB=B$, hence
$\eps\cdot B=0$. Since $\eps$ is arbitrary, the claim follows.
Conversely, let $B':=\prod_{i=1}^kB[g^{-1}_i]$ and 
$I:=\sum_{i=1}^kg_iB$; clearly $IB'=B'$, therefore $I=B$ provided
$B'$ is faithfully flat over $B$.
\end{proof}

\begin{lemma}\label{lem_describe-q-affine}
With the notation of \eqref{subsec_qu-affine-a-sch},
the almost scheme $\Spec\,A\setminus V(J)$ is the subfunctor
of $\Spec\,A$, whose $T$-section are the morphisms
$\sigma:T\to\Spec\,A$ such that $\sigma^\sharp(J)\cO_T=\cO_T$,
for every affine $R^a$-scheme $T$.
\end{lemma}
\begin{proof} Let $\sigma$ be a $T$-section of $\Spec\,A\setminus V(J)$;
by definition, this is the same as saying that the natural morphism
$\amalg_{i=1}^kT\times_{\Spec\,A}\Spec\,A[f^{-1}_i]\to T$
is a fpqc-covering of $T$. Hence the claim is just a restatement
of lemma \ref{lem_crit-f.flatnes}.
\end{proof}

\sset\subsubsection{}
Alternatively, one can regard a quasi-affine $R^a$-scheme as
a difference of two affine $R^a$-schemes. This viewpoint is
elaborated in the following:

\begin{definition}\label{def_some.sheaf.defs}
\index{Closed imbedding, difference of sheaves
on $(R^a\Alg)^o_\fpqc$|indref{def_some.sheaf.defs}}
Let $\phi:\cF\to\cG$ be a morphisms of sheaves on $(R^a\Alg)^o_\fpqc$.
\begin{enumerate}
\item
We say that $\phi$ is a {\em closed imbedding\/} if, for every
affine $R^a$-scheme $T$ and every section $T\to\cG$, the
induced morphism $T\times_\cG\cF\to T$ is a closed imbedding of
affine $R^a$-schemes.
\item
The {\em difference\/} $\cG\setminus\cF$ is the subsheaf of $\cG$
whose $T$-sections are all the $T$-sections $T\to\cG$ of $\cG$
such that $T\times_\cG\cF=\emptyset$, for every affine $R^a$-scheme
$T$. Here $\emptyset$ denotes the initial object of the topos of
sheaves on $(R^a\Alg)^o_\fpqc$.
\end{enumerate}
\end{definition}
Using lemma \ref{lem_describe-q-affine} one checks easily
that the sheaf-theoretic difference $\Spec\,A\setminus\Spec\,A/J$
is the same as the almost scheme considered in
\eqref{subsec_qu-affine-a-sch}.

\begin{lemma}\label{lem_pull-differ}
Let $\phi:\cF\to\cG$ be a morphism of sheaves on
$(R\Alg)^o_\fpqc$. Then:
$$
\sj^*(\cG\setminus\cF)\subset\sj^*\cG\setminus\sj^*\cF.
$$
\end{lemma}
\begin{proof} $\cG\setminus\cF$ is also the largest subobject
of of $\cG$ that has empty intersection with $\Img(\phi)$; in
particular this is defined for every topos. Then the statement
holds for any morphism of topoi (in this case $\sj$), and it
follows easily from the facts that the pull-back functor
(here $\sj^*$) is left exact and sends empty objects to
empty objects.
\end{proof}

\begin{example}\label{ex_proj-Grass}
\index{Almost scheme(s)!$\P_{R^a}(E)$, $\mathrm{Grass}^r_{R^a}(E)$ :
projective space, Grassmannian|indref{ex_proj-Grass}}
Let $E$ be an $R^a$-module.
\begin{enumerate}
\item
We define the presheaf $\P_{R^a}(E)$ on $(R^a\Alg)^o_\fpqc$ as follows.
For every affine $R^a$-scheme $T$, the $T$-sections of $\P(E)$
are the strictly invertible quotients of
$E\otimes_{R^a}\cO_T$, {\em i.e.\/} the equivalence classes
of epimorphisms $E\otimes_{R^a}\cO_T\to L$ where $L$ is a strictly
invertible $\cO_T$-module (see definition \eqref{def_strictly_invertible}).
\item
More generally, for every integer $r\in\N$ we can define
the presheaf $\mathrm{Grass}^r_{R^a}(E)$; its $T$-sections are
the rank $r$ almost $\cO_T$-projective quotients
$E\otimes_{R^a}\cO_T\to P$.
\end{enumerate}
By faithfully flat descent it is easy to verify that both
these presheaves are in fact fpqc-sheaves.
\end{example}

\sset\subsubsection{}
Let $E$ be an $R$-module, $B$ an $R^a$-algebra; a section
$\sigma\in\mathrm{Grass}^r_R(E)(\Spec\,B_*)$ is a quotient map
$\sigma:E\otimes_RB_*\to P$, onto a projective $B_*$-module $P$
of rank $r$. To $\sigma$ we associate the quotient
$\sigma^a:E^a\otimes_RB\to P^a$; after passing to the associated
sheaves we obtain a morphism:
\set\begin{equation}\label{eq_pull-back-Grass}
\sj^*\mathrm{Grass}^r_R(E)\to\mathrm{Grass}^r_{R^a}(E^a)\qquad
\sigma\mapsto\sigma^a.
\end{equation}
Since $P\simeq P^a_*$, it is clear that the rule
$\sigma\mapsto\sigma^a$ defines an injective map of presheaves,
hence \eqref{eq_pull-back-Grass} is a monomorphism.

\begin{lemma}\label{lem_its-a-scheme}
Suppose that $E$ is a finitely generated $R$-module.
Then \eqref{eq_pull-back-Grass} is an isomorphism.
\end{lemma}
\begin{proof} We have to show that \eqref{eq_pull-back-Grass}
is an epimorphism of fpqc-sheaves; hence, let $T$ be any affine
$R^a$-scheme, $\sigma:E^a\otimes_{R^a}\cO_T\to P$ a quotient
morphism, with $P$ almost projective over $\cO_T$ of rank $r$.
The contention is that there exists a faithfully flat morphism
$f:U\to T$ such that $f^*\sigma$ is in the image of
$\mathrm{Grass}^r_R(E)(U_*)$ (where $U_*:=\Spec\,\cO_{U*}$).
By theorem \ref{th_free.in.fpqc} we can then reduce to the case
when $P=\cO_T^r$. We deduce an epimorphism 
$\Lambda^r_{\cO_T}\sigma:
\Lambda^r_{\cO_T}(E^a\otimes_{R^a}\cO_T)\to
\Lambda^r_{\cO_T}(\cO_T^r)\simeq\cO_T$. Let $e_1,...,e_k$
be a finite set of generators of $E$, and for every subset 
$J:=\{j_1,...,j_r\}\subset\{1,...,k\}$ of cardinality $r$, let 
$e_J:=e_{j_1}\wedge...\wedge e_{j_r}\in
\Lambda^r_{\cO_T}(E^a\otimes_{R^a}\cO_T)$; it follows that
the set $(t_J:=\Lambda^r_{\cO_T}\sigma(e_J)~|~J\subset\{1,...,k\})$
generates $\cO_T$. Set $U_J:=\Spec\,\cO_T[t_J^{-1}]$ for every 
$J\subset\{1,...,k\}$ as above; by construction $\sigma$ induces
surjective maps $E\otimes_R\cO_{T*}[t_J^{-1}]\to\cO_{T*}^r[t_J^{-1}]$;
after tensoring by $\cO_T[t_J^{-1}]_*$ we obtain surjections
$E\otimes_R\cO_T[t_J^{-1}]_*\to\cO_T[t_J^{-1}]_*^r$. In other words,
the restriction $\sigma_{|U_J}$ lies in the image of
$\mathrm{Grass}^r_R(E)(U_{J*})$ for every such $J$; on the other
hand, lemma \ref{lem_crit-f.flatnes} says that the natural morphism
$\amalg_JU_J\to T$ is a fpqc covering, whence the claim.
\end{proof}

Lemma \ref{lem_its-a-scheme} says in particular that
$\mathrm{Grass}^r_{R^a}(E^a)$ is an almost scheme associated to
a scheme whenever $E$ is finitely generated. Next we wish to define
the projective $R^a$-scheme associated a graded $R^a$-algebra.

\begin{definition}\label{def_proj-a}
\index{$\Proj(A)$ : projective almost scheme associated to a graded
almost algebra|indref{def_proj-a}}
Let $A:=\oplus_{i\in\N}A_i$ be a graded $R^a$-algebra which is
generated by $A_1$ over $A_0:=R^a$, in the sense that the natural
morphism of graded $R^a$-algebras $\Sym^\bullet_{R^a}A_1\to A$
is an epimorphism. We let $\Proj(A)$ be the sheaf on
$(R^a\Alg)^o_\fpqc$ defined as follows. For every affine
$R^a$-scheme $T$, the $T$-points of $\Proj(A)$ are the
strictly invertible quotients $A_1\otimes_{R^a}\cO_T\to L$ of
$A_1\otimes_{R^a}\cO_T$, inducing morphisms of graded $R^a$-algebras
$A\to\Sym^\bullet_{\cO_T}L$.
\end{definition}

\begin{lemma}\label{lem_closed-imb} 
With the notation of \eqref{def_proj-a}, the natural map of sheaves
$\Proj(A)\to\P_{R^a}(A_1)$ is a closed imbedding.
\end{lemma}
\begin{proof} Let $T$ be any affine $R^a$-scheme and 
$\sigma:T\to\P_{R^a}(A_1)$ a morphism of fpqc-sheaves. We have to
show that $T':=T\times_{\P_{R^a}(A_1)}\Proj(A)$ is a closed
subscheme of $T$. By definition, $\sigma:A_1\otimes_{R^a}\cO_T\to L$
is a strictly invertible quotient; we deduce via
$\Sym^\bullet_{\cO_T}\sigma$ a natural structure of
$\Sym^\bullet_{R^a}A_1$-module on $\Sym^\bullet_{\cO_T}L$.
Let $U$ be an affine $R^a$-scheme and $U\to T'$ a morphism of
fpqc-sheaves; set $J_n:=\Ker\,(\Sym^n_{R^a}A_1\to A)$ for every
$n\in\N$. By definition:
$$
J_n\cdot\Sym^n_{\cO_T}L\subset
\Ker\,(\Sym^n_{\cO_T}L\to\cO_U\otimes_{\cO_T}\Sym^n_{\cO_T}L).
\qquad\text{for every $n\in\N$}
$$
Since $L$ is flat, this is the same as: 
\set\begin{equation}\label{eq_some-Syms}
J_n\cdot\Sym^n_{\cO_T}L\subset
\Ker\,(\cO_T\to\cO_U)\cdot\Sym^n_{\cO_T}L.
\qquad\text{for every $n\in\N$}
\end{equation}
For every $n\in\N$ denote by
$\ev_n:\Sym^n_{\cO_T}L\otimes_{\cO_T}\Sym^n_{\cO_T}L^*\to\cO_T$
the evaluation morphism of the $\cO_T$-module
$\Sym^n_{\cO_T}L$. Then \eqref{eq_some-Syms} is equivalent to:
$$
\sum_{n\in\N}\ev_n
(J_n\cdot\Sym^n_{\cO_T}L\otimes_{\cO_T}\Sym^n_{\cO_T}L^*)
\subset\Ker\,(\cO_T\to\cO_U).
$$
Conversely, let 
\set\begin{equation}\label{eq_def-kernelI}
I:=\sum_{n\in\N}\ev_n
(J_n\cdot\Sym^n_{\cO_T}L\otimes_{\cO_T}\Sym^n_{\cO_T}L^*).
\end{equation}
Then the restriction of $\sigma$ to $\Spec\,\cO_T/I$ is a section
of $\Proj(A)$. All in all, this shows that $T'$ represents the closed
subscheme $\Spec\,\cO_T/I$ of $T$, whence the claim.
\end{proof}

\begin{lemma}\label{lem_a-proj-OK}
Let $S:=\oplus_{i\in\N}S_i$ be a graded $R$-algebra,
with $S_0=R$; suppose that $S$ is generated by $S_1$ and that
$S_1$ is a finitely generated $R$-module. Then the natural
map
$$
\sj^*\Proj(S)\to\Proj(S^a)
$$
is an isomorphism.
\end{lemma}
\begin{proof} Let $B$ be an $R^a$-algebra, and
$\sigma:S_1^a\otimes_{R^a}B\to L$ a strictly invertible $B$-module
quotient. By lemma \ref{lem_its-a-scheme} we know that there is
a faithfully flat $B$-algebra $C$ such that the induced map
$\sigma':S_1\otimes_RC_*\to(L\otimes_BC)_*$ is a rank one projective
quotient. Suppose now that $\sigma$ is a section of
$\Proj(S^a)$ and let $J:=\Ker(\Sym_R^\bullet(S_1)\to S)$;
by definition we have $J^a\subset\Ker\,\Sym_B^\bullet\sigma$,
and since $\Sym^\bullet_{C_*}(L\otimes_BC)_*$ does not contain
$\fm$-torsion, it follows that $J\subset\Ker\,\Sym_{C_*}^\bullet\sigma'$.
\end{proof}

\sset\subsubsection{}\label{subsec_qu-proj-a-sch}
Next we consider quasi-projective almost schemes. Hence, let
$A$ be as in definition \ref{def_proj-a} and 
$\cJ:=\oplus_{i\in\N}\cJ_i\subset A$ a graded ideal. We set
$X:=\Proj(A)\setminus\Proj(A/\cJ)$.

\begin{lemma}\label{lem_Idontknow}
Keep the assumptions of \eqref{subsec_qu-proj-a-sch},
and suppose moreover that $A=S^a$ and $\cJ=J^ a$ for a graded
$R$-algebra $S$ and a finitely generated graded ideal $J\subset S$.
Suppose moreover that $S_1$ is a finitely generated $R$-module and
that $S_1$ generates $S$. Then
$$
X\simeq\sj^*(\Proj(S)\setminus\Proj(S/J)).
$$
\end{lemma}
\begin{proof} Set $Y:=\Proj(S)$, $Y_0:=\Proj(S/J)$,
$X':=Y\setminus Y_0$; in view of lemmata \ref{lem_pull-differ}
and \ref{lem_a-proj-OK}, we have only to show that
$$
\sj^*Y\setminus\sj^*Y_0\subset\sj^*X'.
$$
Hence, let $B$ be an $R^a$-algebra and $\sigma\in\sj^*Y(\Spec\,B)$;
by lemma \ref{lem_a-proj-OK} we can find a faithfully flat morphism
$T\to\Spec\,B$ and a projective $\cO_{T*}$-module $L_*$ of rank one
such that the restriction of $\sigma$ to $T$ is induced by a surjective
map $\tau:S_1\otimes_R{\cO_{T*}}\to L_*$. Suppose now that 
$\sigma$ is a section of $\sj^*Y\setminus\sj^*Y_0$; it follows
that $\tau_0:=\tau\times_YY_0=\emptyset$. From lemma
\ref{lem_closed-imb} and its proof we deduce easily that $\tau_0$
is represented the $R^a$-scheme $\Spec\,\cO_T/I$, where $I$ is
defined as in \eqref{eq_def-kernelI}; consequently $I=\cO_T$.
Since $\Sym^\bullet_{\cO_T}L^*$ is generated by $(L_*)^*$,
which is a finitely generated $\cO_{T*}$-module, and since
$J$ is finitely generated by assumption, we see that $I$ is
generated by the images $x_1,...,x_k$ of finitely many elements
in $J\cdot\Sym^\bullet_{\cO_{T*}}L_*\otimes_{\cO_T}
\Sym^\bullet_{\cO_{T*}}(L_*)^*$.
It then follows from lemma \ref{lem_crit-f.flatnes} that there exists
a covering morphism $U\to T$ such that
$\sum_{i=1}^kx_i\cO_{U*}=\cO_{U*}$, which means that the restriction of
$\sigma$ to $U$ is a section of $\sj^*X'$.
\end{proof}

\sset\subsubsection{}\label{subsec_point-of-view}
In the situation of \eqref{subsec_qu-proj-a-sch}, let 
$Y:=\Spec\,A\setminus\Spec\,(A/\cJ A_+)$, where $A_+:=\oplus_{i>0}A_i$;
then $Y$ represents the functor that assigns to every $R^a$-scheme
$T$ the set of all maps of $R^a$-algebras $\phi:A\to\cO_T$ such that
\set\begin{equation}\label{eq_conditio}
\phi(\cJ A_+)\cdot\cO_T=\cO_T.
\end{equation}
Moreover, we have a map of $R^a$-schemes $\pi:Y\to X$, representing the
natural transformation of functors that assigns to any $T$-section 
$\phi:A\to\cO_T$ of $Y$ the induced map of $\cO_T$-modules
$\phi_1:A_1\otimes_R\cO_T\to\cO_T$. Condition \eqref{eq_conditio}
ensures that $\phi_1$ is surjective and $\phi_1\notin\Proj(A/\cJ)(T)$
whenever $T$ is non-empty, hence this rule yields a well defined
$T$-section of $X$. It is a standard fact that $\pi$ is a $\G_m$-torsor.
The action of $\G_m(T)$ on $Y(T)$ can be described explicitly as follows.
To a given pair $(\phi,u)\in Y(T)\times\cO_T^\times$ one assigns the
unique $R$-algebra map $\phi_u:A\to\cO_T$ such that $\phi_u(x)=u\phi(x)$
for every $x\in A_1$.

\begin{lemma} Keep the notation of \eqref{subsec_point-of-view}.
The subsheaf
$$
Y_\sm:=(\Spec\,A)_\sm\setminus\Spec\,(A/\cJ A_+)\subset Y
$$
is a subtorsor of $Y$ for the natural $\G_m$-action on $Y$.
\end{lemma}
\begin{proof} Clearly $\G_m=(\G_m)_\sm$; hence, by lemma
\ref{lem_smooth-prod}, we have a natural identification:
$$
(\G_m\times_{R^a}Y)_\sm\simeq\G_m\times Y_\sm.
$$
However, the natural transformation of functors 
$\G_m\times_{R^a}Y\to Y\times_XY$ defined by the rule:
$$
(u,\phi)\mapsto(\phi,\phi_u)\quad
\text{for every local section $u$ of $\G_m$ and $\phi$ of $Y$}
$$
is an isomorphism, therefore
$(\G_m\times_{R^a}Y)_\sm\simeq(Y\times_{X^a}Y)_\sm$. Let $T$
be an $R^a$-scheme, $\phi\in Y(T)$ and $u\in\G_m(T)$; it follows
that $(\phi,\phi_u)\in(Y\times_XY)_\sm(T)$ if and only if 
$\phi\in Y_\sm(T)$. Furthermore, it is clear that 
$(\phi,\phi_u)$ is in the smooth locus of $Y\times_XY$ if
and only if $(\phi_u,\phi)$ is (use the $X$-automorphism
of $Y\times_XY$ that swaps the two factors).
However, the pair $(\phi_u,\phi)$ can be written in the form
$(\psi,\psi_v)$ for $\psi=\phi_u$ and $v=u^{-1}$, so that
$\psi$ is in the smooth locus if and only if $(\psi,\psi_v)$ is.
{\em Ergo}, $\phi$ is in the smooth locus if and only if $\phi_u$
is, as required.
\end{proof}

\begin{definition}\label{def_smooth.for.proj}
Let $X$ and $Y$  be as in \eqref{subsec_point-of-view}.
We define the {\em smooth locus of $X$\/} as the quotient 
(in the category of sheaves on $(R^a\Alg)^o_\fpqc$):
$$
X_\sm:=Y_\sm/\G_m.
$$
\end{definition}

The notation $X_\sm$ remains somewhat ambiguous, until
one shows that the smooth locus of $X$ does not depend on
the choice of presentation of $X$ as quotient of a quasi-affine
almost scheme $Y$. While we do not completely elucidate this
issue, we want at least to make the following remark:

\begin{lemma}\label{lem_twist-by-L}
In the situation of  \eqref{subsec_point-of-view},
suppose that $\cL$ is a strictly invertible $R^a$-module, and
set $B:=\oplus_{r\in\N}A_r\otimes_{R^a}\cL^{\otimes r}$,
$\cJ_B:=\oplus_{r\in\N}\cJ_r\otimes_{R^a}\cL^{\otimes r}$,
so $B$ is a graded $R^a$-algebra and $\cJ_B\subset B$ a graded
ideal. Then there are natural isomorphisms
$$
X\simeq X_B:=\Proj(B)\setminus\Proj(B/\cJ_B)\qquad
X_\sm\simeq X_{B,\sm}.
$$
\end{lemma}
\begin{proof} We define as follows a natural isomorphism of
functors $X\stackrel{\sim}{\to}X_B$. To every $R^a$-scheme $T$
and every $T$-section $A_1\otimes_{R^a}\cO_T\to L$, we assign the
strictly invertible quotient
$A_1\otimes_{R^a}\cL\otimes_{R^a}\cO_T\to\cL\otimes_{R^a}L$.
Set $Y_B:=\Spec\,B\setminus\Spec\,(B/\cJ_B B_+)$; it follows
that $Y_B/\G_m\simeq X_B$.
Hence, the identity $X_\sm\simeq X_{B,\sm}$ can be checked
locally on the fpqc topology. However, $\cL$ is locally free
of rank one on $(R^a\Alg)^o_\fpqc$, and the claim is obvious
in case $\cL$ is free.
\end{proof}

\begin{theorem}\label{th_raise-smooth-sects}
Suppose $\tilde\fm$ has homological dimension $\leq 1$, and that
$(R^a,I)$ is a henselian pair. Let $P$ be an almost finitely
generated projective $R^a$-module and 
$\cI,\cJ\subset S_P:=\Sym_{R^a}^\bullet P$ two graded ideals,
with $\cJ$ almost finitely generated. Set
$X:=\Proj(S_P/\cJ)\setminus V(\cI)$. Then, for every finitely
generated subideal $\fm_0\subset\fm$, we have:
$$
\Img(X_\sm(\Spec\,R^a/\fm_0I)\to X_\sm(\Spec\,R^a/I))=
\Img(X_\sm(\Spec\,R^a)\to X_\sm(\Spec\,R^a/I)).
$$
\end{theorem}
\begin{proof} First of all, arguing as in the proof of claim
\ref{cl_reduce-to-tight} we reduce to the case where $I$ is a
tight ideal. Next, set $\cJ:=\oplus_{i\in\N}\cJ_i$,
$\cI:=\oplus_{i\in\N}\cI_i$ and let
$\bar\sigma:(P/\cJ_1)\otimes_{R^a}R^a/\fm_0I\to\bar L$
be an $R^a/\fm_0I$-section of $X$. By theorem
\ref{th_lifts.proj.hensel} we can lift $\bar L$ to an almost
finitely generated projective $R^a$-module $L$. We have 
$\Lambda_{R^a}^2L\otimes_{R^a}R^a/\fm_0I\simeq\Lambda_{R^a}^2\bar L=0$,
whence $\Lambda_{R^a}^2L=0$ by Nakayama's lemma \ref{lem_Naka},
so $L$ is strictly invertible. We set
$S':=\Sym^\bullet_{R^a}(P\otimes_{R^a}L^*)$, 
$\cI':=\oplus_{r\in\N}\cI_r\otimes_{R^a}L^{*\otimes r}$,
$\cJ':=\oplus_{r\in\N}\cJ_r\otimes_{R^a}L^{*\otimes r}$, $B:=S'/\cJ'$,
$Y':=\Spec\,B\setminus\Spec\,(B/\cI'B_+)$ and
$X':=\Proj(B)\setminus\Proj(B/\cI')$. By lemma \ref{lem_twist-by-L}
we have an isomorphism $\beta:X'\simeq X$ preserving the smooth loci;
by inspecting the definition we find that $\beta^{-1}(\bar\sigma)$
lies in the image of the projection
$Y'_\sm(\Spec\,R^a/\fm_0I)\to X'_\sm(\Spec\,R^a/\fm_0I)$.
Hence we can replace $P$ by $P\otimes_{R^a}L^*$ and assume
from start that $\bar\sigma$ lies in the image of the map
$Y_\sm(\Spec\,R^a/\fm_0I)\to X_\sm(\Spec\,R^a/\fm_0I)$.
Let $\bar\tau\in Y_\sm(\Spec\,R^a/\fm_0I)$ that maps to
$\bar\sigma$. By theorem \ref{th_lift-smooth-sects}(i), the image
$\tau_0$ of $\bar\tau$ in $(\Spec\,S_P/\cJ)_\sm$ lifts to a section
$\tau\in(\Spec\,S_P/\cJ)_\sm(\Spec\,R^a)$. To conclude, it suffices
to show the following:
\begin{claim} Suppose that $I$ is tight. Then $\tau\in Y(\Spec\,R^a)$.
\end{claim}
\begin{pfclaim} Let $A:=S_P/\cJ$, so that $\tau:A\to R^a$ is
a morphism of $R^a$-algebras, and set $J:=\tau(A_+\cI)$; by
assumption $\tau_0\in Y(\Spec\,R^a/I)$, which means that
$J\cdot(R^a/I)=(R^a/I)$, {\em i.e.} $J+I=R^a$.
Set $M:=R^a/JR^a$; it then follows that $M=IM$, hence $M=0$
by lemma \ref{lem_Naka}; hence $JR^a=R^a$, which
is the claim.
\end{pfclaim}
\end{proof}

\subsection{Lifting and descent of torsors}\label{sec_lift-tors}
We begin this section with some generalities about group schemes
and their linear representations, that are preliminary to our
later results about liftings of torsors in an almost setting
(theorems \ref{th_lift-iso-tors} and \ref{th_last-one}).

Furthermore, we apply the results of section \ref{sec_elkik} to
derive a descent theorem for $G$-torsors, where $G$ is a group
scheme satisfying some fairly general conditions.

\sset\subsubsection{}
Let $R$ be a ring, $G$ an affine group scheme defined over $R$.
We set $\cO(G):=\Gamma(G,\cO_G)$; so $\cO(G)$ is an $R$-algebra
and the multiplication map of $G$ determines a structure of
co-algebra on $\cO(G)$. In \eqref{subsec_def-G-mods} we have
defined the notion of a (left or right) $G$-action on a quasi-coherent
$\cO_X$-module $M$, where $X$ is a scheme acted on by $G$. We
specialize now to the case where $X=\Spec\,R$, and $G$ acts trivially
on $X$. In such situation, a $G$-action on $M$ is the same as a map
of fpqc-sheaves, from the sheaf represented by $G$ to the sheaf of
automorphisms of $M$; {\em i.e.} the datum, for every $R$-scheme $T$,
of a functorial group homomorphism:
\set\begin{equation}\label{eq_naturally}
G(T)\to\Aut_{\cO_T}(M\otimes_R\cO_T).
\end{equation}
Explicitly, given a map as in \eqref{eq_naturally}, take
$T:=G$; then the identity map $G\to G$ determines an
$\cO(G)$-linear automorphism $\beta$ of $M\otimes_R\cO(G)$,
and one verifies easily that $\beta$ fulfills the conditions
of \eqref{subsec_def-G-mods}. Conversely, an automorphism
$\beta$ as in \eqref{subsec_def-G-mods} extends uniquely to
a well defined map of sheaves \eqref{eq_naturally}.
Furthermore, the $G$-action on $M$ can be prescribed by
choosing, instead of a $\beta$ as above, an $R$-linear
map:
$$
\gamma:M\to M\otimes_R\cO(G)
$$
satisfying certain identities analogous to \eqref{eq_cocycle-action}
and \eqref{eq_rigidify}; of course one can then recover the
corresponding $\beta$ by extension of scalars. One says that the pair
$(M,\gamma)$ is an {\em $\cO(G)$-comodule\/}; more precisely, one
defines right and left $\cO(G)$-comodules, and the bijection just
sketched sets up an equivalence from the category of $R$-modules
with left $G$-actions to the category of right $\cO(G)$-comodules.
One says that an $\cO(G)$-comodule is {\em finitely generated}
(resp. {\em projective}) if its underlying $R$-module has the same
property.

\begin{lemma}\label{lem_Serre}
Suppose that $R$ is a Dedekind domain and that $\cO(G)$ is flat
over $R$. Then:
\begin{enumerate}
\item
Every $\cO(G)$-comodule is the filtered union of its finitely
generated $\cO(G)$-comodules.
\item
Every finitely generated $\cO(G)$-comodule is quotient of a
projective $\cO(G)$-comodule.
\end{enumerate}
\end{lemma}
\begin{proof} (i) is \cite[\S1.5, Cor.]{Ser} and (ii) is
\cite[\S2.2, Prop.3]{Ser}.
\end{proof}

\sset\subsubsection{}\label{subsec_G-torsors}
Let $R$ be a ring, $G$ be a flat group scheme of finite presentation 
over $R$, and suppose that $G$ admits a closed imbedding
as a subgroup scheme $G\subset\GL_n$. We suppose
moreover that the quotient $X_G:=\GL_n/G$ (for the
right action of $G$ on $\GL_n$) is representable
by a quasi-projective $R$-scheme, and that $X_G$ admits 
an ample $\GL_n$-equivariant line bundle. Then $X_G$ is
necessarily of finite presentation by \cite[lemme 17.7.5]{EGA4}.
Furthermore, $X_G$ is smooth; indeed $X_G$ is flat over
$\Spec\,R$, hence smoothness can be checked on the
geometric fibres over the points of $\Spec\,R$, which
reduces to the case where $R$ is an algebraically closed
field. In such case, smoothness over $R$ is the
same as regularity and the latter follows since $\GL_n$
is regular and the quotient map $\GL_n\to X_G$ is faithfully
flat. The following lemma \ref{lem_envisage} provides plenty
of examples of the situation envisaged in this paragraph.

\begin{lemma}\label{lem_envisage}
Suppose that $R$ is a Dedekind domain and $G$ is a flat affine
group scheme of finite type over $R$. Then $G$ fulfills the
conditions of \eqref{subsec_G-torsors}.
\end{lemma}
\begin{proof} The proof is obtained by assembling several
references to the existing literature. To start with, let
$P$ be a finitely generated subrepresentation of the left
(or right) regular representation of $G$ on $\cO(G)$ which
generates $\cO(G)$ as an $R$-algebra. By \cite[1.4.5]{Br-Ti}
$P$ is a finitely generated projective $R$-module and the action
of $G$ on $P$ is faithful, in the sense that it gives a closed
immersion: $G\to\Aut_R(P)$. Taking the direct sum with a trivial
representation one gets a closed immersion of $R$-group schemes:
$G\to\GL_{n,R}$. By a theorem of M.Artin (see \cite[3.1.1]{Ana})
the quotient fppf sheaf $\GL_n/G$ is represented by an algebraic
space of finite presentation over $\Spec\,R$. Then by \cite[Th.4.C]{Ana}
this algebraic space is a scheme. Let $K$ be the fraction field
of $R$; by classical results of Chevalley and Chow, $GL_{n,K}/G_K$
is quasi-projective over $K$ and by \cite[VIII.2]{Ray0} it follows
that $\GL_n/G$ is quasi-projective over $\Spec\,R$.
The ample line bundle $\cL$ on $\GL_n/G$ that one gets always
has a $\GL_n$-linearization. This follows from the proof of
\cite[lemma 1.2]{Su}. For our purpose it suffices to know that
the weaker assertion that some tensor multiple of $\cL$ is
$\GL_n$-linearizable. We may assume that $\cL$ is trivial
on the identity section $S:=\Spec\,R\to\GL_{n,R}/G$.
\begin{claim} 
$\Ker(\mathrm{Pic}\,\GL_{n,R}\stackrel{e^*}{\to}\mathrm{Pic} S)=0$
(where $e:S\to\GL_{n,R}$ is the identity section).
\end{claim}
\begin{pfclaim} $\GL_{n,K}$ is open in an affine space, hence it
has trivial Picard group, hence every divisor on $\GL_{n,R}$ is
linearly equivalent to a divisor whose irreducible components 
do not dominate $S$, {\em i.e.} a pull-back of a divisor on $S$,
which gives the claim.
\end{pfclaim}

The assertion now follows from \cite[VII, Prop.1.5]{Ray0}.
\end{proof}

\begin{lemma}\label{lem_equivariant}
Under the assumptions of \eqref{subsec_G-torsors},
there is a $\GL_n$-equivariant isomorphism of $R$-schemes:
\set\begin{equation}\label{eq_equivariant}
X_G\stackrel{\sim}{\to}\Proj((\Sym^\bullet_R\cP)/\cI)\setminus V(\cJ)
\end{equation}
where $\cP$ is a projective $\GL_n$-comodule, and
$\cI,\cJ\subset\Sym^\bullet_R\cP$ are two finitely generated
graded ideals, that are sub-$\cO(G)$-comodules of $\Sym^\bullet_R\cP$.
\end{lemma}
\begin{proof} Let $\cL$ be an ample $\GL_n$-equivariant
line bundle on $X_G$. Since $X_G$ is quasi-compact, we can find
$n\in\N$ large enough, so that $\cL^{\otimes n}$ is very ample.
We can then replace $\cL$ by $\cL^{\otimes n}$ and therefore
assume that $X_G$ admits a locally closed imbedding in
$\P_R(\Gamma(X_G,\cL))$. Again by compactness argument, we
can find a finitely generated $R$-submodule $W\subset\Gamma(X_G,\cL)$
such that $X_G$ already imbeds into $\P_R(W)$. Let $f_1,...,f_k$ be a
finite set of generators of the $R$-module $W$. By lemma
\ref{lem_Serre}(i) we can find a finitely generated $\Z$-module
$W_0\subset W$ which is a $\cO(\GL_n)$-comodule containing the 
$(f_i~|~i\leq k)$. Up to replacing $W$ by the $R$-module generated
by $W_0$, we can further assume that $W$ is an $\cO(\GL_n)$-comodule.

\begin{claim} $W$ is the quotient of a projective $\cO(G)$-comodule
of finite type $L$.
\end{claim}
\begin{pfclaim} By the foregoing we can assume that $W$ is
generated by a finitely generated $\Z$-submodule $W_0\subset W$
which is an $\cO(G)$-comodule. By lemma \ref{lem_Serre}(ii) we can
write $W_0$ as a quotient of a projective $\cO(G)$-comodule $L_0$;
hence $W$ is a quotient of $L:=L_0\otimes_\Z R$.
\end{pfclaim}

It follows that $X_G$ is a locally closed subscheme of
$\P(L)$. Let $Y$ be the schematic closure of the image of $X_G$;
we can write $Y\simeq\Proj((\Sym_R^\bullet L)/I)$, for some
graded ideal $I$, and $I$ is necessarily an $\cO(G)$-comodule.
Furthermore, we have $Y\setminus X_G=V(J)$, where 
$J\subset\Sym_R^\bullet L$ is another graded ideal, also an
$\cO(G)$-comodule. Arguing as in the foregoing we can write
$J=\bigcup_\alpha J_\alpha$, where $J_\alpha$ runs over the
filtered family of the finitely generated $\cO(G)$-sub-comodules
of $J$. By quasi-compactness one shows easily that
$Y\setminus X_G= V(J_\alpha)$ for $J_\alpha$ large enough;
we set $\cJ:=J_\alpha$. Likewise, we can find a finitely
generated $\cO(G)$-subcomodule $\cI\subset I$ such that
$(\P(L)\setminus V(\cJ))\cap V(I)=(\P(L)\setminus V(\cJ))\cap V(\cI)$.
\end{proof}

\sset\subsubsection{}
Under the assumptions of \eqref{subsec_G-torsors}, let $T$ be
a $\GL_n$-torsor. Then $G$ acts on $T$ via the imbedding
$G\subset\GL_n$, and the functor $T/G$ is representable
by a smooth $R$-scheme.

\begin{lemma}\label{lem_neatsol} 
The functor $T/G:R\text{-}\mathbf{Scheme}\to\Set$
is naturally isomorphic to the functor that assigns to every 
$R$-scheme $S$ the set of all pairs $(H,\omega)$, where $H$
is a $G$-torsor on $S_\fpqc$ and 
$\omega:H\times^G\GL_n\stackrel{\sim}{\to}T$ is an
isomorphism of $\GL_n$-torsors.
\end{lemma}
\begin{proof} Indeed, let $\sigma:S\to T/G$ be any morphism of
$R$-schemes; set $H_\sigma:=S\times_{T/G}T$, which is a 
$G$-torsor on $S_\fpqc$. The natural morphism
$H_\sigma\to T$ induces a well defined isomorphism of 
$\GL_n$-torsors 
$\omega_\sigma:H_\sigma\times^G\GL_n\stackrel{\sim}{\to}T$.
Conversely, given $(H,\omega)$ over $S$, we deduce a $G$-equivariant
morphism $H\to T$; after taking quotients by the left action 
of $G$, we deduce a morphism $S\simeq H/G\to T/G$.
It is easy to verify that these two rules define mutually
inverse natural transformations of functors, whence the claim.
\end{proof}

\sset\subsubsection{}
The right $\GL_n$-torsor $T$ is endowed by a natural left
action by the group $\Aut_{\GL_n}(T)$. The latter is a group
scheme locally isomorphic to $\GL_n$ in the Zariski topology;
especially it is smooth over $R$. After taking the
quotient of $T$ by the right action of $G$ we deduce
a left action of $\Aut_{\GL_n}(T)$ on $T/G$ :
\set\begin{equation}\label{eq_transparent}
\Aut_{\GL_n}(T)\times_RT/G\to T/G.
\end{equation}

\begin{theorem}
Resume the assumptions of proposition {\em\ref{prop_denseimage}}.
Let $G$ be a group scheme over $\Spec\,R[t^{-1}]$ which satisfies
(relative to $R[t^{-1}]$) the conditions of \eqref{subsec_G-torsors}.
Then the natural map
$$
H^1(\Spec\,R[t^{-1}]_\fpqc,G)\to 
H^1(\Spec\,R^\wedge[t^{-1}]_\fpqc,G)
$$
is a bijection.
\end{theorem}
\begin{proof} We begin with the following special case:
\begin{claim}\label{cl_OKforGL} The theorem holds when $G=\GL_n$.
\end{claim}
\begin{pfclaim} Indeed, a $\GL_n$ torsor over
a scheme $X$ is the same as a locally free $\cO_X$-module
of rank $n$. Hence the assertion is just a restatement of 
corollary \ref{cor_elkik}.
\end{pfclaim}

For a given $\GL_n$-torsor $T$, let 
$T^\wedge:=T\times_{R[t^{-1}]}R^\wedge[t^{-1}]$ and denote
by 
$$
H^1_T\subset H^1(\Spec\,R[t^{-1}]_\fpqc,G)
$$
the subset consisting of all classes of $G$-torsors $H$ such that 
$H\times^G\GL_n\simeq T$; define likewise the subset 
$H^1_{T^\wedge}\subset H^1(\Spec\,R^\wedge[t^{-1}]_\fpqc,G)$. 
According to claim \ref{cl_OKforGL} it suffices to show:
\begin{claim} The restriction $H^1_T\to H^1_{T^\wedge}$
is a bijection.
\end{claim}
\begin{pfclaim} The morphism \eqref{eq_transparent} defines
a groupoid $T/G$ of quasi-projective $R[t^{-1}]$-schemes;
$T/G$ is locally isomorphic to $X_G$ in the Zariski topology,
hence it is smooth (see \ref{subsec_G-torsors}) and 
we can therefore apply theorem \ref{th_groupoid}. The assertion
follows directly after one has remarked that, for every 
$R[t^{-1}]$-scheme $S$, the set $\pi_0(T/G(S))$ is in natural
bijection with the set of isomorphism classes of $G$-torsors $H$
on $S_\fpqc$ such that $H\times^G\GL_n\simeq T\times_{\Spec\,R[t^{-1}]}S$.
\end{pfclaim}
\end{proof}

\sset\subsubsection{}\label{subsec_lift-tors}
Finally we come to the lifting problems for $G^a$-torsors. 
Let $A$ be a $V^a$-algebra, set $R:=A_*$ and let $G$ a smooth
affine group scheme of finite type over $\Spec\,R$. We suppose
that $G$ is a closed subgroup scheme of $\GL_{n,R}$, for some $n\in\N$.
Clearly $G^a$ is a group scheme over $S:=\Spec\,A$. Unless
otherwise stated, $G^a$ will be acting on the right on any such
torsor.

\begin{lemma}\label{lem_obvious}
Keep the notation of \eqref{subsec_lift-tors}. The category of
$\GL_n$-torsors over $S_\fpqc$ is naturally equivalent to the
category whose objects are the almost projective $\cO_S$-modules
of constant rank $n$, and whose morphisms are the linear isomorphisms.
\end{lemma}
\begin{proof} Let $\cO_*$ be the structure sheaf of $S_\fpqc$,
{\em i.e.} the sheaf of rings defined by the rule: $T\mapsto\cO_{T*}$.
To a given $\GL_n$-torsor $P$ we assign the $\cO_*$-module 
$F_P:=P\times^{\GL_n}\cO_*^n$ (for the natural left action of
$\GL_n$ on $\cO_*$). 
Locally in $S_\fpqc$ the sheaf $F_P$ is a free $\cO_*$-module of rank $n$,
and the assignment $P\mapsto F_P$ is an equivalence from $\GL_n$-torsors
to the category of all such locally free $\cO_*$-modules $F$ of rank $n$,
whose inverse is the rule: $F\mapsto\Iso_{\cO_*}(\cO_*^n,F)$.
On the other hand, theorem \ref{th_free.in.fpqc} says that any almost
projective $\cO_S$-module $M$ of constant rank $n$ defines a locally free
$\cO_*$-module in $S_\fpqc$, by the rule
$T\mapsto(\cO_T\otimes_{\cO_S}M)_*$  and by faithfully flat descent
it is also clear that the resulting functor is an equivalence as well.
\end{proof}

\begin{theorem}\label{th_lift-iso-tors}
Assume that $\mathrm{hom.dim}_V\tilde\fm\leq 1$.
Let $S$ and $G$ be as in \eqref{subsec_lift-tors},
and $I\subset\cO_S$ an ideal such that the pair $(\cO_S,I)$ is
henselian. Let also $\fm_0\subset\fm$ be a finitely generated
subideal, and set $S_n:=\Spec\,\cO_S/\fm^n_0I$ for $n=0,1$.
Let $P$ and $Q$ be two $G^a$-torsors over $S$, and suppose that
$\bar\beta:P\times_SS_1\stackrel{\sim}{\to}Q\times_SS_1$ is an
isomorphism of $G^a\times_SS_1$-torsors. Then
$\beta\times_S\one_{S_0}$ lifts to an isomorphism
$\beta:P\stackrel{\sim}{\to}Q$ of $G^a$-torsors.
\end{theorem}
\begin{proof} Let $\Iso(P,Q)$ denote the functor that assigns
to every $S$-scheme $T$ the set of isomorphisms of 
$G^a\times_ST$-torsors $P\times_ST\stackrel{\sim}{\to}Q\times_ST$.
It is easy to see that $\Iso(P,Q)$ is a sheaf for the fpqc
topology of $S$. Since, locally in $S_\mathrm{fpqc}$, this
sheaf is represented by the almost scheme $G^a$, it follows by 
faitfully flat descent that $\Iso(P,Q)$ is represented by
an affine $S$-scheme, which we denote by the same name. To the
(right) $G^a$-torsor $P$ we associate a left $G^a$-torsor $P'$, 
whose underlying $S$-scheme is the same as $P$, and whose 
$G^a$-action is defined by the rule: $g\cdot x:=xg^{-1}$
for every $S$-scheme $T$, every $g\in G^a(T)$ and every
$x\in P(T)$. Furthermore, we let $H$ denote the trivial left
and right $G^a$-torsor whose underlying $S$-scheme is $G^a$,
and whose left and right $G^a$-actions are induced by the
multiplication map $G^a\times_S G^a\to G^a$.

\begin{claim}\label{cl_double-act}
There is a natural isomorphism of sheaves on $S_\mathrm{fpqc}$:
$$
\omega:\Iso(P,Q)\simeq Q\times^{G^a}H\times^{G^a}P'.
$$
\end{claim}
\begin{pfclaim} Recall that the meaning of the right-hand
side is as follows. For every $S$-scheme $T$, one consider the 
presheaf whose $T$-sections of are the equivalence
classes of triples $(x,h,y)\in Q(T)\times H(T)\times P'(T)$,
modulo the equivalence relation such that
$(x,g_1\cdot h\cdot g_2,y)\sim(x\cdot g_1,h, g_2\cdot y)$
for every such $(x,h,y)$ and every $g_1,g_2\in G^a(T)$.
Then $Q\times^{G^a}H\times^{G^a}P'$ is the sheaf associated
to this presheaf. The sought isomorphism is defined as follows.
Let $T$ be an $S$-scheme, and $\beta$ a $T$-section of $\Iso(P,Q)$.
We pick a covering morphism $U\to T$ in $S_\mathrm{fpqc}$
such that $P(U)\neq\emptyset$; let $y\in P(U)$; we set
$\omega(\beta,y):=(\beta(y),1,y)$. If now $z\in P(U)$ is any
other section, we have $z=yg$ for some $g\in G^a(U)$,
therefore: 
$\omega(\beta,y)=(\beta(y),g\cdot g^{-1},y)\sim
(\beta(y)\cdot g,1,g^{-1}y)=(\beta(yg),1,yg)=\omega(\beta,z)$,
{\em i.e.} the equivalence class of $\omega(\beta,y)$ does
not depend on the choice of $y$, hence we have a well defined
$U$-section $\omega(\beta)$ of $Q\times^{G^a}H\times^{G^a}P'$.
Furthermore, let $p_i:U\times_TU\to U$, $i=1,2$ be the two
projections; the sections $p_1^*\omega(\beta,y)$ and
$p_2^*\omega(\beta,y)$ differ by an element of $G^a(U\times_TU)$
so, by the same argument, they lie in the same equivalence class,
which means that actually $\omega(\beta)$ comes from $T$, as required.
We leave to the reader the verification that $\omega$ thus defined
is an isomorphism.
\end{pfclaim}

Composing the closed imbedding $G^a\subset\GL_{n,S}$ with
the standard group homomorphism $\GL_{n,S}\subset SL_{n+1,S}$,
we can view $G^a$ as a closed subgroup scheme of $SL_{n,S}$,
whence a closed $G^a$-equivariant imbedding of $S$-schemes 
$G^a\subset M_{n,S}:=\Spec\,\cO_S[x_{ij}~|~i,j\leq n]$.
From claim \ref{cl_double-act} we derive a closed imbedding
of $S$-schemes:
$$
\Iso(P,Q)\subset X:=Q\times^{G^a}M_{n,S}\times^{G^a}P'
$$
(where $G^a$ acts on the left and right on $M_{n,S}$ in the
obvious way). $X$ is locally isomorphic to $M_{n,S}$
in $S_\mathrm{fpqc}$, hence it is of the form 
$\Spec\,(\Sym^\bullet_{\cO_S}M)$, for an almost finitely
generated projective $\cO_S$-module $M$; especially, $X$ is
almost finitely presented over $S$, whence the same holds
for $\Iso(P,Q)$. Now, the given $\bar\beta$ gives a section
of $\Iso(P,Q)$ over $S_1$, so the theorem is an immediate
consequence of theorem \ref{th_lift-smooth-sects}(i).
\end{proof}

\begin{theorem}\label{th_last-one}
Keep the notation and assumptions of theorem
{\em\ref{th_lift-iso-tors}}, and suppose furthermore that $I$ is
tight, and that $G$ fulfills the conditions of 
\eqref{subsec_G-torsors}. Let $P_0$ be any $G^a$-torsor
over $S_0$. Then there exists a $G^a$-torsor $P$ over $S$ with an
isomorphism of $G^a$-torsors: $P_0\simeq P\times_SS_0$.
\end{theorem}
\begin{proof} By assumption there exists a finitely generated
subideal $\fm_0\subset\fm$ and an integer $n>0$ such that 
$I^n\subset\fm_0\cO_S$; then by theorem \ref{th_wrapup}(ii) we
can lift $P_0$ to a $G$-torsor over $S_1:=\cO_S/\fm_0I$.
We form the $\GL_n$-torsor $\bar Q:=\bar P\times^{G^a}\GL_n$.
By lemma \ref{lem_obvious}, $\bar Q$ is the same as an almost
projective $\cO_{S_1}$-module of constant rank $n$; by theorem
\ref{th_lifts.proj.hensel}(i) we can then find an almost
projective $\cO_S$-module $Q$ that lifts $\bar Q$,
and then a standard application of Nakayama's lemma \ref{lem_Naka}
shows that $Q$ has constant rank equal to $n$ (cp. the proof of
theorem \ref{th_raise-smooth-sects}). Thus, $Q$ is a $\GL_n$-torsor
with an isomorphism of $\GL_n$-torsors
$\omega:\bar P\times^{G^a}\GL_n\stackrel{\sim}{\to}Q\times_SS_1$.
By (the almost version of) lemma \ref{lem_neatsol}, the datum of
$(\bar P,\omega)$ is the same as the datum of a morphism 
$\bar\sigma:S_1\to Q/G^a$ of sheaves on
$S_\fpqc$. The contention is that $\bar\sigma\times_SS_0$ lifts
to a morphism $\sigma:S\to Q/G^a$. However, we have a natural
$\GL_n$-equivariant isomorphism
$$
Q/G^a\simeq Q\times^{\GL_n}(\GL_n/G^a)\simeq Q\times^{\GL_n}(\GL_n/G)^a.
$$
By lemmata \ref{lem_Idontknow} and \ref{lem_equivariant}
it follows that there is a $\GL_n$-equivariant isomorphism:
$$
Q/G^a\simeq\Proj(\Sym^\bullet_{\cO_S}\cP')/\cI')\setminus V(\cJ')
$$
with $\cP':=Q\times^{\GL_n}\cP$, $\cI':=\cI\times^{\GL_n}\cP$
and $\cJ':=\cJ\times^{\GL_n}\cP$, where $\cP$, $\cI$ and $\cJ$
are as in \eqref{eq_equivariant}. By lemma
\ref{lem_flat.proj}(ii),(iii) and remark \ref{rem_descent}(i),(ii)
we see that $\cP'$ is almost finitely generated projective and
$\cI'$, $\cJ'$ are almost finitely generated. Finally, $Q/G^a$
is smooth since $X_G$ is (see \ref{subsec_G-torsors}), so the
existence of the section $\sigma$ follows from theorem
\ref{th_raise-smooth-sects}.
\end{proof}


\newpage

\section{Valuation theory}\label{ch_val.theory}

This chapter is an extended detour into valuation 
theory. The first two sections contain nothing new, and are
only meant to gather in a single place some useful material
that is known to experts, but for which satisfactory references
are hard to find. The main theme of sections \ref{sec_val.rings}
through \ref{sec_trasc_exts} is the study of the cotangent
complex of an extension of valuation rings. To give a sample
of our results, suppose that $k$ is a perfect field, and let
$W$ be a valuation ring containing $k$; then we show that
$\Omega_{W/k}$ is a torsion-free $W$-module. Notice that this
assertion would be an easy consequence of the existence of
resolution of singularities for $k$-schemes; our methods enable
us to prove it unconditionally, as well as several other
statements and variants for logarithmic differentials.
Furthermore, consider a finite separable extension $K\subset L$
of henselian valued fields of rank one; it is not difficult to see
that the corresponding extension of valuation rings $K^+\subset L^+$
is almost finite, hence one can define the different ideal
$\cD_{L^+/K^+}$ as in chapter \ref{ch_traces}. In case the
valuation of $K$ is discrete, it is well known that the length
of the module of relative differentials $\Omega_{L^+/K^+}$
equals the length of $L^+/\cD_{L^+/K^+}$; theorem \ref{th_big.deal}
generalizes this identity to the case of arbitrary rank one
valuations; notice that in this case the usual notion of length
won't do, since the modules under considerations are only
almost finitely generated.

Section \ref{sec_deeply.ram} ties up with earlier work of 
Coates and Greenberg \cite{Coat}, in which the notion of
{\em deeply ramified extension\/} of a local field was introduced,
and applied to the study of $p$-divisible groups attached to
abelian varieties defined over such $p$-adic fields. Essentially,
section 2 of \cite{Coat} rediscovers the results of Fresnel and
Matignon \cite{Fr-Ma}, although via a different route, closer to
the original treatment of Tate in \cite{Ta}.
In particular, an algebraic extension $E$ of $K$ is
deeply ramified if and only if $\cD_{E^+/K^+}=(0)$ according to 
the terminology of \cite{Fr-Ma}.  We adopt Coates and Greenberg's 
terminology for our section \ref{sec_deeply.ram}, and we give some 
complements which were not observed in \cite{Coat}; notably,
proposition \ref{prop_deep.ramif}, which we regard as the ultimate
generalization of the one-dimensional case of the almost purity
theorem. Our proof generalizes Faltings' method, which relied on
the above mentioned relationship between differentials and the
different ideal; of course, in the present setting we need to appeal
to our theorem \ref{th_big.deal}, rather than estimating lengths
the way Faltings did. Finally, we extend the definition of
deeply ramified extension to include valued fields of arbitrary rank.

\subsection{Ordered groups and valuations}
In this section we gather some generalities on valuations
and related ordered groups, which will be used in later
sections.

\sset\subsubsection{}\label{subsec_def.values}
\index{Valued field(s)|indref{subsec_def.values}}
\index{Valuation(s)|indref{subsec_def.values}}
\index{Valued field(s)!$\Gamma_K$ : value group of a|indref{subsec_def.values}}
As usual, a {\em valued field\/} $(K,|\cdot|_K)$ consists of 
a field $K$ endowed with a surjective group homomorphism 
$|\cdot|_K:K^\times\to\Gamma_K$ 
onto an ordered abelian group $(\Gamma_K,\leq)$, such that 
\set\begin{equation}\label{eq_ineq.abs.value}
|x+y|_K\leq\max(|x|_K,|y|_K)
\end{equation}
whenever $x+y\neq 0$. We denote by $1$ the neutral element
of $\Gamma_K$, and the composition law of $\Gamma_K$ will be
denoted by: $(x,y)\mapsto x\cdot y$. It is customary to 
extend the map $|\cdot|_K$ to the whole of $K$, by adding a
new element $0$ to the set $\Gamma_K$, and setting 
$|0|:=0$. One can then extend the ordering of $\Gamma_K$ to 
$\Gamma_K\cup\{0\}$ by declaring that $0$ is the smallest 
element of the resulting ordered set. In this way,
\eqref{eq_ineq.abs.value} holds for every $x,y\in K$.
The map $|\cdot|_K$ is called the {\em valuation\/} of $K$
and $\Gamma_K$ is its {\em value group}. 

\sset\subsubsection{}\label{subsec_ext.of.val.fld}
\index{Valued field(s)!extension of a|indref{subsec_ext.of.val.fld}}
An {\em extension of valued fields\/} 
$(K,|\cdot|_K)\subset(E,|\cdot|_E)$ consists of a field extension
$K\subset E$, and a valuation 
$|\cdot|_E:E\to\Gamma_E\cup\{0\}$ together with an imbedding
$j:\Gamma_K\subset\Gamma_E$, such that the restriction to $K$ 
of $|\cdot|_E$ equals $j\circ|\cdot|_K$.

\begin{example}\label{ex_Gauss.norm}
\index{Valuation(s)!Gauss|indref{ex_Gauss.norm}}
Let $|\cdot|:K\to\Gamma\cup\{0\}$ be a valuation on the 
field $K$.

(i) Given a field extension $K\subset E$, it is known that 
there always exist valuations on $E$ which extend $|\cdot|$ 
(cp. \cite[Ch.VI, \S1, n.3, Cor.3]{BouAC}).

(ii) If the field extension $K\subset E$ is algebraic and 
purely inseparable, then the extension of $|\cdot|$ is 
unique. (cp. \cite[Ch.VI, \S8, n.7, Cor.2]{BouAC}).

(iii) We can construct extensions of $|\cdot|$ on the polynomial 
ring $K[X]$, in the following way. Let $\Gamma'$ be an ordered
group with an imbedding of ordered groups $\Gamma\subset\Gamma'$.
For every $x_0\in K$, and every $\rho\in\Gamma'$, we define 
the {\em Gauss valuation\/ $|\cdot|_{(x_0,\rho)}:K[X]\to
\Gamma'\cup\{0\}$ centered at $x_0$ and with radius $\rho$} 
(cp. \cite[Ch.VI, \S10, n.1, Lemma 1]{BouAC}) by the rule:
$$a_0+a_1(X-x_0)+...+a_n(X-x_0)^n\mapsto
\max\{|a_i|\cdot\rho^i~|~i=0,1,...,n\}.$$

(iv) The construction of (iii) can be iterated : for instance,
suppose that we are given a sequence of $k$ elements 
$\rho:=(\rho_1,\rho_2,...,\rho_k)$
of the ordered abelian group $\Gamma'$ of (iv). Then we
can define a Gauss valuation $|\cdot|_{(0,\rho)}$
on the fraction field of $K[X_1,X_2,...,X_k]$, with values 
in $\Gamma'$, by the rule: 
$\sum_{\alpha\in\N^k}a_\alpha X^\alpha\mapsto
\max\{|a_\alpha|\cdot\rho^\alpha~|~\alpha\in\N\}$.

(v) Suppose again it is given an ordered group $\Gamma'$
with an imbedding of ordered groups $\Gamma\subset\Gamma'$.
Let $T\subset\Gamma'/\Gamma$ be a finite torsion subgroup,
say $T\simeq\Z/n_1\Z\oplus...\oplus\Z/n_k\Z$. For every 
$i\leq k$, pick an element $\gamma_i\in\Gamma'$ whose class 
in $\Gamma'/\Gamma$ generates the direct summands $\Z/n_i\Z$
of $T$. Let $x_i\in K$ such that $|x_i|=n_i\cdot\gamma_i$. 
For every $i=1,...,k$ pick an element $y_i$ in a fixed
algebraic closure $E^\mathrm{a}$ of $E$, such that $y_i^{n_i}=x_i$; 
then the field $E:=K(y_1,...,y_k)$ has degree over $K$ equal 
to the order of $T$, and it admits a unique valuation $|\cdot|_E$ 
extending $|\cdot|$. Of course, $|y_i|_E=\gamma_i$ for every 
$i\leq k$.
\end{example}

\sset\subsubsection{}\label{subsec_datum}
We want to explain a construction which is a simultaneous
generalization of examples \ref{ex_Gauss.norm}(iv),(v).
Suppose it is given the datum $\fG:=(G,j,N,\leq)$ 
consisting of:
\begin{enumerate}
\renewcommand{\labelenumi}{(\alph{enumi})}
\item
an abelian group $G$ with an imbedding 
$j:K^\times\hookrightarrow G$ such that $G/j(V^\times)$
is torsion-free;
\item
a subgroup $N$ of $G/j(V^\times)$ such that the 
natural map: 
$$
\Gamma\stackrel{\sim}{\to}K^\times/V^\times
\to\Gamma_\fG:=G/(N+j(V^\times))
$$
is injective;
\item 
an ordering $\leq$ on $\Gamma_\fG$ such that the injective
map: $\Gamma\to\Gamma_\fG$ is order-preserving.
\renewcommand{\labelenumi}{(\roman{enumi})}
\end{enumerate}
Let us denote by $K[G]$ (resp. $K[K^\times]$)
the group $K$-algebra of the abelian group $G$ (resp. 
$K^\times$). Any element of $K[G]$ can be written uniquely
as a formal linear combination $\sum_{g\in G}a_g\cdot[g]$,
where $a_g\in K$ for every $g\in G$, and $a_g=0$ 
for all but a finite number of $g\in G$. We augment 
$K[K^\times]$ over $K$ via the $K$-algebra homomorphism 
\set\begin{equation}\label{eq_augment}
K[K^\times]\to K\quad:\quad [a]\mapsto a\quad
\text{for every $a\in K^\times$}.
\end{equation}
Then we let $K[\fG]:=K[G]\otimes_{K[K^\times]}K$,
where the $K[K^\times]$-algebra structure on $K$
is defined by the augmentation \eqref{eq_augment}.
It is easy to verify that $K[\fG]$ is the maximal
quotient algebra of $K[G]$ that identifies the classes
of $[g\cdot a]$ and $a\cdot[g]$, for every $g\in G$ and 
$a\in K^\times$. 
Pick, for every class $\gamma\in G/K^\times$,
a representative $g_\gamma\in G$. It follows that 
every element of $K[\fG]$ can be written uniquely
as a formal $K$-linear combination 
$\sum_{\gamma\in G/K^\times}a_\gamma\cdot[g_\gamma]$.
We define a map $|\cdot|_\fG:K[\fG]\to\Gamma_\fG\cup\{0\}$ 
by the rule:
\set\begin{equation}\label{eq_inspecting.form}
\sum_{\gamma\in G/K^\times}a_\gamma\cdot[g_\gamma]
\mapsto\max_{\gamma\in G/K^\times}
|a_\gamma|\cdot|g_\gamma|
\end{equation}
where $|g_\gamma|\in\Gamma_\fG$ denotes the class of
$g_\gamma$. One verifies easily that $|\cdot|_\fG$ does
not depend on the choice of representatives $g_\gamma$.
Indeed, if $(h_\gamma~|~\gamma\in G/K^\times)$ is
another choice then, for every $\gamma\in G/K^\times$ 
we have $g_\gamma=j(x_\gamma)\cdot h_\gamma$ for some 
$x_\gamma\in K^\times$; therefore 
$[g_\gamma]=x_\gamma\cdot[h_\gamma]$ and
$|g_\gamma|=|x_\gamma|\cdot|h_\gamma|$. 

\begin{lemma} $K[\fG]$ is an integral domain, and 
$|\cdot|_\fG$ extends to a valuation:
$$|\cdot|_\fG:K(\fG):=\mathrm{Frac}(K[\fG])\to
\Gamma_\fG\cup\{0\}.$$
\end{lemma}
\begin{proof} Let $(G_\alpha~|~\alpha\in I)$ be the
filtered system of the subgroups $G_\alpha$ of $G$ such 
that $K^\times\subset G_\alpha$ and $G_\alpha/K^\times$
is finitely generated. Each $G_\alpha$ defines a datum
$\fG_\alpha:=(G_\alpha,j,N\cap(G_\alpha/j(V^\times)),\leq)$,
and clearly $K[\fG]=\colim{\alpha\in I}K[\fG_\alpha]$.
We can therefore reduce to the case where 
$G/K^\times$ is finitely generated. 
Write $G/K^\times=T\oplus F$, where $T$ is a torsion
group and $F$ is torsion-free. There exist unique subgroups
$\tilde T,\tilde F\supset K^\times$ in $G$ with 
$\tilde T/K^\times=T$ and $\tilde F/K^\times=F$.
Let $\fG_T:=(\tilde T,j,\{0\},\leq)$ and 
$\fG_F:=(\tilde F,j,N,\leq)$ be the corresponding data.
The functor $H\mapsto K[H]$ preserves colimits, since it
is left adjoint to the forgetful functor from $K$-algebras 
to abelian groups; it follows easily that 
$K[\fG]\simeq K[\fG_T]\otimes_KK[\fG_F]$.
By inspecting \eqref{eq_inspecting.form}, one can easily 
show that $K[\fG_T]$ is of the type of example 
\ref{ex_Gauss.norm}(v) and $K[\fG_F]$ is of the type of
example \ref{ex_Gauss.norm}(iv). Especially, $K[\fG]$
is a domain, and $|\cdot|_\fG$ is induced by a Gauss 
valuation of a free algebra over the finite field extension 
$K[\fG_T]$ of $K$.
\end{proof}

The next result shows that an arbitrary valuation is 
always ``close" to some Gauss valuation.

\begin{lemma}\label{lem_approx.sequence}
Let $(E,|\cdot|_E)$ be a valued field extension of 
$(K,|\cdot|)$.
Let $x\in E\setminus K$, and let $(a_i~|~i\in I)$ be 
a net of elements of $K$ (indexed by the directed set 
$(I,\leq)$) with the following property. 
For every $b\in K$ there exists $i_0\in I$ such that 
$|x-a_i|_E\leq|x-b|_E$ for every $i\geq i_0$.
Let $f(X)\in K[X]$ be a polynomial
that splits in $K[X]$ as a product of linear polynomials. 
Then there exists $i_0\in I$ such that 
$|f(x)|_E=|f(X)|_{(a_i,|x-a_i|_E)}$ for every $i\geq i_0$.
\end{lemma}
\begin{proof} To prove the claim, it suffices to consider 
the case when $f(X)=X-b$ for some $b\in K$. However, from 
the definition of the sequence $(a_i~|~i\in I)$ we have 
$\max(|x-a_i|_E,|a_i-b|)\geq|x-b|_E\geq|x-a_i|_E$ for 
every sufficiently large $i\in I$. Therefore, 
$|x-b|_E=\max(|x-a_i|_E,|b-a_i|)=|X-b|_{(a_i,|x-a_i|_E)}$.
\end{proof}

\sset\subsubsection{}\label{subsec_to_see_how}
To see how to apply lemma \ref{lem_approx.sequence},
let us consider the case where $K$ is algebraically
closed and $E=K(X)$, the field of fractions of the
free $K$-algebra in one generator, which we suppose
endowed with some valuation $|\cdot|_E$ with values
in $\Gamma_E$. We apply lemma \ref{lem_approx.sequence} 
to the element $x:=X\in E$. Suppose first that there 
exists an element $a\in K$ that minimizes the function
$K\to\Gamma_E~:~b\mapsto|X-b|_E$. In this case
the trivial net $\{a\}$ fulfills the condition of
the lemma. Since every polynomial of $K[X]$ splits
over $K$, we see that $|\cdot|_E$ {\em is the Gauss 
valuation centered at $a$ and with radius $|X-a|_E$}.
Suppose, on the other hand, that the function
$b\mapsto|X-b|_E$ does not admit a minimum. It will 
still be possible to choose a net of elements 
$\{a_i~|~i\in I\}$ fulfilling the conditions of lemma 
\ref{lem_approx.sequence} (indexed, for instance, by a 
subset of the partially ordered set $\Gamma$).
Then $|\cdot|_E$ is determined by the identity:
$$
|f(X)|_E=\lim_{i\in I}|f(X)|_{(a_i,|X-a_i|_E)}
\quad\text{for every $f(X)\in E$.}
$$
\sset\subsubsection{}\label{subsec_charct.val.rings}
\index{$K^+$ : Valuation ring(s)|indref{subsec_charct.val.rings}}
Given a valuation $|\cdot|$ on a field $K$, the subset 
$K^+:=\{x\in ~|~|x|\leq 1\}$ is a {\em valuation ring\/} 
of $K$, {\em i.e.}, a subring of $K$ such that, for every
$x\in K\setminus\{0\}$, either $x\in K^+$ or $x^{-1}\in K^+$.
The subset $(K^+)^\times$ of units of $K^+$ consists precisely 
of the elements $x\in K$ such that $|x|=1$. Conversely, 
let $V$ be a valuation ring of $K$ with maximal ideal $\fm$;
$V$ induces a valuation $|\cdot|$ on $K$ whose value group is 
$\Gamma_K:=K^\times/V^\times$ (then $|\cdot|$ is just the 
natural projection). The ordering on $\Gamma_K$ is defined
as follows. For given classes $\bar x,\bar y\in\Gamma_K$, 
we declare that $\bar x<\bar y$ if and only if $x/y\in\fm$.

\begin{remark}\label{rem_val.ring.princip}
\index{Pr{\"u}fer Domain|indref{rem_val.ring.princip}}
\index{$K^+$ : Valuation ring(s)!$K^{+\mathrm{h}}$, $K^{+\mathrm{sh}}$,
$K^\mathrm{sh}$ : (strict) henselization of 
a|indref{rem_val.ring.princip}{}, \indref{def_max.tame.ext}}
(i) It follows easily from \eqref{subsec_charct.val.rings}
\index{Valued field(s)!$K^\mathrm{s}$, $K^\mathrm{a}$ : 
separable, algebraic closure of a|indref{def_max.tame.ext}}
that every finitely generated ideal of a valuation ring is 
principal. Indeed, if $a_1,...,a_n$ is a set of generators 
for an ideal $I$, pick $i_0\leq n$ such that 
$|a_{i_0}|=\max_{i\leq n}|a_i|$; then $I=(a_{i_0})$.

(ii) It is also easy to show that any finitely generated 
torsion-free $K^+$-module is free and any torsion-free 
$K^+$-module is flat 
(cp. \cite[Ch.VI, \S3, n.6, Lemma 1]{BouAC}).

(iii) Let $E$ be a field extension of the valued field 
$(K,|\cdot|_K)$. Then the integral closure $W$ of $K^+$ 
in $E$ is the intersection of all the valuation rings 
of $E$ containing $K^+$ 
(cp. \cite[Ch.VI, \S1, n.3, Cor.3]{BouAC}). In particular,
$K^+$ is integrally closed.

(iv) Furthermore, if $E$ is an algebraic extension of
$K$, then $W$ is a {\em Pr{\"u}fer domain}, that is, for 
every prime ideal $\fp\subset W$, the localization $W_\fp$
is a valuation ring. Moreover, the assignment 
$\fm\mapsto W_\fm$ establishes a bijection between the
set of maximal ideals of $W$ and the set of valuation
rings $V$ of $E$ whose associated valuation $|\cdot|_V$ 
extends $|\cdot|_K$ 
(cp. \cite[Ch.VI, \S8, n.6, Prop.6]{BouAC}).

(v) Let $R$ and $S$ be local rings contained in a field
$K$, $\fm_R$ and $\fm_S$ their respective maximal ideals. 
One says that $R$ {\em dominates\/} $S$ if $S\subset R$
and $\fm_S=S\cap\fm_R$. It is clear that the relation
of dominance establishes a partial order structure on
the set of local subrings of $K$. Then a local subring
of $K$ is a valuation ring of $K$ if and only if it is
maximal for the dominance relation (cp. 
\cite[Ch.VI, \S1, n.2, Th.1]{BouAC}).

(vi) Let $K^+$ be a valuation ring of $K$ with maximal
ideal $\fm$, and $K^{+\mathrm{h}}$ a henselization of $K^+$. 
One knows that $K^{+\mathrm{h}}$ is an ind-{\'e}tale local 
$K^+$-algebra (cp. \cite[Ch.VIII, Th.1]{Ray}), hence it 
is integral and integrally closed (cp. 
\cite[Ch.VII, \S2, Prop.2]{Ray}).
Denote by $K^\mathrm{h}$ the field of fractions of 
$K^{+\mathrm{h}}$ and $W$ the integral closure of $K^+$
in $K^\mathrm{h}$. It follows that $W\subset K^{+\mathrm{h}}$.
Let $\fm^\mathrm{h}$ be the maximal ideal of 
$K^{+\mathrm{h}}$; since $\fm^\mathrm{h}\cap K^+=\fm$,
we deduce that $\fq:=\fm^\mathrm{h}\cap W$ is a maximal ideal
of $W$; then by (iv), $W_\fq$ is a valuation ring of 
$K^\mathrm{h}$ dominated by $K^{+\mathrm{h}}$; by (v) it 
follows that $K^{+\mathrm{h}}=W_\fq$, in particular this
shows that the henselization of a valuation ring is again
a valuation ring. The same argument works also for strict
henselizations.
\end{remark}

The following lemma provides a simple method to construct
extensions of valuation rings, which is sometimes useful.

\begin{lemma}\label{lem_simple.meth} Let $(K,|\cdot|)$ be 
a valued field, $\kappa$ the residue field of $K^+$, $R$ a 
$K^+$-algebra which is finitely generated free as a 
$K^+$-module, and suppose that $R\otimes_{K^+}\kappa$ is a 
field. Then $R$ is a valuation ring, and the morphism 
$K^+\to R$ induces an isomorphism of value groups 
$\Gamma_K\stackrel{\sim}{\to}\Gamma_R$.
\end{lemma}
\begin{proof} Let $e_1,...,e_n$ be a $K^+$-basis of $R$.
Let us define a map $|\cdot|_R:R\to\Gamma_R\cup\{0\}$ in 
the following way. Given $x\in R$, write 
$x=\sum_{i=1}^nx_i\cdot e_i$; then 
$|x|_R:=\max\{|x_i|~|~i=1,...,n\}$. If $|x|=1$, then
the image $\bar x$ of $x$ in $R\otimes_{K^+}\kappa$ is not 
zero, hence it is invertible by hypothesis. By Nakayama's
lemma it follows easily that $x$ itself is invertible in
$R$. Hence, every element $y$ of $R$ can be written in the form
$y=u\cdot b$, where $u\in R^\times$ and $b\in K^+$ is an element
such that $|b|=|y|_R$. It follows easily that $R$ is an
integral domain. Moreover, it is also clear that, given
any $x\in\mathrm{Frac}(R)\setminus\{0\}$, either $x\in R$ 
or $x^{-1}\in R$, so $R$ is indeed a valuation ring and
$|\cdot|_R$ is its valuation.
\end{proof}

\begin{lemma}\label{lem_fin.pres.torsion}
Every finitely presented torsion $K^+$-module $M$ is 
isomorphic to a direct sum of the form 
$$(K^+/a_1K^+)\oplus...\oplus(K^+/a_nK^+)$$
where $a_1,...,a_n\in K^+$. More precisely, if 
$F\stackrel{\phi}{\to}M$ is any surjection from a free 
$K^+$-module $F$ of rank $n$, then there is a basis 
$e_1,...,e_n$ of $F$ and elements 
$a_1,...,a_n\in K^+\setminus\{0\}$ such that 
$\Ker\,\phi=(a_1K^+)\oplus...\oplus(a_nK^+)$.
\end{lemma}
\begin{proof} We proceed by induction on the rank $n$ of 
$F$. For $n=1$ the claim follows easily from remark
\ref{rem_val.ring.princip}(i).
Suppose $n>1$; first of all, $S:=\Ker(\phi)$ is finitely 
generated by \cite[Ch.I, \S2, n.8, lemme 9]{BouAC}.
Then $S$ is a free $K^+$-module, in light of remark
\ref{rem_val.ring.princip}(ii); its rank is necessarily
equal to $n$, since $S\otimes_{K^+}K=F\otimes_{K^+}K$.  

The image of the evaluation map $S\otimes_{K^+}F^*\to K^+$ 
given by $f\otimes\alpha\mapsto\alpha(f)$ is a finitely
generated ideal $I\neq 0$ of $K^+$, hence it is principal, 
by remark \ref{rem_val.ring.princip}(i). Let 
$\sum_{i=1}f_i\otimes\alpha_i$ be an element whose image 
generates $I$; this means that $\sum_{i=1}\alpha_i(f_i)$ 
is a generator of $I$, hence one of the terms in the sum, 
say $\alpha_1(f_1)$, is already a generator. The map 
$\alpha_1:S\to I$ is surjective onto a free rank one 
$K^+$-module, therefore it splits, which shows that 
$S=(f_1K^+)\oplus(S\cap\Ker~\alpha_1)$. In particular,
$S':=S\cap\Ker~\alpha_1$ is a finitely generated
torsion-free, hence free $K^+$-module.
Let $e_1,...,e_n$ be a basis of $F$; then 
$f_1=\sum_{i=1}^na_i\cdot e_i$ for some $a_i\in K^+$.
Consider the projection $\pi_i:F\to K^+$ such that
$\pi_i(e_j)=\delta_{ij}$ for $j=1,...,n$; clearly
$\pi_i(f_1)=a_i\in I$. This shows that 
$f_1=\alpha_1(f_1)\cdot g$ for some $g\in F$.
It follows that $\alpha_1(g)=1$, whence 
$F=(gK^+)\oplus\Ker\,\alpha_1$. Set $F':=\Ker\,\alpha_1$; 
we have shown that $M\simeq(K^+/\alpha_1(f_1))\oplus(F'/S')$.
But $F'$ is a free $K^+$-module of rank $n-1$,
hence we conclude by induction.
\end{proof}

\sset\subsubsection{}\label{subsec_standard.setup}
\index{$K^+$ : Valuation ring(s)!standard setup attached
to a|indref{subsec_standard.setup}}
In later sections we will be concerned with almost ring 
theory in the special case where the basic setup $(V,\fm)$
(see \ref{subsec_basic.setup})
consists of a valuation ring $V$. In preparation for
this, we fix the following terminology, which will stand
throughout the rest of this work. If $V$ is a valuation 
ring, then the {\em standard setup\/} attached to $V$ is 
the pair $(V,\fm)$ where $\fm:=V$ in case the value group 
of $V$ is isomorphic to $\Z$ ({\em i.e.} $V$ is a {\em 
discrete valuation ring}), and otherwise $\fm$ is the 
maximal ideal of $V$.

\sset\subsubsection{}\label{subsec_fract.ideals}
\index{Valued field(s)!$\mathrm{Div}(K^{+a})$ : fractional 
ideals of a|indref{subsec_fract.ideals}}
Let $K\to\Gamma_K\cup\{0\}$ be a valuation on the field $K$,
and $K^+$ its valuation ring. We consider the category
$K^{+a}\Mod$ relative to the standard setup $(K^+,\fm)$. 
The topological group $\mathrm{Div}(K^{+a})$ of {\em fractional 
ideals\/} of $K^{+a}$ is the subspace of $\cI_{K^{+a}}(K^a)$ 
which consists of all the submodules $I\neq K^a$ of the 
almost $K^+$-module $K^a$, such that the natural morphism 
$I\otimes_{K^{+a}}K^a\to K^a$ is an isomorphism. The group 
structure is induced by the multiplication of fractional 
ideals.

\begin{remark} One verifies easily that $\mathrm{Div}(K^{+a})$
is isomorphic to the group $D(K^+)$ defined in 
\cite[Ch.VII, \S1, n.1]{BouAC}.
\end{remark}

The structure of the ideals of $K^+$ can be largely 
read off from the value group $\Gamma$. In order to explain 
this, we are led to introduce some notions for general 
ordered abelian groups.

\sset\subsubsection{}
We endow an ordered group $\Gamma$ with the uniform structure 
defined in the following way. For every $\gamma\in\Gamma$ such 
that $\gamma>1$, the subset of $\Gamma\times\Gamma$ given by 
$E(\gamma):=
\{(\alpha,\beta)~|~\gamma^{-1}<\alpha^{-1}\cdot\beta<\gamma\}$
is an entourage for the uniform structure, and the subsets
of this kind form a fundamental system of entourages.
Let $\Gamma^\wedge$ be the completion of $\Gamma$ for
this uniform structure.

\begin{lemma}\label{lem_Picpus} With the notation of
\eqref{subsec_fract.ideals}, 
there exists a natural isomorphism of topological groups:
$\mathrm{Div}(K^{+a})\stackrel{\sim}{\to}\Gamma_K^\wedge$.
\end{lemma}
\begin{proof} We only indicate how to construct the
morphism, and leave the details to the reader.
In light of remark \ref{rem_val.ring.princip}(i), for 
every ideal $I\subset K^{+a}$ we can find a net 
$\{J_i~|~i\in S\}$ of principal ideals converging to 
$I$ (for some filtered ordered set $(S,\leq)$). Let 
$\gamma_i\in\Gamma_K$ be the value of a generator of $J_i$. 
One verifies that the net $\{\gamma_i~|~i\in S\}$ 
converges in $\Gamma_K^\wedge$ to some element $\hat\gamma$. 
Then we assign: $I\mapsto\hat\gamma$. One verifies
that this rule is well-defined and that it extends uniquely
to the whole of $\mathrm{Div}(K^{+a})$. 
\end{proof}

\begin{definition}\label{def_Gamma.plus}
\index{Ordered abelian group(s)!$\Gamma^+$ : submonoid of positive 
elements in an|indref{def_Gamma.plus}}
\index{Ordered abelian group(s)|indref{def_Gamma.plus}}
\index{Ordered abelian group(s)!convex subgroup of an|indref{def_Gamma.plus}}
\index{Ordered abelian group(s)!$\Spec\,\Gamma$ : spectrum of an|indref{def_Gamma.plus}}
\index{Ordered abelian group(s)!$\mathrm{c.rk}(\Gamma)$, $\rk(\Gamma)$ : 
convex and rational rank of an|indref{def_Gamma.plus}}
Let $\Gamma$ be any ordered abelian group with neutral 
element $1$.
\begin{enumerate}
\item
We denote by $\Gamma^+\subset\Gamma$ the subset of
all the $\gamma\in\Gamma$ such that $\gamma\leq 1$.
\item
A subgroup $\Delta$ of $\Gamma$ is said to be 
{\em convex\/} if it satisfies the following property. 
If $x\in\Delta^+$ and $1>y>x$, then $y\in\Delta$. The 
set $\Spec\,\Gamma$ of all the convex subgroups of 
$\Gamma$ will be called the {\em spectrum\/} of $\Gamma$. 
We define the {\em convex rank\/} of $\Gamma$ as the 
supremum $\mathrm{c.rk}(\Gamma)$ over the lengths $r$ of 
the chains  $0\subsetneq\Delta_1\subsetneq...\subsetneq
\Delta_r:=\Gamma$, such that all the $\Delta_i$ are 
convex subgroups. In general 
$\mathrm{c.rk}(\Gamma)\in\N\cup\{\infty\}$, but we will 
mainly encounter situations for which the convex rank 
is a positive integer. It is easy to see that the convex 
rank is always less than or equal to the usual rank, 
defined as $\rk(\Gamma):=\dim_\Q(\Gamma\otimes_\Z\Q)$. 
To keep the two apart, we call {\em rational rank\/} the 
latter.
\end{enumerate}
\end{definition}

\begin{example}\label{ex_tens.with.Q}
\index{Valued field(s)!valuation topology of a|indref{ex_tens.with.Q}}
(i) If $\Gamma$ is an ordered abelian group, there exists 
a unique ordered group structure on $\Gamma\otimes_\Z\Q$ 
such that the natural map $\Gamma\to\Gamma\otimes_\Z\Q$
is order-preserving. Indeed, if $\Gamma$ is the value
group of a valuation $|\cdot|$ on a field $K$, and 
$|\cdot|_{K^\mathrm{a}}$ is any extension of $|\cdot|_K$
to the algebraic closure $K^\mathrm{a}$ of $K$, then
it is easy to see ({\em e.g.} using example 
\ref{ex_Gauss.norm}(v)) that 
$\Gamma_{K^\mathrm{a}}\simeq\Gamma\otimes_\Z\Q$.

(ii) Furthermore, let $K^\sep\subset K^\mathrm{a}$ be the 
separable closure of $K$; we claim that 
$|\cdot|_{K^\mathrm{a}}$ maps $K^\sep$ surjectively
onto $\Gamma_{K^\mathrm{a}}$. Indeed, if $a\in K^\mathrm{a}$
is inseparable over $K^\sep$, then the minimal polynomial
$m(X)\in K^\sep[X]$ of $a$ is of the form $X^{p^m}-b$ for some
$b\in K^\sep$. For $c\in K^\times$, let $m_c(X)\in K^\sep[X]$
be the polynomial $m(X)+cX$; if $a'$ is a root of 
$m_c(X)$, then $a'\in K^\sep$; moreover, 
$|(a-a')^{p^m}|_{K^\mathrm{a}}=|c\cdot a'|_{K^\mathrm{a}}$,
hence for $|c|_{K^\mathrm{a}}$ sufficiently small we have
$|a|_{K^\mathrm{a}}=|a'|_{K^\mathrm{a}}$.

(iii) For any valued field $(K,|\cdot|)$, and every 
$\gamma\in\Gamma_K$, let 
$U_\gamma:=\{x\in K~|~|x|\leq\gamma\}$. One defines
the {\em valuation topology\/} on $K$ as the unique group
topology such that the family $(U_\gamma~|~\gamma\in\Gamma)$ 
is a fundamental system of open neighborhoods of $0$. 
The argument in (ii) shows more precisely that $K^\sep$
is dense in $K^\mathrm{a}$ for the valuation topology 
of $(K^\mathrm{a},|\cdot|_{K^\mathrm{a}})$.

(iv) If $\Delta\subset\Gamma$ is any subgroup, then 
$\mathrm{c.rk}(\Gamma)\leq\mathrm{c.rk}(\Delta)+
\rk(\Gamma/\Delta)$ 
(cp. \cite[Ch.VI, \S10, n.2, Prop.3]{BouAC}).

(v) A subgroup $\Delta\subset\Gamma$ is convex if and only 
if there is an ordered group structure on $\Gamma/\Delta$
such that the natural map $\Gamma\to\Gamma/\Delta$ is
order-preserving. Then the ordered group structure with
this property is unique.

(vi) If $\mathrm{c.rk}(\Gamma)=1$, we can find an 
order-preserving imbedding 
$$
\rho:(\Gamma,\cdot,\leq)\hookrightarrow(\R,+,\leq).
$$ 
Indeed, pick an element $g\in\Gamma$ with $g>1$. For every 
$h\in\Gamma$, and every positive integer $n$, there exists 
a largest integer $k(n)$ such that $g^{k(n)}<h^n$. Then
$(k(n)/n~|~n\in\N)$ is a Cauchy sequence and we
let $\rho(h):=\limdir{n\to\infty}k(n)/n$. One verifies 
easily that $\rho$ is an order-preserving group homomorphism,
and since the convex rank of $\Gamma$ equals one, it follows
that $\rho$ is injective. 
\end{example}

\sset\subsubsection{}\label{subsec_bij.convex}
There is an inclusion-reversing bijection between the
set of convex subgroups of the value group $\Gamma$ of 
a valuation $|\cdot|$ and the set of prime ideals of its 
valuation ring $K^+$. This bijection assigns to a convex 
subgroup $\Delta\subset\Gamma$, the prime ideal 
$\fp_\Delta:=\{x\in K^+~|~\gamma>|x|
\text{ for every $\gamma\in\Delta$}\}$.
Conversely, to a prime ideal $\fp$, there corresponds the
convex subgroup $\Delta_\fp:=
\{\gamma\in\Gamma~|~\gamma>|x|\text{ for all $x\in\fp$}\}$.
Then, the value group of the valuation ring $K^+_\fp$ is
(naturally isomorphic to) $\Gamma/\Delta_\fp$.
Furthermore, $K^+/\fp$ is a valuation ring of its field of
fractions, and its value group is $\Delta_\fp$.

\sset\subsubsection{}\label{subsec_rank.of.val}
\index{Valuation(s)!rank of a|indref{subsec_rank.of.val}}
The {\em rank of a valuation\/} is defined as the convex rank
of its value group. It is clear from \eqref{subsec_bij.convex}
that this is the same as the Krull dimension of the
associated valuation ring.

\sset\subsubsection{}\label{subsec_transc.degree}
For any field extension $F_1\subset F_2$, denote by
$\mathrm{tr.d}(F_2:F_1)$ the transcendence degree of
$F_2$ over $F_1$. Let $E$ be a field extension of the valued 
field $K$, and $|\cdot|_E:E^\times\to\Gamma_E$ an extension 
of the valuation $|\cdot|_K:K^\times\to\Gamma_K$ of $K$ to $E$.
Let $\kappa$ (resp. $\kappa(E)$) be the residue field of
the valuation ring of $(K,|\cdot|)$ (resp. of $(E,|\cdot|)$).
Then we have the inequality: 
$\rk(\Gamma_E/\Gamma_K)+\mathrm{tr.d}(\kappa(E):\kappa)
\leq\mathrm{tr.d}(E:K)$ 
(cp. \cite[Ch.VI, \S10, n.3, Cor.1]{BouAC}).

\sset\subsubsection{}\label{subsec_Mnds.in.Gamma}
The image of $K^+\setminus\{0\}$ in $\Gamma$ is the
monoid $\Gamma^+$. The submonoids of $\Gamma^+$ are
in bijective correspondence with the multiplicative
subsets of $K^+\setminus\{0\}$ which contain $(K^+)^\times$.
The bijection is exhibited by the following "short
exact sequence" of monoids:
$$1\to (K^+)^\times\to K^+\setminus\{0\}\stackrel{\pi}{\to}
\Gamma^+\to 1.$$
Then, to a monoid $M\subset\Gamma^+$ one assigns the 
multiplicative subset $\pi^{-1}(M)$. 

\sset\subsubsection{}\label{subsec_satur.mnds}
\index{Monoid(s)|indref{subsec_satur.mnds}}
\index{Monoid(s)!convex submonoid of a|indref{subsec_satur.mnds}}
Let us say that a submonoid $N$ of a monoid $M$ is 
{\em convex\/} if the following holds. If 
$\gamma,\delta\in M$ and $\gamma\cdot\delta\in N$, then 
$\gamma\in N$ and $\delta\in N$. For every submonoid $N$
there is a smallest convex submonoid $N^\mathrm{con}$
such that $N\subset N^\mathrm{con}$.
One deduces a natural bijection between 
convex submonoids of $\Gamma^+$ and prime ideals of 
$K^+$, by assigning on one hand, to a convex monoid $M$, 
the ideal $\fp(M):=K^+\setminus\pi^{-1}(M)$, and on the 
other hand, to a prime ideal $\fp$, the convex monoid 
$M(\fp):=\pi(K^+\setminus\fp)$. 

\sset\subsubsection{}\label{subsec_ideals.in.mnd}
\index{Monoid(s)!ideal, prime ideal in a|indref{subsec_ideals.in.mnd}}
\index{Monoid(s)!$\Spec\,M$ : spectrum of a|indref{subsec_ideals.in.mnd}}
The subsets of the form $M\setminus N$, where $N$
is a convex submonoid of the monoid $M$, are the first
examples of ideals in a monoid. More generally, one
says that a subset $I\subset M$ is an {\em ideal\/}
of $M$, if $I\cdot M\subset I$. Then we say that $I$
is a {\em prime ideal\/} if $I$ is an ideal such that,
for every $x,y\in M$ with $x\cdot y\in I$, we have either
$x\in I$ or $y\in I$. Equivalently, an ideal $I$ is a
prime ideal if and only if $M\setminus I$ is a submonoid;
in this case $M\setminus I$ is necessarily a convex
submonoid. For a monoid $M$, let 
us denote by $\Spec\,M$ the set of all the prime
ideals of $M$. Taking into account 
\eqref{subsec_bij.convex}, we derive bijections
$$\Spec\,\Gamma\stackrel{\sim}{\to}\Spec\,K^+
\stackrel{\sim}{\to}\Spec\,\Gamma^+
\quad:\quad\Delta\mapsto\fp_\Delta\mapsto
\pi(\fp_\Delta)=\Gamma^+\setminus\Delta^+.$$
Furthermore, the bijection $\Spec\,K^+
\stackrel{\sim}{\to}\Spec\,\Gamma^+$ extends to an
inclusion-preserving bijection between the ideals of 
$K^+$ and the ideals of $\Gamma^+$.

In the sequel, it will be sometimes convenient to study
a monoid via the system of its finitely generated 
submonoids. In preparation for this, we want to delve
a little further into the theory of general commutative
monoids.

\begin{definition}\label{def_delve.into}
\index{Monoid(s)!integral|indref{def_delve.into}}
\index{Monoid(s)!free|indref{def_delve.into}}
\index{Monoid(s)!free!basis of a|indref{def_delve.into}}
Let $M$ be a commutative monoid.
\begin{enumerate}
\item
We say that $M$ is {\em integral\/} if we have $a=b$, 
whenever $a,b,c\in M$ and $a\cdot c=b\cdot c$.
\item
We say that $M$ is {\em free\/} if it isomorphic to 
$\N^{(I)}$ for some index set $I$. In this case, a minimal
set of generators for $M$ will be called a {\em basis}.
\end{enumerate}
\end{definition}

\sset\subsubsection{}\label{subsec_cat.monds}
\index{Monoid(s)!$\mathbf{Mnd}$ : category of commutative|indref{subsec_cat.monds}}
\index{Monoid(s)!$M^\mathrm{gp}$ : group associated to an abelian|indref{subsec_cat.monds}}
Let $\mathbf{Mnd}$ be the category of commutative monoids.
The natural forgetful functor $\Z\Mod\to\mathbf{Mnd}$ 
admits a left adjoint functor $M\mapsto M^\mathrm{gp}$.
Given a monoid $M$, the abelian group $M^\mathrm{gp}$ can 
be realized as the set of equivalence classes of pairs 
$(a,b)\in M\times M$, where $(a,b)\sim(a',b')$ if there 
exists $c\in M$ such that $a\cdot b'\cdot c=a'\cdot b\cdot c$;
the addition is defined termwise, and the unit of the 
adjunction is the map $\phi:M\to M^\mathrm{gp}~:~a\mapsto(a,1)$ 
for every $a\in M$. It is easy to see that $\phi$ is injective 
if and only if $M$ is integral.

\sset\subsubsection{}
The category $\mathbf{Mnd}$ admits arbitrary limits and 
colimits. In particular, it admits direct sums. The functor 
$M\mapsto M^\mathrm{gp}$ commutes with limits and colimits.

\begin{theorem}\label{th_tricky.monds} 
Let $\Delta$ be an ordered abelian group, $N\subset\Delta^+$ 
a finitely generated submonoid. Then there exists a free 
finitely generated submonoid $N'\subset\Delta^+$ such that 
$N\subset N'$.
\end{theorem}
\begin{proof} Since $N$ is a submonoid of a group, it is 
integral, so $N\subset N^\mathrm{gp}$. The group homomorphism 
$N^\mathrm{gp}\subset\Delta$ induced by the imbedding
$N\subset\Delta$ is injective as well. The verification 
is straightforward, using the description of $N^\mathrm{gp}$
in \eqref{subsec_cat.monds}. Then $N^\mathrm{gp}$ inherits
a structure of ordered group from $\Delta$, and we
can replace $\Delta$ by $N^\mathrm{gp}$, thereby reducing
to the case where $\Delta$ is finitely generated and
$N$ spans $\Delta$. Thus, in our situation, the convex
rank $r$ of $\Delta$ is finite; we will argue by induction 
on $r$. Suppose then that $r=1$. In this case we will argue
by induction on the rank $n$ of $\Delta$. If $n=1$, then one 
has only to observe that $\Z^+$ is a free monoid. 
Suppose next that $n=2$; in this case, let $g_1,g_2\in\Delta$ 
be a basis. We can suppose that $g_1<g_2<1$; indeed, if
$g_1>1$, we can replace it by $g_1^{-1}$; then, since $r=1$,
we can find an integer $k$ such that $g_2':=g_2\cdot g_1^k<1$
and $g_2'>g_1$; clearly $g_1,g_2'$ is still a basis of $\Delta$. 
We define inductively a sequence of elements 
$g_i\in\Delta^+$, for every $i>2$, in the following way. 
Suppose that $i>2$ and that the elements $g_3<g_4<...<g_{i-1}$ 
have already been assigned; 
let $k_i:=\sup\{n\in\N~|~g_{i-1}\cdot g_{i-2}^{-n}\leq 1\}$;
notice that, since the convex rank of $\Delta$ equals $1$,
we have $k<\infty$. We set $g_i:=g_{i-1}\cdot g_{i-2}^{-k_i}$.

\begin{claim}\label{cl_minimize} We have 
$g_i^\N\cdot g_{i+1}^\N\subset g_{i+1}^\N\cdot g_{i+2}^\N$ 
for every $i\geq 1$, and 
$\Delta^+=\bigcup_{i\geq 1}(g_i^\N\cdot g_{i+1}^\N)$.
\end{claim}
\begin{pfclaim} The first assertion is obvious. We prove
the second assertion. Let $g\in\Delta^+$; for every $i\geq 1$ 
we can write $g=g_i^{a_i}\cdot g_{i+1}^{b_i}$ for unique 
$a_i,b_i\in\Z$. Notice that $a_i$ and $b_i$ cannot both 
be negative. Suppose that either $a_{i+1}$ or $b_{i+1}$ 
is not in $\N$; we show that in this case 
\set\begin{equation}\label{eq_minimize.sum}
|a_{i+1}|+|b_{i+1}|<|a_i|+|b_i|.
\end{equation}
Indeed, we must have either $a_i<0$ and $b_i>0$, or $a_i>0$ 
and $b_i<0$. However, $a_{i+1}=a_i\cdot k_{i+1}+b_i$ and
$b_{i+1}=a_i$; thus, if $a_i<0$, then $b_{i+1}<0$, and
consequently $a_{i+1}>0$, whence 
\set\begin{equation}\label{eq_ineq.a_i.b_i}
|a_{i+1}|<|b_i|
\end{equation}
and if $a_i>0$, then $a_{i+1}<0$, so again 
\eqref{eq_ineq.a_i.b_i} holds. From \eqref{eq_minimize.sum}
it now follows that eventually $a_i$ and $b_i$ become both
positive.
\end{pfclaim}

Since $N$ is finitely generated, claim \ref{cl_minimize} shows
that $N\subset g_i^\N\cdot g_{i+1}^\N$ for $i>0$ sufficiently
large, so the claim follows in this case.

Next, suppose that the convex rank $r=1$ and $n:=\rk(\Delta)>2$.
Write $\Delta=H\oplus G$ for two subgroups such that 
$\rk(H)=n-1$ and $G=g\Z$ for some $g\in\Delta$. 

\begin{claim}\label{cl_dense} For every 
$\delta\in\Delta^+\setminus\{1\}$ we can find $a,b\in H$ 
such that $\delta<a\cdot g^{-1}<1$ and $\delta<b^{-1}\cdot g<1$.
\end{claim}
\begin{pfclaim}
Let $\rho:\Delta\hookrightarrow\R$ be an order-preserving 
imbedding as in example \ref{ex_tens.with.Q}(vi); since 
$\rk(H)>1$, it is easy to see that $\rho(H)$ is dense in 
$\rho(\Delta)$. The claim is an immediate consequence.
\end{pfclaim}

Let $g_1,...,g_k$ be a set of generators for $N$.
For every $i\leq k$ we can write $g_i=h_i\cdot g^{n_i}$
for unique $h_i\in H$ and $n_i\in\Z$. Suppose that
$n_i\geq 0$; it follows easily from claim \ref{cl_dense}
that there exists $a,b\in H$ such that $a<g<b$ and
$g_i<(b^{-1}\cdot g)^{n_i}<1$ (resp. 
$g_i<(a^{-1}\cdot g)^{n_i}<1$) for every $i$ such that
$n_i\geq 0$ (resp. $n_i<0$). Then for $n_i\geq 0$ set
$h_i':=h_i\cdot b^{n_i}$ and for $n_i<0$ set 
$h_i':=h_i\cdot a^{n_i}$. Notice that $h'_i<1$. 
Set $h_0:=a\cdot b^{-1}$ and let $M$ be the submonoid 
of $H^+$ spanned by $h_0,h'_1,...,h'_k,$; we can imbed
$M$ in a larger submonoid $M'\subset H^+$ such that
$(M')^\mathrm{gp}=H$. Then, by inductive assumption, 
we can imbed $M'$ in a free submonoid $L\subset H^+$. 
Let $l_1,...,l_{n-1}$ be a basis for $L$. 
 
We can write $h_0=\prod_{i=1}^tl_{k_i}$ for some integers
$k_1,...,k_t\in\{1,...,n-1\}$. Notice that 
$h_0\leq b^{-1}\cdot g$ and let $r<t$ be the largest integer 
such that $\prod^r_{i=1}l_{k_i}>g\cdot b^{-1}$; set 
$l':=g\cdot b^{-1}\cdot\prod^r_{i=1}l^{-1}_{k_i}$ and
$l'':=g^{-1}\cdot b\cdot\prod^{r+1}_{i=1}l_{k_i}$. 
\begin{claim}
The submonoid $L'$ generated by 
$\{l_1,...,l_{n-1},l',l''\}\setminus\{l_{k_{r+1}}\}$ contains 
$N$. 
\end{claim}
\begin{pfclaim}
Indeed, since $l'\cdot l''=l_{k_{r+1}}$, it follows
that $L\subset L'$; moreover, 
$g\cdot b^{-1}=l'\cdot\prod^r_{i=1}l_{k_i}$ and 
$g^{-1}\cdot a=l''\cdot\prod^t_{i=r+2}l_{k_i}$ so 
$g\cdot b^{-1},g^{-1}\cdot a\in L'$. Now the claim
follows by remarking that $g_i=h'_i\cdot(g^{-1}\cdot a)^{-n_i}$
if $n_i<0$, and $g_i=h'_i\cdot(g\cdot b^{-1})^{n_i}$
if $n_i\geq 0$.
\end{pfclaim}
 
Now, it is clear that $L\subset\Delta^+$; since moreover 
$L$ spans $\Delta$ and is generated by $n$ elements, it follows
that $L$ is a free monoid, so the proof is concluded in
case $\mathrm{c.rk}(\Delta)=1$.

Finally, suppose $r>1$ and pick a convex subgroup 
$0\neq\Delta_0\subsetneq\Delta$; then the ordering on $\Delta$
induces a unique ordering on $\Delta/\Delta_0$ such that
the projection map $\pi:\Delta\to\Delta/\Delta_0$ is 
order-preserving. Let $N_0:=\pi(N)$. By induction, $N_0$ 
can be imbedded into a finitely generated free submonoid 
$F_0$ of $(\Delta/\Delta_0)^+$. By lifting a minimal set
of generators of $F_0$ to elements $f_1,...,f_n\in\Delta^+$, 
we obtain a free finitely generated monoid 
$F\subset\Delta^+$ with $\pi(F)=F_0$. Now, choose a finite
set $S$ of generators for $N$; we can partition 
$S=S_1\cup S_2$, where $S_1=S\cap\Delta_0$ and 
$S_2=S\setminus\Delta_0$. By construction, for every 
$x\in S_2$ there exist integers $k_{i,x}\geq 0$ ($i=1,...,n$)
such that $y_x:=x\cdot\prod_{i=1}^nf_i^{-k_{i,x}}\in\Delta_0$. 
Let $g:=\max\{y_x~|~x\in S_2\}$; if $g<1$, let $e_i:=f_i$,
otherwise let $e_i:=f_i\cdot g$ for every $i\leq n$.
Since $\Delta_0$ is convex, we have in any case: $e_i<1$
for $i\leq n$. Moreover, the elements 
$z_x:=x\cdot\prod_{i=1}^ne_i^{-k_{i,x}}$ are contained in 
$\Delta_0^+$. By induction, the submonoid of $\Delta^+_0$ 
generated by $S_1\cup\{z_x~|~x\in S_2\}$ is contained in a 
free finitely generated monoid $F'\subset\Delta^+_0$. Using 
the convexity of $\Delta_0$ one verifies easily that 
$N':=F\cdot F'$ is a free monoid. Clearly $N\subset N'$, so 
the assertion follows.
\end{proof}

\begin{remark} Another proof of theorem \ref{th_tricky.monds} 
can be found in \cite[Th.2.2]{Ell}. Moreover, this theorem
can also be deduced from the resolution of singularities of 
toric varieties (\cite[Ch.I, Th.11]{Tors}).
\end{remark}

\subsection{Basic ramification theory}\label{sec_basic.ramif}
This section is a review of some basic ramification theory
in the setting of general valuation rings and their algebraic
extensions.

\sset\subsubsection{}\label{subsec_set.notation}
Throughout this section we fix a valued field $(K,|\cdot|)$. 
Its valuation ring will be denoted $K^+$ and the residue
field of $K^+$ will be denoted by $\kappa$. If $(E,|\cdot|_E)$ 
is any valued field extension of $K$, we will denote by $E^+$ 
the valuation ring of $E$, by $\kappa(E)$ its residue field
and by $\Gamma_E$ its value group. 
Furthermore, we let $K^\mathrm{a}$ be an algebraic closure of 
$K$, and $K^\sep$ the separable closure of $K$ contained 
in $K^\mathrm{a}$.

\sset\subsubsection{}\label{subsec_basic-ramify}
Let $E\subset K^\mathrm{a}$ be a finite extension of $K$. 
Let $W$ be the integral closure of $K^+$ in $E$; by remark 
\ref{rem_val.ring.princip}(iv), to every maximal ideal 
$\fp$ of $W$ we can associate a valuation 
$|\cdot|_\fp:E^\times\to\Gamma_\fp$ extending $|\cdot|$,
and (up to isomorphisms of value groups) every extension
of $|\cdot|$ to $E$ is obtained in this way. Set 
$\kappa(\fp):=W/\fp$; it is known that 
$\sum_{\fp\in\mathrm{Max}(W)}[\Gamma_\fp:\Gamma]\cdot
[\kappa(\fp):\kappa]\leq[E:K]$ (cp. 
\cite[Ch.VI \S8, n.3 Th.1]{BouAC}).

\sset\subsubsection{}\label{subsec_Galois-ramify}
\index{Decomposition subgroup|indref{subsec_Galois-ramify}}
\index{Inertia subgroup|indref{subsec_Galois-ramify}}
Suppose now that $E$ is a Galois extension of $K$.
Then $\Gal(E/K)$ acts transitively on $\mathrm{Max}(W)$.
For a given $\fp\in\mathrm{Max}(W)$, the {\em decomposition
subgroup\/} $D_\fp\subset\Gal(E/K)$ of $\fp$ is the stabilizer
of $\fp$. Then $\kappa(\fp)$ is a normal extension of $\kappa$
and the natural morphism 
$D_\fp\to\Aut(\kappa(\fp)/\kappa)$ is surjective; its
kernel $I_\fp$ is the {\em inertia subgroup\/} at $\fp$
(cp. \cite[Ch.V, \S2, n.2, Th.2]{BouAC} for the case
of a finite Galois extension; the general case is obtained
by passage to the limit over the family of finite Galois
extensions of $K$ contained in $E$). 

\sset\subsubsection{}\label{subsec_fin.Gal-ramif}
If now $E$ is a finite Galois extension of $K$, then
it follows easily from \eqref{subsec_basic-ramify} and
\eqref{subsec_Galois-ramify} that the integers 
$[\Gamma_\fp:\Gamma]$ and $[\kappa(\fp):\kappa]$ are
independent of $\fp$, and therefore, if $W$ admits
$n$ maximal ideals, we have : 
$n\cdot[\Gamma_\fp:\Gamma]\cdot[\kappa(\fp):\kappa]\leq[E:K]$.

\begin{lemma}\label{lem_Gamma-unchanged} Let $K^{+\mathrm{sh}}$ 
be a strict henselization of $K^+$; then $K^{+\mathrm{sh}}$ 
is a valuation ring and $\Gamma_{K^{+\mathrm{sh}}}=\Gamma$.
\end{lemma}
\begin{proof} It was shown in remark 
\ref{rem_val.ring.princip}(vi) that $K^{+\mathrm{sh}}$ is a 
valuation ring. To show the second assertion, let $R$ be
more generally any integrally closed local domain; the (strict) 
henselization of $R$ can be constructed as follows (cp. 
\cite[Ch.X, \S2, Th.2]{Ray}). Let $F:=\mathrm{Frac}(R)$, 
$F^\sep$ a separable closure of $F$, $\fp$ any maximal 
ideal of the integral closure $W$ of $R$ in $F^\mathrm{s}$, 
$D$ and $I$ respectively the decomposition and inertia 
subgroups of $\fp$; let $W^D$ (resp. $W^I$) be the subring 
of elements of $W$ fixed by $D$ (resp. by $I$) and set 
$\fp^D:=W^D\cap\fp$, $\fp^I:=W^I\cap\fp$. 
Then the localization $R^\mathrm{h}:=(W^D)_{\fp^D}$ (resp. 
$R^\mathrm{sh}:=(W^I)_{\fp^I}$) is a henselization (resp. 
strict henselization) of $R$. Now, let us make $R:=K^+$,
so $F:=K$ and $F^\sep:=K^\sep$; let $E\subset K^\sep$ be 
any finite Galois extension of $K$; $W_E:=W\cap E$ is the 
integral closure of $K^+$ in $E$; set $D_E:=D\cap\Gal(E/K)$, 
$I_E:=I\cap\Gal(E/K)$, $E':=E^{D_E}$, $E'':=E^{I_E}$.
Let $\fp':=\fp\cap E'$; it then follows from 
\cite[Ch.VI,\S 12, Th.23]{Za-Sa} that $[\Gamma_{\fp'}:\Gamma]=1$.
Clearly the value group $\Gamma_{K^\mathrm{h}}$ of 
$K^{+\mathrm{h}}$ is the filtered union of all such 
$\Gamma_{\fp'}$, so we deduce $\Gamma_{K^\mathrm{h}}=\Gamma$.
Therefore, in order to prove the lemma, we can assume
that $K=K^\mathrm{h}$. In this case $\Gal(E/K)$ 
coincides with the decomposition subgroup of $\fp'$
and $I_E$ is a normal subgroup of $\Gal(K/E)$
such that $[\Gal(K/E):I_E]$ equals the cardinality
$n$ of $\Aut(\kappa(E)/\kappa)$. By the definition of
$I_E$ it follows that the natural map :
$\Aut(\kappa(E)/\kappa)\to\Aut(\kappa(E'')/\kappa)$
is an isomorphism. We derive :
$[E'':K]=n\leq[\kappa(E''):\kappa]$; then
from \eqref{subsec_fin.Gal-ramif} we obtain
$\Gamma_{E''}=\Gamma$ and the claim follows. 
\end{proof}

\sset\subsubsection{}\label{subsec_henselian}
We suppose now that $K^+$ is a {\em henselian\/} valuation ring, 
with value group $\Gamma$.
Then, on any algebraic extension $E\subset K^\mathrm{a}$ of 
$K$, there is a unique valuation $|\cdot|_E$ extending $|\cdot|$,
and thus a unique inertia subgroup, which we denote simply
by $I$. By remark \ref{rem_val.ring.princip}(iv), 
$E^+$ is the integral closure of $K^+$ in $E$. 

\begin{remark}\label{rem_ineq.Gamma}
(i) In the situation of \eqref{subsec_henselian}, the inequality
of \eqref{subsec_basic-ramify} simplifies to : 
$$[\kappa(E):\kappa]\cdot[\Gamma_E:\Gamma_K]\leq[E:K].$$ 

(ii) Sometimes this inequality is actually an equality; 
this is for instance the case when the valuation of $K$ 
is discrete and the extension $K\subset E$ is finite and 
separable (cp. \cite[Ch.VI, \S8, n.5, Cor.1]{BouAC}).

(iii) However, even when the valuation of $K$ is
discrete, it may happen that the inequality (i) is
strict, if $E$ is inseparable over $K$. As an example,
let $\kappa$ be a perfect field of positive characteristic,
and choose a power series $f(T)\in\kappa[[T]]$
which is transcendental over the subfield $\kappa[T]$.
Endow $F:=\mathrm{Frac}(\kappa(T^{1/p},f(T)))$ with the 
$T$-adic valuation, and let $K$ be the henselization 
of $F$. Then the residue field of $K$ is $\kappa$ and 
the valuation of $K$ is discrete. Let $E:=K[f(T)^{1/p}]$. 
Then $[E:K]=p$, $\Gamma_E=\Gamma_K$ and $\kappa(E)=\kappa$.
\end{remark}

\sset\subsubsection{}\label{subsec_tors.mu}
\index{$\bmu(F)$|indref{subsec_Galois-ramify}}
For a field $F$, we denote by $\bmu(F)$ the torsion subgroup 
of $F^\times$. Let $E$ be a finite Galois extension of $K$
(with $K^+$ still henselian). One defines a pairing 
\set\begin{equation}\label{eq_ram.pairing}
I\times(\Gamma_E/\Gamma_K)\to\bmu(\kappa(E))
\end{equation} 
in the following way. For 
$(\sigma,\gamma)\in I\times\Gamma_E$, let $x\in E^\times$
such that $|x|=\gamma$; then let 
$(\sigma,\gamma)\mapsto\sigma(x)/x$ ($\mod~\fm_E$).
One verifies easily that this definition is independent
of the choice of $x$; moreover, if $x\in K^\times$,
then $\sigma$ acts trivially on $x$, so the definition
is seen to depend only on the class of $\gamma$ in 
$\Gamma_E/\Gamma_K$.

\sset\subsubsection{}\label{subsec_separ.closed}
Suppose furthermore that $\kappa$ is separably closed.
Then the inertia subgroup coincides with the Galois
group $\Gal(E/K)$ and moreover 
$\bmu(\kappa(E))=\bmu(\kappa)$. The pairing 
\eqref{eq_ram.pairing} induces a group homomorphism
\set\begin{equation}\label{eq_ram.mapping}
\Gal(E/K)\to\Hom_\Z(\Gamma_E/\Gamma_K,\bmu(\kappa)).
\end{equation}
Let $p:=\chara(\kappa)$. For a group $G$, let us denote 
by $G^{(p)}$ the maximal abelian quotient of $G$ that does 
not contain $p$-torsion.

\begin{proposition}\label{prop_kernel.p-group} 
Under the assumptions of \eqref{subsec_separ.closed}, 
the map \eqref{eq_ram.mapping} is surjective and its 
kernel is a $p$-group.
\end{proposition}
\begin{proof} Let $n:=[E:K]$. Notice that $\bmu(\kappa)$
does not contain $p$-torsion, hence every homomorphism
$\Gamma_E/\Gamma_K\to\bmu(\kappa)$ factors through
$(\Gamma_E/\Gamma_K)^{(p)}$. Let $m$ be the order of 
$(\Gamma_E/\Gamma_K)^{(p)}$. Let us recall the definition 
of the Kummer pairing: one takes the Galois cohomology 
of the exact sequence of $\Gal(K^\sep/K)$-modules
$$1\to\bmu_m\to(K^\sep)^\times
\stackrel{(-)^m}{\longrightarrow}
(K^\sep)^\times\to 1$$
and applies Hilbert $90$, to derive an isomorphism
$K^\times/(K^\times)^m\simeq H^1(\Gal(K^\sep/K),\bmu_m)$.
Now, since $(m,p)=1$ and $\kappa$ is separably
closed, the equation $X^m=1$ admits $m$ distinct
solutions in $\kappa$. Since $K^+$ is henselian,
these solutions lift to roots of $1$ in $K$, {\em i.e.},
$\bmu_m\subset K^\times$, whence
$H^1(\Gal(K^\sep/K),\bmu_m)\simeq
\Hom_{\mathrm{cont}}(\Gal(K^\sep/K),\bmu_m)$. By working out
the identifications, one checks easily that the resulting
group isomorphism 
$$K^\times/(K^\times)^m\simeq
\Hom_{\mathrm{cont}}(\Gal(K^\sep/K),\bmu_m)$$ 
can be described as follows. To a given $a\in K^\times$, 
we assign the group homomorphism
$$\Gal(K^\sep/K)\to\bmu_m\ :\ 
\sigma\mapsto\sigma(a^{1/m})/a^{1/m}.$$
Notice as well that, since $\kappa$ is separably
closed, more generally every equation of the form 
$X^m=u$ admits $m$ distinct solutions in $\kappa$,
provided $u\neq 0$; again by the henselian property
we deduce that every unit of $K^+$ is an $m$-th power
in $K^\times$; therefore 
$K^\times/(K^\times)^m\simeq\Gamma_K/m\Gamma_K$.

Dualizing, we obtain an isomorphism 
$$\Hom_\Z(\Gamma/m\Gamma,\bmu_m)\simeq
\Hom_\Z(\Hom_{\mathrm{cont}}(\Gal(K^\sep/K),\bmu_m),\bmu_m).$$
However, 
$$\Hom_{\mathrm{cont}}(\Gal(K^\sep/K),\bmu_m)=
\colim{H\subset\Gal(K^\sep/K)}
\Hom_\Z(\Gal(K^\sep/K)/H,\bmu_m)$$
where $H$ runs over the cofiltered system of open normal
subgroups of $\Gal(K^\sep/K)$ such that 
$\Gal(K^\sep/K)/H$ is abelian with exponent dividing $m$. 
It follows that 
$$\Hom_\Z(\Hom_{\mathrm{cont}}
(\Gal(K^\sep/K)/H,\bmu_m),\bmu_m)
\simeq\liminv{H\subset
\Gal(K^\sep/K)}\Gal(K^\sep/K)/H$$
where the right-hand side is a quotient of 
$\Gal(K^\sep/K)^{(p)}$. Hence, we have obtained a surjective
group homomorphism
\set\begin{equation}\label{eq_Kummer.morphism}
\Gal(K^\sep/K)\to
\Hom_\Z(m^{-1}\Gamma_K/\Gamma_K,\bmu_m)\stackrel{\sim}{\to}
\Hom_\Z(m^{-1}\Gamma_K/\Gamma_K,\bmu(\kappa)).
\end{equation}
(Since $\Gamma_K$ is torsion-free, we can identify naturally
$\Gamma_K/m\Gamma_K$ to the subgroup
$m^{-1}\Gamma_K/\Gamma_K\subset
(\Gamma_K\otimes_\Z\Q)/\Gamma_K$). Let 
$j:\Gamma_E/\Gamma_K\hookrightarrow m^{-1}\Gamma_K/\Gamma_K$
be the inclusion map. One verifies directly from the 
definitions, that the maps \eqref{eq_ram.mapping} and 
\eqref{eq_Kummer.morphism} fit into a commutative diagram
$$\xymatrix{\Gal(K^\sep/K) \ar[r] \ar[d] &
\Gal(E/K) \ar[d] \\
\Hom_\Z(m^{-1}\Gamma_K/\Gamma_K,\bmu(\kappa)) \ar[r]^-\rho &
\Hom_\Z(\Gamma_E/\Gamma_K,\bmu(\kappa))
}$$
where the top map is the natural surjection, and
$\rho:=\Hom_\Z(j,\bmu(\kappa))$. Finally, an easy
application of Zorn's lemma shows that $\rho$ is
surjective, and therefore, so is \eqref{eq_ram.mapping}.

It remains to show that the kernel $H$ of 
\eqref{eq_ram.mapping} is a $p$-group. Suppose that 
$\sigma\in H$ and nevertheless $p$ does not divide the 
order $l$ of $\sigma$; then we claim that the $K$-linear
map $\phi:E\to E$ given by 
$x\mapsto\sum_{i=0}^{l-1}\sigma^i(x)$ is an isometry. 
Indeed, $\phi(x)=l\cdot x+\sum_{i=1}^{l-1}(\sigma^i(x)-x)$;
it suffices then to remark that $|l\cdot x|=|x|$ and
$|\sigma^i(x)-x|<|x|$, since $\sigma^i\in H$ for 
$i=0,...,l-1$. Next, for every $x\in E$ we can write
$0=\sigma^l(x)-x=\phi(x-\sigma(x))$; hence $\sigma(x)=x$,
that is, $\sigma$ is the neutral element of $\Gal(E/K)$,
as asserted.
\end{proof}

\begin{corollary}\label{cor_kern.p-group} 
Keep the assumptions of \eqref{subsec_separ.closed}, 
and suppose moreover that $(p,[E:K])=1$. 
Then $\Gamma_E/\Gamma_K\simeq\Hom_\Z(\Gal(E/K),\bmu(K))$.
Moreover, if\/ $\Gamma_E/\Gamma_K\simeq
\Z/q_1\Z\oplus...\oplus\Z/q_k\Z$, then there exist 
$a_1,...,a_k\in K$ with 
$E=K[a_1^{1/q_1},...,a_k^{1/q_k}]$. 
\end{corollary}
\begin{proof} To start out, since $(p,[E:K])=1$,
proposition \ref{prop_kernel.p-group} tells us that
the map \eqref{eq_ram.mapping} is an isomorphism.
In particular, $\Gal(E/K)$ is abelian, and 
$[\Gamma_E:\Gamma_K]\geq[E:K]$, whence 
$[\Gamma_E:\Gamma_K]=[E:K]$ by remark 
\ref{rem_ineq.Gamma}(ii).
Therefore $\Gal(E/K)\simeq\Z/q_1\Z\oplus...\oplus\Z/q_k\Z$
and $E$ is a compositum of cyclic extensions $E_1,...,E_k$ 
of order $q_1,...,q_k$. It follows as well that
$\Gamma_E/\Gamma_K\simeq\Hom_\Z(\Gal(E/K),\bmu(\kappa))$,
so the first assertion holds; furthermore the latter
holds also for every extension of $K$ contained in $E$.
We deduce :
\begin{claim}\label{cl_Galois.ordered} 
The Galois correspondence establishes an inclusion 
preserving bijection between the subgroups of $\Gamma_E$ 
containing $\Gamma_K$, and the subfields of $E$ containing 
$K$.
\end{claim}

To prove the second assertion, we are thus reduced to the 
case where $E$ is a cyclic extension of prime power
order, say $\Gal(E/K)\simeq\Z/q\Z$, with $(q,p)=1$.
Let $\gamma\in\Gamma_E$ be an element whose class
in $\Gamma_E/\Gamma_K$ is a generator; we can find
$a\in K$ such that $|a|=\gamma^q$. Let $E':=K[a^{1/m}]$
and $F:=E\cdot E'$. Since $\Gamma_F$ is torsion-free, 
one sees easily that its subgroups $\Gamma_E$ and 
$\Gamma_{E'}$ coincide. However, $F$ satisfies again 
the assumptions of the corollary, therefore claim 
\ref{cl_Galois.ordered} applies to $F$, and yields $E=E'$.
\end{proof}

\begin{definition}\label{def_max.tame.ext}
\index{Valued field(s)!$K^\tame$ : maximal tame extension 
of a|indref{def_max.tame.ext}}
Let $(K,|\cdot|)$ be a valued field.
We denote by $K^{+\mathrm{sh}}$ be the strict 
henselisation of $K^+$ and set
$K^\mathrm{sh}:=\mathrm{Frac}(K^{+\mathrm{sh}})$. 
The {\em maximal tame extension\/} $K^\tame$ of $K$ 
in its separable closure $K^\sep$ is the union of 
all the finite Galois extensions $E$ of $K^\mathrm{sh}$ 
inside $K^\sep$, such that $([E:K^\mathrm{sh}],p)=1$. 
Notice that, by corollary \ref{cor_kern.p-group}, every
such extension is abelian and the compositum
of two such extensions is again of the same type,
so the family of all such finite extension is 
filtered, and therefore their union is their
colimit, so the definition makes sense.
\end{definition}

\sset\subsubsection{}
Since $\Gamma_{K^\mathrm{sh}}=\Gamma_K$, one verifies 
easily from the foregoing that there is a natural
isomorphism of topological groups 
$\Gal(K^\tame/K^\mathrm{sh})\simeq
\Hom_\Z(\Gamma_K\otimes_\Z\Z_{(p)}/\Gamma_K,\bmu^{(p)})$,
where $\bmu$ denotes the group of roots of $1$
in $K^\mathrm{sh}$ and where we endow the target
with the profinite topology.

\sset\subsubsection{}\label{subsec_stability.K.t}
Let $E\subset K^\sep$ be any separable extension of $K$. 
Then it is easy to check that $E^\tame=E\cdot K^\tame$. 
Indeed, one knows that $E^\mathrm{sh}=E\cdot K^\mathrm{sh}$; then
let $F$ be a finite separable extension of $E$ such
that $([F:E],p)=1$. By taking roots of elements of
$K$ we can find an extension $F'$ of $K$ such that
$([F':K],1)=1$ and 
$(\Gamma_F/\Gamma_K)^{(p)}=(\Gamma_{F'}/\Gamma_K)^{(p)}$
and then $E\cdot F'\cdot K^\mathrm{sh}\supset F$.

\subsection{Algebraic extensions}\label{sec_val.rings}
In this section we return to almost rings: we suppose it 
is given a valued field $(K,|\cdot|)$, and then we will 
study exclusively the almost ring theory relative to the
standard setup attached to $K^+$ (see 
\eqref{subsec_standard.setup}).
For an extension $E$ of $K$, we will use the notation
of \eqref{subsec_set.notation}. Furthermore, we will
denote $W_E$ the integral closure of $K^+$ in $E$.

\begin{lemma}\label{lem_Fitting.vals} Let $R$ be
a ring and $0\to M_1\to M_2\to M_3\to 0$
a short exact sequence of finitely generated torsion
$R$-modules, and suppose that the \/ $\Tor$-dimension 
of $M_3$ is $\leq 1$. Then $F_0(M_2)=F_0(M_1)\cdot F_0(M_3)$.
\end{lemma}
\begin{proof}
We can find epimorphisms $\phi_i:R^{n_i}\to M_i$ for 
$i\leq 3$, with $n_2=n_1+n_3$, fitting into a commutative
diagram with exact rows:
$$
\xymatrix{
0 \ar[r] & R^{n_1} \ar[r] \ar[d]_{\phi_1} & 
R^{n_2} \ar[r] \ar[d]^{\phi_2} & 
R^{n_3} \ar[r] \ar[d]^{\phi_3} & 0 \\
0 \ar[r] & M_1 \ar[r] & M_2 \ar[r] & M_3 \ar[r] & 0.
}
$$
Let $N_i:=\Ker\phi_i$ ($i\leq 3$). By snake lemma we have a 
short exact sequence: $0\to N_1\to N_2\stackrel{\pi}{\to}N_3\to 0$. 
Since the Tor-dimension of the $M_3$ is $\leq 1$, it follows
that $N_3$ is a flat $R$-module. 
\begin{claim}\label{cl_flat.lambda.vanish} 
$\Lambda^{n_3+1}_RN_3=0$.
\end{claim}
\begin{pfclaim} Since $N_3$ is flat, the antisymmetrizer
operator $\bar a_k:\Lambda^{k}_RN_3\to N_3^{\otimes k}$
is injective for every $k\geq 0$ (cp. the proof of proposition
\ref{lem_lambda.of.lambda}). On the other hand, 
$\Lambda^{n_3+1}_RR^{n_3+1}=0$, thus it suffices to
show that the natural map 
$j^{\otimes k}:N_3^{\otimes k}\to(R^{n_3})^{\otimes k}$ 
is injective for every $k\geq 0$. This is clear for $k=0$. 
Suppose that injectivity is known for $j^{\otimes k}$; 
we have $j^{\otimes k+1}=(\one_{R^{\otimes k}}\otimes_Rj)
\circ(j^{\otimes k}\otimes_R\one_{N_3})$. Since $N_3$ is flat, 
we conclude by induction on $k$.
\end{pfclaim}

Next recall that, for every $k\geq 0$ there are 
exact sequences 
\set\begin{equation}\label{eq_alternate.quot}
N_1\otimes_R\Lambda^kN_2\to\Lambda^{k+1}N_2
\stackrel{\pi^{\wedge k+1}}{\longrightarrow}\Lambda^{k+1}N_3.
\end{equation}
(To show that such sequences are exact, one uses the
universality of $\Lambda^{k+1}_RN_3$ for $(k+1)$-linear
alternating maps to $R$-modules).
From \eqref{eq_alternate.quot} and claim 
\ref{cl_flat.lambda.vanish}, a simple argument by
induction on $k$ shows that the natural map
$\psi:\Lambda^{n_1}_RN_1\otimes\Lambda_R^{n_3}N_2\to
\Lambda_R^{n_2}N_2$ is surjective. Finally, by
definition, we have 
$F_0(N_i)=
\Img(\Lambda^{n_i}_RN_i\stackrel{j_i^{\wedge n_i}}
{\longrightarrow}
\Lambda_R^{n_i}R^{n_i}\stackrel{\sim}{\to}R)$. 
To conclude, it suffices therefore to remark that
the diagram:
\set\begin{equation}\label{eq_commute.matrix}
{\diagram
\Lambda^{n_1}_RN_1\otimes\Lambda_R^{n_3}N_2 
\ar[d]_-{\one_{\Lambda^{n_1}_RN_1}\otimes\pi^{\wedge k+1}} 
\ar[r]^-\psi &
\Lambda_R^{n_2}N_2 \ar[d]^-{j_2^{\wedge n_2}} \\
\Lambda^{n_1}_RN_1\otimes\Lambda_R^{n_3}N_3 \ar[r] &
\Lambda_R^{n_2}R^{n_2}
\enddiagram}\end{equation}
commutes. We leave to the reader the task of verifying that
the commutativity of \eqref{eq_commute.matrix} boils down 
to a well-known identity for determinants of matrices.
\end{proof}

\begin{remark}\label{rem_verbatim}
(i) Lemma \ref{lem_Fitting.vals} applies especially to 
a short exact sequence of finitely presented torsion 
$K^+$-modules, since by lemma \ref{lem_fin.pres.torsion}, 
any such module has homological dimension $\leq 1$. 

(ii) By the usual density arguments (cp. the proof of
proposition \ref{prop_Fitting.short}), it then follows that 
lemma \ref{lem_Fitting.vals} holds true {\em verbatim}, 
even when we replace $R$ by $K^+$ and the $R$-modules 
$M_1$, $M_2$, $M_3$ by uniformly almost finitely generated 
torsion $K^{+a}$-modules.
\end{remark}

\begin{proposition}\label{prop_fin.sep.uniform} 
Suppose that $K^+$ is a valuation ring of rank one.
Let $E$ be a finite separable extension of $K$. 
Then $W^a_E$ and $\Omega_{W_E^a/K^{+a}}$ are uniformly 
almost finitely generated $K^{+a}$-modules which admit 
the uniform bounds $[E:K]$ and respectively $[E:K]^2$.
Moreover, $W^a_E$ is an almost projective $K^{+a}$-module.
\end{proposition}
\begin{proof} In view of the presentation 
\eqref{eq_present_Omega}, the assertion for 
$\Omega_{W_E^a/K^{+a}}$ is an immediate consequence of the 
assertion for $W_E^a$. The trace pairing 
$t_{E/K}:E\times E\to K$ is perfect since $E$ is separable 
over $K$. Let $e_1,...,e_n$ be a basis of the $K$-vector 
space $E$ and $e_1^*,...,e^*_n$ the dual basis under the 
trace morphism, so that $t_{E/K}(e_i\otimes e^*_j)=\delta_{ij}$ 
for every $i,j\leq n$. We can assume that $e_i\in W_E$ 
and we can find $a\in K^+\setminus\{0\}$ such that 
$a\cdot e_i^*\in W_E$ for every $i\leq n$. Let $w\in W_E$; 
we can write $w=\sum_{i=1}^n a_i\cdot e_i$ for some 
$a_i\in K$. We have $t_{E/K}(w\otimes a\cdot e^*_j)\in K^+$ for 
every $j\leq n$; on the other hand, 
$t_{E/K}(w\otimes a\cdot e^*_j)=a\cdot a_j$. Thus, if 
we let $\phi:K^n\to E$ be the isomorphism 
$(x_1,...,x_n)\mapsto\sum_{i=1}^nx_i\cdot e_i$, we
see that 
\set\begin{equation}\label{eq_bounding.W}
(K^+)^n\subset\phi^{-1}(W_E)\subset a^{-1}\cdot(K^+)^n.
\end{equation}
We can write $W_E$ as the colimit of the family $\cW$
of all its finitely generated $K^+$-submodules containing 
$e_1,...,e_n$; if $W_0\in\cW$, then $W_0$ is a
free $K^+$-module by remark \ref{rem_val.ring.princip}(ii);
then it is clear from \eqref{eq_bounding.W} that the 
rank of $W_0$ must be equal to $n$. The proof follows 
straightforwardly from the following:
\begin{claim}\label{cl_finish.up} Let $\eps\in\fm$; there 
exists $W_0\in\cW$ such that $\eps\cdot W_E\subset W_0$.
\end{claim}
\begin{pfclaim} Indeed, suppose that this is
not the case. Then we can find an infinite sequence of 
finitely generated submodules 
$\oplus^n_{i=1}e_i\cdot K^+\subset W_0\subset W_1\subset 
W_2\subset ...\subset W_E$ 
such that $\eps\cdot W_{i+1}\nsubseteq W_i$ for every
$i\geq 0$. From \eqref{eq_bounding.W} and lemma 
\ref{lem_Fitting.vals} it follows easily that 
$F_0((K^+)^n/a\cdot(K^+)^n)\subset F_0(W_{k+1}/W_0)=
\prod_{i=0}^kF_0(W_{i+1}/W_i)$ for every $k\geq 0$.
However, 
$F_0(W_{i+1}/W_i)\subset\Ann_{K^+}(W_{i+1}/W_i)\subset
\eps\cdot K^+$ for every $i\geq 0$. We deduce that 
$|a|^n<|\eps|^k$ for every $k\geq 0$, which is absurd,
since the valuation of $K$ has rank one.
\end{pfclaim}
\end{proof}

\sset\subsubsection{}
Suppose that the valuation ring of $K$ has rank one.
Let $K\subset E\subset F$ be a tower of finite separable 
extensions. Let $\fp\subset W_E$ be any prime ideal; then
$W_{E,\fp}$ is a valuation ring (see remark 
\ref{rem_val.ring.princip}(iii)), and $W_{F,\fp}$ is the
integral closure of $W_{E,\fp}$ in $F$. It then follows
from proposition \ref{prop_fin.sep.uniform} and remark
\ref{rem_val.ring.princip}(ii) that $W_{F,\fp}^a$ 
is an almost finitely generated projective 
$W_{E,\fp}^a$-module; we deduce that $W_F^a$ is an almost
finitely generated projective $W_E^a$-module, 
therefore we can define the different ideal of $W_F^a$
over $W_E^a$. To ease notation, we will denote it by
$\cD_{W_F/W_E}$. If $|\cdot|_F$ is a valuation of $F$
extending $|\cdot|$, then $F^+=W_{F,\fp}$ for some prime
ideal $\fp\subset W_F$; moreover, if $|\cdot|_E$ is the
restriction of $|\cdot|_F$ to $E$, then $E^+=W_{E,\fq}$,
where $\fq=\fp\cap W_E$. For this reason, we are led to
define $\cD_{F^+/E^+}:=(\cD_{W_F/W_E})_\fp$.

\begin{lemma}\label{lem_diff.in.towers} 
Let $K\subset E\subset F$ be a tower of
finite separable extensions of $K$. Then: 
\begin{enumerate}
\item
The $W_E^a$-module $(W_E^a)^*$ is invertible.
\item
$\cD_{W_E/K^+}\cdot\cD_{W_F/W_E}=\cD_{W_F/K^+}$.
\end{enumerate}
\end{lemma}
\begin{proof} In view of proposition 
\ref{prop_diff.in.towers}, (ii) follows from (i). 
We show (i): from proposition \ref{prop_fin.sep.uniform}
we can find, for every $\eps\in\fm$, a finitely generated 
$K^+$-submodule $M\subset W_E$ such that 
$\eps\cdot W_E\subset M$.
By remark \ref{rem_val.ring.princip}(ii) it follows that
$M$ is a free $K^+$-module, so the same holds for 
$M^*:=\Hom_{K^+}(M,K^+)$. The scalar multiplication 
$M^*\to M^*~:~\phi\mapsto\eps\cdot\phi$
factors through a map $M^*\to W_E^*$, and if we let $N$ be 
the $W_E$-module generated by image of the latter map, then 
$\eps\cdot W_E^*\subset N$. Furthermore, for every prime ideal
$\fp\subset W_E$, the localization $N_\fp$ is a torsion-free
$W_{E,\fp}$-module; since $W_{E,\fp}$ is a valuation ring,
it follows that $N_\fp$ is free of finite rank, again by
remark \ref{rem_val.ring.princip}(ii). Hence, $N$ is a 
projective $W_E$-module. In particular, this shows that 
$(W_E^a)^*$ is almost finitely generated  projective as 
a $W_E^a$-module. To show that $(W_E^a)^*$ is 
also invertible, it will suffice to show that the rank of 
$N$ equals one. However, the rank of $N$ can be computed as 
$\dim_EN\otimes_{W_E}E$. We have 
$N\otimes_{W_E}E=W_E^*\otimes_{K^+}K=\Hom_K(E,K)$, so the 
assertion follows by comparing the dimensions of the two sides.
\end{proof}

\begin{proposition}\label{prop_firstblood} 
Suppose that $K^+$ has rank one.
Let $K\subset E$ be a finite field extension such that 
$l:=[E:K]$ is a prime. Let $p:=\chara(\kappa)$. Suppose 
that either:
\begin{enumerate}
\renewcommand{\labelenumi}{(\alph{enumi})}
\item
$l\neq p$ and $K=K^\mathrm{sh}$, or
\item
$l=p$ and $K=K^\tame$, or
\item
the valuation of $K$ is discrete and henselian, 
and $E$ is separable over $K$, or
\item
the valuation of $K$ is discrete and henselian,
$\Gamma_E=\Gamma_K$ and $\kappa(E)=\kappa$.
\end{enumerate}
Then : 
\begin{enumerate}
\renewcommand{\labelenumi}{(\roman{enumi})}
\item
In case {\em (a)}, {\em (b)} or {\em (d)} holds,
there exists $x\in E\setminus K$ such that $E^+$ is 
the filtered union of a family of finite $K^+$-subalgebras 
of the form $E_i^+:=K^+[a_ix+b_i]$, ($i\in\N$) where 
$a_i,b_i\in K$ are elements with $|a_ix+b_i|\leq 1$.
\item
In case {\em (c)} holds, there exists an element $x\in E^+$ 
such that $E^+=K^+[x]$.
\item
Furthermore, if $E$ is a separable extension of $K$,
then $H_j(\L_{E^+/K^+})=0$ for every $j>0$.
\item
If $E$ is an inseparable extension of $K$, then
$H_j(\L_{E^+/K^+})=0$ for every $j>1$, and moreover 
$H_1(\L_{E^+/K^+})$ is a torsion-free $E^+$-module.
\end{enumerate}
\end{proposition}
\begin{proof} Let us first show how assertions (iii) and 
(iv) follow from (i) and (ii). Indeed, since the cotangent
complex commutes with colimits of algebras, by (i) and (ii)
we reduce to dealing with an algebra of the form $K^+[w]$ for 
$w\in E^+$. Such an algebra is a complete intersection 
$K^+$-algebra, quotient of the free algebra $K^+[X]$ by the 
ideal $I\subset K^+[X]$ generated by the minimal polynomial
$m(X)$ of $w$. In view of \cite[Ch.III, Cor.3.2.7]{Il}, 
one has a natural isomorphism in $\sD(K^+[w]\Mod)$ 
$$\L_{K^+[w]/K^+}\simeq(0\to I/I^2\stackrel{\delta}{\to}
\Omega_{K^+[X]/K^+}\otimes_{K^+}K^+[w]\to 0).$$
If we identify $\Omega_{K^+[X]/K^+}\otimes_{K^+}K^+[w]$ to the rank 
one free $K^+[w]$-module generated by $dX$, then $\delta$ can 
be given explicitly by the rule: 
$f(X)\mapsto f'(w)dX$, for every $f(X)\in(m(X))$. 
However, $E$ is separable over $K$ if and only if
$m'(w)\neq 0$. It follows that $\delta$ is injective 
if and only if $E$ is separable over $K$, which proves 
(iii). If $E$ is inseparable over $K$, then $\delta$
vanishes identically by the same token. This shows (iv).

We prove (ii). Since the valuation is discrete, we must 
have either $e:=[\Gamma_E:\Gamma]=l$ or 
$f:=[\kappa(E):\kappa]=l$ (see remark \ref{rem_ineq.Gamma}(ii)).
If $e=l$, then pick any uniformizer $a\in E$; every element
of $E$ can be written as a sum $\sum_{i=0}^{l-1}x_i\cdot a^i$ 
with $x_i\in K$ for every $i<l$. Then it is easy to see that
such a sum is in $E^+$ if and only $x_i\in K^+$ for every $i<l$.
In other words, $E^+=K^+[a]$. In case $f=l$, we can write 
$\kappa(E)=\kappa[\bar u]$ for some unit $u\in(E^+)^\times$; 
moreover, $\fm_E=\fm E^+$; then $K^+[u]+\fm E^+=E^+$; 
since in this case $E^+$ is a finite $K^+$-module, we deduce 
$E^+=K^+[u]$ by Nakayama's lemma.

We prove (i). Suppose that (a) holds; then by corollary 
\ref{cor_kern.p-group} it follows that 
$\Gamma_E/\Gamma_K\simeq\Z/l\Z$ and $E=K[a^{1/l}]$
for some $a\in K$. Hence:
\set\begin{equation}\label{eq_while.chant}
|a^{i/l}|\notin\Gamma\text{ for every $i=1,...,l-1$}.
\end{equation}
We can suppose that the valuation of $K$ is not discrete,
otherwise we fall back on case (c); then, for every 
$\eps\in\fm$, there exists $b_\eps\in K$ such that 
$|\eps|<|b_\eps^l\cdot a|<1$. Let $x_0,...,x_{l-1}\in K$
and set $w:=\sum_{i=0}^{l-1}x_i\cdot a^{i/l}$. Clearly
every element of $E$ can be written in this form. From 
\eqref{eq_while.chant} we derive that the values 
$|x_i\cdot a^{i/l}|$ such that $x_i\neq 0$ are all distinct. 
Hence, $|w|=\displaystyle{\max_{0\leq i<l}|x_i\cdot a^{i/l}|}$. 
Suppose now that $w\in E^+$; it follows that 
$|x_i\cdot a^{i/l}|\leq 1$ for $i=0,...,l-1$, and in fact
$|x_i\cdot a^{i/l}|<1$ for $i\neq 0$. Let $\eps\in\fm$
such that $|\eps^{l-1}|>|x_i\cdot a^{i/l}|$ for every 
$i\neq 0$. A simple calculation shows that 
$|x_i\cdot b_\eps^{-i}|<1$ for every $i\neq 0$, in other
words, $w\in K^+[b_\eps\cdot a^{1/l}]$, which proves the
claim in this case.

In order to deal with cases (b) and (d) we need some 
preparation. Let $x\in E\setminus K$ be any element, 
and set:
$$\rho(x):=\displaystyle\inf_{a\in K}|x-a|\in\Gamma^\wedge_E.$$

We consider case (b). Notice that the hypothesis
$K=K^\tame$ implies that the valuation of $K$ is not 
discrete. 
For any $y\in E$ we can write $y=f(x)$ for some  
$f(X):=b_0+b_1X+...+b_dX^d\in K[X]$ with $d:=\deg f(X)<p$. 
The degree of the minimal Galois extension $F$ of $K$ 
containing all the roots of $f(X)$ divides $d!$, hence 
$F\subset K^\tame=K$. In other words, we can write 
$y=a_k\cdot\prod^d_{i=1}(x-\alpha_i)$ for some 
$\alpha_1,...,\alpha_d\in K$.

We distinguish two cases: first, suppose that there exists
$a\in K$ with $|x-a|=\rho(x)$. Replacing $x$ by $x-a$
we may achieve that $|x|\leq|x-a|$ for every $a\in K$.
Then the constant sequence $(a_n:=0~|~n\in\N)$ fulfills 
the condition of lemma \ref{lem_approx.sequence}. Thus, if 
$y=f(x)$ as above is in $E^+$, we must have
$|f(X)|_{(0,\rho(x))}\leq 1$; in other terms: 
\set\begin{equation}\label{eq_gauss.bound}
|b_i|\cdot\rho(x)^i\leq 1\quad\text{for every $i\leq d$}.
\end{equation}
Now, if $\rho(x)\in\Gamma_K$, we can find $c\in K$
such that $x_0:=x\cdot c$ still generates $E$ and
$|x_0|=1$, whence $|b_i/c^i|\leq 1$ for every $i\leq 1$; 
however, 
$y=b_0+(b_1/c)\cdot x_0+(b_2/c^2)x_0^2+...+(b_d/c^d)x^d$,
thus $y\in K^+[x_0]$, so in this case, $E^+$ itself is
one of the $E^+_i$.

In case $\rho(x)\notin\Gamma_K$, since anyway $\Gamma_K$ is
of rank one and not discrete, we can find a sequence
of elements $c_1,c_2,...\in K$ such that, letting
$x_i:=x\cdot c_i$, we have
$$
\text{$|x_j-a|\geq|x_j|$\ \ for every $a\in K$, $j\in\N$;\ \ \ \
$|x_j|<1$\ \ \ \ and\ \ \ \ $|x_j|\to 1$.}
$$
\begin{claim}\label{cl_if.this.fails} 
If $x\notin\Gamma_K$, then $|x^l|\notin\Gamma_K$ for every 
$0<l<p$.
\end{claim}
\begin{pfclaim}
Indeed, suppose that $|x^l|\in\Gamma_K$ for some $0<l<p$; 
since $\Gamma_K$ is $l$-divisible, we can multiply $x$ by 
some $a\in K$ to have $|x^l|=1$, therefore $|x|=1$, a
contradiction.
\end{pfclaim}

From \eqref{eq_gauss.bound} and claim \ref{cl_if.this.fails} 
we deduce that actually $|b_i|\cdot\rho(x)^i<1$ whenever $i>0$.
It follows that, for $j$ sufficiently large, we will 
have $1>|x_j^i|>|b_i|\cdot\rho(x)^i$ for every $i>0$.
Writing 
$y=b_0+(b_1/c_j)x_j+(b_2/c_j^2)x_j^2+...+(b_d/c_j^d)x_j^d$
we deduce $y\in K^+[x_j]$, therefore the sequence of 
$K^+$-subalgebra $K^+[c_i\cdot x]$ will do in this case.

Finally we have to consider the case where the infimum
$\rho(x)$ is not attained for any $x\in E$. In this case, 
since the valuation is not discrete and of rank $1$, we 
can find, for every $x\in E$, a sequence of elements 
$a_0,a_1,a_2,...\in K$ such that
\set\begin{equation}\label{eq_minimizer}
\gamma_j:=|x-a_j|\to\rho(x).
\end{equation}
In particular, for $j$ sufficiently large we will have 
$|x|>|x-a_j|$, therefore $|x|=|a_j|$. This shows:
\set\begin{equation}\label{eq_Gamma.doesnt.grow}
\Gamma_E=\Gamma_K.
\end{equation}
Now, pick $x\in E\setminus K$ and any sequence
of elements $a_i\in K$ such that \eqref{eq_minimizer}
holds; it is clear that $(a_i~|~i\in\N)$ fulfills
the condition of lemma \ref{lem_approx.sequence}.
Consequently 
\set\begin{equation}\label{eq_lunacy}
|y|=|f(X)|_{(a_j,\gamma_j)}\quad\text{for every 
sufficiently large $j$.}
\end{equation}
Let $f(X)=b_{0,j}+b_{1,j}(X-a_j)+...+b_{d,j}(X-a_j)^d$.
\eqref{eq_lunacy} says that 
$|b_{i,j}|\cdot\gamma^i_j\leq 1$ whenever $j$ is sufficiently
large. However, from \eqref{eq_Gamma.doesnt.grow} we know 
that $\gamma_j\in\Gamma_K$. Pick $c_j\in K$ such that
$|c_j|=\gamma_j^{-1}$ and set $x_j:=c_j(x-a_j)$. It follows 
that $|b_{i,j}/c^i_j|\leq 1$ and 
$y=b_{0,j}+(b_{1,j}/c_{1,j})x_j+...+(b_{d,j}/c_j^d)x_j^d$.
Hence $y\in K^+[x_j]$. It is then easy to verify
that the family of all such $K^+$-subalgebras is filtered
by inclusion, and thus conclude the proof of case (b).

At last, we turn to case (d). Notice that, by remark 
\ref{rem_ineq.Gamma}(ii), this case can occur only if
$E$ is inseparable over $K$, and then $l=p$. Let 
$x\in E\setminus K$; let $a\in\fm$ be a uniformizer; 
for given $n\in\N$, suppose that $b_n\in K$ has been 
found such that $|x-b_n|\leq|a^n|$. 
Since $\kappa(E)=\kappa$, we can find an element 
$c\in K^+$ such that $c\equiv(x-b_n)/a^n\pmod{\fm}$. Set 
$b_{n+1}:=b_n+c\cdot a_n$; then $|x-b_{n+1}|\leq|a^{n+1}|$.
This shows that $\rho(x)=0$, and the resulting sequence
$(b_n~|~n\in\N)$ converges to $x$ in the $\fm$-adic
topology. Let $y\in E$; we can write $y=f(x)$ for a
polynomial $f(X)\in K[X]$ of degree $d<p$. Let $F$ be 
the minimal field extension of $K$ that contains all
the roots of $f(X)$. Notice that $[F:K]$ divides $d!$, 
hence $F$ is separable over $K$, and $[E\cdot F:F]=p$.
Let $f(X)=c\cdot\prod^d_{i=0}(X-\alpha_i)$ be the
factorization of $f(X)$ in $F[X]$. By lemma
\ref{lem_approx.sequence} we deduce that, for every
sufficiently large $n\in\N$ we have: 
$|y|=|f(X)|_{(b_n,|x-b_n|)}$, where 
$|\cdot|_{(b_n,|x-b_n|)}$ is the Gauss valuation on
$F(X)$. One then argues as in the proof of case (b),
to show that $y\in E_n^+:=K^+[c_n(x-b_n)]$, with $c_n\in K$
such that $|c_n(x-b_n)|=1$. Again, it is easy to verify
that $E_i^+\subset E_{i+1}^+$ for every $i\in\N$,
so the proof is complete.
\end{proof}

\begin{corollary}\label{cor_firstblood} Let $E$ be a finite
field extension of $K$ of prime degree $l$.
\begin{enumerate}
\item 
If $E$ satisfies condition {\em (a)} of proposition 
{\em\ref{prop_firstblood}}, and the valuation of $K$
is not discrete (but still of rank one), then $\Omega_{E^+/K^+}=0$, 
$\L_{E^+/K^+}\simeq 0$ and $\cD_{E^+/K^+}=E^{+a}$.
\item
If $E$ satisfies condition {\em (c)} of proposition 
{\em\ref{prop_firstblood}}, then we have : 
$F_0(\Omega_{E^+/K^+})=\cD_{E^+/K^+}$ and $H_i(\L_{E^+/K^+})=0$ 
for $i>0$.
\end{enumerate}
\end{corollary}
\begin{proof} (i): Since condition (a) holds, proposition 
\ref{prop_firstblood} and its proof show that there exists
$a\in K$ such that $E^+$ is the increasing union of all 
$K^+$-subalgebras of the form $E_b^+:=K^+[b\cdot a^{1/l}]$, 
where $b\in K^+$ ranges over all elements such that 
$|b^l\cdot a|<1$. Consequently, 
$\Omega_{E^+/K^+}=\colim{b}\Omega_{E_b^+/K^+}$, and 
$\L_{E^+/K^+}=\colim{b}\L_{E_b^+/K^+}$. Then, again from
proposition \ref{prop_firstblood} it follows that 
$H_j(\L_{E^+/K^+})=0$ for every $j>0$.
Hence, in order to show the first two assertions, it suffices 
to show that the filtered system of the $\Omega_{E^+_b/K^+}$ 
is essentially zero. However, the $E_b^+$-module 
$\Omega_{E_b^+/K^+}$ is generated 
by $\omega_b:=d(b\cdot a^{1/l})$, and clearly 
$l\cdot(b^l\cdot a)^{(l-1)/l}\cdot\omega_b=0$. Since
$(l,p)=1$, it follows that 
$(b^l\cdot a)^{(l-1)/l}\cdot\omega_b=0$. On the other 
hand, for $|b|<|c|$ we can write: 
$\omega_b=b\cdot c^{-1}\cdot\omega_c$. Therefore, the
image of $\omega_b$ in $\Omega_{E_c^+/K^+}$ vanishes, 
whenever $|b\cdot c^{-1}|<|c^l\cdot a|^{(l-1)/l}$, 
{\em i.e.\/}, whenever $|b\cdot a^{1/l}|<|c^l\cdot a|<1$. 
Since the valuation of $K$ is of rank one and not discrete, 
such a $c$ can always be found. To show the last 
stated equality, let us recall the following general fact
(for whose proof we refer to \cite[Ch.VII, \S 1]{Ray}).

\begin{claim}\label{cl_dual.basis}
Suppose that $E=K[w]$ for some $w\in E$, and let 
$f(X)\in K[X]$ be its minimal polynomial; the elements 
$1,w,w^2,...,w^{l-1}$ form a basis of the $K$-vector space 
$E$. Let $e_1^*,...,e_n^*$ be the corresponding dual basis
under the trace pairing; then the bases 
$S:=\{e_1^*,...,e_n^*\}$ and $S':=\{w^{l-1}/f'(w),
w^{l-2}/f'(w),...,1/f'(w)\}$ span the same $E^+$-submodule
of $E$.
\end{claim}

Let us take $w=b\cdot a^{1/l}$ for some $b\in V$
such that $|b^l\cdot a|<1$. It follows from claim 
\ref{cl_dual.basis} that 
$(\cD_{E^+/K^+})^{-1}\subset f'(w)^{-1}\cdot E^{+a}$, whence
$f'(w)\in\cD_{E^+/K^+*}$. However, $f'(w)=l\cdot w^{l-1}$,
and from the definition of $w$ we see that 
$|f'(w)|$ can be made arbitrarily close to $1$, by
choosing $|b|$ closer and closer to $|a|^{1/l}$.

(ii): the claim about the cotangent complex is just
a restatement of proposition \ref{prop_firstblood}(iii),(iv).
By proposition \ref{prop_firstblood}(ii) we can
write $V_E=K^+[w]$ for some $w\in E^+$. Let $f(x)\in K^+[X]$
be the minimal polynomial of $w$.
Claim \ref{cl_dual.basis} implies that 
$\cD_{E^+/K^+}=(f'(w))$; a standard calculation yields 
$\Omega_{E^+/K^+}\simeq E^+/(f'(w))$, so the assertion holds.
\end{proof}

\begin{theorem}\label{th_big.deal} 
Let $(E,|\cdot|_E)$ be a finite separable valued field 
extension of $(K,|\cdot|)$ and suppose that $K^+$ has rank 
one. Then $F_0(\Omega_{E^{+a}/K^{+a}})=\cD_{E^+/K^+}$ and 
$H_i(\L_{E^+/K^+})=0$ for $i>0$.
\end{theorem}
\begin{proof} We begin with a few reductions:
\begin{claim}\label{cl_transit.gloria}  
We can assume that $E$ is a Galois extension of $K$.
\end{claim}
\begin{pfclaim} Indeed, let $(L,|\cdot|_L)$ be a Galois 
valued field extension of $K$ extending $(E,|\cdot|_E)$. 
We obtain by transitivity 
(\cite[II.2.1.2]{Il}) a distinguished triangle
\set\begin{equation}\label{eq_transit.gloria}
\sigma^{-1}\L_{L^+/E^+}\to\L_{E^+/K^+}\otimes_{E^+}L^+
\to\L_{L^+/K^+}\to\L_{L^+/E^+}.
\end{equation}
Suppose that the theorem is already known for the Galois 
extensions $K\subset L$ and $E\subset L$. Then
\eqref{eq_transit.gloria} implies that $H_i(\L_{E^+/K^+})=0$
for $i>0$ and moreover provides a short exact sequence
$$0\to\Omega_{E^/K^+}\otimes_{E^+}L^+\to\Omega_{L^+/K^+}
\to\Omega_{L^+/E^+}\to 0.$$
However, on one hand, by lemma \ref{lem_diff.in.towers}(ii) 
the different is multiplicative in towers of 
extensions, and the other hand, the Fitting ideal $F_0$ 
is multiplicative for short exact sequences, by virtue of 
remark \eqref{rem_verbatim}(ii), so the claim follows.
\end{pfclaim}

\begin{claim}\label{cl_reduce.to.hensel} 
We can assume that $K^+$ is strictly henselian.
\end{claim}
\begin{pfclaim} Indeed, let $K^{+\mathrm{sh}}$ be the strict
henselisation of $K^+$ and 
$K^\mathrm{sh}:=\mathrm{Frac}(K^{+\mathrm{sh}})$.
It is known that $K^{+\mathrm{sh}}$ is an ind-{\'e}tale extension 
of $K^+$, therefore $E^+\otimes_{K^+}K^{+\mathrm{sh}}$ is a reduced 
normal semilocal integral and flat $K^{+\mathrm{sh}}$-algebra, whence 
a product of reduced normal local integral and flat 
$K^{+\mathrm{sh}}$-algebras $W_1,...,W_k$. Each such $W_i$ 
is necessarily the integral closure of $K^{+\mathrm{sh}}$
in $E_i:=\mathrm{Frac}(W_i)$. It follows that 
$\L_{E^+/K^+}\otimes_{K^+}K^{+\mathrm{sh}}\simeq
\L_{E^+\otimes_{K^+}K^{+\mathrm{sh}}/K^{+\mathrm{sh}}}\simeq
\oplus^k_{i=1}\L_{W_i\otimes_{K^+}K^{+\mathrm{sh}}}$. 
Furthermore: $\cD_{E^+/K^+}\otimes_{K^{+a}}(K^{+\mathrm{sh}})^a
\simeq\oplus^k_{i=1}\cD_{E_i^+/K^{+\mathrm{sh}}}$ and similarly
for the modules of differentials. We remark as well
that the formation of Fitting ideals commutes with 
arbitrary base changes. In conclusion, it is clear
that the assertions of the theorem hold for the
extension $K\subset E$ if and only if they hold for
each extension $K^\mathrm{sh}\subset E_i$.
\end{pfclaim}

\begin{claim}\label{cl_reduce.to.p.grp} 
Suppose $K=K^\mathrm{sh}$. We can assume that $\Gal(E/K)$ 
is a $p$-group.
\end{claim}
\begin{pfclaim} Indeed, let $P$ be the kernel of
\eqref{eq_ram.mapping}. By proposition 
\ref{prop_kernel.p-group}, $P$ is $p$-group; let $L$
be the fixed field of $P$. Then $L\subset K^\tame$ and,
by virtue of corollary \ref{cor_kern.p-group}, we
see that $L$ admits a chain of subextensions 
$K:=L_0\subset L_1\subset...\subset L_k:=L$ such that each 
$L_i\subset L_{i+1}$ satisfies either condition (a) or (c) 
of proposition \ref{prop_firstblood}. Then, by corollary 
\ref{cor_firstblood} it follows that the assertions
of the theorem are already known for the extensions
$L_i\subset L_{i+1}$. From here, using transitivity
of the cotangent complex and multiplicativity of
the different in towers of extensions, and of the
Fitting ideals for short exact sequences, one shows 
that the assertions hold also for the extension 
$K\subset L$ (cp. the proof of claim \ref{cl_transit.gloria}).
Now, if the assertions are known to hold as well for 
the extension $L\subset E$, again the same argument
proves them for $K\subset E$.
\end{pfclaim}

\begin{claim}\label{cl_discrete.is.ok}
The theorem holds if the valuation of $K$ is discrete.
\end{claim}
\begin{pfclaim} By claim \ref{cl_reduce.to.p.grp}, we 
can suppose that $\Gal(E/K)$ is a $p$-group. Hence, we 
can find a sequence of subextensions 
$E_0:=K\subset E_1\subset E_2\subset...\subset E_n:=E$
with $[E_{i+1}:E_i]=p$, for every $i=0,...,n-1$.
Arguing like in the proof of claim \ref{cl_reduce.to.p.grp}
we see that it suffices to prove the claim for
each of the extensions $E_i\subset E_{i+1}$.
In this case we are left to dealing with an extension 
$K\subset E$ of degree $p$, which is taken care of
by corollary \ref{cor_firstblood}(ii).
\end{pfclaim}

\begin{claim}\label{cl_at.most} 
Suppose $K=K^\mathrm{sh}$, that $\Gal(E:K)$ is a $p$-group 
and that the valuation of $K$ is not discrete. Let $L$ be 
a finite Galois extension of $K$ such that $([L:K],p)=1$. 
Then the natural map $E^+\otimes_{K^+}L^+\to (E\cdot L)^+$ 
is an isomorphism.
\end{claim}
\begin{pfclaim} By corollary \ref{cor_kern.p-group}
we know that $L$ admits a tower of subextensions
of the form $K:=L_0\subset L_1\subset...\subset L_k:=L$,
such that, for each $i\leq k$ we have 
$L_{i+1}=L_i[a^{1/l}]$ for some $a\in L_i$ and some
prime $l\neq p$. By induction on $i$, we can
then reduce to the case where $L=K[a^{1/l}]$
for some $a\in K$ and a prime $l\neq p$. Under the above
assumptions, we must have $E\cap L=K$, hence $a\notin E$. 
Then $E\cdot L=E[a^{1/l}]$ and by proposition 
\ref{prop_firstblood} and its proof, $(E\cdot L)^+$ is the 
filtered union of all its subalgebras of the form 
$E^+[b\cdot a^{1/l}]$, where $b\in E$ ranges over all 
the elements such that $|b^m\cdot a|<1$.
However, since the valuation of $K$ is not discrete
and has rank one, $\Gamma_K$ is dense in $\Gamma_E$, 
and consequently the subfamily consisting of the 
$E^+[b\cdot a^{1/l}]$ with $b\in K$ is cofinal. Finally, 
for $b\in K$ we have 
$E^+[b\cdot a^{1/l}]\simeq E^+\otimes_{K^+}K^+[b\cdot a^{1/l}]$.
By taking colimits, it follows that 
$(E\cdot L)^+\simeq E^+\otimes_{K^+}L^+$.
\end{pfclaim}

\begin{claim} We can assume that $K$ is equal to $K^\tame$.
\end{claim}
\begin{pfclaim} By claim \ref{cl_reduce.to.hensel} we can 
and do assume that $K=K^\mathrm{sh}$, in which case $K^\tame$ 
is the filtered union of all the finite Galois extension $L$ 
of $K$ such that $([L:K],p)=1$. Then 
$K^{\tame+}=\bigcup_LL^+$ and 
$(E\cdot K^\tame)^+=\bigcup_L(E\cdot L)^+$, where $L$ ranges 
over all such extensions. By claim \ref{cl_reduce.to.p.grp} 
we can also assume that $\Gal(E/K)$ is a $p$-group, in which
case, by claim \ref{cl_at.most}, we have 
$E^+\otimes_{K^+}L^+\stackrel{\sim}{\to}(E\cdot L)^+$ for 
every $L$ as above.
Taking colimit, we get $E^+\otimes_{K^+}K^{\tame+}
\stackrel{\sim}{\to}(E\cdot K^\tame)^+$.
Since $K^{\tame+}$ is faithfully flat over $K^+$, this shows
that, in order to prove the theorem, we can replace 
$K$ by $K^\tame$; however, by \eqref{subsec_stability.K.t} 
we have $(K^\tame)^\tame=K^\tame$, whence the claim.
\end{pfclaim}

After this preparation, we are ready to finish the
proof of the theorem. We are reduced to considering
a Galois extension $E$ of $K=K^\tame$ such that $\Gal(E/K)$ 
is a $p$-group; moreover, we can assume that the valuation 
of $K$ is not discrete. Then, arguing as in the proof of claim 
\ref{cl_discrete.is.ok}, we can further reduce to dealing 
with an extension $K\subset E$ of degree $p$; furthermore, 
the condition $K=K^\tame$ still holds, by virtue of 
\eqref{subsec_stability.K.t}. In this situation, 
condition (b) of proposition \ref{prop_firstblood}
is fulfilled, hence $H_j(\L_{E^+/K^+})=0$ for $j>0$, by 
proposition \ref{prop_firstblood}(iii). It remains to show 
the identity $F_0(\Omega_{E^{a+}/K^{+a}})=\cD_{E^+/K^+}$.
By proposition \ref{prop_firstblood}(i), there exists
$x\in E$ such that $E^+$ is the filtered union of a
family of finite $K^+$-subalgebras $E_i^+:=K^+[a_ix+b_i]$ 
($i\in\N$) of $E^+$. Let $f(X)\in K^+[X]$ be the minimal 
polynomial of $x$. By construction of $E_i^+$, it
is clear that they form a Cauchy net in $\cI_{K^{+a}}(E^{+a})$ 
converging to $E^{+a}$. It then follows from lemma 
\ref{lem_Omega.converge}, that the net 
$\{\Omega_{E_i^+/K^+}\otimes_{E_i^+}E^+~|~i\in\N\}$ 
converges to $\Omega_{E^+/K^+}$ in $\cM(E^+)$.
In particular, 
$F_0(\Omega_{E^+/K^+})=\liminv{i\to\infty}
F_0(\Omega_{E_i^+/K^+}\otimes_{E_i^+}E^+)$.
The minimal polynomial of $a_ix_i+b_i$ is 
$f_i(X):=f(a_i^{-1}X-b_i)$,
therefore: $\Omega_{E_i^+/K^+}=E_i^+/(f'_i(a_ix_i+b_i))=
E_i^+/(a_i^{-1}f'(x))$. Consequently, 
$F_0(\Omega_{E^+/K^+})=\liminv{i\to\infty}(a_i^{-1}f'(x))$.
On the other hand, claim \ref{cl_dual.basis} yields:
$\cD_{E_i^+/K^+}=(a_i^{-1}f'(x))$ for every $i\in\N$.
Then the claim follows from lemma \ref{lem_net.of.algs}.
\end{proof}

The final theorem of this section completes and extends
theorem \ref{th_big.deal} to include valuations of
arbitrary rank. 

\begin{theorem}\label{th_gen.alg.case} 
Let $(K,|\cdot|)$ be any valued field and $(E,|\cdot|_E)$ 
any algebraic valued field extension of $(K,|\cdot|)$. 
We have :
\begin{enumerate}
\item
$H_i(\L_{E^+/K^+})=0$ for $i>1$ and $H_1(\L_{E^+/K^+})$ 
is a torsion-free $E^+$-module.
\item
If moreover, $E$ is a separable extension of $K$, then
$H_i(\L_{E^+/K^+})=0$ for $i>0$.
\end{enumerate}
\end{theorem}
\begin{proof} 
Let us show first how to deduce (ii) from (i). Indeed,
suppose that $E$ is separable over $K$. Then 
$\L_{E/K}\simeq 0$. However, by (i), the natural map
$H_1(\L_{E^+/K^+})\to H_1(\L_{E^+/K^+})\otimes_{K^+}K
\simeq H_1(\L_{E/K})$ is injective, so
the assertion follows.

In order to prove (i), we reduce easily to the case of a 
finite algebraic extension. Let us write $K$ as the filtered
union of its subfields $L_\alpha$ that are finitely
generated over the prime field. For each such $L_\alpha$,
let $K_\alpha:=(L_\alpha)^\mathrm{a}\cap K$ and 
$E_\alpha:=(L_\alpha)^\mathrm{a}\cap E$. Then $E_\alpha$
is a finite extension of $K_\alpha$ and $K$ is the
filtered union of the $K_\alpha$. It follows easily
that we can replace the extension $K\subset E$ by
the extension $K_\alpha\subset E_\alpha$, thereby reducing
to the case where the transcendence degree of $K$ over its 
prime field is finite. In this situation, the rank $r$ of $K$ 
is finite (cp. \eqref{subsec_transc.degree}).
We argue by induction on $r$. Suppose first that $r=1$.
We can split into a tower of extensions 
$K\subset K^\sep\cap E\subset E$; then, by using
transitivity (cp. the proof of claim \ref{cl_transit.gloria}), 
we reduce easily to prove the assertion for the subextensions 
$K\subset K^\sep$ and $K^\sep\cap E\subset E$. 
However, the first case is already covered by theorem 
\ref{th_big.deal}, so we can assume that $E$ is purely 
inseparable over $K$. In this case, we can further split 
$E$ into a tower of subextensions of degree equal to $p$; 
thus we reduce to the case where $[E:K]=p$. We apply transitivity
to the tower $K\subset K^\tame\subset E\cdot K^\tame=E^\tame$: 
by proposition \ref{prop_firstblood}(iv) we know that 
$H_i(\L_{E^{\tame+}/K^{\tame+}})$ vanishes for $i>1$ and 
is torsion-free for $i=1$; by theorem \ref{th_big.deal}, 
we have $H_i(\L_{K^{\tame+}/K^+})=0$ for $i>0$, therefore
$H_i(\L_{E^{\tame+}/K^+})$ vanishes for $i>1$ and is torsion-free
for $i=1$. Next we apply transitivity to the tower
$K\subset E\subset E^\tame$ : by theorem \ref{th_big.deal}
we have $H_i(\L_{E^{\tame+}/E^+})=0$ for $i>0$, and the claim
follows easily.

Next suppose that $r>1$, and that the theorem is already
known for ranks $<r$. Arguing as in the proof of claim
\ref{cl_reduce.to.hensel}, we can even reduce to the case 
where $K^+$ is henselian, and then $E^+$ is the integral
closure of $K^+$ in $E$.
Let $\fp_r:=(0)\subset\fp_{r-1}\subset...\subset\fp_0$
be the chain of prime ideals of $K^+$, and for every $i\leq r$
let $\fq_i$ be the unique prime ideal of $E^+$ lying over
$\fp_i$. The valuation ring $E^+_{\fq_1}$ has rank $r-1$,
thus, by inductive assumption, the desired assertions
are known for the extension $K^+_{\fp_1}\subset E^+_{\fq_1}$.
It suffices therefore to show that 
$H_i(\L_{E^+/K^+})\subset H_i(\L_{E^+_{\fq_1}/K^+_{\fp_1}})$ 
for every $i\geq 0$. Pick $a\in\fp_0\setminus\fp_1$.
Then $K^+_{\fp_1}=K^+[a^{-1}]$ and $E^+_{\fq_1}=E^+[a^{-1}]$
and $\L_{E^+_{\fq_1}/K^+_{\fp_1}}=
\L_{E^+/K^+}\otimes_{K^+}K^+[a^{-1}]$. 
Hence, we are reduced to show that multiplication by $a$ 
is injective on the homology of $\L_{E^+/K^+}$. Let 
$R:=K^+/aK^+$
and $R_E:=E^+\otimes_{K^+}R$. We have a short exact sequence 
$0\to K^+\stackrel{a}{\to} K^+\to R\to 0$, therefore,
after tensoring by $\L_{E^+/K^+}$, a distinguished triangle:
$$
\L_{E^+/K^+}\stackrel{a}{\to}\L_{E^+/K^+}\to
\L_{E^+/K^+}\derotimes_{K^+}R\to\sigma\L_{E^+/K^+}.
$$
On the other hand, according to remark
\ref{rem_val.ring.princip}(ii), $E^+$ is flat over $K^+$, 
therefore $\L_{E^+/K^+}\derotimes_{K^+}R\simeq\L_{R_E/R}$ 
(by \cite[II.2.2.1]{Il}). Consequently, it suffices to
show that $H_i(\L_{R_E/R})=0$ for 
$i\geq 2$. However, $R=(K^+/\fp_1)\otimes_{K^+}R$, and 
$R_E=(E^+/\fq_1)\otimes_{K^+}R$; moreover, $E^+/\fq_1$
is the integral closure of the valuation ring $K^+/\fp_1$ 
in the finite field extension $\mathrm{Frac}(E^+/\fq_1)$ of 
$\mathrm{Frac}(K^+/\fp_1)$. Therefore we can replace $K^+$ 
by $K^+/\fp_1$ and $E^+$ by $K^+/\fq_1$. This turns us back
to the case where $r=1$. Then the vanishing of 
$H_i(\L_{E^+/K^+})$ for $i\geq 2$ yields the vanishing
of $H_i(\L_{R_E/R})$ for $i>2$. Moreover, since 
$H_1(\L_{E^+/K^+})$ is torsion-free, multiplication by $a$
on $H_1(\L_{E^+/K^+})$ is injective, therefore 
$H_2(\L_{R_E/R})$ vanishes as well.
\end{proof}

\subsection{Logarithmic differentials}\label{sec_log.stuff}
In this section $K^+$ is a valuation ring of arbitrary rank.
We keep the notation of \eqref{subsec_set.notation}.
We start by reviewing some facts on logarithmic structures,
for which the general reference is \cite{Ka}.

\sset\subsubsection{}\label{subsec_asso.sh.groups}
\index{$\mathbf{Mnd}_X$, $\Z\Mod_X$|indref{subsec_asso.sh.groups}}
Let $\mathbf{Mnd}_X$ (reps. $\Z\Mod_X$) be the category 
of sheaves of commutative monoids (resp. of abelian groups) 
on a topological space $X$. The forgetful functor 
$\Z\Mod_X\to\mathbf{Mnd}_X$ admits a left adjoint
functor $\underline M\mapsto\underline M^\mathrm{gp}$. 
If $\underline M$ is a sheaf of monoids,
$\underline M^\mathrm{gp}$ is the sheaf associated to the
presheaf defined by : $U\mapsto\underline M(U)^\mathrm{gp}$
for every open subset $U\subset X$. 

The functor $\Gamma:\mathbf{Mnd}_X\to\mathbf{Mnd}$ that 
associates to every sheaf of monoids its global sections, 
admits a left adjoint 
$\mathbf{Mnd}\to\mathbf{Mnd}_X~:~M\mapsto M_X$.
For a monoid $M$, $M_X$ is the sheaf associated to the
constant presheaf with value $M$.

\sset\subsubsection{}\label{subsec_recall.prelog}
\index{Pre-log structure(s)|indref{subsec_recall.prelog}}
\index{Pre-log structure(s)!$\prelog_X$ : category of|indref{subsec_recall.prelog}}
Recall that a {\em pre-log structure\/} on a scheme $X$ 
is a morphism of sheaves of commutative monoids : 
$\alpha:\underline M\to\cO_X$, where the monoid structure 
of $\cO_X$ is induced by multiplication of local sections.
We denote by $\prelog_X$ the category of pre-log
structures on $X$. 

To a monoid $M$ and a morphism of monoids 
$\phi:M\to\Gamma(X,\cO_X)$, one can associate a pre-log
structure $\phi_X:M_X\to\cO_X$ by composing the
induced morphism of constant sheaves 
$M_X\to\Gamma(X,\cO_X)_X$ with the counit of the adjunction 
$\Gamma(X,\cO_X)_X\to\cO_X$.

\sset\subsubsection{}\label{eq_prelog.and.maps}
To a morphism $\phi:Y\to X$ of schemes, one can associate
a pair of adjoint functors $\phi^*:\prelog_X\to\prelog_Y$
and $\phi_*:\prelog_Y\to\prelog_X$. Let 
$(\underline M,\alpha:\underline M\to\cO_X)$ (resp.
$(\underline N,\beta:\underline N\to\cO_Y)$) be a pre-log 
structure on $X$ (resp. on $Y$) and 
$\phi^\flat:\cO_X\to\phi_*\cO_Y$
$\phi^\sharp:\phi^{-1}\cO_X\to\cO_Y$ the natural morphisms 
(unit and counit of the adjunction $(\phi^{-1},\phi_*)$
on sheaves of $\Z$-modules); then 
$\phi^{-1}\underline M\stackrel{\phi^{-1}\alpha}{\longrightarrow}
\phi^{-1}\cO_X\stackrel{\phi^\sharp}{\longrightarrow}\cO_Y$ 
defines $\phi^*(\underline M,\alpha:\underline M\to\cO_X)$
and $\phi_*(\underline N,\beta:\underline N\to\cO_Y)$
is the morphism of sheaves of monoids
$\gamma:\phi_*\underline N\times_{\phi_*\cO_Y}\cO_X\to\cO_X$
which makes commute the cartesian diagram
$$\xymatrix{
\phi_*\underline N\times_{\phi_*\cO_Y}\cO_X
\ar[r]^-\gamma \ar[d] & \cO_X \ar[d]^-{\phi^\flat} \\
\phi_*\underline N \ar[r]^-{\phi_*\beta} & \phi_*\cO_Y.
}$$
\sset\subsubsection{}\label{subsec_def.logstruct}
\index{Log structure(s)|indref{subsec_def.logstruct}}
\index{Log structure(s)!$\mathbf{log}_X$ : category 
of|indref{subsec_def.logstruct}}
A pre-log structure $\alpha$ is said to be a 
{\em log structure\/} if 
$\alpha^{-1}(\cO^\times_X)\simeq\cO^\times_X$. We denote 
by $\mathbf{log}_X$ the category of log structures on $X$. 
The forgetful functor $\mathbf{log}_X\to\prelog_X~:~
\underline M\mapsto\underline M^\mathrm{pre\text{-}log}$ 
admits a left adjoint 
\set\begin{equation}\label{eq_functor.log.assoc}
\prelog_X\to\mathbf{log}_X\quad:\quad
\underline M\mapsto\underline M^{\log}
\end{equation}
and the resulting the diagram: 
\set\begin{equation}\label{eq_resulting.diag}
{\diagram{
\alpha^{-1}(\cO^\times_X) \ar[r] \ar[d] & 
\underline M \ar[d] \\
\cO^\times_X \ar[r] & \underline M^{\log}
}\enddiagram}\end{equation}
is cocartesian in the category of pre-log structures.
From this, one can easily verify that the unit of the
adjunction : $\underline M\mapsto
(\underline M^\mathrm{pre\text{-}log})^{\log}$ is an 
isomorphism for every log structure $\underline M$.

\sset\subsubsection{}
The category $\mathbf{log}_X$ admits arbitrary colimits;
indeed, since the unit of the adjunction 
\eqref{eq_functor.log.assoc} is an isomorphism,
it suffices to construct such colimits in the category 
of pre-log structures, and then apply the functor
$(-)\mapsto(-)^{\log}$ which preserves colimits, since it 
is a left adjoint. In particular, $\mathbf{log}_X$ admits 
arbitrary direct sums, and for any family 
$(\underline M_i~|~i\in I)$ of pre-log structures 
we have 
$(\oplus_{i\in I}\underline M_i)^{\log}\simeq
\oplus_{i\in I}\underline M_i^{\log}$.

\sset\subsubsection{}
For any morphism of schemes $Y\to X$ we remark that, if
$(\underline M,\alpha)$ is a log structure on $Y$, then the 
pre-log structure $\phi_*(\underline M,\alpha)$ is actually
a log structure (this can be checked on the stalks). We
deduce a pair of adjoint functors $(\phi^*,\phi_*)$ for log
structures, as in \eqref{eq_prelog.and.maps}. These are formed 
by composing the corresponding functors for pre-log structures 
with the functor \eqref{eq_functor.log.assoc}.

\sset\subsubsection{}\label{subsec_def.reg.log}
\index{Log structure(s)!regular|indref{subsec_def.reg.log}}
We say that a log structure $\underline M$ is {\em regular\/}
if $\underline M=(M_X)^{\log}$ for some free monoid $M$, and
the associated morphism of monoids $\phi:M\to\Gamma(X,\cO_X)$
maps $M$ into the set of non-zero-divisors of $\Gamma(X,\cO_X)$.

\sset\subsubsection{}\label{subsec_aux.categ}
\index{$\mathbf{Hom}_{\cO_X}(\cF,*)$|indref{subsec_aux.categ}}
For an $\cO_X$-module $\cF$, denote by 
$\mathbf{Hom}_{\cO_X}(\cF,*)$ the category of all homomorphisms 
of $\cO_X$-modules $\cF\to\cA$ (for any $\cO_X$-module $\cA$). A 
morphism from $\cF\to\cA$ to $\cF\to\cB$ is a morphism $\cA\to\cB$
of $\cO_X$-modules which induces the identity on $\cF$.
This category admits arbitrary colimits.

\sset\subsubsection{}\label{subsec_log.differ}
\index{Pre-log structure(s)!$\Omega_{X/\Z}(\log\underline M)$ : logarithmic differentials of a|indref{subsec_log.differ}}
Given a pre-log structure $\alpha:\underline M\to\cO_X$,
one defines the sheaf of {\em logarithmic differentials\/}
$\Omega_{X/\Z}(\log\underline M)$ as the quotient of
the $\cO_X$-module $\Omega_{X/\Z}\oplus
(\cO_X\otimes_{\Z_X}\underline M^\mathrm{gp})$ by the
$\cO_X$-submodule generated by the local sections of
the form $(d\alpha(m),-\alpha(m)\otimes m)$, for every 
local section $m$ of $\underline M$. (The meaning
of this is, that one adds to $\Omega_{X/\Z}$ the
logarithmic differentials $\alpha(m)^{-1}d\alpha(m)$).
For every local section $m$ of $\underline M$, we
denote by $d\log(m)$ the image of $1\otimes m$ in
$\Omega_{X/\Z}(\log\underline M)$. The assignment 
$\underline M\mapsto(\Omega_{X/\Z}\to
\Omega_{X/\Z}(\log\underline M))$
defines a (covariant) functor :
$$
\Omega:\prelog_X\to\mathbf{Hom}_{\cO_X}(\Omega_{X/\Z},*).
$$

\begin{lemma}\label{lem_regular.mnd}
Let $X$ be a scheme.
\begin{enumerate}
\item 
The functor\/ $\Omega$ commutes with all colimits.
\item
The functor\/ $\Omega$ factors through the functor 
\eqref{eq_functor.log.assoc}.
\item
Let $j:U\to X$ be a formally {\'e}tale morphism of 
schemes and $\underline M$ a log structure on $X$. 
Then the natural morphism:
$j^*\Omega_{X/\Z}(\log\underline M)\to
\Omega_{U/\Z}(\log j^*\underline M)$
is an isomorphism.
\item
If $\underline M$ is a regular log structure, then 
$\Omega(\underline M)$ is a monomorphism of $\cO_X$-modules.
\end{enumerate}
\end{lemma}
\begin{proof} (i): It is clear that $\Omega$ commutes with
filtered colimits. Thus, to show that it commutes with all
colimits, it suffices to show that it commutes with finite
direct sums and with coequalizers. We consider first direct
sums. We have to show that, for any
two pre-log structures $\underline M_1$ and $\underline M_2$, 
the natural morphism
$$\Omega_{X/\Z}(\log\underline M_1)\mathop\amalg_{\Omega_{X/\Z}}
\Omega_{X/\Z}(\log\underline M_2)\to
\Omega_{X/\Z}(\log\underline M_1\oplus\underline M_2)$$
is an isomorphism.
Notice that the functor 
$(-)\mapsto(-)^\mathrm{gp}$ of \eqref{subsec_asso.sh.groups}
commutes with colimits, since it is a left adjoint.
It follows that the diagram
$$\xymatrix{
\Omega_{X/\Z} \ar[r] \ar[d] & 
\Omega_{X/\Z}\oplus
(\cO_X\otimes_{\Z_X}\underline M_1^\mathrm{gp}) \ar[d] \\
\Omega_{X/\Z}\oplus
(\cO_X\otimes_{\Z_X}\underline M_2^\mathrm{gp}) \ar[r] & 
\Omega_{X/\Z}\oplus
(\cO_X\otimes_{\Z_X}
(\underline M_1\oplus\underline M_2)^\mathrm{gp})
}$$
is cocartesian. Thus, we are reduced to show that the kernel 
of the map 
$$\Omega_{X/\Z}\oplus(\cO_X\otimes_{\Z_X}\underline M^\mathrm{gp})
\to\Omega_{X/\Z}(\log\underline M)$$
is generated by the images of the kernels of the corresponding
maps relative to $\underline M_1$ and $\underline M_2$.
However, any section of $\underline M_1\oplus\underline M_2$ 
can be written locally in the form $x\cdot y$ for two local 
sections $x$ of $\underline M_1$ and $y$ of $\underline M_2$. 
Then we have :
$$\begin{array}{r@{\:=\:}l}
(d\alpha(x\cdot y),
-\alpha(x\cdot y)\otimes(x\cdot y)) &
(\alpha(x)\cdot d\alpha(y)+
\alpha(y)\cdot d\alpha(x),
-\alpha(x)\cdot\alpha(y)\otimes(x\cdot y)) \\
& \alpha(x)\cdot(d\alpha(y),-\alpha(y)\otimes y)+
\alpha(y)\cdot(d\alpha(x),-\alpha(x)\otimes x)
\end{array}$$
so the claim is clear. Next, suppose that 
$\phi,\psi:\underline M\to\underline N$ are two morphisms of
pre-log structures. Let $:\alpha:\underline Q\to\cO_X$ be 
the coequalizer of $\phi$ and $\psi$. Clearly, $\underline Q$ 
is the coequalizer of $\phi$ and $\psi$
in the category of sheaves of monoids. The functor
$\underline M\mapsto\underline M^\mathrm{gp}$ preserves
colimits, so we are reduced to consider the cokernel
of $\beta:=\phi^\mathrm{gp}-\psi^\mathrm{gp}$.
Moreover, clearly we have 
$\Coker(\beta)\otimes_{\Z_X}\cO_X\simeq
\Coker(\beta\otimes_{\Z_X}\cO_X)$; the claim follows easily.

(ii): Let us apply the functor $\Omega$ to the cocartesian
diagram \eqref{eq_resulting.diag}. In view of (i), the
resulting diagram of $\cO_X$-modules is cocartesian.
However, it is easy to check that 
$\Omega_{X/\Z}(\log\alpha^{-1}(\cO^\times_X))\simeq
\Omega_{X/\Z}(\log\cO^\times_X)\simeq\Omega_{X/\Z}$.
The assertion follows directly.

(iii): one uses \cite[Ch.IV, Cor. 17.2.4]{EGA4};
the details will be left to the reader.

(iv): By (ii), the functor $\Omega$ descends to a
functor 
\set\begin{equation}\label{eq_new.omega}
\mathbf{log}_X\to
\mathbf{Hom}_{\cO_X}(\Omega_{X/\Z},*).
\end{equation}
Since the unit of the adjunction \eqref{eq_functor.log.assoc}
is an isomorphism, it follows easily that \eqref{eq_new.omega}
commutes with all colimits of log structures.
Hence, to verify that $\Omega(\underline M)$ is a
monomorphism when $\underline M$ is regular, we are 
immediately reduced to the case when $\underline M$ 
is the regular log structure associated to a morphism of 
monoids $\phi:\N\to\Gamma(X,\cO_X)$. Let $f:=\phi(1)$. 
It is easy to check that in this case, the diagram
$$\xymatrix{
\cO_X \ar[r]^-f \ar[d]_{df} & \cO_X \ar[d]^{d\log f}\\
\Omega_{X/\Z} \ar[r]^-\beta & \Omega_{X/\Z}(\log\underline M)
}$$
is cocartesian. By assumption, $f$ is a non-zero-divisor,
thus multiplication by $f$ is a monomorphism of
$\cO_X$-modules, so the assertion follows. 
\end{proof}
\sset\subsubsection{}\label{subsec_this.general}
\index{$K^+$ : Valuation ring(s)!total log structure on a|indref{subsec_this.general}}
This general formalism will be applied here to the
following situation. We consider the submonoid 
$M:=K^+\setminus\{0\}$ of $K^+$. 
The imbedding $M\subset K^+$ induces a log structure on
$\Spec\,K^+$, which we call the {\em total log structure\/}
on $K^+$. More generally, we consider the natural
projection $\pi:M\to\Gamma^+$ (see 
\eqref{subsec_Mnds.in.Gamma}); then for every submonoid 
$N\subset\Gamma^+$, we have a log structure $\underline N$ 
corresponding to the imbedding $\pi^{-1}(N)\subset K^+$. 
To ease notation, we will denote by $\Omega_{K^+/\Z}(\log N)$ 
the corresponding $K^+$-module of logarithmic differentials.

\begin{proposition}\label{prop_boom} In the situation of 
\eqref{subsec_this.general}, let $\Delta\subset\Gamma$ 
be any subgroup, $N$ a prime ideal of $\Delta^+$ (cp. 
\eqref{subsec_ideals.in.mnd}) and suppose that the convex rank 
of $\Sigma:=\Delta/(\Delta^+\setminus N)^\mathrm{gp}$ equals 
one. Then we have a short exact sequence
$$0\to\Omega_{K^+/\Z}(\log\Delta^+\setminus N)
\stackrel{j}{\to}\Omega_{K^+/\Z}(\log\Delta^+)
\stackrel{\rho}{\to}
\Sigma\otimes_\Z(K^+/\pi^{-1}(N)\cdot K^+)\to 0.$$
\end{proposition}
\begin{proof} Let us first remark that the assumptions and
the notation make sense : indeed, since $N$ is a prime
ideal of $\Delta^+$, it follows that $M:=\Delta^+\setminus N$
is a convex submonoid of $\Delta^+$, hence $M=(M^\mathrm{gp})^+$
and $M^\mathrm{gp}$ is a convex subgroup of $\Delta$ 
(cp. \eqref{subsec_satur.mnds}), therefore $\Sigma$ is an 
ordered group (cp. example \eqref{ex_tens.with.Q}(v)), 
and hence it makes sense to say that its convex rank equals 
one.

Let us show that $j$ is injective. We can write $\Delta^+$
as the colimit of the filtered family of its finitely
generated submonoids $F_\alpha$. For each such $F_\alpha$,
theorem \ref{th_tricky.monds} gives us a free finitely
generated submonoid $L_\alpha\subset\Delta^+$ such that
$F_\alpha\subset L_\alpha$. Clearly $\Delta^+$ is the colimit 
of the $L_\alpha$, and $M$ is the colimit of the 
$M_\alpha:=M\cap L_\alpha$. Thus $\underline\Delta^+$
is the colimit of the $\underline L_\alpha$ and $\underline M$
is the colimit of the $\underline M_\alpha$. Let $S_\alpha$ 
be a basis of $L_\alpha$. Since $M$ is convex in $\Delta$, 
we see that $M_\alpha$ is free with basis $S_\alpha\cap M$ 
and $L_\alpha=M_\alpha\oplus N_\alpha$, where $N_\alpha$
is the free submonoid spanned by $S_\alpha\setminus M$.
For each $e\in S_\alpha\setminus M$, pick arbitrarily an 
element $x_e\in K^+$ such that $|x_e|=e$. The map $e\mapsto x_e$
can be extended to a map of monoids $N_\alpha\to K^+$,
and then to a pre-log structure $\nu_\alpha:
(N_\alpha)_{\Spec V}\to\cO_{\Spec K^+}$. 
Clearly we have an isomorphism of pre-log structures:
$\underline L_\alpha=\underline M_\alpha\oplus\nu_\alpha$.
Since the formation of logarithmic differentials commutes
with colimits of monoids, we are reduced  to showing that 
the analogous map 
$$j_\alpha:\Omega_{K^+/\Z}(\log M_\alpha)
\to\Omega_{K^+/\Z}(\log L_\alpha)$$
is injective. By lemma \ref{lem_regular.mnd}(i), we have 
$\Ker(j_\alpha)\simeq\Ker(\Omega_{K^+/\Z}\to
\Omega_{K^+/\Z}(\log\nu_\alpha))$.
By lemma \ref{lem_regular.mnd}(iv), the latter map
is injective, whence the assertion.

Next we proceed to show how to construct $\rho$.
Define a map 
$$\tilde\rho:X:=\Omega_{K^+/\Z}\oplus
(\pi^{-1}\Delta\otimes_\Z K^+)\to
\Sigma\otimes_\Z(K^+/\pi^{-1}(N)\cdot K^+)$$
by the rule : $(\omega,a\otimes b)\mapsto
\bar{\pi(a)}\otimes\bar b$,
for any $\omega\in\Omega_{K^+/\Z}$, $a\in\pi^{-1}\Delta$,
$b\in K^+$. 
\begin{claim}\label{cl_rho_descends}
$\Ker\tilde\rho$ contains the kernel of 
the surjection $X\to\Omega_{K^+/\Z}(\log\Delta^+)$. 
\end{claim}
\begin{pfclaim} It suffices to show that 
$a\otimes|a|\in\Ker\tilde\rho$ whenever 
$\pi(a)\notin(\Delta^+\setminus N)^\mathrm{gp}$.
However, $\pi(a)\notin(\Delta^+\setminus N)^\mathrm{gp}
\Leftrightarrow\pi(a)\notin\Delta^+\setminus N
\Leftrightarrow\pi(a)\in N\Leftrightarrow a\in\pi^{-1}(N)$,
so the claim follows.
\end{pfclaim}

By claim \ref{cl_rho_descends} we deduce that $\tilde\rho$
descends to the map $\rho$ as desired. It is now obvious
that $\rho$ is surjective and that its kernel contains the 
image of $j$. To conclude the proof, it suffices to show 
that the cokernel of $j$ is annihilated by $\pi^{-1}N$.
We are thus reduced to showing that $\pi^{-1}(N)$
annihilates the classes in $\Coker(j)$ of the
elements $d\log(e)$, for every $e\in\pi^{-1}N$.
Let $a\in\pi^{-1}(N)$. Since the convex rank of $\Sigma$ 
equals one, and $\pi(e)\in N$, there exists $k\geq 0$ 
and $b\in K^+$ such that $e=a^k\cdot b$ and $|b|<|a|$. 
In particular, $|b|\in\Delta^+$, and we can write:
$d\log(e)=d\log(a^k\cdot b)=b\cdot k\cdot d\log(a)+
a^k\cdot d\log(b)$, and it is clear that $a$ annihilates
each of the terms of this expression.
\end{proof}

\begin{corollary}\label{cor_boom}
In the situation of \eqref{subsec_this.general}, we have :
\begin{enumerate}
\item
The natural map : 
$\beta_{K^+}:\Omega_{K^+/\Z}\to\Omega_{K^+/\Z}(\log\Gamma^+)$ 
is injective.
\item
Suppose moreover that $K^+$ has finite rank. Let 
$\fp_r:=0\subset\fp_{r-1}\subset...\subset\fp_0:=\fm_K$
be the chain of all the prime ideals of $K^+$.
Denote by $\Delta_r:=\Gamma_K\supset\Delta_{r-1}\supset...
\supset\Delta_0:=0$ the corresponding ascending chain
of convex subgroups of $\Gamma_K$ (see 
\eqref{subsec_bij.convex}). Then $\Coker\,\beta_{K^+}$
admits a finite filtration 
$\mathrm{Fil}^\bullet(\Coker\,\beta_{K^+})$ indexed by the 
totally ordered set $\Spec(K^+)$, such that :
$$\gr^{\fp_i}(\Coker\,\beta_{K^+})\simeq
(\Delta_{i+1}/\Delta_i)\otimes_\Z(K^+/\fp_i)\quad
\text{for every $\fp_i\in\Spec\,K^+$}.$$
\end{enumerate}
\end{corollary}
\begin{proof} (i): Since the formation of differentials and
logarithmic differentials commutes with colimits of 
$\Z$-algebras and log structures, we can reduce to the
case where $K$ is a field of finite type over its prime
field. In this case the convex rank of $\Gamma$ is finite,
so the assertion follows from proposition \ref{prop_boom} 
and an easy induction.

(ii): is a straightforward consequence of proposition 
\ref{prop_boom}.
\end{proof}

\subsection{Transcendental extensions}\label{sec_trasc_exts}
In this section we extend the results of section 
\ref{sec_val.rings} to the case of arbitrary extensions of 
valued fields.

\sset\subsubsection{}\label{subsec_fix-notation_rho}
\index{$\rho_{E^+/K^+}$|indref{subsec_fix-notation_rho}}
We fix the following notation throughout this section. 
For a valued field extension $(K,|\cdot|)\subset(E,|\cdot|_E)$,
we let
$$
\rho_{E^+/K^+}:E^+\otimes_{K^+}\Omega_{K^+/\Z}(\log\Gamma^+)\to
\Omega_{E^+/\Z}(\log\Gamma^+_E)
$$
be the natural morphism. One of the main results of this
section states that $\Coker(\rho_{E^+/K^+})$ is injective
with torsion-free cokernel when $K$ is algebraically closed
(theorem \ref{th_finally}) or when $K$ has characteristic zero
(lemma \ref{lem_in.char.zero}). 
Lurking behind the results of this sections there should 
be some notion of "logarithmic cotangent complex", which 
however is not currently available.

\sset\subsubsection{}
Let $\fG:=(G,j,N,\leq)$ be a datum as in  
\eqref{subsec_datum}. We wish to study the total
log structure of the valued field $(K(\fG),|\cdot|_\fG)$.
We consider the morphism of monoids 
\set\begin{equation}\label{eq_morph.of.monds}
G\to K[\fG]\quad:\quad g\mapsto[g].
\end{equation} 
Let $K[\fG]^+$ be the subring of the elements $x\in K[\fG]$
such that $|x|_\fG\leq 1$. Let $\pi:G\to\Gamma_\fG$ be the 
projection; for every submonoid $M\subset\Gamma_\fG^+$, 
the preimage $\pi^{-1}M$ is a submonoid of $G$, and the 
restriction of \eqref{eq_morph.of.monds} induces a morphism 
of monoids $\pi^{-1}(M)\to K[\fG]^+$, whence a pre-log 
structure $\pi^{-1}M_X$ on $X:=\Spec\,K[\fG]^+$ (see 
\eqref{subsec_recall.prelog}). To ease notation, we
set $\underline M:=(\pi^{-1}M_X)^{\log}$ and we will
write $\Omega_{X/\Z}(\log M)$ for the associated sheaf
of log differentials.

\begin{lemma}\label{lem_K.fG} 
Resume the notation of \eqref{subsec_datum}. Then the 
natural diagram
$$\xymatrix{
K^\times\otimes_\Z K[\fG]^+\ar[r]^-\alpha \ar[d]_\beta & 
G\otimes_\Z K[\fG]^+ \ar[d]^\eta \\
\Omega_{K^+/\Z}(\log\Gamma^+)\otimes_{K^+}K[\fG]^+ 
\ar[r] & \Omega_{K[\fG]^+/\Z}(\log\Gamma_\fG^+)
}$$
is cocartesian.
\end{lemma}
\begin{proof} 
Let $P$ be the push out of $\alpha$ and $\beta$. We already 
have a map $\phi:P\to\Omega_{K[\fG]^+/\Z}(\log\Gamma_\fG^+)$, 
and by inspecting the definition of $K[\fG]^+$ one verifies 
easily that $\phi$ is surjective; thus we need only find a 
left inverse for $\phi$. Let us remark also that $\beta$, and 
consequently $\eta$, is surjective, hence it suffices to
 exhibit:
\begin{enumerate}
\renewcommand{\labelenumi}{(\alph{enumi})}
\item
a derivation $\delta:K[\fG]^+\to P$ such that 
$\eta(\delta a)=da$ for every $a\in K[\fG]^+$;
\item
a $\Z$-linear map 
$\psi:(\pi^{-1}\Gamma^+_\fG)^\mathrm{gp}\to G$ such that
$\eta\circ\psi(\gamma)=d\log(g)$ for every 
$g\in\pi^{-1}\Gamma^+_\fG$.
\end{enumerate}
Of course we can take for $\psi$ the natural identification
$(\pi^{-1}\Gamma^+_\fG)^\mathrm{gp}\stackrel{\sim}{\to}G$.
To define $\delta$, choose arbitrarily a set of representatives
$(g_\gamma~|~\gamma\in G/K^\times)$ for the classes of 
$G/K^\times$. Then every $a\in K[\fG]$ can be written in 
a unique way as a $K$-linear combination 
$a=\sum_{\gamma\in G/K^\times}a_\gamma\cdot[g_\gamma]$; we 
define $\delta':K[\fG]^+\to G\otimes_\Z K[\fG]^+$ by the 
rule: $a\mapsto
\sum_{\gamma\in G/K^\times}g_\gamma\otimes a_\gamma$.
It is easy to check that the image of $\delta'(a)$ in $P$ does
not depend on the choices of representatives, and this defines
our sought derivation $\delta$.
\end{proof}

\sset\subsubsection{}\label{subsec_val.ringK(fG)}
Let $K(\fG)^+$ be the valuation ring
of the valuation $|\cdot|_\fG$. It is easy to see that 
$K(\fG)^+=K[\fG]^+_\fp$, where $\fp$ is the ideal of elements 
$x\in K[\fG]$ such that $|x|_\fG<1$. It then follows from lemma
\ref{lem_regular.mnd}(i),(iii) that the diagram of lemma
\ref{lem_K.fG} remains cocartesian when we replace everywhere
$K[\fG]^+$ by $K(\fG)^+$.

\begin{proposition}\label{prop_case.pure.transc} 
Suppose that $K$ is algebraically closed, let $(E,|\cdot|_E)$
be a purely transcendental valued field extension of 
$(K,|\cdot|)$ with $\mathrm{tr.deg}(E:K)=1$. Then:
\begin{enumerate}
\item
$\Omega_{E^+/K^+}$ is a torsion-free $E^+$-module and 
$H_i(\L_{E^+/K^+})=0$ for every $i>0$.
\item
The natural map of $E^+$-modules:
$\Omega_{K^+/\Z}(\log\Gamma^+)\otimes_{K^+}E^+\to
\Omega_{E^+/\Z}(\log\Gamma_E^+)$
is injective with torsion-free cokernel.
\item
Suppose that $\Gamma_E=\Gamma$. Then the natural diagram
$$
\xymatrix{
\Omega_{K^+/\Z}\otimes_{K^+}E^+ \ar[rr] \ar[d] & &
\Omega_{E^+/\Z} \ar[d] \\
\Omega_{K^+/\Z}(\log\Gamma^+)\otimes_{K^+}E^+ 
\ar[rr]^-{\rho_{E^+/K^+}} & &
\Omega_{E^+/\Z}(\log\Gamma_E^+)
}$$
is cocartesian.
\end{enumerate}
\end{proposition}
\begin{proof}Let $\fm_E$ be the maximal ideal of $E^+$
and $X\in E$ such that $E=K(X)$.
Following \eqref{subsec_to_see_how}, we distinguish two cases,
according to whether there exists or there does not exist an 
element $a\in K$ which minimizes the function 
$K\to\Gamma_E~:~b\mapsto|X-b|_E$. Suppose first that such
an element does not exist. We pick a net $(a_i~|~i\in I)$
satisfying the conditions of lemma \ref{lem_approx.sequence}, 
relative to the element $x:=X$. For a given $b\in K$,
choose $i\in I$ such that $|X-a_i|_E<|X-b|_E$; then we have
$|a_i-b|=|X-b|_E$ and it follows easily that 
$\Gamma_E=\Gamma_K$ in this case. Then, for every $i\in I$ 
we can find $b_i\in K$ such that $|X-a_i|_E=|b_i|$. Let 
$f(X)/g(X)\in E^+$ be the quotient of two elements 
$f(X),g(X)\in K[X]$. By lemma \ref{lem_approx.sequence}, 
we have $|f(X)/g(X)|_{(a_i,|b_i|)}\leq 1$ and 
$\gamma:=|g(X)|_E=|g(X)|_{(a_i,|b_i|)}$ for every sufficiently
large $i\in I$.  Pick $a\in K$ such that $|a|=\gamma$.
Arguing as in the proof of case (b) of proposition 
\ref{prop_firstblood}(i), we deduce that
$a^{-1}\cdot g(X),a^{-1}\cdot f(X)\in A_i:=K^+[(X-a_i)/b_i]$, 
and if we let $\fp_i:=A_i\cap\fm_{E^+}$,
then $f(X)/g(X)\in E^+_i:=A_{i,\fp_i}$. 
It is also easy to see that the family of the $K^+$-algebras 
$E^+_i$ is filtered by inclusion. Clearly $\Omega_{E^+_i/K^+}$
is a free $E^+_i$-module of rank one, and 
$H_i(\L_{E^+_i/K^+})=0$ for every $i>0$, so (i) follows
easily in this case. Notice that, since $\Gamma=\Gamma_E$,
the log structure $\underline\Gamma^+_E$ on $\Spec\,E^+$
(notation of \eqref{subsec_this.general}) is the log structure
associated to the morphism of monoids 
$(K^+)\setminus\{0\}\to E^+$.
It follows easily that, for every $i\in I$, we have a 
cocartesian diagram
\set\begin{equation}\label{eq_analogue.for.A_i}
{\diagram{
\Omega_{K^+/\Z}\otimes_{K^+}A_i \ar[r]^-{\alpha_i} \ar[d] &
\Omega_{A_i/\Z} \ar[d] \\
\Omega_{K^+/\Z}(\log\Gamma^+)\otimes_{K^+}A_i \ar[r] &
\Omega_{A_i/\Z}(\log\Gamma^+)
}\enddiagram}
\end{equation}
where moreover, $\alpha_i$ is split injective; the diagram
of (iii) is obtained from \eqref{eq_analogue.for.A_i}, by 
localizing at $\fp$ and taking colimits over the family $I$; 
since both operations preserve colimits, we get (ii) and (iii) 
in this case. Finally, suppose that there exists an
element $a\in K$ such that $|X-a|$ is minimal; we can
replace $X$ by $X-a$, and thus assume that $a=0$. By 
\eqref{subsec_to_see_how} it follows that $|\cdot|_E$
is a Gauss valuation; then this case can be realized
as the valuation $|\cdot|_\fG$ associated to the
datum $\fG:=(K^\times\oplus\Z,j,N,\leq)$, where $j$
is the obvious imbedding, and $N$ is either $\Z$ or
$\{0\}$, depending on whether $|X|_E\in\Gamma$ or 
otherwise. In either case, \eqref{subsec_val.ringK(fG)} 
tells us that the map of (ii) is split injective, with 
cokernel isomorphic to $E^+$, so (ii) holds. Suppose first
that $|X|_E\in\Gamma$. Then we can find $b\in K$ such that 
$|X/b|_E=1$, and one verifies easily that $E^+$ is the 
localization of $A:=K^+[X/b]$ at the prime ideal 
$\fm_K\cdot E^+$. Clearly \eqref{eq_analogue.for.A_i}
remains cocartesian when we replace $A_i$ by $A$ and
$\alpha_i$ by the corresponding map $\alpha$; the latter
is still split injective, so (iii) follows easily.
(i) is likewise obvious in this case. In case 
$|X|_E\notin\Gamma$, we distinguish three cases. First,
suppose that $|X|_E<|b|$ for every $b\in K^\times$.
Then $K^+[X/b]\subset E^+$ for every $a\in K^\times$, and
indeed it is easy to check that $E^+$ is the filtered
union of its $K^+$-subalgebras of the form 
$K^+[X/b]_{\fp_b}$, where $\fp_b$ is the prime ideal 
generated by $\fm_K$ and $X/b$. Again (i) follows.
The second case, when $|X|_E>|b|$ for every $b\in K^\times$,
is reduced to the former, by replacing $X$ with $X^{-1}$.
It remains only to consider the case where there exist
$a_0,b_0\in K$ such that $|a_0|<|X|_E<|b_0|$; then we can find
a net $(a_i,b_i~|~i\in I)$ consisting of pairs of elements
of $K^\times$, such that $|a_i|<|X|_E<|b_i|$ for every 
$i\in I$, and moreover, for every $a,b\in K^\times$
such that $|a|<|X|<|b|$, there exists $i_0\in I$ with
$|a|<|a_i|$ and $|b_i|<|b|$ whenever $i\geq i_0$.
In such a situation, one verifies easily that $E^+$
is the filtered union of its $K^+$-subalgebras of
the form $E^+_i:=K^+[a_i/X,X/b_i]_{\fp_i}$, where $\fp_i$
is the prime ideal generated by $\fm_K$ and the elements
$a_i/X$, $X/b_i$. Each ring $E^+_i$ is a complete intersection
$K^+$-algebra, isomorphic to $K^+[X,Y]/(X\cdot Y-a_i/b_i)$.
It follows that $\L_{E^+_i/K^+}$ is acyclic in degrees
$>0$, and $\Omega_{E^+_i/K^+}\simeq 
X\cdot E^+_i\oplus Y\oplus E^+_i/(XdY+YdX)$. We leave
to the reader the verification that this $E^+_i$-module is 
torsion-free.
\end{proof}

\begin{theorem}\label{th_was.theor.three} 
Let $(K,|\cdot|)\subset(E,|\cdot|_E)$ be any extension 
of valued fields. Then 
\begin{enumerate}
\item
$H_i(\L_{E^+/K^+})=0$ for every $i>1$ and $H_1(\L_{E^+/K^+})$ 
is a torsion-free $E^+$-module.
\item
If $K$ is perfect, then $H_i(\L_{E^+/K^+})=0$ for every $i>0$.
\end{enumerate}
\end{theorem}
\begin{proof} Let us show first how to deduce (ii) from (i).
Indeed, we reduce easily to the case where $E$ is finitely
generated over $K$. Then, if $K$ is perfect, we can find
a subextension $F\subset E$ which is purely transcendental
over $K$, and such that $E$ is separable over $F$; by 
transitivity, we deduce that 
$\L_{E/K}\simeq E\otimes_F\L_{F/K}$; moreover 
$H_i(\L_{F/K})=0$ for $i>0$; by (i) we know that 
$H_1(\L_{E^+/K^+})$ imbeds into 
$H_1(\L_{E^+/K^+})\otimes_{E^+}E\simeq H_1(\L_{E/K})$,
so the assertion follows.

To show (i), let $|\cdot|_{E^\mathrm{a}}$ be a valuation
on the algebraic closure $E^\mathrm{a}$ of $E$, which
extends $|\cdot|_E$; recall that $E^{\mathrm{a}+}$ is
a faithfully flat $K^+$-module by remark
\ref{rem_val.ring.princip}(ii). We apply transitivity to 
the tower $K^+\subset E^+\subset E^{\mathrm{a}+}$ to see 
that the theorem holds for the extension 
$(K,|\cdot|)\subset(E,|\cdot|_E)$
if and only if it holds for 
$(K,|\cdot|)\subset(E^\mathrm{a},|\cdot|_{E^\mathrm{a}})$
and for 
$(E,|\cdot|_E)\subset(E^\mathrm{a},|\cdot|_{E^\mathrm{a}})$.
For the latter extension the assertion is already known
by theorem \ref{th_gen.alg.case}(i), so we are reduced
to prove the theorem for the case 
$(K,|\cdot|)\subset(E^\mathrm{a},|\cdot|_{E^\mathrm{a}})$.
Similarly, we apply transitivity to the tower
$K^+\subset K^{\mathrm{a}+}\subset E^{\mathrm{a}+}$
to reduce to the case where both $K$ and $E$ are
algebraically closed. Then we can write $E$ as the
filtered union of the algebraic closures $E_i^\mathrm{a}$ 
of its finitely generated subfields $E_i$, thereby reducing 
to prove the theorem for the extensions 
$K\subset E_i^\mathrm{a}$; hence we can assume that 
$\mathrm{tr.d}(E:K)$ is finite.
Again, by transitivity, we further reduce to the case
where the transcendence degree of $E$ over $K$ equals
one. In this case, we can pick an element $X\in E$
transcendent over $K$, and write $E=K(X)^\mathrm{a}$.
Using once more transitivity, we reduce to show the
assertion for the purely transcendental extension
$K\subset K(X)$, in which case proposition 
\ref{prop_case.pure.transc}(i) applies, and concludes 
the proof.
\end{proof}

\begin{lemma}\label{lem_oldstyle} 
Let $R\to S$ be a ring homomorphism.
\begin{enumerate}
\item
Suppose that $\F_p\subset R$, denote by $\Phi_R:R\to R$ the 
Frobenius endomorphism of $R$, and define similarly $\Phi_S$. 
Let $R_{(\Phi)}:=\Phi_R^*R$ and $S_{(\Phi)}:=\Phi_S^*S$ (cp. 
\eqref{subsec_pull_by_Frob}). Suppose moreover that the 
natural morphism :
$$
R_{(\Phi)}\derotimes_RS\to S_{(\Phi)}
$$ 
is an isomorphism in $\sD(R\Mod)$. Then 
$\L_{S/R}\simeq 0$ in $\sD(s.S\Mod)$.
\item
Suppose that $S$ is a flat $R$-algebra and let $p$ be a
prime integer, $b\in R$ a non-zero-divisor such that 
$p\cdot R\subset b^p\cdot R$. Suppose moreover that the
Frobenius endomorphisms of $R':=R/b^p\cdot R$ and
$S':=S/b^p\cdot S$ are surjective. Then the natural 
morphism :
$$
\L_{S/R}\to\L_{S[b^{-1}]/R[b^{-1}]}
$$
is an isomorphism in $\sD(s.S\Mod)$.
\end{enumerate}
\end{lemma}
\begin{proof} (i): Let $P^\bullet:=P^\bullet_R(S)$ be the 
standard simplicial resolution of $S$ by free $R$-algebras. 
Then $\L_{S/R}\simeq\Omega_{P^\bullet/R}\otimes_{P^\bullet}S$.
Let $\Phi_{P^\bullet}:P^\bullet\to P^\bullet_{(\Phi)}$ 
be the termwise Frobenius endomorphism of the simplicial 
algebra $P^\bullet$. As usual, we can write $\Phi_{P^\bullet}=
(\Phi_R\otimes_R\one_{P^\bullet})\circ\Phi_{P^\bullet/R}$,
where the relative Frobenius 
$\Phi_{P^\bullet/R}:R_{(\Phi)}\otimes_RP^\bullet\to 
P^\bullet_{(\Phi)}$ is a morphism of simplicial 
$R_{(\Phi)}$-algebras.
Concretely, if $P^k=R[X_i~|~i\in I]$ is a free algebra
on generators $(X_i~|~i\in I)$, then 
\set\begin{equation}\label{eq_P^bullet.oldstyle}
\Phi_{P^k/R}(X_i)=X_i^p\quad 
\text{for every $i\in I$.}
\end{equation}
Under the assumption of the lemma, $\Phi_{P^\bullet/R}$
is a quasi-isomorphism of simplicial $R_{(\Phi)}$-algebras.
It then follows from \cite[Ch.II, Prop.1.2.6.2]{Il}
that $\Phi_{P^\bullet/R}$ induces an isomorphism
\set\begin{equation}\label{eq_the_latter}
R_{(\Phi)}\otimes_R\L_{S/R}\stackrel{\sim}{\to}
\L_{S_{(\Phi)}/R_{(\Phi)}}\simeq\L_{S/R}.
\end{equation}
However, \eqref{eq_P^bullet.oldstyle} shows that 
\eqref{eq_the_latter} is represented by a map of simplicial
complexes which is termwise the zero map, so the 
claim follows.

(ii): Under the stated assumptions, the Frobenius map induces
an isomorphism of $R$-algebras: 
$R/b\cdot R\stackrel{\sim}{\to}R'_{(\Phi)}$ (resp. of
$S$-algebras: $S/b\cdot S\stackrel{\sim}{\to}S'_{(\Phi)}$).
Thus the map
$$
S'\otimes_{R'}(R')_{(\Phi)}\to S'_{(\Phi)}\quad:\quad
x\otimes y\mapsto\Phi_{S'}(x)\cdot y
$$
is an isomorphism. Since moreover $S'$ is a flat $R'$-algebra,
we see that the assumption of (i) is satisfied, whence 
$\L_{S'/R'}\simeq 0$. If we now tensor the short exact sequence
$0\to R\to R\to R'\to 0$ by $\L_{S/R}$, we obtain a distinguished
triangle
$$
\L_{S/R}\stackrel{b^p}{\longrightarrow}\L_{S/R}
\to\L_{S/R}\derotimes_RR'\to\sigma\L_{S/R}.
$$
However, $\L_{S/R}\derotimes_RR'\simeq\L_{S'/R'}$ by 
\cite[II.2.2.1]{Il}, so we have shown that $b^p$ acts as
an isomorphism on $\L_{S/R}$. In other words,
$\L_{S/R}\simeq\L_{S/R}\otimes_SS[b^{-1}]\simeq
\L_{S[b^{-1}]/R[b^{-1}]}$,
as claimed.
\end{proof}

\begin{lemma}\label{lem_in.char.zero} 
Let $(K,|\cdot|)\subset(E,|\cdot|_E)$ be an extension 
of valued fields and suppose that $\Q\subset\kappa(K)$. 
Then the map $\rho_{E^+/K^+}$ of \eqref{subsec_fix-notation_rho} 
is injective with torsion-free cokernel.
\end{lemma}
\begin{proof} To begin with, let $(F,|\cdot|_F)$ be any 
valued field extension of $(E,|\cdot|_E)$. We remark that: 
\set\begin{equation}\label{eq_triple.extes}
\rho_{F^+/E^+}\circ
(\rho_{E^+/K^+}\otimes_{E^+}\one_{F^+})=\rho_{F^+/K^+}.
\end{equation}
\begin{claim}\label{cl_algebrr-case} 
Suppose moreover that $E$ is an algebraic extension of $K$. 
Then $\rho_{E^+/K^+}$ is an isomorphism.
\end{claim}
\begin{pfclaim} Applying \eqref{eq_triple.extes}, with
$F:=E^\mathrm{a}$, we reduce easily to prove the claim 
in case $E$ is algebraically closed. Let $K^\mathrm{sh}$
be the field of fractions of the strict henselization of
$K^+$ (which we see as imbedded in $E^+$). Let 
$j:\Spec\,K^{\mathrm{sh}+}\to\Spec\,K^+$ be the morphism
induces by the imbedding $K\subset K^\mathrm{sh}$; In view 
of lemma \ref{lem_Gamma-unchanged}, the log structure 
$\underline\Gamma_{K^\mathrm{sh}}^+$ on 
$\Spec\,K^{\mathrm{sh}+}$ (notation of 
\eqref{subsec_this.general}) equals $j^*\underline\Gamma$.
Since moreover $K^{\mathrm{sh}+}$ is local ind-{\'e}tale over
$K^+$, we deduce from \ref{lem_regular.mnd}(iii) that 
$\rho_{K^{\mathrm{sh}+}/K^+}$ is an isomorphism. Then 
arguing as in the foregoing, we see that it suffices to 
prove the claim for the case when $K=K^\mathrm{sh}$. 
Since everything in sight commutes with filtered unions 
of field extensions, we can even reduce to the case where 
$E$ is a finite (Galois) extension of $K$. Then, by corollary 
\ref{cor_kern.p-group}, this case can be realized as the 
extension associated to some datum $\fG:=(G,j,N,\leq)$, 
where moreover $G/K^\times$ is a finite torsion group. 
Since by assumption $\Q\subset K[\fG]^+$, the claim follows 
by lemma \ref{lem_K.fG} and \eqref{subsec_val.ringK(fG)}.
\end{pfclaim}

Now, if $K\subset E$ is an arbitrary extension, we
can apply \eqref{eq_triple.extes} with $F:=E^\mathrm{a}$
and claim \ref{cl_algebrr-case} to the extension 
$E\subset E^\mathrm{a}$ to reduce to the case where
$E$ is algebraically closed. Then we can apply again
claim \ref{cl_algebrr-case} to the extension 
$K\subset K^\mathrm{a}$ and \eqref{eq_triple.extes}
with $E:=K^\mathrm{a}$ and $F:=E$, to reduce to the
case where also $K$ is algebraically closed. Then,
by the usual argument we reduce to the case of an
extension of finite transcendence degree, and even
to the case of transcendence degree equal to one.
We factor the latter as a tower of extensions
$K\subset K(X)\subset E$, where $X$ is transcendental
over $K$, hence $E$ algebraic over $K(X)$.
So finally we are reduced to the case $E=K(X)$,
in which case we conclude by proposition 
\ref{prop_case.pure.transc}(ii).
\end{proof}

\begin{theorem}\label{th_with.K-alg.cl} 
Let $(K,|\cdot|)\subset(E,|\cdot|_E)$ be an extension 
of valued fields, with $K$ algebraically closed. Then 
$\Omega_{E^+/K^+}$ is a torsion-free $E^+$-module.
\end{theorem}
\begin{proof} Pick a valuation $|\cdot|_{E^\mathrm{a}}$
of the field $E^\mathrm{a}$ extending $|\cdot|_E$.
We have an exact sequence: 
$H_1(\L_{E^{\mathrm{a}+}/E^+})\to
E^{\mathrm{a}+}\otimes_{E^+}\Omega_{E^+/K^+}\to
\Omega_{E^{\mathrm{a}+}/K^+}$, where the leftmost term
is torsion-free by theorem \ref{th_gen.alg.case},
so it suffices to show that $\Omega_{E^{\mathrm{a}+}/K^+}$
is torsion-free, and we can therefore assume that $E$
is algebraically closed. In this case, if now $\chara(K)>0$, 
it follows that the Frobenius endomorphism of $E^+$ is 
surjective; then, for any $a\in E^+$ we can write 
$da=d(a^{1/p})^p=p\cdot a^{(p-1)/p}\cdot da^{1/p}=0$,
so actually $\Omega_{E^+/K^+}=0$. In case $\chara(K)=0$
and $\chara(\kappa(K))=p$, let us pick an element $b\in K^+$
such that $|b^p|\geq|p|$.  
Since $K$ and $E$ are algebraically closed, the Frobenius 
endomorphisms on $K^+/b^pK^+$ and $E^+/b^pE^+$ are surjective,
so $\L_{E^+/K^+}\simeq\L_{E^+[b^{-1}]/K^+[b^{-1}]}$ by
lemma \ref{lem_oldstyle}(ii).
Now, $K^+[b^{-1}]$ is the valuation ring of a valuation
$|\cdot|'$ on $K$, which extends to a valuation $|\cdot|'_E$
on $E^+$ whose valuation ring is $E^+[b^{-1}]$. Furthermore,
the residue fields of these valuations are fields of
characteristic zero. Hence, we have reduced the proof
of the theorem to the case where $\kappa(K)\supset\Q$.
We can further reduce to the case where $\mathrm{tr.d}(E:K)$
is finite and $K$ is the algebraic closure of an extension
of finite type of its prime field. By lemma 
\ref{lem_in.char.zero} we have a commutative diagram with 
exact rows:
$$
\xymatrix{
& E^+\otimes_{K^+}\Omega_{K^+/\Z} \ar[r] 
\ar[d]^{\one_{E^+}\otimes\beta_{K^+}} &
\Omega_{E^+/\Z} \ar[r] \ar[d]^{\beta_{E^+}} & 
\Omega_{E^+/K^+} \ar[r] \ar[d]^\gamma & 0 \\
0 \ar[r] & 
E^+\otimes_{K^+}\Omega_{K^+/\Z}(\log\Gamma^+) \ar[r] &
\Omega_{E^+/\Z}(\log\Gamma^+_E) \ar[r] &
\Coker(\rho_{E^+/K^+}) \ar[r] & 0
}$$
where $\beta_{K^+}$ and $\beta_{E^+}$ are the maps of 
corollary \ref{cor_boom}(i).
By virtue of lemma \ref{lem_in.char.zero}, it suffices 
to show that $\gamma$ is injective.  Since $\beta_{E^+}$ is
injective by corollary \ref{cor_boom}(i), the snake lemma
reduces us to prove :
\begin{claim}\label{cl_coke-betas}
The induced map 
$E^+\otimes_{K^+}\Coker\,\beta_{K^+}\to\Coker\,\beta_{E^+}$ 
is injective.
\end{claim}
\begin{pfclaim} Under our current assumptions, $K^+$ 
and $E^+$ are valuation rings of finite Krull dimension,
by \eqref{subsec_transc.degree}.
Let $\fp_r:=0\subset\fp_{r-1}\subset...\subset\fp_0:=\fm_E$
be the chain of all the prime ideals of $E^+$.
Denote by $\Delta_r:=\Gamma_E\supset\Delta_{r-1}\supset...
\supset\Delta_0:=0$ the corresponding ascending chain
of convex subgroups of $\Gamma_E$ (see 
\eqref{subsec_bij.convex}).
Let $\mathrm{Fil}^\bullet(\Coker\,\beta_{E^+})$ (resp.
$\mathrm{Fil}^\bullet(\Coker\,\beta_{K^+})$) be the finite 
filtration indexed by the totally ordered set $\Spec\,E^+$
(resp. $\Spec\,K^+$), provided by corollary \ref{cor_boom}(ii). 
Since it is preferable to work with a single indexing set, 
we use the surjection 
$\Spec\,E^+\to\Spec\,K^+$, to replace by $\Spec\,E^+$ 
the indexing of the filtration on $\Coker\,\beta_{K^+}$;
of course in this way some of the graded subquotients
become trivial, but we do not mind. With this notation 
we can write down the identities:
$$
E^+\otimes_{K^+}\gr^{\fp_i}(\Coker\,\beta_{K^+})\simeq
E^+\otimes_{K^+}
((\Delta_{i+1}\cap\Gamma)/(\Delta_i\cap\Gamma))
\otimes_\Z(K^+/\fp_i\cap K^+)
$$
for every $\fp_i\in\Spec(E^+)$. Furthermore, our map 
$\phi:E^+\otimes_{K^+}\Coker\,\beta_{K^+}\to
\Coker\,\beta_{E^+}$ 
respects these filtrations. If now $\fp_i\in\Spec(E^+)$ is 
the radical of the extension of a prime ideal of $K^+$, then
clearly the map $\gr^{\fp_i}(\phi):
E^+\otimes_{K^+}\gr^{\fp_i}(\Coker\,\beta_{K^+})
\to\gr^{\fp_i}(\Coker\,\beta_{E^+})$ is induced by
the imbeddings 
$(\Delta_{i+1}\cap\Gamma)/(\Delta_i\cap\Gamma)\subset
\Delta_{i+1}/\Delta_i$ and 
$K^+/\fp_i\cap K^+\subset E^+/\fp_i$, and it is therefore
injective. On the other hand, if $\fp_i$ is not the radical
of an ideal extended from $K^+$, we have 
$\gr^{\fp_i}(\Coker\,\beta_{K^+})=0$, so
$\gr^{\fp_i}(\phi)$ is trivially injective in this case
as well. Since the map $\gr^\bullet(\phi)$ is injective,
the same holds for $\phi$, which concludes the proof
of the claim and of the theorem.
\end{pfclaim}
\end{proof}

\begin{theorem}\label{th_for_sep.closure} 
Let $|\cdot|_{K^\sep}$ be a valuation on the separable 
closure $K^\sep$ of $K$, extending the valuation of $K$. 
Then the map $\rho:=\rho_{K^{\sep+}/K^+}$ is injective.
\end{theorem}
\begin{proof} Suppose first that $\Gamma$ is divisible. 
By the usual reductions, we can assume that $K$ is finitely 
generated over its prime field, hence that the convex rank 
of $\Gamma_{K^\sep}$ is finite.
We consider the commutative diagram :
$$\xymatrix{
\Omega_{K^+/\Z}\otimes_{K^+}K^{\sep+} 
\ar[r]^-\alpha 
\ar[d]_{\beta_{K^+}\otimes\one_{K^{\sep+}}} &
\Omega_{K^{\sep+}/\Z} \ar[d]^{\beta_{K^{\sep+}}} \\
\Omega_{K^+/\Z}(\log\Gamma^+)\otimes_{K^+}K^{\sep+}
\ar[r]^-\rho & 
\Omega_{K^{\sep+}/\Z}(\log\Gamma^+_{K^\sep})
}$$
where $\beta_{K^+}$ and $\beta_{K^{\sep+}}$ are the maps
of corollary \ref{cor_boom}(i).
By theorem \ref{th_gen.alg.case}(ii), $\alpha$ is injective,
and the same holds for $\beta_{K^{\sep+}}$, by corollary 
\ref{cor_boom}(i). It follows that 
$\Img(\beta_{K^+}\otimes_{K^+}\one_{K^{\sep+}})
\cap\Ker\,\rho=\{0\}$, in other words, the induced map 
\set\begin{equation}\label{eq_maps-of-kers}
\Ker\,\rho\to K^{\sep+}\otimes_{K^+}\Coker\,\beta_{K^+}
\end{equation}
is injective. By corollary \ref{cor_boom}(ii), there is
a filtration $\mathrm{Fil}^\bullet(\Coker\,\beta_{K^+})$ on 
$\Coker\,\beta_{K^+}$, indexed by $\Spec\,K^+$, such that
$\gr^{\fp_i}(\Coker\,\beta_{K^+})\simeq
(\Delta_{i+1}/\Delta_i)\otimes_\Z(K^+/\fp_i)$, where
$\Delta_i,\Delta_{i+1}$ are two convex subgroups of
$\Gamma_K$. However, since we assume that $\Gamma$ is 
divisible, the same holds for $\Delta_{i+1}/\Delta_i$; we 
deduce that $\gr^{\fp_i}(\Coker\,\beta_{K^+})$ vanishes 
whenever $\mathrm{Frac}(K^+/\fp_i)$ is a field of positive 
characteristic. What this means is that the filtration
$\mathrm{Fil}^\bullet(\Coker\,\beta_{K^+})$ is actually 
indexed by $\Spec\,K^+\otimes_\Z\Q\subset\Spec\,K^+$, and 
the natural map : 
\set\begin{equation}\label{eq_maps-of-cokers}
\Coker\beta_{K^+}\to
\Q\otimes_\Z\Coker\,\beta_{K^+}
\end{equation}
is an isomorphism. The same holds also for 
$\Coker\,\beta_{K^{\sep+}}$.
If $K^+_\Q:=K^+\otimes_\Z\Q=\{0\}$, then 
$\Coker\,\beta_{K^+}=\Coker\,\beta_{K^{\sep+}}=\{0\}$, 
hence $\Ker\,\rho=0$, which is what we had to show.

In case $K^+_\Q\neq\{0\}$, then $K^+_\Q$ is a valuation 
ring of $K$ with residue field of characteristic 
zero. However, by lemma \ref{lem_regular.mnd}(iii) it follows 
easily that $\Q\otimes_\Z\Coker\,\beta_{K^+}\simeq
\Coker\,\beta_{K^+_\Q}$, and likewise
$\rho_{K^{\sep+}/K^+}\otimes_\Z\one_\Q=
\rho_{K^{\sep+}_\Q/K^+_\Q}$, where 
$K^{\sep+}_\Q:=K^{\sep+}\otimes_\Z\Q$ 
is a valuation ring of $K^\sep$ whose valuation 
extends that of $K^+_\Q$. Since \eqref{eq_maps-of-cokers}
is an isomorphism, \eqref{eq_maps-of-kers} factors
through $\Ker\rho_{K^{\sep+}_\Q/K^+_\Q}$;
however, the latter vanishes by lemma \ref{lem_in.char.zero}.
Since \eqref{eq_maps-of-kers} is injective, we derive
$\Ker\,\rho=0$, so the theorem holds in this case.

In case $\Gamma$ is not necessarily divisible, let us
choose a datum $\fG:=(G,j,N,\leq)$ as in \eqref{subsec_datum}, 
such that $G:=(K^\sep)^\times\oplus F$, where $F$ is a 
torsion-free abelian group (whose composition law we write 
in multiplicative notation) and $N$ is the graph of a 
surjective group homomorphism 
$\phi:F\to(K^\sep)^\times$. Notice that in this case 
$\Gamma_\fG\simeq\Gamma_{K^\sep}\simeq
\Gamma\otimes_\Z\Q$, and the restriction to $F$ of the
projection $G\to\Gamma_\fG$ is the map 
$x\mapsto|\phi(x^{-1})|_{K^\sep}$.
Let now $H:=K^\times\oplus F$ and define a new datum
$\fH:=(H,j,H\cap N,\leq)$; since $\phi$ is surjective,
clearly we still have $\Gamma_\fH\simeq\Gamma\otimes_\Z\Q$.
Notice as well that $K^\sep(\fG)$ is separable 
over $K(\fH)$.
Set $\rho_\fH:=\rho_{K(\fH)^+/K^+}$, 
$\rho_\fG:=\rho_{K^\sep(\fG)^+/K^{\sep+}}$
and $\rho_{\fG/\fH}:=\rho_{K^\sep(\fG)^+/K(\fH)^+}$.
We consider the diagram :
$$\xymatrix{
\Omega_{K^+/\Z}(\log\Gamma^+)\otimes_{K^+}K^\sep(\fG)^+
\ar[d]_{\rho_\fH\otimes\one_{K^\sep(\fG)^+}}
\ar[rr]^-{\rho\otimes\one_{K^\sep(\fG)^+}} & & 
\Omega_{K^{\sep+}/\Z}(\log\Gamma^+_{K^\sep})
\otimes_{K^{\sep+}}K^\sep(\fG)^+ 
\ar[d]^-{\rho_\fG} \\
\Omega_{K(\fH)^+/\Z}(\log\Gamma_\fH^+)
\otimes_{K(\fH)^+}K^\sep(\fG)^+ 
\ar[rr]^-{\rho_{\fG/\fH}} & &
\Omega_{K^\sep(\fG)^+/\Z}(\log\Gamma^+_\fG).
}$$
Since $F$ is torsion-free, 
it follows easily from \eqref{subsec_val.ringK(fG)} and lemma 
\ref{lem_K.fG} that $\rho_\fH$ and $\rho_\fG$ are injective 
with torsion-free cokernels. 
Hence, in order to prove that $\rho$ is injective, it suffices 
to show that $\rho_{K^\sep(\fG)^+/K(\fH)^+}$ is.
Finally, let $E$ be the separable closure of $K(\fH)$
and choose a valuation on $E$ which extends the valuation of
$K^\sep(\fG)$; we notice that 
$\Ker\rho_{K^\sep(\fG)^+/K(\fH)^+}\subset
\Ker\rho_{E^+/K(\fH)^+}$. Therefore, we can replace 
$K$ by $K(\fH)$ and reduce to the case where $\Gamma$ is
divisible, which has already been dealt with.
\end{proof}

\begin{theorem}\label{th_finally}
Let $(K,|\cdot|)\subset(E,|\subset|_E)$ be an extension 
of valued fields, with $K$ algebraically closed. Then 
$\rho_{E^+/K^+}$ is injective with torsion-free cokernel.
\end{theorem}
\begin{proof} By the usual arguments, we can suppose that
$E$ is finitely generated over $K$. Let $|\cdot|_{E^\sep}$ 
be an extension of the valuation $|\cdot|_E$ to the separable 
closure $E^\sep$ of $E$. We have 
$\Ker\rho_{E^+/K^+}\subset\Ker\rho_{E^{\sep+}/K^+}$,
and even 
$\Coker\rho_{E^+/K^+}\subset\Coker\rho_{E^{\sep+}/K^+}$, 
by theorem \ref{th_for_sep.closure}. Thus we can replace
$E$ by $E^\sep$ and suppose that $E$ is separably closed,
hence $\Gamma_E$ divisible, by example \ref{ex_tens.with.Q}(ii).
By corollary \ref{cor_boom}(i) we have a commutative diagram 
with exact rows :
$$\xymatrix{
0 \ar[r] & E^+\otimes_{K^+}\Omega_{K^+/\Z} \ar[d]_\alpha 
\ar[r] &
E^+\otimes_{K^+}\Omega_{K^+/\Z}(\log\Gamma^+) \ar[r] 
\ar[d]^{\rho_{E^+/K^+}} &
E^+\otimes_{K^+}\Coker\,\beta_{K^+} \ar[r] \ar[d]^\gamma & 0 \\
0 \ar[r] & \Omega_{E^+/\Z} \ar[r] & 
\Omega_{E^+/\Z}(\log\Gamma_E^+) \ar[r] &
\Coker\,\beta_{E^+} \ar[r] & 0.
}$$
By theorem \ref{th_was.theor.three}(ii), the map $\alpha$ 
is injective. The same holds for $\gamma$, in view of claim
\ref{cl_coke-betas}. It follows already that $\rho_{E^+/K^+}$
is injective. Moreover, by theorem \ref{th_with.K-alg.cl},
$\Coker\,\alpha$ is a torsion-free $E^+$-module. Since both
$\Gamma$ and $\Gamma_E$ are divisible, it follows easily 
from corollary \ref{cor_boom}(ii) that $\Coker\beta_K$
and $\Coker\beta_E$ are $\Q$-vector spaces (cp. the proof
of theorem \ref{th_for_sep.closure}), hence the same holds 
for $\Coker\,\gamma$. Consequently, $\Coker\rho_{E^+/K^+}$
is a torsion-free $\Z$-module, and thus we are reduced to 
show that $\Q\otimes_\Z\Coker\rho_{E^+/K^+}$ is a torsion-free
$E^+$-module. However, $\Q\otimes_\Z\Coker\rho_{E^+/K^+}\simeq
\Coker\rho_{E^+_\Q/K^+_\Q}$, where $E^+_\Q:=E^+\otimes_\Z\Q$
and $K^+_\Q:=K^+_\Q$ are valuation rings with residue fields
of characteristic zero (or else they vanish, in which case
we are done). But the assertion to prove is already known in
this case, by lemma \ref{lem_in.char.zero}.
\end{proof}

\begin{corollary}\label{cor_over-perfect-field} Let 
$(K,|\cdot|)$ be a valued field, and $k$ a perfect field 
such that $k\subset K^+$. Then $\Omega_{K^+/k}$ and 
$\Coker\rho_{K^+/k}$ are torsion-free $K^+$-modules.
\end{corollary}
\begin{proof} We have $k^\mathrm{a}\subset K^{\mathrm{sh}+}$; let 
$E:=k^\mathrm{a}\cdot K\subset K^\mathrm{sh}$ and denote by 
$j:\Spec\,E^+\to\Spec\,K^+$ the morphism induced by the 
imbedding $K\subset E$. By lemma \ref{lem_Gamma-unchanged} the 
natural map : $j^*\underline\Gamma^+\to\underline\Gamma_E^+$
is an isomorphism of log structures; moreover 
$\Omega_{k^\mathrm{a}/k}=0$, since $k$ is perfect. Hence
$\Omega_{K^+/k}\subset\Omega_{E^+/k^\mathrm{a}}$, and furthermore,
by lemma \ref{lem_regular.mnd}(iii) we have
$\Coker\rho_{K^+/k}\subset\Coker\rho_{E^+/k^\mathrm{a}}$.
Then the assertion follows from theorems \ref{th_with.K-alg.cl} 
and \ref{th_finally}.
\end{proof}

\begin{remark} Notice that corollary \ref{cor_over-perfect-field}
is a straightforward consequence of a standard (as yet unproven)
conjecture on the existence of resolution of singularities over 
perfect fields.
\end{remark}

\subsection{Deeply ramified extensions}\label{sec_deeply.ram}
We keep the notation of section \ref{sec_val.rings}. We borrow
the notion of deeply ramified extension of valuation rings from
the paper \cite{Coat}, even though our definition applies more
generally to valuations of arbitrary rank.

\begin{definition}\label{def_deeply-ramif}
\index{Valued field(s)!deeply ramified|indref{def_deeply-ramif}}
Let $(K,|\cdot|)$ be a valued field, $|\cdot|_{K^\sep}$ a 
valuation on $K^\sep$ which extends $|\cdot|$.
We say that $(K,|\cdot|)$ is {\em deeply ramified\/} if 
$\Omega_{K^{\sep +}/K^+}=0$. Notice that the definition
does not depend on the choice of the extension 
$|\cdot|_{K^\sep}$.
\end{definition}

\begin{proposition}\label{prop_deep.ramif} 
Let $(K,|\cdot|)$ be a valued field whose valuation
has rank one, $|\cdot|_{K^\sep}$ an extension of $|\cdot|$ 
to $K^\sep$. Then the following conditions are equivalent :
\begin{enumerate}
\item
$(K,|\cdot|)$ is deeply ramified;
\item
The morphism of almost algebras $(K^+)^a\to(K^{\sep +})^a$ 
is weakly {\'e}tale;
\item
$(\Omega_{K^{\sep +}/K^+})^a=0$.
\end{enumerate}
Moreover, the above equivalent conditions imply that the
valuation of $K$ is not discrete.
\end{proposition}
\begin{proof} We leave to the reader the verification that
(iii) (and, {\em a fortiori}, (i)) can hold only in case
the valuation of $K$ is not discrete.
 
Let $K^{+\mathrm{sh}}$ be a strict henselization of $K^+$ 
contained in $K^{\sep +}$, and $K^\mathrm{sh}$ its fraction 
field. It is easy to check that $(K,|\cdot|)$ is deeply
ramified if and only if 
$(K^\mathrm{sh},|\cdot|_{K^\mathrm{sh}})$ is. Moreover, in 
view of lemma \ref{lem_itoiv}(iv), condition (ii) holds
for $(K,|\cdot|)$ if and only if it holds for 
$(K^\mathrm{sh},|\cdot|_{K^\mathrm{sh}})$, and similarly
for condition (iii). Hence we can assume that $K^+$ is strictly
henselian. It is also clear that 
(i)$\Rightarrow$(iii)$\Leftarrow$(ii). To show that 
(iii)$\Rightarrow$(ii), let $E\subset K^\sep$ be a finite 
separable extension of $K$ and set $E^+:=K^{\sep +}\cap E$; it
follows from theorem \ref{th_big.deal} that the natural map
$\Omega_{E^+/K^+}\otimes_{E^+}K^{\sep +}\to
\Omega_{K^{\sep +}/K^+}$ is injective; thus, if (iii) holds, 
we deduce that $(\Omega_{E^+/K^+})^a=0$ for every finite 
separable extension $E$ of $K$. Again by theorem 
\ref{th_big.deal} we derive that 
$\cD_{E^+/K^+}=E^{+a}$ for every such $E$. Finally, lemmata 
\ref{lem_diff.in.towers}(i) and \ref{lem_etale.different} show
that $E^{+a}$ is {\'e}tale over $K^{+a}$, whence (ii).
Suppose next that the residue characteristic of $K^+$ is zero;
then every finite extension $E$ of $K$ factors as a tower of 
Kummer extensions of prime degree, therefore 
$\Omega_{E^+/K^+}=0$ by corollary \ref{cor_firstblood}(i), 
which implies that (i)$\Leftrightarrow$(iii) in this case. 
Finally, suppose that the residue characteristic is $p>0$. 
Let us choose $a\in K^+$ such that $|a|\geq|p|$. It follows 
easily from example \ref{ex_tens.with.Q}(iii) that every 
element $x\in K^{\sep+}$ can be written in the form 
$x=y^p+a\cdot z$ for some $y,z\in K^{\sep+}$.
Hence $dx=p\cdot dy^{p-1}+a\cdot dz$, which means that 
$\Omega_{K^{\sep+}/K^+}=a\cdot\Omega_{K^{\sep+}/K^+}$. 
Therefore, even in this case we deduce (iii)$\Rightarrow$(i).
\end{proof}

\sset\subsubsection{}
Let us say that a $K^{+a}$-module $M$ is {\em $K^{+a}$-divisible\/}
if, for every $x\in K^+\setminus\{0\}$ we have $M=x\cdot M$. 

\begin{lemma}\label{lem_struct.T} Let $(K,|\cdot|)$ be a valued field 
such that $\Q\subset K$, and let $p:=\chara(\kappa)>0$. Let 
$(K^\sep,|\cdot|_{K^\sep})$ be an extension of the
valuation $|\cdot|$ to a separable closure of $K$.
Denote by $T$ the $K^{\sep +}$-torsion submodule of\/
$\Omega_{K^{\sep +}/\Z}$. Then $T\simeq K^\sep/K^{\sep +}$.
\end{lemma}
\begin{proof} Let $(\Q^\mathrm{a},|\cdot|_{\Q^{\mathrm{a}}})$
be the restriction of $|\cdot|_{K^\sep}$ to the
algebraic closure of $\Q$ in $K^\sep$. From theorem \ref{th_finally} 
it follows easily that $T\simeq\Omega_{\Q^{\mathrm{a}+}/\Z}
\otimes_{\Q^{\mathrm{a}+}}K^{\sep+}$, hence we can suppose
that $K=\Q$. From theorem \ref{th_big.deal} we deduce
that the natural map 
$\Q^{\mathrm{a}+}\otimes_{E^+}\Omega_{E^+/\Z}\to
\Omega_{\Q^{\mathrm{a}+}/\Z}$
is injective for every subextension $E\subset\Q^\mathrm{a}$.  
For every $n\in\N$ and every subextension $E\subset\Q^\mathrm{a}$, 
let $E_n:=E(\zeta_{p^n})$, where $\zeta_{p^n}$ is any primitive 
$p^n$-th root of $1$ and set $E_\infty:=\bigcup_{n>0}E_n$.

\begin{claim}\label{cl_disc.image} For every finite subextension 
$E\subset\Q^\mathrm{a}$, there exists $n\in\N$ such that the image of 
$E_n^+\otimes_{E^+}\Omega_{E^+/\Z}$ in $\Omega_{E_n^+/\Z}$
is included in the image of $E_n^+\otimes_{\Q_n^+}\Omega_{\Q_n^+/\Z}$.
\end{claim}
\begin{pfclaim} For every $n\in\N$, $E^+_n$ is a discrete valuation 
ring and $\kappa(E_n)$ is a finite separable extension of $\kappa(\Q)=\F_p$; 
from the exact sequence 
$\fm_{E_n}/\fm^2_{E_n}\to\Omega_{E^+_n/\Z}\to
\Omega_{\kappa(E_n)/\F_p}=0$
we deduce that $\Omega_{E_n^+/\Z}$ is a (torsion) cyclic $E_n^+$-module.
By comparing the annihilators of the modules under consideration, one
obtains easily the claim.
\end{pfclaim}

A standard calculation shows that 
$\Omega_{\Q^+_\infty/\Z}\simeq\Q_\infty/\Q^+_\infty$. This, together
with claim \ref{cl_disc.image} implies the lemma.
\end{proof}

\begin{proposition}\label{prop_deep.ramif-II}
Keep the notation and assumptions of proposition
{\em\ref{prop_deep.ramif}} and suppose moreover that the 
characteristic $p$ of the residue field $\kappa$ of $K^+$ 
is positive and that the valuation on $K$ is not discrete. 
Let $(K^\wedge,|\cdot|^\wedge)$ be the completion of 
$(K,|\cdot|)$ for the valuation topology. 
Then the following conditions are equivalent:
\begin{enumerate}
\item
$(K,|\cdot|)$ is deeply ramified;
\item
The Frobenius endomorphism of $K^{\wedge+}/pK^{\wedge+}$ 
is surjective;
\item
For some $b\in K^+\setminus\{0\}$ such that $1>|b|\geq|p|$, 
the Frobenius endomorphism on $(K^+/bK^+)^a$ is an 
epimorphism;
\item
$\Omega_{K^+/\Z}(\log\Gamma^+)$ is a $K^+$-divisible 
$K^+$-module;
\item
$\Omega_{K^+/\Z}$ is a $K^+$-divisible $K^+$-module;
\item
$(\Omega_{K^+/\Z})^a$ is a $K^{+a}$-divisible $K^{+a}$-module;
\item
$\Coker\,\rho_{K^{\sep+}/K^+}=0$ (notation of
\eqref{subsec_fix-notation_rho});
\item
$\Coker(\rho_{K^{\sep+}/K^+})^a=0$.
\end{enumerate}
\end{proposition}
\begin{proof} Suppose that (i) holds; then by proposition 
\ref{prop_deep.ramif} it follows that the morphism
$(K^+)^a\to(K^{\sep+})^a$ is weakly {\'e}tale, so the same holds
for the morphism $(K^+/bK^+)^a\to(K^{\sep+}/bK^{\sep+})^a$,
for every $b\in K^+$ with $|b|\geq|p|$.
In view of example \ref{ex_tens.with.Q}(iii), one sees that the 
Frobenius endomorphism on $K^{\sep+}/bK^{\sep+}$ is 
an epimorphism. Using theorem \ref{th_Frobenius}(ii) we deduce 
that the Frobenius endomorphism on $(K^+/bK^+)^a$ is an 
epimorphism as well.
This shows that (i)$\Rightarrow$(iii). To show that 
(iii)$\Rightarrow$(ii), let us choose $\eps\in\fm\setminus\{0\}$
such that $|\eps^p|>|b|$; by hypothesis, for every $x\in K^+$
there exists $y\in K^+$ such that 
$\eps^p\cdot x-y^p\in bK^+$. 
It follows easily that the Frobenius endomorphism is surjective
on $K^+/(b\cdot\eps^{-p})K^+$. Replacing $b$ by $b\cdot\eps^{-p}$
we can assume that the Frobenius endomorphism is surjective on
$K^+/bK^+$. Let $b_1\in K^+$ such that $1>|b_1^p|\geq|b|$;
we let $\mathrm{Fil_1^\bullet}(K^+/pK^+)$ (resp. 
$\mathrm{Fil_2^\bullet}(K^+/pK^+)$) be the $b_1$-adic (resp. 
$b_1^p$-adic) filtration on $K^+/pK^+$. The group topology on 
$K^+/pK^+$ defined by the filtrations 
$\mathrm{Fil_i^\bullet}(K^+/pK^+)$ ($i=1,2$) is the same
as the one induced by the valuation topology of $K^+$;
moreover, one verifies easily that the Frobenius endomorphism
defines a morphism of filtered abelian groups 
$\mathrm{Fil_1^\bullet}(K^+/pK^+)\to
\mathrm{Fil_2^\bullet}(K^+/pK^+)$
and that the associated morphism of graded abelian groups is
surjective. It then follows from 
\cite[Ch.III, \S2, n.8, Cor.2]{BouAC} that (ii) holds.
Next suppose that (ii) holds; choose $b\in K^+$ such that
$1>|b|>|b^p|\geq|p|$; by hypothesis, the Frobenius endomorphism
on $K^+/b^pK^+$ is surjective; the same holds for the
Frobenius map on $K^{\sep+}/b^pK^{\sep+}$, in view of 
example \ref{ex_tens.with.Q}(iii). Hence, the assumptions
of lemma \ref{lem_oldstyle}(ii) are fulfilled, and we
deduce that 
$\L_{K^{\sep+}/K^+}\simeq\L_{K^{\sep+}[b^{-1}]/K^+[b^{-1}]}$.
Now, if $\chara(K)=p$, this implies already that 
$\Omega_{K^{\sep+}/K^+}\simeq\Omega_{K^\sep/K}=0$, which is 
(i). In case $\chara(K)=0$, we only deduce that
$\Omega_{K^{\sep+}/K^+}\simeq \Omega_{K^{\sep+}[1/p]/K^+[1/p]}$; 
however, $R:=K^+[1/p]$ is a valuation ring of residue 
characteristic $0$. We are therefore reduced to showing that $R$
is deeply ramified. Arguing as in the proof of proposition
\ref{prop_deep.ramif} we can even assume that $R$ is
strictly henselian, in which case the assertion follows
from corollary \ref{cor_firstblood}(i).

Furthermore, (ii) implies easily that 
$\Omega_{(K^+/b^pK^+)/\Z}=0$ (since $dx^p=p\cdot dx^{p-1}=0$). 
Let $I:=b^pK^+$; it follows that the natural map 
$I/I^2\to(K^+/b^pK^+)\otimes_\Z\Omega_{K^+/\Z}$
is surjective, {\em i.e.\/} 
$\Omega_{K^+/\Z}=b^p\cdot\Omega_{K^+/\Z}+K^+\cdot db^p\subset
b^p\cdot\Omega_{K^+/\Z}$, which implies (v). Next, by corollary
\ref{cor_boom}(ii), we have 
$\Omega_{K^+/\Z}(\log\Gamma^+)/\Omega_{K^+/\Z}\simeq
\kappa\otimes_\Z\Gamma$, and this last term vanishes since 
(ii) implies that $\Gamma=\Gamma^p$. This shows that 
(ii)$\Rightarrow$(iv) as well. Clearly (v)$\Rightarrow$(vi).
Suppose that (vi) holds. We will need the following :
\begin{claim}\label{cl_several.cls} $\Omega_{K^{\sep+}/\Z}$ is 
a $K^{\sep+}$-divisible module and 
$C:=\Coker(\rho_{K^{\sep+}/K^+})$ is a $K^+$-torsion module 
(notation of \eqref{subsec_fix-notation_rho}). Furthermore, 
$C^a\simeq(\Omega_{K^{\sep+}/K^+})^a$.
\end{claim}
\begin{pfclaim} In view of example \ref{ex_tens.with.Q}(iii), 
$(K^\sep,|\cdot|_{K^\sep})$ satisfies condition (iii), hence
the first assertion follows from the implications 
(iii)$\Leftrightarrow$(ii)$\Rightarrow$(v), which have already 
been shown. Furthermore, it is clear that 
$\Omega_{K^{\sep+}/K^+}$ is a torsion $K^+$-module and
therefore the second assertion follows easily from corollary
\ref{cor_boom}. The latter corollary also implies that 
$(\Omega_{K^+/\Z})^a\simeq\Omega_{K^+/\Z}(\log\Gamma^+)^a$, 
and similarly for $\Omega_{K^{\sep+}/\Z}$, whence the third 
assertion.
\end{pfclaim}

Now, suppose first that $\F_p\subset K^+$; in this case 
$\Omega_{K^{\sep+}/\Z}$ is a torsion-free $K^{\sep+}$-module
according to corollary \ref{cor_over-perfect-field}. Let 
$b\in K^+$ be any element; by theorem \ref{th_for_sep.closure} 
and snake lemma we deduce that the $b$-torsion submodule 
$C[b]^a:=\Ker(C^a\to C^a:x\mapsto b\cdot x)$ is isomorphic to 
the cokernel of the scalar multiplication by $b$ on the module 
$K^{\sep+}\otimes_{K^+}\Omega_{K^+/\Z}(\log\Gamma^+)$; the 
latter vanishes by assumption, and by claim \ref{cl_several.cls}
we have $C^a=\bigcup_{b\in K^+}C[b]^a$, whence $C^a=0$, which
is equivalent to (i) by claim \ref{cl_several.cls} and 
proposition \ref{prop_deep.ramif}.

Finally, in case $K^+$ is of mixed characteristic, denote by
$T$ (resp. $T'$) the $K^+$-torsion submodule of
$\Omega_{K^{\sep+}/\Z}$ (resp. of 
$K^{\sep+}\otimes_{K^+}\Omega_{K^+/\Z}$) and define
$T[b]$ (resp. $T'[b]$) as its $b$-torsion submodule, 
for any $b\in K^+$. The foregoing argument shows that $T^a$ is 
isomorphic to the $K^{+a}$-torsion submodule of 
$(\Omega_{K^{\sep+}/\Z}(\log\Gamma^+_{K^\sep}))^a$, and 
similarly for $(T')^a$; moreover, by snake lemma we obtain a
short exact sequence 
$0\to T'[b]^a\to T[b]^a\to C[b]^a\to 0$  
for every $b\in K^+$, whence a short exact sequence
$0\to (T')^a\to T^a\to C^a\to 0$.
Under (vi), $(T')^a$ is a divisible module; however,
it is clear from lemma \ref{lem_struct.T} that the only divisible
$(K^{\sep+})^a$-submodules of $T^a$ are $0$ and $T^a$.
Consequently, in light of claim \ref{cl_several.cls} and
proposition \ref{prop_deep.ramif},  in order to prove that 
(vi)$\Rightarrow$(i), it suffices to show that $(T')^a\neq 0$. 
In turns, this is implied by the following :
\begin{claim} The image in $\Omega_{K^+/\Z}(\log\Gamma)$ of 
$d\log(p)\in\Omega_{\Q^+/\Z}(\log\Gamma^+_\Q)$, has annihilator 
$pK^+$.
\end{claim}
\begin{pfclaim} (Of course, $\Q^+:=\bar\Q^+\cap\Q$). By theorem 
\ref{th_for_sep.closure} it suffices to consider the image of 
$d\log(p)$ in $\Omega_{K^{\sep+}/\Z}(\log\Gamma^+_{K^\sep})$.
Then, by theorem \ref{th_finally}, we reduce to consider the 
case $K^\sep=\bar\Q$. Then, once more by theorem
\ref{th_for_sep.closure}, it suffices to look at the 
annihilator of $d\log(p)$ in $\Omega_{\Q^+/\Z}(\log\Gamma^+_\Q)$ 
itself, and the claim follows.
\end{pfclaim}

Since (vi) is implied by both (iv) and (v), we deduce at once 
that all the conditions (i)-(vi) are equivalent. Furthermore, 
it is clear from claim \ref{cl_several.cls} and proposition 
\ref{prop_deep.ramif} that both (vii) and (viii) are equivalent
to (i), so the proposition follows.
\end{proof}

\begin{remark}\label{rem_already.done} 
By inspecting the proof of proposition
\ref{prop_deep.ramif-II}, we see that the argument for 
(ii)$\Rightarrow$(i) still goes through for valued 
fields $(K,|\cdot|)$ of arbitrary rank and characteristic 
$p>0$.
\end{remark}

\begin{lemma}\label{lem_consecutive} 
Let $(K,|\cdot|)$ be a valued field and $b\in K$ an element 
with $0<|b|<1$. Denote by $\fq(b)$ the radical of the ideal 
$bK^+$ and set $\fp(b):=\bigcap_{r>0}b^r\cdot K^+$. Then
$\fp(b)$ and $\fq(b)$ are {\em consecutive\/} prime
ideals, {\em i.e.} there are no prime ideals strictly
contained between $\fp(b)$ and $\fq(b)$. Equivalently,
the ring $W(b):=(K^+/\fp(b))_{\fq(b)}$ is a valuation ring 
of rank one and the image of $b$ is topologically nilpotent
in the valuation topology of $W(b)$.
\end{lemma}
\begin{proof} It is easy to verify that $\fp(b)$ and $\fq(b)$ 
are prime ideals, and using \eqref{subsec_bij.convex} one 
deduces that $W(b)$ is a valuation ring of rank one, 
which means that $\fp(b)$ and $\fq(b)$ are consecutive,
\end{proof}

\begin{theorem}\label{th_deep.ramificat} 
Let $(K,|\cdot|)$ be a valued field,
$(K^\wedge,|\cdot|^\wedge)$ its completion. The following 
conditions are equivalent :
\begin{enumerate}
\item
$(K,|\cdot|)$ is deeply ramified.
\item
For every valued extension $(E,|\cdot|_E)$ of $(K,|\cdot|)$,
for every $b\in K^+\setminus\{0\}$ and for every $i>0$ we have
$H_i(\L_{(E^+/bE^+)/(K^+/bK^+)})=0$.
\item
For every pair of consecutive prime ideals 
$\fp\subset\fq\subset K^+$, the valuation ring $(K^+/\fp)_\fq$ 
is deeply ramified.
\item
For every pair of convex subgroups 
$H_1\subset H_2\subset\Gamma$, the quotient $H_2/H_1$ is not 
isomorphic to $\Z$, and moreover, if $p:=\chara(\kappa)>0$, the 
Frobenius endomorphism on $K^{\wedge+}/pK^{\wedge+}$ is 
surjective.
\end{enumerate}
\end{theorem}
\begin{proof} To show that (ii)$\Rightarrow$(i), we take 
$E:=K^\sep$ and we choose a valuation on $K^\sep$ extending 
$|\cdot|$. Then, by arguing as in the proof of lemma
\ref{lem_oldstyle}(ii), we deduce from (ii) that the scalar 
multiplication by $b$ on $\Omega_{K^{\sep+}/K^+}$ is injective.
Since the latter is a torsion $K^{\sep+}$-module, we deduce 
(i). To show (i)$\Rightarrow$(ii), 
we reduce first to the case where $E=E^\mathrm{a}$; indeed, let 
$|\cdot|_{E^\mathrm{a}}$ be a valuation on $E^\mathrm{a}$
extending $|\cdot|_E$ and suppose
that the sought vanishing is known for the extension 
$(K,|\cdot|)\subset(E^\mathrm{a},|\cdot|_{E^\mathrm{a}})$; by 
transitivity, it then suffices to show :

\begin{claim}
$H_i(\L_{(E^{\mathrm{a}+}/bE^{\mathrm{a}+})/(E^+/bE^+)})=0$ 
for every $i>1$.
\end{claim}
\begin{pfclaim} By \cite[II.2.2.1]{Il} we have
$\L_{(E^{\mathrm{a}+}/bE^{\mathrm{a}+})/(E^+/bE^+)}\simeq
\L_{E^{\mathrm{a}+}/E^+}
\derotimes_{E^{\mathrm{a}+}}E^{\mathrm{a}+}/bE^{\mathrm{a}+}$, 
whence a spectral sequence
$$
E^2_{pq}:=
\Tor_p^{E^{\mathrm{a}+}}
(H_q(\L_{E^{\mathrm{a}+}/E^+}),E^{\mathrm{a}+}/bE^{\mathrm{a}+})
\Rightarrow 
H_{p+q}(\L_{(E^{\mathrm{a}+}/bE^{\mathrm{a}+})/(E^+/bE^+)}).
$$
Since $E^{\mathrm{a}+}/bE^{\mathrm{a}+}$ is an 
$E^{\mathrm{a}+}$-module of Tor-dimension 
$\leq 1$, we see that $E^2_{pq}=0$ for every $p>1$; furthermore,
by theorem \ref{th_gen.alg.case}(i), it follows that 
$E^2_{pq}=0$ whenever $p,q>0$, so the claim follows.
\end{pfclaim}

Thus, we can suppose that  $E$ is algebraically closed.
A spectral sequence analogous to the foregoing computes 
$H_i(\L_{(E^+/bE^+)/(K^{\mathrm{a}+}/bK^{\mathrm{a}+})})$,
and using theorems \ref{th_was.theor.three}(ii) and 
\ref{th_with.K-alg.cl}, we find that the latter vanishes 
for $i>0$. Consequently, by applying transitivity to the tower
of extensions $K^+/bK^+\subset K^{\mathrm{a}+}/bK^{\mathrm{a}+}
\subset E^+/bE^+$, we reduce to show the assertion for the case
$E=K^\mathrm{a}$. However, by example \ref{ex_tens.with.Q}(iii)
we have 
$K^{\mathrm{a}+}/bK^{\mathrm{a}+}\simeq K^{\sep+}/bK^{\sep+}$
for every $b\in K^+\setminus\{0\}$, so we can further reduce
to the case $E=K^\sep$. In this case, one concludes the proof
by another spectral sequence argument, this time using 
assumption (i) and theorem \ref{th_gen.alg.case}(ii) to show 
that the relevant terms $E^2_{pq}$ vanish.

To show that (iii)$\Rightarrow$(iv), we consider two subgroups
$H_1\subset H_2$ as in (iv); if $\mathrm{c.rk}(H_2/H_1)>1$,
then clearly $H_2/H_1$ cannot be isomorphic to $\Z$, so we can
assume that $H_1$ and $H_2$ are consecutive, so that the 
corresponding prime ideals are too (see \eqref{subsec_bij.convex}).
In this case, (iii) and proposition \ref{prop_deep.ramif} show 
that $H_2/H_1$ is not isomorphic to $\Z$, which is the first 
assertion of (iv). To prove the second assertion, it will suffice
to show the following :
\begin{claim}\label{cl_imagine} 
Suppose that $p:=\chara(\kappa)>0$ and that (iii) holds. 
Then, for every $b\in K^+\setminus\{0\}$ with $1>|b|\geq|p|$, 
the Frobenius endomorphism on $K^+/bK^+$ is surjective.
\end{claim}
\begin{pfclaim} For such a $b$ as above, define $\fp(b)$, 
$\fq(b)$ and $W(b)$ as in lemma \ref{lem_consecutive}; then 
$W(b)$ is a valuation ring of rank one, so it is deeply 
ramified by assumption (iii). Then, by proposition 
\ref{prop_deep.ramif-II} it follows that the Frobenius 
endomorphism is surjective on 
$W(b)/b\cdot W(b)\simeq K^+_{\fq(b)}/bK^+_{\fq(b)}$. 
We remark that $bK^+_{\fq(b)}\subset K^+$; there 
follows a natural imbedding: 
$K^+/bK^+_{\fq(b)}\subset W(b)/b\cdot W(b)$, commuting
with the Frobenius maps. It is then easy to deduce 
that the Frobenius endomorphism is surjective on 
$K^+/bK^+_{\fq(b)}$. Moreover, by proposition 
\ref{prop_deep.ramif}, the valuation of $W(b)$ is not 
discrete, hence its value group is isomorphic to a 
dense subgroup of $(\R,\geq)$ (see example 
\ref{ex_tens.with.Q}(vi)); therefore, by 
\eqref{subsec_bij.convex} and example 
\ref{ex_tens.with.Q}(v), we deduce that there exists 
an element $c\in K^+$ such that $|b|>|c^{3p}|$ and 
$|b|<|c^{2p}|$. These inequalities have been chosen 
so that $c^p K^+_{\fq(b)}\subset K^+$ and 
$bK^+_{\fq(b)}\subset c^{2p} K^+_{\fq(b)}$, 
whence $bK^+_{\fq(b)}\subset c^pK^+$, and finally
we conclude that the Frobenius endomorphism induces
a surjection: $K^+/cK^+\to K^+/c^pK^+$.
We let $\mathrm{Fil}^\bullet_1(K^+/bK^+)$ 
(resp. $\mathrm{Fil}^\bullet_2(K^+/bK^+)$)
be the $c$-adic (resp. $c^p$-adic) filtration on
$K^+/bK^+$. The foregoing implies that the
Frobenius endomorphism induces a morphism of filtered
modules $\mathrm{Fil}^\bullet_1(K^+/bK^+)\to
\mathrm{Fil}^\bullet_2(K^+/bK^+)$ which is
surjective on the associated graded modules; by
\cite[Ch.III, \S2, n.8, Cor.2]{BouAC} the claim follows.
\end{pfclaim}

Next, assume (iv) and let $W:=(K^+/\fp)_\fq$, for
two consecutive prime ideals $\fp\subset\fq\subset K^+$.
By assumption the Frobenius map is surjective on 
$K^+/bK^+$, whenever $b\in K^+\setminus\{0\}$ 
and $|b|\geq|p|$; we deduce easily that the Frobenius
endomorphism is surjective on $W/bW$, which implies (iii), 
in view of proposition \ref{prop_deep.ramif-II}.

(i)$\Rightarrow$(iii): indeed, let $\fp\subset\fq$
be as in (iii); we need to show that $(K^+/\fp)_\fq$
is deeply ramified. After replacing $K^+$ by $K^+_\fq$
we can assume that $\fq$ is the maximal ideal of $K^+$.
The ring $k^+:=K^+/\fp$ is a valuation ring; let 
$k:=\mathrm{Frac}(k^+)$, $|\cdot|_k$ the valuation on $k$
corresponding to $k^+$, and $(k',|\cdot|_{k'})$ a finite 
separable valued extension of $(k,|\cdot|_k)$.
It suffices to show that $\Omega_{k^{\prime+}/k^+}=0$.
We have $k'\simeq k[X]/(f(X))$ for some 
irreducible monic polynomial $f(X)\in k[X]$; let 
$\tilde f(X)\in K^+_\fp[X]$ be a lifting of $f(X)$ to
a monic polynomial. Then 
$E^+:=K^+_\fp[X]/(\tilde f(X))$ is the
integral closure of $K^+$ in the finite separable
extension $E:=\mathrm{Frac}(E^+)$ of $K$, and 
$E^+/\fp E^+\simeq k'$, so that $E^+$ is a valuation
ring, by lemma \ref{lem_simple.meth}. Furthermore, the 
preimage of $k^{\prime+}$ in $E^+$ is a valuation ring 
$R$ of $E$ with $R\cap K=K^+$ and 
$R/\fp R\simeq k^{\prime+}$. From (i) and
theorem \ref{th_gen.alg.case}(ii) we deduce that
$\Omega_{R/K^+}=0$, whence $\Omega_{k^{\prime+}/k^+}=0$
as required.

Finally we show that (iv) implies (i). We distinguish
several cases. The case when $p:=\chara(K)>0$ has already 
been dealt with, in view of remark \ref{rem_already.done}.
Next suppose that $\chara(\kappa)=0$; we will adapt the
argument given for the rank one case to prove corollary
\ref{cor_firstblood}(i). As usual, we reduce to the case
where $K$ is strictly henselian; it suffices to show that
$\Omega_{E^+/K^+}=0$ for every finite extension 
$(E,|\cdot|_E)$ of $K$. Then $E$ factors as a tower of
subextensions $E_0:=K\subset E_1\subset E_2\subset...\subset E_n:=E$
such that $E_{i+1}=E_i[b_i^{1/l_i}]$ for every $i=0,...,n-1$,
where $l_i:=[E_{i+1}:E_i]$ is a prime number and 
$b_i\in E_i$ such that $|b_i|\notin l_i\cdot\Gamma_{E_i}$.
It is easy to see that assumption (iv) is inherited
by every finite algebraic extension of $K$, hence we
can reduce to the case $E=K[b^{1/l}]$, with $b:=b_1$, $l:=l_1$.
One verifies as in the proof of proposition 
\ref{prop_firstblood}(i) that 
$E^+$ consists of the elements of the form 
$\sum_{i=0}^{l-1}x_i\cdot b^{i/l}$ such that 
$x_i\in K$ and $|x_i\cdot b^{i/l}|_E\leq 1$ for every 
$i=0,...,l-1$ and we have to show that 
$d(x_i\cdot b^{i/l})=0$ for every $i\leq l-1$. We may assume that
$i>0$, and up to replacing $b$ by $b^i\cdot x_i^l$, we can obtain 
that $b\in E^+$; we have then to verify that $db^{1/l}=0$. Define 
$\fp(b)$, $\fq(b)$ as in lemma \ref{lem_consecutive}, so that 
$\fp(b)$ and $\fq(b)$ are consecutive prime ideals, therefore 
$(E^+/\fp(b))_{\fq(b)}$ is a deeply ramified rank one valuation 
ring; in particular, its value group is not discrete. Then, using 
\eqref{subsec_bij.convex} and example 
\ref{ex_tens.with.Q}(v), we deduce that there exists 
an element $c\in K^+$ such that $|b|>|c^{l+1}|$ and
$|b|<|c^l|$. We can write $b=x\cdot c^l$ for some $x\in K^+$,
whence $db^{1/l}=c\cdot dx^{1/l}$. However, 
$|c|\leq|b|^{(l-1)/l(l+1)}\leq|a^{1/l}\cdot b^{-1}|^{l-1}=
|x|^{(l-1)/l}$. Since $x^{(l-1)/l}\cdot dx^{1/l}=0$, the claim
follows. Finally, suppose that $p:=\chara(\kappa)>0$ and
$\chara(K)=0$. Arguing as in the previous case, we produce
an element $b\in K^+$ such that $|b^p|>|p|$ and 
$|b^{p+1}|<|p|$. The Frobenius map is surjective on 
$K^+/b^pK^+$ by assumption, and on 
$K^{\sep+}/b^pK^{\sep+}$ by example \ref{ex_tens.with.Q}(iii), 
hence 
$\Omega_{K^{\sep+}/K^+}\simeq\Omega_{K^{\sep+}[1/p]/K^+[1/p]}$,
by lemma \ref{lem_oldstyle}(ii).
Now it suffices to remark that $K^+[1/p]$ is a valuation
ring with residue field of characteristic zero, so we are
reduced to the previous case, and the proof is concluded.
\end{proof}

\begin{remark}\label{rem_deeply.ramificat} 
By inspection of the proof, it is easy to check that condition 
(ii) of theorem \ref{th_deep.ramificat} is equivalent to the 
following. There exists a subset $S\subset K^+\setminus\{0\}$ 
such that the convex subgroup generated by $|S|:=\{|s|~|~s\in S\}$ 
equals $\Gamma_K$ and 
$H_i(\L_{E^+/K^+}\otimes_{K^+}K^+/s\cdot K^+)=0$ for every valued 
field extension $(E,|\cdot|_E)$ of $(K,|\cdot|)$, every $s\in S$ 
and every $i>0$.
\end{remark}


\newpage


\def\sA{{\text{\sf A}}}
\def\sT{{\text{\sf T}}}
\def\fV{{\mathfrak V}}
\def\fU{{\mathfrak U}}
\def\fX{{\mathfrak X}}
\def\fY{{\mathfrak Y}}
\def\fW{{\mathfrak W}}
\def\fZ{{\mathfrak Z}}
\def\fd{\mathfrak d}
\def\ff{\mathfrak f}
\def\fg{\mathfrak g}
\def\fh{\mathfrak h}
\def\fn{\mathfrak n}
\def\Spa{\mathrm{Spa}}
\def\Cont{\mathrm{Cont}}
\def\an{\mathrm{an}}
\def\na{\mathrm{na}}
\def\Exan{\mathrm{Exan}}


\section{Analytic geometry}\label{ch_analytic}
In this final chapter we bring into the picture $p$-adic
analytic geometry and formal schemes. The first three
sections develop a theory of the {\em analytic cotangent
complex\/} : we show how to attach a complex $\L^\an_{X/Y}$
to any morphism of locally finite type $\phi:X\to Y$ of formal
schemes or of R.Huber's adic spaces. This complex is obtained
via {\em derived completion\/} from the usual cotangent complex
of the morphism of ringed spaces underlying $\phi$. We prove
that $\L^\an_{X/Y}$ controls the {\em analytic\/} deformation
theory of the morphism $\phi$, in the same way as the usual
cotangent complex computes the deformations of the map of
ringed spaces underlying $\phi$. We hope the reader will agree
with us that these sections - though largely independent from the
rest of the monograph - are not misplaced, in view of the prominence
of the cotangent complex construction throughout our work.
Some of what we do here had been already anticipated in
an appendix of Andr\'e's treatise \cite[Suppl.(c)]{An}.

The main result of the remaining sections \ref{sec_Analytic_over_deep}
and \ref{sec_semicont} is a kind of weak purity statement valid for
affinoid varieties over a deeply ramified valued field of rank one
(theorem \ref{th_weak.purity}). The occurence of analytic
geometry in purity issues (and in $p$-adic Hodge theory at large) 
is rather natural; indeed, the literature on the subject is
littered with indications of the relevance of analytic varieties,
and already in \cite{Ta}, Tate explicitly asked for a $p$-adic 
Hodge theory in the framework of rigid varieties. We elect instead
to use the language of adic spaces, introduced by R.Huber in
\cite{Hu1}. Adic spaces are generalizations of Zariski-Riemann
spaces, that had already made a few cameo appearences in earlier
works on rigid analytic geometry. Recall that Zariski-Riemann
spaces were introduced originally by Zariski in his quest for
the resolution of singularities of algebraic varieties. Zariski's
idea was to attach to any (singular, reduced and irreducible) variety
$X$ defined over a field $k$, the ringed space
$(\tilde X,\cO_{\tilde X})$ defined as the projective limit of
the cofiltered system of all blow-up maps $X_\alpha\to X$; so
a point of $\tilde X$ is a compatible system
$\tilde x:=(x_\alpha\in X_\alpha~|~X_\alpha\to X)$. It is easy to
verify that, for any such $\tilde x$, the stalk $\cO_{\tilde X,\tilde x}$
is a valuation ring dominating the image of $\tilde x$ in $X$.
The strategy to construct a regular model for $X$ was broken up
in two stages : first one sought to show that for any point
$\tilde x\in\tilde X$ one can find some $x_\alpha$ under $\tilde x$
which is non-singular in the blow-up $X_\alpha$. This first stage goes
under the name of local uniformization; translated in algebraic terms, this
means that for every valuation ring $v$ of the field $k(X)$ of rational
functions on $X$, there is a model of $k(X)$ on which the center of $v$
is a non-singular point (a model is a reduced irreducible $k$-scheme
$Y$ of finite type such that $k(Y)=k(X)$). 
Local uniformization and the quasi-compactness of the Zariski-Riemann
space imply that there is a ``finite resolving system'', {\em i.e.}
a finite number of models such that every valuation of $k(X)$ has a
non-singular center on one of them (notice that the local uniformization
of a point $x_\alpha$ ``spreads around'' to an open neighborhood $U_\alpha$
of $x_\alpha$ in $X_\alpha$, so we achieve uniformization not just
for $\tilde x$, but for all the valuations contained in the preimage
$\tilde U$ of $U_\alpha$, which is open in $\tilde X$).

The second step is to try to reduce the number of models in a finite 
resolving system; restricting to open subvarieties one reduces this to
the question of going from a resolving system of cardinality two to
a non-singular model. This is what Zariski called the ``fundamental
theorem'' (see \cite[\S2.5]{Zar}). Zariski showed local uniformization
for all valuation rings containing the field $\Q$, but for the fundamental
theorem he could only find proofs in dimensions $\leq 3$.

Our theorem \ref{th_weak.purity} is directly inspired by Zariski's
strategy : rather than looking at the singularities of a variety,
we try to resolve the singularities of an \'etale covering $Y\to X$
of smooth affinoid adic spaces over a deeply ramified non-archimedean
field $K$ (so the map is singular only on the special fibre of a given
integral model defined on the valuation ring of $K$).
We assume that $X$ admits generically \'etale coordinates
$t_1,...,t_d\in\cO_X(X)$, and the role of Zariski's $\tilde X$ is played
by the projective system of all finite coverings of the form
$X_n:=X[t_1^{1/n},...,t_d^{1/n}]$, where $n$ ranges over the positive
integers. We have to show that the \'etale covering
$Y_n:=X_n\times_XY\to X_n$ becomes ``less and less'' singular (that is,
on its special fibre) as $n$ grows; so we need a numerical invariant
that quantifies the singularity of the covering : this is the
{\em discriminant\/} $\fd_{Y/X}$. The analogue of local uniformization
is our proposition \ref{prop_to.extremes}, whose proof uses the
results of section \ref{sec_deeply.ram}. Next, in order to exploit the
quasi-compactness of the affinoid adic spaces $Y_n$, we have to show
that the estimates furnished by proposition \ref{prop_to.extremes}
``spread around'' in an open subset; to this aim we prove that the
discriminant function is semicontinuous : this is the purpose of section
\ref{sec_semicont}. Notice that in our case we do not need an analogue
of Zariski's difficult fundamental theorem; this is because the discriminant
function always decreases when $n$ grows : what makes the problem of
resolution of singularities much harder is that a non-singular point
$x_\alpha$ on a blow up model $X_\alpha$ may be dominated by singular
points $x_\beta$ on some further blow up $X_\beta\to X_\alpha$
(so, for instance one cannot trivially reduce the cardinality of
a resolving system by forming joins of the various models of the
system). Besides, our spaces are not varieties, but adic spaces,
{\em i.e.} morally we work only ``at the level of the Zariski-Riemann
space'' and we do not need - as for Zariski's problem - to descend
to (integral) models of our spaces; so we are dealing exclusively
with valuation rings (or mild extensions thereof), rather than more
complicated local rings.

In essence, this is the complete outline
of the method; however, our theorem is - alas - much weaker than Faltings'
and does not yield by itself the kind of Galois cohomology vanishings
that are required to deduce comparison theorems for the cohomology of
algebraic varieties; we explain more precisely the current status of
the question in \eqref{subsec_explain}.

Throughout this chapter we fix a valued field $(K,|\cdot|)$ with
valuation of rank one, complete for its valuation topology.
As usual, $\fm$ denotes the maximal ideal of $K^+$.
We also let $a$ be a topologically nilpotent element
in $K^\times$.

\subsection{Derived completion functor}\label{subsec_der.compl}
Let $A$ be a complete $K^+$-algebra of topologically finite 
presentation. For any $A$-module $M$, we denote by $M^\wedge$ 
the (separated) $a$-adic completion of $M$.

\begin{proposition}\label{prop_general.coherence}
Let $A$ be as in \eqref{subsec_der.compl}.
\begin{enumerate}
\item 
Every finitely generated $A$-module which is torsion-free
as a $K^+$-module, is finitely presented.
\item
$A$ is a coherent ring.
\item
Let $N$ be a finitely presented $A$-module, $N'\subset N$
a submodule. Then there exists an integer $c\geq 0$ such that 
\set\begin{equation}\label{eq_Artin-Rees}
a^k N\cap N'\subset a^{k-c} N'
\end{equation}
for every $k\geq c$. In particular, the topology on $N'$
induced by the $a$-adic topology on $N$, agrees with the
$a$-adic topology of $N'$.
\item
Every finitely generated $A$-module is $a$-adically complete
and separated.
\item
Every submodule of a free $A$-module $F$ of finite type is 
closed for the $a$-adic topology of $F$.
\item
Every $A$-algebra of topologically finite type is separated.
\end{enumerate}
\end{proposition}
\begin{proof}
(i) is an easy consequence of \cite[Lemma 1.2]{Bosch}. To show 
(ii), one chooses a presentation  $A:=K^+\langle T_1,...,T_n\rangle/I$ 
for some finitely generated ideal $I$, and then reduces to prove
the statement for $K^+\langle T_1,...,T_n\rangle$, in 
which case it follows from (i). Next, let $N$, $N'$ be 
as in (iii) and define $T$ to be the $K^+$-torsion
submodule of $N'':=N/N'$; clearly $T$ is an $A$-submodule,
and the $A$-module $N''/T$ is $K^+$-torsion-free, therefore
is finitely presented by (i). Since $N$ is finitely generated,
this implies that $M:=\Ker(N\to N''/T)$ is finitely
generated. Hence, there exists an integer $c\geq 0$ such
that $a^c M\subset N'$. If now $k\geq c$, we have
$a^k N\cap N'\subset a^k N\cap M=a^k M\subset a^{k-c}N'$, 
which shows (iii). Next let us show:

\begin{claim}\label{cl_foregoing} 
Assertion (iv) holds for every finitely presented 
$A$-module.
\end{claim}
\begin{pfclaim}
Let $N$ be a finitely presented $A$-module and choose
a presentation $0\to K\to A^n\to N\to 0$. By (ii), $K$
is again finitely presented, and by (iii), the topology
on $K$ induced by the $a$-adic topology on $A^n$ coincides 
with the $a$-adic topology of $K$. Hence, after taking
$a$-adic completion, we obtain a short exact sequence :
$0\to K^\wedge\to A^n\to N^\wedge\to 0$ (see \cite[Th.8.1]{Mat}). 
It follows that the natural map $K\to K^\wedge$ is injective, 
which shows that the map $N\to N^\wedge$ is surjective for 
every finitely presented $A$-module $N$. In particular, this 
holds for $K$, whence $K\simeq K^\wedge$, and 
$N\simeq N^\wedge$, as claimed.
\end{pfclaim}

Finally, let $M$ be a submodule of $A^n$. By (iii), the topology
on $M$ induced by $A^n$ coincides with the $a$-adic topology.
Consequently, if $M$ is finitely presented, then $M$ is
complete for the $a$-adic topology by claim \ref{cl_foregoing}, 
hence complete as a subspace of $A^n$, hence closed in $A^n$. 
For an arbitrary $M$, define $\bar M:=\bigcup_{n>0}(M:a^n)$; 
then $\bar M$ is a submodule of $A^n$ and $A/\bar M$ is torsion-free 
as a $K^+$-module, so it is finitely presented by (i), therefore 
$\bar M$ is finitely presented by (ii). It follows that 
$a^c\bar M\subset M$ for some $c\geq 0$, whence
$a^kA^n\cap\bar M=a^kA^n\cap M$ for every $k\geq c$.
By the foregoing, $\bar M$ is complete, so $a^kA^n\cap\bar M$
is also complete, and finally $M$ is complete, hence closed.
This settles (v) and (iv) follows as well.
(vi) is an immediate consequence of (v).
\end{proof}

\begin{lemma}\label{lem_finite.pres.B} 
Let $A\to B$ be a map of $K^+$-algebras of topologically finite 
presentation. Then $B$ is of topologically finite presentation 
as an $A$-algebra.
More precisely, if $\phi:A\langle T_1,...,T_n\rangle\to B$
is any surjective map, $\Ker\,\phi$ is finitely generated.
\end{lemma}
\begin{proof} By proposition \ref{prop_general.coherence}(vi),
$B$ is complete and separated, hence we can find a surjective
map $\phi:A\langle T_1,...,T_n\rangle\to B$. It remains to
show that $\Ker\,\phi$ is finitely generated for any such $\phi$.
We can write $A:=K^+\langle T_{n+1},...,T_m\rangle/I$ for
some finitely generated ideal $I$, and thus reduce to the
case where $A=K^+$ and $\phi:K^+\langle T_1,...,T_n\rangle\to B$.
We will need the following :
\begin{claim} Let $\alpha:K^+\langle Y_1,...,Y_{r+s}\rangle\to B$ 
be a surjective map and $\beta:K^+\langle Y_1,...,Y_r\rangle\to
K^+\langle Y_1,...,Y_{r+s}\rangle$ the natural imbedding. Suppose
that $\gamma:=\alpha\circ\beta$ is surjective as well. Then
$\Ker\,\alpha$ is finitely generated if and only if $\Ker\,\gamma$
is finitely generated. 
\end{claim} 
\begin{pfclaim}\label{cl_ingarbuglia} For $i=r+1,...,r+s$, choose 
$f_i\in K^+\langle Y_1,...,Y_r\rangle$ such that 
$\gamma(f_i)=\alpha(Y_i)$. We define a surjective map 
$\delta:K^+\langle Y_1,...,Y_{r+s}\rangle\to 
K^+\langle Y_1,...,Y_r\rangle$ by setting $\delta(Y_i):=Y_i$ for
$i\leq r$ and $\delta(Y_i):=f_i$ for $i>r$. Clearly 
$\gamma\circ\delta=\alpha$. There follows a short exact sequence
$0\to\Ker\,\delta\to\Ker\,\alpha\to\Ker\,\gamma\to 0$. However,
$\Ker\,\delta$ is the closure of the ideal $I$ generated by 
$Y_i-f_i$ for $i=r+1,...,r+s$. By proposition 
\ref{prop_general.coherence}(v), we deduce that $\Ker\,\delta=I$,
and the claim follows easily.
\end{pfclaim}

By hypothesis there is at least one surjection 
$\psi:K^+\langle Y_1,...,Y_r\rangle\to B$ with finitely generated
kernel. Let $\mu:B\hat\otimes_{K^+}B\to B$ be the multiplication map
and set $\theta:=\mu\circ(\phi\hat\otimes_{K^+}\psi):
K^+\langle T_1,...,T_n,X_1,...,X_k\rangle\to B$. Applying twice claim 
\ref{cl_ingarbuglia} we deduce first that $\Ker\,\theta$ is finitely
generated, and then that $\Ker\,\phi$ is too, as required. 
\end{proof}

\begin{lemma}\label{lem_fat.compl.mods} 
Let $F$ be a flat $A$-module. Then:
\begin{enumerate}
\item
$F^\wedge$ is a flat $A$-module.
\item
For every finitely presented $A$-module $M$, the natural
map 
\set\begin{equation}\label{eq_comp.complet}
M\otimes_AF^\wedge\to(M\otimes_AF)^\wedge
\end{equation} 
is an isomorphism.
\end{enumerate}
\end{lemma}
\begin{proof} To begin with, we claim that the functor 
$N\mapsto(N\otimes_AF)^\wedge$ is exact on the abelian 
category of finitely presented $A$-modules. Indeed, let 
$\underline E:=(0\to N'\to N\to N''\to 0)$ be an exact 
sequence of finitely presented $A$-modules; we have to
show that $(\underline{E}\otimes_AF)^\wedge$ is still
exact. Obviously $\underline{E}\otimes_AF$ is exact, so
the assertion will follow by \cite[Th.8.1(ii)]{Mat}, once 
we know:
\begin{claim} 
The topology on $N'\otimes_AF$ induced by the 
imbedding into $N\otimes_AF$ agrees with the $a$-adic 
topology.
\end{claim} 
\begin{pfclaim} By proposition \ref{prop_general.coherence}(iii),
we can find $c\geq 0$ such that \eqref{eq_Artin-Rees} holds. 
Since $F$ is flat, we derive
$$
a^k(N\otimes_AF)\cap(N'\otimes_AF)\subset 
a^{k-c}(N'\otimes_AF)
$$
which implies the claim.
\end{pfclaim}

(ii): clearly \eqref{eq_comp.complet} is an 
isomorphism in case $M$ is a free module of finite type. 
For a general $M$, one chooses a resolution 
$\underline R:=(A^n\to A^m\to M\to 0)$; by the foregoing, 
the sequence $(\underline{R}\otimes_AF)^\wedge$ is still
exact, so one concludes by applying the 5-lemma to
the map of complexes $\underline{R}\otimes_AF^\wedge\to
(\underline{R}\otimes_AF)^\wedge$.

(i): we have to show that, for every injective map
of $A$-modules $f:N'\to N$, $f\otimes_A\one_{F^\wedge}$
is still injective. By the usual reductions, we can 
assume that both $N$ and $N'$ are finitely presented.
In view of (ii), this is equivalent to showing that the 
induced map 
$(N'\otimes_AF)^\wedge\to(N\otimes_AF)^\wedge$ is 
injective, which is already known.
\end{proof}

\sset\subsubsection{}
We will need to consider the left derived functor of
the $a$-adic completion functor, which we denote:
\set\begin{equation}\label{eq_notat.der.compl}
\sD^-(A\Mod)\to\sD^-(A\Mod)\quad:\quad
(K^\bullet)\mapsto(K^\bullet)^\wedge.
\end{equation}
As usual, it can be defined by completing termwise 
bounded above complexes of projective $A$-modules. 
However, the following lemma shows that it can also 
be computed by arbitrary flat resolutions.

\begin{lemma}\label{lem_fat.compl.cplx} 
Let $\phi:K^\bullet_1\to K^\bullet_2$ be 
a quasi-isomorphism of bounded above complexes of flat 
$A$-modules and denote by $(K_i^\bullet)^\wedge$ the 
termwise $a$-adic completion of $K_i^\bullet$ ($i=1,2$). 
Then the induced morphism 
\set\begin{equation}\label{eq_casotto}
(K_1^\bullet)^\wedge\to 
(K_2^\bullet)^\wedge
\end{equation}
is a quasi-isomorphism.
\end{lemma}
\begin{proof} Since $K^\bullet_1$ and $K^\bullet_2$ 
are termwise flat, we deduce quasi-isomorphisms
$$
\phi_n:K^\bullet_{1,n}:=K^\bullet_1\otimes_AA/a^n A\to 
K_{2,n}:=K^\bullet_2\otimes_AA/a^n A
$$
for every $n\in\N$. The map of inverse system of complexes 
$(K^\bullet_{1,n})_{n\in\N}\to(K^\bullet_{2,n})_{n\in\N}$ 
can be viewed as a morphism of complexes of objects of the 
abelian category $(A\Mod)^\N$ of inverse systems of 
$A$-modules. As such, it induces a morphism $(\phi_n)_{n\in\N}$
in the derived category $\sD((A\Mod)^\N)$, and it is clear 
that $(\phi_n)_{n\in\N}$ is a quasi-isomorphism.
Let 
$$
R\lim:\sD((A\Mod)^\N)\to\sD(A\Mod)
$$ 
be the right derived functor of the inverse limit functor
$\lim:(A\Mod)^\N\to A\Mod$. We remark that, for every
$j\in\Z$, the inverse systems $(K^j_{i,n})_{n\in\N}$ 
($i=1,2$) are acyclic for the functor $\lim$, since 
their transition maps are surjective. We derive that
$R\lim(K^\bullet_{i,n})_{n\in\N}\simeq(K_i^\bullet)^\wedge$, 
and, under this identification, the morphism \eqref{eq_casotto} 
is the same as $R\lim(\phi_n)_{n\in\N}$. Since the latter
preserves quasi-isomorphisms, the claim follows.
\end{proof}

\sset\subsubsection{}\label{subsec_essential.img}
\index{$\sD^-(A\Mod)^\wedge$|indref{subsec_essential.img}}
We denote by $\sD^-(A\Mod)^\wedge$ the essential image
of the functor \eqref{eq_notat.der.compl}.

\begin{corollary} {\em (i)\ } For any object $K^\bullet$ of\/
$\sD^-(A\Mod)$, the natural morphism 
$$
(K^\bullet)^\wedge\to((K^\bullet)^\wedge)^\wedge
$$
is a quasi-isomorphism.
\begin{enumerate}
\addtocounter{enumi}{1}
\item
$\sD^-(A\Mod)^\wedge$ is a full triangulated subcategory 
of $\sD^-(A\Mod)$.
\end{enumerate}
\end{corollary}
\begin{proof} Notice first that there are two natural morphisms
as in (i), which coincide : namely, for any complex $E$ in\/
$\sD^-(A\Mod)$ one has a natural morphism $u_E:E\to E^\wedge$;
then one can take either $(u_K)^\wedge$ or $u_{K^\wedge}$. Now,
(i) is an immediate consequence of lemmata \ref{lem_fat.compl.cplx} 
and \ref{lem_fat.compl.mods}. Clearly $\sD^-(A\Mod)^\wedge$ is 
preserved by shift and by taking cones of arbitrary morphisms; 
furthermore, it follows from (i) that it is a full 
subcategory of\/ $\sD^-(A\Mod)$.
\end{proof}

We will need some generalities on pseudo-coherent
complexes of $R$-modules (for an arbitrary ring $R$), 
which we borrow from \cite[Exp.I]{Bert}. In our situation, 
the definitions can be simplified somewhat, since we are 
only concerned with sheaves over the one-point site that 
are pseudo-coherent relative to the subcategory of free 
$A$-modules of finite type.

\sset\subsubsection{}\label{subsec_pseudo-coh}
\index{($n$-)pseudo-coherent complex(es)|indref{subsec_pseudo-coh}}
For given $n\in\Z$, one says that a complex $K^\bullet$ 
of $R$-modules is {\em $n$-pseudo-coherent\/} if there 
exists a quasi-isomorphism $E^\bullet\to K^\bullet$ where
$E^\bullet$ is a complex bounded above such that
$E^i$ is a free $R$-module of finite type for every
$i\geq n$. One says that $K^\bullet$ is {\em pseudo-coherent\/}
if it is $n$-pseudo-coherent for every $n\in\Z$.

\sset\subsubsection{}\label{subsec_der.cat.pseudo-coh}
\index{($n$-)pseudo-coherent complex(es)!$\sD(R\Mod)_{\mathrm{coh}}$ :
derived category of|indref{subsec_der.cat.pseudo-coh}}
Let $K^\bullet$ be a $n$-pseudo-coherent (resp. 
pseudo-coherent) complex of $R$-modules, and 
$F^\bullet\to K^\bullet$ a quasi-isomorphism. Then
$F^\bullet$ is $n$-pseudo-coherent (resp. pseudo-coherent)
(\cite[Exp.I, Prop.2.2(b)]{Bert}). It follows that
that the pseudo-coherent complexes form a (full) 
subcategory $\sD(R\Mod)_{\mathrm{coh}}$ of $\sD(R\Mod)$.

\sset\subsubsection{}\label{subsec_furthermore}
Furthermore, let $X\to Y\to Z\to X[1]$ be a distinguished
triangle in $\sD^-(R\Mod)$. If $X$ and $Z$ are 
$n$-pseudo-coherent (resp. pseudo-coherent), then the same
holds for $Y$ (\cite[Exp.I, Prop.2.5(b)]{Bert}).

\begin{lemma}\label{lem_one.of.those} Let $n,p\in\N$, 
$K^\bullet$ a $n$-pseudo-coherent complexes in 
$\sD^{\leq 0}(R\Mod)$, and $\cF_p$ one of the functors 
$\otimes^p_R$, $\Sym^p_R$, $\Lambda^p_R$, $\Gamma^p_R$ 
defined in \cite[I.4.2.2.6]{Il}. Then $L\cF_p(K^\bullet)$ 
is an $n$-pseudo-coherent complex.
\end{lemma}
\begin{proof} It is well known that $\cF_p$ sends 
free $R$-modules of finite type to free $R$-modules of finite 
type. It follows easily that the assertion of the lemma can be 
checked by inspecting the definition of the unnormalized chain 
complex associated to a simplicial complex, and of the simplicial 
complex associated to a chain complex via the Dold-Kan
correspondence. We omit the details.
\end{proof}

\sset\subsubsection{}
Let $K^\bullet$ be a pseudo-coherent complex. By
(\cite[Exp.I, Prop.2.7]{Bert}) there exists a quasi-isomorphism 
$E^\bullet\to K^\bullet$ where $E^\bullet$ is a bounded above 
complex of free $R$-modules of finite type.

\sset\subsubsection{}
Suppose now that $R$ is coherent; then we deduce easily
that a complex $K^\bullet$ of $R$-modules is pseudo-coherent 
if and only if $H^i(K^\bullet)$ is a coherent $R$-module
for every $i\in\Z$ and $H^i(K^\bullet)=0$ for every sufficiently
large $i\in\Z$ (\cite[Exp.I, Cor.3.5]{Bert}).
By proposition \ref{prop_general.coherence}(iv) it follows 
also that $\sD^-(A\Mod)_\mathrm{coh}\subset\sD^-(A\Mod)^\wedge$
for every $K^+$-algebra $A$ as in \eqref{subsec_der.compl}.

\sset\subsubsection{}
Let $A$ be as in \eqref{subsec_der.compl} and $M$ an 
$A$-module of finite presentation. We denote by $M[0]$ the 
complex consisting of the module $M$ placed in degree zero. 
Any finite presentation of $M$ can be extended to
a quasi-isomorphism $E^\bullet\to M[0]$, where $E^\bullet$
is a complex of free $A$-modules of finite type 
and $E^i=0$ for $i>0$ (\cite[Exp.I, Cor.3.5(a)]{Bert}). 
Together with proposition \ref{prop_general.coherence}(iv), 
it follows easily that the natural morphism 
$M[0]\to M[0]^\wedge$ is a quasi-isomorphism.

\begin{lemma}\label{lem_booga.booga} 
Let $A\to B$ be a map of complete $K^+$-algebras of
topologically finite presentation. Then, for every 
object $K^\bullet$ of $\sD^-(A\Mod)$, the natural morphism
$$
(K^\bullet\derotimes_AB)^\wedge\to
(K^{\bullet\wedge}\derotimes_AB)^\wedge
$$
is a quasi-isomorphism.
\end{lemma}
\begin{proof} We can suppose that $K^\bullet$ is a complex
of free $A$-modules. Then we are reduced to showing that, for
every free $A$-module $F$, the natural map 
$(F\otimes_AB)^\wedge\to(F^\wedge\otimes_AB)^\wedge$ is
an isomorphism. We leave this task to the reader.
\end{proof}

\begin{definition}\label{def_anal.cot.cplx}
\index{$\L_{B/A}^\an$ : Analytic cotangent complex for topological
algebras|indref{def_anal.cot.cplx}}
\index{$\Omega_{B/A}^\an$ : module of analytic
differentials for topological algebras|indref{def_anal.cot.cplx}{},
\indref{subsec_anal.omega}}
Let $\phi:A\to B$ be a map of complete $K^+$-algebras
of topologically finite presentation. The $B$-module of 
{\em analytic differentials relative to $\phi$} is
defined as $\Omega_{B/A}^\an:=\Omega_{B/A}^\wedge$.
The {\em analytic cotangent complex\/} of $\phi$ is the
complex $\L_{B/A}^\an:=(\L_{B/A})^\wedge$. Directly
on the definition we derive a natural isomorphism
\set\begin{equation}\label{eq_omega.an}
H_0(\L_{B/A}^\an)\simeq\Omega_{B/A}^\an
\end{equation}
Notice that $\L_{B/A}^\an$ is defined here via the standard resolution 
$P_A(B)\to B$, and it is therefore well defined as a {\em complex\/} 
of $B$-modules, not just as an object in the derived category
$\sD^-(B\Mod)^\wedge$. This will be essential in order to 
globalize the construction to formal schemes, in 
\eqref{subsec_global-form.sch}, and to adic spaces, in 
definition \ref{def_cot.cpxadicsp}.
\end{definition}

The following lemma will be useful in section 
\ref{sec_Analytic_over_deep}.

\begin{lemma}\label{lem_compare-oms} 
Let $A\to B$ be a continuous map of $K^+$-algebras
of topologically finite presentation. The natural map
$\phi_{B/A}:\Omega_{B/A}\to\Omega_{B/A}^\an$
is surjective with $K^+$-divisible kernel.
\end{lemma}
\begin{proof} One writes $B=B_0/I$ with 
$B_0:=A\langle T_1,...,T_n\rangle$ and $I\subset B_0$. 
Directly from the construction of $\Omega_{B_0/A}^\an$
one checks that $\phi_{B_0/A}$ is onto. Then,
a little diagram chasing shows that $\phi_{B/A}$ is
onto as well, and yields a surjective map 
$B\otimes_{B_0}\Ker\,\phi_{B_0/A}\to\Ker\,\phi_{B/A}$.
This allows to reduce to the case where $B=B_0$.
In this case, $\Ker\,\phi_{B/A}$ is generated by the
terms of the form 
$\delta(f):=df-\sum_{i=1}^n(\partial f/\partial T_i)\cdot dT_i$,
where $f$ ranges over all the elements of $B_0$.
For given $f\in B_0$, we can write $f=f_0+af_1$,
with $f_0\in A[T_1,...,T_n]$, $f_1\in B_0$. 
It follows easily that 
$\delta(f)=\delta(af_1)=a\cdot\delta(f_1)$, whence 
the claim.
\end{proof}
\begin{remark} In view of \eqref{eq_omega.an}, lemma 
\ref{lem_compare-oms} is also implied by the following 
more general observation. Let
$K_{B/A}:=\Cone(\psi_{B/A}:\L_{B/A}\to\L_{B/A}^\an)[1]$;
one has: $K_{B/A}\derotimes_BB/aB\simeq 0$. Indeed,
directly on the definition of $\L^\an_{B/A}$ one
sees that $\psi_{B/A}\derotimes_B\one_{B/aB}$ is an isomorphism.
\end{remark}

\begin{proposition}\label{prop_form.smooth} 
Let $\phi:A\to B$ be a map of complete $K^+$-algebras 
of topologically finite presentation, and suppose that 
$\phi$ is formally smooth for the $a$-adic topology. 
Then there is a natural quasi-isomorphism
$$
\L_{B/A}^\an\simeq\Omega_{B/A}^\an[0].
$$
\end{proposition}
\begin{proof} For every $n\in\N$, set $A_n:=A/a^n\cdot A$
and $B_n:=B/a^n\cdot B$. The hypothesis on $\phi$ implies
that $\phi_n:=\phi\otimes_A\one_{A_n}$ is of finite presentation
and formally smooth for the discrete topology, therefore
\set\begin{equation}\label{eq_Omega_B_n}
\L_{B_n/A_n}\simeq\Omega_{B_n/A_n}[0]\simeq
\Omega_{B/A}\otimes_AA_n[0]
\end{equation} 
for every $n\in\N$. Moreover, $\phi$ is flat by \cite[Lemma 1.6]{Bosch},
hence $\L_{B_n/A_n}\simeq\L_{B/A}\otimes_AA_n$.
On the other hand, for every $i\in\Z$ there is a short
exact sequence (cp. \cite[Th.3.5.8]{We})
$$
0\to\lim_{n\in\N}{}^1H^{i-1}(\L_{B/A}\otimes_AA_n)\to
H^i(\L_{B/A}^\an)\to
\lim_{n\in\N}H^i(\L_{B/A}\otimes_AA_n)\to 0.
$$
In view of \eqref{eq_Omega_B_n}, the inverse system
$(H^{i-1}(\L_{B/A}\otimes_AA_n))_{n\in\N}$
vanishes for $i\neq 1$ and has surjective transition
maps for $i=1$, hence its $\lim^1$ vanishes for every 
$i\in\Z$, and the claim follows easily.
\end{proof}

\begin{proposition}\label{prop_surj.pseudo.coh} 
Let $\phi:A\to B$ be a surjective map of complete 
$K^+$-algebras of topologically finite presentation.
Then $\L_{B/A}$ is a pseudo-coherent complex, in 
particular it lies in $\sD^-(B\Mod)^\wedge$ and 
$\L_{B/A}\simeq\L_{B/A}^\an$.
\end{proposition}
\begin{proof} First of all, notice that by lemma 
\ref{lem_finite.pres.B}, $B$ is of finite presentation, 
hence it is coherent as an $A$-module. Let $P:=P_A(B)$ 
be the standard simplicial resolution of $B$ by free 
$A$-algebras.  We obtain a morphism of simplicial 
$B$-algebras $\phi:B\otimes_AP\to B$ by tensoring with 
$B$ the augmentation $P\to B$ (here 
$B$ is regarded as a constant simplicial algebra).
By the foregoing, $P$ is pseudo-coherent, hence $P\otimes_AB$
lies in $\sD(B\Mod)_\mathrm{coh}$. Let $J:=\Ker\,\phi$. 
The short exact sequence of complexes
$0\to J\to P\otimes_AB\to B\to 0$ is split, therefore
$J$ is also pseudo-coherent.
Recall that we have natural isomorphisms: 
$J^i/J^{i+1}\stackrel{\sim}{\to}\Sym_B^i(\L_{B/A})$
for every $i\in\N$ (where $J^0:=B\otimes_AP$ and 
$\Sym^0_B(\L_{B/A}):=B$) (\cite[Ch.III,\S 3.3]{Il}). Furthermore, 
we have (see {\em loc.cit.}) :
\set\begin{equation}\label{eq_th.Quillen}
H_n(J^i)=0\quad\text{for every $n,i\in\N$ such that $i>n$}.
\end{equation}
We prove by induction on $n$ that $\L_{B/A}$
is $n$-pseudo-coherent for every $n\leq 1$. If $n=1$
there is nothing to prove. Suppose that the claim is
known for the integer $n$. It then follows by lemma
\ref{lem_one.of.those} that $J^i/J^{i+1}$ is 
$n$-pseudo-coherent for every $i>0$. However, it follows
from \eqref{eq_th.Quillen} that $J^i$ is $n$-pseudo-coherent 
as soon as $i>-n$. Hence, by \eqref{subsec_furthermore} (and 
an easy induction), we deduce that $J^i$ is 
$n$-pseudo-coherent for every $i\in\N$. Hence $J^2[1]$
is $(n-1)$-pseudo-coherent; if we now apply \eqref{subsec_furthermore}
to the distinguished triangle $J\to\L_{B/A}\to J^2[1]\to J[1]$,
we deduce that $\L_{B/A}$ is $(n-1)$-pseudo-coherent.
\end{proof}

\begin{theorem}\label{th_cot.is.coh} 
Let $A\to B\to C$ be maps of complete $K^+$-algebras of 
topologically finite presentation. Then:
\begin{enumerate}
\item
$\L_{B/A}^\an$ lies in $\sD^-(B\Mod)_\mathrm{coh}$.
\item
There is a natural distinguished triangle in $\sD^-(C\Mod)$ :
\set\begin{equation}\label{eq_dist.tri.analy}
C\otimes_B\L_{B/A}^\an\to\L_{C/A}^\an
\to\L_{C/B}^\an\to C\otimes_B\L_{B/A}^\an[1].
\end{equation}
\end{enumerate}
\end{theorem}
\begin{proof} (i): by lemma \ref{lem_finite.pres.B} we can find 
a surjection $B_0:=A\langle T_1,...,T_n\rangle\to B$ from a 
topologically free $A$-algebra onto $B$. If we apply transitivity 
to the sequence of maps $A\to B_0\to B$ and take the (derived) 
completion of the resulting distinguished triangle, we end 
up with the triangle:
$$
(B\otimes_{B_0}\L_{B_0/A})^\wedge\to
\L_{B/A}^\an\to\L_{B/B_0}^\an\to
(B\otimes_{B_0}\L_{B_0/A})^\wedge[1].
$$
We know already from proposition \ref{prop_surj.pseudo.coh} 
that $\L_{B/B_0}$ is pseudo-coherent, hence it coincides 
with $\L_{B/B_0}^\an$. Lemma \ref{lem_booga.booga}
yields a quasi-isomorphism: 
$(B\otimes_{B_0}\L_{B_0/A})^\wedge\stackrel{\sim}{\to}
(B\otimes_{B_0}\L_{B_0/A}^\an)^\wedge$; in view 
of proposition \ref{prop_form.smooth}, $\L_{B_0/A}^\an$
is a free $B_0$-module of finite rank in degree zero, in 
particular it is pseudo-coherent, so the same holds
for $(B\otimes_{B_0}\L_{B_0/A})^\wedge$, and taking into
account \eqref{subsec_furthermore}, the claim follows.

(ii): if we apply transitivity to the sequence of maps 
$A\to B\to C$, and then we complete the distinguished
triangle thus obtained, we obtain \eqref{eq_dist.tri.analy},
except that the first term is replaced by 
$(C\otimes_B\L_{B/A})^\wedge$, which we can also write
as $(C\otimes_B\L_{B/A}^\an)^\wedge$, in view
of lemma \ref{lem_booga.booga}. However, by (i), 
$\L_{B/A}^\an$ is pseudo-coherent, so 
it remains such after tensoring by $C$; in particular
$C\otimes_B\L_{B/A}^\an$ is already complete,
and the claim follows.
\end{proof}

\subsection{Cotangent complex for formal schemes and adic spaces}
In this section we show how to globalize the definition
of the analytic cotangent complex introduced in section 
\ref{subsec_der.compl}. We consider two kinds of globalization :
first we define the cotangent complex of a morphism $\ff:\fX\to\fY$
of formal schemes locally of finite presentation over $\Spf\,K^+$; 
then we will define the cotangent complex for a morphism
of adic spaces locally of finite type over $\Spa(K,K^+)$.

\begin{lemma}\label{lem_flat-stalks}
Let $\fX:=\Spf A$ be an affine formal scheme finitely
presented over $\Spf\,K^+$. For every $f\in A$, let
$\fD(f):=\{x\in\fX~|~f\notin\fm_x\}$. The natural map 
$A\to\Gamma(\fD(f),\cO_\fX)$ is flat.
\end{lemma}
\begin{proof} Since $\Gamma(\fD(f),\cO_\fX)$ is the
$a$-adic completion of $A_f$, the lemma follows from
lemma \ref{lem_fat.compl.mods}(i).
\end{proof}

\sset\subsubsection{}\label{subsec_global-form.sch}
Let $\ff:\fX\to\fY$ be a morphism of formal schemes locally 
of finite presentation over $\Spf\,K^+$, and suppose that 
$\fY$ is separated. For every affine open subset $U\subset\fX$, 
the small category $F_U$ of all affine open subsets $V\subset\fY$ 
with $\ff(U)\subset V$, is cofiltered under inclusion (or else 
it is empty). For every $V\in F_U$, $\cO_\fY(V)$ is a 
$K^+$-algebra of topologically finite presentation, hence the 
induced morphism $\cO_\fY(V)\to\cO_\fX(U)$ is of the kind 
considered in definition \ref{def_anal.cot.cplx}. We set 
$$
L(U/\fY):=\colim{V\in F_U^o}
\L^\an_{\cO_\fX(U)/\cO_\fY(V)}.
$$
\begin{definition}\label{def_cot.cplx.f-schemes}
\index{Formal scheme(s)!$\L^\an_{\fX/\fY}$ : Cotangent complex
of a morphism of|indref{def_cot.cplx.f-schemes}}
The mapping $U\mapsto L(U/\fY)$ defines a complex of presheaves
on a cofinal family of affine open subsets of $\fX$. By applying 
degreewise the construction of \cite[Ch.0, \S 3.2.1]{EGAI}, we 
can extend the latter to a complex of presheaves of $\cO_\fX$-modules 
on $\fX$. We define the {\em analytic cotangent complex\/} 
$\L^\an_{\fX/\fY}$ of the morphism $\ff:\fX\to\fY$ as the complex 
of sheaves associated to this complex of presheaves (this means 
that we form degreewise the associated sheaf, and we consider the 
resulting complex).
\end{definition}

\sset\subsubsection{}\label{subsec_notseparated}
More generally, if $\fY$ is not necessarily 
separated, we can choose an affinoid covering 
$\fY=\bigcup_{i\in I}\fU_i$ and the construction above
applies to the restrictions $\fV_i:=\ff^{-1}(\fU_i)\to\fU_i$;
since the definition of $\L^\an_{\fV_i/\fU_i}$ is local on
$\fV_i$, one can then glue them into a single cotangent complex 
$\L^\an_{\fX/\fY}$.

\begin{lemma}\label{lem_compcotcompl}
Let $\L_{\fX/\fY}$ denote the (usual) cotangent complex
of the morphism 
$$
(f,f^\sharp):(\fX,\cO_\fX)\to(\fY,\cO_\fY)
$$
of ringed spaces; there exists a natural map of complexes
\set\begin{equation}\label{eq_mapfromusual}
\L_{\fX/\fY}\to\L^\an_{\fX/\fY}
\end{equation}
inducing an isomorphism
\set\begin{equation}\label{eq_aftertime}
\L_{\fX/\fY}\otimes_{\cO_\fX}\cO_\fX/a^n\cO_\fX\stackrel{\sim}{\to}
\L^\an_{\fX/\fY}\otimes_{\cO_\fX}\cO_\fX/a^n\cO_\fX
\end{equation}
for every $n\in\N$.
\end{lemma}
\begin{proof} It suffices to construct \eqref{eq_mapfromusual}
in case $\fY$ is affine. According to \cite[Ch.II, (1.2.3.6)]{Il}
and \cite[Ch.II, (1.2.3.4)]{Il}, the complex $\L_{\fX/\fY}$ is
naturally isomorphic to the sheafification of the complex of
presheaves defined by the rule:
$$
U\mapsto\colim{V\in F_U^o}\L_{\cO_\fX(U)/\cO_\fY(V)}
$$
and then it is clear how to define \ref{eq_mapfromusual}.
From the construction it is obvious that \eqref{eq_aftertime}
is an isomorphism.
\end{proof}

It is occasionally important to know that
both $\L^\an_{\fX/\fY}$ and the morphism \eqref{eq_mapfromusual}
are well defined in the category of complexes of $\cO_\fX$-modules
(not just in its derived category).

\sset\subsubsection{}\label{subsec_anal-diffs}
\index{Formal scheme(s)!$\Omega_{\fX/\fY}^\an$ : sheaf of analytic relative
differentials of a morphism of|indref{subsec_anal-diffs}}
Furthermore, denote by $\Omega_{\fX/\fY}^\an$ the sheaf of {\em analytic
relative differentials} for the morphism $\ff$, which is defined as
$(\cI/\cI^2)_{|\fX}$, where $\cI$ is the ideal defining the
diagonal imbedding $\fX\to\fX\times_\fY\fX$. We see easily that 
there is a natural isomorphism
\set\begin{equation}\label{eq_compwithOmega}
H_0(\L^\an_{\fX/\fY})\simeq\Omega_{\fX/\fY}^\an.
\end{equation}

\begin{proposition}\label{prop_glob.coh.form} 
Let $\ff:\fX\to\fY$ be a morphism of formal schemes locally of finite 
type over $\Spf(K^+)$. Then:
\begin{enumerate}
\item
$\L^\an_{\fX/\fY}$ is a pseudo-coherent complex of $\cO_\fX$-modules.
\item
If\/ $\ff$ is a closed imbedding, then \eqref{eq_mapfromusual} is a 
quasi-isomorphism.
\item
If\/ $\ff$ is formally smooth, then \eqref{eq_compwithOmega} induces
a quasi-isomorphism: $\L^\an_{\fX/\fY}\stackrel{\sim}{\to}
\Omega_{\fX/\fY}^\an[0]$.
\end{enumerate}
\end{proposition}
\begin{proof} (i): according to \cite[Exp.I, Prop.2.1(b)]{Bert},
it suffices to show that $H_i(\L^\an_{\fX/\fY})$ is a 
coherent sheaf of $\cO_\fX$-modules for every $i\in\N$.
To this aim, let $U\subset\fX$ be an affine open subset
such that the family $F_U$ (notation of 
\eqref{subsec_global-form.sch}) is not empty; pick any 
$V\in F_U$. After replacing $\fX$ by $U$, we can suppose that
$U=\fX$. Set $A:=\Gamma(\fX,\cO_\fX)$ and let $L_i^\triangle$ be the 
sheaf of coherent $\cO_\fX$-modules associated to the coherent 
$A$-module $L_i:=H_i(\L^\an_{\cO_\fX(\fX)/\cO_\fY(V)})$ 
(cp. \cite[Ch.I, \S 10.10.1]{EGAI}, where this concept is discussed 
in the case of locally noetherian formal schemes).
The assertion will be an immediate consequence 
of theorem \ref{th_cot.is.coh}(i) and the following :
\begin{claim} There is a natural isomorphism of
$\cO_\fX$-modules:
$L_i^\triangle\stackrel{\sim}{\to}H_i(\L^\an_{\fX/\fY})$.
\end{claim}
\begin{pfclaim}  By the definition 
of $\L^\an_{\fX/\fY}$ we deduce a natural morphism of $\cO_\fX$-modules: 
$\alpha:L_i^\triangle\to H_i(\L^\an_{\fX/\fY})$. It therefore suffices to 
show that $\alpha$ induces an isomorphism on the stalks. To this aim, 
we remark first that the natural map $L_i\to H_i(L(\fX/\fY))$ is an 
isomorphism. Indeed, it suffices to consider another open subset 
$V'\in F_\fX$ with $V'\subset V$; we have 
$\L_{\cO_\fY(V')/\cO_\fY(V)}^\an\simeq 0$ by proposition 
\ref{prop_form.smooth}, and then it follows by transitivity 
(theorem \ref{th_cot.is.coh}(ii)) that the map
$L_i\to H_i(\L_{\cO_\fX(\fX)/\cO_\fY(V')}^\an)$ is an 
isomorphism.
More generally, this argument shows that, for every affine open 
subset $U'\subset\fX$, the natural map 
$H_i(\L^\an_{\cO_\fX(U')/\cO_\fY(V)})\to H_i(L(U'/\fY))$ 
is an isomorphism. However, on one hand we have 
$(L_i^\triangle)_x\simeq L_i\otimes_A\cO_{\fX,x}$.
On the other hand, we have (\cite[Ch.0, \S 3.2.4]{EGAI}) :
\set\begin{equation}\label{eq_H_i}
H_i(\L^\an_{\fX/\fY})_x\simeq\colim{x\in U'}H_i(L(U'/\fY))
\end{equation}
where the colimit ranges over the set $S$ of all affine open 
neighborhoods of $x$ in $\fX$. We can replace $S$ by the
cofinal subset of all open neighborhoods of the form
$\fD(f)$ (for $f\in A$ such that $f\notin\fm_x$).
Then, lemma \ref{lem_flat-stalks}, together with another easy 
application of transitivity allows to identify the right-hand 
side of \eqref{eq_H_i} with 
$H_i(L(\fX/\fY))\otimes_A\cO_{\fX,x}$, so (i) follows.
(ii) and (iii) are immediate consequences of proposition 
\ref{prop_surj.pseudo.coh} and respectively proposition 
\ref{prop_form.smooth}.
\end{pfclaim}
\end{proof}

\begin{proposition}\label{prop_transitformal} Let 
$\fX\stackrel{\ff}{\to}\fY\stackrel{\fg}{\to}\fZ$ be two 
morphisms of formal schemes locally of finite presentation 
over $\Spf\,K^+$. There is a natural distinguished triangle 
in $\sD^-(\cO_\fX\Mod)$
\set\begin{equation}\label{prop_global.transit}
L\ff^*\L^\an_{\fY/\fZ}\to\L^\an_{\fX/\fZ}\to\L^\an_{\fX/\fY}\to
L\ff^*\L^\an_{\fY/\fZ}[1].
\end{equation}
\end{proposition}
\begin{proof} As explained in \cite[Ch.II, \S 2.1]{Il}, for 
every sequence of ring homomorphisms $A\to B\to C$, the 
transitivity triangle is induced by a functorial exact sequence
of complexes $\L_{C/B/A}$ of flat $C$-modules (a "true triangle" 
in {\em loc. cit.}). Suppose now that $A,B,C$ are complete
$K^+$-algebras of topologically finite type; upon $a$-adic 
completion, one deduces a true triangle $\L^\wedge_{C/B/A}$.
Then, to every sequence of affine open subsets $U\subset\fX$,
$V\subset\fY$, $W\subset\fZ$ such that $\ff(U)\subset V$
and $\fg(V)\subset W$, one can associate the true triangle
$\L^\wedge_{\cO_\fX(U)/\cO_\fY(V)/\cO_\fZ(W)}$; 
Since the construction is functorial in all arguments, one
derives a presheaf of true triangles on a cofinal family
of open subsets of $\fX$, which we can then sheafify in the
usual manner. The resulting distinguished triangle in
$\sD^-(\cO_\fX\Mod)$ gives rise to \eqref{prop_global.transit}.
\end{proof} 

In the following we wish to define the cotangent complex
of a morphism of adic spaces (studied in \cite{Hu1} and \cite{Hu2}).
For simplicity, we will restrict to adic spaces of topologically
finite type over $\Spa(K,K^+)$, which suffice for our applications.
For the convenience of the reader we recall a few basic
definitions from \cite{Hu2}.

\sset\subsubsection{}\label{subsec_global-adic.spc}
\index{f-adic ring(s)|indref{subsec_global-adic.spc}}
\index{f-adic ring(s)!ring of definition of an|indref{subsec_global-adic.spc}}
\index{f-adic ring(s)!ideal of definition of an|indref{subsec_global-adic.spc}}
An {\em f-adic ring\/} is a topological ring $A$ that admits 
an open subring $A_0$ such that the induced topology on
$A_0$ is pre-adic and defined by a finitely generated ideal 
$I\subset A_0$. As an example, every $K$-algebra of topologically 
finite type is an f-adic ring. A subring $A_0$ with the 
above properties is called a {\em ring of definition\/} for 
$A$, and $I$ is an {\em ideal of definition}.
One denotes by $A^\circ$ the open subring of power-bounded 
elements of $A$.

\sset\subsubsection{}\label{subsec_top.fin.type}
\index{f-adic ring(s)! topologically finite type map of|indref{subsec_top.fin.type}}
Let $A\to B$ be complete f-adic rings and $\phi:A\to B$ a
ring homomorphism. One says that $\phi$ is {\em of topologically
finite type\/} if there exist rings of definition $A_0\subset A$
and $B_0\subset B$ such that $\phi(A_0)\subset B_0$, the restriction 
$A_0\to B_0$ factors through a quotient map ({\em i.e.\/} open
and surjective)  $A_0\langle T_1,...,T_n\rangle\to B_0$ and
$B$ is finitely generated over $A\cdot B_0$.

\sset\subsubsection{}\label{subsec_affinoid-ring}
\index{Affinoid ring(s)|indref{subsec_affinoid-ring}}
An {\em affinoid ring\/} is a pair $A=(A^\triangleright,A^+)$
consisting of an f-adic ring $A^\triangleright$ and a
subring $A^+\subset A^\triangleright$ which is open,
integrally closed in $A^\triangleright$ and contained
in the subring $A^\circ$. $A^+$ is called the {\em subring
of integral elements\/} of $A$. 

\sset\subsubsection{}\label{subsec_compl.aff.rng}
\index{Affinoid ring(s)!completion of|indref{subsec_compl.aff.rng}}
The {\em completion\/} $A^\wedge$ of an affinoid ring 
$A=(A^\triangleright,A^+)$ is the pair 
$((A^\triangleright)^\wedge,(A^+)^\wedge)$ (it turns out
that $(A^+)^\wedge$ is integrally closed in 
$(A^\triangleright)^\wedge$).

A homomorphism 
$\phi:(A^\triangleright,A^+)\to(B^\triangleright,B^+)$
of affinoid rings is a ring homomorphism 
$\phi^\triangleright:A^\triangleright\to B^\triangleright$ 
such that $\phi(A^+)\subset B^+$. One says that $\phi$ is 
{\em of topologically finite type\/} if $\phi^\triangleright$ 
is of topologically finite type and there exists an open
subring $C\subset B^+$ such that $B^+$ is the integral closure
of $C$, $\phi(A^+)\subset C$ and the induced map $A^+\to C$
is of topologically finite type. 
(cp. \eqref{subsec_top.fin.type}).

\sset\subsubsection{}\label{subsec_for.instance}
For instance, $(K,K^+)$ is an affinoid ring, complete for its
valuation topology; notice that in this case we have $K^+=K^\circ$.
A complete f-adic ring of topologically finite type over
$K$ is the same as a $K$-algebra of topologically finite
type. Furthermore, suppose that $A$ and $B$ are complete
f-adic rings of topologically finite type over $K$, and
let $f:A\to B$ be a continuous ring homomorphism. Then
there is a unique subring of integral elements $B^+$ such
that $f:(A,A^\circ)\to(B,B^+)$ is a morphism of affinoid 
rings of topologically finite type; namely one must take 
$B^+:=B^\circ$. Especially, $A^+:=A^\circ$ is the only
ring of integral elements of $A$ such that the affinoid
ring $(A,A^+)$ is of topologically finite type over $(K,K^+)$
(\cite[Prop.2.4.15]{Hu1}).

\sset\subsubsection{}
Given an arbitrary ring $A$, a {\em valuation\/} on $A$ is
a map $|\cdot|:A\to\Gamma\cup\{0\}$ where $\Gamma$
is an ordered abelian group whose composition law 
we denote multiplicatively, and the ordering is extended
to $\Gamma\cup\{0\}$ as usual. Then $|\cdot|$ is required
to satisfy the usual conditions, namely: 
$|x\cdot y|=|x|\cdot|y|$ and $|x+y|\leq\max(|x|,|y|)$
for every $x,y\in A$, and $|0|=0$, $|1|=1$.

\sset\subsubsection{}\label{subsec_def-Cont}
\index{f-adic ring(s)!$\Cont(A)$ : continuous valuations 
of an|indref{subsec_def-Cont}}
Now, let $A$ be an f-adic ring, and $|\cdot|:A\to\Gamma\cup\{0\}$
a valuation on $A$. For every $\gamma\in\Gamma$,
let $U_\gamma:=\{\alpha\in\Gamma~|~\alpha<\gamma\}\cup \{0\}$.
We endow $\Gamma\cup\{0\}$ with the topology which restricts 
to the discrete topology on $\Gamma$, and which admits 
$(U_\gamma~|~\gamma\in\Gamma)$ as a fundamental system of open 
neighborhoods of $0$. We say that $|\cdot|$ is {\em continuous\/}
if it is continuous with respect to the above topology on
$\Gamma\cup\{0\}$. One denotes by $\Cont(A)$ the
set of all (equivalence classes of) continuous valuations 
on $A$. Given $a,b\in A$, let $U(a/b)\subset\Cont(A)$ be the 
subset of all valuations $|\cdot|$ such that $|a|\leq|b|\neq 0$.
$\Cont(A)$ is endowed with the topology which admits 
the collection $(U(a/b)~|~a,b\in A)$ as a sub-basis.
With this topology, $\Cont(A)$ is a spectral
topological space (see \cite[1.1.13]{Hu2} for the definition
of spectral space). In particular, this implies that
$\Cont(A)$ admits a basis of quasi-compact open subsets.
Such a basis is provided by the {\em rational subsets},
defined as follows. A subset $U\subset\Cont(A)$ is called
rational if there exist $f_1,...,f_n,g\in A$ such that
the ideal $J:=f_1 A+...+f_n A$ is open in $A$ 
and $U$ consists of all $|\cdot|\in\Cont(A)$ such that 
$|f_i|\leq|g|\neq 0$ for every $i=1,...,n$. (Notice that, 
since we have chosen to restrict to f-adic rings containing 
$K$, asking for $J$ to be an open ideal is the same as 
requiring that $J=A$). Given $f_1,...,f_n,g\in A$ with
the above property, we denote by $R(f_1/g,...,f_n/g)$ 
the corresponding rational subset.

\sset\subsubsection{}\label{subsec_adic-spectrum}
\index{Affinoid ring(s)!$\Spa\,A$ : adic spectrum of 
an|indref{subsec_adic-spectrum}}
If $A:=(A^\triangleright,A^+)$, then one defines the
subset $\Spa\,A:=\{|\cdot|\in\Cont(A^\triangleright)~|~|a|\leq 1
\text{ for every $a\in A^+$}\}\subset\Cont(A^\triangleright)$. 
$\Spa\,A$, endowed with the subspace topology, is called 
the {\em adic spectrum\/} of the affinoid ring $A$.
$\Spa\,A$ is a pro-constructible subset of $\Cont(A)$,
hence it is a spectral space too.
Any continuous map $A\to B$ of affinoid rings induces
in the obvious way a continuous map on adic spectra: 
$\Spa\,B\to\Spa\,A$. 

\sset\subsubsection{}\label{subsec-adic.spaces}
\index{Adic space(s)|indref{subsec-adic.spaces}}
\index{Adic space(s)!affinoid|indref{subsec-adic.spaces}}
For any affinoid ring $A$, one can endow $X:=\Spa\,A$ with
a presheaf $\cO_X$ of topological rings, as follows. First of
all, for any $f_1,...,f_n,g\in A^\triangleright$ as in
\eqref{subsec_def-Cont}, one defines an affinoid ring 
$A(f_1/g,...,f_n/g)$, such that 
$A(f_1/g,...,f_n/g)^\triangleright:=(A^\triangleright)_g$ and 
$A(f_1/g,...,f_n/g)^+$ is the integral closure of the
subring $A[f_1/g,...,f_n/g]$ in $A(f_1/g,...,f_n/g)^\triangleright$. 
If $B\subset A^\triangleright$ is a ring of definition and 
$I\subset B$ an ideal of definition, let $B(f_1/g,...,f_n/g)$
be the subring of $(A^\triangleright)_g$ generated by $B$ and
$f_1/g$ ,..., $f_n/g$; we endow $B(f_1/g,...,f_n/g)$ with
the pre-adic topology defined by the ideal 
$I\cdot B(f_1/g,...,f_n/g)$; then the f-adic
topology on $A(f_1/g,...,f_n/g)^\triangleright$ is defined
to be the unique ring topology for which $B(f_1/g,...,f_n/g)$ 
is a ring of definition. Next, let 
$A\langle f_1/g,...,f_n/g\rangle:=A(f_1/g,...,f_n/g)^\wedge$
(cp. \eqref{subsec_compl.aff.rng}). With this preliminaries, 
one sets: 
$$
\cO_X(R(f_1/g,...,f_n/g)):=
A\langle f_1/g,...,f_n/g\rangle^\triangleright.
$$
In this way, $\cO_X$ is well defined on every rational
subset. One can then extend the definition to an arbitrary
open subset of $\Spa\,A$, following \cite[Ch.0, \S 3.2.1]{EGAI}.
It is not difficult to check that, for every open subset
$U\subset\Spa\,A$, and every $x\in U$, any valuation $|\cdot|_x$
in the equivalence class $x$ extends to the whole of $\cO_X(U)$, 
hence to the stalk $\cO_{X,x}$. One denotes by $\cO^+_X$ the
sub-presheaf defined by the rule: $\cO^+_X(U):=
\{f\in\cO_X(U)~|~|f|_x\leq 1\text{ for every $x\in U$}\}$.
In the cases of interest, the presheaf $\cO_X$ is a sheaf
(and $\cO^+_X$ is therefore a subsheaf). In such cases, one 
can show that, for every rational subset $R(f_1/g,...,f_n/g)$,
the natural map 
$A\langle f_1/g,...,f_n/g\rangle^+\to\cO^+_X(R(f_1/g,...,f_n/g))$
is an isomorphism of topological rings.

This holds notably when $A^\triangleright$ is a $K$-algebra of 
topologically finite type. One calls the datum 
$(\Spa\,A,\cO_{\Spa\,A},\cO^+_{\Spa\,A})$ an {\em affinoid adic 
space}.
General adic spaces are obtained as usual, by gluing affinoids.
Adic spaces form a category, whose morphisms $f:X\to Y$ are the 
morphisms of topologically locally ringed spaces 
$(X,\cO_X)\to(Y,\cO_Y)$ which induce morphisms of sheaves
$f^*\cO_Y^+\to\cO_X^+$.

\sset\subsubsection{}\label{subsec_fin.type.adic}
\index{Adic space(s)!morphism locally of finite type 
of|indref{subsec_fin.type.adic}}
Let $f:X\to Y$ be a morphism of adic spaces. One says that
$f$ is {\em locally of finite type\/} if for every 
$x\in X$ there exist open affinoid subspaces $U\subset X$,
$V\subset Y$ such that $x\in U$, $f(U)\subset V$ and the
induced morphism of affinoid rings 
$(\cO_Y(V),\cO^+_Y(V))\to(\cO_X(U),\cO^+_X(U))$ is of
topologically finite type.

\sset\subsubsection{}\label{subsec_smooth.etc-adic}
\index{Adic space(s)!{\'e}tale morphism of|indref{subsec_smooth.etc-adic}}
\index{Adic space(s)!smooth, unramified morphism of|indref{subsec_smooth.etc-adic}}
A morphism $f:X\to Y$ between adic spaces (defined over 
$\Spa(K,K^+)$) is called {\em smooth\/} (resp. {\em unramified},
resp. {\em {\'e}tale\/}) if $f$ is locally of finite type and if, 
for any affinoid ring $A$, any ideal $I$ of $A^\triangleright$ 
with $I^2=\{0\}$ and any morphism $\Spa\,A\to Y$, the mapping 
$\Hom_Y(\Spa\,A,X)\to\Hom_Y(\Spa\,A/I,X)$ is surjective (resp.
injective, resp. bijective).

\sset\subsubsection{}\label{subsec_functor.d}
\index{Adic space(s)!$d(\fX)$ : adic space associated to a formal 
scheme|indref{subsec_functor.d}}
In \cite[\S 1.9]{Hu2} it is shown how to associate functorially 
to every formal scheme $\fX$ (say locally of finite presentation 
over $\Spf\,K^+$) an adic space $d(\fX)$, together with a 
morphism of topologically ringed spaces $\lambda:d(\fX)\to\fX$, 
characterized by a certain universal property which we won't 
spell out here, but that includes the condition that 
$\Img(\cO_\fX\to\lambda_*\cO_{d(\fX)})\subset\cO_{d(\fX)}^+$.
If $\fX=\Spf A_0$ for a $K^+$-algebra $A_0$ of topologically finite
type, then $d(\fX)=\Spa\,A$, where $A$ is the affinoid ring
$(A_0\otimes_{K^+}K,A^+)$, with $A^+$ defined as the integral
closure of the image of $A_0$ in $A_0\otimes_{K^+}K$.
Moreover, $\fX$ is quasi-compact if and only if $d(\fX)$ is.

\sset\subsubsection{}\label{subsec_projlim.is.adic}
Let $\fX$ be a formal scheme of finite presentation over 
$\Spf\,K^+$. The collection $\cC_\fX$ of all morphisms 
$\ff:\fX'\to\fX$ of formal schemes of finite presentation 
over $\Spf\,K^+$ such that $d(\ff)$ is an isomorphism, 
forms a small cofiltered category (with morphisms given 
as usual by the commutative diagrams). It is shown in 
\cite[\S 3.9]{Hu1} that there is a natural isomorphism of 
topologically ringed spaces
\set\begin{equation}\label{eq_natural.ident}
(d(\fX),\cO^+_{d(\fX)})\stackrel{\sim}{\to}
\liminv{(\fX'\to\fX)\in\cC}(\fX',\cO_{\fX'}).
\end{equation}
(Actually, the argument in {\em loc.cit.} is worked out only
in the case of noetherian formal schemes, but it is not
difficult to adapt it to the present situation). 

\sset\subsubsection{}\label{subsec_cot.cplx.affrings}
Let $f:A\to B$ be a morphism of affinoid rings of topologically
finite type over $(K,K^+)$ (especially $A^+=A^\circ$ and
$B^+=B^\circ$, see \eqref{subsec_for.instance}). We let $\cC_f$ 
be the filtered family consisting of all the pairs $(A_0,B_0)$ of 
$K^+$-algebras of topologically finite presentation, such that
$A_0$ (resp. $B_0$) is an open subalgebra of $A^\circ$ 
(resp. of $B^\circ$) and $f(A_0)\subset B_0$. The {\em analytic
cotangent complex\/} of the morphism $f$ is the complex of
$B^+$-modules
$$
\L^+_{B/A}:=\colim{(A_0,B_0)\in\cC_f}\L^\an_{B_0/A_0}.
$$

\sset\subsubsection{}\label{}
Let $f:X\to Y$ be a morphism of adic spaces locally 
of finite type over $\Spa(K,K^+)$, and suppose that 
$Y$ is separated. For every affinoid open subset $U\subset X$, 
the small category $F_U$ of all affinoid open subsets $V\subset Y$ 
with $f(U)\subset V$, is cofiltered under inclusion (or else 
it is empty). For every $V\in F_U$, $(\cO_Y(V),\cO_Y(U)^+)$
is an affinoid $(K,K^+)$-algebra of topologically finite type,
hence the induced morphism $\cO_Y(V)\to\cO_X(U)$ is of the kind 
considered in \eqref{subsec_cot.cplx.affrings}. We set 
$$
L(U/Y):=\colim{V\in F_U^o}
\L^+_{\cO_X(U)/\cO_Y(V)}.
$$

\begin{definition}\label{def_cot.cpxadicsp}
\index{Adic space(s)!$\L^+_{X/Y}$, $\L^\an_{X/Y}$ : analytic
cotangent complex of a morphism of|indref{def_cot.cpxadicsp}}
The mapping $U\mapsto L(U/Y)$ defines a complex of presheaves
on a cofinal family of affinoid open subsets of $X$. By applying 
degreewise the construction of \cite[Ch.0, \S 3.2.1]{EGAI}, we 
can extend the latter to a complex of presheaves of $\cO_X^+$-modules 
on $X$. We define the {\em analytic cotangent complex\/} 
$\L^+_{X/Y}$ of the morphism $f:X\to Y$ as the complex 
of sheaves associated to this complex of presheaves (cp.
definition \eqref{def_cot.cplx.f-schemes}). We also set
$\L^\an_{X/Y}:=\L^+_{X/Y}\otimes_{K^+}K$.
The definition can be extended to the case of a morphism
$f:X\to Y$ where $Y$ is not necessarily separated: one argues 
as in \eqref{subsec_notseparated}.
\end{definition}

\sset\subsubsection{}
Let $X^+$ denote the ringed space $(X,\cO_X^+)$ and define likewise
$Y^+$. Just as in the proof of lemma \ref{lem_compcotcompl}, by 
inspecting the construction we obtain natural maps of complexes
\set\begin{equation}\label{eq_onXplus}
\L_{X^+/Y^+}\to\L^+_{X/Y}\qquad \L_{X/Y}\to\L^\an_{X/Y}.
\end{equation}

\sset\subsubsection{}\label{subsec_anal.omega}
Let $A$, $B$ be complete f-adic rings and $A\to B$ a ring
homomorphism of topologically finite presentation. We refer to 
\cite[\S 1.6]{Hu2} for the construction of a {\em universal 
$A$-derivation of $B$}, which is a continuous $A$-derivation
$d:B\to\Omega_{B/A}^\an$ from $B$ to a complete
topological $B$-module $\Omega^\an_{B/A}$, universal for 
$A$-derivations $B\to M$ to complete topological $B$-modules 
$M$. The construction of $\Omega^\an_{B/A}$ can be 
globalized to a sheaf of relative differentials $\Omega_{X/Y}^\an$ 
for any morphism of adic spaces $X\to Y$ locally of finite 
presentation. Then one checks easily that:
\set\begin{equation}\label{eq_omegas}
H_0(\L^\an_{X/Y})\simeq\Omega_{X/Y}^\an.
\end{equation}

\begin{proposition}\label{prop_transit.adic} Let 
$X\stackrel{f}{\to}Y\stackrel{g}{\to}Z$ be two 
morphisms of adic spaces locally of finite type over 
$\Spa(K,K^+)$. There is a natural distinguished triangle 
in $\sD^-(\cO_X\Mod)$
\set\begin{equation}\label{eq_transit.adic}
Lf^*\L^\an_{Y/Z}\to\L^\an_{X/Z}\to\L^\an_{X/Y}\to
Lf^*\L^\an_{Y/Z}[1].
\end{equation}
\end{proposition}
\begin{proof} {\em Mutatis mutandis}, this is the same as the
proof of proposition \ref{prop_transitformal}, so we can leave
the details to the reader.
\end{proof}

\begin{remark}\label{rem_transit.adic} In fact, the proof of proposition
\ref{prop_transit.adic} shows that \eqref{eq_transit.adic}
is represented by a {\em functorial\/} true triangle
(see \cite[Ch.I, \S3.2.4]{Il}). This will be important
in the sequel, when we will need to compute the truncation
$\tau_{[-1}\L_{X/Z}$ of the analytic cotangent in the situation
contemplated in proposition \ref{prop_truncateL}.
\end{remark}

\begin{theorem}\label{th_cot.smooth.vanish} 
Let $\ff:\fX\to\fY$ (resp. $f:X\to Y$) be a morphism of 
formal schemes (resp. of adic spaces) locally of finite 
presentation over $\Spf\,K^+$ (resp. locally of finite 
type over $\Spa(K,K^+)$) such that the 
induced morphism $d(\ff):d(\fX)\to d(\fY)$ is smooth
(resp. such that $f$ is smooth). Then:
\begin{enumerate}
\item
$\L^\an_{\fX/\fY}\otimes_{K^+}K\simeq
\Omega_{\fX/\fY}^\an[0]\otimes_{K^+}K$ 
in $\sD^-((\cO_\fX\otimes_{K^+}K)\Mod)$.
\item
$\L^\an_{X/Y}\simeq\Omega_{X/Y}^\an[0]$ in $\sD^-(\cO_X\Mod)$.
\end{enumerate}
\end{theorem}
\begin{proof} After the usual reductions, both assertions come
down to the following situation. We have a map of $K^+$-algebras 
of topologically finite presentation $\phi:A_0\to B_0$, such that
$d(\Spf\,\phi)$ is a smooth morphism of affinoid adic spaces.
We have to show that 
$\L_{B_0/A_0}^\an\otimes_{K^+}K\simeq
\Omega^\an_{B_0/A_0}\otimes_{K^+}K[0]$ in 
$\sD(B\Mod)$. We can write $B_0=C_0/I_0$, where 
$C_0:=A_0\langle T_1,...,T_n\rangle$ and $I_0$ is some 
finitely generated ideal. Set $I:=I_0\otimes_{K^+}K$,
$A:=A_0\otimes_{K^+}K$ and let $\fn$ be a maximal ideal 
of $C:=C_0\otimes_{K^+}K$ with $I\subset\fn$. 

\begin{claim}\label{cl_regular.seq} 
$I_\fn$ is generated by a regular sequence of elements 
of the local ring $C_\fn$.
\end{claim}
\begin{pfclaim} Let $\fp:=\fn\cap A$; $\fp$ is a maximal
ideal in $A$, and its residue field $K'$ is a finite extension
of $K$. Let $\bar\fn$ be the image of $\fn$ in 
$C\otimes_AK'$ and $\cI\subset\cO_{d(\Spf\,C_0)}$ the 
sheaf of ideals corresponding to $I$; the maximal ideal 
$\fn$ yields a point in $d(\Spf\,C_0)$, which we denote by 
$x(\fn)$. We have an isomorphism on the $\fn$-adic completions:
\set\begin{equation}\label{eq_on.n-adic.compl}
(\cI_{x(\fn)})^\wedge\simeq(I_\fn)^\wedge.
\end{equation}
Moreover, there are natural maps :
$$
\cI_{x(\fn)}/\cI_{x(\fn)}^2\to\bar\fn/\bar\fn{}^2\to
\Omega_{C/A}^\an\otimes_CC/\fn
$$
and, by \cite[Prop.2.5]{BoschIII}, there exists a set of
generators $g_1,...,g_k$ for $\cI_{x(\fn)}$ such that 
the images $dg_1,...,dg_k$ in 
$\Omega_{C/A}^\an\otimes_CC/\fn$ 
are linearly independent; it follows that the images
$\bar g_1,...,\bar g_k$ in $\bar\fn/\bar\fn{}^2$ are also linearly
independent.  Due to \eqref{eq_on.n-adic.compl}, we can
assume that $g_1,...,g_k\in I_\fn$, and then it follows that
$g_1\otimes 1,...,g_k\otimes 1$ are the first $k$ elements 
of a regular system of parameters for the regular local ring 
$C_\fn\otimes_AK'$. From lemma \ref{lem_fat.compl.mods} 
it follows that $C_\fn$ is a flat $A_\fp$-module; then by 
\cite[Ch.0, Prop.15.1.16]{EGAIV} we deduce that $g_1,...,g_n$ 
is a regular sequence of elements of $C_\fn$, as required.
\end{pfclaim}

Set $B:=B_0\otimes_{K^+}K$; it follows from claim 
\ref{cl_regular.seq} and \cite[Ch.III, Prop.3.2.4]{Il} that 
$C_\fn\otimes_C\L_{C/B}\simeq\L_{C_\fn/B_\fn}\simeq 
I_\fn/I_\fn^2[1]$ for every maximal ideal $\fn$, hence 
$\L_{C/B}\simeq I/I^2[1]$. By proposition \ref{prop_surj.pseudo.coh},
we derive $\L_{C_0/B_0}^\an\otimes_{K^+}K\simeq I/I^2[1]$.
Finally, by theorem \ref{th_cot.is.coh}(ii) and proposition 
\ref{prop_form.smooth} we deduce an isomorphism in 
$\sD^-(B\Mod)$ :
$$
\L_{B_0/A_0}^\an\otimes_{K^+}K\simeq(0\to I/I^2\to
B\otimes_{C_0}\Omega_{C_0/A_0}^\an\to 0)
$$ 
and the latter complex is quasi-isomorphic to 
$\Omega_{B/A}^\an[0]$ by \cite[Prop.2.5]{BoschIII}.
\end{proof}

\begin{lemma}\label{lem_closed.imbadic} 
Let $f:X\to Y$ be a closed imbedding of adic spaces locally
of finite type over $\Spa(K,K^+)$. Then the natural morphism
$$
\L_{X/Y}\to\L^\an_{X/Y}
$$
is an isomorphism in $\sD(\cO_X\Mod)$.
\end{lemma}
\begin{proof} 
The question is local on $X$, so we can reduce to the case
where $Y=\Spa(A,A^\circ)$ for some $K$-algebra of topologically
finite type and $X=\Spa(B,B^\circ)$, where $B=A/I$ for some 
ideal $I\subset A$. Let $\cC$ be the filtered family of all triples
of the form $(U,A_0,B_0)$ where $U\subset Y$ is an affinoid open 
neighborhood of $X$, $A_0\subset\cO_Y(U)^\circ$, $B_0\subset B^\circ$
are two open $K^+$-subalgebras of topologically finite presentation
with $\Img(A_0\to B)\subset B_0$.

\begin{claim} The family $\cC'\subset\cC_\pi$ of all
$(U,A_0,B_0)$ such that $\Img(A_0\to B)=B_0$ is cofinal.
\end{claim}
\begin{pfclaim} Let $(U,A_0,B_0)$ be any triple in $\cC$.
We can find finitely many elements $\bar b_1,...,\bar b_k\in B_0$
such that $A_0[\bar b_1,...,\bar b_k]$ is dense in $B_0$.
Choose elements $b_1,...,b_k\in A$ that lift the $\bar b_i$;
the subset $U':=\{x\in U~|~|b_i(x)|\leq 1;~i=1,...,k\}$
is an affinoid open neighborhood of $X$ and the topological
closure $A_0'$ of $A_0[b_1,...,b_k]$ is a $K^+$-algebra
of topologically finite presentation contained in $\cO_X(U')^\circ$;
moreover $\Img(A_0'\to B)=B_0$, hence $(U',A_0',B_0)\in\cC$.
\end{pfclaim}

Finally, let $(U,A_0,B_0)\in\cC'$; by proposition 
\ref{prop_surj.pseudo.coh} we have 
$\L^\an_{B_0/A_0}\simeq\L_{B_0/A_0}$. After taking colimits
we deduce:
$$
\colim{X\subset U\subset Y}\L_{B/\cO_Y(U)}^\an\simeq
\colim{X\subset U\subset Y}\L_{B/\cO_Y(U)}
$$
where $U$ runs over the cofiltered family of all open affinoid
neighborhoods of $X$ in $Y$. The lemma follows easily.
\end{proof}

\begin{proposition}\label{prop_truncateL}
Let $X\stackrel{f}{\to}Y\stackrel{g}{\to}Z$ be morphisms of 
adic spaces locally of finite type over $\Spa(K,K^+)$, where
$f$ is a closed imbedding with defining ideal $\cI\subset\cO_Y$, 
and $g$ is smooth. Then there is a natural isomorphism in 
$\sD(\cO_X\Mod)$
\set\begin{equation}\label{eq_truncateL}
\tau_{[-1}\L_{X/Z}^\an\stackrel{\sim}{\to}
(0\to(\cI/\cI^2)_{|X}\stackrel{d}{\to}
f^*\Omega_{Y/Z}^\an\to 0)
\end{equation}
where $d$ is the natural map.
\end{proposition}
\begin{proof} From lemma \ref{lem_closed.imbadic} and
\cite[Ch.III, Cor.1.2.8.1]{Il} we deduce a natural isomorphism
$$
\tau_{[-1}\L^\an_{X/Y}\simeq(\cI/\cI^2)_{|X}[1].
$$
Taking into account proposition \ref{prop_transit.adic},
theorem \ref{th_cot.smooth.vanish}, one can repeat the proof
of \cite[Ch.III, Cor.1.2.9.1]{Il}, which yields a distinguished
triangle:
\set\begin{equation}\label{eq_Illuyields}
f^*H_0(\L_{Y/Z})[0]\to\tau_{[-1}\L_{X/Z}\to 
H_1(\L_{X/Y})[1]\to
\end{equation}
such that the connecting morphism 
$(\cI/\cI^2)_{|X}\simeq H_1(\L_{X/Y})\to f^*H_0(\L_{Y/Z})
\simeq f^*\Omega_{Y/Z}^\an$ is naturally identified (up to
sign) with the differential map $f\mapsto df$. There follows
an isomorphism such as \eqref{eq_truncateL}; however,
the naturality of \eqref{eq_truncateL} is not explicitly
verified in {\em loc.cit.} : in general, this kind of manipulations,
when carried out in the derived category, do not lead to functorial
identifications. The problem is that we can have morphisms of
distinguished triangles:
$$
\xymatrix{ A^\bullet \ar[r] \ar[d]_f & B^\bullet \ar[r] \ar[d]_g & 
          C^\bullet \ar[r] \ar[d]^h & A^\bullet[1] \ar[d]^{f[1]}  \\
           A'{}^\bullet \ar[r] & B'{}^\bullet \ar[r] & C'{}^\bullet
          \ar[r] & A'{}^\bullet[1] 
}$$
such that both $f$ and $g$ are the zero maps, and yet $h$
is not. The issue can be resolved if one remarks that
\eqref{eq_Illuyields} is deduced from a transitivity triangle
via proposition \ref{prop_transit.adic}, and thus it is actually
well defined in the derived category of true triangles
$\sT(\cO_X\Mod)$ (cp. remark \ref{rem_transit.adic} and 
\cite[Ch.I, \S3.2.4]{Il}). It suffices then to apply the
following 
\begin{claim} Let $\cC$ be any abelian category,
$\underline T:=(0\to A^\bullet\to B^\bullet\to C^\bullet\to 0)$
a true triangle of $\sT(\cC)$ and $n\in\Z$ an integer such that
$H^{i+1}(A^\bullet)=0=H^i(C^\bullet)$ for every $i\neq n$. Then
there is a natural isomorphism in $\sT(\cC)$:
$$
\underline T\simeq(0\to H^{n+1}(A^\bullet)[n+1] \to
\Cone(H^n(C^\bullet)[n]\stackrel{d}{\to}H^{n+1}(A^\bullet)[n+1]) \to
H^n(C^\bullet)[n]\to 0)
$$
where $d$ is the connecting morphism of the long exact homology
sequence associated to $\underline T$.
\end{claim}
\begin{pfclaim} Under the stated assumptions, the natural
maps of complexes $\phi:\tau_{[n}C^\bullet\to C$ and 
$\psi:A^\bullet\to\tau_{n+1]}A^\bullet$ are quasi-isomorphisms;
we set $\underline T':=\psi*\underline T*\phi$ (notation of
\eqref{subsec_more.ext.stuff}. Clearly the pullback and push out
maps define a natural isomorphism 
$\underline T'\stackrel{\sim}{\to}\underline T$ in $\sT(\cC)$.
Say that 
$\underline T'=(0\to A'{}^\bullet\to B'{}^\bullet\to C'{}^\bullet\to 0)$;
since $H^{n+1}(C'{}^\bullet)=0$, it follows that the complex
$\underline T'':=\tau_{[n+1}\underline T':=(0\to\tau_{[n+1} A'{}^\bullet\to
\tau_{[n+1} B'{}^\bullet\to\tau_{[n+1} C'{}^\bullet\to 0)$ is
again a true triangle of $\sT(\cC)$, naturally isomorphic
to $\underline T'$. Likewise, $\underline T''$ is naturally 
isomorphic to the true triangle $\tau_{n+1]}\underline T''$, and
by inspection, one sees that the latter has the expected shape.
\end{pfclaim}
\end{proof}

\subsection{Deformations of formal schemes and adic spaces}
\label{sec_deformformsch}
\begin{lemma}\label{lem_meaning} Let $(\fX,\cO_\fX)$ be a formal
scheme locally of finite presentation over $\Spf\,K^+$, let $\cI$ be 
a coherent $\cO_\fX$-module and $(0\to\cI\to\cO_1\to\cO_\fX\to 0)$
an extension of sheaves of $K^+$-algebras on $\fX$.
Then $(\fX,\cO_1)$ is a formal scheme locally of finite
presentation over $\Spf\,K^+$.
\end{lemma}
\begin{proof} We may and do assume that $\fX$ is an affine 
formal scheme over $\Spf\,K^+$; then the assertion follows 
from claims \ref{cl_ext.topfin} and \ref{cl_locwell} below.
\begin{claim}\label{cl_ext.topfin}
$\cO_1(\fX)$ is a complete $K^+$-algebra of topologically 
finite presentation.
\end{claim}
\begin{pfclaim} Since $\fX$ is affine, we derive a short
exact sequence $0\to\cI(\fX)\to\cO_1(\fX)\to\cO_\fX(\fX)\to 0$.
Since $\cI$ is coherent, the $\cO_\fX(\fX)$-module $\cI(\fX)$
is finitely presented. Let $T\subset\cI(\fX)$ be the $K^+$-torsion 
submodule; by proposition \ref{prop_general.coherence}(i),
it follows that $\cI(\fX)/T$ is finitely presented, hence
$T$ is finitely generated over $\cO_\fX(\fX)$; likewise, the 
$K^+$-torsion submodule $T'$ of $\cO_\fX(\fX)$ is a finitely 
generated $\cO_\fX(\fX)$-module, and then the $K^+$-torsion 
submodule $T_1$ of $\cO_1(\fX)$ is a finitely generated 
$\cO_1(\fX)$-module. Let $N:=\Ker(\cO_1(\fX)\to\cO_\fX(\fX)/T')$; 
in the usual way we derive that there exists $k_0\in\N$ such that
$$
a^k\cO_1(\fX)\cap\cI(\fX)=a^kN\cap\cI(\fX)\subset 
a^{k-k_0}\cI(\fX)\qquad\text{for every $k\geq k_0$.}
$$
We deduce by \cite[Th.8.1(ii)]{Mat} a short exact sequence 
of complete $K^+$-algebras
$0\to\cI(\fX)\to\cO_1(\fX)^\wedge\to\cO_\fX(\fX)\to 0$,
which shows that $\cO_1(\fX)$ is complete. Since $\cI(\fX)$
and $\cO_\fX(\fX)$ are topologically of finite type, we
can then find a continuous surjection 
$A:=K^+\langle T_1,...,T_n\rangle\to\cO_1(\fX)$; by lemma 
\ref{lem_finite.pres.B} it follows that $\cO_1(\fX)$ is a 
finitely presented $A$-algebra, hence $\cI(\fX)$ is a 
finitely presented $A$-module, and then the same holds for
$\cO_1(\fX)$, so the claim is proved.
\end{pfclaim}

\begin{claim}\label{cl_locwell} For  $f\in\cO_\fX(\fX)$, let 
$\fD(f)\subset\fX$ be as in lemma \ref{lem_flat-stalks}. Then
the natural map $(\cO_1(\fX)_f)^\wedge\to\cO_1(\fD(f))$
is an isomorphism.
\end{claim}
\begin{pfclaim} The existence of the said map is a consequence
of claim \ref{cl_ext.topfin}; however, we have as well a natural
short exact sequence 
$0\to\cI(\fX)_f^\wedge\to\cO_1(\fD(f))\to\cO_\fX(\fX)_f^\wedge\to 0$,
so the assertion is immediate.
\end{pfclaim}
\end{proof}

The meaning of lemma \ref{lem_meaning} is that the square
zero deformations of $\fX$ by $\cI$ in the category of ringed
spaces are the same as those in the category of formal 
$\Spf\,K^+$-schemes locally of finite presentation. Especially,
the latter are classified by the appropriate $\Ext$-group of
the (usual) cotangent complex of the map of ringed spaces 
$\fX\to\Spf\,K^+$; we aim to show that the same computation 
can be carried out with the analytic cotangent complex 
$\L^\an_{\fX/K^+}$ introduced in definition \ref{def_cot.cplx.f-schemes}.

\sset\subsubsection{}\label{subsec_sometopoi}
Let $T$ be a topos, $I$ a small category; the category $T^I$
of all functors from $I$ to $T$ is a topos in a natural way,
and we define a functor $c:T\to T^I$ by assigning to an object
$X$ of $T$ the constant functor $c_X:I\to T$ of value $X$
(so $c_X(i)=X$ for every object $i$ of $I$, and $c_X(\phi)=\one_X$
for every morphism $\phi$ of $I$). The functor $c$ admits a
right adjoint $\liminv{I}:T^I\to T$, and the adjoint pair 
$(c,\liminv{I})$ defines a morphism of topoi $\pi_T:T^I\to T$.
If $\phi:S\to T$ is any morphism of topoi, we obtain a
commutative diagram of topoi:
$$
\xymatrix{ S^I \ar[r]^{\pi_S} \ar[d]_{\phi^I} & 
           S \ar[d]^\phi \\
           T^I \ar[r]^{\pi_T} & T
}$$
whence two spectral sequences:
$$
\begin{array}{r@{\;:=\;}l}
E^{pq}_2 & R^p\phi_*\, \liminv{i\in I}^q\cF_i 
          \Rightarrow R^{p+q}(\phi\circ\pi_S)_*(\cF_i~|~i\in I)\\
F^{pq}_2 & \liminv{i\in I}^p\, R^q\phi_*\cF_i
          \Rightarrow R^{p+q}(\phi\circ\pi_S)_*(\cF_i~|~i\in I)
\end{array}
$$
for every abelian sheaf $(\cF_i~|~i\in I)$ on $S^I$.

\begin{lemma}\label{lem_hilims0} 
Let $\fX$ be an adic formal scheme over $\Spf\,K^+$,
$\cF$ a coherent $\cO_\fX$-module and set $\cF_n:=\cF/a^n\cF$ for
every $n\in\N$. Then $(\cF_n~|~n\in\N)$ defines an abelian
sheaf on $\fX^\N$, and we have $\liminv{n\in\N}^q\cF_n=0$ on
$\fX$, for every $q>0$.
\end{lemma}
\begin{proof} We apply the spectral sequences of 
\eqref{subsec_sometopoi} to the indexing category $\N$ and
the morphism of topoi $\phi:\fU\to\Spf\,K^+$, where $\fU\subset\fX$
is any affine open subset. Hence $F^{pq}_2=0$ whenever $q>0$, 
therefore the abutment is isomorphic to 
$F^{p0}_2=\liminv{n\in\N}^p\,\cF_n(\fU)$; however, the inverse 
system $(\cF_n(\fU)~|~n\in\N)$ is surjective, hence 
$F^{p0}=0$ for $p>0$, and it follows that $E_\infty^{pq}=0$
whenever $p+q>0$ and
$E_\infty^{00}\simeq\liminv{n\in\N}\,\cF_n(\fU)\simeq\cF(\fU)$.
Since $\fU$ is arbitrary, the claim follows easily.
\end{proof}

\sset\subsubsection{}\label{subsec_trivialdual}
With the notation of \eqref{subsec_sometopoi}, let 
$\cO_\bullet:=(\cO_i~|~i\in I)$ be a ring object of the topos 
$T^I$ (briefly: a $T^I$-ring), or which is the same, a functor 
$\cO$ from $I$ to the category of $T$-rings. Let also $\cO_T$ be 
a $T$-ring, and $\pi_T^\sharp:\cO_T\to\pi_{T*}\cO_\bullet$
a morphism of $T$-rings. Then the pair $(\pi_T,\pi_T^\sharp)$ defines
a morphism of ringed topoi $\pi_T:(T^I,\cO_\bullet)\to(T,\cO_T)$. The
corresponding functor $\pi_{T*}=\liminv{I}:\cO_\bullet\Mod\to\cO_T\Mod$ 
admits a left adjoint $\pi_T^*$, defined by the rule: 
$$
\cF\mapsto\pi_T^{-1}\cF\otimes_{\pi_T^{-1}\cO_T}\cO_\bullet=
(\cF\otimes_{\cO_T}\cO_i~|~i\in I).
$$ 
(Of course, $\pi_T^{-1}$ is the same as the functor $c:T\to T^I$
of \eqref{subsec_sometopoi}).
The adjoint pair
$(\pi_{T*},\pi_T^*)$ extends to an adjoint pair of derived functors
$(R\pi_{T*},L\pi_T^*)$; more precisely, for every complex 
$K_1^\bullet\in\sD^-(\cO_T\Mod)$ and
$K_2^\bullet\in\sD^+(\cO_\bullet\Mod)$ there is a natural isomorphism
(``trivial duality")
$$
\Hom_{\sD(\cO_\bullet\Mod)}
(\pi_T^{-1}K_1^\bullet\derotimes_{\pi_T^{-1}\cO_T}\cO_\bullet,
K_2^\bullet)\simeq
\Hom_{\sD(\cO_T\Mod)}(K_1^\bullet,R\liminv{I}K_2^\bullet)
$$
(see \cite[Ch.III, Prop.4.6]{Il}).

\begin{proposition}\label{prop_topos.vs.anal}
Let $\fX$ be formal scheme locally of finite
presentation over $\Spf\,K^+$ and $\cF$ any coherent $\cO_\fX$-module. 
Then the natural morphism $\L_{\fX/K^+}\to\L^\mathrm{an}_{\fX/K^+}$
induces isomorphisms
$$
\hExt^i_{\cO_\fX}(\L^\an_{\fX/K^+},\cF)\stackrel{\sim}{\to}
\hExt^i_{\cO_\fX}(\L_{\fX/K^+},\cF)\qquad\text{for every $i\in\N$}.
$$
\end{proposition}
\begin{proof} We endow the topos $\fX^\N$ with the ring object 
$\cO_\bullet:=(\cO_\fX/a^n\cO_\fX~|~n\in\N)$; the natural morphism
$\cO_\fX\to\liminv{\N}\,\cO_\bullet$ determines a morphism of ringed
topoi $\pi:(\fX^\N,\cO_\bullet)\to(\fX,\cO_\fX)$ as in 
\eqref{subsec_trivialdual}. By lemma \ref{lem_hilims0}, the natural 
map $\cF\to R\liminv{\N}\,\pi^*\cF$ is an isomorphism and by lemma
\ref{eq_comp.complet}, both $\L_{\fX/K^+}$ and $\L^\an_{\fX/K^+}$ 
are complexes of flat $\cO_\fX$-modules, hence the trivial duality
isomorphism reads
$$
\hExt^i_{\cO_\bullet}(\pi^*\L^\an_{\fX/K^+},\pi^*\cF)\simeq
\hExt^i_{\cO_\fX}(\L^\an_{\fX/K^+},\cF)\qquad
\text{for every $i\in\N$}
$$
and likewise for $\L_{\fX/K^+}$. To conclude it remains only
to remark that $\pi^*\L_{\fX/K^+}\simeq\pi^*\L^\an_{\fX/K^+}$.
\end{proof}

Lemma \ref{lem_meaning} and proposition \ref{prop_topos.vs.anal}
enable one to derive the standard results on nilpotent deformations,
along the lines of section \ref{sec_lift} : the statements are 
unchanged, except that the topos-theoretic cotangent complex is
replaced by the (much more manageable) analytic one. The details
will be left to the reader.

\sset\subsubsection{}\label{subsec_andnowadic}
Let $f:X\to Y$ be a morphism of adic spaces locally of finite type
over $\Spa(K,K^+)$. We would like to show that the infinitesimal 
deformation theory of $f$ is captured by the analytic cotangent 
complex $\L^\an_{X/Y}$, in analogy with the usual topos-theoretic 
situation, and with the treatment for formal schemes already presented.
This turns out to be indeed the case, however the proofs are rather 
more delicate than those for formal schemes, due to the existence 
of square zero deformations of affinoid algebras that are not 
themselves affinoid. In homological terms, this reflects the fact
that the natural map
\set\begin{equation}\label{eq_compare.Ext1}
\hExt^1_{\cO_X}(\L_{X/Y}^\an,\cF)\to
\hExt^1_{\cO_X}(\L_{X/Y},\cF)
\end{equation}
is not in general an isomorphism for arbitrary coherent
$\cO_X$-modules $\cF$. Nevertheless, we have the following:

\begin{proposition}\label{prop_atleastinj} 
Let $f:X\to Y$ be as in \eqref{subsec_andnowadic}.
Then the map \eqref{eq_compare.Ext1} is injective for every 
coherent $\cO_X$-module $\cF$.
\end{proposition}
\begin{proof} Using the long exact $\hExt$ sequence, we reduce
to showing: 
\begin{claim}
$\hExt^0_{\cO_X}(\Cone\:(\L_{X/Y}\to\L_{X/Y}^\an),\cF)=0$.
\end{claim}
\begin{pfclaim} Indeed, we have
$$
\begin{array}{r@{\,\simeq\,}l}
\hExt^0_{\cO_X}(\Cone\:(\L_{X/Y}\to\L_{X/Y}^\an),\cF) &
\Ext^0_{\cO_X}(H_0(\Cone\:(\L_{X/Y}\to\L_{X/Y}^\an)),\cF) \\
& \Ext^0_{\cO_X}(\Coker(H_0(\L_{X/Y})\to H_0(\L_{X/Y}^\an)),\cF) \\
& \Ext^0_{\cO_X}(\Coker(\Omega_{X/Y}\to\Omega_{X/Y}^\an),\cF)
\end{array}
$$
and the latter group vanishes, due to the universal property
of $\Omega_{X/Y}^\an$ (see \eqref{subsec_anal.omega}). 
\end{pfclaim}
\end{proof}

\sset\subsubsection{}\label{subsec_analdefor}
\index{Adic space(s)!$\Exan_Y(X,\cF)$ : analytic deformations
of an|indref{subsec_analdefor}}
Let $f:X\to Y$ be a morphism of adic spaces locally of finite type 
over $\Spa(K,K^+)$. An {\em analytic deformation\/} of $X$ over $Y$ 
is a datum of the form $(j,\cF,\beta)$ consisting of :
\begin{enumerate}
\renewcommand{\labelenumi}{(\alph{enumi})}
\item
a closed imbedding $j:X\to X'$ of $Y$-adic spaces, such that 
the ideal $\cI\subset\cO_{X'}$ defining $j$ satisfies $\cI^2=0$, and
\item
a coherent $\cO_X$-module $\cF$ with an isomorphism of $\cO_X$-modules 
$\beta:j^*\cI\stackrel{\sim}{\to}\cF$.
\renewcommand{\labelenumi}{(\roman{enumi})}
\end{enumerate}
One defines in the obvious way a morphism of analytic deformations,
and we let $\Exan_Y(X,\cF)$ denote the set of isomorphism classes of
analytic deformations of $X$ by $\cF$ over $Y$.

\begin{proposition}\label{prop_deform.affin} 
Let $X$ be an affinoid adic space of finite type over $\Spa(K,K^+)$
and $(j:X\to X',\cF,\beta)$ an analytic deformation. Then $X'$ is 
affinoid of finite type over $\Spa(K,K^+)$.
\end{proposition}
\begin{proof} To start out, we notice that $j$ is a homeomorphism
on the topological spaces underlying $X$ and $X'$, hence 
$\cO_{X'}(X')=\Gamma(X,j^*\cO_{X'})$; we deduce a short exact
sequence of continuous maps:
$$
0\to\cF(X)\to\cO_{X'}(X')\to\cO_X(X)\to 0.
$$
\begin{claim}\label{cl_liftbdds} 
Let $\phi:B\to A$ be a surjective map of complete 
$K$-algebras of topologically finite type, such that $\Ker\,\phi$
is a square zero ideal of $B$. Then $B^\circ=\phi^{-1}(A^\circ)$.
\end{claim}
\begin{pfclaim} Quite generally, let $C$ be any $K$-algebra of
topologically finite type; according to \cite[\S6.2.3, Prop.1]{Bo-Gun}
one has $C^\circ=\{f\in C~|~|f(x)|\leq 1\text{ for every }
x\in\mathrm{Max}\, C\}$. In our situation, it is clear that 
$\mathrm{Max}\, A=\mathrm{Max}\, B$, whence the claim.
\end{pfclaim}

\begin{claim}\label{cl_globals}
$\cO_{X'}(X')$ is a complete $K$-algebra of topologically
finite type over $K$.
\end{claim}
\begin{pfclaim} Let $\bar g_1,...,\bar g_m$ be a set of topological
generators for $\cO_X(X)$, and choose arbitrary liftings 
$g_i,...,g_m\in\cO_{X'}(X')$. Pick also a finite set of
generators $g_{m+1},...,g_n$ for the $\cO_X(X)$-module $\cF(X)$.
We define a map $\phi:K[T_1,...,T_n]\to\cO_{X'}(X')$ by the rule
$T_i\mapsto g_i$ ($i=1,...,n$). Let $X'=\bigcup_kU_k$ be a
covering of $X'$ by finitely many of its affinoid domains.
We deduce a map $\psi:K[T_1,...,T_n]\to\prod_k\cO_{X'}(U_k)$.
By viewing $0\to j_*\cF\to\cO_{X'}\to j_*\cO_X\to 0$ as a short
exact sequence of coherent $\cO_{X'}$-modules on $X'$,
we deduce short exact sequences
$0\to\cF(U_k)\to\cO_{X'}(U_k)\stackrel{\pi_k}{\to}\cO_X(U_k)\to 0$ 
for every $k$. By the open mapping theorem (see
\cite[\S2.8.1]{Bo-Gun}) one deduces easily that the topology of
$\cO_X(U_k)$ is the same as the quotient topology deduced from
the surjection $\pi_k$, especially $\cO_X(U_k)$ is a $K$-algebra
of topologically finite type. Since the images of 
$\bar g_1,...,\bar g_m$ in $\cO_X(U_k)$ are power bounded for every
$k$, it then follows from claim \ref{cl_liftbdds} that the images
of $g_1,...,g_m$ in $\cO_{X'}(U_k)$ are power bounded for every $k$.
Hence $\psi$ extends to a map 
$\psi^\wedge:K\langle T_1,...,T_n\rangle\to\prod_k\cO_{X'}(U_k)$,
and by construction, $\psi^\wedge$ factors through a continuous map
$\phi^\wedge:K\langle T_1,...,T_n\rangle\to\cO_{X'}(X')$.
It is easy to check that $\phi^\wedge$ is surjective, so
the claim follows, again by the open mapping theorem.
\end{pfclaim}

Set $A:=\cO_{X'}(X')$; by claim \ref{cl_globals} we have the 
affinoid adic space $X'':=\Spa(A,A^\circ)$ of finite type
over $\Spa(K,K^+)$. By construction we get morphisms 
$X\stackrel{j}{\to}X'\stackrel{\alpha}{\to}X''$ of adic spaces
inducing homeomorphisms on the underlying topologies.
To conclude it suffices to show:
\begin{claim} The morphism 
$\alpha^\sharp:\cO_{X''}\to\alpha_*\cO_{X'}$ is an isomorphism
of sheaves of topological algebras.
\end{claim}
\begin{pfclaim} Using claim \ref{cl_HerrBosch} we see that
it suffices to show that $\alpha^\sharp$ is an isomorphism
of sheaves of $K$-algebras. However, for any affinoid open
domain in $X''$ we have a commutative diagram with exact
rows:
$$
\xymatrix{
0 \ar[r] & \cF(U) \ar[r] \ddouble & \cO_{X''}(U) \ar[r] \ar[d] &
\cO_X(U) \ar[r] \ddouble & 0 \\
0 \ar[r] & \cF(U) \ar[r] & \cO_{X'}(U) \ar[r] &
\cO_X(U) \ar[r] & 0.
}$$
Since the affinoid open subsets of $X''$ form a basis of 
the topology for both $X''$ and $X'$, the claim follows.
\end{pfclaim}
\end{proof}

\sset\subsubsection{}
In the situation of \eqref{subsec_analdefor}, let $(j:X\to X',\cF,\beta)$
be an analytic deformation of $X$ over $Y$. Clearly $j$ is a
homeomorphism on the underlying topological spaces, hence we
can view the given deformation as the datum of a map of
$f^{-1}\cO_Y$-algebras $\cO_{X'}\to\cO_X$, whence a transitivity
distinguished triangle (proposition \ref{prop_transit.adic})
$$
Lj^*\L^\an_{X'/S}\to\L^\an_{X/S}\to\L^\an_{X/X'}\to
Lj^*\L^\an_{X'/S}[1]
$$
which in turns yields a distinguished triangle:
$$
R\Hom_{\cO_X}(\L^\an_{X/X'},\cF)\to
R\Hom_{\cO_X}(\L^\an_{X/S},\cF)\to
R\Hom_{\cO_{X'}}(\L^an_{X'/S},j_*\cF)\to.
$$
Especially, we get a map
\set\begin{equation}\label{eq_combine1}
\Ext^1_{\cO_X}(\L_{X/X'}^\an,\cF)\to
\Ext^1_{\cO_X}(\L_{X/S}^\an,\cF).
\end{equation}
Let $\cI\subset\cO_{X'}$ be the ideal that defines the imbedding
$j$; by proposition \ref{prop_truncateL} we have a natural
isomorphism 
$\L_{X/X'}^\an\stackrel{\sim}{\to} j^*\cI[1]$, whence
an isomorphism
\set\begin{equation}\label{eq_combine2}
\Ext^1_{\cO_X}(\L_{X/X'}^\an,\cF)\stackrel{\sim}{\to}
\Hom_{\cO_X}(j^*\cI,\cF).
\end{equation}
Combining \eqref{eq_combine1} and \eqref{eq_combine2} we see
that the given isomorphism $\beta:j^*\cI\to\cF$ determines
a unique element $e^\an(X',\beta)\in\Ext^1_{\cO_X}(\L_{X/S}^\an,\cF)$.
One verifies easily that $e^\an(X',\beta)$ depends only on
the isomorphism class of the analytic deformation $(j,\cF,\beta)$,
therefore it defines a map
$$
e^\an:\Exan_Y(X,\cF)\to\Ext^1_{\cO_X}(\L_{X/S}^\an,\cF).
$$
By inspecting the construction, it is easy to check that
$e^\an$ fits into a commutative diagram
\set\begin{equation}\label{eq_fitsinto}
{\diagram
\Exan_Y(X,\cF) \ar[r]^{e^\an} \ar[d] & 
\hExt^1_{\cO_X}(\L_{X/S}^\an,\cF) \ar[d] \\
\Exal_{f^{-1}\cO_Y}(\cO_X,\cF) \ar[r]^e & 
\hExt^1_{\cO_X}(\L_{X/Y},\cF)
\enddiagram}
\end{equation}
where $e$ is the isomorphism of \cite[Ch.III, Th.1.2.3]{Il} and
the right vertical arrow is \eqref{eq_compare.Ext1}.

\begin{theorem}
For every coherent $\cO_X$-module $\cF$, the map $e^\an$
is a natural bijection.
\end{theorem}
\begin{proof} In view of \eqref{eq_fitsinto}, in order to show 
that $e^\an$ is injective, it suffices to verify that if
$(j_i:X\to X_i,\beta_i)$ for $i=1,2$ are two analytic deformations
of $X$ by $\cF$ over $Y$, and 
$(j_1,\beta_1)\stackrel{\sim}{\to}(j_2,\beta_2)$ is an
isomorphism in the category of extensions of $f^{-1}\cO_Y$-algebras,
then the corresponding map $\cO_{X_1}\to\cO_{X_2}$ is 
{\em continuous}, {\em i.e.} it is an isomorphism of sheaves of 
topological algebras (since in this case the map restricts to an
isomorphism $\cO_{X_1}^+\stackrel{\sim}{\to}\cO_{X_2}^+$). This 
can be checked locally, so we may assume that $X$ is affinoid;
in this case both $X_1$ and $X_2$ are affinoid as well, by proposition 
\ref{prop_deform.affin}. Then the assertion is a straightforward
consequence of the following well known:
\begin{claim}\label{cl_HerrBosch} 
Every map $A\to B$ from a noetherian $K$-Banach algebra
$A$ to $K$-algebra $B$ of topologically finite type, is continuous.
\end{claim}
\begin{pfclaim} This is \cite[\S6.1.3, Th.1]{Bo-Gun}.
\end{pfclaim} 

The surjectivity is a local issue as well; indeed, any class in 
the target of $e^\an$ represents an extension of $f^{-1}\cO_Y$-algebras
$0\to\cF\to\cE\to\cO_X\to 0$, and the question amounts to
showing that $\cE$ represents an analytic deformation, which
can be checked locally. Thus, suppose that $X=\Spa(B,B^\circ)$,
$Y=\Spa(A,A^\circ)$. We can write $B=P/I$, where 
$P:=A\langle T_1,...,T_n\rangle$, and by proposition 
\ref{prop_truncateL}, the complex $\tau_{[-1}\L_{X/Y}^\an$ is 
naturally isomorphic to the complex of sheaves associated to 
the complex of $B$-modules 
$0\to I/I^2\to\Omega_{P/A}^\an\otimes_PB\to 0$. Whence, a natural
isomorphism
$$
\begin{array}{r@{\;\simeq\;}l}
\hExt_{\cO_X}^1(\L_{X/Y}^\an,\cF) &
\Coker(\Hom_{\cO_X}(j^*\Omega_{Z/Y}^\an,\cF)\to
\Hom_{\cO_X}((I/I^2)^\sim,\cF)) \\
& \Coker(\Hom_B(\Omega_{P/A}^\an\otimes_PB,\cF(X))\to
\Hom_B(I/I^2,\cF(X)))
\end{array}
$$
(where we have denoted by $j:X\to Z:=\Spa(P,P^\circ)$ the induced
closed imbedding). However, given $\alpha:I/I^2\to\cF(X)$,
the corresponding $f^{-1}\cO_Y$ algebra $\cE$ is the sheaf of
algebras associated to the $A$-algebra $E$ defined as follows.
Let $\beta:I/I^2\to E_0:=(P/I^2)\oplus\cF(X)$ be the map
given by the rule: $x\mapsto(x,\alpha(x))$ for every $x\in I/I^2$;
$E_0$ is endowed with an $A$-algebra structure given by the rule 
$(x,f)\cdot(y,g)=(xy,fy+gx)$ and one verifies easily that the
image of $\beta$ is an ideal in $E_0$; then set $E:=E_0/\Img\,\beta$.
It is clear that $E_0$ is an $A$-algebra of topologically finite
type, and the projection $E\to B$ determines a closed imbedding
$j:X\to X':=\Spa(E,E^\circ)$ that identifies $j^*\cO_{X'}$ with
$\cE$. 
\end{proof}

With the foregoing results, one can derive the usual results
on existence of deformations and thereof obstructions; again
we leave the details to the industrious reader.
To conclude this section we want to show that \eqref{eq_compare.Ext1}
fails to be surjective already in the simplest situations.
To carry out this analysis requires the use of more refined
commutative algebra : the following proposition \ref{prop_regularity}
collects all the information that we shall be needing.

\begin{proposition}\label{prop_regularity} 
Let $A$ be a complete $K$-algebra of topologically finite type,
$n\in\N$ an integer, and set $P:=A\langle T_1,...,T_n\rangle$.
Then: 
\begin{enumerate}
\item
the natural morphism $\Spec\,P\to\Spec\,A$ is regular.
\item
$H_i(\L_{P/A})=0$ for every $i>0$ and $H_0(\L_{P/A})$ is a flat
$P$-module.
\end{enumerate}
\end{proposition}
\begin{proof} We begin with the following observation:
\begin{claim}\label{cl_doesnotexceed} 
Let $R\to S$ be a faithfully flat morphism
of noetherian local rings. If $S$ is regular, then $R$ is regular.
\end{claim}
\begin{pfclaim} One verifies with no trouble that the homological
dimension of $R$ does not exceed the homological dimension of $S$.
\end{pfclaim}

\begin{claim}\label{cl_Andreth} Let $R\to S$ be a local morphism 
of local noetherian rings; let $\fn\subset S$ be the maximal ideal.
Suppose that $R$ is quasi-excellent and that $S$ is formally
smooth over $R$ for its $\fn$-adic topology. Then the induced
morphism $\Spec\,S\to\Spec\,R$ is regular.
\end{claim}
\begin{pfclaim} This is \cite[Th.]{AnII}.
\end{pfclaim}

In view of claim \ref{cl_Andreth}, the proof of (i) is reduced to 
the following:

\begin{claim}\label{cl_excellence} 
Let $\fn\subset P$ be any maximal ideal and set $\fq:=\fn\cap A$. Then:
\begin{enumerate}
\item
$P_\fn$ is formally smooth over $A_\fq$ for its $\fn$-adic topology.
\item
$A_\fq$ is an excellent local ring.
\end{enumerate}
\end{claim} 
\begin{pfclaim} (i): It is well known that $\fq$ is a maximal
ideal and the residue field $K':=A/\fq$ is a finite extension
of $K$. It follows easily from lemma 
\ref{lem_fat.compl.mods}(i) that $P_\fn$ is flat over $A$,
hence by \cite[Ch.IV, Th.19.7.1]{EGAIV} it suffices to show
that $P_\fn\otimes_AK'$ is geometrically regular over $K'$, 
which allows to reduce to the case where $A=K$. Then, in view
of claim \ref{cl_doesnotexceed}, it suffices to show that,
for every finite field extension $K\subset K'$ there exists
a larger finite extension $K'\subset K''$ such that the semilocal
ring $P_\fn\otimes_KK''$ is regular. However, the residue
field $\kappa:=P/\fn$ is a finite extension of $K$, and if
$K\subset K'$ is any finite normal extension with $\kappa\subset K'$, 
every maximal ideal of $P\otimes_KK'$ containing $\fn\otimes_KK'$ 
is of the form $\fp:=(T_1-a_1,...,T_n-a_n)$ for certain $a_1,...,a_n\in K'$.
We are therefore reduced to the case where $\fn=(T_1-a_1,...,T_n-a_n)$
for some $a_i\in K$, $i=1,...,n$; in this case the $\fn$-adic
completion $P_\fn^\wedge$ of $P_\fn$ is isomorphic to 
$K[[X_1,...,X_n]]$, which is regular. (ii) is a theorem due to Kiehl, 
whose complete proof is reproduced in \cite[Th.1.1.3]{Conr}.
\end{pfclaim}

Finally, (ii) follows from (i) and from the following:
\begin{claim}\label{cl_AndreSupp} 
Let $R\to S$ be a map of noetherian rings. Then the following 
are equivalent:
\begin{enumerate}
\renewcommand{\labelenumi}{(\alph{enumi})}
\item
$H_i(\L_{S/R})=0$ for every $i>0$ and $H_0(\L_{S/R})$
is a flat $S$-module.
\item
The induced morphism $\Spec\,S\to\Spec\,R$ is regular.
\renewcommand{\labelenumi}{(\roman{enumi})}
\end{enumerate}
\end{claim}
\begin{pfclaim} Taking into account \cite[Th.28.7]{Mat}, this 
is seen to be a paraphrase of \cite[Suppl.(c), Th.30, p.331]{An}.
\end{pfclaim}
\end{proof}

\sset\subsubsection{}\label{subsec_def.Omega.na}
\index{$\Omega^\na_{P/K}$ : non-analytic
differentials|indref{subsec_def.Omega.na}}
Let us now specialize to the case where $P:=K\langle T\rangle$;
it is well known that $P$ is a principal ideal domain; since
$\Omega_{P/K}$ is a free $P$-module, we can choose a splitting: 
\set\begin{equation}\label{eq_choose.split}
\Omega_{P/K}\simeq\Omega_{P/K}^\an\oplus\Omega^\na_{P/K}
\end{equation}
where $\Omega^\na_{P/K}:=\Ker(\Omega_{P/K}\to\Omega_{P/K}^\an)$
(``na'' stays for ``not analytic''). Our first aim is to show that $P$
admits nilpotent extensions (by finitely generated $P$-modules) that
are not affinoid algebras. In view of proposition \ref{prop_regularity}(ii),
the square zero extensions of $P$ by a $P$-module $F$ are
classified by $\Ext^1_P(\Omega_{P/K},F)$ (and then it is
clear that any non-trivial element of this group cannot represent
an analytic deformation, since $P$ admits none). Hence, we come 
down to computing $\Ext^1_P(\Omega^\na_{P/K},F)$.

\begin{lemma}\label{lem_vectorspace}{\em (i) }
$\Omega^\na_{P/K}$ is a vector space over the fraction field 
$E$ of $P$.
\begin{enumerate}
\addtocounter{enumi}{1}
\item
If either $\chara(K)=0$ or $\chara(K)=p>0$ and $[K:K^p]=\infty$,
then $\Omega_{P/K}^\na\neq 0$.
\end{enumerate}
\end{lemma}
\begin{proof} (i): by proposition \ref{prop_regularity}(ii)
we know that $\Omega^\na_{P/K}$ is a flat $P$-module, especially
it is torsion-free. Let $f\in P$ be any irreducible element;
then $\kappa:=P/fP$ is a finite field extension of $K$,
and we have 
$\Hom_P(\Omega_{P/K},\kappa)=\Hom_P(\Omega_{P/K}^\an,\kappa)$
by the universal property of $\Omega_{P/K}^\an$. In other
words, $\Omega^\na\otimes_P\kappa=0$, so multiplication
by $f$ is a bijection on $\Omega^\na_{P/K}$, and the claim
follows.

(ii): suppose first that $\chara(K)=0$;
by (i) we have: $\Omega_{E/K}\simeq E\otimes_P\Omega_{P/K}\simeq
\Omega_{P/K}^\na\oplus(E\otimes_P\Omega_{P/K}^\an)$. Since
$\Omega_{P/K}^\an$ is a free $P$-module of rank one with generator
$dT$, and since $E$ is a separable extension of $K$, it suffices
to show that $\mathrm{tr.deg}(E:K(T))>0$. This can be checked
explicitly; for instance, let $\pi\in K$ be a non-zero element
such that $\log(1+\pi T)\in P$; it is well known that this power
series is transcendental over $K(T)$. Finally, suppose that
$\chara(K)=p>0$ and $[K:K^p]=\infty$. We construct a splitting
for the surjection $\Omega_{E/K}\to E\otimes_P\Omega_{P/K}^\an$,
as follows. The tower of field extensions 
$K\subset F:=E^p\cdot K(T)\subset E$ yields an exact sequence
$$
E\otimes_F\Omega_{F/K}\stackrel{\alpha}{\to}\Omega_{E/K}
\stackrel{\beta}{\to}\Omega_{E/F}\to 0
$$
and it is easy to check that $\alpha$ factors through a surjection
$E\otimes_F\Omega_{F/K}\to E\otimes_P\Omega_{P/K}^\an$ and the
induced map $E\otimes_P\Omega_{P/K}^\an\to\Omega_{E/K}$
is the sought splitting. Hence, it suffices to exhibit an
element $f\in E$ such that $\beta(df)\neq 0$. To this aim, we
will use the following general remark:
\begin{claim}\label{cl_Kerdiff} Let $E$ be a field with 
$\chara(E)=p>0$ and $F\subset E$ a subfield with $E^p\subset F$.
Then: $\Ker(d:E\to\Omega_{E/F})=F$.
\end{claim}
\begin{pfclaim} Let $(x_i~|~i\in I)$ be a $p$-basis for the
extension $F\subset E$, and for every $i\in I$ set 
$E_i:=F[x_i~|~i\in I\setminus\{i\}]$. Clearly $E=E_i[x_i]$
and the minimal polynomial of $x_i$ over $E_i$ is $m_i(X):=X^p-x^p_i$.
A simple calculation shows that $\Ker(d:E\to\Omega_{E/E_i})=E_i$.
Since $F=\bigcap_{i\in I}E_i$, it follows that 
$(d:E\to\Omega_{E/F})=\bigcap_{i\in I}\Ker(d:E\to\Omega_{E/E_i})=F$.
\end{pfclaim}

In view of claim \ref{cl_Kerdiff}, we need to exhibit an element
$f\in E$ such that $f\notin E^p\cdot K(T)$. However, under the
standing assumptions, we can find a countable sequence $(b_n~|~n\in\N)$
of elements $b_n\in K^\circ$ that are linearly independent over
$K^p$. Set $f:=\sum_{n\in\N}a^{pn}b_nT^n$, so that $f\in P$, and
suppose by way of contradiction, that $f\in E^p\cdot K(T)$; then
we could find $g\in E\cdot K^{1/p}(T^{1/p})$ such that $g^p=f$.
In turns, this means that $g\in K'\cdot E(T^{1/p})$ for some finite
$K$-extension $K'\subset K^{1/p}$. However, any such $g$ can be
written uniquely in the form $g=\sum_{n\in\N}c_nT^{n/p}$ for
a sequence $(c_n~|~n\in\N)$ of elements $c_n\in K'$. Then 
$g^p=\sum_{n\in\N}c_n^pT^n$, and consequently $c_n^p=a^{pn}b_n$
for every $n\in\N$, in other words, $K'$ contains all the elements
$b_n^{1/p}$ for every $n\in\N$, which is absurd.
\end{proof}

\sset\subsubsection{}\label{subsec_adeles}
We shall complete our calculation for $F:=P$. To this aim,
recall that $E$ is an injective $P$-module, therefore
$$
\Ext^1_P(\Omega^\na_{P/K},P)\simeq
\Hom_P(\Omega^\na_{P/K},E/P)/\Hom_P(\Omega^\na_{P/K},E).
$$
Let $\sA_E$ be the (finite) adele ring of $E$, {\em i.e.}
the restricted product of all the completions 
$$
\sA_E:=\sideset{}{'}\prod_{\fp\in\mathrm{Max}P}E_\fp.
$$
By the strong approximation theorem (cp. 
\cite[Ch.II, \S15, Th.]{Ca-Fr}) the natural imbedding 
$E\subset\sA_E$ has everywhere dense image, whence an isomorphism:
$$
E/P\simeq\mathop\oplus_{\fp\in\mathrm{Max}P}E_\fp/P_\fp.
$$
A simple calculation shows that the map
$E_\fp\to\Hom_P(E,E_\fp/P_\fp)$ given by the rule: 
$f\mapsto(x\mapsto xf)$ is an isomorphism; whence a natural
isomorphism:
$$
\sA_E\stackrel{\sim}{\to}\Hom_P(E,E/P)\qquad a\mapsto(x\mapsto ax).
$$
Finally, $\Ext^1_P(E,P)\simeq\sA_E/E$; 
since $\Omega^\na_{P/K}$ is a direct sum of copies of $E$, this
achieves our aim of producing non-trivial square zero extensions
of $K$-algebras: $0\to P\to\cE\to P\to 0$, whenever the hypotheses
of lemma \ref{lem_vectorspace}(ii) are fulfilled.

\sset\subsubsection{}\label{subsec_local-counterp}
Now we want to carry out the local counterpart of the calculations of
\eqref{subsec_adeles}. Namely, let $X:=\Spa(P,P^\circ)$; we will
show that there exist non-trivial deformations of $(X,\cO_X)$
by the $\cO_X$-module $\cO_X$, in the category of locally
ringed spaces. To this aim, choose any $K$-rational point $p\in X$,
and let $\cO_X(\infty)$ denote the quasi-coherent $\cO_X$-module
whose local sections on any open subset $U\subset X$ are the 
meromorphic functions on $U$ with poles (of arbitrary finite order)
only at $p$. We deduce an exact sequence
$$
\Hom_{\cO_X}(\Omega_{X/K}^\na,\cO_X(\infty))\stackrel{f}{\to}
\Hom_{\cO_X}(\Omega_{X/K}^\na,\cO_X(\infty)/\cO_X)\stackrel{g}{\to}
\Ext^1_{\cO_X}(\Omega_{X/K}^\na,\cO_X)
$$
where $\Omega_{X/K}^\na$ denotes the quasi-coherent sheaf obtained
by sheafifying the modules $\Omega^\na$ defined as in 
\eqref{subsec_def.Omega.na}. After choosing a global splitting
\eqref{eq_choose.split}, we can write 
$\Omega_{X/K}\simeq\Omega_{X/K}^\an\oplus\Omega_{X/K}^\na$.

\begin{lemma}\label{lem_local_counterp}
With the notation of \eqref{subsec_local-counterp}, we have:
$\Hom_{\cO_X}(\Omega_{X/K}^\na,\cO_X(\infty))=0$.
\end{lemma}
\begin{proof} Indeed, let $\phi:\Omega_{X/K}^\na\to\cO_X(\infty)$
be any $\cO_X$-linear map; we extend $\phi$ to a map
$\phi':\Omega_{X/K}\to\cO_X(\infty)$ by prescribing that $\phi'$
restricts to the zero map on the direct factor $\Omega_{X/K}^\an$.
Then $\phi'$ corresponds to a $K$-linear derivation 
$\partial:\cO_X\to\cO_X(\infty)$; the restriction of $\partial$ to
$U:=X\setminus\{p\}$ is a derivation of $\cO_U$ with values in the
coherent $\cO_U$-module $\cO_U$; hence the restriction of 
$\phi'_{|U}$ to $\Omega_{U/K}^\na$ vanishes identically.
It follows that $\partial_{|U}$ vanishes identically. Let $j:U\to X$
be the imbedding; since the natural map $\cO_X(\infty)\to j_*\cO_U$
is injective, it follows that $\partial$ must vanish as well,
whence $\phi=0$.
\end{proof}

\begin{lemma}\label{lem_stalk-notzero} 
With the notation of \eqref{subsec_local-counterp}, let $x\in X$
be a $K$-rational point and suppose that the hypotheses of lemma 
{\em\ref{lem_vectorspace}(ii)} hold for $K$. 
Then $\Omega_{X/K,x}^\na\neq 0$.
\end{lemma}
\begin{proof} By lemma \ref{lem_vectorspace}(ii) we know already
that the global sections of $\Omega_{X/K}^\na$ do not vanish.
It suffices therefore to show that, for every affinoid subdomain
$U:=\Spa(A,A^\circ)\subset X$ containing $x$, the natural map
$\Omega_{P/K}\otimes_PA\to\Omega_{A/K}$ is injective. However, the
kernel of this map is a quotient of $H_1(\L_{A/P})$, hence it
suffices to show the following
\begin{claim} $H_i(\L_{A/P})=0$ for every $i>0$.
\end{claim}
\begin{pfclaim} In light of claim \ref{cl_AndreSupp}, it suffices to
show that the map $P\to A$ is regular. Then, in view of claims
\ref{cl_Andreth} and \ref{cl_excellence}(ii), it suffices to show
that, for every maximal ideal $\fn\subset A$, the ring $A$ is formally
smooth over $P$ for its $\fn$-adic topology. Let $\fq:=\fn\cap P$;
since $A$ is flat over $P$, \cite[Ch.IV, Th.19.7.1]{EGAIV} reduces
to showing that the induced morphism $P/\fq\to A/\fq$ is formally
smooth, which is trivial, since the latter is an isomorphism.
\end{pfclaim}
\end{proof}

\sset\subsubsection{} It follows easily from proposition
\ref{prop_regularity}(ii) and lemma \ref{lem_flat-stalks} 
that $\L_{X/K}\simeq\Omega_{X/K}[0]$, so the extensions 
$0\to\cO_X\to\cE\to\cO_X\to 0$ are classified by 
$\Ext^1_{\cO_X}(\Omega_{X/K}^\na,\cO_X)$. In order
to show that the latter is not trivial, it suffices, 
in view of lemma \ref{lem_local_counterp}, to exhibit a nonzero 
map $\Omega_{X/K}^\na\to Q:=\cO_X(\infty)/\cO_X$.
However, $Q$ is a skyscaper sheaf sitting at the point
$p$, with stalk equal to $F/\cO_{X,p}$, where 
$F:=\mathrm{Frac}(\cO_{X,p})$. We are therefore reduced
to showing the existence of a nonzero map 
$\Omega_{X/K,p}^\na\to F$. However, we deduce easily from
lemma \ref{lem_vectorspace} that $\Omega_{X/K,p}^\na$ is an
$F$-vector space, and by lemma \ref{lem_stalk-notzero} the latter
does not vanish, provided that either $\chara(K)=0$ or
$K$ has characteristic $p>0$ and $[K:K^p]=\infty$.

\sset\subsubsection{}
The foregoing results notwithstanding, there is at least one
situation in which the abstract topos deformation theory of an
adic space is the same as its analytic deformation theory.
This is explained by the following final proposition.

\begin{proposition} Suppose that $\chara(K)=p>0$ and that
$[K:K^p]<\infty$. Then, for every morphism $X\to Y$ of adic
spaces locally of finite type over $\Spa(K,K^+)$, the natural
map $\L_{X/Y}\to\L^\an_{X/Y}$ is an isomorphism in $\sD(\cO_X\Mod)$.
\end{proposition}
\begin{proof} It suffices to show the corresponding statement
for maps of affinoid rings $A\to B$, hence choose a presentation
$B\simeq A\langle T_1,...,T_n\rangle/I$. We have to verify that
the induced maps $H_i(\L_{B/A})\to H_i(\L^\an_{B/A})$ are
isomorphisms for every $i\in\N$. The sequence 
$A\to P:=A\langle T_1,...,T_n\rangle\to B$ yields by transitivity
two distinguished triangles : one relative to the usual
cotangent complex and one relative to the analytic cotangent
complex. Hence, by the five lemma, we are reduced to showing
that the maps 
$H_i(B\otimes_P\L_{P/A})\to H_i(B\otimes_P\L^\an_{P/A})$
and $H_i(\L_{B/P})\to H_i(\L_{B/P}^\an)$ are isomorphisms for every
$i\in\N$. For the latter we appeal to lemma \ref{lem_closed.imbadic}.
Also, it follows easily from proposition \ref{prop_regularity}
and theorem \ref{th_cot.smooth.vanish} that both
$H_i(B\otimes_P\L_{P/A})$ and $H_i(B\otimes_P\L^\an_{P/A})$
vanish for $i>1$, hence we are reduced to showing that
the map 
\set\begin{equation}\label{eq_truncateall}
\tau_{[-1}\L_{B/A}\to\tau_{[-1}\L_{B/A}^\an
\end{equation}
is an isomorphism in $\sD(\cO_X\Mod)$. By proposition 
\ref{prop_truncateL} we have: 
$\tau_{[-1}\L_{B/A}^\an\simeq
(0\to I/I^2\to B\otimes_P\Omega_{P/A}^\an\to 0)$; on the other
hand, proposition \ref{prop_regularity} and 
\cite[Ch.III, Cor.1.2.9.1]{Il} give: 
$\tau_{[-1}\L_{B/A}\simeq(0\to I/I^2\to B\otimes_P\Omega_{P/A}\to 0)$.
Both these isomorphisms are natural, so \eqref{eq_truncateall}
is represented by the obvious map between the latter
complexes. Thus, we are reduced to showing that the
natural map $\Omega_{P/A}\to\Omega_{P/A}^\an$ is an
isomorphism. Let $C\subset P$ be the $A$-subalgebra generated
by the image of the Frobenius endomorphism $\Phi:P\to P$; a 
standard calculation shows that $\Omega_{P/A}\simeq\Omega_{P/C}$.
However, under the standing assumptions, $P$ is finite over its
subalgebra $C$, therefore $\Omega_{P/A}$ is a finite $P$-module,
and the claim follows.
\end{proof}

\subsection{Analytic geometry over a deeply ramified 
base}\label{sec_Analytic_over_deep}
\index{Adic space(s)!$\Omega^+_{X/K}$ : sheaf of analytic
differentials of an|indref{sec_Analytic_over_deep}}
In this section  we assume throughout that $(K,|\cdot|)$ 
is a deeply ramified complete valued field, with valuation 
of rank one. Recall that $a\in K^\times$ denotes a topologically 
nilpotent element of $K$.

If $\fX$ is a formal scheme of finite type over $\Spf\,K^+$, 
we will sometimes write $\L^\an_{\fX/K^+}$ instead of 
$\L^\an_{\fX/\Spf\,K^+}$. Similarly we define $\L^+_{X/K}$ 
for an adic space $X$ of finite type over $\Spa(K,K^+)$, and set 
$\Omega^+_{X/K}:=H_0(\L^+_{X/K})$.

\begin{theorem} Let $X$ be a smooth adic space over $\Spa(K,K^+)$. 
Then $\L^+_{X/K}\simeq\Omega^+_{X/K}[0]$ in $\sD^-(\cO_X^+\Mod)$, 
and $\Omega^+_{X/K}$ is a flat sheaf of $\cO^+_X$-modules.
\end{theorem}
\begin{proof} Both assertions can be checked on the stalks,
therefore let $x\in X$ be any point. The stalk $\cO_{X,x}$
is a local ring and its residue field $\kappa(x)$ carries
a natural valuation; the preimage in $\cO_{X,x}$ of the
corresponding valuation ring $\kappa(x)^+$ is the subring 
$\cO_{X,x}^+$. Let $I:=\bigcap_{n\in\N}a^n\cO^+_{X,x}$; it 
follows from this description that $\kappa(x)^+=\cO^+_{X,x}/I$.
Especially, we have
\set\begin{equation}\label{eq_separate-stalk}
\cO^+_{X,x}/a\cO^+_{X,x}\simeq\kappa(x)^+/a\cdot\kappa(x)^+.
\end{equation}
\begin{claim}\label{cl_when.a.regul} 
Let $M$ be an $\cO^+_{X,x}$-module, and suppose that $a$ is 
regular on $M$. Then $M/IM$ is a flat $\kappa(x)^+$-module. 
\end{claim}
\begin{pfclaim} By snake lemma we derive 
$\Ker(M/IM\stackrel{\cdot a}{\to}M/IM)\subset
\Coker(IM\stackrel{\cdot a}{\to}IM)$. However, it is clear that
$I=aI$, so $a$ is regular on $M/IM$, whence
the latter is a torsion-free $\kappa(x)^+$-module and
the claim follows.
\end{pfclaim}

Let $\cU$ be the cofiltered system of all affinoid open
neighborhoods of $x$ in $X$; for $U\in\cU$ let $F_U$ be 
the filtered system of all $K^+$-subalgebras of $\cO^+_X(U)$
of topologically finite presentation. We derive 
$$
\begin{array}{r@{\:\simeq\:}l}
(\L_{X/K}^+)_x\derotimes_{K^+}K^+/aK^+ &
\colim{U\in\cU}\colim{A\in F_U}\L_{A/K^+}^\an
\derotimes_{K^+}K^+/aK^+ \\
& \colim{U\in\cU}\colim{A\in F_U}\L_{A/K^+}
\derotimes_{K^+}K^+/aK^+ \\
& \L_{\cO^+_{X,x}/K^+}\derotimes_{K^+}K^+/aK^+ \\
& \L_{\kappa(x)^+/K^+}\derotimes_{K^+}K^+/aK^+.
\end{array}
$$
Together with theorem \ref{th_deep.ramificat}, this implies 
already that scalar multiplication by $a$ is an automorphism 
of $H_i(\L^+_{X/K})$, for every $i>0$. However, according to
theorem \ref{th_cot.smooth.vanish}(ii), $H_i(\L^+_{X/K})$
is a $K^+$-torsion sheaf of $\cO_X^+$-modules, for $i>0$, whence
the first assertion. It also follows that $(\Omega_{X/K}^+)_x$
is a torsion-free, hence flat, $K^+$-module. To prove that
$(\Omega_{X/K}^+)_x$ is a flat $\cO^+_{X,x}$-module, we 
remark first that 
$(\Omega_{X/K}^\an)_x\simeq(\Omega_{X/K}^+)_x\otimes_{K^+}K$, and
the latter is a flat $\cO_{X,x}$-module, since $X$ is smooth
over $\Spa(K,K^+)$. By lemma \ref{lem_Tor.orgy} it suffices 
therefore to show 
\begin{claim}
$(\Omega_{X/K}^+)_x\otimes_{K^+}K^+/aK^+$ is a flat 
$\cO^+_{X,x}\otimes_{K^+}K^+/aK^+$-module. 
\end{claim}
\begin{pfclaim} By claim \ref{cl_when.a.regul} we know
that $(\Omega^+_{X/K})_x\otimes_{\cO^+_{X,x}}\kappa(x)^+$
is a flat $\kappa(x)^+$-module.
In view of \eqref{eq_separate-stalk}, the claim follows after 
base change to $\kappa(x)^+/a\cdot\kappa(x)^+$.
\end{pfclaim}
\end{proof}

\begin{definition}\label{def_diff.divisible}
\index{Formal scheme(s)!Deeply ramified system of|indref{def_diff.divisible}}
Let $(\fX_\alpha~|~\alpha\in I)$ be a system of formal schemes 
of finite presentation over $\Spf\,K^+$, indexed by a small 
cofiltered category $I$.
\begin{enumerate}
\item
Let $\fX_\infty:=\liminv{\alpha\in I}~\fX_\alpha$, where
the limit is taken in the category of locally ringed spaces.
For every $\alpha\in I$, let $\pi_\alpha:\fX_\infty\to\fX_\alpha$
be the natural morphism of locally ringed spaces. We define 
$\Omega_{\fX_\infty/K^+}^\an:=
\colim{\alpha\in I^o}\pi_\alpha^*(\Omega_{\fX_\alpha/K^+}^\an)$, 
which is a sheaf of $\cO_{\fX_\infty}$-modules.

More generally, we let $\L^\an_{\fX_\infty/K^+}:=
\colim{\alpha\in I^o}\pi_\alpha^*(\L^\an_{\fX_\alpha/K^+})$.
\item
We say that the cofiltered system $(\fX_\alpha~|~\alpha\in I)$
is {\em deeply ramified\/} if the natural morphism
$\Omega_{\fX_\infty/K^+}^\an\to\Omega_{\fX_\infty/K^+}^\an\otimes_{K^+}K$
is an epimorphism.
\end{enumerate}
\end{definition}

\begin{lemma}\label{lem_equiv.on.tower} 
Let $(\fX_\alpha~|~\alpha\in I)$ be a cofiltered system as 
in definition {\em\ref{def_diff.divisible}}. For any morphism 
$\beta\to\alpha$ of $I$, let 
$\ff_{\alpha\beta}:\fX_\beta\to\fX_\alpha$ be the corresponding
morphism of formal $\Spf\,K^+$-schemes. Moreover, for every 
$\alpha\in I$, let $\Omega_{\fX_\alpha/K^+}^\mathrm{tf}$ be 
the image of the morphism 
$\Omega_{\fX_\alpha/K^+}^\an\to
\Omega_{\fX_\alpha/K^+}^\an\otimes_{K^+}K$
(``$\mathrm{tf}$" stays for torsion-free). The following two 
conditions are equivalent :
\begin{enumerate}
\item
The system $(\fX_\alpha~|~\alpha\in I)$ is deeply ramified.
\item
For every $\alpha\in I$ there is a morphism $\beta\to\alpha$
of $I$, such that the image of the natural morphism 
$\ff_{\alpha\beta}^*(\Omega_{\fX_\alpha/K^+}^\mathrm{tf})\to
\Omega_{\fX_\beta/K^+}^\mathrm{tf}$ is contained 
in the subsheaf $a\cdot\Omega_{\fX_\beta/K^+}^\mathrm{tf}$.
\end{enumerate}
\end{lemma}
\begin{proof} It is clear that (ii)$\Rightarrow$(i). We show 
that (i)$\Rightarrow$(ii). Under the above assumptions,
every $\fX_\alpha$ is quasi-compact, hence we can cover it
by finitely many affine formal schemes $\fU_i:=\Spf\,A_i$ 
($i=1,...,n$) of finite type over $\Spf\,K^+$. Then, for every 
$i=1,...,n$, the restriction of 
$\Omega_{\fX_\alpha/K^+}^\an$ to $\fU_i$ is the coherent
sheaf $(\Omega_{A_i/K^+}^\an)^\triangle$ (notation
of \cite[Ch.I, \S 10.10.1]{EGAI}). Hence 
$\Omega^\mathrm{tf}_{\fX_\alpha/K^+}$ is a coherent sheaf of 
$\cO_{\fX_\alpha}$-modules. For every morphism $\beta\to\alpha$, 
let 
$$
U_{\alpha\beta}:=\{x\in\fX_\beta~|~
\Img(\ff^*_{\alpha\beta}(\Omega^\mathrm{tf}_{\fX_\alpha/K^+})_x\to
(\Omega^\mathrm{tf}_{\fX_\beta/K^+})_x)\subset 
a\cdot(\Omega^\mathrm{tf}_{\fX_\beta/K^+})_x\}.
$$
$U_{\alpha\beta}$ is therefore a constructible open subset 
of $\fX_\beta$, and we denote its complement by $Z_{\alpha\beta}$. 
By assumption (i) we know that 
$$
\lim_{\beta\to\alpha}Z_{\alpha\beta}=
\bigcap_{\beta\to\alpha}\pi_\alpha^{-1}(Z_{\alpha\beta})=
\emptyset.
$$
If we retopologize the reduced schemes $Z_{\alpha\beta}$ by
their constructible topologies, we get an inverse system of 
compact spaces, and deduce that some $Z_{\alpha\beta}$ is empty
by \cite[Ch.I, \S 9, n.6, Prop.8(b)]{Bou}.
\end{proof}

\begin{example}\label{ex_prototype} 
The prototype of deeply ramified systems is  given by the 
tower of morphisms 
\set\begin{equation}\label{eq_tower-diff.div}
...\to\B^d_{K^+}(0,\rho^{1/p^n})\stackrel{\phi_n}{\to}
\B^d_{K^+}(0,\rho^{1/p^{n-1}})\stackrel{\phi_{n-1}}{\to}...
\stackrel{\phi_1}{\to}\B^d_{K^+}(0,\rho)
\end{equation}
where, for any $r=(r_1,...,r_d)\in(K^\times)^d$, we have denoted 
$$\B^d_{K^+}(0,|r|):=
\Spf\,K^+\langle r_1^{-1}T_1,...,r_d^{-1}T_d\rangle
$$ 
({\em i.e.}, the formal $d$-dimensional polydisc defined by the 
equations $|T_i|\leq|r_i|$, $i=1,...,d$). The morphisms $\phi_n$ are 
induced by the ring homomorphisms $T_i\mapsto T_i^p$ ($i=1,...,d$).
Notice that the tower \eqref{eq_tower-diff.div} is defined 
whenever 
$\rho_i\in\Gamma_K^{p^\infty}:=\bigcap_{n\in\N}\Gamma_K^{p^n}$
for every $i=1,...,d$. 
We leave to the reader the verification that condition (ii) of lemma 
\ref{lem_equiv.on.tower} is indeed satisfied.
\end{example}

\begin{lemma}\label{lem_cofilt.syst} 
Let $\fX:=(\fX_\alpha~|~\alpha\in I)$ be a cofiltered 
system as in definition {\em\ref{def_diff.divisible}}.
\begin{enumerate}
\item
If $\fY:=(\fY_\alpha~|~\alpha\in I)\to(\fX_\alpha~|~\alpha\in I)$
is a morphism of cofiltered systems such that the induced
morphisms of adic spaces $d(\fY_\alpha)\to d(\fX_\alpha)$
are unramified for every $\alpha\in I$ (cp. \eqref{subsec_functor.d}), 
then $\fY$ is deeply ramified if $\fX$ is.
\item
Let $\fZ:=(\fZ_\beta~|~\beta\in J)$ be another such cofiltered
system, and suppose that $\fX$ and $\fZ$ are isomorphic as
pro-objects of the category of formal schemes. Then $\fX$
is deeply ramified if and only if $\fZ$ is.
\item
If $\fX$ and $\fZ:=(\fZ_\beta~|~\beta\in J)$ are two 
deeply ramified cofiltered systems, then the
fibred product 
$\fX\times\fZ:=(\fX_\alpha\times_{\Spf(K^+)}\fZ_\beta~|~
(\alpha,\beta)\in I\times J)$ is deeply ramified.
\end{enumerate}
\end{lemma}
\begin{proof} (i): by \cite[Prop.2.2]{BoschIII} 
the natural morphism $\Omega_{\fX_\infty/K^+}^\an\otimes_{K^+}K\to
\Omega_{\fY_\infty/K^+}^\an\otimes_{K^+}K$ is an epimorphism;
the claim follows easily.
(ii) and (iii) are easy and shall be left to the reader.
\end{proof}

The counterpart of the above definitions for adic schemes
is given in the following:

\begin{definition}\label{def_diff.div.adic}
\index{Adic space(s)!deeply ramified system of|indref{def_diff.div.adic}}
Let $(X_\alpha~|~\alpha\in I)$ be a system of adic spaces 
of finite type over $\Spa(K,K^+)$, indexed by a small 
cofiltered category $I$.
\begin{enumerate}
\item
Let $(X_\infty,\cO_{X_\infty},\cO_{X_\infty}^+):=
\liminv{\alpha\in I}~(X_\alpha,\cO_{X_\alpha},\cO_{X_\alpha}^+)$, 
where the limit is taken in the category of locally
ringed spaces. For every $\alpha\in I$, let 
$\pi_\alpha:X_\infty\to X_\alpha$ be the natural 
morphism of locally ringed spaces. We define 
$\Omega_{X_\infty/K}^+:=
\colim{\alpha\in I^o}\pi_\alpha^*(\Omega_{X_\alpha/K}^+)$, 
which is a sheaf of $\cO_{X_\infty}^+$-modules.
More generally, we let $\L^+_{X_\infty/K}:=
\colim{\alpha\in I^o}\pi_\alpha^*(\L^+_{X_\alpha/K})$.
\item
We say that the cofiltered system $(X_\alpha~|~\alpha\in I)$
is {\em deeply ramified\/} if the natural morphism
$\Omega_{X_\infty/K}^+\to\Omega_{X_\infty/K}^+\otimes_{K^+}K$
is an epimorphism.
\end{enumerate}
\end{definition}

\sset\subsubsection{}
Let $(X_\alpha~|~\alpha\in I)$ be as in definition
\ref{def_diff.div.adic}. For every $x\in X_\infty$, we let 
$\kappa(x)^+:=\cO^+_{X_\infty,x}/
\bigcap_{n\in\N}a^n\cO^+_{X_\infty,x}$. The ring
$\kappa(x)^+$ is a filtered colimit of valuation rings, hence
it is a valuation ring. Moreover, the image of $a$ in
$\kappa(x)^+$ is topologically nilpotent for the valuation
topology of $\kappa(x)^+$.

\sset\subsubsection{}\label{subsec_go.to.extreme}
Given a cofiltered system $\underline\fX:=(\fX_\alpha~|~\alpha\in I)$
of formal schemes, one obtains a cofiltered system of
adic spaces $\underline X:=(X_\alpha:=d(\fX_\alpha)~|~\alpha\in I)$,
and using \eqref{subsec_projlim.is.adic} and lemma
\ref{lem_cofilt.syst}(i) one sees easily that $\underline X$
is deeply ramified whenever $\underline\fX$ is.
Together with example \ref{ex_prototype}, this yields
plenty of examples of deeply ramified systems of adic spaces.

\begin{proposition}\label{prop_to.extremes} 
Let $\underline X:=(X_\alpha~|~\alpha\in I)$ be a deeply ramified 
cofiltered system of adic spaces. Then, for every point 
$x\in X_\infty$, the valuation ring $\kappa(x)^+$ is deeply ramified.
\end{proposition}
\begin{proof}
For every $K^+$-module $M$, let us denote by $T_n(M)$ the
submodule of $M$ annihilated by $a^n$. Furthermore, let
$T(M):=\bigcup_{n\in\N}T_n(M)$. Since the cofiltered
system $\underline X$ is deeply ramified, we have: 
$(\Omega_{X_\infty/K}^+)_x=
T(\Omega_{X_\infty/K}^+)_x+a\cdot(\Omega_{X_\infty/K}^+)_x$.
To lighten notation, let $\cO^+_x:=\cO^+_{X_\infty,x}$.
From lemma \ref{lem_compare-oms} one deduces easily that 
the natural map 
$\Omega_{\cO^+_x/K^+}\to(\Omega_{X_\infty/K}^+)_x$
is surjective with $a$-divisible kernel.
Hence, by snake lemma, the induced map
$T_n(\Omega_{\cO^+_x/K^+})\to T_n(\Omega_{X_\infty/K}^+)_x$
is surjective for every $n$, and {\em a fortiori\/} the
map $T(\Omega_{\cO^+_x/K^+})\to T(\Omega_{X_\infty/K}^+)_x$
is onto. It follows easily that 
$\Omega_{\cO^+_x/K^+}=T(\Omega_{\cO^+_x/K^+})+
a\cdot\Omega_{\cO^+_x/K^+}$, and consequently:
$\Omega_{\kappa(x)^+/K^+}=T(\Omega_{\kappa(x)^+/K^+})+
a\cdot\Omega_{\kappa(x)^+/K^+}$. However, it follows easily
from theorem \ref{th_deep.ramificat} that 
$T(\Omega_{\kappa(x)^+/K^+})=0$, so finally 
$\Omega_{\kappa(x)^+/K^+}=a\cdot\Omega_{\kappa(x)^+/K^+}$ and 
\set\begin{equation}\label{eq_final.countdown}
\L_{\kappa(x)^+/K^+}\derotimes_{K^+}K^+/aK^+\simeq 0.
\end{equation}
On the other hand, if $(E,|\cdot|_E)$ is any valued field 
extension of $\kappa(x)^+$, we have 
\set\begin{equation}\label{eq_final.countdown.go}
H_i(\L_{E^+/K^+}\otimes_{K^+}K^+/aK^+)=0\quad
\text{for all $i>0$}
\end{equation}
by theorem \ref{th_deep.ramificat}. From \eqref{eq_final.countdown}, 
\eqref{eq_final.countdown.go} and transitivity for the tower 
$K^+\subset\kappa(x)^+\subset E^+$, we derive 
$H_i(\L_{E^+/\kappa(x)^+}\otimes_{K^+}K^+/aK^+)=0$ for all
$i>0$.
Again by theorem \ref{th_deep.ramificat} and remark 
\ref{rem_deeply.ramificat}
we conclude.
\end{proof}

\sset\subsubsection{}
Let $\underline X$ be a cofiltered system as in definition 
\ref{def_diff.div.adic}. Let $\cA$ be a sheaf of 
$\cO^+_{X_\infty}$-algebras. We say that $\cA$ is a 
{\em weakly {\'e}tale\/ $\cO^+_{X_\infty}$-algebra\/} if, for every
$x\in X_\infty$, the stalk $\cA_x^a$ is a weakly
{\'e}tale $\cO^{+a}_{X_\infty,x}$-algebra.

\begin{theorem}\label{th_weak.purity} Suppose that $K$ is 
deeply ramified, and let $\underline X:=(X_\alpha~|~\alpha\in I)$ 
be a deeply ramified cofiltered system. Let also
$\underline f:\underline Y:=(Y_\alpha~|~\alpha\in I)\to\underline X$ 
be a morphism of cofiltered systems, such that the morphisms
$Y_\alpha\to X_\alpha$ are finite {\'e}tale for every 
$\alpha\in I$. Then $f_{\infty*}\cO_{Y_\infty}^+$
is a weakly {\'e}tale $\cO_{X_\infty}^+$ algebra.
\end{theorem}
\begin{proof} To lighten notation, let us write $\cO^+$
(resp. $\cA^+$) instead of $\cO_{X_\infty}^+$
(resp. $f_{\infty*}\cO_{Y_\infty}^+$). For every
$\alpha\in I$ consider the cofiltered system 
$\underline Z(\alpha):=\underline X\times_{X_\alpha}Y_\alpha$ 
indexed by $I/\alpha$ (the category of morphisms $\beta\to\alpha$),
which is defined by setting $\underline Z(\alpha)_{\beta\to\alpha}:=
X_\beta\times_{X_\alpha}Y_\alpha$. We have obvious
morphisms of cofiltered systems 
$\underline f_{/\alpha}:\underline Z(\alpha)\to\underline X$;
the sheaf $\cA^+$ is the colimit of the sheaves 
$f_{/\alpha\infty *}\cO^+_{Z(\alpha)_\infty}$.
Hence it suffices to prove the assertion for the latter
sheaves, and therefore in order to show the theorem we can 
and do assume that there exists $\alpha\in I$ such that, for 
every $\beta\to\alpha$, the induced commutative diagram
$$
\xymatrix{
Y_\beta \ar[r] \ar[d] & Y_\alpha \ar[d] \\
X_\beta \ar[r] & X_\alpha
}
$$
is cartesian. Let $x\in X_\infty$.

\begin{claim}\label{cl_halfway} 
$\cA_x^+$ is a flat $\cO^+_x$-algebra. 
\end{claim}
\begin{pfclaim} 
On one hand, by assumption $\cA^+_x[1/a]$ is a flat 
$\cO^+_x[1/a]$-algebra; on the other hand, 
$\cA_x^+\otimes_{\cO^+_x}\kappa(x)^+$ is a flat 
$\kappa(x)^+$-module by claim \ref{cl_when.a.regul}, so 
that $\cA_x^+/a\cA^+_x$ is a flat 
$\cO^+_x/a\cO^+_x$-module; thus the claim follows
from lemma \ref{lem_Tor.orgy}.
\end{pfclaim}

Let $e\in C:=\cA^+_x\otimes_{\cO^+_x}\cA^+_x[1/a]$ be the
idempotent provided by lemma \ref{prop_idemp}. In view of
claim \ref{cl_halfway}, we only have to show that 
$\eps\cdot e\in C^+:=\cA^+_x\otimes_{\cO^+_x}\cA^+_x$ for
every $\eps\in\fm$. Let $\bar e$ be the image of $e$ in 
$C\otimes_{\cO_x}\kappa(x)$.
Set $I:=\bigcap_{n\geq 0}a^n\cO^+_x$.
\begin{claim}\label{cl_int.closure}
$\cA_x^+/I\cA_x^+$ is the integral closure of of $\kappa(x)^+$
in $\cA_x/I\cA_x$.
\end{claim}
\begin{pfclaim} First of all, since $I=aI$, we have 
$I\cA_x=I\cA_x^+$, and thus the natural map 
$\cA_x^+/I\cA_x^+\to\cA_x/I\cA_x$
is injective. Next, we remark that $\cA_x^+$
is the integral closure of $\cO_x^+$ in $\cA_x$; indeed,
this follows from \cite[\S 6.2.2, Lemma 3, Prop.2]{Bo-Gun},
after taking colimits. This already shows that $\cA_x^+/I\cA_x^+$
is a subalgebra of $\cA_x/I\cA_x$ integral over $\kappa(x)^+$.
To conclude, suppose that $\bar f\in\cA_x/I\cA_x$ satisfies
an integral equation: 
$\bar f{}^n+\bar b_1\cdot\bar f{}^{n-1}+...+\bar b_n=0$, for
certain $\bar b_1,...,\bar b_n\in\cA_x^+/I\cA_x^+$; pick arbitrary
representatives $f\in\cA_x$, $b_i\in\cA_x^+$ of these elements.
Then $f^n+b_1\cdot f^{n-1}+...+b_n\in I\cA_x$.
Since $I\cA_x=I\cA_x^+$, we deduce that $f$ is integral
over $\cA_x^+$, so $f\in\cA_x^+$ and the claim follows.
\end{pfclaim}

\begin{claim}\label{cl_etale.closure}
$(\cA_x^+/I\cA_x^+)^a$ is an {\'e}tale $\kappa(x)^{+a}$-algebra.
\end{claim}
\begin{pfclaim} In view of propositions \ref{prop_converse}(ii) 
and \ref{prop_crit.projectivity}, it suffices to show that
$(\cA_x^+/I\cA_x^+)^a$ is weakly {\'e}tale over $\kappa(x)^{+a}$.
Let $\fp$ be the height one prime ideal of $\kappa(x)^+$; then 
$\kappa(x)^+_\fp$ is a rank one valuation ring and the localization 
map induces isomorphisms 
$\kappa(x)^{+a}\stackrel{\sim}{\to}\kappa(x)^{+a}_\fp$,
$(\cA_x^+/I\cA_x^+)^a\stackrel{\sim}{\to}(\cA_x^+/I\cA_x^+)_\fp^a$
(recall that the standing basic setup is the standard setup
of $K^+$). Taking claim \ref{cl_int.closure} into account, it 
suffices then to show that $(\cA_x^+/I\cA_x^+)_\fq^a$ is weakly 
{\'e}tale over $\kappa(x)^{+a}_\fp$, for every prime ideal 
$\fq\subset\cA_x^+/I\cA_x^+$ of height one.
By proposition \ref{prop_to.extremes}, $\kappa(x)^+$
is deeply ramified, hence the same holds for $\kappa(x)^+_\fp$.
Let $\kappa(x)^\sep$ be a separable closure of $\kappa(x)$,
and $\kappa(x)^{\sep +}$ a valuation ring of $\kappa(x)^\sep$
dominating $\kappa(x)^+_\fp$; we can assume that 
$(\cA_x^+/I\cA_x^+)_\fq\subset\kappa(x)^{\sep +}$. 
From claim \ref{cl_int.closure} we deduce that
$(\cA_x^+/I\cA_x^+)_\fq$ is a valuation ring, hence 
$\kappa(x)^{\sep +}$ is a faithfully flat
$(\cA_x^+/I\cA_x^+)_\fq$-algebra;
then the claim follows in view of lemma \ref{lem_itoiv}(viii) 
and proposition \ref{prop_deep.ramif}.
\end{pfclaim}

From claim \ref{cl_etale.closure} we deduce that 
$\eps\cdot\bar e\in C^+\otimes_{\cO^+_x}\kappa(x)^+$ 
for every $\eps\in\fm$. To conclude, it suffices therefore
to remark:
\begin{claim}\label{cl_liftup} For every $\eps\in\fm$, we have
$\eps\cdot e\in C^+$ if and only if
$\eps\cdot\bar e\in C^+\otimes_{\cO^+_x}\kappa(x)^+$.
\end{claim}
\begin{pfclaim}
Let $e_\eps\in C^+$ be any lifting
of $\eps\cdot\bar e$; then $\eps\cdot e-e_\eps$ is in the
kernel of the projection $C\to C\otimes_{\cO^+_x}\kappa(x)$.
Let $n\in\N$ be a large enough integer, so that 
$a^n\cdot e\in C^+$; it follows that 
$a^n\cdot(\eps\cdot e-e_\eps)$ is in the kernel of the
projection $C^+\to C^+\otimes_{\cO^+_x}\kappa(x)^+$, consequently
$a^n\cdot(\eps\cdot e-e_\eps)=a^n\cdot c$ for some $c\in C^+$.
Since $a$ is a non-zero-divisor in $C^+$, it follows that
$\eps\cdot e=c+e_\eps\in C^+$, as required.
\end{pfclaim}
\end{proof}

\subsection{Semicontinuity of the discriminant}\label{sec_semicont}
\begin{definition}\label{def_discriminant}
\index{Almost module(s)!almost projective!$\fd_A(P,b)$ : discriminant of a bilinear form on a|indref{def_discriminant}}
Let $(V,\fm)$ be a basic setup, $A$ a $V^a$-algebra and 
$P$ an almost projective $A$-module of constant rank $r\in\N$. 
Suppose moreover that $P$ is endowed with a bilinear form 
$b:P\otimes_AP\to A$. 
We let $\beta:P\to P^*$ be the $A$-linear morphism defined by 
the rule: $\beta(x)(y):=b(x\otimes y)$ for every $x,y\in P_*$.
The {\em discriminant\/} of the pair $(P,b)$ is the ideal
$$
\fd_A(P,b):=
\Ann_A\,\Coker(\Lambda_A^r\beta:\Lambda^r_AP\to\Lambda^r_AP^*).
$$
\end{definition}

\sset\subsubsection{}\label{subsec_disc-of-algebra}
\index{Almost algebra(s)!almost projective!$\fd_{B/A}$ : discriminant 
of an|indref{subsec_disc-of-algebra}}
As a special case, we can consider the pair $(B,t_{B/A})$
consisting of an $A$-algebra $B$ which is almost projective 
of constant rank $r$ over $A$, and its trace form $t_{B/A}$. 
In this situation, we let $\fd_{B/A}:=\fd_A(B,t_{B/A})$,
and we call this ideal the {\em discriminant of the $A$-algebra} 
$B$.

\begin{lemma}\label{lem_discr.criterio} 
Let $B$ be an almost projective $A$-algebra of constant rank 
$r$ as an $A$-module. Then $B$ is {\'e}tale over $A$ if and only 
if\/ $\fd_{B/A}=A$.
\end{lemma}
\begin{proof} By theorem \ref{th_proj.etale}, it is clear that 
$\fd_{B/A}=A$ when $B$ is {\'e}tale over $A$. Suppose therefore 
that $\fd_{B/A}=A$; it follows that $\Lambda^r_A\tau_{B/A}$ is 
an epimorphism.
However, by proposition \ref{prop_decomp.fin.rank}, 
$\Lambda^r_AB$ and $\Lambda^r_AB^*$ are invertible $A$-modules.
It then follows by lemma \ref{lem_explain.invert}(iv) that
$\Lambda^r_A\tau_{B/A}$ is an isomorphism, hence $\tau_{B/A}$ 
is an isomorphism, by virtue of proposition \ref{prop_Cramer}.
One concludes again by theorem \ref{th_proj.etale}.
\end{proof}

\begin{lemma}\label{lem_Fitting.discr} 
Let $(V,\fm)$ be the standard setup associated
to a valued field $(K,|\cdot|)$ (cp. \eqref{subsec_standard.setup},
especially, $V:=K^+$).
Let $P'\subset P$ be two almost projective $V^a$-modules
of constant rank equal to $r$. Let $b:P\otimes_{V^a}P\to V^a$ 
be a bilinear form, such that $b\otimes_{V^a}\one_{K^a}$ is a 
perfect pairing, and denote by $b'$ the restriction of $b$ to 
$P'\otimes_{V^a}P'$. Then we have:
$$
\fd_{V^a}(P',b')=F_0(P/P')^2\cdot\fd_{V^a}(P,b).
$$
\end{lemma}
\begin{proof} Let $j:P'\to P$ be the imbedding, $\beta:P\to P^*$
(resp. $\beta':P'\to P^{\prime *}$) the $V^a$-linear morphism 
associated to $b$ (resp. to $b'$).
The assumptions implies that $\Lambda^r_{V^a}j$, $\Lambda^r_{V^a}j^*$,
$\Lambda^r_{V^a}\beta$ and $\Lambda^r_{V^a}\beta'$ are all
injective, and clearly we have
$$
\Lambda^r_{V^a}j^*\circ\Lambda^r_{V^a}\beta\circ\Lambda^r_{V^a}j=
\Lambda^r_{V^a}\beta'.
$$
There follow short exact sequences:
$$
\begin{array}{c}
0\to\Coker\,\Lambda^r_{V^a}\beta\to\Coker\,\Lambda^r_{V^a}(j^*\circ\beta)\to
\Coker\,\Lambda^r_{V^a}j^*\to 0 \\
0\to\Coker\,\Lambda^r_{V^a}j\to\Coker\,\Lambda^r_{V^a}\beta'\to
\Coker\,\Lambda^r_{V^a}(j^*\circ\beta)\to 0.
\end{array}
$$
Using lemma \ref{lem_Fitting.vals} and remark \ref{rem_verbatim}(ii),
we deduce: 
$$
F_0(\Coker\,\Lambda^r_{V^a}\beta')=
F_0(\Coker\,\Lambda^r_{V^a}j)\cdot
F_0(\Coker\,\Lambda^r_{V^a}\beta)\cdot
F_0(\Coker\,\Lambda^r_{V^a}j^*).
$$
\begin{claim}\label{cl_calculate-Fit}
Let $Q$ be an almost finitely generated projective
$A$-module of constant rank $r$, which is also uniformly almost
finitely generated. Let $Q'\subset Q$ any submodule. Then:
$$
F_0(Q/Q')=\Ann_A(\Coker(\Lambda^r_AQ'\to\Lambda^r_AQ)).
$$
\end{claim}
\begin{pfclaim} By lemma \ref{lem_annihilate.and.forget}
and proposition \ref{prop_Fitting.short}(ii) we can check
the identity after a faithfully flat base change, hence
we can suppose that $Q\simeq A^r$, in which case the identity
holds by definition of the Fitting ideal $F_0$.
\end{pfclaim}

Notice that $\Lambda^r_{V^a}P$ and $\Lambda^r_{V^a}P'$ are
invertible $V^a$-modules by virtue of proposition 
\ref{prop_decomp.fin.rank}, consequently claim \ref{cl_calculate-Fit}
applies and yields:
$F_0(\Coker\,\Lambda^r_{V^a}\beta')=\fd_{V^a}(P',b')$,
$F_0(\Coker\,\Lambda^r_{V^a}\beta)=\fd_{V^a}(P,b)$ and
$F_0(\Coker\,\Lambda^r_{V^a}j)=F_0(\Coker\,\Lambda^r_{V^a}j^*)=
F_0(P/P')$.
\end{proof}

\sset\subsubsection{}\label{subsec_discrim.sep.ext}
\index{Valued field(s)!extension of a!$\fd^+_{L/K}$ : discriminant of an 
algebraic|indref{subsec_discrim.sep.ext}}
After these generalities, we return to the standard
setup $(K^+,\fm)$ of this chapter, associated to a valued
field $(K,|\cdot|)$ of rank one (cp. \eqref{subsec_standard.setup}).
Consider an {\'e}tale $K$-algebra $L$; we denote by $W_L$ the
integral closure of $K^+$ in $L$. $L$ is the product of
finitely many separable field extensions of $K$, therefore
$W_L^a$ is an almost projective $K^{+a}$-module of constant
rank $n$, by proposition \ref{prop_fin.sep.uniform}.
Hence, the discriminant of $W_L^a$ over $K^{+a}$ is defined,
and to lighten notation, we will denote it by $\fd^+_{L/K}$. 
Furthermore, since $L$ is {\'e}tale over $K$, it is clear that 
$\fd^+_{L/K}$ is a fractional ideal of $K^{+a}$ (cp. 
\eqref{subsec_fract.ideals}). Let 
$|\cdot|:\mathrm{Div}(K^{+a})\to\Gamma_K^\wedge$ be the isomorphism
provided by lemma \ref{lem_Picpus}. We obtain an element
$|\fd^+_{L/K}|\in\Gamma_K^\wedge$; after choosing (cp. example 
\eqref{ex_tens.with.Q}(vi)) an order preserving isomorphism
\set\begin{equation}\label{eq_choice.of.iso}
((\Gamma_K\otimes_\Z\Q)^\wedge,\leq)\stackrel{\sim}{\to}
(\R_{>0},\leq)
\end{equation}
on the multiplicative group of positive real
numbers, we can then view $|\fd^+_{L/K}|\in(0,1]$.

\begin{lemma}\label{lem.Huber.again} 
Let $K$, $L$ be as in \eqref{subsec_discrim.sep.ext} and denote 
by $K^\wedge$ the completion of $K$ for the valuation topology. 
Set $L^\wedge:=K^\wedge\otimes_KL$. Then 
$K^{\wedge+}\otimes_{K^+}\fd^+_{L/K}=\fd^+_{L^\wedge/K^\wedge}$.
\end{lemma}
\begin{proof} Since the base change $K\to K^\wedge$ is
faithfully flat, everything is clear from the definitions, once
we have established that 
$W_{L^\wedge}\simeq K^{\wedge+}\otimes_{K^+}W_L$. However, both
rings can be identified with the $a$-adic completion $(W_L)^\wedge$
of $W_L$, so the assertion follows.
\end{proof}

\sset\subsubsection{}\label{subsec_approach.semicont}
Let $X$ be an adic space locally of finite type over 
$\Spa(K,K^+)$. $X$ is a locally spectral space, and every 
point $x\in X$ admits a unique maximal generisation
$r(x)\in X$. The valuation ring $\kappa(r(x))^+$ has
rank one, and admits a natural imbedding 
$K^+\subset\kappa(r(x))^+$, continuous for the valuation 
topologies; especially, the image of the topologically 
nilpotent element $a$ is topologically nilpotent in 
$\kappa(r(x))^+$. This imbedding induces therefore a
natural isomorphism of completed value groups
$$
(\Gamma_K\otimes_\Z\Q)^\wedge\stackrel{\sim}{\to}
(\Gamma_{\kappa(r(x))}\otimes_\Z\Q)^\wedge.
$$
In particular, our original choice of isomorphism 
\eqref{eq_choice.of.iso} fixes univocally a similar isomorphism
for every point $r(x)$.
We denote by ${\mathcal M}(X)$ the set $r(X)$ endowed with the quotient 
topology induced by the mapping $X\to r(X)$ : $x\mapsto r(x)$.
This topology is coarser than the subspace topology induced by 
the imbedding into $X$. The mapping $x\mapsto r(x)$ is a
retraction of $X$ onto the subset ${\mathcal M}(X)$ of its maximal 
points. If $X$ is a quasi-separated quasi-compact adic space, 
${\mathcal M}(X)$ is a compact Hausdorff topological space
(\cite[8.1.8]{Hu2}).

\sset\subsubsection{}\label{subsec_let.it.be}
\index{Adic space(s)!{\'e}tale morphism of!$\fd^+_{Y/X}$ : discriminant 
of an|indref{subsec_let.it.be}}
Let $X$ be as in \eqref{subsec_approach.semicont}, and let
$f:Y\to X$ be a finite {\'e}tale morphism of adic spaces. For
every point $x\in X$, the fibre 
$E(x):=(f_*\cO_Y)_x\otimes_{\cO_{X,x}}\kappa(x)$ is a finite
{\'e}tale $\kappa(x)$-algebra. If now $x\in{\mathcal M}(X)$, we
can consider the discriminant $\fd^+_{E(x)/\kappa(x)}$ defined
as in \eqref{subsec_discrim.sep.ext} (warning: notice that the 
definition makes sense when we choose the standard setup associated 
to the valuation ring $\kappa(x)^+$; since it may happen that the
valuation of $K$ is discrete and that of $\kappa(x)$ is not
discrete, the setups relative to $K$ and to $\kappa(x)$
may not agree in general). Upon passing to absolute 
values, we finally obtain a real valued function:
$$
\fd^+_{Y/X}:{\mathcal M}(X)\to(0,1] \qquad 
x\mapsto|\fd^+_{E(x)/\kappa(x)}|.
$$
The study of the function $\fd^+_{Y/X}$ is reduced easily
to the case where $X$ (hence $Y$) are affinoid. In such
case, one can state the main result in a more general form,
as follows.

\begin{definition}\label{def_spectrum} Let $A$ be any
(commutative unitary) ring.
\index{${\mathcal N}(A)$ : Spectrum of seminorms of a 
ring|indref{def_spectrum}}
\begin{enumerate}
\item
We denote by ${\mathcal N}(A)$ the set consisting of all 
multiplicative ultrametric seminorms $|\cdot|:A\to\R$. For every
$x\in{\mathcal N}(A)$ and $f\in A$ we write usually $|f(x)|$
in place of $x(f)$. ${\mathcal N}(A)$ is endowed with the 
coarsest topology such that, for every $f\in A$,
the real-valued map $|f|:{\mathcal N}(A)\to\R$ given by the
rule: $x\mapsto|f(x)|$, is continuous. 
\item
For every $x\in{\mathcal N}(A)$, we let 
$\mathrm{Supp}(x):=\{f\in A~|~|f(x)|=0\}$. Then 
$\mathrm{Supp}(x)$ is a prime ideal and we set
$\kappa(x):=\mathrm{Frac}(A/\mathrm{Supp}(x))$. The seminorm
$x$ induces a valuation on the residue field $\kappa(x)$,
and as usual we denote by $\kappa(x)^+$ its valuation ring.
\item
Let $A\to B$ be a finite {\'e}tale morphism.
For every $x\in{\mathcal N}(A)$, we let $E(x):=B\otimes_A\kappa(x)$.
Then $E(x)$ is an {\'e}tale $\kappa(x)$-algebra, so we can 
define $\fd_{B/A}^+(x):=\fd^+_{E(x)/\kappa(x)}$ (cp. the warning
in \eqref{subsec_let.it.be}). By setting $x\mapsto|\fd_{B/A}^+(x)|$ 
we obtain a well-defined function
$$
|\fd_{B/A}^+|:{\mathcal N}(A)\to(0,1].
$$
\end{enumerate}
\end{definition}

\sset\subsubsection{}
If $X=\Spa(A,A^+)$, with $A$ a complete $K$-algebra of
topologically finite type, then ${\mathcal M}(X)$ is naturally
homeomorphic to the subspace ${\mathcal M}(A)$ of ${\mathcal N}(A)$
consisting of the continuous seminorms that extend the absolute value
of $K$ given by \eqref{eq_choice.of.iso}. 
It is shown in \cite[\S 1.2]{Berk} that ${\mathcal M}(A)$
is a compact Hausdorff space, for every Banach $K$-algebra $A$.

\begin{proposition}\label{prop_ring.version} 
Let $A$ be a ring, $B$ a finite {\'e}tale $A$-algebra. 
Then the function $|\fd^+_{B/A}|$ is lower semi-continuous 
({\em i.e.} it is continuous for the topology of $(0,1]$ 
whose open subsets are of the form $(c,1]$, $c\in[0,1]$).
\end{proposition}
\begin{proof} Let $f\in A$ be any element; notice that
${\mathcal N}(A[1/f])$ is naturally homeomorphic to the open
subset $U(f):=\{x\in{\mathcal N}(A)~|~|f(x)|\neq 0\}$.
Hence, after replacing $A$ by some localization, we can assume
that $B$ is a free $A$-module, say of rank $n$. For every $b\in B$,
let $\chi(b,T):=T^n+s_1(b)\cdot T^{n-1}+...+s_n(b)$ be the 
characteristic polynomial of the $A$-linear endomorphism $B\to B$ 
given by the rule $b'\mapsto b'\cdot b$.
\begin{claim}\label{cl_charact.integrality} 
For every point $x\in{\mathcal N}(A)$ and every $b\in B$, the 
following are equivalent:
\begin{enumerate}
\item
$b\otimes 1\in W_{E(x)}$.
\item
$|s_i(b)(x)|\leq 1$ for $i=1,...,n$.
\end{enumerate}
\end{claim}
\begin{pfclaim} Indeed, if (ii) holds, then the image of $\chi(b,T)$ 
in $\kappa(x)[T]$ is a monic polynomial with coefficients in
$\kappa(x)^+$ and $b\otimes 1$ is one of its roots (Cayley-Hamilton), 
hence $b\otimes 1$ is integral over $\kappa(x)^+$, which is (i).
Conversely, if (i) holds, let $E(x)=\prod_{j=1}^kE_j$ be the
decomposition of $E(x)$ as product of finite separable extensions
of $\kappa(x)$, and let $b_j\in E_j$ be the image of $b$. It
follows that $b_j\in W_{E_j}$ for every $j=1,...,k$, and moreover
the image $\bar\chi(b,T)$ of $\chi(b,T)$ in $\kappa(x)[T]$ decomposes 
as a product $\prod_{j=1}^k\chi(b_j,T)$. It suffices therefore to 
show that the coefficients $s_{ij}(b_j)$ satisfy (ii) for every 
$i\leq n$ and $j\leq k$, so we can assume that $E(x)$ is a field.
Let $m_b(T)\in\kappa(x)[T]$ be the minimal polynomial of $b\otimes 1$; 
it is well known that $\bar\chi(b,T)$ divides $m_b(T)^n$, hence 
the roots of $\bar\chi(b,T)$ are conjugates of $b$ under the action 
of $G:=\Gal(\kappa(x)^\mathrm{a}/\kappa(x))$. Let $C$ be the integral
closure of $\kappa(x)^+$ in a finite Galois extension of $\kappa(x)$
containing $E(x)$; $C$ is an integral $\kappa(x)^+$-algebra and the 
Galois conjugates of $b\otimes 1$ are all contained in $C$. Since the
latter are the roots of $\bar\chi(b,T)$, the elements
$s_i(b)\otimes 1$ are symmetric polynomials of the elements 
$\sigma(b\otimes 1)$ ($\sigma\in G$), so $s_i(b)\otimes 1\in C\cap
\kappa(x)=\kappa(x)^+$, which is (ii).
\end{pfclaim}

Let $t_{B/A}$ be the trace morphism of the $A$-algebra $B$, and 
let $x\in{\mathcal N}(A)$. Then the trace morphism 
$t_x:=t_{E(x)/\kappa(x)}$ equals $t_{B/A}\otimes_A\one_{\kappa(x)}$.

\begin{claim} 
For every real number $\eps>0$ we can find a free 
$\kappa(x)^+$-submodule $W_\eps$ of $W_{E(x)}$, such that 
$W_\eps\otimes_{\kappa(x)^+}\kappa(x)=E(x)$ and 
\set\begin{equation}\label{eq_firstpunch}
|\fd_{\kappa(x)^+}(W_\eps,t_x)|+\eps>|\fd^+_{B/A}(x)|.
\end{equation}
\end{claim}
\begin{pfclaim} From claim \ref{cl_finish.up}, we derive that, 
for every positive real number $\delta<1$ there exists a free 
finitely generated $\kappa(x)^+$-submodule $W_\delta\subset W_{E(x)}$ 
such that $|F_0(W^a_{E(x)}/W^a_\delta)|>\delta$; in view of lemma 
\ref{lem_Fitting.discr}, the claim follows easily.
\end{pfclaim}

Let $w_1,...,w_n$ be a basis of $W_\eps$; up to replacing $A$
by a localization $A[1/g]$ for some $g\in A$, we
can write $w_i=b_i\otimes 1$, for some $b_i\in B$ ($i=1,...,n$).
Consequently:
\set\begin{equation}\label{eq_secondpunch}
|\fd_{\kappa(x)^+}(W_\eps,t_x)|=|\det(t_{B/A}(b_i\otimes b_j))(x)|.
\end{equation}
By claim \ref{cl_charact.integrality} we have $|s_j(b_i)(x)|\leq 1$
for $i,j=1,...,n$. Let $1>\delta>0$ be a real number; 
for every $i,j\leq n$, we define an open neighborhood $U_{ij}$ 
of $x$ in ${\mathcal N}(A)$ as follows. Suppose first that $|s_j(b_i)(x)|<1$; 
since the real-valued function $y\mapsto|s_j(b_i)(y)|$ is continuous 
on ${\mathcal N}(A)$ (for the standard topology of $\R$), we can 
find $U_{ij}$ such that $|s_j(b_i)(y)|\leq 1$ for all $y\in U_{ij}$.

Suppose next that $|s_j(b_i)(x)|=1$; then, up to replacing
$A$ by $A[1/s_j(b_i)]$, we can assume that $s_j(b_i)$ is
invertible in $A$. We pick $U_{ij}$ such that 
\set\begin{equation}\label{eq_ineq.deltas}
|s_j(b_i)(y)|\leq 1+\delta\quad\text{for every $y\in U_{ij}$}.
\end{equation}
We set $U:=\bigcap_{1\leq i,j\leq n}U_{ij}$.
Next, we define, for every $y\in U$, an element $c_y\in A$,
as follows. Choose $\alpha,\beta\leq n$ such that 
$|s_\beta(b_\alpha)(y)|=\max_{1\leq i,j\leq n}|s_j(b_i)(y)|$. 
If $|s_\beta(b_\alpha)(y)|\leq 1$, then set $c_y:=1$;
if $|s_\beta(b_\alpha)(y)|>1$, set $c_y:=s_\beta(b_\alpha)^{-1}$.
Then $|s_j(c_y\cdot b_i)(y)|\leq 1$ for every $i,j=1,...,n$ 
and every $y\in U$. Let $\bar W_y$ be the $\kappa(y)^+$-submodule 
of $B\otimes_A\kappa(y)$ spanned by the images of
$c_y\cdot b_1,...,c_y\cdot b_n$. It follows that 
$\bar W_y\subset W_{E(y)}$ for every $y\in U$. 
We compute:
$$
\begin{array}{r@{\:=\:}l}
|\fd_{Y/X}^+(y)|\geq|\fd_{\kappa(y)^+}(\bar W_y,t_y)| &
|\det(t_{B/A}(c_y\cdot b_i\otimes c(y)\cdot b_j))(y)| \\
& |c_y(y)|^{2n}\cdot|\det(t_{B/A}(b_i\otimes b_j))(y)|.
\end{array}
$$
However, the real-valued function 
$y\mapsto|\det(t_{B/A}(b_i\otimes b_j))(y)|$ is continuous
on ${\mathcal N}(A)$, therefore, combining \eqref{eq_firstpunch}
and \eqref{eq_secondpunch}, we see that, up to shrinking
further the open neighborhood $U$, we can assume that
$|\det(t_{B/A}(b_i\otimes b_j))(y)|+\eps>|\fd^+_{B/A}(x)|$
for all $y\in U$, so finally:
$$
|\fd^+_{B/A}(y)|\geq(1-\delta)^{2n}\cdot(|\fd_{B/A}^+(x)|-\eps)
\quad\text{for every $y\in U$}
$$
which implies the claim.
\end{proof}

\begin{theorem}\label{th_disc.upper.semi.cont} 
Let $Y\to X$ be as in \eqref{subsec_let.it.be}. Then the map 
$\fd^+_{Y/X}$ is lower semi-continuous.
\end{theorem}
\begin{proof} We can assume that $X=\Spa(A,A^+)$, where $A$
is a complete $K$-algebra of topologically finite type, and 
therefore $Y=\Spa(B,B^+)$, for a finite {\'e}tale $A$-algebra $B$. 
Then there is a natural homeomorphism 
$\omega:{\mathcal M}(X)\stackrel{\sim}{\to}{\mathcal M}(A)$, so 
the theorem follows from proposition \ref{prop_ring.version} and :
\begin{claim} $\fd^+_{Y/X}=\fd^+_{B/A}\circ\omega$.
\end{claim}
\begin{pfclaim} Let $x\in{\mathcal M}(X)$; $x$ corresponds to
a rank one valuation of $A$, whose value group we identify with
(a subgroup of) $\R_{>0}$ according to \eqref{subsec_approach.semicont}.
The resulting multiplicative seminorm is $\omega(x)$.
We derive easily a natural imbedding 
$\iota:\kappa(\omega(x))\subset\kappa(x)$, compatible with the
identifications of value groups. One knows moreover that $\iota$
induces an isomorphism on completions 
$\iota^\wedge:\kappa(\omega(x))^\wedge
\stackrel{\sim}{\to}\kappa(x)^\wedge$, so the claim follows from
lemma \ref{lem.Huber.again}.
\end{pfclaim}
\end{proof}

\begin{lemma}\label{lem_down.to.finite} 
Let $(K_\alpha,|\cdot|_\alpha~|~\alpha\in I)$ be a system of 
valued field extensions of $(K,|\cdot|)$, indexed by a filtered 
small category $I$ and such that $K^+_\alpha$ is a valuation 
ring of rank one for every $\alpha\in I$. Let 
moreover $L$ be a finite {\'e}tale $K_\beta$-algebra, for some 
$\beta\in I$. Set $L_\alpha:=L\otimes_{K_\beta}K_\alpha$ for 
every morphism $\beta\to\alpha$ in $I$. Then :
\begin{enumerate}
\item
$(K_\infty,|\cdot|_\infty):=
\colim{\alpha\in I}(K_\alpha,|\cdot|_\alpha)$
is a valued field extension of $(K,|\cdot|)$, with valuation
ring of rank one.
\item
Set $L_\infty:=L\otimes_{K_\beta}K_\infty$. Then, for every
sequence of morphisms $\beta\to\gamma\to\alpha$ in $I$, we have 
$|\fd^+_{L_\alpha/K_\alpha}|\geq|\fd^+_{L_\gamma/K_\gamma}|$, 
and moreover:
$\displaystyle{\lim_{(\beta\to\alpha)\in I}}|\fd^+_{L_\alpha/K_\alpha}|=
|\fd^+_{L_\infty/K_\infty}|$.
\end{enumerate}
\end{lemma}
\begin{proof} (i) is obvious. The proof of the first assertion 
in (ii) is easy and shall be left to the reader. For the second
assertion in (ii) we remark that, due to
claim \ref{cl_finish.up}, for every $\eps<1$ there exists a
free $K_\infty^+$-submodule $W_\eps\subset W_{L_\infty}$ of finite 
type, such that $|F_0(W_{L_\infty}^a/W_\eps^a)|>\eps$. We can find 
$\alpha\in I$ such that $W_\eps=W_0\otimes_{K_\alpha^+}K^+_\infty$ 
for some free $K_\alpha^+$ submodule $W_0\subset L^+_\alpha$.
It then follows from lemma \ref{lem_Fitting.discr} that 
$$
|\fd^+_{L_\alpha/K_\alpha}|\geq
|\fd_{K^+_\alpha}(W_0,t_{L_\alpha/K_\alpha})|=
|\fd_{K_\infty^+}(W_\eps,t_{L_\infty/K_\infty})|>\eps^2\cdot
|\fd^+_{L_\infty/K_\infty}|
$$ 
for every morphism $\alpha\to\beta$ in $I$.
\end{proof}

\sset\subsubsection{}\label{subsec_situation}
Suppose now that $(K,|\cdot|)$ is deeply ramified and let 
$X:=(X_\alpha~|~\alpha\in I)$ be a deeply ramified 
cofiltered system of adic spaces of finite type over 
$\Spa(K,K^+)$. Suppose furthermore that it is given, for some 
$\beta\in I$, a finite {\'e}tale morphism  $Y_\beta\to X_\beta$ 
of adic $\Spa(K,K^+)$-spaces of finite type.
For every morphism $\alpha\to\beta$ of $I$ we set 
$Y_\alpha:=Y_\beta\times_{X_\beta}X_\alpha$ and
denote by $f_\alpha:Y_\alpha\to X_\alpha$ the 
induced morphism of adic spaces.

\begin{theorem}\label{th_defect} In the situation of 
\eqref{subsec_situation}, for every positive real number 
$\eps<1$ there exists a morphism $\alpha\to\beta$ in $I$ 
such that, for every morphism $\gamma\to\alpha$ we have 
$$
|\fd^+_{Y_\gamma/X_\gamma}(x)|>\eps\quad
\text{for every $x\in{\mathcal M}(X_\gamma)$}.
$$
\end{theorem}
\begin{proof} Notice that $\underline Y:=(Y_\alpha~|~\alpha\to\beta)$
is a cofiltered system, hence we can define $X_\infty$ and
$Y_\infty$ as in \eqref{subsec_go.to.extreme}, and we
obtain a morphism of locally ringed spaces 
$f_\infty:Y_\infty\to X_\infty$. For every $\alpha\in I$, 
let $\pi_\alpha:X_\infty\to X_\alpha$ be the natural
morphism. Moreover, let 
${\mathcal M}(X_\infty):=\liminv{\alpha\in I}~{\mathcal M}(X_\alpha)$; 
as a topological space, it is compact, by Tychonoff's theorem and the
fact that ${\mathcal M}(X_\alpha)$ is compact
\eqref{subsec_approach.semicont}); as a set, it admits an injective
(usually non-continuous) map ${\mathcal M}(X_\infty)\to X_\infty$, so
we can identify it as a subset of the latter.

Let $x\in{\mathcal M}(X_\infty)$; by proposition 
\ref{prop_to.extremes}, the valuation ring $\kappa(x)^+$
is deeply ramified. Set $\kappa(x):=\kappa(x)^+\otimes_{K^+}K$; 
it is clear that the morphism
$$
\kappa(x)\to 
E(x):=
(f_{\infty *}\cO_{Y_\infty}^+)_x\otimes_{\kappa(x)^+}\kappa(x)
$$
is finite and {\'e}tale. Let $W_x$ be the integral closure of
$\kappa(x)^+$ in $E(x)$. By proposition \ref{prop_deep.ramif} 
we deduce easily that the induced morphism of $K^{+a}$-algebras 
$\kappa(x)^{+a}\to W_x^a$ is weakly {\'e}tale, hence {\'e}tale
by proposition \ref{prop_fin.sep.uniform}. Consequently
$|\fd_{E(x)/\kappa(x)}^+|=1$, in light of lemma
\ref{lem_discr.criterio}. For every $\alpha\to\beta$, let
$x_\alpha:=\pi_\alpha(x)$. Then $\kappa(x)$ is the colimit 
of the filtered system $(\kappa(x_\alpha)~|~\alpha\to\beta)$, 
and similarly $E(x)$ is the colimit of the finite {\'e}tale 
$\kappa(x_\alpha)$-algebras $(f_{\beta *}\cO_{Y_\beta})_{x_\beta}
\otimes_{\cO^+_{X_\beta,x_\beta}}\kappa(x_\alpha)$
(for all $\alpha\to\beta$). In this situation, lemma 
\ref{lem_down.to.finite} applies and shows that, for every
$\eps<1$ there exists $\alpha(\eps,x)$ such that 
\set\begin{equation}\label{eq_yet.again.cpct}
|\fd^+_{Y_\alpha/X_\alpha}(x_\alpha)|>\eps\quad
\text{for every $\alpha\to\alpha(\eps,x)$}.
\end{equation}
In light of theorem \ref{th_disc.upper.semi.cont}, for every
$\alpha\to\beta$, the subset
$$
X_\alpha(\eps):=\{y\in{\mathcal M}(X_\alpha)~|~
|\fd^+_{Y_\alpha/X_\alpha}(y)|\leq\eps\}
$$
is closed in ${\mathcal M}(X_\alpha)$, hence compact.
From \eqref{eq_yet.again.cpct} we see that 
$\liminv{\alpha\to\beta}X_\alpha(\eps)=\emptyset$, therefore
one of the $X_\alpha(\eps)$ must be empty 
(\cite[Ch.I, \S 9, n.6, Prop.8(b)]{Bou}), and the claim
follows.
\end{proof}

\sset\subsubsection{}\label{subsec_defect}
\index{Adic space(s)!{\'e}tale morphism of!$\mathrm{def}(f)$ : defect of 
an|indref{subsec_defect}}
Let us choose an imbedding 
$\rho:(\Gamma_K,\cdot,\leq)\hookrightarrow(\R,+,\leq)$ as 
in example \ref{ex_tens.with.Q}(vi).
Let $f:Y\to X$ be a finite {\'e}tale morphism of adic spaces
of finite type over $\Spa(K,K^+)$. For every $x\in X$,
set $A(x):=(f_*\cO^+_Y)_x\otimes_{\cO^+_{X,x}}(f_*\cO^+_Y)_x$; 
we denote by $e_x\in A(x)\otimes_{K^+}K$ the unique idempotent 
characterized by the conditions of proposition \ref{prop_idemp}.
We define the {\em defect\/} of the morphism $f$ as the
real number 
$$
\mathrm{def}(f):=
\inf\{r\in\R_{\geq 0}~|~\eps\cdot e_x\in A(x)
\text{ for every $x\in X$ and every $\eps\in K^+$ with 
$\rho(|\eps|)\leq -r$ }\}.
$$
Clearly $\mathrm{def}(f)\geq 0$ and $\mathrm{def}(f)=0$
if and only if $(f_*\cO_Y^+)^a_x$ is an {\'e}tale 
$\cO_{X,x}^{+a}$-algebra for every $x\in X$. Furthermore
we remark that, by proposition \ref{prop_crit.projectivity}, 
the map $\cO_{X,x}^{+a}\to(f_*\cO_Y^+)^a_x$ is weakly {\'e}tale if 
and only if it is {\'e}tale.

\begin{corollary} In the situation of \eqref{subsec_situation},
for every real number $r>0$ there exists $\alpha\in I$ such 
that, for every morphism $\gamma\to\alpha$, we have:
$\mathrm{def}(f_\gamma)<r$.
\end{corollary}
\begin{proof} Let $r>0$; according to theorem \ref{th_defect}
and claim \ref{cl_liftup}, there exists $\alpha\in I$ such that
$\eps\cdot e_x\in A(x)$ for every $\gamma\to\alpha$, every
$x\in{\mathcal M}(X_\gamma)$ and every $\eps\in K^+$ with
$\rho(|\eps|)\leq -r$. If now $y\in X_\gamma$ is any point, there
is a unique generisation $x$ of $y$ in ${\mathcal M}(X_\gamma)$.
Let $\phi:A(y)\to A(x)$ be the induced specialisation map, and
set $\phi_K:=\phi\otimes_{K^+}\one_K$. One verifies easily that 
$$
\phi_K^{-1}(\fm\cdot A(x))\subset A(y).
$$
Since $\phi_K(e_y)=e_x$, the claim follows easily.
\end{proof}

\sset\subsubsection{}\label{subsec_explain}
To conclude, we want to explain briefly what kind of Galois 
cohomology calculations are enabled by the results of this 
section. Let $f:Y\to X$ be a finite {\'e}tale {\em Galois\/}
morphism of $\Spa(K,K^+)$-adic spaces of finite type, and
let $G$ denote the group of $X$-automorphisms of $Y$.
Denote by $f_*\cO_Y^+[G]\Mod$ the category of 
$f_*\cO_Y^+$-modules on $X$, endowed with a semilinear action 
of $G$. Let $\underline\Gamma^G:f_*\cO_Y^+[G]\Mod\to\cO_X^+\Mod$ 
be the functor that associates to an $f_*\cO_Y^+[G]$-module 
the sheaf of its $G$-invariant local sections. A standard 
argument shows that, for every $f_*\cO^+_Y[G]$-modules $\cF$, 
the cone of the natural morphism in $\sD(f_*\cO^+_Y\Mod)$
\set\begin{equation}\label{eq_climax}
\cO^+_Y\otimes_{\cO^+_X}R\underline\Gamma^G\cF\to\cF
\end{equation}
is annihilated by all $\eps\in\fm$ such 
that $\rho(\eps)<-\mathrm{def}(f)$. 
However, for applications one is rather more interested in 
understanding the Galois cohomology groups $H^i:=H^i(G,H^0(X,\cF))$.
One can try to study $H^i$ via \eqref{eq_climax}; indeed,
a bridge between these two objects is provided by the
higher derived functors of the related functor
$\Gamma^G:f_*\cO_Y^+[G]\Mod\to\Gamma(X,\cO^+)\Mod$, defined
by $\cF\mapsto\Gamma(X,\underline\Gamma^G\cF)=\Gamma(X,\cF)^G$.
We have two spectral sequences converging to $R\Gamma^G\cF$,
namely
$$
\begin{array}{r@{\:\Rightarrow\:}l}
E_2^{pq}:H^p(X,R^q\underline\Gamma^G\cF) & R^{p+q}\Gamma^G\cF \\
F_2^{pq}:H^p(G,H^q(X,\cF)) & R^{p+q}\Gamma^G\cF.
\end{array}
$$
Using \eqref{eq_climax} one deduces that $E_2^{pq}$ degenerates
up to some torsion, which can be estimated precisely in terms
of the defect of the morphism $f$. However, the spectral sequence 
$F_2^{pq}$ contains the terms $H^q(X,\cF)$, about which not much 
is currently known. In this direction, the only results
that we could found in the literature concern the calculation
of $H^i(Y,\cO^+_Y)$, for an affinoid space, under some very
restrictive assumptions : in \cite{Bart} these groups are shown
to be almost zero modules for $i>0$, in case $Y$ admits a smooth
formal model over $K^+$; in \cite{Put} the case of generalized
polydiscs is taken up, and the same kind of almost vanishing is
proven.


\newpage


\def\ic{\mathrm{i.c.}}
\def\Cov{\mathbf{Cov}}

\section{Appendix}

In this appendix we have gathered a few miscellaneous results
that were found in the course of our investigation, and which
may be useful for other applications.

\subsection{Simplicial almost algebras}\label{sec_simplicius}
\index{Simplicial almost algebra(s)|indref{sec_simplicius}}
\index{Simplicial almost algebra(s)!modules over a|indref{sec_simplicius}}
We need some preliminaries on simplicial objects : first of all,
a {\em simplicial almost algebra\/} is just an object in the category 
$s.(V^a\Alg)$. Then for a given simplicial almost algebra $A$ 
we have the category $A\Mod$ of $A$-modules : it 
consists of all simplicial almost $V$-modules $M$ such that 
$M[n]$ is an $A[n]$-module and such that the face and 
degeneracy morphisms $d_i:M[n]\to M[n-1]$ and 
$s_i:M[n]\to M[n+1]$ $(i=0,1,...,n)$ are $A[n]$-linear.  

\sset\subsubsection{}\label{subsec_simpl.der.catgory}
\index{Simplicial almost algebra(s)!modules over a!derived category of
|indref{subsec_simpl.der.catgory}}
We will need also the derived category of $A$-modules; it 
is defined as follows.
A bit more generally, let $\cC$\ \ be any abelian category. 
For an object $X$ of $s.\cC$\ \ let $N(X)$ be the normalized 
chain complex (defined as in \cite[I.1.3]{Il}). By the 
theorem of Dold-Kan (\cite[Th. 8.4.1]{We}) $X\mapsto N(X)$ 
induces an equivalence $N:s.\cC\to\sC_\bullet(\cC)$ with the 
category $\sC_\bullet(\cC)$ of chain complexes of object of 
$\cC$ that vanish in positive degrees.
Now we say that a morphism $X\to Y$ in $s.\cC$\ \ is a 
{\em quasi-isomorphism\/} if the induced morphism 
$N(X)\to N(Y)$ is a quasi-isomorphism of chain complexes.
\medskip
 
\sset\subsubsection{}\label{subsec_simpl_exact}
\index{Simplicial almost algebra(s)!exact|indref{subsec_simpl_exact}}
\index{Simplicial almost algebra(s)!modules over a!quasi-isomorphism of
|indref{subsec_simpl_exact}}
In the following we fix a simplicial almost algebra $A$.
\begin{definition} We say that $A$ is {\em exact\/} if the 
almost algebras $A[n]$ are exact for all $n\in\N$. A morphism 
$\phi:M\to N$ of $A$-modules (or $A$-algebras) 
is a {\em quasi-isomorphism\/} if the morphism $\phi$ of underlying 
simplicial almost $V$-modules is a quasi-isomorphism. We define the 
category $\sD(A\Mod)$ (resp. the category $\sD(A\Alg)$) 
as the localization of the category $A\Mod$ (resp. $A\Alg$) with 
respect to the class of quasi-isomorphisms. 
\end{definition}

\sset\subsubsection{}
As usual, the morphisms in $\sD(A\Mod)$ can be computed via
a calculus of fraction on the category $\Hot(A)$
of simplicial complexes up to homotopy. Moreover, if $A_1$
and $A_2$ are two simplicial almost algebras, then the 
``extension of scalars'' functors define equivalences of
categories
$$\begin{array}{l}
\sD(A_1\times A_2\Mod)\stackrel{\sim}{\longrightarrow}
\sD(A_1\Mod)\times\sD(A_2\Mod) \\
\sD(A_1\times A_2\Alg)\stackrel{\sim}{\longrightarrow}
\sD(A_1\Alg)\times\sD(A_2\Alg).
\end{array}$$
\begin{proposition}\label{prop_preserve} 
Let $A$ be a simplicial $V^a$-algebra.
\begin{enumerate}
\item
The functor on 
$A$-algebras given by $B\mapsto(s.V^a\times B)_{!!}$ preserves 
quasi-isomorphisms and therefore induces a functor 
$\sD(A\Alg)\to\sD((s.V^a\times A)_{!!}\Alg)$.
\item
The localisation functor $R\mapsto R^a$ followed by ``extension
of scalars'' via $s.V^a\times A\to A$ induces a functor
$\sD((s.V^a\times A)_{!!}\Alg)\to\sD(A\Alg)$
and the composition of this and the above functor is naturally
isomorphic to the identity functor on $\sD(A\Alg)$.
\end{enumerate}
\end{proposition}
\begin{proof} (i) : let $B\to C$ be a quasi-isomorphism 
of $A$-algebras. Clearly the induced morphism 
$s.V^a\times B\to s.V^a\times C$ is still a quasi-isomorphism 
of $V$-algebras. But by remark \ref{rem_exact.alg}, 
$s.V^a\times B$ and $s.V^a\times C$ are exact simplicial 
almost $V$-algebras; moreover, it follows from corollary 
\ref{cor_left.exact} that $(s.V^a\times B)_!\to(s.V^a\times C)_!$ 
is a quasi-isomorphism of $V$-modules. Then the claim 
follows easily from the exactness of the sequence 
\eqref{eq_left.adj.algebras}. Now (ii) is clear.
\end{proof}

\begin{remark} In case $\fm$ is flat, then all
$A$-algebras are exact, and the same argument shows that
the functor $B\mapsto B_{!!}$ induces a functor 
$\sD(A\Alg)\to\sD(A_{!!}\Alg)$.
In this case, composition with localisation is naturally
isomorphic to the identity functor on $\sD(A\Alg)$.
\end{remark}

\begin{proposition}\label{prop_Gabber} 
Let $f:R\to S$ be a map of $V$-algebras. We have:
\begin{enumerate}
\item
If $f^a:R^a\to S^a$ is an isomorphism, then $\L_{S/R}^a\simeq 0$
in $\sD(s.S^a\Mod)$.
\item
The natural map $\L^a_{S^a/R^a}\to\L^a_{S/R}$ is an isomorphism
in $\sD(s.S^a\Mod)$.
\end{enumerate}
\end{proposition}
\begin{proof} (ii) is an easy consequence of (i) and of transitivity.
To show (i) we prove by induction on $q$ that 

\medskip

\noindent{\bf VAN}($q;S/R$){\hskip 2.5cm} $H_q(\L_{S/R}^a)=0.$
\medskip

For $q=0$ the claim follows immediately from \cite[II.1.2.4.2]{Il}.
Therefore suppose that $q>0$ and that {\bf VAN}($j;D/C$) is known for
all almost isomorphisms of $V$-algebras $C\to D$ and all $j<q$.
Let $\bar R:=f(R)$. Then by transitivity (\cite[II.2.1.2]{Il}) 
we have a distinguished triangle in $\sD(s.S^a\Mod)$
$$
(S\otimes_{\bar R}\L_{\bar R/R})^a \stackrel{u}{\to} 
\L^a_{S/R} \stackrel{v}{\to} \L^a_{S/\bar R} \to  
\sigma(S\otimes_{\bar R}\L_{\bar R/R})^a.
$$
We deduce that {\bf VAN}($q;\bar R/R$) and {\bf VAN}($q;S/\bar R$) imply
{\bf VAN}($q;S/R$), thus we can assume that $f$ is either injective or
surjective. Let $S_\bullet\to S$ be the simplicial $V$-algebra
augmented over $S$ defined by $S_\bullet:=P_V(S)$. It is a simplicial
resolution of $S$ by free $V$-algebras, in particular the augmentation
is a quasi-isomorphism of simplicial $V$-algebras. Set 
$R_\bullet:=S_\bullet\times_SR$. This is a simplicial $V$-algebra
augmented over $R$ via a quasi-isomorphism. Moreover, the induced
morphisms $R[n]^a\to S[n]^a$ are isomorphisms. By \cite[II.1.2.6.2]{Il}
there is a quasi-isomorphism 
$\L_{S/R}\simeq\L^\Delta_{S_\bullet/R_\bullet}$.
On the other hand we have a spectral sequence
$$E^1_{ij}:=H_j(\L_{S[i]/R[i]})\Rightarrow 
H_{i+j}(\L^\Delta_{S_\bullet/R_\bullet}).$$
It follows easily that 
{\bf VAN}$(j;S[i]/R[i])$ for all $i\ge 0,j\le q$
implies {\bf VAN}$(q;S/R)$.
Therefore we are reduced to the case where
$S$ is a free $V$-algebra and $f$ is either injective or surjective.
We examine separately these two cases. If $f:R\to V[T]$ is surjective,
then we can find a right inverse $s:V[T]\to R$ for $f$. By applying
transitivity to the sequence $V[T]\to R\to V[T]$ we get a distinguished
triangle
$$(V[T]\otimes_R\L_{R/V[T]})^a \stackrel{u}{\to}\L^a_{V[T]/V[T]}
\stackrel{v}{\to}\L^a_{V[T]/R}\to\sigma(V[T]\otimes_R\L_{R/V[T]})^a.$$
Since $\L^a_{V[T]/V[T]}\simeq 0$ there follows an isomorphism :
$H_q(\L_{V[T]/R})^a\simeq H_{q-1}(V[T]\otimes_R\L_{R/V[T]})^a$.
Furthermore, since $f^a$ is an isomorphism, $s^a$ is an isomorphism as well,
hence by induction (and by a spectral sequence of the type
\cite[I.3.3.3.2]{Il}) $H_{q-1}(V[T]\otimes_R\L_{R/V[T]})^a\simeq 0$.
The claim follows in this case.

Finally suppose that $f:R\to V[T]$ is injective. Write $V[T]=\Sym(F)$,
for a free $V$-module $F$ and set $\tilde F=\tilde\fm\otimes_VF$; since 
$f^a$ is an isomorphism, 
$\Img(\Sym(\tilde F)\to\Sym(F))\subset R$.
We apply transitivity to the sequence 
$\Sym(\tilde F)\to R\to\Sym(F)$. By arguing as above 
we are reduced to showing that 
$\L^a_{\Sym(F)/\Sym(\tilde F)}\simeq 0.$ We know that 
$H_0(\L^a_{\Sym(F)/\Sym(\tilde F)})\simeq 0$ and we will 
show that  
$H_q(\L^a_{\Sym(F)/\Sym(\tilde F)})\simeq 0$ for $q>0$. 
To this purpose we apply transitivity to the sequence 
$V\to\Sym(\tilde F)\to\Sym(F)$. As $F$ and 
$\tilde F$ are flat $V$-modules, \cite[II.1.2.4.4]{Il} 
yields $H_q(\L_{\Sym(F)/V})\simeq
H_q(\L_{\Sym(\tilde F)/V})\simeq 0$ for $q>0$ and 
$H_0(\L_{\Sym(\tilde F)/V})$ is a flat 
$\Sym(\tilde F)$-module.
In particular 
$H_j(\Sym(F)\otimes_{\Sym(\tilde F)}\L_{\Sym(\tilde F)/V})
\simeq 0$ for all $j>0$. 
Consequently $H_{j+1}(\L_{\Sym(F)/\Sym(\tilde F)})\simeq 0$ 
for all $j>0$ and 
$H_1(\L_{\Sym(F)/\Sym(\tilde F)})\simeq
\Ker(\Sym(F)\otimes_{\Sym(\tilde F)}
\Omega_{\Sym(\tilde F)/V}\to\Omega_{\Sym(F)/V})$.
The latter module is easily seen to be almost zero.
\end{proof}

\begin{theorem}\label{th_Gabber} Let $\phi:R\to S$ be a map 
of simplicial $V$-algebras inducing an isomorphism 
$R^a\stackrel{\sim}{\to}S^a$ in $\sD(R^a\Mod)$. 
Then $(\L^\Delta_{S/R})^a\simeq 0$ in $\sD(S^a\Mod)$.
\end{theorem}
\begin{proof} Apply the base change theorem \cite[II.2.2.1]{Il}
to the (flat) projections of $s.V\times R$ onto $R$ and respectively
$s.V$ to deduce that the natural map 
$\L^\Delta_{s.V\times S/s.V\times R}\to
\L^\Delta_{S/R}\oplus\L^\Delta_{s.V/s.V}\to\L^\Delta_{S/R}$ is a 
quasi-isomorphism in $\sD(s.V\times S\Mod)$. By proposition 
\ref{prop_preserve} the induced morphism 
$(s.V\times R)^a_{!!}\to(s.V\times S)^a_{!!}$ is still a 
quasi-isomorphism. There are spectral sequences
$$\begin{array}{c}
E^1_{ij}:=H_j(\L_{(V\times R[i])/(V\times R[i])_{!!}^a})\Rightarrow 
H_{i+j}(\L^\Delta_{(s.V\times R)/(s.V\times R)_{!!}^a}) \\
F^1_{ij}:=H_j(\L_{(V\times S[i])/(V\times S[i])_{!!}^a})\Rightarrow 
H_{i+j}(\L^\Delta_{(s.V\times S)/(s.V\times S)_{!!}^a}).
\end{array}$$
On the other hand, by proposition \ref{prop_Gabber}(i) we
have
$\L^a_{(V\times R[i])/(V\times R[i])_{!!}^a}\simeq 0
\simeq\L^a_{(V\times S[i])/(V\times S[i])_{!!}^a}$
for all $i\in\N$. Then the theorem follows directly from 
\cite[II.1.2.6.2(b)]{Il} and transitivity. 
\end{proof}

\sset\subsubsection{}
In view of proposition \ref{prop_Gabber}(i) we have  
$\L^a_{(V^a\times A)_{!!}/V\times A_{!!}}\simeq 0$
in $\sD(V^a\times A\Mod)$. By this, transitivity 
and localisation (\cite[II.2.3.1.1]{Il}) we derive that 
$\L^a_{B/A}\to\L^a_{B_{!!}/A_{!!}}$ is a quasi-isomorphism 
for all $A$-algebras $B$. If $A$ and $B$ are exact 
({\em e.g.\/} if $\fm$ is flat), we conclude from proposition 
\ref{prop_quasi.exact} that the natural map 
$\L_{B/A}\to\L_{B_{!!}/A_{!!}}$ is a quasi-isomorphism.

\sset\subsubsection{}
Finally we want to discuss left derived functors of (the almost version
of) some notable non-additive functors that play a role in deformation
theory. Let $R$ be a simplicial $V$-algebra. Then we have an obvious 
functor $G:\sD(R\Mod)\to\sD(R^a\Mod)$ obtained by applying 
dimension-wise the localisation functor. Let $\Sigma$ be the 
multiplicative set of morphisms of $\sD(R\Mod)$ that induce almost 
isomorphisms on the cohomology modules. An argument as in section 
\ref{sec_homol} shows that $G$ induces an equivalence of categories
$\Sigma^{-1}\sD(R\Mod)\to\sD(R^a\Mod)$.

\sset\subsubsection{}
Now let $R$ be a $V$-algebra and $\cF_p$ one of the functors 
$\otimes^p$, $\Lambda^p$, $\Sym^p$, $\Gamma^p$ defined in 
\cite[I.4.2.2.6]{Il}.
\begin{lemma}\label{lem_cF} Let $\phi:M\to N$ be an almost 
isomorphism of $R$-modules. Then $\cF_p(\phi):\cF_p(M)\to\cF_p(N)$ 
is an almost isomorphism.
\end{lemma}
\begin{proof} Let $\psi:\tilde\fm\otimes_VN\to M$ be the 
map corresponding to $(\phi^a)^{-1}$ under the bijection
\eqref{eq_alm.morph}. By inspection, the compositions 
$\phi\circ\psi:\tilde\fm\otimes_VN\to N$ and 
$\psi\circ(\one_{\tilde\fm}\otimes\phi):\tilde\fm\otimes_VM\to M$ 
are induced by scalar multiplication. Pick any $s\in\fm$ and lift 
it to an element $\tilde s\in\tilde\fm$; define $\psi_s:N\to M$
by $n\mapsto\psi(\tilde s\otimes n)$ for all $n\in N$.
Then $\phi\circ\psi_s=s\cdot\one_N$ and 
$\psi_s\circ\phi=s\cdot\one_M$. This easily implies that 
$s^p$ annihilates $\Ker\,\cF_p(\phi)$ and $\Coker\,\cF_p(\phi)$.
In light of proposition \ref{prop_less.obv}(ii), the claim 
follows. 
\end{proof}

\sset\subsubsection{}\label{subsec_app.gammas.etc}
Let $B$ be an almost $V$-algebra. We define a functor 
$\cF_p^a$ on $B\Mod$ by letting $M\mapsto(\cF_p(M_!))^a$, 
where $M_!$ is viewed as a $B_{!!}$-module or a $B_*$-module 
(to show that these choices define the same functor it 
suffices to observe that $B_*\otimes_{B_{!!}}N\simeq N$ 
for all $B_*$-modules $N$ such that $N=\fm\cdot N$).
For all $p>0$ we have diagrams :
\set\begin{equation}\label{eq_derotimes}{
\diagram
R\Mod \ar[r]^{\cF_p} \ar@<.5ex>[d] & R\Mod \ar@<.5ex>[d] \\
R^a\Mod \ar[r]^{\cF_p^a} \ar@<.5ex>[u] & 
R^a\Mod \ar@<.5ex>[u]
\enddiagram}\end{equation}
where the downward arrows are localisation and the upward
arrows are the functors $M\mapsto M_!$.
Lemma \ref{lem_cF} implies that the downward arrows in the
diagram commute (up to a natural isomorphism) with the 
horizontal ones. It will follow from the following proposition 
\ref{prop_upward} that the diagram commutes also going upward. 

\sset\subsubsection{}
For any $V$-module $N$ we have an exact sequence 
$\Gamma^2N\to\otimes^2N\to\Lambda^2N\to 0$. As observed in
the proof of proposition \ref{prop_less.obv}, the symmetric
group $S_2$ acts trivially on $\otimes^2\tilde\fm$ and
$\Gamma^2\tilde\fm\simeq\otimes^2\tilde\fm$, so 
$\Lambda^2\tilde\fm=0$. Also we have natural isomorphisms
$\Gamma^p\tilde\fm\simeq\tilde\fm$ for all $p>0$.
\begin{proposition}\label{prop_upward}
Let $R$ be a commutative ring and $L$ a flat $R$-module 
with $\Lambda^2L=0$. Then for $p>0$ and for all $R$-modules
$N$ we have natural isomorphisms
$$\Gamma^p(L)\otimes_R\cF_p(N)
\stackrel{\sim}{\to}{\cF_p}(L\otimes_RN).$$
\end{proposition}
\begin{proof} Fix an element $x\in\cF_p(N)$. 
For each $R$-algebra $R'$ and each element $l\in R'\otimes_RL$
we get a map $\phi_l:R'\otimes_RN\to R'\otimes_RL\otimes_RN$
by $y\mapsto l\otimes y$, hence a map 
$\cF_p(\phi_l):R'\otimes_R\cF_p(N)\simeq\cF_p(R'\otimes_RN)\to
\cF_p(R'\otimes_RL\otimes_RN)\simeq 
R'\otimes_R\cF_p(L\otimes_RN)$. For varying $l$ we obtain
a map of sets $\psi_{R',x}:R'\otimes_RL\to 
R'\otimes_R\cF_p(L\otimes_RN)$ : 
$l\mapsto\cF_p(\phi_l)(1\otimes x)$. According to the 
terminology of \cite{Ro}, the system of maps $\psi_{R',x}$ 
for $R'$ ranging over all $R$-algebras forms a homogeneous 
polynomial law of degree $p$ from $L$ to $\cF_p(L\otimes_RN)$, 
so it factors through the universal homogeneous degree $p$ 
polynomial law $\gamma_p:L\to\Gamma^p(L)$ . The resulting 
$R$-linear map $\bar\psi_x:\Gamma^p(L)\to\cF_p(L\otimes_RN)$ 
depends $R$-linearly on $x$, hence we derive an $R$-linear 
map $\psi:\Gamma^p(L)\otimes_R\cF_p(N)\to{\cF_p}(L\otimes_RN)$.
Next notice that by hypothesis $S_2$ acts trivially on $\otimes^2L$ 
so $S_p$ acts trivially on 
$\otimes^pL$ and we get an isomorphism
$\beta:\Gamma^p(L)\stackrel{\sim}{\longrightarrow}\otimes^pL$.
We deduce a natural map 
$(\otimes^pL)\otimes_R\cF_p(N)\to\cF_p(L\otimes_RN)$.
Now, in order to prove the proposition for the case 
$\cF_p=\otimes^p$, it suffices to show that this last 
map is just the natural isomorphism that ``reorders
the factors''. Indeed, let $x_1,...,x_n\in L$ and 
$q:=(q_1,...,q_n)\in\N^n$ such that $|q|:=\sum_iq_i:=p$;
then $\beta$ sends the generator 
$x_1^{[q_1]}\cdot...\cdot x_n^{[q_n]}$ to 
$\binom{p}{q_1,...,q_n}\cdot 
x_1^{\otimes q_1}\otimes...\otimes x_n^{\otimes q_n}$.
On the other hand, pick any $y\in\otimes^pN$ and let 
$R[T]:=R[T_1,...,T_n]$ be the polynomial $R$-algebra in $n$ 
variables; write 
$(T_1\otimes x_1+...+T_n\otimes x_n)^{\otimes p}\otimes y=
\psi_{R[T],y}(T_1\otimes x_1+...+T_n\otimes x_n)=
\sum_{r\in\N^n}T^r\otimes w_r$ with $w_r\in\otimes^p(L\otimes_RN)$.
Then $\psi((x_1^{[q_1]}\cdot...\cdot x_n^{[q_n]})\otimes y)=w_q$
(see \cite[pp.266-267]{Ro})
and the claim follows easily. Next notice 
that $\Gamma^p(L)$ is flat, so that tensoring with $\Gamma^p(L)$ 
commutes with taking coinvariants (resp. invariants) under the 
action of the symmetric group; this implies the assertion for 
$\cF_p:=\Sym^p$ (resp. $\cF_p:=\mathrm{TS}^p$). 
To deal with $\cF_p:=\Lambda^p$ recall that for any $V$-module 
$M$ and $p>0$ we have the antisymmetrizer operator 
$a_M:=\sum_{\sigma\in S_p}\mathrm{sgn}(\sigma)\cdot\sigma:
\otimes^pM\to\otimes^pM$ and a surjection $\Lambda^p(M)\to\Img(a_M)$ 
which is an isomorphism for $M$ free, hence for $M$ flat.
The result for $\cF_p=\otimes^p$ (and again the flatness of 
$\Gamma^p(L)$) then gives 
$\Gamma^p(L)\otimes\Img(a_N)\simeq\Img(a_{L\otimes_RN})$,
hence the assertion for $\cF_p=\Lambda^p$ and $N$ flat.
For general $N$ let $F_1\stackrel{\partial}{\to}
F_0\stackrel{\eps}{\to}N\to 0$ be a presentation with $F_i$
free. Define $j_0,j_1:F_0\oplus F_1\to F_0$ by 
$j_0(x,y):=x+\partial(y)$ and $j_1(x,y):=x$. By functoriality
we derive an exact sequence 
$$\xymatrix{
\Lambda^p(F_0\oplus F_1) \ar@<.5ex>[r] \ar@<-.5ex>[r] & 
\Lambda^p(F_0) \ar[r] & \Lambda^p(N) \ar[r] & 0
}$$
which reduces the assertion to the flat case. For $\cF_p:=\Gamma^p$
the same reduction argument works as well (cp. \cite[p.284]{Ro})
and for flat modules the assertion for $\Gamma^p$ follows
from the corresponding assertion for $\mathrm{TS}^p$.
\end{proof}
\begin{lemma} Let $A$ be a simplicial almost algebra, $L,E$ 
and $F$ three $A$-modules, $f:E\to F$ a quasi-isomorphism. 
If $L$ is flat or $E,F$ are flat, then 
$L\otimes_Af:L\otimes_AE\to L\otimes_AF$ is a quasi-isomorphism.
\end{lemma}
\begin{proof} It is deduced directly from \cite[I.3.3.2.1]{Il}
by applying $M\mapsto M_!$.
\end{proof}

\sset\subsubsection{}
As usual, this allows one to show that 
$\otimes:\Hot(A)\times\Hot(A)\to
\Hot(A)$ admits a left derived functor 
$\derotimes:\sD(A\Mod)\times\sD(A\Mod)\to\sD(A\Mod)$.
If $R$ is a simplicial $V$-algebra then we have essentially
commutative diagrams
$$\xymatrix{
\sD(R\Mod)\times\sD(R\Mod) \ar[r]^-{\derotimes} \ar@<.5ex>[d] &
\sD(R\Mod) \ar@<.5ex>[d] \\
\sD(R^a\Mod)\times\sD(R^a\Mod)\ar[r]^-{\derotimes}\ar@<.5ex>[u]&
\sD(R^a\Mod) \ar@<.5ex>[u]
}$$
where again the downward (resp. upward) functors are induced by 
localisation (resp. by $M\mapsto M_!$).

\sset\subsubsection{}
We mention the derived functors of the non-additive functor $\cF_p$ 
defined above in the simplest case of modules over a constant
simplicial ring. Let $A$ be a (commutative) $V^a$-algebra.
\begin{lemma} If $u:X\to Y$ is a quasi-isomorphism of flat 
$s.A$-modules then $\cF_p^a(u):\cF_p^a(X)\to\cF_p^a(Y)$ is 
a quasi-isomorphism.
\end{lemma}
\begin{proof} This is deduced from \cite[I.4.2.2.1]{Il} applied
to $N(X_!)\to N(Y_!)$ which is a quasi-isomorphism of chain complexes
of flat $A_{!!}$-modules. We note that {\em loc. cit.\/} deals 
with a more general mixed simplicial construction of $\cF_p$ which
applies to bounded above complexes, but one can check that it
reduces to the simplicial definition for complexes in 
$\cC_\bullet(A_{!!})$.
\end{proof}

\sset\subsubsection{}
Using the lemma one can construct 
$L\cF_p^a:\sD(s.A\Mod)\to\sD(s.A\Mod)$. If $R$ is a
$V$-algebra we have the derived category version of the 
essentially commutative squares \eqref{eq_derotimes}, 
relating $L\cF_p:\sD(s.R\Mod)\to\sD(s.R\Mod)$
and $L\cF_p^a:\sD(s.R^a\Mod)\to\sD(s.R^a\Mod)$.

\subsection{Fundamental group of an almost algebra}\label{sec_Galois}
\index{Topos|indref{sec_Galois}}
We will need some generalities from \cite[Exp.V]{SGA1} and 
\cite[Exp.VI]{SGA4-2}. In the following we fix a universe $U$
and suppose that all our categories are $U$-categories and all
our topoi are $U$-topoi in the sense of 
\cite[Exp.I, D{\'e}f.1.1 and Exp.IV, D{\'e}f.1.1]{SGA4-1}. No further 
mention of universes will be necessary.

\sset\subsubsection{}\label{subsec_qc.on.site}
\index{Topos!quasi-compact object of a|indref{subsec_qc.on.site}}
Let $\cC$ be a site. Recall (\cite[Exp.VI, D{\'e}f.1.1]{SGA4-2})
that an object $X$ of $\cC$ is called {\em quasi-compact\/} if, 
for every covering family $(X_i\to X~|~i\in I)$ there is a 
finite subset $J\subset I$ such that the subfamily 
$(X_j\to X~|~j\in J)$ is still covering.

\sset\subsubsection{}\label{subsec_qc.on.topos}
\index{Topos!quasi-compact, quasi-separated, 
coherent morphism in a|indref{subsec_qc.on.topos}}
Let $E$ be a topos; in the following we always endow $E$ with 
its canonical topology (\cite[Exp.II, D{\'e}f 2.5]{SGA4-1}), 
so $E$ is a site in a natural way and the terminology of 
\eqref{subsec_qc.on.site} applies to the objects of $E$.
Moreover, if $\cC$ is any site and $\eps:\cC\to\cC^\sim$ 
the natural functor to the category $\cC^\sim$ of sheaves 
on $\cC$, then an object $X$ of $\cC$ is quasi-compact in 
$\cC$ if and only if $\eps(X)$ is quasi-compact in $\cC^\sim$ 
(\cite[Exp.VI, Prop.1.2]{SGA4-2}).

Furthermore, since in $E$ all finite limits are representable, 
we can make the following further definitions
(\cite[Exp.VI, D{\'e}f.1.7]{SGA4-2}). A morphism $f:X\to Y$ 
in $E$ is called {\em quasi-compact\/} if, for every morphism 
$Y'\to Y$ in $E$ with quasi-compact $Y'$, the object 
$X\times_YY'$ is quasi-compact. We say that $f$ is 
{\em quasi-separated\/} if the diagonal morphism 
$X\to X\times_YX$ is quasi-compact. We say that $f$ is 
{\em coherent\/} if it is quasi-compact and quasi-separated.

\sset\subsubsection{}\label{subsec_construct}
\index{Topos!quasi-separated, coherent object of a|indref{subsec_construct}}
Let $X$ be an object of a topos $E$. We say that $X$ is
{\em quasi-separated\/} if, for every quasi-compact object
$S$ of $E$, every morphism $S\to X$ is quasi-compact.
We say that $X$ is {\em coherent\/} if it is quasi-compact 
and quasi-separated (\cite[Exp.VI, D{\'e}f.1.13]{SGA4-2}).
We denote by $E_\mathrm{coh}$ the full subcategory of $E$
consisting of all the coherent objects.

Suppose that the object $Y$ of $E$ is coherent and let
$f:X\to Y$ be a coherent morphism; by 
\cite[Exp.VI, Prop.1.14(ii)]{SGA4-2}, it follows that $X$ 
is coherent.

\begin{definition}\label{def_coh.topos} (cp. \cite[Exp.VI, D{\'e}f.2.3]{SGA4-2})
\index{Topos!coherent|indref{def_coh.topos}}
We say that a topos $E$ is {\em coherent\/} if it satisfies the
following conditions:
\begin{enumerate}
\item
$E$ admits a full generating subcategory $\cC$ consisting of
coherent objects.
\item
Every object $X$ of $\cC$ is quasi-separated over the final
object of $E$, {\em i.e.\/} the diagonal morphism 
$X\to X\times X$ is quasi-compact.
\item
The final object of $E$ is quasi-separated.
\end{enumerate}
\end{definition}

\sset\subsubsection{}
If $E$ is a coherent topos, then $E_\mathrm{coh}$ is 
stable under arbitrary finite limits (of $E$) 
(\cite[Exp.VI, 2.2.4]{SGA4-2}).
Moreover, a topos $E$ is coherent if and only if it is 
equivalent to a topos of the form $\cC^\sim$, where $\cC$ 
is a site whose objects are quasi-compact and whose finite 
limits are representable.

It is possible to characterize the categories of the form
$E_\mathrm{coh}$ arising from a coherent topos : this leads
to the following definition.

\begin{definition}\label{def_pretopos}
\index{Pretopos|indref{def_pretopos}}
A {\em pretopos\/} is a small category $\cC$ satisfying the
following conditions (\cite[Exp.VI, Exerc.3.11]{SGA4-2}).
\begin{itemize}
\item[(PT1)]
All finite limits are representable in $\cC$.
\item[(PT2)]
All finite sums are representable in $\cC$ and they are
universal and disjoint.
\item[(PT3)]
Every equivalence relation in $\cC$ is effective, and
every epimorphism is universally effective.
\end{itemize}
\end{definition}

\sset\subsubsection{}\label{subsec_pretopos}
As in \cite[Exp.VI, Exerc.3.11]{SGA4-2}, we leave to the reader
the verification that, for every coherent topos $E$, the
subcategory $E_\mathrm{coh}$ is a pretopos, and $E$ induces 
on $E_\mathrm{coh}$ the {\em precanonical\/} topology,
{\em i.e.} the topology whose covering families 
$(X_i\to X~|~i\in I)$
are those admitting a finite subfamily which is covering
for the canonical topology of $E_\mathrm{coh}$. One deduces
that $E$ is equivalent to $(E_\mathrm{coh})^\sim$,
the topos of sheaves on the precanonical topology
of $E_\mathrm{coh}$. 

Conversely, if $\cC$ is a pretopos,
let $E:=\cC^\sim$ be the topos of sheaves on the precanonical
topology of $\cC$; then $E$ is a coherent topos and the 
natural functor $\eps:\cC\to E$ induces an equivalence
of $\cC$ with $E_\mathrm{coh}$.

\sset\subsubsection{}\label{subsec_precanon}
Furthermore, if $\cC$ is a pretopos (endowed with the precanonical
topology), the natural functor $\cC\to\cC^\sim$ commutes with finite 
sums, with quotients under equivalence relations, and it is left exact 
({\em i.e.} commutes with finite limits) (\cite[Exp.VI, Exerc.3.11]{SGA4-2}).

\sset\subsubsection{}\label{subsec_point-topos}
\index{Topos!point of a|indref{subsec_point-topos}}
Recall (\cite[Exp.IV, D{\'e}f.6.1]{SGA4-1}) that a {\em point\/}
of a topos $E$ is a morphism of topoi $p:\Set\to E$ (where one
views the category $\Set$ as the topos of sheaves on the
one-point topological space). By \cite[Exp.IV, Cor.1.5]{SGA4-1},
the assignment $p\mapsto p^*$ defines an equivalence from
the category of points of $E$ to the opposite of the category 
of all functors $F:E\to\Set$ that commute with all colimits and
are left exact.
A functor $F:E\to\Set$ with these properties is called a 
{\em fibre functor\/} for $E$. By \cite[Exp.IV, Cor.1.7]{SGA4-1},
a functor $E\to\Set$ is a fibre functor if and only if it
is left exact and it takes covering families to surjective
families.

\sset\subsubsection{}\label{subsec_Deligne}
By a theorem of Deligne (\cite[Exp.VI, Prop.9.0]{SGA4-2})
every coherent non-empty topos admits at least a fibre functor. 
(Actually Deligne's result is both more precise and more general, 
but for our purposes, the foregoing statement will suffice).

In several contexts, it is useful to attach fibre functors to
categories that are not quite topoi. These situations are
axiomatized in the following definition.

\begin{definition}\label{def_Gal-category} A {\em Galois category\/}
\index{$\cC$ : Category(ies)!Galois|indref{def_Gal-category}}
(\cite[Exp.V, \S5]{SGA1}) is the datum of a 
category $\cC$ and a functor $F$ from $\cC$ to the category
$\fSet$ of finite sets, satisfying the following conditions:
\begin{itemize}
\item[(G1)]
all finite limits exist in $\cC$ (in particular
$\cC$ has a final object).
\item[(G2)]
Finite sums exist in $\cC$ (in particular $\cC$ has an initial
object). Also, for every object $X$ of $\cC$ and every
finite group $G$ of automorphisms of $X$, the quotient $X/G$ 
exists in $\cC$.
\item[(G3)]
Every morphism $u:X\to Y$ in $\cC$ factors as a composition
$X\stackrel{u'}{\to}Y'\stackrel{u''}{\to}Y$, where $u'$ is a
strict epimorphism and $u''$ is both a monomorphism and an 
isomorphism on a direct summand of $Y$.
\item[(G4)]
The functor $F$ is left exact.
\item[(G5)]
$F$ commutes with finite direct sums and with 
quotients under actions of finite groups of automorphisms. 
Moreover $F$ takes strict epimorphisms to epimorphisms.
\item[(G6)]
Let $u:X\to Y$ be a morphism of $\cC$. Then $u$ is an 
isomorphism if and only if $F(u)$ is.
\end{itemize}
\end{definition}

\sset\subsubsection{}\label{subsec_Galois.grp}
\index{$\cC$ : Category(ies)!Galois!fibre functor of 
a|indref{subsec_Galois.grp}}
Given a Galois category $(\cC,F)$, one says that $F$ is a 
{\em fibre functor\/} of $\cC$.
It is shown in \cite[Exp.V, \S4]{SGA1} that for any Galois
category $(\cC,F)$ the functor $F$ is pro-representable
and its automorphism group is a profinite group $\pi$ in
a natural way. Furthermore, $F$ factors naturally through 
the forgetful functor $\pi\text{-}\fSet\to\fSet$ from the 
category $\pi\text{-}\fSet$ of (discrete) finite sets with 
continuous $\pi$-action, and the resulting functor 
$\cC\to\pi\text{-}\fSet$ is an equivalence.

\sset\subsubsection{}
The category of {\'e}tale coverings $\Cov(A)$ of an 
almost algebra $A$ (to be defined in \eqref{subsec_covers})
is not directly presented as a Galois 
category, since it does not afford an a priori choice of 
fibre functor; rather, the existence of a fibre functor is 
deduced from Deligne's theorem. The argument only appeals to 
some general properties of the category $\Cov(A)$, 
which are abstracted in the following definition 
\ref{def_pregalois} and lemma \ref{lem_pregalois}.

\begin{definition}\label{def_pregalois}
\index{$\cC$ : Category(ies)!Pregalois|indref{def_pregalois}}
A {\em pregalois category\/} is a category $\cC$ satisfying
the following conditions.
\begin{itemize}
\item[(PG1)]
Every monomorphism $X\to Y$ in $\cC$ induces an isomorphism
of $X$ onto a direct summand of $Y$.
\item[(PG2)]
$\cC$ admits a final object $e$ which is connected and non-empty 
(that is, $e$ is not an initial object).
\item[(PG3)]
For every object $X$ of $\cC$, there exists $n\in\N$ such
that, for every non-empty object $Y$ of $\cC$, the product 
$X\times Y$ exists in $\cC$ and is not $Y$-isomorphic to an 
object of the form $\displaystyle{(\mathop{\amalg}^n Y)\amalg Z}$ 
(where $Z$ is any other object of $\cC$).
\end{itemize}
\end{definition}

\begin{lemma}\label{lem_pregalois}
Let $\cC$ be a pregalois pretopos. Then there exists a 
functor $F:\cC\to\fSet$ such that $(\cC,F)$ is a Galois category.
\end{lemma}
\begin{proof} (G1) holds because it is the same as (PT1).
(G2) follows easily from (PT2) and (PT3). In order to show 
(G3) we will need the following:
\begin{claim}\label{cl_mono.epi} A morphism $u:X\to Y$ of 
$\cC$ is an isomorphism if and only if it is both a 
monomorphism and an epimorphism.
\end{claim}
\begin{pfclaim} One direction is obvious, so we can suppose 
that $u$ is both a monomorphism and an epimorphism. By (PG1),
it follows that, up to composing with an isomorphism, 
$Y=X\amalg Z$ for some object $Z$ of $\cC$, and $u$ can
be identified with the natural morphism $X\to X\amalg Z$.
Let $v:Z\to X\amalg Z$ be the natural morphism; by (PT3)
the induced morphism $u\times_YZ:X\times_YZ\to Z$ is
an epimorphism and by (PT2) we have $X\times_YZ\simeq\emptyset$,
the initial object of $\cC$. Since in $\cC$ all epimorphisms
are effective, one derives that $Z\simeq\emptyset$,
and the claim follows.
\end{pfclaim}

Now, let $u:X\to Y$ be a morphism in $\cC$; the induced 
morphisms
$$
\pr_1,\pr_2:\xymatrix{X\times_YX\ar@<.5ex>[r]
\ar@<-.5ex>[r] & X}
$$
define an equivalence relation; by (PT3) there is a corresponding
quotient morphism $u':X\to Y'$ and moreover $u'$ is a strict
epimorphism. Clearly $u$ factors via a morphism $u'':Y'\to Y$.
We need to show that $u''$ is a monomorphism, or equivalently,
that the induced diagonal morphism $\delta:Y'\to Y'\times_YY'$ 
is an isomorphism. However, there is a natural commutative 
diagram
$$
\xymatrix{
X\times_{Y'}X \ar[r]^\alpha \ar[d] & X\times_YX \ar[d] \\ 
Y' \ar[r]^-\delta &  Y'\times_YY'}
$$
where $\alpha$ is an isomorphism by construction and both
vertical arrows are epimorphisms. It follows that $\delta$
is an epimorphism; since it is obviously a monomorphism
as well, we deduce (G3) in view of (PG1) and claim 
\ref{cl_mono.epi}. Let $\cC^\sim$ be the topos of sheaves
on the precanonical topology of $\cC$; by \eqref{subsec_pretopos} 
and Deligne's theorem \eqref{subsec_Deligne}, there exists
a fibre functor $\cC^\sim\to\Set$. Composing with the
natural functor $\cC\to\cC^\sim$ we obtain a functor 
$F:\cC\to\Set$ which is left exact, commutes with finite sums 
and quotients under equivalence relations, in view of 
\eqref{subsec_precanon}, so (G4) and (G5) hold for $F$.
\begin{claim}\label{cl_F.non-empty} 
Let $X$ be a non-empty object of $\cC$. Then $F(X)\neq\emptyset$.
\end{claim}
\begin{pfclaim} Since we know already that (G3) holds in
$\cC$, we deduce using (PG2) that the unique morphism $X\to e$
is an epimorphism. Then (PT3) says that $e$ can be written
as the quotient of $X$ under the induced equivalence relation
$\pr_1,\pr_2:\xymatrix{X\times_eX\ar@<.5ex>[r]\ar@<-.5ex>[r] & X}$.
Since $F$ commutes with quotients under such equivalence relations, 
the claim follows after remarking that $F(e)\neq\emptyset$, 
\end{pfclaim}

\begin{claim}\label{cl_reflect.epi} 
Let $u:X\to Y$ be a morphism in $\cC$ such that
$F(u)$ is surjective. Then $u$ is an epimorphism.
\end{claim}
\begin{pfclaim} We use (G3) to factor $u$ as an epimorphism
$u':X\to Y'$ followed by a monomorphism of the form 
$Y'\to Y'\amalg Z$. We need to show that $Z=\emptyset$ or
equivalently, in view of claim \ref{cl_F.non-empty}, that
$F(Z)=\emptyset$. However, the assumption implies that 
$F(Y')$ maps onto $F(Y'\amalg Z)$; on the other hand, 
$F$ commutes with finite sums, so the claim holds.
\end{pfclaim}

\begin{claim}\label{cl_reflect.mono} 
Let $u:X\to Y$ be a morphism in $\cC$ such that
$F(u)$ is injective. Then $u$ is a monomorphism.
\end{claim}
\begin{pfclaim} The assumption means that the induced
diagonal map $F(X)\to F(X)\times_{F(Y)}F(X)$ is bijective.
Then claim \ref{cl_reflect.epi} implies that the diagonal
morphism $X\to X\times_YX$ is an epimorphism. The latter is 
also obviously a monomorphism, hence an isomorphism, in view
of claim \ref{cl_mono.epi}; but this means that $u$ is a 
monomorphism.
\end{pfclaim}

Now, taking into account claims \ref{cl_mono.epi},
\ref{cl_reflect.epi} and \ref{cl_reflect.mono} we deduce
that (G6) holds for $F$.
It remains only to show that $F$ takes values in finite sets.
So, suppose by contradiction that $F(X)$ is an infinite
set for some object $X$ of $\cC$. We define inductively
a sequence of objects $(Y_i~|~i\in\N)$, with morphisms
$\phi_{i+1}:Y_{i+1}\to Y_i$ for every $i\in\N$, as follows. Let 
$Y_0:=e$, $Y_1:=X$ and $\phi_1$ the unique morphism.
Let then $i>0$ and suppose that $Y_i$ and $\phi_i$ have 
already been given. Using the diagonal morphism, we can write 
$Y_i\times_{Y_{i-1}}Y_i\simeq Y_i\amalg Z$ for some
object $Z$; we set $Y_{i+1}:=Z$ and let $\phi_{i+1}$ be the
restriction of the projection 
$\pr_1:Y_i\times_{Y_{i-1}}Y_i\to Y_i$. Notice that 
$F(Y_i)\neq\emptyset$ for every $i\in\N$ (indeed all the fibers
of the induced map $F(Y_i)\to F(Y_i)$ are infinite whenever $i>0$); 
in particular $Y_i\neq\emptyset$ for every $i\in\N$. 
On the other hand, for every $n>0$, $X\times Y_n$ admits a 
decomposition of the form 
$\displaystyle{(\mathop{\amalg}^nY_n)}\amalg Y_{n+1}$, which
is against (PG3). The contradiction concludes the proof of
the lemma.
\end{proof}

\sset\subsubsection{}\label{subsec_Fundgrp}
\index{Topos!constant, locally constant, bounded object 
in a|indref{subsec_Fundgrp}}
Let $E$ be a topos. Recall that an object $X$ of $E$ is
said to be {\em constant\/} if it is a direct sum of
copies of the final object $e$ of $E$ 
(\cite[Exp.IX, \S2.0]{SGA4-3}). The object $X$ is 
{\em locally constant\/} if there exists a covering
$(Y_i\to e~|~i\in I)$ of $e$, such that, for every $i\in I$, 
the restriction of $X\times Y_i$ is constant on the 
induced topos $E_{/Y_i}$ (\cite[Exp.III, \S5.1]{SGA4-1}).
If additionally there exists an integer $n$, such that
so that each $(X\times Y_i)_{|Y_i}$ is a direct sum of 
at most $n$ copies of $Y_i$, then we say that $X$ 
is {\em bounded}.

\begin{lemma}\label{lem_lcc} 
Let $E$ be a topos. Denote by $E_\mathrm{lcb}$
the full subcategory of all locally constant bounded
objects of $E$ (see \eqref{subsec_construct}). Then:
\begin{enumerate}
\item
$E_\mathrm{lcb}$ is a pretopos.
\item
If the final object of $e$ is connected and 
non-empty, $E_\mathrm{lcb}$ is a pregalois pretopos.
\end{enumerate}
\end{lemma}
\begin{proof} Let $X$ be an object of $E_\mathrm{lcb}$, and
$(Y_i\to e)$ a covering of the final object of $E$ by non-empty
objects, such that $X\times Y_i$ is constant on $E_{/Y_i}$ for 
every $i\in I$, say $(X\times Y_i)_{|Y_i}\simeq Y_i\times S_i$, 
where $S_i$ is some set. Since $X$ is bounded, there exists $n\geq 0$
such that the cardinality of every $S_i$ is no greater than
$n$. Since all finite limits are representable in $E$, axiom 
(PT1) can be checked locally on $E$, so we can reduce to the 
case of a finite inverse system $\underline{Z}:=(Z_j~|~j\in J)$ 
of constant bounded objects $Z_j:=e\times S_j$. Furthermore, 
thanks to \cite[Exp.IX, Lemme 2.1(i)]{SGA4-3} we can assume 
that the inverse system is induced from an inverse system 
$\underline{S}:=(S_j~|~j\in J)$ of maps of sets, in which case 
it is easy to check that 
$\liminv{J}\:\underline{Z}\simeq\liminv{J}\:\underline{S}$, which
implies (PT1).

(PT2) is immediate. Similarly, if $R$ is
a locally constant bounded equivalence relation on an
object $X$ of $E_\mathrm{lcb}$, then $X/R$ is again in 
$E_\mathrm{lcb}$; indeed, since equivalence relations in $E$ 
are universally effective, this can be checked locally on
a covering $(Y_i\to e~|~i\in I)$. Then again, by 
\cite[Exp.IX, Lemme 2.1(i)]{SGA4-3}
we can reduce to the case of an equivalence relation on
sets, where everything is obvious. This shows that (PT3)
holds, and proves (i). Suppose next that $e$ is connected
and non-empty; since $e$ is in $E_\mathrm{lcb}$, it follows 
that (PG2) holds in $E_\mathrm{lcb}$. To show (PG1), consider
a morphism $u:X\to Y$ in $E_\mathrm{lcb}$. As in the foregoing,
we can find a covering $(Z_i\to e~|~i\in I)$ such that
$(X\times Z_i)_{|Z_i}\simeq Z_i\times S_i$, 
$(Y\times Z_i)_{|Z_i}\simeq Z_i\times T_i$ for some sets
$S_i$, $T_i$, and $u_{|Z_i}$ is induced by a map $u_i:S_i\to T_i$.
Let $S'_i:=T_i\setminus u_i(S_i)$. Clearly we have an
isomorphism 
$(Y\times Z_i)_{|Z_i}\simeq(X\times Z_i)_{|Z_i}\amalg
(Z_i\times S'_i)$ for every $i\in I$. Since the induced
decompositions agree on $Z_i\times Z_j$ for every $i,j\in I$,
the constant objects $Z_i\times S'_i$ glue to a locally
constant object $X'$ of $E$, and $u$ induces an isomorphism
$Y\simeq X\amalg X'$. 
Finally, if $X$ is in $E_\mathrm{lcb}$, find sets 
$(S_i~|~i\in I)$ and a covering $(Y_i\to e~|~i\in I)$ such 
that $(X\times Y_i)_{|Y_i}\simeq S_i\times Y_i$, and let 
$m$ be the maximum of the cardinalities of the sets $S_i$.
Clearly (PG3) holds for $X$, if one chooses  $n:=m+1$.
\end{proof}

\sset\subsubsection{}\label{subsec_covers}
\index{Almost algebra(s)!$\Cov(\Spec\,A)$ : 
{\'e}tale covering of an|indref{subsec_covers}}
Let $A$ be an almost algebra. We consider the site 
$\cS_A:=(A\Alg)^o_\mathrm{fpqc}$ obtained by endowing the 
category of affine $A$-schemes with the fpqc topology (see 
\eqref{subsec_fpqc}). Moreover, the category of 
{\em {\'e}tale coverings\/} of $\Spec\,A$ is defined as the 
full subcategory $\Cov(\Spec\,A)$ of the category of
affine $A$-schemes consisting of all {\'e}tale $A$-schemes of
finite rank.

\begin{proposition}\label{prop_covers} 
The natural functor $\eps:\cS_A\to\cS_A^\sim$ induces an
equivalence of the category $\Cov(\Spec\,A)$ onto
the category of locally constant bounded sheaves on $\cS_A$. 
\end{proposition}
\begin{proof} Let $B$ be an {\'e}tale $A$-algebra of rank $r$.
By proposition \ref{prop_decomp.fin.rank} there is an isomorphism
of almost algebras $A\simeq\prod^r_{i=0}A_i$ such that 
$B_i:=B\otimes_AA_i$ is an $A_i$-algebra of constant rank $i$
for every $i=0,...,r$. In particular, $B_0=0$ and $B_i$ is
faithfully flat {\'e}tale $A_i$-algebra for every $i>0$. We use the 
diagonal morphism to obtain a decomposition 
$B_i\otimes_{A_i}B_i\simeq B_i\times C_i$, where $C_i$ is
again an {\'e}tale $B_i$-algebra of constant rank $i-1$. Iterating
this procedure we find faithfully flat $A_i$-algebras $D_i$
such that $B\otimes_AD_i$ is $D_i$-isomorphic to a direct
product of $i$ copies of $D_i$, for every $i>0$. Setting 
$D_0:=A_0$, we obtain a covering 
$(\eps(\Spec\,D_i\to\Spec\,A)~|~i=0,...,r)$
of $\eps(\Spec\,A)$ in $\cS_A^\sim$ such that the restriction of
$\eps(\Spec\,B)$ to each $\eps(\Spec\,D_i)$ is a bounded constant 
sheaf. This show that the restriction of $\eps$ to 
$\Cov(\Spec\,A)$ lands in the category of locally constant
bounded sheaves.
Since the fpqc topology is coarser than the canonical topology,
$\eps$ is fully faithful. To show that $\eps$ is essentially
surjective amounts to an exercise in faithfully flat descent :
clearly every constant sheaf is represented by an {\'e}tale 
$A$-scheme; then one uses \cite[Exp.IX, Lemme 2.1(i)]{SGA4-3} 
to show that any descent datum of bounded constant sheaves
is induced by a cocycle system of morphisms for the corresponding
representing algebras, and one can descend the latter. We leave
the details to the reader.
\end{proof}

\sset\subsubsection{}\label{subsec_connected}
\index{Almost algebra(s)!connected|indref{subsec_connected}}
Let us say that an affine almost scheme $S$ is {\em connected}
if $\eps(S)$ is a connected object of the
category $\cS_{\cO_S}^\sim$, which simply means that the only
non-zero idempotent of $\cO_{S*}$ is the identity. In this case, 
proposition \ref{prop_covers} and lemma \ref{lem_lcc}
show that $\Cov(S)$ is a pregalois pretopos,
hence it admits a fibre functor $F:\Cov(S)\to\fSet$
by lemma \ref{lem_pregalois}.

\begin{definition}\label{def_Fund-grp}
\index{Almost algebra(s)!connected!$\pi_1(\Spec\,A)$ : fundamental group 
of a|indref{def_Fund-grp}}
Suppose that the affine almost scheme $S$ is connected. The 
{\em fundamental group\/} of $S$ is the group $\pi_1(S)$ 
defined as the automorphism group of any fibre functor 
$F:\Cov(S)\to\fSet$, endowed with its natural 
profinite topology (see \eqref{subsec_Galois.grp}).
\end{definition}

\sset\subsubsection{}
It results from the general theory (\cite[Exp.V, \S5]{SGA1}) 
that $\pi_1(S)$ is independent (up to isomorphism) on
the choice of fibre functor. We refer to {\em loc.cit.}
for a general study of fundamental groups of Galois
categories. In essence, several of the standard results
for schemes admit adequate almost counterpart. We conclude
this section with a sample of such statements.

\begin{definition}\label{def_ic}
\index{Almost algebra(s)!integral closure of an|indref{def_ic}}
Let $A\subset B$ be a pair of $V^a$-algebras;
the {\em integral closure\/} of $A$ in $B$ is the subalgebra
$\ic(A,B):=\ic(A_*,B_*)^a$, where for a pair of rings $R\subset S$
we let $\ic(R,S)$ be the integral closure of $R$ in $S$.
\end{definition}

\begin{lemma}\label{lem_ic}
If $R\subset S$ is any pair of $V$-algebras, then
$\ic(R^a,S^a)=\ic(R,S)^a$.
\end{lemma}
\begin{proof} It suffices to show the following:
\begin{claim} Given a commutative diagram of $V$-algebras
$$
\xymatrix{ R \ar[r] \ar[d] & S \ar[d] \\
           R_1 \ar[r] & S_1
}
$$
whose vertical arrows are almost isomorphism, the induced map
$\ic(R,S)\to\ic(R_1,S_1)$ is an almost isomorphism.
\end{claim}
\begin{pfclaim} Clearly the kernel of the map is almost zero;
it remains to show that for every $b\in\ic(R_1,S_1)$ and $\eps\in\fm$,
the element $\eps b$ lifts to $\ic(R,S)$. By assumption we have a
relation $b^N+\sum_{i=1}^Na_ib^{N-i}=0$, with $a_i\in R_1$, so
$(\eps b)^N+\sum_{i=1}^N\eps^i a_i(\eps b)^{N-i}=0$. By lifting
$\eps^ia_i$ to some $\bar a_i\in R$, we deduce a monic polynomial
$P(T)$ over $R$ such that $P(\eps b)=0$, so if $\bar b$ is a lifting
of $\eps b$, we have $\fm\cdot P(\bar b)=0$. Since the restriction
$\fm S\to\fm S_1$ is surjective, we can choose $\bar b\in\fm S$,
so $\bar b\cdot P(\bar b)=0$.
\end{pfclaim}
\end{proof}

\begin{remark}\label{rem_ic} (i) If $A\subset B$ are $V^a$-algebras,
then $A$ is integrally closed in $B$ if and only if $A_*$ is integrally
closed in $B_*$. Indeed, by lemma \ref{lem_ic} we know that the
integral closure of $A_*$ in $B_*$ is almost isomorphic to $A_*$
and any such $V$-algebra must be contained in $A_*$.

(ii) If $(A_i\subset B_i~|~i\in I)$ is a (possibly infinite)
family of pairs of $V^a$-algebras, then $\prod_{i\in I}A_i$
is integrally closed in $\prod_{i\in I}B_i$ if and only if
$A_i$ is integrally closed in $B_i$ for every $i\in I$.
\end{remark}

The following proposition is an analogue of \cite[Prop.18.12.15]{EGA4}.

\begin{proposition}\label{prop_analogue}
Let $A\subset B$ be a pair of $V^a$-algebras such that $A=\ic(A,B)$.
\begin{enumerate}
\item
For any \'etale almost finite projective $A$-algebra $A_1$ of almost
finite rank we have $A_1=\ic(A_1,A_1\otimes_AB)$.
\item
Suppose that $A$ and $B$ are connected, and choose a fibre functor
$F$ for the category $\Cov(B)$. Then the functor $\Cov(A)\to\fSet$ :
$C\mapsto F(C\otimes_AB)$ is a fibre functor, and the induced group
homomorphism: $\pi_1(\Spec\,B)\to\pi_1(\Spec\,A)$ is surjective.
\end{enumerate}
\end{proposition}
\begin{proof} (i): using remark \ref{rem_ic}(ii) we reduce to the
case where $A_1$ has constant finite rank over $A$. Set
$B_1:=A_1\otimes_AB$, and suppose that $x\in B_{1*}$ is integral
over $A_{1*}$, or equivalently, over $A_*$. Consider the element
$e\in(A_1\otimes_AA_1)_*$ provided by proposition \ref{prop_idemp};
for given $\eps\in\fm$ write $\eps\cdot e=\sum^k_ic_i\otimes d_i$
for some $c_i,d_i\in A_{1*}$.
According to remark \ref{rem_zeta.identity} we have:
$\sum^k_ic_i\cdot\Tr_{B_1/B}(xd_i)=\eps\cdot x$ for every $x\in B_{1*}$.
Corollary \ref{cor_Cayley} implies that $\Tr_{B_1/B}(xd_i)$ is integral
over $A_*$ for every $i\leq k$, hence it lies in $A_*$, by remark
\ref{rem_ic}(i). Hence $\eps\cdot x\in A_{1*}$, as claimed.

(ii): is an immediate consequence of (i) and of the general
theory of \cite[Exp.V, \S5]{SGA1}.
\end{proof}

\sset\subsubsection{}\label{subsec_fund.grp}
As a special case, let $R$ be a $V$-algebra whose spectrum is
connected; we suppose that there exists an element $\eps\in\fm$
which is regular in $R$. Suppose moreover that $R$ is integrally 
closed in $R[\eps^{-1}]$, consequently $\Spec\,R[\eps^{-1}]$ 
is connected and $\pi_1(\Spec\,R[\eps^{-1}])$ is well defined. 
It follows as well that $\Spec\,R^a$ is connected. 
Indeed, if $\Spec\,R^a$ were not connected, neither would be 
$\Spec\,R^a_*$; but since $R^a_*\subset R[\eps^{-1}]$, this 
is absurd. Then $\pi_1(\Spec\,R^a)$ is also well defined and,
after a fibre functor for $\Cov(R[\eps^{-1}])$
is chosen, the functors
$\Cov(\Spec\,R)\to\Cov(\Spec\,R^a)\to\Cov(\Spec\,R[\eps^{-1}])$ : 
$B\mapsto B^a\mapsto B^a_*[\eps^{-1}]$ induce continuous group
homomorphisms (\cite[Exp.V, Cor.6.2]{SGA1})
\set\begin{equation}\label{eq_cont.gp.hom}
\pi_1(\Spec\,R[\eps^{-1}])\to\pi_1(\Spec\,R^a)\to\pi_1(\Spec\,R).
\end{equation}

\begin{proposition}\label{prop_fund.grp} 
Under the assumptions of \eqref{subsec_fund.grp}, we have:
\begin{enumerate}
\item
$R^a_*$ is integrally closed in $R[\eps^{-1}]$.
\item
The maps \eqref{eq_cont.gp.hom} are surjective.
\end{enumerate}
\end{proposition}
\begin{proof} (i): let $x\in R[\eps^{-1}]$ such that
$x^n+a_1\cdot x^{n-1}+...+a_n=0$ for some elements 
$a_1,...,a_n\in R_*^a$. For every $\delta\in\fm$ we have
$(\delta\cdot x)^n+
(\delta\cdot x)\cdot(\delta\cdot x)^{n-1}+...+\delta^n\cdot a_n=0$,
which shows that $\delta\cdot x\in R$, since by assumption
$R$ is integrally closed in $R[\eps^{-1}]$ and 
$\delta^i\cdot a_i\in R$ for every $i=1,...,n$. The 
assertion follows.

(ii): it suffices to show the assertion for the leftmost 
map and for the composition of the two maps. However, the
composition of the two maps is actually a special case of
the leftmost map (for the classical limit $V=\fm$), so we need
only consider the leftmost map. Then the assertion follows
from (i) and proposition \ref{prop_analogue}(ii).
\end{proof}

\printindex


\newpage

\begin{thebibliography}{99}

\bibitem[1]{Ana}
{\scshape S.Anantharaman},
Sch\'emas en groupes, espaces homog\`enes et espaces alg\'ebriques
sur une base de dimension 1.
{\it Bull. Soc. Math. France} 33 (1973) pp.5-79. 

\bibitem[2]{An}
{\scshape M.Andr{\'e}},
Homologie des alg{\`e}bres commutatives.
{\it Springer Grundl. Math. Wiss.} 206 (1974).

\bibitem[3]{AnII}
{\scshape M.Andr{\'e}},
Localisation de la lissit{\'e} formelle.
{\it Manuscr. Math.} 13 (1974) pp.297-307.

\bibitem[4]{SGA4-1}
{\scshape M.Artin et al.},
Th{\'e}orie des topos et cohomologie {\'e}tale des sch{\'e}mas - tome 1.
{\it Springer Lect. Notes Math.} 269 (1972).

\bibitem[5]{SGA4-2}
{\scshape M.Artin et al.},
Th{\'e}orie des topos et cohomologie {\'e}tale des sch{\'e}mas - tome 2.
{\it Springer Lect. Notes Math.} 270 (1972).

\bibitem[6]{SGA4-3}
{\scshape M.Artin et al.},
Th{\'e}orie des topos et cohomologie {\'e}tale des sch{\'e}mas - tome 3.
{\it Springer Lect. Notes Math.} 305 (1973).

\bibitem[7]{Bart}
{\scshape W.Bartenwerfer},
Die h{\"o}heren metrischen Kohomologiegruppen affinoider R{\"a}ume.
{\it Math. Ann.} 241 (1979) pp.11-34.

\bibitem[8]{Bass}
{\scshape H.Bass},
Algebraic $K$-theory.
{\it W.A. Benjamin} (1968).

\bibitem[9]{Be-La}
{\scshape A.Beauville, Y.Laszlo},
Un lemme de descente.
{\it C.R. Acad. Sc. Paris} 320 (1995) pp.335-340.

\bibitem[10]{Berk}
{\scshape V.Berkovich},
Spectral theory and analytic geometry over non-archimedean fields.
{\it AMS Math. Surveys and Monographs} 33 (1990).

\bibitem[11]{Bert}
{\scshape P.Berthelot et al.},
Th{\'e}orie des Intersection et Th{\'e}or{\`e}mes de Riemann-Roch.
{\it Springer Lect. Notes Math.} 225 (1971).

\bibitem[12]{Bo-Gun}
{\scshape S.Bosch, U.G{\"u}ntzer, R.Remmert},
Non-Archimedean analysis.
{\it Springer Grundl. Math. Wiss.} 261 (1984).

\bibitem[13]{Bosch}
{\scshape S.Bosch, W.L{\"u}tkebohmert},
Formal and rigid geometry I. Rigid spaces.
{\it Math. Ann.} 295 (1993) pp.291-317.

\bibitem[14]{BoschIII}
{\scshape S.Bosch, W.L{\"u}tkebohmert, M.Raynaud},
Formal and rigid geometry III. The relative maximum principle.
{\it Math. Ann.} 302 (1995) pp.1-29.

\bibitem[15]{Bourbaki}
{\scshape N.Bourbaki}
Alg{\`e}bre.
{\it Hermann} (1970).

\bibitem[16]{BouAC}
{\scshape N.Bourbaki},
Alg{\`e}bre Commutative.
{\it Hermann} (1961).

\bibitem[17]{BouAH}
{\scshape N.Bourbaki},
Alg{\`e}bre Homologique.
{\it Masson} (1980).

\bibitem[18]{Bou}
{\scshape N.Bourbaki},
Topologie G{\'e}n{\'e}rale.
{\it Hermann} (1971).

\bibitem[19]{Br-Ti}
{\scshape F.Bruhat, J.Tits},
Groupes r\'eductifs sur un corps local - Chapitre II. 
Sch\'emas en groupes. Existence d'une donn\'ee radicielle valu\'ee. 
{\it Publ. Math. IHES} 60 (1984) pp.1-194.

\bibitem[20]{Ca-Fr}
{\scshape J.W.S.Cassels, A.Fr{\"o}hlich},
Algebraic number theory.
{\it Academic Press} (1967).

\bibitem[21]{Coat}
{\scshape J.Coates, R.Greenberg},
Kummer theory for abelian varieties over local fields. 
{\it Invent. Math.} 124 (1996) pp.129--174. 

\bibitem[22]{Conr}
{\scshape B.Conrad},
Irreducible components of rigid spaces.
{\it Ann. Inst. Fourier} 49 (1999) pp.473-541.

\bibitem[23]{De}
{\scshape P.Deligne},
Cat{\'e}gories tannakiennes.
{\it Grothendieck Festchrift vol.II, Birkhauser Progress in Math.} 87
(1990) pp.111-195.

\bibitem[24]{DeM}
{\scshape P.Deligne, J.Milne},
Tannakian categories.
{\it Springer Lect. Notes Math.} 900 (1982) pp.101-228.

\bibitem[25]{SGA3}
{\scshape M.Demazure, A.Grothendieck et al.},
Sch{\'e}mas en Groupes I.
{\it Springer Lect. Notes Math.} 151 (1970).

\bibitem[26]{EGAI}
{\scshape J.Dieudonn{\'e}, A.Grothendieck},
{\'E}l{\'e}ments de G{\'e}om{\'e}trie Alg{\'e}brique - Chapitre I.
{\it Publ. Math. IHES} 4 (1960).

\bibitem[27]{EGA}
{\scshape J.Dieudonn{\'e}, A.Grothendieck},
{\'E}l{\'e}ments de G{\'e}om{\'e}trie Alg{\'e}brique - Chapitre II.
{\it Publ. Math. IHES} 8 (1961).

\bibitem[28]{EGAIV}
{\scshape J.Dieudonn{\'e}, A.Grothendieck},
{\'E}l{\'e}ments de G{\'e}om{\'e}trie Alg{\'e}brique - Chapitre IV, partie 1.
{\it Publ. Math. IHES} 20 (1964).

\bibitem[29]{EGAIV-3}
{\scshape J.Dieudonn{\'e}, A.Grothendieck},
{\'E}l{\'e}ments de G{\'e}om{\'e}trie Alg{\'e}brique - Chapitre IV, partie 3.
{\it Publ. Math. IHES} 28 (1966).

\bibitem[30]{EGA4}
{\scshape J.Dieudonn{\'e}, A.Grothendieck},
{\'E}l{\'e}ments de G{\'e}om{\'e}trie Alg{\'e}brique - 
Chapitre IV, partie 4.
{\it Publ. Math. IHES} 32 (1967).

\bibitem[31]{Elk}
{\scshape R.Elkik},
Solutions d'{\'e}quations {\`a} coefficients dans un anneau hens{\'e}lien.
{\it Ann. Sci. ENS} 6 (1973) pp.553-604.

\bibitem[32]{Ell}
{\scshape G.Elliott},
On totally ordered groups and $K_0$.
{\it Springer Lect. Notes Math.} 734 (1979) pp.1-49.

\bibitem[33]{Fa1}
{\scshape G.Faltings},
$p$-adic Hodge theory.
{\it J. Amer. Math. Soc.} 1 (1988) pp.255--299. 

\bibitem[34]{Fa2}
{\scshape G.Faltings},
Almost {\'e}tale extensions.
{\it Preprint Max-Planck-Institut f{\"u}r Mathematik} 3 (1998).

\bibitem[35]{Fer}
{\scshape D.Ferrand},
Descente de la platitude par un homomorphisme fini.
{\it C.R. Acad. Sc. Paris} 269 (1969) pp.946-949.

\bibitem[36]{Fe-Ra}
{\scshape D.Ferrand, M.Raynaud},
Fibres formelles d'un anneau local noeth{\'e}rien.
{\it Ann. Sci.E.N.S.} 3 (1970) pp.295--311.

\bibitem[37]{Fr-Ma}
{\scshape J.Fresnel, M.Matignon},
Produit tensoriel topologique de corps valu{\'e}s.
{\it Canadian J. Math.} 35 (1983) pp.218-273.

\bibitem[38]{Ga}
{\scshape P.Gabriel},
Des cat{\'e}gories ab{\'e}liennes.
{\it Bull. Soc. Math. France} 90 (1962) pp.323-449.

\bibitem[39]{Gi}
{\scshape J.Giraud},
Cohomologie non ab{\'e}lienne.
{\it Springer Grundl. Math. Wiss.} 179 (1971).

\bibitem[40]{SGA1}
{\scshape A.Grothendieck et al.},
Rev{\^e}tements {\'E}tales et Groupe Fondamental.
{\it Springer Lect. Notes Math.} 224 (1971).

\bibitem[41]{Gru}
{\scshape L.Gruson},
Dimension homologique des modules plats sur un anneau commutatif 
noetherien.
{\it Symposia Mathematica} Vol. XI;
Academic Press, London (1973) pp. 243--254.

\bibitem[42]{Gr-Ra}
{\scshape L.Gruson, M.Raynaud},
Crit{\`e}res de platitude et de projectivit{\'e}.
{\it Invent. Math.} 13 (1971) pp.1-89.

\bibitem[43]{Hak}
{\scshape M.Hakim},
Topos annel{\'e}s et sch{\'e}mas relatifs.
{\it Springer Ergebnisse Math. Grenz.} 64 (1972).

\bibitem[44]{Hu1}
{\scshape R.Huber},
Bewertungsspektrum und rigide Geometrie.
{\it Regensburger Math. Schriften} 23 (1993).

\bibitem[45]{Hu2}
{\scshape R.Huber},
{\'E}tale cohomology of rigid analytic varieties
and adic spaces.
{\it Vieweg Aspects of Math.} 30 (1996).

\bibitem[46]{Il}
{\scshape L.Illusie},
Complexe cotangent et d{\'e}formations I.
{\it Springer Lect. Notes Math.} 239 (1971).

\bibitem[47]{Il2}
{\scshape L.Illusie},
Complexe cotangent et d{\'e}formations II.
{\it Springer Lect. Notes Math.} 283 (1972).

\bibitem[48]{Ka}
{\scshape K.Kato},
Logarithmic structures of Fontaine-Illusie.
{\it Algebraic analysis, geometry, and number theory - 
Johns Hopkins Univ. Press} (1989) pp.191--224.

\bibitem[49]{Tors}
{\scshape G.Kempf et al.},
Toroidal embeddings I.
{\it Springer Lect. Notes Math.} 339 (1973).

\bibitem[50]{Lan}
{\scshape S.Lang},
Algebra - Third edition.
{\it Addison-Wesley} (1993).

\bibitem[51]{La}
{\scshape D.Lazard},
Autour de la platitude.
{\it Bull. Soc. Math. France} 97 (1969) pp.81-128.

\bibitem[52]{MLazard}
{\scshape M.Lazard},
Commutative formal groups.
{\it Springer Lect. Notes Math.} 443 (1975).

\bibitem[53]{Ma}
{\scshape S.Mac Lane},
Categories for the working mathematician.
{\it Springer Grad. Text Math.} 5 (1971).

\bibitem[54]{Mat}
{\scshape H.Matsumura},
Commutative ring theory.
{\it Cambridge Univ. Press} (1986).

\bibitem[55]{Mit}
{\scshape B.Mitchell},
Rings with several objects.
{\it Advances in Math.} 8 (1972) pp.1-161.

\bibitem[56]{Mumford}
{\scshape D.Mumford},
Abelian varieties.
{\it Oxford U.Press} (1970).

\bibitem[57]{Nee}
{\scshape A.Neeman},
A counterexample to a 1961 ``theorem'' in homological algebra.
{\it Invent. Math.} 148 (2002) pp.397-420.

\bibitem[58]{Oli}
{\scshape J.-P.Olivier},
Descente par morphismes purs.
{\it C.R. Acad. Sc. Paris} 271 (1970) pp.821-823. 

\bibitem[59]{Ray0}
{\scshape M.Raynaud},
Faisceaux amples sur les sch\'emas en groupes et les espaces homog\`enes.
{\it Springer Lect. Notes Math.} 119 (1970).

\bibitem[60]{Ray}
{\scshape M.Raynaud},
Anneaux locaux hens{\'e}liens.
{\it Springer Lect. Notes Math.} 169 (1970).

\bibitem[61]{Ro}
{\scshape N.Roby},
Lois polyn{\^o}mes et lois formelles en th{\'e}orie des modules.
{\it Ann. Sci.E.N.S.} 80 (1963) pp.213-348.

\bibitem[62]{Za-Sa}
{\scshape P.Samuel, O.Zariski},
Commutative algebra vol.II.
{\it Springer Grad. Text Math.} 29 (1975).

\bibitem[63]{Ser}
{\scshape J.-P.Serre},
Groupes de Grothendieck des sch\'emas en groupes r\'eductifs
d\'eployes.
{\it Publ. Math. IHES} 34 (1968) pp.37-52.

\bibitem[64]{Ta}
{\scshape J.Tate},
$p$-divisible groups.
{\it Proc. conf. local fields, Driebergen} (1967) pp.158-183.

\bibitem[65]{Put}
{\scshape M.van der Put},
Cohomology on affinoid spaces.
{\it Compositio. Math.}, 45 (1982) pp.165-198.

\bibitem[66]{Su}
{\scshape H.Sumihiro},
Equivariant completion II.
{\it J. Math. Kyoto}, 15 (1975) pp.573-605.

\bibitem[67]{We}
{\scshape C.Weibel},
An introduction to homological algebra.
{\it Cambridge Univ. Press} (1994).

\bibitem[68]{Zar}
{\scshape O.Zariski},
Reduction of the singularities of algebraic 3-dimensional
varieties.
{\it Ann. of Math.} (1944) pp.472-542.

\end{thebibliography}
\end{document}